%% file: Geometric_topology.tex
\documentclass [11pt]{amsbook} 

\usepackage{amssymb} \usepackage{amsfonts} \usepackage{amsmath}
\usepackage{amsthm} \usepackage{epsfig} 
\usepackage{mathrsfs} 
\usepackage{color}
\usepackage{amscd}
\usepackage[all]{xy}
\usepackage{graphicx}
\usepackage{etoolbox}

\usepackage[english]{babel}
\usepackage[utf8]{inputenc} 	
\usepackage[T1]{fontenc}		
\usepackage{mathtools}		

\usepackage{cmbright}	

\newcommand*{\transp}[2][-3mu]{\ensuremath{\mskip1mu\prescript{\smash{\mathrm t\mkern#1}}{}{\mathstrut#2}}}%

\newtheorem{lemma}{Lemma}[section]
\newtheorem{teo}[lemma]{Theorem}
\newtheorem{prop}[lemma]{Proposition}
\newtheorem{cor}[lemma]{Corollary} 
\newtheorem{conj}[lemma]{Conjecture}

\theoremstyle{definition}
\newtheorem{defn}[lemma]{Definition}

\newtheorem{example}[lemma]{Example}
\newtheorem{ex}[lemma]{Exercise} 

\theoremstyle{remark}
\newtheorem{oss}[lemma]{Remark} 
\newtheorem{rem}[lemma]{Remark} 
\newtheorem{warning}[lemma]{Warning} 

\numberwithin{section}{chapter}  
\numberwithin{subsection}{section}  

\numberwithin{figure}{chapter}
\numberwithin{table}{chapter}

\newcommand{\matK}{\ensuremath {\mathbb{K}}}
\newcommand{\matN}{\ensuremath {\mathbb{N}}}
\newcommand{\matR} {\ensuremath {\mathbb{R}}}
\newcommand{\matQ} {\ensuremath {\mathbb{Q}}}
\newcommand{\matZ} {\ensuremath {\mathbb{Z}}}
\newcommand{\matC} {\ensuremath {\mathbb{C}}}
\newcommand{\matP} {\ensuremath {\mathbb{P}}}
\newcommand{\matH} {\ensuremath {\mathbb{H}}}
\newcommand{\matS} {\ensuremath {\mathbb{S}}}

\newcommand{\matRP} {\ensuremath {\mathbb{RP}}}
\newcommand{\matCP} {\ensuremath {\mathbb{CP}}}
\newcommand{\matX}{\ensuremath {\mathbb{X}}}

\newcommand{\calI} {\ensuremath {\mathscr{I}}}
\newcommand{\calL} {\ensuremath {\mathcal{L}}}
\newcommand{\calS} {\ensuremath {\mathscr{S}}}
\newcommand{\calM} {\ensuremath {\mathscr{M}}}
\newcommand{\calC} {\ensuremath {\mathscr{C}}}
\newcommand{\calA} {\ensuremath {\mathcal{A}}}

\newcommand{\calF} {\ensuremath {\mathscr{F}}}

\newcommand{\calT} {\ensuremath {\mathcal{T}}}

\newcommand{\calG}{\ensuremath {\mathscr{G}}}
\newcommand{\calML}{\calM\calL}
\newcommand{\PML}{\matP\calM\calL}
\newcommand{\frakg}{\ensuremath {\mathfrak{g}}}
\newcommand{\frakh}{\ensuremath {\mathfrak{h}}}

\newcommand{\frakgl}{\ensuremath {\mathfrak{gl}}}
\newcommand{\fraksl}{\ensuremath {\mathfrak{sl}}}
\newcommand{\frakon}{\ensuremath {\mathfrak{o}}}
\newcommand{\frakso}{\ensuremath {\mathfrak{so}}}

\newcommand{\nota} [1] {\caption{\footnotesize{#1}}}
\newcommand{\matr} [4] {\big({\tiny\begin{array}{@{}c@{\ }c@{}} #1 & #2 \\ #3 & #4 \\ \end{array}} \big)}

\newcommand{\interior}[1]{{\rm int}(#1)}
\newfont{\Got}{eufm10 scaled 1200}

\newcommand{\GL}{{\rm GL}}
\newcommand{\SL}{{\rm SL}}
\newcommand{\SO}{{\rm SO}}
\newcommand{\SU}{{\rm SU}}
\newcommand{\On}{{\rm O}}
\newcommand{\Un}{{\rm U}}

\newcommand{\Fix}{{\rm Fix}}
\newcommand{\Span}{{\rm Span}}

\newcommand{\st}{{\rm st}}

\newcommand{\img}{{\rm Im\,}}
\newcommand{\Img}{{\rm Im\,}}

\newcommand{\id}{{\rm id}}
\newcommand{\rk}{{\rm rk}}
\newcommand{\PSLC}{\ensuremath{\matP{\rm SL}_2(\matC)}}
\newcommand{\PGLC}{\ensuremath{\matP{\rm GL}_2(\matC)}}
\newcommand{\GLC}{\ensuremath{{\rm GL}_2(\matC)}}

\newcommand{\SLC}{\ensuremath{{\rm SL}_2(\matC)}}
\newcommand{\SLR}{\ensuremath{{\rm SL}_2(\matR)}}
\newcommand{\SLZ}{\ensuremath{{\rm SL}_2(\matZ)}}
\newcommand{\PSLR}{\ensuremath{\matP{\rm SL}_2(\matR)}}

\newcommand{\PSLZ}{\ensuremath{\matP{\rm SL}_2(\matZ)}}

\newcommand{\FN}{{\rm FN}}
\newcommand{\tr}{{\rm tr}}
\newcommand{\ML}{\ensuremath{\calM\calL}}
\newcommand{\Conf}{{\rm Conf}}
\newcommand{\sh}{{\rm sh}}
\newcommand{\Rot}{{\rm Rot}}
\newcommand{\Def}{{\rm Def}}
\newcommand{\Stab}{{\rm Stab}}

\newcommand{\dimo}[1]{\vspace{2pt}\noindent\textit{Proof of \ref{#1}}.\ }
\newcommand{\finedimo}{{\hfill\hbox{$\square$}\vspace{2pt}}}
\newcommand{\cref}{c^{\rm ref}}
\newcommand{\isom}{\cong}
\newcommand{\Isom}{\cong}
\newcommand{\Hom}{{\rm Hom}}
\newcommand{\Vol}{{\rm Vol}}
\newcommand{\Area}{{\rm Area}}
\newcommand{\Iso}{{\rm Isom}}

\newcommand{\Isomet}{{\rm Isom}}
\newcommand{\Aut}{{\rm Aut}}
\newcommand{\Int}{{\rm Int}}
\newcommand{\Out}{{\rm Out}}
\newcommand{\Omeo}{{\rm Homeo}}
\newcommand{\Diffeo}{{\rm Diffeo}}
\newcommand{\inj}{{\rm inj}}
\newcommand{\Teich}{{\rm Teich}}
\newcommand{\MCG}{{\rm MCG}}

\newcommand{\Nil}{{\rm Nil}}
\newcommand{\Sol}{{\rm Sol}}
\newcommand{\St}{{\rm St}}

\renewcommand{\emptyset}{\varnothing}

\newcommand{\timtil}{
\begin{picture}(16,8)
\put(4,0){$\times$}
\put(5,5){$\scriptstyle\sim$}
\end{picture}
}

\newcommand{\timforsetil}{
\begin{picture}(16,8)
\put(4,0){$\times$}
\put(2,5){${\scriptscriptstyle(}{\scriptstyle\sim}{\scriptscriptstyle)}$}
\end{picture}
}

\newtoggle{BW}
\togglefalse{BW}

\newtoggle{printed}
\toggletrue{printed}

\makeatletter

\def\l@subsection{\@tocline{2}{0pt}{2.5pc}{5pc}{}}

\makeatother

\makeindex

\begin{document}

\frontmatter

\author{Bruno Martelli \\ \vspace{2 cm} \iftoggle{printed}{\vspace{10 cm}}{\includegraphics[width = 7 cm]{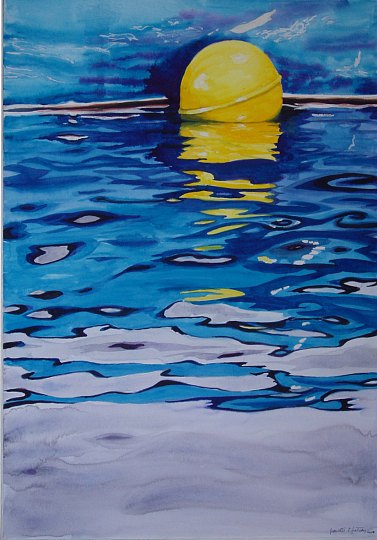} \vspace{3 cm}}}
\title{\iftoggle{printed}{\vspace{3 cm}}{} An Introduction to Geometric Topology}

\maketitle

\begin{center}
Bruno Martelli \\
University of Pisa, Italy \\

\vspace{1 cm}
{\tt http://people.dm.unipi.it/martelli/}

\vspace{5 cm}

Version 2.0 of April 2022 \\

\vspace{1 cm}

\end{center}

\newpage
\thispagestyle{empty}
    \null\vspace{\stretch {1}}
        \begin{flushright}
                \emph{To Ada, Sara, and Ylenia.}
        \end{flushright}
\vspace{\stretch{2}}\null

\tableofcontents

\mainmatter

\chapter*{Introduction}
The aim of this book is to introduce the reader to an area of mathematics called \emph{geometric topology}. The text should be suitable to a master or PhD student in mathematics interested in geometry, and more generally to any curious mathematician with a standard background in topology and analysis.

We interpret here the term ``geometric topology'' in a quite restrictive sense: for us, this topic is the study of topological manifolds via the assignment of a preferred ``geometric structure,'' that is some Riemannian metric of a particular nice kind. It is one of the most important discoveries in modern geometry that such preferred metrics exist on every compact manifold of dimension two and three, and this is exactly what this book is about. 

In other words, this book is an introduction to surfaces and three-manifolds, and to their \emph{geometrisation}, due to Poincar\'e and Koebe in 1907 in dimension two and to Thurston and Perelmann in 2002 in dimension three. Therefore this is also a textbook on \emph{low-dimensional topology}, except that we completely neglect four-manifolds, that form a relevant part of this area but which do not (yet?) fit in any geometrisation perspective.

There are already many good textbooks on surfaces, so our main new contribution is probably to furnish a complete introduction to Thurston's geometrisation of three-manifolds, that includes both the topological side of the story (the decomposition of three-manifolds along spheres and tori, the classification of Seifert manifolds) and the geometric side (hyperbolic geometry, Thurston's Dehn filling Theorem, Mostow rigidity, the eight geometries). 

This book is essentially self-contained: in the first chapter we recall all the relevant background material in differential topology and geometry, and from the second chapter on we never mention a theorem without furnishing a proof -- with only one important but unavoidable exception: Perelman's solution of the geometrisation conjecture, stated in Chapter \ref{eight:chapter}. 

\subsection*{Outline}
Here is a brief outline of the material contained in each chapter.
The book is divided into three parts. In the first, we introduce hyperbolic geometry and its relatives, the elliptic and flat geometries. Hyperbolic geometry is by far the richest, the most beautiful, and also the most important geometry in dimensions two and three: its relevance is witnessed by the folk sentence that a ``generic'' manifold of dimension two or three has a hyperbolic structure, that is a Riemannian structure locally isometric to the hyperbolic plane or space.

Chapter \ref{preliminaries:chapter} includes all the preliminaries that we will need: in particular, we quickly review various general notions of differential topology and geometry that can be found in many excellent books.  

Chapter \ref{space:chapter} introduces the reader to the hyperbolic space $\matH^n$: as opposite to the sphere $S^n$ and to the euclidean space $\matR^n$, the hyperbolic space may be represented in at least four different ways, none of which is prevalent: each representation is called a \emph{model}. We introduce the models, and then study the isometries, the compactification of $\matH^n$, and its subspaces. 

Chapter \ref{manifold:chapter} deals with hyperbolic manifolds. Maybe unexpectedly, these objects are related to a number of different beautiful mathematical concepts like discrete subgroups of Lie groups, polyhedra, and tessellations. We also present some generalisations like manifolds with geodesic boundary, cone manifolds, and orbifolds.

Chapter \ref{thick-thin:chapter} introduces the \emph{thick-thin decomposition}, a general structural theorem on hyperbolic manifolds that allows us to understand in particular the topology of the complete hyperbolic manifolds that have finite volume but are not compact. It also contains various information on flat and elliptic manifolds, including Bieberbach's Theorem.

Chapter \ref{infinity:chapter} introduces a variety of notions and results on hyperbolic manifolds that are connected with the points at infinity of $\matH^n$. 

The second part of the book deals with surfaces. We classify and geometrise every surface of finite type in Chapter \ref{surfaces:chapter}, where we also study simple closed curves in surfaces and define the \emph{mapping class group}.

Chapter \ref{Teichmuller:chapter} presents the \emph{Teichm\"uller space} of a genus-$g$ surface, as the space of its hyperbolic structures. We use the Fenchel--Nielsen coordinates to show that this space is in fact homeomorphic to $\matR^{6g-6}$.

Chapter \ref{automorphisms:chapter} introduces Thurston's beautiful theory of diffeomorphisms of surfaces. We introduce Thurston's compactification of the Teichm\"uller space, and the consequent classification of the elements of the mapping class group into three classes. We define and study some intriguing objects called \emph{geodesic currents} and \emph{laminations}.

The last (and longest) part of the book is devoted to three-manifolds.  In Chapter \ref{Three:topology:chapter} we introduce the first topological facts, including the prime decomposition, incompressible surfaces, and Haken manifolds.

Chapter \ref{Seifert:chapter} is entirely devoted to Seifert manifolds, a class of three-manifolds that contains many interesting examples. We classify these manifolds completely.

In Chapter \ref{construction:chapter} we present various techniques that topologists use every day to construct three-manifolds: Heegaard splittings, knots and links, Dehn surgery, and surface bundles. We end by stating and proving the canonical torus decomposition.

In Chapter \ref{eight:chapter} we finally move from topology to geometry: there are eight relevant geometries in dimension three, and we introduce them with some detail. We show in particular that the Seifert manifolds occupy precisely six of them.

In Chapter \ref{Mostow:chapter} it is due time to start investigating the most interesting of the three-dimensional geometries, the hyperbolic one: we prove that in dimension three every manifold has at most one hyperbolic structure, and this important fact is known as \emph{Mostow's rigidity Theorem}.

In Chapter \ref{three:chapter} we construct many examples of hyperbolic three-manifolds, by introducing ideal triangulations and \emph{Thurston's equations}.

Chapter \ref{Hyperbolic:Dehn:chapter} contains a complete proof of Thurston's \emph{hyperbolic Dehn filling theorem} and a discussion on the volumes of hyperbolic three-manifolds.

As we mentioned above, this book is almost entirely self-contained and the bibliography is minimised to the strict necessary: each chapter ends with a short section containing the pertaining references, that consist essentially in the sources that we have consulted for that chapter. Many of the topics presented here have their origin in Thurston's notes and papers and are of course already covered by other books, that we have widely used, so our bibliography consists mainly of secondary sources.

\subsection*{Acknowledgements}
This book has been written in the years 2011-16, as a slowly growing textbook that has accompanied the author and his students during the master courses in Pisa on hyperbolic geometry, surfaces, and three-manifolds. The author has largely profited of many discussions with his students and colleagues, and in particular he would like to warmly thank Giovanni Alberti, Vinicius Ambrosi, Ludovico Battista, Riccardo Benedetti, Andrea Bianchi, Francesco Bonsante, Alessio Carrega, Jacopo Guoyi Chen, Gemma Di Petrillo, Xiaoming Du, Leonardo Ferrari, Irene Filoscia, Stefano Francaviglia, Roberto Frigerio, John Hubbard, Wolfgang Keller, Patrick Lin, Filippo Mazzoli, Alice Merz, Gabriele Mondello, Andrea Monti, Matteo Novaga, Carlo Petronio, Nicola Picenni, Stefano Riolo, Federico Salmoiraghi, Leone Slavich, Chaitanya Tappu, Giacomo Tendas, Dylan Thurston, and Gabriele Viaggi for their help.

\subsection*{Copyright notices}
The text is released under the  \emph{Creative Commons}-BY-NC-SA license. 
You are allowed to distribute, modify, adapt, and use of this work for non-commercial purposes, as long as you correctly attribute its authorship and keep using the same license. 
 
The pictures used here are all in the public domain (both those created by the author and those downloaded from Wikipedia, which were already in the public domain), except the following ones that are released with a CC-BY-SA license and can be downloaded from Wikipedia: 
\begin{itemize}
\item Figure \ref{Connected_sum:fig} created by Oleg Alexandrov,
\item Figure \ref{Torus-triang:fig} (triangulated torus) created by Ag2gaeh,
\item Figure \ref{square_tiling:fig} (tessellation of hyperbolic plane), Figure \ref{Horocycle_normals:fig} (horosphere), Figure \ref{pseudosphere:fig} (pseudosphere), and Figure \ref{triangular2:fig} (tessellations of hyperbolic plane) created by Claudio Rocchini,
\item Figure \ref{icosaedro:fig}-right (dodecahedron) created by DTR,
\item Figure \ref{honeycomb_dodecahedron_534:fig} and \ref{honeycomb_dodecahedron_536:fig} (dodecahedral tessellations) created by Roice3,
\item Figure \ref{Voronoi:fig} (Voronoi tessellation) created by Mysid and Cyp,
\item Figure \ref{modular_group:fig} (fundamental domain of the modular group) and Figure \ref{Seifert_surface:fig} (Seifert surface for the trefoil knot) created by Kilom691,
\item Figure \ref{Klein_bottle:fig} created by Tttrung,
\item Figure \ref{lamination:fig} created by Adam Majewski,
\item Figure \ref{Hopf:fig} (The Hopf fibration) created by Niles Johnson,
\item Figure \ref{sum_of_knots:fig} (Knot sum) originally created by Maksim and tPbrocks13,
\item Figure \ref{satellites:fig} (Satellite knots) originally created by RyBu.
\end{itemize}
The painting in the front page is a courtesy of Mariette Michelle Egreteau. The Figures
\ref{horoballs:fig}, \ref{convex-hull:fig}, \ref{hyperboloid:fig} and \ref{3ps:fig} were made by Jacopo Guoyi Chen.

\part{Hyperbolic geometry}
\include{Preliminaries}

\include{Space}

\include{Manifold}

\include{Thick-thin}

\include{Infinity}

\part{Surfaces}
\include{Surfaces}

\include{Teichmuller}

\include{Automorphisms}

\part{Three-manifolds}
\include{Three_topology}

\include{Seifert}
\include{Construction}

\include{Eight}

\include{Mostow}
\include{Three}

\include{Hyperbolic_Dehn}

\include{Biblio}
\printindex

\end{document}

%% file: Preliminaries.tex
\chapter{Preliminaries} \label{richiami:chapter} \label{preliminaries:chapter}
We expect the reader to be familiar with the mathematics usually taught to undergraduates, like multivariable differential calculus, group theory, topological spaces, and fundamental groups. 

Some knowledge on differentiable and Riemannian manifolds would also help, at least on an intuitive way: in any case, this chapter introduces from scratch everything we need from differential topology and geometry. We also include some important information on groups (like Lie groups and group actions), and a few basic notions of measure theory and algebraic topology, with a very quick overview of homology theory. Finally, we use these tools to introduce cell complexes, handle decompositions, and triangulations: these are the main practical instruments that we have to build manifolds concretely.

Most results stated here are given without a proof: details can be found in various excellent books, some of which will be cited below.

\section{Differential topology} \label{differential:topology:section}
We introduce manifolds, bundles, embeddings, tubular neighbourhoods, isotopies, and connected sums. Throughout this book, we will always work in the smooth category. 
The material contained in this section is carefully explained in \emph{Differentiable manifolds} of Kosinksi \cite{K}.

\subsection{Differentiable manifolds}
A \emph{topological manifold} of dimension $n$ is a paracompact Hausdorff topological space $M$ locally homeomorphic to $\matR^n$. In other words, there is a covering $\{U_i\}$ of $M$ consisting of open sets $U_i$ homeomorphic to open sets $V_i$ in $\matR^n$.
\index{manifold}

\begin{figure}
\begin{center}
\includegraphics[width = 9 cm] {\iftoggle{BW}{test-BW}{test}}
\nota{The \emph{Alexander horned sphere} is a subset of $\matR^3$ homeomorphic to the 2-sphere $S^2$. It divides $\matR^3$ into two connected components, none of which is homeomorphic to an open ball. It was constructed by Alexander as a counterexample to a natural three-dimensional generalisation of Jordan's curve theorem. The natural generalisation would be the following: does every 2-sphere in $\matR^3$ bound a ball? If the 2-sphere is only topological, the answer is negative
as this counterexample shows. If the sphere is a differentiable submanifold, the answer is however positive as proved by Alexander himself.}
\label{Alexander:fig}
\end{center}
\end{figure}

Topological manifolds are difficult to investigate, their definition is too general and allows to directly define and prove only few things. Even the notion of dimension is non-trivial: to prove that an open set of $\matR^k$ is not homeomorphic to an open set of $\matR^h$ for different $k$ and $h$ we need to use non-trivial constructions like homology. It is also difficult to treat topological subspaces: for instance, the \emph{Alexander horned sphere} shown in Figure \ref{Alexander:fig} is a subspace of $\matR^3$ topologically homeomorphic to a 2-sphere. It is a complicated object that has many points that are not ``smooth'' and that cannot be ``smoothened'' in any reasonable way.

We need to define some ``smoother'' objects, and for that purpose we can luckily invoke the powerful multivariable infinitesimal calculus. Let $U\subset \matR^n$ be an open set: a map $f\colon U\to \matR^k$ is  \emph{smooth} if it is $C^\infty$, \emph{i.e.}~it has partial derivatives of any order.

\begin{defn} Let $M$ be a topological manifold.
A \emph{chart} is a fixed homeomorphism $\varphi_i\colon U_i \to V_i$ between an open set $U_i$ of $M$ and an open set $V_i$ of $\matR^n$. An \emph{atlas} is a set of charts $\big\{(U_i, \varphi_i)\big\}$ where the open sets $U_i$ cover $M$. 

If $U_i\cap U_j\neq\varnothing$ there is a \emph{transition map} $\varphi_{ji} = \varphi_j\circ \varphi_i^{-1}$ that sends homeomorphically the open set $\varphi_i({U_i\cap U_j})$ onto the open set $\varphi_j({U_i\cap U_j})$. Since these two open sets are in $\matR^n$, it makes sense to require $\varphi_{ij}$ to be smooth. The atlas is \emph{differentiable} if all transition maps are smooth. 
\end{defn}

\index{manifold!differentiable manifold}
\begin{defn}
A \emph{differentiable manifold} is a topological manifold that is equipped with a differentiable atlas.
\end{defn}

We will often use the word \emph{manifold} to indicate a differentiable manifold.
The integer $n$ is the \emph{dimension} of the manifold. We have defined the objects, so we now turn to their morphisms. 

\begin{defn}
A map $f\colon M \to M'$ between differentiable manifolds is \emph{smooth} if it is smooth when read locally through charts. This means that for every $p \in M$ there are two charts $(U_i, \varphi_i)$ of $M$ and $(U_j', \varphi_j')$ of $N$ with $p\in U_i$ and $f(p)\in U_j'$ such that the composition
$\varphi_j'\circ f \circ \varphi_i^{-1}$ is smooth wherever it is defined. 
\end{defn}

A \emph{diffeomorphism} is a smooth map $f\colon M \to M'$ that admits a smooth inverse $g\colon M'\to M$. A \emph{curve} in $M$ is a smooth map $\gamma \colon I \to M$ defined on some open interval $I$ of the real line, which may be bounded or unbounded. 
\index{diffeomorphism}

\begin{defn} A differentiable manifold is \emph{oriented} if it is equipped with an atlas where all transition functions are orientation-preserving (that is, the determinant of their differential at any point is positive).
\end{defn}

A manifold which can be oriented is called \emph{orientable}.
\index{manifold!oriented manifold}

\subsection{Tangent space} \label{tangente:subsection}
Let $M$ be a differentiable manifold of dimension $n$. We may define for every point $p\in M$ a $n$-dimensional vector space $T_pM$ called the \emph{tangent space}. 

The space $T_pM$ may be defined briefly as the set of all curves $\gamma\colon (-a,a) \to M$ such that $\gamma(0) = p$ and $a>0$ is arbitrary, considered up to some equivalence relation. The relation is the following: we identify two curves that, read on some chart $(U_i,\varphi_i)$, have the same tangent vector at $\varphi_i(p)$. The definition does not depend on the chart chosen.
\index{tangent space}

\begin{figure}
\begin{center}
\includegraphics[width = 6 cm] {\iftoggle{BW}{Tangentialvektor_new-BW}{Tangentialvektor_new}}
\nota{The tangent space in $x$ may be defined as the set of all curves $\gamma$ with $\gamma(0)=x$ seen up to an equivalent relation that identifies two curves having (in some chart) the same tangent vector at $x$. This condition is chart-independent.}
\label{Tangentialvektor:fig}
\end{center}
\end{figure}

A chart identifies $T_pM$ with the usual tangent space $\matR^n$ at $\varphi_i(p)$ in the open set $V_i = \varphi_i(U_i)$. Two distinct charts $\varphi_i$ and $\varphi_j$ provide different identifications with $\matR^n$, which differ by a linear isomorphism: the differential $d\varphi_{ji}$ of the transition map $\varphi_{ij}$. The structure of $T_pM$ as a vector space is then well-defined, while its identification with $\matR^n$ is not.

Every smooth map $f\colon M \to N$ between differentiable manifolds induces at each point $p\in M$ a linear map $df_p\colon T_pM \to T_{f(p)} N$ between tangent spaces in the following simple way: the curve $\gamma$ is sent to the curve $f\circ\gamma$.

\begin{defn} A smooth map $f\colon M \to N$ is a \emph{local diffeomorphism} at a point $p \in M$ if there are two open sets $U\subset M$ and $V\subset N$ containing respectively $p$ and $f(p)$ such that $f|_U\colon U \to V$ is a diffeomorphism.\index{diffeomorphism!local diffeomorphism}
\end{defn}

The inverse function theorem in $\matR^n$ implies easily the following fact, that demonstrates the importance of tangent spaces.

\begin{teo}
Let $f\colon M \to N$ be a smooth map between manifolds of the same dimension. The map is a local diffeomorphism at $p\in M$ if and only if the differential $df_p\colon T_pM \to T_{f(p)}N$ is invertible.
\end{teo}

In this theorem a condition satisfied at a single point (differential invertible at $p$) implies a local property (local diffeomorphism). Later on, we will see that in Riemannian geometry a condition satisfied at a single point may also imply a global property.

If $\gamma\colon I \to M$ is a curve, its \emph{velocity} $\gamma'(t)$ in $t\in I$ is the tangent vector $\gamma'(t) = d\gamma_t(1)$. Here ``1'' means the vector 1 in the tangent space $T_tI = \matR$. We note that the velocity is a vector and not a number: the modulus of a tangent vector is not defined in a differentiable manifold (because the tangent space is just a real vector space, without a norm).

\subsection{Differentiable submanifolds}
Let $N$ be a differentiable manifold of dimension $n$. 
\begin{defn}
A subset $M\subset N$ is a \emph{differentiable submanifold} of dimension $m\leqslant n$ if every $p\in M$ has an open neighbourhood $U\subset N$ diffeomorphic to an open set $V\subset \matR^n$ via a map $\varphi\colon U \to V$ that sends $U\cap M$ onto $V\cap L$, where $L$ is a linear subspace of dimension $m$.
\end{defn}
\index{submanifold!differentiable submanifold}

The pairs $\{(U\cap M, \varphi|_{U\cap M})\}$ form an atlas for $M$, which inherits a structure of $m$-dimensional differentiable manifold. At every point $p\in M$ the tangent space $T_p M$ is a linear subspace of $T_p N$. 

\subsection{Fibre bundles}
We introduce an important type of smooth maps. 
\begin{defn}
A \emph{smooth fibre bundle} is a smooth map
$$\pi\colon E \longrightarrow B$$
such that every fibre $\pi^{-1}(p)$ is diffeomorphic to a fixed manifold $F$ and $\pi$ looks locally like a projection. This means that $B$ is covered by open sets $U_i$ equipped with diffeomorphisms $\psi_i \colon U_i\times F \to \pi^{-1}(U_i) $ such that $\pi\circ\psi_i$ is the projection on the first factor.
\end{defn}
\index{fibre bundle}

The manifolds $E$ and $B$ are called the \emph{total} and \emph{base manifold}, respectively. The manifold $F$ is the \emph{fibre} of the bundle. A \emph{section} of the bundle is a smooth map $s\colon B\to E$ such that $\pi\circ s = \id_B$. Two fibre bundles $\pi\colon E \to B$ and $\pi'\colon E' \to B$ are \emph{isomorphic} if there is a diffeomorphism $\psi\colon E \to E'$ such that $\pi = \pi'\circ\psi$.

\subsection{Vector bundles}
A \emph{smooth vector bundle} is a smooth fibre bundle where every fibre $\pi^{-1}(p)$ has the structure of a $n$-dimensional real vector space which varies smoothly with $p$. Formally, we require that $F=\matR^n$ and $\psi_i (p,\cdot)\colon F \to \pi^{-1}(p) $ be an isomorphism of vector spaces for every $\psi_i$ as above.
\index{vector bundle}

The \emph{zero-section} of a smooth vector bundle is the section $s\colon B \to E$ that sends $p$ to $s(p)=0$, the zero in the vector space $\pi^{-1}(p)$. The image $s(B)$ of the zero-section is typically identified with $B$ via $s$.

Two vector bundles are \emph{isomorphic} if there is a diffeomorphism $\psi$ as above, which restricts to an isomorphism of vector spaces on each fibre.
As every manifold here is differentiable, likewise every bundle will be smooth and we will hence often omit this word.

\subsection{Tangent and normal bundle}
Let $M$ be a differentiable manifold of dimension $n$. The union of all tangent spaces 
$$TM = \bigcup_{p\in M} T_pM$$ 
is naturally a differentiable manifold of double dimension $2n$, called the \emph{tangent bundle}. The tangent bundle $TM$ is naturally a vector bundle over $M$, the fibre over $p\in M$ being the tangent space $T_pM$. 

\index{normal space}
Let $M\subset N$ be a smooth submanifold of $N$. The \emph{normal space} at a point $p\in M$ is the quotient vector space $\nu_pM = T_pN/_{T_pM}$. The \emph{normal bundle} $\nu M$ is the union
$$\nu M = \bigcup_{p\in M} \nu_p M$$ 
and is also naturally a smooth vector bundle over $M$. The normal bundle is not  canonically contained in $TN$ like the tangent bundle, but (even more usefully) it may be embedded directly in $N$, as we will soon see.

\subsection{Vector fields}
A \emph{vector field} $X$ on a smooth manifold $M$ is a section of its tangent bundle. A point $p$ and a vector $v\in T_pM$ determine an \emph{integral curve} $\alpha\colon I\to M$ starting from $v$, that is a smooth curve with $\alpha(0) = p, \alpha'(0) = v$, and $\alpha'(t) = X(\alpha(t))$ for all $t\in I$. The curve $\alpha$ is unique if we require the interval $I$ to be maximal. It depends smoothly on $p$ and $v$.
\index{vector field}

\subsection{Immersions and embeddings}
A smooth map $f\colon M \to N$ between manifolds is an \emph{immersion} if its differential is everywhere injective: note that this does not imply that $f$ is injective. 
The map is an \emph{embedding} if it is a diffeomorphism onto its image: this means that $f$ is injective, its image is a submanifold, and $f\colon M \to f(M)$ is a diffeomorphism. 
\index{immersion}
\index{embedding}

\begin{teo} If $M$ is compact, every injective immersion $f\colon M \to N$ is an embedding.
\end{teo}

\subsection{Isotopy and ambient isotopy} \label{isotopy:subsection}
Let $X$ and $Y$ be topological spaces. We recall that a \emph{homotopy} between two continuous maps $\varphi, \psi\colon X \to Y$ is a continuous map $F\colon X\times [0,1]\to Y$  such that $F_0 = \varphi$ and $F_1 = \psi$, where $F_t = F(\cdot, t)$. 

\index{isotopy}
Let $M$ and $N$ be differentiable manifolds. A \emph{smooth isotopy} between two embeddings $\varphi, \psi\colon M \to N$ is a smooth homotopy $F$ between them, such that every map $F_t$ is an embedding. Again, we will shortly use the word \emph{isotopy} to mean a smooth isotopy. We note that many authors do not require an isotopy to be smooth, but we do.

Being isotopic is an equivalence relation for smooth maps $M\to N$: two isotopies $F_t$ and $G_u$ can be glued if $F_1 = G_0$, and since we want a smooth map we priorly modify $F$ and $G$ so that they are constant near $t=1$ and $u=0$. This can be done easily by reparametrising $t$ and $u$.

An \emph{ambient isotopy} on $N$ is an isotopy between $\id_N$ and some other diffeomorphism $\varphi\colon N \to N$, such that every level is a diffeomorphism. Two embeddings $\varphi, \psi\colon M \to N$
are \emph{ambiently isotopic} if there is an ambient isotopy $F$ of $N$ such that $\psi = F_1\circ \varphi$. 

\begin{teo} 
If $M$ is compact, two embeddings $\varphi, \psi\colon M \to N$ are isotopic if and only if they are ambiently isotopic.
\end{teo}

\subsection{Tubular neighbourhood}
Let $M\subset N$ be a differentiable submanifold. A \emph{tubular neighbourhood} of $M$ is an open subset $U\subset N$ with a diffeomorphism $\nu M \to U$ sending identically the zero-section onto $M$. 
\index{tubular neighbourhood}

\begin{teo} Let $M\subset N$ be a closed differentiable submanifold. A tubular neighbourhood for $M$ exists and is unique up to an isotopy fixing $M$ and up to pre-composing with a bundle isomorphism of $\nu M$.
\end{teo}

If we are only interested in the open set $U$ and not its parametrisation, we can of course disregard the bundle isomorphisms of $\nu M$.

Vector bundles are hence useful (among other things) to understand neighbourhoods of submanifolds. Since we will be interested essentially in manifolds of dimension at most 3, two simple cases will be important.

\begin{prop} \label{unique:line:bundle:prop}
A connected compact manifold $M$ has a unique line bundle $E\to M$ up to isomorphism with orientable total space $E$.
\end{prop}

The orientable line bundle on $M$ is a product $M\times \matR$ precisely when $M$ is also orientable. If $M$ is not orientable, the unique orientable line bundle is indicated by $M \timtil \matR$.

\begin{prop} \label{S1:bundle:prop}
For every $n$, there are exactly two vector bundles of dimension $n$ over $S^1$ up to isomorphism, one of which is orientable.
\end{prop}

Again, the orientable vector bundle is $S^1 \times \matR^n $ and the non-orientable one is denoted by $S^1 \timtil \matR^n$. These simple facts allow us to fully understand the possible neighbourhoods of curves in surfaces, and of curves and surfaces inside orientable 3-manifolds. 

\subsection{Manifolds with boundary}
Let a \emph{differentiable manifold $M$ with boundary} be a topological space with charts on a fixed half-space of $\matR^n$ instead of $\matR^n$, forming a smooth atlas. (By definition, maps from subsets of $\matR^n$ are smooth if they locally admit extensions to smooth functions defined on open domains.) 
\index{manifold!manifold with boundary}

The points corresponding to the boundary of the half-space form a subset of $M$ denoted by $\partial M$ and called \emph{boundary}. The boundary of a $n$-manifold is naturally a $(n-1)$-dimensional manifold without boundary. The \emph{interior} of $M$ is $M\setminus \partial M$.

We can define the tangent space $T_pM$ at a point $p\in \partial M$ as the set of all curves in $M$ starting from or arriving to $p$, with the same equivalence relation as above. The space $T_pM$ is a vector space that contains the hyperplane $T_p \partial M$. Most of the notions introduced for manifolds extend in an appropriate way to manifolds with boundary. A \emph{submanifold} of a manifold with boundary $M$ is the image of an embedding $N \hookrightarrow M$, where $N$ is another manifold with boundary.

\index{disc}
The most important manifold with boundary is certainly the \emph{disc} 
$$D^n = \big\{x\ \big|\ \|x\|\leqslant 1\big\}\subset \matR^n.$$ 
More generally, a \emph{disc} in a $n$-manifold $N$ is a submanifold $D\subset N$ with boundary diffeomorphic to $D^n$. Since a disc is in fact a (closed) tubular neighbourhood of any point in its interior, the uniqueness of tubular neighbourhoods implies the following.

\begin{teo}[Cerf -- Palais] \label{dischi:teo}
Let $N$ be connected and oriented. Two orientation-preserving embeddings $f,g\colon D^n \to \interior {N}$ are ambiently isotopic.
\end{teo}

A \emph{boundary component} $N$ of $M$ is a connected component of $\partial M$. A \emph{collar} for $N$ is an open neighbourhood diffeomorphic to $N\times [0,1)$. As for tubular neighbourhoods, every compact boundary component has a collar, unique up to isotopy.
\index{collar}

A \emph{closed manifold} is a compact manifold without boundary.
\index{manifold!closed manifold}

\subsection{Cut and paste} \label{cut:subsection}
If $M\subset N$ is a closed and orientable $(n-1)$-submanifold in the interior of an orientable $n$-manifold $N$, it has a tubular neighbourhood diffeomorphic to $M \times \matR$. The operation of \emph{cutting} $N$ along $M$ consists of the removal of the portion $M\times (-1,1)$. The resulting manifold $N'$ has two new boundary components $M\times \{-1\}$ and $M\times \{1\}$, both diffeomorphic to $M$. By the uniqueness of the tubular neighbourhood, the manifold $N'$ depends (up to diffeomorphisms) only on the isotopy class of $M\subset N$.

\index{cut and paste}
Let $M$ and $N$ be connected. The submanifold $M$ is \emph{separating} if its complement consists of two connected components. The cut manifold has two or one components depending on whether $M$ is separating or not.

Conversely, let $M$ and $N$ be $n$-manifolds with boundary, and let $\varphi \colon \partial M \to \partial N$ be a diffeomorphism. It is possible to \emph{glue} $M$ and $N$ along $\varphi$ and obtain a new $n$-manifold as follows.

A na\"\i ve approach would consist in taking the topological space $M\sqcup N$ and identify $p$ with $\varphi (p)$ for all $p\in \partial M$. The resulting quotient space is indeed a topological manifold, but the construction of a smooth atlas is not immediate. A quicker method consists of taking two collars $\partial M\times [0,1)$ and $\partial N\times [0,1)$ of the boundaries and then consider the topological space
$$(M\setminus \partial M) \sqcup (N\setminus \partial N).$$
Now we identify the points $(p,t)$ and $(\varphi (p), 1-t)$ of the open collars, for all $p\in \partial M$ and all $t\in (0,1)$. Having now identified two \emph{open} subsets of $M\setminus \partial M$ and $N\setminus \partial N$, a differentiable atlas for the new manifold is immediately derived from the atlases of $M$ and $N$.

\begin{prop} \label{depend:diffeo:prop}
The resulting smooth manifold depends (up to diffeomorphism) only on the isotopy class of $\varphi$.
\end{prop}

\subsection{Connected sum} \label{sum:subsection}
\begin{figure}
\centering
\includegraphics[width = 7 cm] {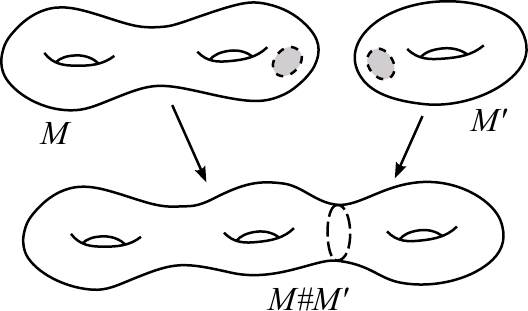}
\nota{A connected sum of closed surfaces.}
\label{Connected_sum:fig}
\end{figure}

We introduce an important cut-and-paste operation. The \emph{connected sum} of two connected oriented $n$-manifolds $M$ and $M'$ is a new $n$-manifold obtained by choosing two $n$-discs $D\subset \interior{M}$ and $D'\subset \interior{M'}$ and an orientation-reversing diffeomorphism $\varphi\colon D \to D'$. The new manifold is constructed in two steps: we first remove the interiors of $D$ and $D'$ from $M\sqcup M'$, thus creating two new boundary components $\partial D$ and $\partial D'$, and then glue these boundary components along $\varphi|_{\partial D}$. See Figure \ref{Connected_sum:fig}.\index{connected sum}

We denote the resulting manifold by $M\# M'$. Since $\varphi$ is orientation-reversing, the manifold $M\# M'$ is oriented coherently with $M$ and $M'$. Theorem \ref{dischi:teo} implies the following.

\begin{prop}
The oriented smooth manifold $M\# M'$ depends (up to diffeomorphism) only on $M$ and $M'$. The connected sum operation $\#$ is commutative, associative, and $S^n$ serves as the identity element.
\end{prop}

We can invert a connected sum as follows. Suppose a $n$-manifold $M$ contains a separating $(n-1)$-sphere $N\subset \interior{M}$. We fix a diffeomorphism $\varphi\colon N \to \partial D^n$. By cutting $M$ along $N$ we get two new boundary components diffeomorphic to $N$ and we glue two discs to them via the map $\varphi$. We get two manifolds $N_1, N_2$ such that $M=N_1\# N_2$.

\subsection{Transversality}
Let $f\colon M \to N$ be a smooth map between manifolds and $X\subset N$ be a submanifold. We say that $f$ is \emph{transverse} to $X$ if for any $p\in f^{-1}X$ the following condition holds:
$$\img (d f_p) + T_{f(p)}X = T_{f(p)}N.$$
The maps transverse to a fixed $X$ are generic, that is they form an open dense subset in the space of all smooth maps from $X$ to $Y$, with respect to an appropriate topology. In particular the following holds.
\begin{teo}
Let $f\colon M \to N$ be a continuous map and $d$ a distance on $N$ compatible with the topology of $N$. For every $\varepsilon >0$ there is a smooth map $g$ transverse to $X$, homotopic to $f$, with $d(f(p), g(p))<\varepsilon$ for all $p\in M$.
\end{teo}
\index{transverse map}

\section{Riemannian geometry}

We briefly introduce Riemannian manifolds and their geometric properties: distance, geodesics, volume, exponential map, injectivity radius, completeness, curvature, and isometries. These notions are carefully explained in do Carmo's \emph{Riemannian Geometry} \cite{Do}.

\subsection{Metric tensor}
A differentiable manifold lacks many natural geometric notions, such as distance between points, angle between vectors, path lengths, geodesics, volumes, etc. It is a quite remarkable fact that the introduction of a single additional mathematical entity suffices to recover all these geometric notions: this miraculous object is the metric tensor.

A \emph{metric tensor} for a differentiable manifold $M$ is the datum of a scalar product on each tangent space $T_p$ of $M$, which varies smoothly on $p$. More specifically, on a chart the scalar product may be expressed as a matrix, and we require that its coefficients vary smoothly on $p$.
\index{metric tensor}

\begin{defn}
A \emph{Riemannian manifold} is a differentiable manifold with a metric tensor that is positive definite at every point.
\end{defn}
\index{manifold!Riemannian manifold}

A Riemannian manifold is usually denoted as a pair $(M,g)$, where $M$ is the manifold and $g$ is the metric tensor. 

\begin{example}
The \emph{Euclidean space} is the manifold $\matR^n$ equipped with the Euclidean metric tensor $g(x,y) = \sum_{i=1}^n x_iy_i$ at every tangent space $T_p = \matR^n$.
\end{example}

\begin{example}
Every differential submanifold $N$ in a Riemannian manifold $M$ is also Riemannian: it suffices to restrict at every $p\in N$ the metric tensor on $T_pM$ to the linear subspace $T_pN$.

In particular, the \emph{sphere} 
$$S^n = \big\{x\in \matR^{n+1}\ \big| \ \|x\| = 1 \big\}$$
is a submanifold of $\matR^{n+1}$ and is hence Riemannian.
\end{example}
\index{sphere}

The metric tensor $g$ defines in particular a norm for every tangent vector, and an angle between tangent vectors with the same basepoint. The velocity $\gamma'(t)$ of a curve $\gamma \colon I \to M$ at time $t\in I$ now has a norm $\|\gamma'(t)\|\geqslant 0$ called \emph{speed}, and two curves that meet at a point with non-zero velocities form a well-defined angle. The \emph{length} of $\gamma$ may be defined as
$$L(\gamma) = \int_I \|\gamma'(t)\| dt$$
and can be finite or infinite. A \emph{reparametrisation} of $\gamma$ is the curve $\eta \colon J \to M$ obtained as $\eta = \gamma\circ\varphi$ where $\varphi\colon J\to I$ is a diffeomorphism of intervals. The length is invariant under reparametrisations, that is $L(\gamma) = L(\eta)$.

\subsection{Distance and geodesics}
Let $(M,g)$ be a connected Riemannian manifold. The curves in $M$ now have a length and hence may be used to define a distance on $M$. 

\begin{defn}
The \emph{distance} $d(p,q)$ between two points $p, q\in M$ is 
$$d(p,q) = \inf_\gamma L(\gamma)$$
where $\gamma$ varies among all curves $\gamma\colon [0,1] \to M$ with $\gamma(0)=p$ and $\gamma(1)=q$. 
\end{defn}
The manifold $M$ equipped with the distance $d$ is a metric space (compatible with the initial topology of $M$).

\begin{defn}
A \emph{geodesic} is a curve $\gamma\colon I \to M$ having constant speed $k$ that realises locally the distance. This means that every $t \in I$ has a closed neighbourhood $[t_0,t_1]\subset I$ such that $d(\gamma(t_0), \gamma(t_1)) = L(\gamma|_{[t_0,t_1]}) = k(t_1-t_0)$.
\end{defn}
\index{geodesic}

Note that with this definition the constant curve $\gamma(t) = p_0$ is a geodesic with constant speed $k=0$. Such a geodesic is called \emph{trivial}. A curve that realises the distance locally may not realise it globally. 

\begin{example}
The non-trivial geodesics in Euclidean space $\matR^n$ are affine lines run at constant speed. 
The non-trivial geodesics in the sphere $S^n$ are portions of great circles, run at constant speed. 
\end{example}

A geodesic $\gamma:I \to M$ is \emph{maximal} if it cannot be extended to a geodesic on a strictly larger interval $J \supset I$. Maximal geodesics are determined by some first-order conditions: 

\begin{teo}  \label{geod:teo}
Let $p\in M$ be a point and $v\in T_pM$ a tangent vector. There exists a unique maximal geodesic $\gamma\colon I \to M$ with $\gamma(0) = p$ and $\gamma'(0) = v$. The interval $I$ is open and contains $0$.
\end{teo}

\subsection{Volume}
If the differentiable manifold $M$ is oriented, the metric tensor also induces a \emph{volume form}.
\index{volume form}

The best method to define volumes on a $n$-manifold $M$ is to construct an appropriate \emph{$n$-form}. A $n$-form $\omega$ is an alternating multilinear form
$$\omega_p: \underbrace{T_p \times \ldots \times T_p}_n \longrightarrow \matR$$
at each point $p\in M$, which varies smoothly with $p$. The alternating condition means that if we swap two vectors the result changes by a sign. 
The $n$-forms are designed to be integrated: it makes sense to write 
$$\int_D \omega$$
on any open set $D$. A \emph{volume form} on an oriented manifold $M$ is a $n$-form $\omega$ such that $\omega_p(v_1,\ldots, v_n)> 0$ for every positive basis $v_1,\ldots, v_n$ of $T_p$ and every $p\in M$. 

The metric tensor defines a volume form as follows: it suffices to set $\omega_p(e_1,\ldots, e_n)=1$ on every positive orthonormal basis $e_1,\ldots, e_n$. With this definition every open set $D$ of $M$ has a well-defined \emph{volume} 
$$\Vol(D) = \int_D \omega$$
which is a positive number or infinity. If $D$ has compact closure the volume is necessarily finite; in particular, a compact Riemannian manifold $M$ has finite volume $\Vol (M)$.
\index{volume}

On a chart, the volume form can be calculated from the metric tensor $g$ via the following formula: 
$$\omega = \sqrt{\det g} \cdot dx_1\cdots dx_n.$$

\subsection{Completeness}
A Riemannian manifold $(M,g)$ is a metric space, and as such it may be complete or not. For instance,
a compact Riemannian manifold is always complete; on the other hand, by removing a point from a Riemannian manifold we always get a non-complete space. Non-compact manifolds like $\matR^n$ typically admit both complete and non-complete Riemannian structures. 
The completeness of a Riemannian manifold may be expressed in various ways:

\index{Hopf -- Rinow theorem}
\begin{teo}[Hopf -- Rinow] \label{HR:teo}
Let $(M,g)$ be a connected Riemannian manifold. The following are equivalent: 
\begin{enumerate}
\item $M$ is complete,
\item a subset of $M$ is compact if and only if it is closed and bounded,
\item every geodesic can be extended on the whole $\matR$.
\end{enumerate}
If $M$ is complete, any two points $p,q \in M$ are joined by a \emph{minimising} geodesic $\gamma$, \emph{i.e.}~a curve such that $L(\gamma) = d(p,q)$.
\end{teo}

\subsection{Exponential map} \label{esponenziale:subsection}
Let $M$ be a complete Riemannian manifold. The fact that geodesics are determined by first-order conditions allows us to introduce the following.

\index{exponential map}
\begin{defn}
Let $p \in M$ be a point. The \emph{exponential map}
$$\exp_p\colon T_pM \longrightarrow M$$
is defined as follows.
A vector $v\in T_pM$ determines a maximal geodesic $\gamma_v\colon \matR \to M$ with $\gamma_v(0) = p$ and $\gamma_v'(0) = v$. We set $\exp_p(v) = \gamma_v(1)$.
\end{defn}

\begin{teo} \label{exp:teo}
The differential of the exponential map $\exp_p$ at the origin is the identity and hence $\exp_p$ is a local diffeomorphism at the origin.
\end{teo}

Via the exponential map, a small neighborhood of the origin in  $T_pM$ can be used as a chart near $p$: we recover here the intuitive idea that the tangent space approximates the manifold near $p$.

When $M$ is not complete, the exponential map is only defined in some open star-shaped neighbourhood $V_p \subset T_pM$ of the origin, and Theorem \ref{exp:teo} holds also in this case.

\subsection{Injectivity radius}
Let $M$ be a complete Riemannian manifold. We introduce the following.

\begin{defn} 
The \emph{injectivity radius} $\inj_pM>0$ of $M$ at a point $p\in M$ is the supremum of all $r>0$ such that $\exp_p|_{B_0(r)}$ is a diffeomorphism onto its image.
\end{defn}
\index{injectivity radius}

Here $B_0(r)$ is the open ball with centre $0$ and radius $r$ in the tangent space $T_pM$. The injectivity radius is always positive by Theorem \ref{exp:teo}. For every $r< \inj_pM$ the exponential map transforms the ball of radius $r$ in $T_p(M)$ into the metric ball of radius $r$ in $M$. That is, the following equality holds: 
$$\exp_p(B_0(r)) = B_p(r)$$
and the metric ball $B_p(r)$ is indeed diffeomorphic to an open ball in $\matR^n$. When $r\geq \inj_p M$ this may not be true: for instance if $M$ is compact there is a $R>0$ such that $B_p(R) = M$.

The injectivity radius $\inj_p(M)$ varies continuously with respect to $p\in M$; 
the injectivity radius $\inj (M)$ of $M$ is defined as
$$\inj (M) = \inf_{p\in M} \inj_p M.$$

\begin{prop} \label{iniettivita:prop}
A compact Riemannian manifold has positive injectivity radius.
\end{prop}
\begin{proof}
The injectivity radius $\inj_p M$ is positive and varies continuously with $p$, hence it has a positive minimum on the compact $M$.
\end{proof}

Finally we note the following. A \emph{loop} is a curve $\gamma\colon [a,b] \to M$ with $\gamma(a) = \gamma(b)$.

\begin{prop}
Let $M$ be a complete Riemannian manifold. Every loop $\gamma$ in $M$ shorter than $2\cdot \inj (M)$ is homotopically trivial.
\end{prop}
\begin{proof}
Set $p=\gamma(a) = \gamma (b)$. Since $\gamma$ is shorter than $2\cdot\inj (M)$, it cannot escape the ball $B_p(r)$ for some $r<\inj (M)\leqslant \inj_p M$. This ball is diffeomorphic to a ball in $\matR^n$, hence in particular it is contractible, so $\gamma$ is homotopically trivial.
\end{proof}

A proof of the continuity of $\inj_p(M)$ can be found in \cite{Lee}.

\subsection{Curvature} \label{curvatura:subsection}
The \emph{curvature} of a Riemannian manifold $(M,g)$ is some mathematical entity that measures how distorted $g$ is when compared to the familiar Euclidean structure on $\matR^n$. The curvature is encoded by various kinds of mathematical objects, and some of them may be frustratingly complicated when $n=\dim M \geqslant 3$. Luckily, most of the Riemannian manifolds that we will encounter in this book have ``constant curvature'' and therefore these objects will simplify dramatically. We briefly recall them here.
\index{curvature}

\index{Levi-Civita connection}
\index{Christoffel symbol}
The metric tensor $g$ induces the \emph{Levi-Civita connection} $\nabla$, which can be used to differentiate any vector field with respect to a tangent vector at any point of $M$. We concentrate on an open chart and hence suppose that $M=U$ is an open set of $\matR^n$; let $e_1, \ldots, e_n$ be the coordinate vector fields on $U$ and $\nabla_i$ be the Levi-Civita differentiation with respect to $e_i$. We have
$$\nabla_i e_j = \Gamma^k_{ij} e_k$$
for some smooth functions $\Gamma^k_{ij} \colon U \to \matR$ called the \emph{Chistoffel symbols}. Here we use the Einstein notation: repeated indices should be added from $1$ to $n$. The Christoffel symbols are determined by $g$ via the formula:
$$\Gamma_{ij}^k = \frac{1}{2}g^{k\ell} \left(\frac{\partial g_{\ell i}}{\partial x^j} + \frac{\partial g_{\ell j}}{\partial x^i} - \frac{\partial g_{ij}}{\partial x^\ell} \right).$$
Here $g^{k\ell}$ is just the inverse matrix of $g = g_{k\ell}$.
Christoffel symbols are very useful but depend heavily on the chosen chart; a more intrinsic object is the \emph{Riemann tensor} that may be defined (quite obscurely) as
\index{Riemann tensor}
$$R^i{}_{jk\ell} = \partial_k\Gamma^i_{\ell j}
    - \partial_\ell\Gamma^i_{k j}
    + \Gamma^i_{k m}\Gamma^m_{\ell j}
    - \Gamma^i_{\ell m}\Gamma^m_{k j}.$$ 
By contracting two indices we get the \emph{Ricci tensor} \index{Ricci tensor}
$$R_{ij} = {R^k}_{ik j} =
\partial_{k}{\Gamma^k_{ji}} - \partial_{j}\Gamma^k_{k i}
+ \Gamma^k_{km} \Gamma^m_{ji}
- \Gamma^k_{jm}\Gamma^m_{k i}.
$$ 
The Ricci tensor is symmetric, and by further contracting we get the \emph{scalar curvature}
$$R = g^{ij}R_{ij}.$$
\index{scalar curvature}

\subsection{Scalar curvature} \label{gauss:subsection}
What kind of geometric information is contained in the tensors introduced above? The scalar curvature $R$ is the simplest curvature detector, and it measures the local discrepancy between volumes in $M$ and in $\matR^n$. Let $B_p(\varepsilon)\subset M$ be the $\varepsilon$-ball in $M$ centred at $p$ and $B_0(\varepsilon)\subset \matR^n$ be the $\varepsilon$-ball in the Euclidean $\matR^n$. We have
$$\Vol\big(B_p(\varepsilon)\big) = \Vol\big(B_0(\varepsilon)\big) \cdot \left( 1- \frac R{6(n+2)}\varepsilon^2 \right) + {\rm O}(\varepsilon^4).$$
We note in particular that $R$ is positive (negative) if $B_p(\varepsilon)$ has smaller (bigger) volume than the usual Euclidean volume.

If $M$ has dimension 2, that is it is a surface, the Riemann and Ricci tensors simplify dramatically and are fully determined by the scalar curvature $R$, which is in turn equal to twice the \emph{gaussian curvature} $K$: if $M$ is contained in $\matR^3$ the gaussian curvature is defined as the product of its two principal curvatures, but when $M$ is abstract principal curvatures make no sense. On surfaces, the formula above reads as
$$\Area(B_p(\varepsilon)) = \pi\varepsilon^2 - \frac{\pi \varepsilon^4}{12}K + {\rm O}(\varepsilon^4).$$
\index{gaussian curvature}

\begin{figure}
\begin{center}
\includegraphics[width = 6 cm] {\iftoggle{BW}{Gaussian_curvature-BW}{640px-Gaussian_curvature}}
\nota{Three surfaces in space (hyperboloid with one sheet, cylinder, sphere) whose gaussian curvature is respectively negative, null, and positive at each point. The curvature on the sphere is constant, while the curvature on the hyperboloid varies: a complete surface in $\matR^3$ cannot have constant negative curvature.}
\label{Gaussian:fig}
\end{center}
\end{figure}

\subsection{Sectional curvature}
If $(M,g)$ has dimension $n\geqslant 3$ the scalar curvature is a weak curvature detector when compared to the Ricci and Riemann tensors. Moreover, there is yet another curvature detector which encodes the same amount of information of the full Riemann tensor, but in a more geometric way: this is the \emph{sectional curvature} and is defined as follows. 
\index{sectional curvature}

\begin{defn}
Let $(M,g)$ be a Riemannian manifold. Let $p \in M$ be a point and $W \subset T_pM$ be a 2-dimensional vector subspace. By Theorem \ref{exp:teo} there exists an open set $U_p \subset T_pM$ containing the origin where $\exp_p$ is a diffeomorphism onto its image. In particular $S = \exp_p(U_p \cap W)$ is a small smooth surface in $M$ passing through $p$, with tangent plane $W$. As a submanifold of $M$, the surface $S$ has a Riemannian structure induced by $g$.

The \emph{sectional curvature} of $(M,g)$ along $W$ is defined as the gaussian curvature of $S$ in $p$.
\end{defn}

We can use the sectional curvature to unveil the geometric nature of the Ricci tensor:
the Ricci tensor $R_{ij}$ in $p$ measures the average sectional curvature along axes, that is for every unit vector $v\in T_p(M)$ the number $R_{ij} v^iv^j$ is $(n-1)$ times the average sectional curvature of the planes $W\subset T_p(M)$ containing $v$.

The Riemann tensor is determined by the sectional curvatures and vice-versa. In dimension $n=3$ the Ricci tensor fully determines the sectional curvatures and hence also the Riemann tensor. This is not true in dimension $n\geqslant 4$.

\subsection{Constant sectional curvature}
A Riemannian manifold $(M,g)$ has \emph{constant sectional curvature} $K$ if the sectional curvature of every 2-dimensional vector space $W\subset T_pM$ at every point $p\in M$ is always $K$.

\begin{oss}
On a Riemannian manifold $(M,g)$ one may \emph{rescale} the metric by some factor $\lambda >0$ substituting $g$ with the tensor $\lambda g$. At every point the scalar product is rescaled by $\lambda$. Consequently, lengths of curves are rescaled by $\sqrt \lambda$ and volumes are rescaled by $\lambda^{\frac n2}$. The sectional curvature is rescaled by $1/\lambda$. 
\end{oss}

By rescaling the metric we may transform every Riemannian manifold with constant sectional curvature $K$ into one with constant sectional curvature $-1$, $0$, or $1$.

\begin{example}
Euclidean space $\matR^n$ has constant curvature zero. A sphere of radius $R$ has constant curvature $1/R^2$.
\end{example}

\subsection{Isometries}
Every honest category has its morphisms. Riemannian manifolds are so rigid, that in fact one typically introduces only isomorphisms: these are called \emph{isometries}.
\index{isometry}

\begin{defn}
A diffeomorphism $f\colon M\to N$ between two Riemannian manifolds $(M,g)$ and $(N,h)$ is an \emph{isometry} if it preserves the scalar product. That is, the equality
$$\langle v,w \rangle = \langle df_p(v), df_p(w) \rangle$$
holds for all $p\in M$ and every pair of vectors $v,w \in T_pM$. The symbols $\langle,\rangle$ indicate the scalar products in $T_pM$ and $T_{f(p)}N$.
\end{defn}

As we said, isometries are extremely rigid. These are determined by their first-order behaviour at any single point.

\begin{teo} \label{isom:teo}
Let $f,g\colon M \to N$ be two isometries between two connected Riemannian manifolds. If there is a point $p\in M$ such that $f(p) = g(p)$ and $df_p = dg_p$, then $f  =g$ everywhere.
\end{teo}
\begin{proof}
Let us show that the subset $S\subset M$ of the points $p$ such that $f(p) = g(p)$ and $df_p = dg_p$ is open and closed.

The locus where two functions coincide is typically closed, and this holds also here (to prove it, take a chart). We prove that it is open: pick $p\in S$. 
By Theorem \ref{exp:teo} there is an open neighbourhood $U_p \subset T_pM$ of the origin where the exponential map is a diffeomorphism onto its image. We show that the open set $\exp_p(U_p)$ is entirely contained in $S$.

A point $q\in\exp_p(U_p)$ is the image $q=\exp(v)$ of a vector $v\in U_p$ and hence $q=\gamma(1)$ for the geodesic $\gamma$ determined by the data $\gamma(0)=p, \gamma'(0)=v$. The maps $f$ and $g$ are isometries and hence send geodesics to geodesics: here $f\circ \gamma$ and $g\circ\gamma$ are geodesics starting from $f(p)=g(p)$ with the same initial velocities and thus they coincide. This implies that $f(q)=g(q)$. Since $f$ and $g$ coincide on the open set $\exp_p(U_p)$, also their differentials do.
\end{proof}

\subsection{Local isometries}
A \emph{local isometry} $f\colon M \to N$ between Riemannian manifolds is a map where every $p\in M$ has an open neighbourhood $U$ such that $f|_U$ is an isometry onto its image. Theorem \ref{isom:teo} applies with the same proof to local isometries.
\index{isometry!local isometry}

The following proposition relates nicely the notions of local isometry, topological covering, and completeness.

\begin{prop} \label{local:isometry:prop}
Let $f\colon M \to N$ be a local isometry.
\begin{enumerate}
\item If $M$ is complete, the map $f$ is a covering.
\item If $f$ is a covering, then $M$ is complete $\Longleftrightarrow$ $N$ is complete.
\end{enumerate}
\end{prop}
\begin{proof}
Since $f$ is a local isometry, every geodesic in $M$ projects to a geodesic in $N$. If $f$ is also a covering, the converse holds: every geodesic in $N$ lifts to a geodesic in $M$ (at any starting point). 

If $f$ is a covering we can thus project and lift geodesics via $f$: therefore every geodesic in $M$ can be extended to $\matR$ if and only if every geodesic in $N$ can; this proves (2) using the Hopf -- Rinow Theorem \ref{HR:teo}.

We prove (1) by showing that the ball $B=B(p,\inj_pN)$ is well-covered for all $p\in N$. Since $M$ is complete, every geodesic in $N$ can be lifted to a geodesic in $M$ (at any starting point).  For every $\tilde p \in f^{-1}(p)$ the map $f$ sends the geodesics exiting from $\tilde p$ to geodesics exiting from $p$ and hence sends isometrically $B(\tilde p, \inj_p N)$ onto $B$. On the other hand, given a point $q\in f^{-1}(B)$, the geodesic in $B$ connecting $f(q)$ to $p$ lifts to a geodesic connecting $q$ to some point $\tilde p \in f^{-1}(p)$. Therefore 
$$f^{-1}\big(B(p,\inj_p N)\big) = \bigsqcup_{\tilde p \in f^{-1}(p)} B(\tilde p, \inj_p N)$$
and $f$ is a covering.
\end{proof}

\begin{prop} \label{volume:covering:prop}
Let $f\colon M \to N$ be a local isometry and a degree-$d$ covering. We have
$$\Vol(M) = d \cdot \Vol (N).$$
\end{prop}
\begin{proof}[Sketch of the proof]
We may find a disjoint union of well-covered open sets in $N$ whose complement has zero measure. Every such open set lifts to $d$ copies of it in $M$, and the zero-measure set lifts to a zero-measure set.
\end{proof}

The formula makes sense also when some of the quantities $\Vol(M)$, $\Vol(N)$, and $d$ are infinite.

\subsection{Totally geodesic submanifolds}
A differentiable submanifold $M$ in a Riemannian manifold $N$ is \emph{totally geodesic} if 
every geodesic in $M$ with the induced metric is also a geodesic in $N$. 
\index{submanifold!totally geodesic submanifold}

When $\dim M=1$ this notion is equivalent to that of an unparametrized embedded geodesic; if $\dim M\geqslant 2$ then $M$ is a quite peculiar object: generic Riemannian manifolds do not contain totally geodesic surfaces at all. An equivalent condition is that, for every $p\in M$ and every $v \in T_p(M)$, the unique geodesic in $N$ passing through $p$ with velocity $v$ stays in $M$ for some interval $(-\varepsilon, \varepsilon)$.

\subsection{Riemannian manifolds with boundary}
Many geometric notions in Riemannian geometry extend easily to manifolds $M$ with boundary. The boundary $\partial M$ of a Riemannian manifold $M$ is naturally a Riemannian manifold without boundary. 
\index{submanifold!Riemannian submanifold}
A particularly nice (and exceptional) case is when $\partial M$ is totally geodesic. 

\section{Measure theory} \label{measure:theory:section}
We will use some basic measure theory only in Chapter \ref{automorfismi:chapter}.

\subsection{Borel measure}
A \emph{Borel set} in a topological space $X$ is any subset obtained from open sets through the operations of countable union, countable intersection, and relative complement. Let $\calF$ denote the set of all Borel sets.
A \emph{Borel measure} on $X$ is a function $\mu\colon \calF \to [0,+\infty]$ which is additive on any countable collection of disjoint sets. 
\index{Borel set}
\index{Borel measure}

The measure is \emph{locally finite} if every point has a neighbourhood of finite measure and is \emph{trivial} if $\mu(S)=0$ for all $S\in\calF$. 

\begin{ex} If $\mu$ is a locally finite Borel measure then $\mu(X)<+\infty$ for any compact Borel set $K\subset X$.
\end{ex}

\begin{example}
Let $D\subset X$ be a discrete set. The \emph{Dirac measure} $\delta_D$ concentrated in $D$ is the measure
$$\delta_D(S) = \#(S\cap D).$$
Since $D$ is discrete, the measure $\delta_D$ is locally finite.
\end{example}

The \emph{support} of a measure is the set of all points $x\in X$ such that $\mu(U)>0$ for any open set $U$ containing $x$. The support is a closed subset of $X$. The measure is \emph{fully supported} if its support is $X$. The support of $\delta_D$ is of course $D$. A point $x\in X$ is an \emph{atom} for $\mu$ if $\mu(\{x\})>0$.

\subsection{Construction by local data}
A measure can be defined using local data in the following way.

\begin{prop} \label{ricoprimento:prop}
Let $\{U_i\}_{i\in I}$ be a countable, locally finite open covering of $X$ and for any $i\in I$ let $\mu_i$ be a locally finite Borel measure on $U_i$. If $\mu_i|_{U_i\cap U_j} = \mu_j|_{U_i\cap U_j}$ for all $i,j\in I$ there is a unique locally finite Borel measure  $\mu$ on $X$ whose restriction to $U_i$ is $\mu_i$ for all $i$.
\end{prop}
\begin{proof}
For every finite subset $J\subset I$ we define 
$$X_J = \big(\cap_{j\in J} U_j\big) \setminus \big(\cup_{i\in I\setminus J} U_i\big).$$ 
The sets $X_J$ form a countable partition of $X$ into Borel sets and every $X_J$ is equipped with a measure $\mu_J = \mu_j|_{X_j}$ for any $j\in J$. Define $\mu$ by setting
$$\mu(S)= \sum_{j\in J} \mu(S\cap X_j)$$
on any Borel $S\subset X$. 
\end{proof}

When $X$ is a reasonable space some hypothesis may be dropped.

\begin{prop} \label{ricoprimento2:prop}
If $X$ is paracompact and separable, Proposition \ref{ricoprimento:prop} holds for any open covering $\{U_i\}_{i\in I}$.
\end{prop}
\begin{proof}
By paracompactness and separability the open covering $\{U_i\}$ has a refinement that is locally finite and countable: apply Proposition \ref{ricoprimento:prop} to the refinement to get a measure $\mu$. To prove that indeed $\mu|_{U_i} = \mu_i$ apply Proposition \ref{ricoprimento:prop} again to the covering of $U_i$ given by the refinement.
\end{proof}

\subsection{Topology on the measures space} \label{topology:measure:subsection}
In what follows we suppose for simplicity that $X$ is a finite-dimensional topological manifold, although everything is valid in a much wider generality. We indicate by $\calM(X)$ the space of all locally finite Borel measures on $X$ and by $C_c(X)$ the space of all continuous functions $X\to \matR$ with compact support: the space $C_c(X)$ is not a Banach space, but it is a topological vector space.

Recall that the \emph{topological dual} of a topological vector space $V$ is the vector space $V^*$ formed by all continuous linear functionals $V\to \matR$.
A measure $\mu\in\calM(X)$ acts like a continuous functional on $C_c(X)$ as follows
$$\mu\colon f \longmapsto \int_\mu f$$
and hence defines an element of $C_c^*(X)$. A functional in $C_c^*(X)$ is \emph{positive} if it assumes non-negative values on non-negative functions. 

\begin{teo}[Riesz representation] The space $\calM(X)$ may be identified in this way to the subset in $C_c(X)^*$ of all positive functionals. 
\end{teo}

The space $\calM(X)$ in $C_c(X)^*$ is closed with respect to sum and product with a positive scalar.
We now use the embedding of $\calM(X)$ into $C_c(X)^*$ to define a natural topology on $\calM(X)$.

\begin{defn}
Let $V$ be a real topological vector space. Every vector $v\in V$ defines a functional in $V^*$ as $f\mapsto f(v)$. The \emph{weak-* topology} on $V^*$ is the weakest topology among those where these functionals are continuous.
\end{defn}

We give $C_c(X)^*$, and hence $\calM(X)$, the weak-* topology. 

\subsection{Sequences of measures} \label{mu:sequence:subsection}
By definition, a sequence of measures  $\mu_i$ converges to $\mu$ if and only if $\int_{\mu_i} f \to \int_\mu f$ for any $f\in C_c(X)$. This type of weak convergence is usually denoted with the symbol
$\mu_i \rightharpoonup \mu$.
 
\begin{ex} 
Let $x_i$ be a sequence of points in $X$ that tend to $x\in X$: we get $\delta_{x_i} \rightharpoonup \delta_x$. 
\end{ex}

It is important to note that $\mu_i \rightharpoonup \mu$ does \emph{not} imply $\mu_i(U) \to \mu(U)$ for any open (or closed) set $U$: consider for instance a sequence $\mu_i = \delta_{x_i}$ with $x_i$ exiting from (or entering into) the set $U$. We can get this convergence on compact sets if we can control their topological boundary.

\begin{prop} \label{delta:zero:prop}
Let $K\subset X$ be a Borel compact subset. If $\mu_i \rightharpoonup \mu$ and $\mu(\partial K)=0$ then $\mu_i(K) \to \mu(K)$.
\end{prop}

On Banach spaces, the unitary ball is compact in the weak-* topology. Here $\calM(X)$ is not a Banach space, but we have an analogous compactness theorem.

\begin{teo} \label{compact:measure:teo}
A sequence of measures $\mu_i$ such that $\mu_i(K)$ is bounded on every Borel compact set $K\subset X$ converges on a subsequence.
\end{teo}

\section{Groups}
We recall some basic definitions and properties of groups. 

\subsection{Presentations}
Recall that a \emph{finite presentation} of a group $G$ is a description of $G$ as
$$\langle g_1,\ldots, g_k\ |\ r_1,\ldots, r_s \rangle$$
where $g_1,\ldots, g_k\in G$ are the \emph{generators} and $r_1,\ldots, r_s$ are words in $g_i^{\pm 1}$ called \emph{relations}, such that 
$$G \Isom F(g_i)/_{N(r_j)}$$
where $F(g_i)$ is the free group generated by the $g_i$'s and $N(r_j)\triangleleft F(g_i)$ is the \emph{normal closure} of the $r_j$'s, the smallest normal subgroup containing them. 
\index{finite presentation}

Not every group $G$ has a finite presentation: a necessary (but not sufficient) condition is that $G$ must be finitely generated. 

\subsection{Commutators}
Let $G$ be a group. The \emph{commutator} of two elements $h,k \in G$ is the element \index{commutator}
$$[h,k] = hkh^{-1}k^{-1}.$$
The commutator $[h,k]$ is trivial if and only if $h$ and $k$ commute.
Commutators do not form a subgroup in general, but we can use them to generate one. 

More generally, let $H, K<G$ be two subgroups. We define $[H, K]$ as the subgroup of $G$ generated by all commutators $[h,k]$
where $h$ and $k$ vary in $H$ and $K$, respectively. Every element in $[H, K]$ is a product of commutators $[h,k]$ and of their inverses $[h,k]^{-1}=[k,h]$, and we have $[H, K] = [K, H]$.

\begin{prop} If $H$ and $K$ are normal subgroups of $G$ then $[H,K] < H\cap K$ and $[H,K]$ is a normal subgroup of $G$.
\end{prop}
\begin{proof}
If $H$ and $K$ are normal we get $[h,k] = hkh^{-1}k^{-1} \in H\cap K$. Moreover, for every $g\in G$ we have 
$$g^{-1}\cdot[h_1,k_1]^{\pm 1} \cdots [h_i, k_i]^{\pm 1} \cdot g = [g^{-1}h_1g,g^{-1}k_1g]^{\pm 1} \cdots
[g^{-1}h_ig, g^{-1}k_i g]^{\pm 1}$$
hence $[H,K]$ is normal in $G$ if $H$ and $K$ are.
\end{proof}

The group $[G,G]$ is the \emph{commutator subgroup} of $G$. It is trivial if and only if $G$ is abelian.
\index{commutator!commutator subgroup}

\subsection{Series}
We can use commutators iteratively to create some characteristic subgroups of $G$. There are two natural ways to do this, and they produce two nested sequences of subgroups.

The \emph{lower central series} of $G$ is the sequence of normal subgroups
\index{lower central series}
$$ G = G_0 > G_1 > \ldots > G_n > \ldots $$
defined iteratively by setting
$$G_{n+1} = [G_n, G].$$ The \emph{derived series} is the sequence of normal subgroups
\index{derived series}
$$G = G^{(0)} \triangleright G^{(1)} \triangleright \ldots \triangleright G^{(n)} \triangleright \ldots $$
defined by setting
$$G^{(n+1)} = [G^{(n)}, G^{(n)}].$$
We clearly have $G^{(n)} < G_{n}$ for all $n$. 
\index{group!nilpotent group} \index{group!solvable group}
\begin{defn} A group $G$ is \emph{nilpotent} if $G_n$ is trivial for some $n$. It is \emph{solvable} if $G^{(n)}$ is trivial for some $n$.
\end{defn}
The following implications are obvious:
\begin{center}
$G$ abelian $\Longrightarrow$ $G$ nilpotent $\Longrightarrow$ $G$ solvable.
\end{center}
The converses are false, as we will see.

\begin{prop} \label{subgroup:prop}
Subgroups and quotients of abelian (nilpotent, solvable) groups are also abelian (nilpotent, solvable).
\end{prop}
\begin{proof}
If $H<G$ then $H_n<G_n$ and $H^{(n)} < G^{(n)}$. If $H=G/_N$, then $H_n = G_n/_N$ and $H^{(n)} = G^{(n)}/_N$. Both equalities are proved by induction on $n$.
\end{proof}

\subsection{Nilpotent groups} \label{nilpotent:subsection}
Every abelian group $G$ is obviously nilpotent since $G_1$ is trivial. 

\index{group!Heisenberg group}
\begin{ex} \label{Heisenberg:ex}
Let the \emph{Heisenberg group} $\Nil$ consist of all matrices 
$$\begin{pmatrix} 1 & x & z \\ 0 & 1 & y \\ 0 & 0 & 1 \end{pmatrix}$$
where $x, y,$ and $z$ vary in $\matR$, with the multiplication operation. Prove that the Heisenberg group is non-abelian and nilpotent. Indeed the matrices with $x=y=0$ form an abelian subgroup $\matR < \Nil$ and we have
$$[\Nil, \Nil] = \matR, \qquad [\matR, \Nil] = \{e\}.$$
\end{ex}

We will use at some point the following criterion.

\begin{prop} \label{criterion:nilpotent:prop}
Let $G$ be a group generated by some set $S$ and $n>0$ a number. Suppose that
$$[a_1, \ldots [a_{n-1}, [a_n, b]] \cdots ]$$
is trivial for all $a_1,\ldots,a_n,b \in S$. Then $G_n=\{e\}$ and thus $G$ is nilpotent.
\end{prop}
\begin{proof}
We claim that $G_n$ is generated by some elements of type
$$[a_1, \ldots [a_{m-1}, [a_m, b]] \cdots ]$$
with $m\geqslant n$ and $a_1,\ldots,a_m,b\in S$: this clearly implies the proposition.
The claim is proved by induction on $n$ using the formula
$$[a,bc] = [a,b]\cdot [b,[a,c]]\cdot [a,c]$$
which holds in every group.
\end{proof}

We note the following. 
\begin{prop} \label{center:prop}
A nilpotent non-trivial group has non-trivial centre.
\end{prop}
\begin{proof}
Let $G_n$ be the last non-trivial group in the lower central series. Since $G_{n+1} = [G_n,G]$ is trivial, the centre of $G$ contains $G_n$.
\end{proof}

\subsection{Solvable groups}
Solvable groups form a strictly larger class than nilpotent groups.
\begin{ex}
The permutation group $S_3$ is solvable but not nilpotent.
\end{ex}
\begin{ex}
The group ${\rm Aff}(\matR)=\{x\mapsto ax+b\ |\ a\in\matR^*, b\in\matR\}$ of affine transformations is solvable but not nilpotent. Indeed $[{\rm Aff}(\matR), {\rm Aff}(\matR)] = \matR$ consists of all translations and
$$[\matR, {\rm Aff}(\matR)] = \matR, \qquad [\matR, \matR] = \{e\}.$$
\end{ex}

\begin{ex} The permutation group $S_n$ is not solvable for $n\geqslant 5$.
\end{ex}
The latter fact is related to the existence of polynomials of degree $n\geqslant 5$ that are not solvable by radicals. 

Solvable groups are farther from being abelian than nilpotent groups: for instance, they may have trivial centre (like $S_3$). However, they still share some  nice properties with the abelian world:

\begin{prop} \label{solvable:normal:prop}
A solvable group contains a non-trivial normal abelian subgroup.
\end{prop}
\begin{proof}
Let $G^{(n)}$ be the last non-trivial group in the derived series. It is normal in $G$ and abelian, since $[G^{(n)}, G^{(n)}]$ is trivial.
\end{proof}

\subsection{Lie groups}
A \emph{Lie group} is a smooth manifold $G$ which is also a group, such that the operations \index{group!Lie group}
\begin{align*}
G\times G & \to G, \quad (a,b) \mapsto ab \\
G & \to G, \quad a \mapsto a^{-1}
\end{align*}
are smooth. 

\begin{example}
The basic examples are $\GL(n,\matR)$ and $\GL(n,\matC)$ consisting of invertible real and complex $n\times n$ matrices, with the multiplication operation. These contain many interesting Lie subgroups:
\begin{itemize}
\item $\SL(n,\matR)$ and $\SL(n,\matC)$, the matrices with determinant 1;
\item $\On(n)$ and $\Un(n)$, the orthogonal and unitary matrices;
\item $\SO(n) = \SL(n,\matR) \cap \On(n)$ and $\SU(n) = \SL(n,\matC)\cap \Un(n)$.
\end{itemize}
\end{example}

The Lie groups $\On(n)$ and U$(n)$ are compact.
\begin{ex}
The following Lie groups are isomorphic:
$$S^1 \isom \Un(1) \isom \SO(2). $$
\end{ex}
Here are some ways to construct more examples:
\begin{itemize}
\item the product of two Lie groups is a Lie group,
\item a closed subgroup of a Lie group is a Lie group,
\item the universal cover of a Lie group is a Lie group.
\end{itemize}

A Lie group $G$ is often not connected, and we denote by $G^\circ$ the connected component of $G$ containing the identity. The subset $G^\circ$ is a normal Lie subgroup of $G$. For instance $\On(n)$ has two components and $\On(n)^\circ = \SO(n)$. 

Another important example is the Lie group
$$\On(m,n) = \big\{A \in \GL(m+n,\matR) \ \big| \ \transp A I_{m,n} A = I_{m,n} \big\}$$
defined for any pair $m,n \geqslant 0$ of integers, where $I_{m,n}$ is the diagonal matrix
$$I_{m,n} = \begin{pmatrix} I_m & 0 \\ 0 & -I_n \end{pmatrix}.$$
The Lie group $\On(m,n)$ has two connected components if either $m$ or $n$ is zero (in that case we recover the orthogonal group) and has four connected components otherwise.

\begin{oss}
Not all the Lie groups are subgroups of $\GL(n,\matC)$: for instance, the universal cover of $\SL(2,\matR)$ is not.
\end{oss}

A \emph{homomorphism} $\varphi \colon G \to H$ of Lie groups is a smooth group homomorphism. The kernel of $\varphi$ is a closed subgroup and hence a Lie subgroup of $G$. An isomorphism of Lie groups is a group isomorphism that is also a diffeomorphism. 

\subsection{Vector fields, metrics, and differential forms}
A Lie group $G$ acts on itself in two ways: by left and right multiplication. Both actions are smooth, free, and transitive.

A geometric object on $G$ is \emph{left-invariant} if it is invariant by left multiplication. Right-invariance is defined analogously. In general, it is easy to find objects that are either left- or right-invariant, but not necessarily both. 

There is a natural 1-1 correspondence between:
\begin{itemize}
\item vectors in $T_eG$ and left-invariant vector fields on $G$, 
\item scalar products on $T_eG$ and left-invariant Riemannian metrics on $G$,
\item $n$-forms on $T_eG$ and left-invariant $n$-forms on $G$,
\item orientations on $T_eG$ and left-invariant orientations on $G$.
\end{itemize}

This holds because left-multiplication by an element $g\in G$ is a diffeomorphism that transports everything from $T_eG$ to $T_gG$. The same holds for right-invariant objects.

Every basis of $T_eG$ extends in this way to $n$ independent left-invariant vector fields: this shows that $G$ is \emph{parallelizable}, \emph{i.e.}~the tangent bundle of $G$ is trivial.

\subsection{Simple Lie groups} \label{simple:Lie:subsection}
A \emph{simple Lie group} $G$ is a connected, non-abelian Lie group $G$ that does not contain any non-trivial connected normal subgroup (the trivial cases being of course $\{e\}$ and $G$ itself).

The classification of simple Lie groups is due to E.~Cartan. The following theorem furnishes many examples.

\begin{teo}
The groups $\On(m,n)^\circ$ are all simple when $m+n\geqslant 3$, except the cases $(m,n) = (4,0), (2,2), (0,4)$.
\end{teo}

\subsection{Haar measures} \label{Haar:subsection}
Let $G$ be a Lie group. As we just said, non-vanishing left-invariant $n$-forms $\omega$ on $G$ are in 1-1 correspondence with non-trivial $n$-forms on $T_eG$, and since the latter are all proportional, the form $\omega$ is unique up to multiplication by a non-zero scalar.

The form $\omega$ defines a left-invariant orientation (a basis $v_1,\ldots, v_n$ is positive when $\omega(v_1,\ldots, v_n)>0$); it is a volume form and hence defines a left-invariant locally finite Borel measure on $G$, called the \emph{Haar measure} of $G$. Summing up, the Haar measure depends only on $G$ up to rescaling.

Since left- and right-multiplications commute, right-multiplication by an element $g\in G$ transforms $\omega$ into another non-vanishing left-invariant form which must be equal to $\lambda(g) \omega$ for some positive real number $\lambda(g)$. This defines a Lie group homomorphism $\lambda \colon G \to \matR_{>0}$ to the multiplicative group of positive real numbers, called the \emph{modular function}.

The Lie group $G$ is \emph{unimodular} if the left-invariant Haar measures are also right-invariant, that is if the modular function $\lambda$ is trivial. 
Recall that a group $G$ is \emph{simple} if its normal subgroups are $\{e\}$ and $G$.
\index{group!Lie group!unimodular Lie group}

\begin{prop} \label{compact:unimodular:prop}
Compact, abelian, or simple Lie groups are unimodular.
\end{prop}
\begin{proof}
If $G$ is compact, the image of $\lambda$ is compact in $\matR_{>0}$ and hence trivial. If $G$ is abelian, left- and right-multiplications coincide. If $G$ is simple, the connected normal Lie subgroup $(\ker \lambda)^\circ$ must be $G$. 
\end{proof}

\begin{cor} \label{O:unimodular:cor}
The group $\On(m,n)$ is unimodular if $(m,n)$ is distinct from $(4,0), (2,2), (0,4)$.
\end{cor}
\begin{proof}
The identity component $\On(m,n)^\circ$ is simple and hence unimodular. The modular function $\lambda\colon \On(m,n) \to \matR$ is trivial on the finite-index subgroup $\On(m,n)^\circ$ and is hence trivial.
\end{proof}

\begin{oss}
When $(m,n) = (4,0), (2,2), (0,4)$ the group $\On(m,n)$ is \emph{semisimple}, which means that it looks roughly like a product of simple groups, and is in fact unimodular also in this case. We will not need to introduce this concept rigorously here. 

The solvable Lie group ${\rm Aff}(\matR)$ is not unimodular.
\end{oss}

\subsection{Discrete subgroups}
Let $G$ be a Lie group. A closed subgroup $H< G$ is \emph{discrete} if it forms a discrete topological subset, that is if every point in $H$ is isolated.

\begin{ex} \label{isolated:ex}
A subgroup $H<G$ is discrete if and only if $e\in G$ is an isolated point in $G$.
\end{ex}

\subsection{The Selberg Lemma}
The group $\GL(n,\matC)$ is the prototypical Lie group, since it contains many Lie groups. A couple of purely algebraic facts about this group will have some important geometric consequences in this book. The first is the Selberg Lemma. Recall that a group \emph{has no torsion} if every non-trivial element in it has infinite order.

\begin{lemma}[Selberg's Lemma] \label{Selberg:lemma}
 Let $G$ be a finitely generated subgroup of $\GL(n,\matC)$. There is a finite-index normal subgroup $H\triangleleft G$ without torsion.
\end{lemma}
\index{Selberg lemma}

The second fact is quite related. A group $G$ is \emph{residually finite} if one of the following equivalent conditions holds:
\begin{itemize}
\item for every non-trivial element $g\in G$ there is a finite group $F$ and a surjective homomorphism $\varphi\colon G \to F$ with $\varphi(g) \neq e$;
\item for every non-trivial element $g \in G$ there is a finite-index normal subgroup $H \triangleleft G$ which does not contain $g$;
\item the intersection of all finite-index normal subgroups in $G$ is trivial. 
\end{itemize}
\index{group!residually finite group}

\begin{lemma} \label{RF:lemma}
Every finitely generated subgroup of $\GL(n,\matC)$ is residually finite.
\end{lemma}

The proof of these two lemmas is not particularly hard (see for instance \cite[Chapter 7.6]{R}) but it employs some purely algebraic techniques that are distant from the scope of this book. Their geometric consequences, as we will see, are quite remarkable.

\subsection{Lie algebras}
A \emph{Lie algebra} $\frakg$ is a real vector space equipped with an alternate bilinear product \index{Lie algebra}
$$[,] \colon \frakg \times \frakg \to \frakg$$ 
called \emph{Lie bracket} that satisfies the Jacobi identity
$$[x,[y,z]]+[y,[z,x]] + [z,[x,y]]=0$$
for all $x,y,z \in \frakg$. 

The tangent space $T_eG$ at $e\in G$ of a Lie group $G$ has a natural Lie algebra structure and is denoted by $\frakg$. It is equipped with a natural \emph{exponential map} $\exp\colon \frakg \to G$ that sends $0$ to $e$ and is a local diffeomorphism at $0$.

The differential $d\varphi_e\colon \frakg \to \frakh$ of a Lie group homomorphism $\varphi \colon G\to H$ is a Lie algebra homomorphism and is denoted by $\varphi_*$. The diagram
$$
\xymatrix{ 
\frakg \ar[r]^{\varphi_*} \ar[d]_\exp & \frakh \ar[d]^\exp \\
G\ar[r]_\varphi & H
}
$$
commutes. The basic example is the following: the Lie algebra $\frakgl(n,\matC)$ of $\GL(n,\matC)$ is the vector space $M(n,\matC)$ of all $n\times n$ matrices; here the map $\exp$ is the usual matrix exponential and $[A,B] = AB-BA$ is the usual commutator bracket. 
The Lie algebras of $\GL(n,\matR)$, $\SL(n,\matR)$, $\On(n)$, $\SO(n)$ are: 
\begin{align*}
\frakgl(n,\matR) & = M(n,\matR), \\
\fraksl(n,\matR) & = \{A\in M(n,\matR)\ |\ \tr A = 0 \}, \\
\frakon(n) & = \{A \in M(n, \matR) \ |\ A + \transp A = 0 \}, \\
\frakso(n) & = \{A \in M(n, \matR) \ |\ A + \transp A = 0,\ \tr A = 0 \}.
\end{align*}

\begin{prop}  \label{Lie:covering:prop}
Let $G$ and $H$ be connected. A Lie group homomorphism $\varphi\colon G \to H$ is a topological covering if and only if $\varphi_*$ is an isomorphism.
\end{prop}

\section{Group actions}
Groups acting of spaces are so important in geometry, that they are sometimes used as a definition of  ``geometry'' itself.

\subsection{Definitions}
The \emph{action} of a group $G$ on a topological space $X$ is a homomorphism 
$$G\to \Omeo (X)$$ 
where $\Omeo(X)$ is the group of all self-homeomorphisms of $X$. The quotient set $X/_G$ is the set of all orbits in $X$ and we give it the usual quotient topology. We denote by $g(x)$ the image of $x\in X$ along the homeomorphism determined by $g\in G$. The action is:
\begin{itemize}
\item \emph{free} if $g(x) \neq x$ for all non-trivial $g\in G$ and all $x\in X$;
\item \emph{properly discontinuous} if any two points $x,y \in X$ have neighbourhoods $U_x$ and $U_y$ such that the set
$$\big\{g \in G\ \big|\ g(U_x) \cap U_y \neq \emptyset \big\}$$
is finite.
\end{itemize}
The relevance of these definitions is due to the following.
\begin{prop} \label{equivalent:prop}
Let $G$ act on a Hausdorff connected space $X$. The following are equivalent:
\begin{enumerate}
\item $G$ acts freely and properly discontinuously;
\item the quotient $X/_G$ is Hausdorff and the map $X \to X/_G$ is a covering.
\end{enumerate}
\end{prop}

A covering of type $X \to X/_G$ is called \emph{regular}. A covering $X \to Y$ is regular if and only if the image of $\pi_1(X)$ in $\pi_1(Y)$ is normal, and in that case $G$ is the quotient of the two groups. In particular every universal cover is regular.
Summing up:

\begin{cor} \label{XG:cor}
Every path-connected locally contractible Hausdorff topological space $X$ is the quotient $\tilde X/_G$ of its universal cover by the action of some group $G$ acting freely and properly discontinuously.
\end{cor}

The group $G$ is isomorphic to $\pi_1(X)$.

\begin{ex} \label{proper:ex}
Let a discrete group $G$ act on a locally compact space $X$. The following are equivalent:
\begin{itemize}
\item the action is properly discontinuous; 
\item for every compact $K\subset X$, the set $\big\{g\ |\ g(K) \cap K \neq \emptyset\big\}$ is finite; \item the map 
$ G \times X \to X \times X$ that sends $(g,x)$ to $(g(x), x)$ is proper.
\end{itemize}
\end{ex}

\subsection{Isometry group}
We now want to consider the case where a group $G$ acts by isometries on a Riemannian manifold $M$. In this case we have a homomorphism
$$G \to \Iso(M)$$
with values in the \emph{isometry group} $\Iso(M)$ of $M$, \emph{i.e.}~the group of all isometries $f \colon M \to M$. The group $\Iso(M)$ is much smaller than $\Omeo (M)$. \index{isometry!isometry group} \index{Myers -- Steenrod theorem}

\begin{teo}[Myers -- Steenrod] \label{MS:teo}
The group $\Iso (M)$ has a natural Lie group structure compatible with the compact-open topology. The map $F\colon \Iso(M) \times M \to M \times M$ that sends
$(\varphi, p)$ to $(\varphi(p), p)$ is proper.
\end{teo}
\begin{cor} If $M$ is compact then $\Iso(M)$ is compact.
\end{cor}

\begin{prop} \label{stabilizer:compact:prop}
The stabilizer of any point $x\in M$ is a compact Lie subgroup of $\Iso(M)$.
\end{prop}
\begin{proof}
The stabilizer of $x$ is closed and hence is a Lie subgroup of $\Iso(M)$. The isometries that fix $x$ are determined by their orthogonal action on $T_xM$ and therefore form a compact Lie subgroup of $O(n)$. 
\end{proof}

The hard part of the Myers -- Steenrod Theorem is to endow $\Iso(M)$ with a Lie group structure. In fact, in all the concrete cases that we will encounter, the Lie group structure of $\Iso(M)$ will be evident from the context, so we will not need the full strength of Theorem \ref{MS:teo}. The rest of the theorem is not difficult to prove and we can leave it as an exercise.

\begin{ex} Prove that $F$ is proper.
\end{ex}

We denote by $\Iso^+(M)$ the subgroup of $\Iso(M)$ consisting of all the orientation-preserving isometries. It has index one or two in $\Iso(M)$.

\subsection{Discrete groups}
Let $M$ be a Riemannian manifold. The group $\Iso(M)$ is a Lie group, so it makes sense to consider discrete subgroups.

\begin{prop} \label{iff:discrete:prop}
A group $\Gamma<\Iso(M)$ acts properly discontinuously on $M$ if and only if it is discrete.
\end{prop}
\begin{proof}
If $\Gamma$ is discrete, Theorem \ref{MS:teo} and Exercise \ref{proper:ex} imply that it acts properly discontinuously. Conversely, if $\Gamma$ acts properly discontinuously then $e\in \Gamma$ is easily seen to be isolated, and we apply Exercise \ref{isolated:ex}.
\end{proof}

If $\Gamma < \Iso(M)$ is discrete and acts freely, the quotient map $M \to M/_\Gamma$ is a covering. Moreover, the Riemannian structure projects from $M$ to $M/_\Gamma$.

\begin{prop} \label{quotient:prop}
Let $\Gamma <\Iso(M)$ act freely and properly discontinuously on $M$. There is a unique Riemannian structure on the manifold $M/_\Gamma$ such that the covering $\pi\colon M \to M/_\Gamma$ is a local isometry.
\end{prop}
\begin{proof}
Let $U\subset M/_\Gamma$ be a well-covered set: we have $\pi^{-1}(U) = \sqcup_{i\in I} U_i$ and $\pi|_{U_i}\colon U_i \to U$ is a homeomorphism. Pick an $i\in I$ and assign to $U$ the smooth and Riemannian structure of $U_i$ transported along $\pi$. The resulting structure on $U$ does not depend on $i$ since the open sets $U_i$ are related by isometries in $\Gamma$. We get a Riemannian structure on $M/_\Gamma$, determined by the fact that $\pi$ is a local isometry.
\end{proof}

\subsection{Measures} \label{measure:descends:subsection}
The case where a group $G$ acts by preserving some measure instead of a Riemannian metric is also very interesting, although quite different. This situation will occur only in Chapter \ref{automorfismi:chapter}.

Let a group $G$ act on a manifold $M$. A Borel measure $\mu$ on $M$ is \emph{$G$-invariant} if $\mu(S) = \mu(g(S))$ for every Borel set $S\subset M$. 

If $G$ acts on $M$ freely and properly discontinuously, then $M \to M/_G$ is a covering and every $G$-invariant measure $\mu$ on $M$ descends to a natural measure on the quotient $M/_G$ which we still indicate by $\mu$. The measure on $M/_G$ is defined as follows: for every well-covered open set $U\subset M/_G$ we have $\pi^{-1}(U) = \sqcup_{i\in I}U_i$ and we assign to $U$ the measure of $U_i$ for any $i\in I$. This assignment extends to a unique Borel measure on $M/_G$ by Proposition \ref{ricoprimento2:prop}.

Note that the measure on $M/_G$ is \emph{not} the push-forward of $\mu$, namely it is not true that $\mu(U) = \mu(\pi^{-1}(U))$.

\section{Homology}
The singular homology theory needed in this book is not very deep: all the homology groups of the manifolds that we consider are boringly determined by their fundamental group. The theory is quickly reviewed in this section, a standard introduction is Hatcher's \emph{Algebraic Topology} \cite{HAT}.

\subsection{Definition}
Let $X$ be a topological space and $R$ be a ring. The cases $R = \matZ$, $\matR$, or $\matZ/_{2\matZ}$ are typically the most interesting ones.

A \emph{singular $k$-simplex} is a continuous map $\alpha\colon \Delta_k \to X$ from the standard  $k$-dimensional simplex $\Delta_k$ into $X$. A \emph{$k$-chain} is an abstract linear combination
$$\lambda_1\alpha_1+\ldots +\lambda_h\alpha_h$$
of singular $k$-simplexes $\alpha_1,\ldots,\alpha_h$ with coefficients $\lambda_1,\ldots,\lambda_h \in R$. The set $C_k(X,R)$ of all $k$-chains is a $R$-module. There is a linear boundary map $\partial_k\colon C_k(X,R) \to C_{k-1}(X,R)$ such that $\partial_{k-1}\circ \partial_k=0$. The \emph{cycles} and \emph{boundaries} are the elements of the submodules  
$$Z_k(X,R) = \ker \partial_k, \quad  B_k(X,R) = \img \partial_{k+1}.$$
The $k$-th \emph{homology group} is the quotient
$$H_k(X,R) = Z_k(X,R)/_{B_k(X,R)}.$$
We sometimes omit $R$ and write $H_k(X)$ instead of $H_k(X,R)$. \index{homology group}

By taking the dual spaces $C^k(X, R) = \Hom(C_k(X,R), R)$ we define analogously the \emph{cohomology group} $H^k(X,R)$. (Co-)homology groups are also defined for pairs $(X,Y)$ with $Y\subset X$. When $X=U\cup V$ and $U,V$ are open we get the exact \emph{Mayer -- Vietoris sequence}:
\index{Mayer -- Vietoris sequence}
$$\ldots \longrightarrow H_{n+1}(X) \longrightarrow H_n(U\cap V) \longrightarrow H_n(U) \oplus H_n(V) \longrightarrow H_n(X) \longrightarrow\ldots$$
(Co-)homology groups are \emph{functorial} in the sense that continuous maps induce natural homomorphisms of groups. Homotopic maps induce the same homomorphisms.

If $X$ is path-connected, then $H_0(X,R)=R$ and there is a canonical homomorphism $$\pi_1(X) \longrightarrow H_1(X,\matZ).$$
The homomorphism is surjective and its kernel is generated by the commutators of $\pi_1(X)$: in other words $H_1(X,\matZ)$ is the \emph{abelianization} of $\pi_1(X)$. Concerning cohomology, we have
$$H^1(X,\matZ) = \Hom\big(H_1(X,\matZ),\matZ\big)  = \Hom(\pi_1(X), \matZ).$$
The relation between $H^i(X,\matZ)$ and $H_i(X,\matZ)$ for $i>1$ is not so immediate.

\subsection{Dualities}
Let $M$ be a compact oriented connected $n$-manifold with (possibly empty) boundary. The abelian group $H_k(M,\matZ)$ is finitely generated and hence decomposes as  \index{Betti number}
$$H_k(M,\matZ) \isom F_k \oplus T_k$$ 
where $F_k = \matZ^{b_k}$ is free and $T_k$ is finite. The \emph{torsion} subgroup $T_k$ consists of all finite-order elements in $H_k(M,\matZ)$. The rank $b_k$ of $F_k$ is the $k$-th \emph{Betti number} of $M$. In cohomology things change only a little:
$$H^k(M,\matZ) \isom F_k \oplus T_{k-1}.$$
All these groups vanish when $k>n$. Even when the torsion vanishes, there is no canonical isomorphism between $H_k(M)$ and $H^k(M)$. On the other hand,
the \emph{Lefschetz duality} provides a canonical identification  \index{Lefschetz duality}
$$H^k(M) = H_{n-k}(M,\partial M)$$
for any ring $R$. When $\partial M=\emptyset$ this is the \emph{Poincar\'e duality} $H^k(M) = H_{n-k}(M)$. In particular we get  \index{Poincar\'e duality}
$$H_n(M, \partial M,\matZ) = H^0(M,\matZ) \isom \matZ$$ 
and the choice of an orientation for $M$ is equivalent to a choice of a generator $[M] \in H_n(M,\partial M,\matZ)$ called the \emph{fundamental class} of $M$. \index{fundamental class}

An important exact sequence for $M$ is the following:
$$\ldots \longrightarrow H_n(M) \longrightarrow H_n(M,\partial M) \longrightarrow H_{n-1}(\partial M) \longrightarrow H_{n-1}(M) \longrightarrow \ldots$$

\subsection{Intersection form} \label{intersection:subsection}
Let $G$ and $H$ be finitely generated abelian groups, seen as $\matZ$-modules. A bilinear form
$$\omega\colon G\times H \longrightarrow \matZ$$
is \emph{non-degenerate} if for every infinite-order element $g\in G$ there is a  $h\in H$ such that $\omega(g,h)\neq 0$. If $G=H$, we say that $\omega$ is \emph{symmetric} (resp.~\emph{skew-symmetric}) if $\omega(g_1,g_2)$ equals $\omega(g_2,g_1)$ (resp.~$-\omega(g_2,g_1)$) for all $g_1,g_2\in G$. A skew-symmetric non-degenerate form is called \emph{symplectic}.  \index{symplectic form}

Let $M$ be a compact oriented connected $n$-manifold with (possibly empty) boundary. The Lefschetz duality furnishes a non-degenerate bilinear form
$$\omega\colon H_k(M,\matZ) \times H_{n-k}(M,\partial M, \matZ) \longrightarrow \matZ$$
called the \emph{intersection form}. It has the following geometric interpretation.
\index{intersection form in homology}

An oriented closed $k$-submanifold $S\subset M$ defines a class $[S] \in H_k(M)$ as the image of its fundamental class via the map $i_*\colon H_k(S) \to H_k(M)$ induced by the inclusion. If $S$ has boundary and is properly embedded (that is, $\partial S = \partial M \cap S)$, it defines a class $[S] \in H_k(M,\partial M)$.

Suppose two oriented submanifolds $S$ and $S'$ have complementary dimensions $k$ and $n-k$ and intersect transversely: every intersection point $x$ is isolated and has a sign $\pm 1$, defined by comparing the orientations of $T_x S \oplus T_xS'$ and $T_x M$. The \emph{algebraic intersection} $S\cdot S'$ of $S$ and $S'$ is the sum of these signs.  \index{algebraic intersection}

\begin{teo} \label{pairing:teo}
Let $S, S'$ be transverse and represent two classes $[S]\in H_k(M)$ and $[S']\in H_{n-k}(M,\partial M)$. We have
$$\omega ([S], [S']) = S\cdot S'.$$
\end{teo}

\begin{cor} The intersection number of two transverse oriented submanifolds of complementary dimension depends only on their homology classes.
\end{cor}

When $M$ is closed and has even dimension $2n$, the central form
$$\omega\colon H_n(M,\matZ) \times H_n(M, \matZ) \longrightarrow \matZ$$
is symmetric or skew-symmetric, depending on whether $n$ is even or odd. 

Everything we said about fundamental classes, Lefschetz duality, and intersection forms holds for non-orientable manifolds as well, provided that we pick $R = \matZ/_{2\matZ}$ and consider bilinear forms and intersections in $\matZ/_{2\matZ}$.

\section{Cells and handle decompositions} 
Many nice topological spaces can be constructed iteratively starting from finitely many points and then attaching some discs of increasing dimension called \emph{cells}. When the topological space is a differentiable manifold one typically thickens the cells to \emph{handles}. 

Cells and handles are beautifully introduced by Hatcher \cite{HAT} and Kosinksi \cite{K}, respectively.

\subsection{Cell complexes}
A \emph{finite cell complex} of dimension $n$ (briefly, a $n$-complex) is a topological space obtained iteratively in the following manner:
\begin{itemize}
\item a $0$-complex $X^0$ is a finite set of points, 
\item a $n$-complex $X^n$ is obtained from a $(n-1)$-complex $X^{n-1}$ by attaching finitely many \emph{$n$-cells}, that is copies of $D^n$ glued along continuous maps $\varphi\colon \partial D^n \to X^{n-1}$.
\end{itemize}
The closed subset $X^k\subset X^n$ is the \emph{$k$-skeleton} of $X^n$, for all $k<n$.  \index{cell complex}

\begin{prop}
Let $X$ be a $n$-complex. The inclusion map $X^k \hookrightarrow X$ induces an isomorphism $\pi_{j}(X^k) \to \pi_j(X)$ for all $j<k$.
\end{prop}
\begin{proof}
By transversality, maps $S^j\to X$ and homotopies between them can be homotoped away from cells of dimension $\geqslant j+2$. 
\end{proof}

In particular, the space $X$ is connected if and only if $X^1$ is, and its fundamental group is captured by $X^2$.

\begin{teo} Every differentiable compact $n$-manifold may be realised topologically as a finite $n$-complex.
\end{teo}

A presentation for the fundamental group of a cell complex $X$ can be constructed as follows.
If $x_0\in X^0$, we fix a maximal tree $T\subset X^1$ containing $x_0$ and equip the $k$ arcs in $X^1\setminus T$ with some arbitrary orientations. These arcs determine some generators $g_1,\ldots, g_k\in\pi_1(X,x_0)$. The boundary of a $2$-cell makes a circular path in $X^1$: every time it crosses an arc $g_i$ in some direction (entering from one side and exiting from the other) we write the corresponding letter $g_i^{\pm 1}$ and get a word. The $s$ two-cells produce $s$ word relations. We have constructed a presentation for $\pi_1(X)$.

\subsection{Euler characteristic}
The \emph{Euler characteristic} of a $n$-complex $X$ is the integer  \index{Euler characteristic}
$$\chi (X) = \sum_{i=0}^n (-1)^i C_i$$
where $C_i$ is the number of $i$-cells in $X$. It is also equal to 
$$\chi(X) = \sum_{i=0}^n (-1)^i b_i(X)$$
where $b_i(X)$ is the $i$-th Betti number of $X$. Therefore $\chi(X)$ is a number which depends only on the topology of $X$, that can be easily calculated from any cell decomposition of $X$.

\begin{prop} If $\tilde X \to X$ is a degree-$d$ covering of finite complexes then $\chi(\tilde X) = d\cdot \chi(X)$.
\end{prop}
\begin{proof}
A $k$-cell in $X$ is simply-connected and hence lifts to $d$ distinct $k$-cells in $\tilde X$. A cell decomposition of $X$ thus induces one of $\tilde X$ where the numbers $C_i$ are all multiplied by $d$.
\end{proof}

It is much harder to control homology under coverings. The Euler characteristic of a closed manifold measures the obstruction of constructing a nowhere-vanishing vector field.

\begin{teo} A closed differentiable manifold $M$ has a nowhere-vanishing vector field if and only if $\chi(M)=0$.
\end{teo}

\subsection{Aspherical cell-complexes}
A finite cell complex is locally contractible and hence has a universal cover $\tilde X$. If $\tilde X$ is contractible the complex  $X$ is called \emph{aspherical}. \index{cell complex!aspherical cell complex}

\begin{example}
If a closed manifold $M$ is covered by $\matR^n$, it is aspherical. 
\end{example}

By a theorem of Whitehead, a finite cell complex $X$ is aspherical if and only if all its higher homotopy groups $\pi_i(X)$ with $i\geqslant 2$ vanish.

The following theorem says that maps to aspherical spaces are determined (up to homotopy) by homomorphisms between fundamental groups.

\begin{teo} \label{complessi:omotopia:teo}
Let $(X,x_0), (Y,y_0)$ be pointed connected finite cellular complexes. If $Y$ is aspherical, every homomorphism $\pi_1(X,x_0) \to \pi_1(Y,y_0)$ is induced by a continuous map $(X,x_0)\to (Y,y_0)$, unique up to homotopy.
\end{teo}
\begin{proof}
We construct a continuous map $f\colon X^k \to Y^k$ iteratively on the $k$-skeleta, starting from $k=1$. 

Let $T$ be a maximal tree in $X^1$. The oriented 1-cells $g_1,\ldots, g_k$ in $X^1\setminus T$ define generators in $\pi_1(X,x_0)$. We define $f\colon X^1 \to Y^1$ by sending $T$ to $y_0$ and each $g_i$ to any loop in $Y$ representing the image of $g_i$ along the given homomorphism $\pi_1(X,x_0) \to \pi_1(Y,y_0)$.

The map $f$ sends the boundary of each 2-cell to a homotopically trivial loop in $Y$ and hence extends to a map $f\colon X^2 \to Y^2$. Since $Y$ is aspherical, the higher homotopy groups $\pi_i(Y)$ with $i\geqslant 2$ vanish and hence $f$ extends to a map $f\colon X^k \to Y^k$ iteratively for all $k\geqslant 3$.

We prove that $f$ is unique up to homotopy. Let $f'$ be a map that realises the given homomorphism on fundamental groups. We can construct a homotopy between $f$ and $f'$ iteratively on $X^k$ as follows. 

For $k=1$, we can suppose that both $f$ and $f'$ send $T$ to $y_0$. By hypothesis they send the generators $g_i$ to homotopic loops, hence we can homotope $f'$ to $f$ on $X^1$. For $k\geqslant 2$, the maps $f$ and $f'$ on each $k$-cell are homotopic because they glue to a map $S^k \to Y$, which is null-homotopic because $\pi_k(Y)$ is trivial.
\end{proof}

\begin{cor} \label{aspherical:cor}
Let $X$ and $Y$ be connected finite aspherical complexes. Every isomorphism $\pi_1(X) \to \pi_1(Y)$ is realised by a homotopic equivalence $X \to Y$, unique up to homotopy.
\end{cor}

In particular the homotopy type of an aspherical manifold is fully determined by its fundamental group.

\begin{cor}
Two aspherical closed manifolds of distinct dimension have non-isomorphic fundamental groups.
\end{cor}
\begin{proof}
Two closed manifolds of different dimension cannot be homotopically equivalent because they have non-isomorphic homology groups.
\end{proof}

We cite for completeness this result, although we will never use it. \index{Cartan--Hadamard theorem}

\begin{teo}[Cartan--Hadamard] A complete Riemannian manifold $M$ with sectional curvature everywhere $\leqslant 0$ has a universal covering diffeomorphic to $\matR^n$ and is hence aspherical.
\end{teo}
\begin{proof}[Sketch of the proof]
Pick a point $x\in M$. Since $M$ is complete, the exponential map $\exp_x\colon T_xM \to M$ is defined on $T_xM$. The fact that the sectional curvatures are $\leqslant 0$ imply that $(d\exp_x)_y$ is invertible for any $y\in T_x$ and $\exp_x$ is a covering.
\end{proof}

At a single point in this book we will need the following theorem.
\begin{teo} \label{no:torsion:aspherical:teo}
The fundamental group of an aspherical manifold $M$ has no torsion.
\end{teo}
\begin{proof}[Sketch of the proof]
Up to passing to a cover, it suffices to consider the case $\pi_1(M) \isom \matZ/_{n\matZ}$ for some $n\geqslant 2$. This case is excluded because the cohomology of $M$ is isomorphic to the (suitably defined) cohomology of $\matZ/_{n\matZ}$, which has however infinite dimension.
\end{proof}

\subsection{Gluing portion of boundaries}
Every compact differentiable $n$-manifold $M$ can be obtained topologically as a finite complex. The finite complex structure is however not designed to describe the smooth structure of $M$, and for that purpose it is better to replace $k$-cells with some thickened objects called \emph{$k$-handles}. These handles are $n$-discs glued iteratively along portions of their boundaries. Before describing them, we briefly explain how smooth manifolds can be glued along portions of their boundaries.

Let $M$ and $N$ be two $n$-manifolds with boundary and $X\subset \partial M$, $Y\subset \partial N$ be two compact $(n-1)$-submanifolds with boundary. A diffeomorphism $\varphi\colon X\to Y$ defines a topological space
$$M \cup_{\varphi} N$$
that may be promoted to a differentiable manifold: it suffices to use a collar as in Section \ref{cut:subsection}, see Figure \ref{boundary_glue:fig}.
The glued manifold depends (up to diffeomorphism) only on the isotopy class of $\varphi$.

\begin{figure}
\begin{center}
\includegraphics[width = 12.5 cm] {\iftoggle{BW}{boundary_glue-BW}{boundary_glue}}
\nota{To glue two smooth manifolds along portions of their boundaries (left) we pick a collar for these portions (centre), then we remove $X,Y$ and identify the interiors of the collars (right).}
\label{boundary_glue:fig}
\end{center}
\end{figure}

\subsection{Handles} \label{handles:subsection}
Let $M$ be a (possibly empty or disconnected) $n$-manifold with boundary and $0\leqslant k \leqslant n$. A \emph{$k$-handle} is a manifold $D^k\times D^{n-k}$ attached to $M$ along some diffeomorphism $\varphi\colon \partial D^k \times D^{n-k}\to Y\subset \partial M$, hence producing a new manifold $M'$. 

\begin{figure}
\begin{center}
\includegraphics[width = 12 cm] {\iftoggle{BW}{handles-BW}{handles}}
\nota{Two 0-handles\iftoggle{BW}{}{ (yellow)}, two 1-handles\iftoggle{BW}{}{ (orange)}, and one 2-handle\iftoggle{BW}{}{ (red)} in dimension 2 (left). Two 0-handles\iftoggle{BW}{}{ (yellow)} and one 1-handle\iftoggle{BW}{}{ (orange)} in dimension 3 (right).}
\label{handles:fig}
\end{center}
\end{figure}

For instance, a $0$-handle is a $D^0 \times D^n = D^n$ attached to nothing, since $\partial D^0 = \varnothing$. This means that attaching a 0-handle to $M$ consists of adding a disjoint disc $D^n$ to it. For instance, by attaching a 0-handle to the empty set we create a disc $D^n$ out of nothing. 

A $1$-handle is a $D^1\times D^{n-1}$ attached along $\partial D^1 \times D^{n-1} = S^0 \times D^{n-1}$, that is two copies of $D^{n-1}$. Some examples in dimension $n=2,3$ are shown in Figure \ref{handles:fig}. A $2$-handle is a $D^2\times D^{n-2}$ attached along $\partial D^2 \times D^{n-2} = S^1\times D^{n-2}$. When $n=2$ this is a disc attached along its boundary, see Figure \ref{handles:fig}-(left).

\subsection{Handle decompositions}
A sequence of handle attachments \index{handle decomposition}
$$\emptyset \rightsquigarrow M_1 \rightsquigarrow \ldots \rightsquigarrow M_k = M$$
starting from the empty set and producing a compact manifold $M$ with (possibly empty) boundary is called a \emph{handle decomposition} for $M$. Using Morse theory one proves the following.

\begin{teo} Every compact manifold (possibly with boundary) can be obtained from a handle decomposition.
\end{teo}

\begin{example}
The disc $D^n$ has an obvious handle decomposition consisting of a single $0$-handle. It also has more complicate handle decompositions, as Figure \ref{handles:fig} shows.
\end{example}

By transversality, handles may always be reordered so that the lower index handles are attached first, and handles of the same index are attached simultaneously. So we can think of a decomposition as the appearing of some 0-handles, then the simultaneous attaching of some 1-handles, then of some 2-handles, and so on. 

A handle decomposition of a closed manifold may be turned upside down, by reversing all arrows and interpreting every $k$-handle as a $(n-k)$-handle. 

\begin{prop} \label{0:handle:prop}
Every compact connected manifold $M$ has a handle decomposition with one $0$-handle and at most one $n$-handle.
\end{prop}

If $M$ has a handle decomposition with $n_i$ handles of index $i$ then
$$\chi(M) = \sum_{i=0}^n (-1)^i n_i.$$

\subsection{Triangulations} \label{triangulations:subsection}
\begin{figure}
\begin{center}
\includegraphics[width = 7 cm] {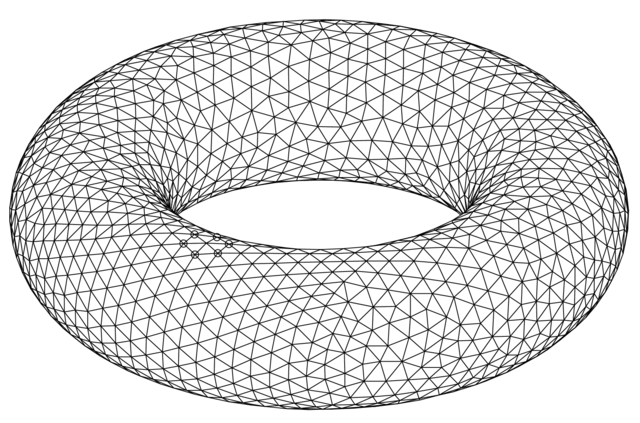}
\nota{A triangulated torus.}
\label{Torus-triang:fig}
\end{center}
\end{figure}

Instead of decomposing a manifold into handles, one may decide to decompose it into simplexes as in Figure \ref{Torus-triang:fig}. One such decomposition is called a triangulation: we now give a formal definition.

An (abstract and finite) \emph{simplicial complex} is a finite set $S$ of \emph{vertices} together with a set $X$ of subsets of $S$, called \emph{faces}, which contains all singletons and such that every subset of a face is also a face. A face of cardinality $k+1$ has \emph{dimension} $k$. \index{simplicial complex}

A simplicial complex $X$ has a topological realisation $|X|$, in which every face $F\in X$ of dimension $k$ transforms into a $k$-simplex with vertices in $F$ (we get a cell complex). The complex $X$ is \emph{purely $n$-dimensional} if every face is contained in a face of dimension $n$. This implies that $|X|$ is the union of its $n$-simplexes.

Let now $M$ be a differentiable $n$-manifold. A \emph{smooth $k$-simplex} in $M$ is a subset $\Delta \subset M$ diffeomorphic to a $k$-simplex in $\matR^n$ (we mean that there is a diffeomorphism between open neighbourhoods of $\Delta$ and of the simplex sending the first to the second). A \emph{smooth triangulation} of $M$ is a pure $n$-dimensional simplicial complex $X$ together with a homeomorphism between $|X|$ and $M$ that transforms every $k$-simplex in $X$ into a smooth $k$-simplex in $M$. \index{triangulation}

\begin{teo} Every compact manifold with (possibly empty) boundary has a smooth triangulation.
\end{teo}

Handle decompositions are often more flexible and efficient than triangulations, but triangulations have a more combinatorial nature and can for instance be handled by a computer. Both constructions are widely used.

On a simplicial complex $X$, the \emph{star} of a vertex $v$ is the set of all simplexes containing $v$ plus all their sub-simplexes; the \emph{link} of $v$ is the subset of the star consisting of the simplexes that do not contain $v$. Topologically, the star of $v$ is a closed neighbourhood of $v$ homeomorphic to a cone on its link.
\index{star of a vertex} \index{link of a vertex}

\subsection{Hypersurfaces and cohomology}
The techniques introduced in this chapter furnish another characterisation of the first cohomology group of a manifold. Let $[X,Y]$ denote the homotopy classes of continuous maps $X \to Y$.

\begin{prop}
Let $X$ be a cell complex. There is a canonical identification
$$H^1(X;\matZ) = [X,S^1].$$
\end{prop}
\begin{proof}
We know that $H^1(X;\matZ) = \Hom(\pi_1(X),\matZ)$.
The circle $S^1$ is aspherical with $\pi_1(S^1)=\matZ$, so by Theorem \ref{complessi:omotopia:teo} every homomorphism $\pi_1(X) \to \pi_1(S^1)$ is realised by a map $X\to S^1$, unique up to homotopy.
\end{proof}

We have noted that an oriented $k$-submanifold $S\subset M$ determines a class $[S]\in H_k(M,\matZ)$. It is natural to ask whether every homology class may be represented by an oriented submanifold: this is false in general, but it is true in codimension one.

\begin{prop} \label{homology:prop}
Let $M$ be a compact oriented $n$-manifold with (possibly empty) boundary. Every class in $H_{n-1}(M,\partial M;\matZ)$ is represented by an oriented properly embedded hypersurface $S\subset M$.
\end{prop}
\begin{proof}
We have
$$H_{n-1}(M,\partial M;\matZ) = H^1(M;\matZ) = \Hom (\pi_1(M), \matZ) = [M, S^1].$$
Every map $M\to S^1$ is homotopic to a smooth map, hence each class $\alpha \in H^1(M,\matZ)$ is represented by a smooth map $f\colon M\to S^1$. By Sard's lemma there is a regular value $x\in S^1$, whose counter-image $S = f^{-1}(x)$ is a hypersurface, transversely oriented by $f$ and hence oriented (since $M$ is).

To prove that $[S] = \alpha$ we verify that $[S]$ and $\alpha$ act on $\pi_1(M)$ in the same way. Pick a loop $\gamma \in \pi_1(M)$ transverse to $S$. The number $\alpha(\gamma)$ is the degree of $f\circ\gamma$, that is the number of times $f\circ \gamma$ crosses $x$ counted with signs, which equals the number of times $\gamma$ crosses $S$ counted with signs.
\end{proof}

\begin{prop}
An oriented connected properly embedded hypersurface $S\subset M$ is separating if and only if $[S]\in H_{n-1}(M,\partial M; \matZ)$ is trivial.
\end{prop}
\begin{proof}
The surface $S$ is separating if and only if $S \cdot \alpha = 0$ for every loop $\alpha$ transverse to $S$ (exercise). Since $\omega([S], [\alpha]) = S \cdot \alpha$ and $\omega$ is non-degenerate (see Section \ref{intersection:subsection}), this holds precisely when $[S]=0$.
\end{proof}

\subsection{Topological discs}
Finally, we mention a purely topological theorem that will be used only at one point in this book, in Chapter \ref{automorfismi:chapter}.

\begin{teo} \label{topological:disc:teo}
Let $M$ be a compact topological manifold, whose interior and boundary are homeomorphic respectively to $B^n$ and $S^{n-1}$. Then $M$ is homeomorphic to $D^n$.
\end{teo}
\begin{proof}
We use two important theorems on topological manifolds, whose proofs can be found in the first chapter of the book \emph{Topological embeddings} of Rushing \cite{Ru}.

As in every compact topological manifold, the boundary $\partial M$ has a topological collar by a theorem of Brown (see Theorem 1.7.4 in that book). The interior boundary sphere of the collar is contained in $\interior M \cong B^n$ and hence bounds a closed disc by the Generalised Sch\"onflies Theorem (see Theorem 1.8.2 there). The manifold $M$ is obtained by collaring a closed disc and is hence a closed disc.
\end{proof}

%% file: Space.tex
\chapter{Hyperbolic space} \label{spazio:chapter} \label{space:chapter}
In every dimension $n\geqslant 2$ there exists a unique simply connected complete Riemannian manifold with constant sectional curvature $1$, $0$, or $-1$. These are the sphere $S^n$, the Euclidean space $\matR^n$, and the hyperbolic space $\matH^n$. 

These manifolds are the three most important spaces in Riemannian geometry. We introduce in this chapter the least familiar and the most interesting of the three: hyperbolic space. 

\section{The models of hyperbolic space}
In contrast with $S^n$ and $\matR^n$, the hyperbolic space $\matH^n$ may be described in various different ways, no-one of which is prevalent in the literature. Each description is a \emph{model} for $\matH^n$. We first introduce the hyperboloid model, which has a more algebraic flavour, and then we turn to the disc and half-space models that are somehow more geometric (and easier to visualise in dimensions $n=2$ and $3$).

\subsection{Hyperboloid}
The sphere $S^n$ is the set of all points with norm 1 in $\matR^{n+1}$, equipped with the Euclidean scalar product. Analogously, we may define $\matH^n$ as the set of all points of norm $-1$ in $\matR^{n+1}$, equipped with the usual \emph{Lorentzian} scalar product. This set forms a hyperboloid with two sheets, and we choose one. 
\index{Lorentzian scalar product}

\begin{defn} The \emph{Lorentzian scalar product} on $\matR^{n+1}$ is:
$$\langle x, y \rangle = \sum_{i=1}^n x_iy_i - x_{n+1}y_{n+1}.$$
It has signature $(n,1)$. A vector $x \in \matR^{n+1}$ is \emph{time-like}, \emph{light-like}, or \emph{space-like} if $\langle x, x \rangle$ is negative, null, or positive respectively.
The \emph{hyperboloid model} $I^n$ is defined as follows:
\index{vector!time-like, light-like, and space-like vector} \index{hyperbolic space} \index{hyperbolic space!hyperboloid model}
$$I^n = \big\{x \in \matR^{n+1}\ \big| \ \langle x, x \rangle = -1, \ x_{n+1} > 0 \big\}.$$
\end{defn}

\begin{figure}
\begin{center}
\includegraphics[width = 6 cm] {\iftoggle{BW}{Hyperboloid2-BW}{640px-Hyperboloid2}}
\nota{The hyperboloid with two sheets defined by the equation $\langle x, x \rangle = -1$. The model $I^n$ for $\matH^n$ is the upper sheet.}
\label{Hyperboloid2:fig}
\end{center}
\end{figure}

The set of points $x$ with $\langle x,x \rangle = -1$ is a \emph{hyperboloid with two sheets}, and $I^n$ is the connected component (sheet) with $x_{n+1}>0$. Let us prove a general fact. For us, a scalar product is a real non-degenerate symmetric bilinear form.\index{$I^n$}

\begin{prop}
Let $\langle, \rangle$ be a scalar product on $\matR^{n+1}$. The function $f\colon \matR^n \to \matR$ given by
$$f(x) = \langle x,x\rangle$$
is everywhere smooth and has differential  
$$df_x(y) = 2\langle x,y\rangle.$$
\end{prop}
\begin{proof}
The following equality holds:
$$\langle x+y, x+y\rangle = \langle x,x\rangle + 2 \langle x,y \rangle + \langle y, y \rangle.$$
The component $\langle x,y\rangle $ is linear in $y$ while $\langle y,y \rangle$ is quadratic. 
\end{proof}
\begin{cor}
The hyperboloid $I^n$ is a Riemannian manifold.
\end{cor}
\begin{proof}
The hyperboloid is the set of points with $f(x)=\langle x,x\rangle = -1$. For all $x\in I^n$ the differential $y\mapsto 2\langle x,y\rangle $ is surjective and hence the hyperboloid is a differentiable submanifold of codimension 1. 

The tangent space $T_xI^n$ at $x\in I^n$ is the hyperplane
$$T_x = \ker df_x = \big\{ y \ \big| \ \langle x, y \rangle = 0 \big\} = x^\bot$$
orthogonal to $x$ in the Lorentzian scalar product. Since $x$ is time-like, the restriction of the Lorentzian scalar product to $x^\bot$ is positive definite and hence defines a metric tensor on $I^n$. \end{proof}

The hyperboloid $I^n$ is a model for hyperbolic space $\matH^n$. We will soon prove that it is simply connected, complete, and has constant curvature $-1$. 

\subsection{Isometries of the hyperboloid}
The isometries of $I^n$ are easily classified using linear algebra. 

Let $O(n,1)$ be the group of linear isomorphisms $f$ of $\matR^{n+1}$ that preserve the Lorentzian scalar product, \emph{i.e.}~such that $\langle v,w\rangle = \langle f(v),  f(w) \rangle$ for any $v,w \in \matR^n$. An element in $O(n,1)$ preserves the hyperboloid with two sheets, and the elements preserving the upper sheet $I^n$ form a subgroup of index two in $O(n,1)$ that we indicate with $O^+(n,1)$.\index{$O(n,1)$}\index{$O^+(n,1)$}

\begin{prop} \label{isom:prop}
The following equality holds:
$$\Iso (I^n) = O^+(n,1).$$
\end{prop}
\begin{proof}
Pick $f\in O^+(n,1)$. If $x\in I^n$ then $f(x) \in I^n$ and $f$ sends $x^\bot$ to $f(x)^\bot$ isometrically, hence $f \in\Iso (I^n)$. Therefore $O^+(n,1) \subseteq \Iso (I^n)$.

To prove the converse inclusion we show that for every pair $x,y \in I^n$ and every
linear isometry $g\colon x^\bot \to y^\bot$ there is an element $f\in O^+(n,1)$ such that $f(x)=y$ and $f|_{x^\bot} = g$. Since isometries are determined by their first-order behaviour at a point $x$, this implies $\Iso(I^n) \subseteq O^+(n,1)$.

Via elementary linear algebra we prove that $O^+(n,1)$ acts transitively on $I^n$ and hence we may suppose that $x=y=(0,\ldots,0,1)$. Now $x^\bot = y^\bot$ is the horizontal hyperplane and $g\in O(n)$.  To define $f$ simply take
$$f = \begin{pmatrix} g & 0 \\ 0 & 1 \end{pmatrix}.$$
The proof is complete.
\end{proof}

The isometry groups of $S^n$ and $\matR^n$ are described analogously:
\begin{prop} \label{isom2:prop}
The following equalities hold:
\begin{align*}
\Iso(S^n) & = O(n+1), \\
\Iso(\matR^n) & = \big\{x \mapsto Ax+b\ \big| \ A\in O(n), b \in \matR^n\big\}.
\end{align*}
\end{prop}
\begin{proof}
The proof is analogous to the one above.
\end{proof}

We have also proved the following fact. A \emph{frame} at a point $p$ in a Riemannian manifold $M$ is an orthonormal basis for $T_pM$.

\begin{cor} \label{riferimento:cor}
Let $M= S^n$, $\matR^n$, or $\matH^n$. Given two points $p,q \in M$ and two frames at $p$ and $q$, there is a unique isometry that carries the first frame to the second.
\end{cor}

\subsection{Subspaces}
We introduce the following natural objects. \index{hyperbolic space!subspace of hyperbolic space}

\begin{defn}
A $k$-dimensional \emph{subspace} of $\matR^n$, $S^n$, $I^n$ is respectively:
\begin{itemize}
\item an affine $k$-dimensional space in $\matR^n$,
\item the intersection of a $(k+1)$-dimensional vector subspace of $\matR^{n+1}$ with $S^n$,
\item the intersection of a $(k+1)$-dimensional vector subspace of $\matR^{n+1}$ with $I^n$, when it is non-empty.
\end{itemize}
\end{defn}
\begin{oss}
Elementary linear algebra shows that the following conditions are equivalent for any $(k+1)$-dimensional vector subspace $W\subset\matR^{n+1}$:
\begin{enumerate}
\item $W\cap I^n \neq \emptyset$,
\item $W$ contains at least a time-like vector,
\item the signature of $\langle,\rangle|_W$ is $(k,1)$.
\end{enumerate}
\end{oss}

A $k$-subspace in $\matR^n, S^n, \matH^n$ is itself isometric to $\matR^k, S^k, \matH^k$. 
The non-empty intersection of two subspaces is always a subspace. An isometry of $\matR^n, S^n, \matH^n$ sends $k$-subspaces to $k$-subspaces. 

\begin{ex}
Let $S$ be a $k$-subspace in $\matR^n, S^n$, or $\matH^n$ and $p\in S$ a point. There is a unique $(n-k)$-subspace $S'$ intersecting $S$ orthogonally in $p$.
\end{ex}

\subsection{Reflections} \label{reflections:subsection}
We now introduce a basic kind of isometry.
The \emph{reflection} $r_S$ along a subspace $S$ in $I^n$ is an isometry of $I^n$ defined as follows. By definition $S=I^n \cap W$ with $\langle,\rangle|_W$ non-degenerate, hence $\matR^{n+1} = W \oplus W^\bot$ and we set $r_S|_W = \id$ and $r_S|_{W^\bot} = -\id$.
Analogous definitions work for subspaces of $S^n$ and $\matR^n$. 
\index{reflection}

\begin{ex}
The reflection $r_S$ has fixed set $S$ and preserves all the subspaces orthogonal to $S$. It is orientation-preserving if and only if $S$ has even codimension.
\end{ex}

\begin{prop} \label{reflections:generate:prop}
Reflections along hyperplanes generate the isometry groups of $S^n$, $\matR^n$, and $\matH^n$.
\end{prop}
\begin{proof}
It is a standard linear algebra fact that orthogonal reflections along vector hyperspaces
generate $O(n)$. This proves the case $S^n$ and shows that reflections generate the stabiliser of any point in $\matR^n$ and $\matH^n$. To conclude it suffices to check that reflections act transitively on points: to send $x$ to $y$, reflect along the hyperplane orthogonal to the segment connecting $x$ to $y$ in its midpoint.
\end{proof}

\subsection{Lines}
A $1$-subspace is a \emph{line}. We now show that lines and geodesics are the same thing. We recall the hyperbolic trigonometric functions: \index{line}
$$\sinh (t) = \frac{e^t-e^{-t}}2, \quad \cosh (t) = \frac{e^t+e^{-t}}2.$$

\begin{prop} \label{lines:prop}
A non-trivial complete geodesic in $S^n$, $\matR^n$, or $\matH^n$ is a line run at constant speed. Concretely, let $p\in M$ be a point and $v\in T_pM$ a unit vector. The geodesic $\gamma$ exiting from $p$ with velocity $v$ is:
\begin{itemize}
\item $\gamma(t) = \cos (t) \cdot p + \sin (t) \cdot v$ if $M=S^n$,
\item $\gamma(t) = p + tv$ if $M=\matR^n$,
\item $\gamma(t) = \cosh (t) \cdot p + \sinh (t) \cdot v$ if $M=I^n$.
\end{itemize}
\end{prop}
\begin{proof}
Let $p\in I^n$ be a point, $v\in T_pM$ a unit vector, and $\gamma$ the geodesic exiting from $p$ with velocity $v$.
The plane $W\subset \matR^{n+1}$ generated by $p$ and $v$ intersects $I^n$ into a line $l$ containing $p$ and tangent to $v$. To prove that $l$ is the support of $\gamma$ we use a symmetry argument: the reflection $r_l$ fixes $p$ and $v$ and hence $\gamma$, therefore $\gamma$ is forced to be contained in its fixed locus, which is $l$. This shows that non-trivial geodesics are lines run at constant speed.

We now consider the curve $\alpha(t) = \cosh (t) \cdot p + \sinh(t) \cdot v$. We have $\alpha(0) = p$ and $\alpha'(0) = v$. It remains to prove that $\alpha$ parametrizes $l$ with unit speed, and from this we deduce that $\gamma = \alpha$. We note that 
\begin{align*}
\langle \alpha (t), \alpha(t) \rangle & = \cosh^2(t) \langle p,p \rangle + 2\cosh(t)\sinh(t) \langle p,v \rangle + \sinh^2(t) \langle v,v\rangle \\
 & = -\cosh^2(t) + \sinh^2(t) = -1.
 \end{align*}
Therefore $\alpha$ parametrizes $l$. Its velocity is 
$$\alpha'(t) = \cosh'(t)\cdot p + \sinh'(t)\cdot v = \sinh(t)\cdot p + \cosh(t) \cdot v$$
whose squared norm is $-\sinh^2(t) + \cosh^2(t) = 1$. Therefore $\gamma=\alpha$. The proofs for $S^n$ and $\matR^n$ are analogous.
\end{proof}

\begin{cor}
The spaces $S^n$, $\matR^n$, and $\matH^n$ are complete.
\end{cor}
\begin{proof}
The previous proposition shows that geodesics are defined on $\matR$, hence the space is complete by the Hopf -- Rinow Theorem \ref{HR:teo}.
\end{proof}

It is easy to show that two points in $\matH^n$ are contained in a unique line.
\index{Euclid's V postulate}
\begin{oss}
Euclid's fifth postulate holds only in $\matR^2$. Given a line $r$ and a point $P\not\in r$, there is exactly one line passing through $P$ and disjoint from $r$ in $\matR^2$, there is no-one in $S^2$, and there are infinitely many in $\matH^2$.
\end{oss}

We can easily calculate the distance between two points.

\begin{prop}
Let $p, q\in M$ be two points. We have
\begin{itemize}
\item $\cos (d(p,q)) = \langle p,q \rangle $ if $M=S^n$,
\item $\cosh (d(p,q)) = -\langle p,q \rangle $ if $M=I^n$.
\end{itemize}
\end{prop}
\begin{proof}
Let $\gamma$ be a geodesic that goes from $p$ to $q$ at unit speed. We have $\gamma(0)=p$ and $\gamma(t)=q$ with $t=d(p,q)$. Set $v=\gamma'(0)$ and use Proposition \ref{lines:prop} to conclude.
\end{proof}

\subsection{The Poincar\'e disc} \label{Poincare:subsection}
We introduce two models of $\matH^n$, the \emph{disc} and the \emph{half-space}, that are easier to visualise especially in the dimensions $n=2,3$ we are interested in. The first model is the \emph{Poincar\'e disc}\index{$D^n$}
$$D^n = \big\{x\in \matR^n\ \big| \ \|x\|<1 \big\}.$$

\begin{figure}
\begin{center}
\includegraphics[width = 5 cm] {\iftoggle{BW}{bigezione-BW}{bigezione}}
\nota{The projection towards $P=(0,\ldots,0,-1)$ induces a bijection between the hyperboloid model $I^n$ and the disc model $D^n$.}
\label{bigezione:fig}
\end{center}
\end{figure}

The metric tensor on $D^n$ is obviously not the Euclidean one of $\matR^n$, but instead is induced by a particular diffeomorphism between $I^n$ and $D^n$ that we construct now. We identify $\matR^n$ with the horizontal hyperplane $x_{n+1}=0$ in $\matR^{n+1}$ and note that the projection towards $P=(0,\ldots,0,-1)$ depicted in Figure \ref{bigezione:fig} induces a bijection between $I^n$ and the horizontal disc $D^n\subset \matR^n$. 
The projection $p$ may be written as:
$$p(x_1,\ldots,x_{n+1}) = \frac {(x_1,\ldots, x_n)}{x_{n+1} +1}$$
and is indeed a diffeomorphism $p\colon I^n \to D^n$ that transports the metric tensor on $I^n$ to some metric tensor $g$ on $D^n$. \index{hyperbolic space!disc model}

\begin{prop}
The metric tensor $g$ at $x\in D^n$ is:
$$g_x = \left( \frac 2 {1-\|x\|^2}\right)^2 \cdot g^E_x$$
where $g^E$ is the Euclidean metric tensor on $D^n \subset \matR^n$.
\end{prop}
\begin{proof}
The inverse $q\colon D^n \to I^n$ of $p$ is:
$$q(x) = \left( \frac{2x_1}{1-\|x\|^2}, \ldots, \frac{2x_n}{1-\|x\|^2}, \frac{1+\|x\|^2}{1-\|x\|^2}\right).$$
Pick $x\in D^n$. Rotations around the $x_{n+1}$ axis are isometries of $I^n$ and commute with $p$, therefore they are isometries of $(D^n, g)$ too. Up to rotating we may take
$x = (x_1,0,\ldots, 0)$ and
find
$$dq_x = 
\frac {2}{1-x_1^2} \cdot
\begin{pmatrix}
\frac{1+x_1^2}{1-x_1^2} & 0 & \cdots & 0 \\
0 & 1 & \cdots & 0 \\
\vdots & \vdots & \ddots & \vdots \\
0 & 0 & \cdots & 1 \\
2\frac{x_1}{1-x_1^2} & 0 & \cdots & 0
\end{pmatrix}.
$$
The column vectors form an orthonormal basis of $T_{q(x)}I^n$. Hence $dq_x$  stretches all vectors of a constant $\frac 2{1-x_1^2}$. Therefore $g_x =\frac 4{(1-x_1^2)^2} g^E_x$. 
\end{proof}

The Poincar\'e disc is a \emph{conformal} model of $\matH^n$: it is a model where the metric differs from the Euclidean metric only by multiplication by a positive scalar $\big(\frac 2{1-\|x\|^2}\big)^2$ that depends smoothly on $x$. We note that the scalar tends to infinity when $x$ tends to $\partial D^n$. On a conformal model the lengths of vectors are different from the Euclidean lengths, but the angles formed by two adjacent vectors coincide with the Euclidean ones. Shortly: lengths are distorted but angles are preserved. \index{hyperbolic space!conformal model}

Let us see how we can easily visualise $k$-subspaces in the disc model.

\begin{prop}
The $k$-subspaces in $D^n$ are the intersections of $D^n$ with $k$-spheres and $k$-planes of $\matR^n$ orthogonal to $\partial D^n$.
\end{prop}
\begin{proof}
Since every $k$-subspace is an intersection of hyperplanes, we easily restrict to the case $k=n-1$. A hyperplane in $I^n$ is $I^n\cap v^\bot$ for some space-like vector $v$. If $v$ is horizontal (\emph{i.e.}~its last coordinate is zero) then $v^\bot$ is vertical and $p(I^n\cap v^\bot) = D^n \cap v^\bot$, a hyperplane orthogonal to $\partial D^n$.

If $v$ is not horizontal, up to rescaling and rotating around $x_{n+1}$ we may suppose $v=(\alpha,0,\ldots, 0, 1)$ with $\alpha>1$. The hyperplane is
$$\big\{x_1^2 +\ldots + x_n^2 - x_{n+1}^2 = -1 \big\} \cap \big\{x_{n+1} = \alpha x_1 \big\}.$$
On the other hand the sphere in $\matR^n$ of centre $(\alpha,0,\ldots,0)$ and radius $\sqrt{\alpha^2 -1}$ is orthogonal to $\partial D^n$ and is described by the equation
$$\left\{(y_1-\alpha)^2 +y_2^2 + \ldots + y_n^2 = \alpha^2 - 1\right\} =
\left\{y_1^2 + \ldots + y_n^2 - 2\alpha y_1 = -1 \right\}$$
which is equivalent to $||y||^2 = -1 + 2\alpha y_1$. If $y = p(x)$ the relations
$$y_1 = \frac {x_1}{x_{n+1}+1}, \quad \|y\|^2 = \frac{x_{n+1}-1}{x_{n+1}+1}$$
transform the latter equation in $x_{n+1}=\alpha x_1$.
\end{proof}

Three lines in $D^2$ delimiting a hyperbolic triangle are drawn in Figure \ref{Triangolo_iperbolico:fig}. Since the disc is a conformal model, the inner angles $\alpha, \beta, \gamma$ are the ones we see on the picture. In particular we verify easily that $\alpha+\beta+\gamma<\pi$. 

\begin{figure}
\begin{center}
\includegraphics[width = 4 cm] {\iftoggle{BW}{Triangolo_iperbolico-BW}{Triangolo_iperbolico}}
\nota{Three lines that determine a hyperbolic triangle in the Poincar\'e disc. The inner angles $\alpha, \beta,$ and $\gamma$ coincide with the Euclidean ones, and we have $\alpha + \beta + \gamma < \pi$.}
\label{Triangolo_iperbolico:fig}
\end{center}
\end{figure}

\begin{figure}
\begin{center}
\includegraphics[width = 5 cm] {\iftoggle{BW}{Order-5_tetrakis_square_tiling-BW}{Order-5_tetrakis_square_tiling}}
\nota{A \emph{tessellation} of $S^2$, $\matR^2$ o $\matH^2$ is a subdivision of the plane into polygons. The tessellation of $\matH^2$ shown here is obtained by drawing infinitely many lines in the plane. The triangles have inner angles $\frac \pi 2, \frac \pi {5}, \frac \pi {5}$ and are all isometric.}
\label{square_tiling:fig}
\end{center}
\end{figure}

\begin{ex} For any triple of positive angles $\alpha,\beta,\gamma$ with $\alpha+\beta+\gamma < \pi$ there is a triangle with inner angles $\alpha,\beta,\gamma$. This triangle is unique up to isometry.
\end{ex}

\subsection{Inversions} \label{inversions:subsection}
To construct our second conformal model we need to introduce a geometric transformation called inversion. \index{sphere inversion}

\begin{defn}
Let $S=S(x_0,r)$ be the sphere in $\matR^n$ centred in $x_0$ and with radius $r$. The \emph{inversion} along $S$ is the map $\varphi\colon \matR^n\setminus\{x_0\} \to \matR^n\setminus \{x_0\}$ defined as follows:
$$\varphi(x) = x_0 + r^2 \frac {x-x_0}{\|x-x_0\|^2}.$$
\end{defn}
The map may be extended continuously on the whole sphere $S^n$, identified with $\matR^n \cup \{\infty\}$ through the stereographic projection, by setting $\varphi(x_0) = \infty$ and $\varphi(\infty) = x_0$. A geometric description is shown in Figure \ref{Inversion_illustration1:fig}. 

We have already talked about conformal models. More generally, a diffeomorphism $f\colon M\to N$ between two oriented Riemannian manifolds is \emph{conformal} (respectively, \emph{anticonformal}) if for any $p\in M$ the differential $df_p$  
is the product of a scalar dilation $\lambda_p >0$ and an isometry that preserves (respectively, inverts) the orientation. 
The scalar $\lambda_p$ depends smoothly on $p$. A conformal map preserves the angle between two tangent vectors but multiplies their lengths by $\sqrt{\lambda_p}$.

\begin{ex} The stereographic projection is conformal.
\end{ex}

\begin{figure}
\begin{center}
\includegraphics[width = 4 cm] {\iftoggle{BW}{Inversion_illustration1-BW}{Inversion_illustration1}}
\includegraphics[width = 4 cm] {\iftoggle{BW}{Inversion_illustration2-BW}{Inversion_illustration2}}
\includegraphics[width = 4 cm] {\iftoggle{BW}{Inversion_illustration3-BW}{Inversion_illustration3}}
\nota{The inversion trough a sphere of centre $O$ and radius $r$ moves $P$ to $P'$ so that $OP \times OP' = r^2$ (left). It transforms a $k$-sphere $S$\iftoggle{BW}{}{ (blue)} into a $k$-plane\iftoggle{BW}{}{ (green)} if $O\in S$ (centre) or into a $k$-sphere\iftoggle{BW}{}{ (green)} if $O\not\in S$ (right).}
\label{Inversion_illustration1:fig}
\end{center}
\end{figure}

\begin{prop} \label{inversione:prop}
The following hold:
\begin{enumerate}
\item every inversion is a smooth and anticonformal map;
\item inversions send $k$-spheres and $k$-planes to $k$-spheres and $k$-planes.
\end{enumerate}
\end{prop}
\begin{proof}
Up to conjugating with translations and dilations we may suppose $x_0=0$ and $r=1$. The inversion is $\varphi(x) = \frac{x}{\|x\|^2}$ and we show that $d\varphi_x$ is $\frac 1{\|x\|^2}$ times a reflection with respect to the hyperplane orthogonal to $x$. We may suppose $x=(x_1,0,\ldots,0)$ and calculate the partial derivatives:
$$ 
\varphi(x_1,\ldots,x_n) = \frac{(x_1,\ldots, x_n)}{\|x\|^2}, \qquad
\frac{\partial \varphi_i}{\partial x_j} = \frac{\delta_{ij} \|x\|^2 - 2x_ix_j}{\|x\|^4}.
$$
The partial derivatives at $x=(x_1,0,\ldots, 0)$ are
$$\frac{\partial \varphi_1}{\partial x_1} = -\frac {1}{x_1^2}, \quad  \frac{\partial \varphi_i}{\partial x_i} = \frac {1}{x_1^2}, \quad \frac{\partial \varphi_j}{\partial x_k} = 0$$
for all $i>1$ and $j\neq k$. This proves our assertion.

The fact that an inversion preserves sphere and planes may be easily reduced to the two-dimensional case (with circles and lines), a classical fact of Euclidean geometry.
\end{proof}

\subsection{The half-space model}
We introduce another conformal model. The \emph{half-space model} is the space
\index{hyperbolic space!half-space model}\index{$H^n$}
$$H^n = \big\{(x_1,\ldots,x_n)\in \matR^n \ \big|\ x_n>0\big\}.$$
It is obtained from the disc model $D^n$ by an inversion in $\matR^n$ with centre $C=(0,\ldots,0,-1)$ and radius $\sqrt 2$ as shown in Figure \ref{inversion:fig}. The boundary $\partial H^n$ is the horizontal hyperplane $\{x_n=0\}$, to which we add the point $\infty$ at infinity to get a bijective correspondence between $\partial H^n$ and $\partial D^n$ through the inversion.

\begin{prop} The half-space $H^n$ is a conformal model for $\matH^n$. Its $k$-subspaces are the $k$-planes and $k$-spheres in $\matR^n$ orthogonal to $\partial H^n$. 
\end{prop}
\begin{proof}
The inversion is anticonformal and hence preserve angles, moreover it transforms $k$-spheres and $k$-planes in $D^n$ orthogonal to $\partial D^n$ into $k$-spheres and $k$-planes in $H^n$ orthogonal to $\partial H^n$.
\end{proof}

\begin{figure}
\begin{center}
\includegraphics[width = 5 cm] {\iftoggle{BW}{inversion-BW}{inversion}}
\nota{The inversion along the sphere with centre $C=(0,\ldots,0,-1)$ and radius $\sqrt 2$ transforms the Poincar\'e disc $D^n$ into the half-space model $H^n$. Here $n=2$.}
\label{inversion:fig}
\end{center}
\end{figure}

\begin{figure}
\begin{center}
\includegraphics[width = 10 cm] {\iftoggle{BW}{semispazio-BW}{semispazio}}
\nota{Lines and planes in the half-space model for $\matH^3$.}
\label{semispazio:fig}
\end{center}
\end{figure}

Some lines and planes in $H^3$ are drawn in Figure \ref{semispazio:fig}. The metric tensor $g$ on $H^n$ has a particularly nice form. 

\begin{prop}
The metric tensor on $H^n$ is:
$$g_x = \frac 1 {x_n^2} \cdot g^E$$
where $g^E$ is the Euclidean metric tensor on $H^n \subset \matR^n$.
\end{prop}
\begin{proof}
The inversion $\varphi\colon D^n \to H^n$ is the function
\begin{align*}
\varphi(x_1,\ldots, x_n) & = (0,\ldots,0,-1) + 2\frac{(x_1,\ldots,x_{n-1},x_n+1)}{\|(x_1,\ldots,x_{n-1},x_n+1)\|^2} \\
& = \frac{(2x_1,\ldots,2x_{n-1},1-\|x\|^2)}{\|x\|^2 +2x_n+1}.
\end{align*}
As seen in the proof of Proposition \ref{inversione:prop}, the inversion $\varphi$ is anticonformal with scalar dilation
$$\frac{2}{\|(x_1,\ldots,x_{n-1},x_n+1)\|^2} = \frac{2}{\|x\|^2+2x_n+1}.$$
The map $\varphi$ hence transforms the metric tensor $\left(\frac 2 {1-\|x\|^2}\right)^2\cdot g^E$ in $x\in D^n$ into the following metric tensor in $\varphi(x) \in H^n$:
$$\left( \frac 2 {1-\|x\|^2}\right)^2 \cdot \left(\frac{\|x\|^2+2x_n+1}2\right)^2
\cdot g^E = 
\frac 1 {\varphi_n(x)^2}\cdot g^E.$$
The proof is complete.
\end{proof}

\subsection{Geodesics in the conformal models}
The disc $D^n$ and the half-space $H^n$ are conformal models for $\matH^n$. In both models the hyperbolic metric differs from the Euclidean one only by a multiplication with some function.

In the half-space $H^n$ the lines are Euclidean vertical half-lines or half-circles orthogonal to $\partial H^n$ as in Figure \ref{semispazio:fig}. Vertical geodesics have a particularly simple form.
\begin{prop} \label{verticale:prop}
The vertical geodesic in $H^n$ passing through the point $(x_1,\ldots,x_{n-1},1)$ at time $t=0$ and pointing upward with unit speed is:
$$\gamma(t) = (x_1,\ldots,x_{n-1},e^t).$$
\end{prop}
\begin{proof}
A tangent vector $v$ at $(x_1,\ldots,x_n) \in H^n$ has norm $\frac{\|v\|^E}{x_n}$ where $\|v\|^E$ indicates the Euclidean norm. We get 
$$\|\gamma'(t)\| = \|(0,\ldots,0,e^t)\| = \frac {e^t}{e^t} = 1.$$
Therefore $\gamma(t)$ runs on a line at unit speed.
\end{proof}
We can easily deduce a parametrisation for the geodesics in $D^n$ passing through the origin. Recall the hyperbolic tangent:
$$\tanh (t) = \frac {\sinh (t)}{\cosh (t)} = \frac{e^t-e^{-t}}{e^t+e^{-t}} = \frac{e^{2t} -1}{e^{2t}+1}.$$
\begin{prop}
The geodesic in $D^n$ passing through the origin at time $t=0$ and pointing towards $x \in S^{n-1}$ at unit speed is:
$$\gamma(t) = \frac{e^{t}-1}{e^{t}+1}\cdot x = \left(\tanh \tfrac {t}2\right)\cdot x.$$
\end{prop}
\begin{proof}
We can suppose $x=(0,\ldots,0,1)$ and obtain this parametrisation from that of the vertical line in $H^n$ through inversion.
\end{proof}

We obtain in particular:
\begin{cor} \label{exp:cor}
The exponential map $\exp_0\colon T_0D^n \to D^n$ at the origin $0 \in D^n$ is the diffeomorphism:
$$\exp_0(x) = \frac{e^{\|x\|}-1}{e^{\|x\|}+1}\cdot \frac{x}{\|x\|} = \left(\tanh\tfrac{\|x\|}2 \right)\cdot \frac x{\|x\|}.$$
\end{cor}

Since the isometries of $\matH^n$ act transitively on points, we deduce that the exponential map at any $p\in \matH^n$ is a diffeomorphism. As a consequence, the injectivity radius of $\matH^n$ is $\infty$, as in the Euclidean $\matR^n$.

\subsection{Isometries of the conformal models}
In the half-space model it is easy to identify some isometries:

\begin{prop} \label{traslazioni:prop}
The following are isometries of $H^n$:
\begin{itemize}
\item horizontal translations $x \mapsto x+b$ with $b=(b_1,\ldots, b_{n-1}, 0)$,
\item dilations $x \mapsto \lambda x$ with $\lambda>0$,
\item inversions with respect to spheres orthogonal to $\partial H^n$.
\end{itemize}
\end{prop}
\begin{proof}
Horizontal translations obviously preserve the metric tensor $g= \frac 1{x_n^2}\cdot g^E$. 
We indicate by $\|\cdot \|$ and $\|\cdot\|^E$ the hyperbolic and Euclidean norm of tangent vectors. On a dilation $\varphi( x)= \lambda x$ we get 
$$\|d\varphi_x(v)\| = \frac{\|d\varphi_x(v)\|^E}{\varphi(x)_n} =
\frac{\lambda \|v\|^E}{\lambda x_n} = 
\frac{\|v\|^E}{x_n} = \|v\|.
$$

Concerning inversions, up to conjugating by translations and dilations it suffices to consider the map $\varphi(x) = \frac{x}{\|x\|^2}$. We have already seen that $d\varphi_x$ is $\frac {1}{\|x\|^2}$ times a linear reflection. Therefore
$$\|d\varphi_x(v)\| = \frac{\|d\varphi_x(v)\|^E}{\varphi(x)_n} =
\frac{\|v\|^E/\|x\|^2 }{x_n/\|x\|^2} = 
\frac{\|v\|^E}{x_n} = \|v\|.
$$
This completes the proof.
\end{proof}

In the disc model we can easily write the isometries that fix the origin:

\begin{prop}
The group $O(n)$ acts isometrically on $D^n$.
\end{prop}
\begin{proof}
The metric tensor on $D^n$ has a spherical symmetry.
\end{proof}

It is harder to write the isometries that fix another point of $D^n$ or $H^n$. 

\begin{prop} \label{generate:Hn:prop}
On the conformal models $D^n$ and $H^n$, the isometry group is generated by inversions along spheres and reflections along Euclidean planes orthogonal to the boundary.
\end{prop}
\begin{proof}
We know from Proposition \ref{traslazioni:prop} that these maps are isometries of $H^n$. One such isometry fixes a hyperplane $S$, hence it is the hyperbolic reflection $r_S$ (which is the unique non-trivial isometry fixing $S$). Hyperbolic reflections generate the isometry group by Proposition \ref{reflections:generate:prop}. 

In $D^n$ the proof is the same, we leave as an exercise to prove that a sphere inversion preserves the metric tensor and is hence an isometry.
\end{proof}

\subsection{Balls in the conformal models}
How does a metric ball $B(x_0,r)$ in the hyperbolic space look like? The answer in the conformal models is surprisingly  simple. 

\begin{prop} \label{sfere:prop}
In the conformal models balls are Euclidean balls.
\end{prop}
\begin{proof}
In the disc model $B(0,r)$ is a ball centred at $0$ by symmetry. The ball $B(x_0,r)$ at another point $x_0$ is obtained from $B(0,r)$ by composing inversions, which send spheres to spheres and hence balls to balls. The inversion $D^n \to H^n$ also sends balls to balls. 
\end{proof}

The centre of the ball is not its Euclidean centre in general!

\subsection{The Klein model} \label{Klein;subsection}
There is a fourth model for the hyperbolic space that is some kind of intermediate version between the hyperboloid and the Poincar\'e disc model.\index{hyperbolic space!Klein model}\index{hyperbolic space!projective model}

The \emph{Klein} or \emph{projective model} for $\matH^n$ is obtained by embedding the hyperboloid $I^n$ inside $\matRP^n$ via the projection $\matR^{n+1}\setminus \{0\} \to \matRP^n$. The image of this embedding is an open disc $K^n\subset \matRP^n$ bounded by the quadric $x_1^2+\ldots +x_n^2-x_{n+1}^2 = 0$. We equip $K^n$ with the metric tensor transported from $I^n$ , so that $K^n$ is indeed a model for $\matH^n$. 

When read in the chart $x_{n+1} = 1$, the Klein model $K^n$ becomes the open disc $x_1^2 + \ldots + x_n^2 < 1$. It is like the Poincar\'e disc $D^n$, but with a different metric tensor! 

The Klein model $K^n$ is not conformal and its metric tensor is a bit more complicated than that of $D^n$ or $H^n$. On the other hand, the subspaces in $K^n$ are easier to identify: by definition, these are just the projective subspaces of $\matRP^n$ intersected with the open disc $K^n$. Therefore lines are straight lines, but the angles and distances are distorted. Note also that the isometries of $K^n$ are projective transformations. The distance function is particularly nice:

\begin{ex}
The distance of two distinct points $p, q$ in $K^n$ is 
$$d(p,q) = \frac 12 | \log \beta (p,q,r,s) |$$
where $r,s$ are the intersections of the projective line $l$ containing $p$ and $q$ with $\partial K^n$, and $\beta(p,q,r,s)$ is the cross-ratio of the four points.\index{cross-ratio}
\end{ex}

\begin{figure}
\begin{center}
\includegraphics[width = 4.5 cm] {\iftoggle{BW}{Klein-BW}{Klein}}
\nota{A tessellation of $\matH^2$ into regular triangles and heptagons in the Klein model.}
\label{uniform_tiling:fig}
\end{center}
\end{figure}

\section{Compactification and isometries of hyperbolic space}
In this section we compactify the hyperbolic space $\matH^n$ by adding some ``points at infinity''. The compactification will then be used to classify the isometries of $\matH^n$ into three types. We also study the mutual position of two subspaces and define some peculiar hyperbolic objects: the horospheres. 

\subsection{Points at infinity}
Let a \emph{geodesic half-line} in $\matH^n$ be a geodesic $\gamma\colon [0,+\infty) \to \matH^n$ with constant unit speed. \index{geodesic!geodesic half-line}
\begin{defn} The set $\partial \matH^n$ of the \emph{points at infinity} in $\matH^n$ is the set of all geodesic half-lines, considered up to the following equivalence relation:
$$\gamma_1 \sim \gamma_2 \ \Longleftrightarrow \sup_{t\in [0,+\infty)}\big\{d\big(\gamma_1(t), \gamma_2(t)\big) \big\} < + \infty.$$
\end{defn}

We add its points at infinity to $\matH^n$ by defining
$$\overline {\matH^n} = \matH^n \cup \partial \matH^n.$$

\begin{prop} \label{bordo:prop}
On the disc model there is a natural 1-1 correspondence between $\partial D^n$ and $\partial \matH^n$ and hence between the closed disc $\overline{D^n}$ and $\overline{\matH^n}$.
\end{prop}
\begin{proof}
A geodesic half-line $\gamma$ in $D^n$ is a circle or line arc orthogonal to $\partial D^n$ and hence the Euclidean limit $\lim_{t\to\infty}\gamma(t)$ is a point in $\partial D^n$. We now prove that two half-lines converge to the same point if and only if they lie in the same equivalence class.

Suppose that two half-geodesics $\gamma_1$, $\gamma_2$ converge to the same point $p\in\partial D^n$. Up to isometries and inversions, we can change the conformal model to $H^n$ and take $p= \infty$. In this nicer setting $\gamma_1$ and $\gamma_2$ are vertical lines:
$$\gamma_1(t) = (x_1,\ldots,x_{n-1},x_ne^t), \quad
\gamma_2(t) = (y_1,\ldots,y_{n-1},y_ne^t). $$ 
The geodesic
$$\gamma_3(t) = (y_1,\ldots,y_{n-1},x_ne^t) $$ 
is equivalent to $\gamma_2$ since $d(\gamma_2(t),\gamma_3(t))=|\ln \frac{y_n}{x_n}|$ for all $t$ and is also equivalent to $\gamma_1$ because $d(\gamma_1(t),\gamma_3(t))\to 0$ as shown in Figure \ref{verticali:fig}.

\begin{figure}
\begin{center}
\includegraphics[width = 6 cm] {\iftoggle{BW}{verticali-BW}{verticali}}
\nota{Two vertical lines $\gamma_1$ and $\gamma_3$ in the half-space model $H^n$ at \emph{Euclidean} distance $d$. The \emph{hyperbolic} length of the horizontal segment between them at height $x_n$ is $\frac d{x_n}$ and hence tends to zero as $x_n\to\infty$ (left). 
Using as a height parameter the more intrinsic hyperbolic arc-length, we see that the two vertical geodesics $\gamma_1$ and $\gamma_3$ approach at exponential rate, since $d(\gamma_1(t), \gamma_3(t))\leqslant de^{-t}$ (right).}
\label{verticali:fig}
\end{center}
\end{figure}

Suppose that $\gamma_1$ and $\gamma_2$ converge to distinct points in $\partial D^n$. We can use the half-space model again and suppose that $\gamma_1$ is upwards vertical and $\gamma_2$ tends to some other point in $\{x_n=0\}$. In that case we easily see that $d(\gamma_1(t), \gamma_2(t))\to \infty$: for any $M>0$ there is a $t_0>0$ such that $\gamma_1(t)$ and $\gamma_2(t)$ lie respectively in $\{x_{n+1}>M\}$ and $\left\{x_n<\frac 1M\right\}$ for all $t>t_0$. Whatever curve connects these two open sets, it has length at least $\ln{M^2}$, hence $d(\gamma_1(t), \gamma_2(t)) > \ln{M^2}$ for all $t>t_0$.
\end{proof}

\subsection{The compactification}
We can give $\overline{\matH^n}$ the topology of $\overline{D^n}$, and in that way we have \emph{compactified} $\matH^n$ by adding its points at infinity. The interior of $\overline{\matH^n}$ is $\matH^n$, and the points at infinity form a sphere $\partial \matH^n$. \index{hyperbolic space!compactification of hyperbolic space}

The topology on $\overline {\matH^n}$ may also be defined intrinsically: for any $p\in\partial \matH^n$ we define a system of open neighbourhoods of $p$ in $\overline{\matH^n}$ as follows. Let $\gamma$ be a half-line with $[\gamma] = p$ and $V$ be an open neighbourhood of the vector $\gamma'(0)$ in the unitary sphere in $T_{\gamma(0)}\matH^n$. Pick $r>0$ and define the subset of $\overline{\matH^n}$:
\begin{align*}
U(\gamma, V, r) & = \left\{\alpha(t)\ \big|\ \alpha(0) = \gamma(0), \ \alpha'(0) \in V, \ t>r \right\} \\
& \bigcup \left\{[\alpha] \ \big|\ \alpha(0) = \gamma(0), \ \alpha'(0) \in V \right\}
\end{align*}
where $\alpha$ indicates a half-line in $\matH^n$ and $[\alpha]\in\partial\matH^n$ its class, see Figure \ref{intorno:fig}. We define an open neighbourhoods system $\{U(\gamma, V, r)\}$ for $p$ by letting $\gamma$, $V$, and $r$ vary. The resulting topology on $\overline{\matH^n}$ coincides with that induced by $\overline{D^n}$.

\begin{figure}
\begin{center}
\includegraphics[width = 5 cm] {\iftoggle{BW}{intorno-BW}{intorno}}
\nota{An open neighbourhood $U(\gamma, V, r)$ of $p\in\partial\matH^n$\iftoggle{BW}{ (in light grey)}{ (in yellow)}. We use the Klein model here.}
\label{intorno:fig}
\end{center}
\end{figure}

Note that $\matH^n$ is a complete Riemannian manifold (and hence a complete metric space), while its compactification $\overline{\matH^n}$ is only a topological space: a point in $\partial \matH^n$ has infinite distance from any other point in $\overline{\matH^n}$.

\subsection{Klein model} \label{projective:compactification:subsection}
In the Klein model $K^n\subset \matRP^n$, the compactification is obtained by adding the quadric $\partial K^n = \{x_1^2 + \ldots + x_n^2 - x_{n+1}^2=0\}$, which is the image of the light cone in $\matR^{n+1}$. Analogously, in the hyperboloid model $I^n$ we may represent $\partial I^n$ as the set of rays in the light cone.

A nice feature of the Klein model is that the points $x\in\matRP^n$ that lie ``beyond the infinity'', that is outside $\overline {K^n}$, can also be given a natural geometric interpretation. One such $x\in \matRP^n \setminus \overline{K^n}$ is a space-like ray in $\matR^{n+1}$ and as such it defines a hyperplane $x^\bot$ of signature $(n-1,1)$ in $\matR^{n+1}$, which projects to a hyperbolic hyperplane in $K^n$. The points beyond the infinity are in 1-1 correspondence with the hyperbolic hyperplanes in $K^n$.

\subsection{Incident, parallel, and ultraparallel subspaces}
A $k$-subspace $S\subset \matH^n$ has a topological closure $\overline S$ in the compactification $\overline{\matH^n}$. The boundary $\partial S = \overline S\cap \partial\matH^n$ of $S$ is a $(k-1)$-sphere.

For instance, the boundary of a line $l$ consists of two distinct points, the \emph{endpoints} of $l$. The boundary of a plane is a circle. Two distinct points in $\partial \matH^n$ are the endpoints of a unique line. A circle in the sphere $\partial D^3$ is the boundary of a unique plane in the disc model.

The usual \emph{distance} $d(A,B)$ between two subsets $A,B$ in a metric space is defined as 
$$d(A,B) = \inf_{x\in A, y\in B} \big\{d(x,y)\big\}.$$

There are three types of configurations for two subspaces in $\matH^n$, depicted in Figure \ref{lines:fig}.

\begin{figure}
\begin{center}
\includegraphics[width = 11 cm] {\iftoggle{BW}{lines-BW}{lines}}
\nota{Two incident, asymptotic parallel, and ultraparallel lines.}
\label{lines:fig}
\end{center}
\end{figure}

\begin{prop} \label{sottospazi:prop}
Let $S$ and $S'$ be subspaces in $\matH^n$ of arbitrary dimension. Precisely one of the following holds: 
\begin{enumerate}
\item $S\cap S' \neq \emptyset$,
\item $S\cap S' = \emptyset$ and $\overline S \cap \overline {S'}$ is a point in $\partial \matH^n$; moreover $d(S,S')=0$ and there is no geodesic orthogonal to both $S$ and $S'$,
\item $\overline S\cap \overline{S'} = \emptyset$; moreover $d=d(S,S')>0$ and there is a unique geodesic $\gamma$ orthogonal to both $S$ and $S'$: the segment of $\gamma$ between $S$ and $S'$ is the unique arc connecting them of length $d$.
\end{enumerate}
\end{prop}
\begin{proof}
If $\overline S \cap \overline{S'}$ contains two points then it contains the line connecting them and hence $S\cap S'\neq \emptyset$.

In (2) we use the half-space model and send $\overline S\cap \overline {S'}$ at infinity. Then $S$ and $S'$ are Euclidean vertical subspaces and Figure \ref{verticali:fig} shows that $d(S,S')=0$. Geodesics are vertical or half-circles and cannot be orthogonal to both $S$ and $S'$.

In (3), let $x_i \in S$ and $x'_i \in S'$ be such that $d(x_i,x_i')\to d$. Since $\overline{\matH^n}$ is compact, on a subsequence $x_i\to x\in \overline S$ and $x_i'\to x'\in\overline {S'}$. By hypothesis $x\neq x'$ and hence $x,x' \in \matH^n$ since $d<\infty$. Therefore $d>0$.

\begin{figure}
\begin{center}
\includegraphics[width = 4 cm] {\iftoggle{BW}{unica_geodetica-BW}{unica_geodetica}}
\nota{Two ultraparallel subspaces $S, S'$ and a line $\gamma$ orthogonal to both.}
\label{unica_geodetica:fig}
\end{center}
\end{figure}

Let $\gamma$ be the line passing through $x$ and $x'$. The segment between $x$ and $x'$ has length $d(x,x') = d$. The line is orthogonal to $S$ and $S'$: if it had an angle smaller than $\frac\pi 2$ with $S'$ we could find another point $x''\in S'$ near $x'$ with $d(x,x'')<d$. We can draw $S, S', \gamma$ as in Figure \ref{unica_geodetica:fig} by placing the origin between $x$ and $x'$: no other line can be orthogonal to both $S$ and $S'$.
\end{proof}

Two subspaces of type (1), (2) or (3) are called respectively \emph{incident}, \emph{asymptotically parallel}, and \emph{ultra-parallel}.
\index{subspace!incident, asymptotically parallel, ultra-parallel subspaces}

\subsection{The conformal sphere at infinity}
The sphere at infinity $\partial \matH^n$ has no metric structure, but it has instead a \emph{conformal structure}, that is a Riemannian structure considered up to conformal transformations. Before defining it we note the following.

\begin{prop}
Every isometry $\varphi \colon \matH^n \to \matH^n$ extends to a unique homeomorphism $\varphi:\overline{\matH^n} \to \overline{\matH^n}$. An isometry $\varphi$ is determined by its trace $\varphi|_{\partial \matH^n}$ at the boundary.
\end{prop}
\begin{proof}
The extension of $\varphi$ to $\partial \matH^n$ is defined intrinsically: a boundary point is a class $[\gamma]$ of geodesic half-lines and we set $\varphi ([\gamma]) = [\varphi(\gamma)]$. 

To prove the second assertion we show that an isometry $\varphi$ that fixes the points at infinity is the identity. The isometry $\varphi$ fixes every line as a set (because it fixes its endpoints), and since every point is the intersection of two lines it fixes also every point.
\end{proof}

We give $\partial \matH^n$ the conformal structure of the sphere $\partial D^n$. The group $\Iso(\matH^n)$ is generated by sphere inversions, which act conformally: although the metric tensor of $\partial D^n$ is not preserved by this action, its conformal class is preserved and hence the conformal structure on $\partial \matH^n$ is well-defined.

\subsection{Elliptic, parabolic, and hyperbolic isometries}
It is convenient to classify the isometries of $\matH^n$ into three types.
\index{isometry!elliptic, parabolic, and hyperbolic isometry}

\begin{prop} \label{fixed:prop}
Let $\varphi$ be a non-trivial isometry of $\matH^n$. Precisely one of the following holds: 
\begin{enumerate}
\item $\varphi$ has at least one fixed point in $\matH^n$, 
\item $\varphi$ has no fixed points in $\matH^n$ and has exactly one in $\partial \matH^n$,
\item $\varphi$ has no fixed points in $\matH^n$ and has exactly two in $\partial \matH^n$.
\end{enumerate}
\end{prop}
\begin{proof}
The extension $\varphi:\overline{\matH^n} \to \overline{\matH^n}$ is continuous and has a fixed point by Brouwer's Theorem. We only need to prove that if $\varphi$ has three fixed points $P_1, P_2, P_3$ at the boundary then it has some fixed point in the interior. The isometry $\varphi$ fixes the line $r$ with endpoints $P_1$ and $P_2$.
There is only one line $s$ with endpoint $P_3$ and orthogonal to $r$ (exercise): the isometry $\varphi$ must also fix $s$ and hence fixes the point $r\cap s$.
\end{proof}
Isometries of type (1), (2), and (3) are called respectively \emph{elliptic}, \emph{parabolic}, and \emph{hyperbolic}. A hyperbolic isometry fixes two points $p,q\in \partial \matH^n$ and hence preserves the unique line $l$ with endpoints $p$ and $q$. The line $l$ is the \emph{axis} of the hyperbolic isometry, which acts on $l$ as a translation.

\subsection{Horospheres} \label{orosfere:subsection}
Parabolic transformations are related to some objects in $\matH^n$ called horospheres. \index{horosphere}

\begin{defn} Let $p$ be a point in $\partial\matH^n$. A \emph{horosphere} centred in $p$ is a connected complete hypersurface orthogonal to all the lines exiting from $p$.
\end{defn}

Horospheres may be easily visualised in the half-space model $H^n$ by sending $p$ at infinity. The lines exiting from $p$ are the Euclidean vertical lines and the horospheres centred at $p$ are precisely the horizontal hyperplanes $\{x_n = k\}$ with $k>0$. 

\begin{oss}
Since the metric tensor $g = \frac 1{x_n^2}g^E$ is constant on each hyperplane $\{x_n = k\}$, each horosphere is \emph{isometric} to the Euclidean $\matR^n$. 
\end{oss}

The horospheres centred at $p\neq \infty$ in $\partial H^n$ or at any point $p\in \partial D^n$ are precisely the Euclidean spheres tangent in $p$ to the sphere at infinity. The horospheres in $\matH^2$ are circles and are called \emph{horocycles}, see Figure \ref{Horocycle_normals:fig}. The portion of $\matH^n$ delimited by a horosphere is called a \emph{horoball}.

\begin{figure}
\begin{center}
\includegraphics[width = 5 cm] {\iftoggle{BW}{Horocycle_normals-BW}{Horocycle_normals}}
\nota{A horocycle in $\matH^2$ centred in $p\in\partial\matH^2$ is a circle tangent to $p$. It is orthogonal to all the lines exiting from $p$.}
\label{Horocycle_normals:fig}
\end{center}
\end{figure}

Let us go back to the isometries of $\matH^n$. We sometimes write a point in the half-space $H^n$ as a pair $(x,t)$ with $x\in \matR^{n-1}$ and $t>0$. Isometries with nice fixed points have nice expressions in the conformal models. 

\begin{prop} \label{tre:prop}
Let $\varphi$ be an isometry of $\matH^n$:
\begin{enumerate}
\item if $\varphi$ is elliptic with fixed point $0\in D^n$ then
$$\varphi(x) = Ax$$
for some matrix $A\in \On(n)$; 
\item if $\varphi$ is parabolic with fixed point $\infty$ in $H^n$ then
$$\varphi (x,t) = (Ax+b, t)$$ 
for some matrix $A\in \On(n-1)$ and some vector $b$;
\item if $\varphi$ is hyperbolic with fixed points $0$ and $\infty$ in $H^n$ then
$$\varphi (x,t) = \lambda(Ax,t)$$
for some matrix $A\in \On(n-1)$ and some positive scalar $\lambda\neq 1$.
\end{enumerate}
\end{prop}
\begin{proof}
Point (1) is obvious. In (2) the isometry $\varphi$ fixes $\infty$ and hence permutes the horospheres centred at $\infty$: we first prove that this permutation is trivial. The map $\varphi$ sends a horosphere $O_0$ at height $t=t_0$ to a horosphere $O_1$ at some height $t=t_1$. If  $t_1\neq t_0$, up to changing $\varphi$ with its inverse we may suppose that $t_1<t_0$. 

We know that the map $\psi\colon O_1 \to O_0$ sending $(x,t_1)$ to $(x,t_0)$ is a contraction: hence $\varphi\circ\psi \colon O_1 \to O_1$ is a contraction and thus has a fixed point $(x,t_1)$. Therefore $\varphi(x,t_0) = (x,t_1)$. Since $\varphi(\infty)=\infty$, the vertical geodesic passing through $(x,t_0)$ and $(x,t_1)$ is preserved by $\varphi$, and hence we have found another fixed point $(x,0)\in\partial \matH^n$, a contradiction.

We now know that $\varphi$ preserves the horosphere $O$ at height $t$, for all $t$.
The metric tensor on $O$ is Euclidean (rescaled by $\frac 1{t^2}$), hence $\varphi$ acts on $O$ like an isometry $x\mapsto Ax+b$.
Since $\varphi$ sends vertical geodesics to vertical geodesics, it acts with the same formula on each horosphere and we are done.

Concerning (3), the axis $l$ of $\varphi$ is the vertical line with endpoints $0=(0,0)$ and $\infty$, and $\varphi$ acts on $l$ by translations: hence it sends $(0,1)$ to some $(0,\lambda)$. The differential $d\varphi$ at $(0,1)$ is necessarily $\matr A 0 0 \lambda$ for some $A\in \On(n-1)$ and hence $\varphi$ is globally as stated. The case $\lambda = 1$ is excluded because $(0,1)$ would be a fixed point in $\matH^n$.
\end{proof}

The \emph{minimum displacement} $d=d(\varphi)$ of an isometry $\varphi$ of $\matH^n$ is \index{minimum displacement}
$$d(\varphi) = \inf_{x\in\matH^n} d\big(x,\varphi(x)\big).$$
A point $x$ realises the minimum displacement if $d(x,\varphi(x)) = d(\varphi)$.
\begin{cor} \label{spostamento:cor}
The following hold:
\begin{enumerate}
\item an elliptic transformation $\varphi$ has $d=0$ realised on its fixed points;
\item a parabolic transformation $\varphi$ with fixed point $p\in \partial\matH^n$ has $d=0$ realised nowhere and fixes every horosphere centred in $p$; 
\item a hyperbolic transformation $\varphi$ with fixed points $p,q\in\partial\matH^n$ has $d>0$ realised on its axis.
\end{enumerate}
\end{cor}
\begin{proof}
Point (1) is obvious. Point (2) was already noticed while proving Proposition \ref{tre:prop}. Concerning (3), let $l$ be the axis of the hyperbolic transformation $\varphi$. The hyperplane orthogonal to $l$ at a point $x\in l$ is sent to the hyperplane orthogonal to $l$ in $\varphi(x)$. The two hyerplanes are ultraparallel and by Proposition \ref{sottospazi:prop} their minimum distance is realised at the points $x$ and $\varphi(x)$. Hence the points on $l$ realise the minimum displacement for $\varphi$.
\end{proof}

\section{Isometry groups in dimensions two and three} \label{isometry:2:3:section}
With the hyperboloid model the isometry group $\Iso(\matH^n)$ is the matrix group $O^+(n,1)$. We now see that in dimensions $n=2$ and $3$ the group $\Iso^+(\matH^n)$ is also isomorphic to some familiar groups of $2\times 2$ matrices. 

We start with a concise tour on M\"obius transformations. 
\index{M\"obius transformation and anti-transformation}

\subsection{M\"obius transformations}
Let the \emph{Riemann sphere} be $S=\matC \cup \{\infty\}$, homeomorphic to $S^2$. Consider the group
$$ \PSLC = \SLC/_{\pm I} = \GLC/_{\{ \lambda I\}} = \PGLC $$
of all $2\times 2$ invertible complex matrices considered up to scalar multiplication. The group $\PSLC$ acts on $S$ as follows: a matrix $\matr abcd \in \PSLC$ determines the \emph{M\"obius transformation} 
$$z \mapsto \frac{az+b}{cz+d}$$
which is an orientation-preserving self-diffeomorphism of $S$, and in fact also a  biolomorphism.

\begin{ex} \label{Mobius:freely:ex}
M\"obius transformations act freely and transitively on triples of distinct points in $S$.
\end{ex}

A matrix $ \matr abcd\in \PSLC$ also determines a \emph{M\"obius anti-transformation}
$$z \mapsto \frac{a\bar z+b}{c\bar z+d}$$
which is an orientation-reversing self-diffeomorphism of $S$. The composition of two anti-transformations is a M\"obius transformation. Transformations and anti-transformations together form a group $\Conf(S)$ which contains the M\"obius transformations as an index-two subgroup.

\begin{prop}
Circle inversions and line reflections are both M\"obius anti-transformations and generate $\Conf(S)$.
\end{prop}
\begin{proof}
By conjugating with translations $z\mapsto z+b$ and complex dilations $z \mapsto az$ every circle inversion transforms into the inversion along the unit circle $z\mapsto \frac{1}{\bar z}$, and every line reflection transforms into $z\mapsto \bar z$.

By composing line reflections we get all translations and rotations, and by composing circle inversions we get all dilations. With these operations and the inversion $z\mapsto \frac 1 {\bar z}$ one can easily act transitively on triples of points. They generate $\Conf(S)$ by Exercise \ref{Mobius:freely:ex}.
\end{proof}

\subsection{M\"obius transformations of $H^2$}
We consider the half-plane $H^2\subset \matC$ as $H^2 = \{z\ |\ \Im z >0\}$ and denote by $\Conf(H^2)$ the subgroup of $\Conf(S)$ consisting of all maps that preserve $H^2$.

\begin{rem}
By standard results in complex analysis, the two groups $\Conf(S)$ and $\Conf(H^2)$ contain precisely all the conformal diffeomorphisms of $S$ and $H^2$, whence their names. We will not use this fact here.
\end{rem}

\begin{prop}
The maps in $\Conf(H^2)$ are of the form
$$z \mapsto \frac{az+b}{cz+d} \qquad {\rm and}
\qquad z \mapsto \frac{a\bar z+b}{c\bar z+d}$$
with $a,b,c,d \in \matR$ and having $ad-bc$ equal to $1$ and $-1$, respectively. 
\end{prop}
\begin{proof}
The transformations listed have real coefficients and hence preserve the line $\matR\cup \infty$ and permute the two half-planes in $\matC\setminus \matR$. The sign condition on $ad-bc = \pm 1$ ensures precisely that $i$ is sent to some point in $H^2$ and hence $H^2$ is preserved.

On the other hand, a transformation that preserves $H^2$ must preserve $\partial H^2 = \matR\cup \infty$ and it is easy to see that since the images of $0,1,\infty$ are real all the coefficients $a,b,c,d$ can be taken in $\matR$.
\end{proof}

The M\"obius transformations in $\Conf(H^2)$ form a subgroup of index two which is naturally isomorphic to
$$\PSLR = \SLR/_{\pm I}.$$

An ordered triple of distinct points in $\matR\cup \infty$ is \emph{positive} if they are oriented counterclockwise, like $0,1, \infty$.

\begin{ex} \label{triple:ex}
The group $\PSLR$ acts freely and transitively on positive triples of points in $\matR\cup \infty$. 
\end{ex}

Let $C\subset \matC$ be a circle or line orthogonal to $\matR$. The inversion or reflection along $C$ preserves $H^2$ and is hence an element of $\Conf(H^2)$.

\begin{prop} \label{H2:generate}
Inversions along circles and reflections along lines orthogonal to $\matR$ generate $\Conf(H^2)$.
\end{prop}
\begin{proof}
Composing reflections we obtain all horizontal translations $z \mapsto z+b$ with $b\in\matR$, composing inversions we obtain all dilations $z \mapsto \lambda z$ with $\lambda \in \matR^*$. These maps together with the inversion $z \mapsto \frac {1}{\bar z}$ act transitively on positive triples of points in $\matR\cup \infty$.
\end{proof}

\begin{ex}
The inversion sending $H^2$ to $D^2$ is 
\begin{equation*} \label{anti:eq}
z \mapsto \frac{\bar z+i}{i\bar z+1}.
\end{equation*}
\end{ex}

\subsection{Isometries of $H^2$}
After this short detour on M\"obius transformations, we turn back to our hyperbolic spaces. We can characterise the isometry group of $H^2$. 
\begin{prop} \label{H2:prop}
We have $\Iso(H^2) = \Conf(H^2)$.
\end{prop}
\begin{proof}
Both groups are generated by inversions along circles and reflections along lines orthogonal to $\partial H^2 = \matR$ by Propositions \ref{generate:Hn:prop} and \ref{H2:generate}.
\end{proof}

In particular we have
$$\Iso^+(H^2) = \PSLR.$$
We will henceforth identify these two groups.
The trace of an element in $\PSLR$ is well-defined up to sign and carries some relevant information:
\begin{prop} \label{PSLR:prop}
A non-trivial isometry $A\in \PSLR$ is elliptic, parabolic, hyperbolic $\Longleftrightarrow$ respectively $|\tr A|<2$, $|\tr A|=2$, $|\tr A|>2$.
\end{prop}
\begin{proof}
Take $A=\big(\begin{smallmatrix} a & b \\ c & d \end{smallmatrix}\big)$ with $\det A = ad-bc = 1$. The M\"obius transformation $z\mapsto \frac {az+b}{cz+d}$ has a fixed point $z\in\matC$ if and only if
$$\frac{az+b}{cz+d} = z \Longleftrightarrow cz^2+(d-a)z-b = 0.$$
We find
$$\Delta = (d-a)^2 +4bc = (d+a)^2 - 4 = \tr^2 A - 4.$$
There is a fixed point in $H^2$ if and only if $\Delta <0$; if $\Delta > 0$ we find two fixed points in $\matR\cup\{\infty\}$ and if $\Delta = 0$ only one.
\end{proof}

\subsection{Isometries of $H^3$}
We have proved that $\Iso^+(H^2) = \PSLR$. Quite surprisingly, the isometry group $\Iso^+(H^3)$ is also isomorphic to a group of $2\times 2$ matrices! To prove this, we make the following identifications:
$$\matR^3 = \matC \times \matR = \{(z,t)\ |\ z\in\matC, t\in\matR\}$$
hence $H^3 = \{(z,t)\ |\ t >0\}$.
We also write $\matC$ for $\matC \times \{0\}$. The boundary trace of an isometry of $H^3$ is a homeomorphism of the Riemann sphere
$$\partial H^3 = \matC \cup\{\infty\} = S.$$
\begin{prop} \label{H3:prop}
The boundary trace induces an identification
$$\Iso(H^3) = \Conf(S).$$
\end{prop}
\begin{proof}
The group $\Iso(H^3)$ is generated by inversions along spheres and reflections along planes orthogonal to $\partial H^3$. Their traces are inversions along circles and reflections along lines in $S$. These generate $\Conf(S)$.
\end{proof}
In particular we have
$$ \Iso^+(H^3) = \PSLC.$$
We will also henceforth identify these two groups. As above, the trace of an element in $\PSLC$ is well-defined up to sign and carries some information:
\begin{prop}
A non-trivial isometry $A\in \PSLC$ is elliptic, parabolic, hyperbolic if and only if respectively $\tr A\in (-2,2)$, $\tr A=\pm 2$, $\tr A\in\matC \setminus [-2,2]$.
\end{prop}
\begin{proof}
Every non-trivial matrix $A\in\SLC$ is conjugate to one of:
$$\pm\begin{pmatrix} 1 & 1 \\ 0 & 1 \end{pmatrix}, \quad 
\begin{pmatrix} \lambda & 0 \\ 0 & \lambda^{-1} \end{pmatrix} $$
for some $\lambda \in \matC^*$, and these represent the following isometries:
$$(z,t) \longmapsto (z + 1,t), \quad (z,t) \longmapsto (\lambda^2 z, |\lambda|^2 t).$$
In the first case $\tr A = \pm 2$ and $A$ is parabolic with fixed point $\infty$, in the second case $A$ has a fixed point in $H^3$ if and only if $|\lambda|=1$, \emph{i.e.}~$\tr A = \lambda + \lambda^{-1}\in (-2,2)$, the fixed point being $(0,1)$. If $|\lambda|\neq 1$ there are two fixed points $0$ and $\infty$ at infinity and hence $A$ is hyperbolic.
\end{proof}

Summing up, we have
$$\Iso^+(H^2) = \PSLR, \qquad \Iso^+(H^3) = \PSLC.$$
The group $\PSLR$ acts directly on $H^2$, while $\PSLC$ acts on the boundary sphere of $H^3$.

\section{Geometry of hyperbolic space}
We study the geometry of $\matH^n$. We prove that $\matH^n$ has constant sectional curvature $-1$, that the distance function along lines is convex, we define convex combinations and barycenters, and study parallel transport along lines. Finally, we prove that $\Iso(\matH^n)$ is a unimodular Lie group.

\subsection{Area and curvature}
We can verify that $\matH^n$ has constant sectional curvature $-1$. It should be no surprise that $\matH^n$ has constant curvature, since it has many symmetries (\emph{i.e.}~isometries). To calculate its sectional curvature we compute the area of a disc. 

\begin{prop}
The disc of radius $r$ in $\matH^2$ has area 
$$A(r) = \pi \left(e^{\frac r2} -e ^{-\frac r2}\right)^2 =  4\pi \sinh^2 {\tfrac r2} = 2\pi (\cosh r -1).$$
\end{prop}
\begin{proof}
Recall that the volume form is
$$\omega = \sqrt{\det g} \cdot dx_1\cdots dx_n.$$
Let $D(r)$ be a disc in $\matH^2$ of radius $r$. If we centre it in $0$ in the disc model, its Euclidean radius is $\tanh \frac r2$ by Corollary \ref{exp:cor} and we get
\begin{align*}
A(r) & = \int_{D(r)} \sqrt{\det g} \cdot dxdy = \int_{D(r)} \left(\frac 2{1-x^2-y^2}\right)^2dxdy \\
 & = \int_0^{2\pi} \int_0^{\tanh\frac r2} \left(\frac {2}{1-\rho^2}\right)^2 \rho \cdot d\rho d\theta = 2\pi\left[ \frac 2{1-\rho^2}\right]_0^{\tanh \frac r2} \\
 & = 4\pi \left(\frac 1{1-{\tanh^2 \frac r2}} - 1 \right) = 4\pi \sinh^2 \tfrac r2.
\end{align*}
The proof is complete.
\end{proof}

\begin{cor}
The hyperbolic space $\matH^n$ has sectional curvature $-1$.
\end{cor}
\begin{proof}
Pick $p\in \matH^n$ and $W\subset T_p$ a 2-dimensional subspace. The image $\exp_p(W)$ is the hyperbolic plane tangent to $W$ in $p$. On a hyperbolic plane 
$$A(r) = 2\pi(\cosh r - 1) = 2\pi \left(\frac {r^2}{2!} + \frac {r^4}{4!} + o(r^4)\right) = \pi r^2 + \frac{\pi r^4}{12} + o(r^4)$$
and hence $K=-1$ following the area formula in Section \ref{gauss:subsection}.
\end{proof}

\subsection{Convexity of the distance function}
We recall that a function $f\colon \matR^n \to \matR$ is \emph{strictly convex} if 
$$f(tv+(1-t)w)< tf(v) + (1-t)f(w)$$
for any pair $v,w \in \matR^n$ of distinct points and any $t\in (0,1)$. 
\begin{ex} A positive strictly convex function is continuous and admits a minimum if and only if it is proper.
\end{ex}

We now prove that the distance function is strictly convex on disjoint lines of $\matH^n$. Given two lines $l, l'\subset \matH^n$, we fix an isometry of each line with $\matR$ and we get an identification of $l\times l'$ with the Euclidean plane $\matR\times\matR$.\index{convexity of the distance function}

\begin{figure}
\begin{center}
\includegraphics[width = 9cm] {\iftoggle{BW}{convex-BW}{convex}}
\nota{Distance between points in disjoint lines is a strictly convex function in hyperbolic space.}
\label{convex:fig}
\end{center}
\end{figure}

\begin{prop} \label{distanza:convessa:prop}
Let $l,l'\subset \matH^n$ be two disjoint lines. The map
\begin{align*}
l \times l' & \longrightarrow \matR_{\geqslant 0} \\
(x,y) & \longmapsto d(x, y)
\end{align*}
is strictly convex; it is proper if and only if the lines are ultraparallel. 
\end{prop}
\begin{proof}
With our identifications we have $(x,y) \in \matR \times \matR$. The function $d$ is clearly continuous, hence to prove its convexity it suffices to show that  
$$d\left(\frac{x_1+x_2}2, \frac{y_1+y_2}2\right) < \frac{d(x_1,y_1) + d(x_2,y_2)}2$$
for any pair of distinct points $(x_1,y_1), (x_2,y_2) \in l\times l'$. Suppose $x_1\neq x_2$ and denote by $m$ and $n$ the midpoints $\frac{x_1+x_2}2$ and $\frac{y_1+y_2}2$ as in Figure \ref{convex:fig}. 

Let $\sigma_p$ be the reflection at the point $p\in\matH^n$. The isometry $\tau = \sigma_n\circ \sigma_m$ translates the line $r$ containing the segment $mn$ by the quantity $2d(m,n)$: hence it is a hyperbolic transformation with axis $r$. We draw the points  $o = \tau(m)$ and $z_i = \tau(x_i)$ in the figure and note that $z_1 = \sigma_n(x_2)$, hence $d(x_2,y_2) = d(z_1,y_1)$. The triangular inequality implies that
$$d(x_1, z_1) \leqslant d(x_1,y_1) + d(y_1, z_1) = d(x_1,y_1) + d(x_2, y_2).$$
A hyperbolic transformation has minimum displacement on its axis $r$ and $x_1 \neq m$ is not in $r$, hence
$$2d(m,n) = d(m,o) = d(m,\tau(m)) < d(x_1, \tau(x_1)) = d(x_1,z_1).$$
Finally we get $2d(m,n) <  d(x_1,y_1) + d(x_2, y_2)$ and hence $d$ is convex. 

The function $d$ is proper, that is it has minimum, if and only if the two lines are ultraparallel by Proposition \ref{sottospazi:prop}.
\end{proof}

\begin{ex}
The distance function on parallel lines in $\matR^n$ is not strictly convex (it is only convex).
\end{ex}

\subsection{Convex combinations} \label{convex:subsection}
Let $p_1,\ldots, p_k$ be $k$ points in $\matH^n, \matR^n$, or $S^n$ and $t_1,\ldots, t_k$ be non-negative numbers with $t_1+\ldots + t_k=1$. The \emph{convex combination}  \index{convex combination} 
$$p = t_1p_1+ \ldots + t_kp_k$$ 
is another point in the space defined as follows:
\begin{align*}
{\rm in\ } \matR^n: \quad p & =  t_1p_1+ \ldots + t_kp_k \\
{\rm in\ } I^n,S^n: \quad p & = \frac{t_1p_1+ \ldots + t_kp_k}{\| t_1p_1+ \ldots + t_kp_k\|} 
\end{align*}
where $\|v\| = \sqrt{-\langle v, v \rangle}$ on $I^n$. Using convex combination we may define the \emph{barycenter} of the points as $\frac 1k p_1 + \ldots + \frac 1k p_k$. The barycenter may in turn be used to prove the following.
\index{barycenter}

\begin{prop} \label{iterate:prop}
Let $\varphi\colon \matH^n \to \matH^n$ be a non-trivial isometry and $k\geqslant 2$:
\begin{itemize}
\item if $\varphi$ is elliptic then $\varphi^k$ is elliptic or trivial; 
\item if $\varphi$ is parabolic then $\varphi^k$ is parabolic;
\item if $\varphi$ is hyperbolic then $\varphi^k$ is hyperbolic.
\end{itemize}
\end{prop}
\begin{proof}
If $\varphi\colon \matH^n \to\matH^n$ is an isometry with no fixed points, then $\varphi^k$ also is: if $\varphi^k(x)=x$ then $\varphi$ fixes the finite set $\{x, \varphi(x), \ldots, \varphi^{k-1}(x)\}$ and hence also its barycenter.

If $\varphi$ is parabolic then it fixes the horospheres centred at some point $p\in\partial\matH^n$ and also $\varphi^k$ does, hence it is still parabolic (it cannot be hyperbolic). If $\varphi$ is hyperbolic it has two fixed points at infinity, and $\varphi^k$ too.
\end{proof}

\begin{prop} \label{finite:fixed:prop}
Every finite group $\Gamma<\Iso(\matH^n)$ fixes a point in $\matH^n$.
\end{prop}
\begin{proof}
The barycenter of any orbit is fixed by $\Gamma$.
\end{proof}

\subsection{Parallel transport} \label{parallel:subsection}
On Riemannian manifolds, the parallel transport is a way to slide frames along geodesics. On $\matH^n$ we can do this simply as follows: for every geodesic $\gamma$, we put $\gamma$ into vertical position in the half-space model $H^n$, and slide the frames vertically in the obvious way.  \index{parallel transport}

This construction furnishes in particular, for every pair $x,y \in \matH^n$ of points, a canonical isometry between $T_x\matH^n$ and $T_y\matH^n$, obtained by sliding frames along the unique geodesic $\gamma$ containing $x$ and $y$.

This canonical identification is of course not transitive on a triple of non-collinear points $x,y,z$: the curvature of $\matH^n$ is responsible for that.

\subsection{Unimodularity}
The following is a consequence of Corollary \ref{O:unimodular:cor}.

\begin{cor} \label{unimodular:cor}
The isometry group $\Isomet(\matH^n)$ is unimodular.
\end{cor}

\begin{oss} \label{Haar:oss}
A Haar measure for $\Iso(\matH^n)$ may be constructed concretely as follows: fix a point $x\in \matH^n$ and define the measure of a Borel set $S\subset \Iso(\matH^n)$ as the measure of $S(x) = \cup_{\varphi\in S} \varphi(x)\subset\matH^n$. This measure is obviously left-invariant, and is hence also right-invariant since $\Iso(\matH^n)$ is unimodular. As a consequence, it does not depend on the choice of $x$.
\end{oss}

\subsection{References}
The hyperbolic space is introduced in various books: two standard references are Benedetti -- Petronio \cite{BP} and Ratcliffe \cite{R}, and most of the arguments presented here were borrowed from these two sources. The proof of Proposition \ref{distanza:convessa:prop} was taken from Farb -- Margalit \cite{FM}.
Thurston's notes contain some useful trigonometric formulae that we have omitted, see \cite[Chapter 2]{Th}.

%% file: Manifold.tex
\chapter{Hyperbolic manifolds} \label{varieta:chapter} \label{manifold:chapter} 
 \label{orbifold:chapter} \label{orbifolds:chapter}

A hyperbolic manifold is a Riemannian manifold locally isometric to the hyperbolic space $\matH^n$. Maybe the most striking aspect of geometric topology is that, despite this quite restrictive definition, there are plenty of hyperbolic manifolds around, especially in the dimensions $n=2$ and $3$. For that reason hyperbolic manifolds (and hence hyperbolic geometry) play a central role in the topology of surfaces and three-manifolds.

The study of complete hyperbolic manifolds is tightly connected to that of discrete subgroups in the Lie group $\Iso(\matH^n)$ and of polyhedra in $\matH^n$, so it has both an algebraic and geometric flavour. We start this chapter by describing these connections; then we show some examples and discuss some important variations: non-complete hyperbolic manifolds, hyperbolic manifolds with geodesic boundary, cone manifolds, and orbifolds.

\section{Discrete groups of isometries}
We define hyperbolic manifolds and prove a crucial theorem, that says that every complete hyperbolic manifold is isometric to a quotient $\matH^n/_\Gamma$ for some discrete group $\Gamma <\Iso(\matH^n)$ acting freely on $\matH^n$.

\subsection{Hyperbolic, flat, and elliptic manifolds}
We introduce three important classes of Riemannian manifolds.
\index{hyperbolic manifold} \index{flat manifold} \index{elliptic manifold}

\begin{defn} A \emph{hyperbolic} (resp.~\emph{flat} or \emph{elliptic}) \emph{manifold} is a connected Riemannian $n$-manifold that may be covered by open sets isometric to open sets of $\matH^n$ (resp.~$\matR^n$ o $S^n$). 
\end{defn}

A hyperbolic (resp.~flat or elliptic) manifold has constant sectional curvature $-1$ (resp.~$0$ or $+1$). We show that the model $\matH^n$ is indeed unique.

\begin{teo} \label{simply:teo}
Every complete simply connected hyperbolic $n$-manifold $M$ is isometric to $\matH^n$.
\end{teo}
\begin{proof}
Pick a point $x\in M$ and choose an isometry $D\colon U \to V$ between an open ball $U$ containing $x$ and an open ball $V\subset \matH^n$. 
We show that $D$ extends (uniquely) to an isometry $D\colon M\to\matH^n$. 

For every $y \in M$, choose an arc $\alpha\colon [0,1]\to M$ from $x$ to $y$. By compactness there is a partition $0=t_0 < t_1 < \ldots < t_k =1$ and for each $i=0,\ldots, k-1$ an isometry $D_i\colon U_i \to V_i$ from an open ball $U_i$ in $M$ containing $\alpha([t_i,t_{i+1}])$ to an open ball $V_i\subset \matH^n$. 

We may suppose that $U_0 \subset U$ and $D_0 = D|_{U_0}$. Inductively on $i$, we now modify $D_i$ so that $D_{i-1}$ and $D_i$ coincide on the component $C$ of $U_{i-1}\cap U_i$ containing $\alpha(t_i)$. To do so, note that
$$D_{i-1}\circ D_i^{-1}\colon D_i(C) \longrightarrow D_{i-1}(C)$$
is an isometry of open connected sets in $\matH^n$ and hence extends to an isometry of $\matH^n$. Then it makes sense to compose $D_i$ with $D_{i-1}\circ D_i^{-1}$, so that the new maps $D_{i-1}$ and $D_i$ coincide on $C$. Finally, we define $D(y) = D_{k-1}(y)$.

The proof that $D(y)$ is well-defined is a standard argument. First, it is easy to check that different partitions $0=t_0 <\ldots <t_k=1$ do not vary $D(y)$, just by considering a common refinement. Then we consider another path $\beta$ connecting $x$ to $y$. Since $M$ is simply-connected, there is a homotopy connecting $\alpha$ and $\beta$. The image of the homotopy is compact and is hence covered by finitely many open balls $U_i$ isometric to open balls $V_i\subset \matH^n$ via some maps $D_i$. By the Lebesgue number theorem, there is a $N>0$ such that in the grid in $[0,1]\times [0,1]$ of $\frac 1N \times \frac 1N$ squares, the image of every square is entirely contained in at least one $U_i$. We can now modify as above the isometries $D_i$ inductively on the grid, starting from the bottom-left square, so that they all glue up and show that $D(y)$ does not depend on $\alpha$ or $\beta$.

The resulting map $D\colon M \to \matH^n$ is a local isometry by construction. Since $M$ is complete, the map $D$ is a covering by Proposition \ref{local:isometry:prop}. Since $\matH^n$ is simply connected, the covering $D$ is a homeomorphism and $D$ is actually an isometry.
\end{proof}

The isometry $D\colon M \to \matH^n$ constructed in the proof is called a \emph{developing map}.
The same proof shows that every complete simply connected flat (or elliptic) $n$-manifold is isometric to $\matR^n$ (or $S^n$). \index{developing map}

\subsection{Complete hyperbolic manifolds}
We have determined the unique complete simply connected hyperbolic $n$-manifold, and we now look at complete hyperbolic manifolds with arbitrary fundamental group. We first note that if $\Gamma < \Iso(\matH^n)$ is a group of isometries that acts freely and properly discontinuously on $\matH^n$, the quotient manifold $\matH^n/_\Gamma$ has a natural Riemannian structure that promotes the covering
$$\pi\colon\matH^n \longrightarrow \matH^n/_\Gamma$$
to a local isometry, see Proposition \ref{quotient:prop}. The quotient $\matH^n/_\Gamma$ is a complete hyperbolic manifold. 
We now show that every complete hyperbolic manifold is realised in this way:

\begin{prop} \label{corto:prop}
Every complete hyperbolic $n$-manifold $M$ is isometric to $\matH^n/_\Gamma$ for some subgroup $\Gamma<\Iso(\matH^n)$ acting freely and properly discontinuously.
\end{prop}
\begin{proof}
The universal cover of $M$ inherits a Riemannian structure that is complete (by Proposition \ref{local:isometry:prop}), hyperbolic, and simply connected: hence it is isometric to $\matH^n$ by Theorem \ref{simply:teo}. The deck transformations $\Gamma$ of the covering $\matH^n\to M$ are necessarily locally isometries, hence isometries. We conclude that $M=\matH^n/_\Gamma$ and $\Gamma$ acts freely and properly discontinuously using Proposition \ref{equivalent:prop}.
\end{proof}

Note that $\Gamma$ is isomorphic to the fundamental group $\pi_1(M)$.

\begin{oss}
A group $\Gamma<\Iso(\matH^n)$ acts freely if and only if it does not contain elliptic isometries: that is, every non-trivial isometry in $\Gamma$ is either hyperbolic or parabolic.
\end{oss}
Note also that $\Gamma$ acts properly discontinuously if and only if it is discrete, see Proposition \ref{iff:discrete:prop}. 

\begin{oss}
The same proofs show that every complete flat or spherical $n$-manifold is isometric to $\matR^n/_\Gamma$ or $S^n/_\Gamma$ for some discrete group $\Gamma$ of isometries acting freely on $\matR^n$ or $S^n$. 
\end{oss}

\begin{cor} There is a natural 1-1 correspondence
$$\left\{\!
\begin{array}{c}{\rm complete\ hyperbolic} \\ {\rm manifolds\ } M {\rm\ up\ to\ isometry}\end{array}\!\right\} \longleftrightarrow
\left\{\!
\begin{array}{c}{\rm discrete\ subgroups\ } \Gamma < \Iso(\matH^n) \\ {\rm without\ elliptics} \\ {\rm up\ to\ conjugation}\end{array}\!\right\}$$
\end{cor}
\begin{proof}
When passing from the complete hyperbolic manifold $M$ to the group $\Gamma$, the only choice we made is an isometry between the universal cover of $M$ and $\matH^n$. Different choices produce conjugate groups $\Gamma$.
\end{proof}

\subsection{Discrete groups}
We investigate some basic properties of discrete groups $\Gamma$ of isometries of $\matH^n$.

\begin{ex} If $\Gamma<\Iso(\matH^n)$ is discrete then it is countable.
\end{ex}

Note that $\Gamma$ is not necessarily finitely generated. We denote by $\Gamma_p < \Gamma$ the stabiliser of a point $p\in \matH^n$.

\begin{prop} \label{discrete:prop}
Let $\Gamma <\Iso(\matH^n)$ be discrete and $p\in\matH^n$ be a point. The stabiliser $\Gamma_p$ is finite and the orbit $\Gamma(p) = \{g(p)\ |\ g\in \Gamma\}$ is discrete.
\end{prop}
\begin{proof}
Both are obvious consequence of the fact that $\Gamma$
acts properly discontinuously.
\end{proof}

Of course $\Gamma$ acts freely on $\matH^n$ if and only if $\Gamma_p$ is trivial for all $p\in\matH^n$. A set of subspaces in $\matH^n$ is \emph{locally finite} if every compact subset in $\matH^n$ intersects only finitely many of them.

\begin{prop} \label{dense:prop}
Let $\Gamma<\Iso(\matH^n)$ be discrete. The points $p\in\matH^n$ with trivial stabiliser $\Gamma_p$ form an open dense set in $\matH^n$.
\end{prop}
\begin{proof}
The fixed-points set $\Fix(g)$ of a non-trivial isometry $g$ is a proper subspace of $\matH^n$. The subspaces $\Fix(g)$ with $g\in\Gamma$ are locally finite: if infinitely many of them intersect a compact set they accumulate and $\Gamma$ does not act properly discontinuously. The complement of a locally finite set of proper subspaces is open and dense.
\end{proof}

Recall that a group \emph{has no torsion} if every non-trivial element has infinite order.

\begin{prop} \label{no:torsion:prop}
A discrete group $\Gamma < \Iso(\matH^n)$ acts freely on $\matH^n$ if and only if it has no torsion.
\end{prop}
\begin{proof}
By Proposition \ref{iterate:prop} parabolic and hyperbolic elements have infinite order. On the other hand, an elliptic element $g\in\Gamma$ has finite order since $\Gamma_p$ is finite for $p\in\Fix(g)$.
\end{proof}

\begin{cor}
The fundamental group of a complete hyperbolic manifold has no torsion.
\end{cor}

It is now time to exhibit some examples. Recall that $H^2$ is the half-plane model, see Section \ref{isometry:2:3:section}.

\begin{example}
The \emph{modular group} 
$$\Gamma = \PSLZ < \PSLR = \Iso^+(H^2)$$ 
consists of all matrices in $\PSLR$ having integer entries and is clearly a discrete subgroup. It does not act freely on $H^2$, however: the matrix $\matr 01{-1}0$ represents the elliptic transformation $z\mapsto -\frac 1z$ with fixed point $i$. \index{modular group}
\end{example} 

\subsection{Coverings}
We now make a simple but crucial observation: if $\Gamma < \Iso(\matH^n)$ acts freely and properly discontinuously, then also every subgroup $\Gamma' < \Gamma$ does; we get a manifolds covering
$$\matH^n/_{\Gamma'} \longrightarrow \matH^n/_\Gamma$$
whose degree $d$ is precisely the index of $\Gamma'$ in $\Gamma$.  Recall from Proposition \ref{volume:covering:prop} that we get
$$\Vol\big(\matH^n/_{\Gamma'}\big) = d \cdot \Vol\big(\matH^n/_{\Gamma}\big)$$
where some of the terms in the formula may be infinite. Moreover, every covering of a hyperbolic complete $M= \matH^n/_\Gamma$ is constructed in this way: there is a nice bijective correspondence
$$\big\{{\rm coverings\ of\ }M \big\} \longleftrightarrow
\big\{{\rm subgroups\ of\ }\Gamma \big\} .$$
This holds of course also for flat and spherical manifolds.

\subsection{Congruence subgroups}
We can now exhibit a family of two-dimensional hyperbolic manifolds. 

Pick an integer $m\geqslant 2$. Let $\SL_2(\matZ/_{m\matZ})$ be the group of $2\times 2$ matrices with coefficients in $\matZ/_{m\matZ}$ and determinant 1. We define the quotient 
$$\matP\SL_2(\matZ/_{m\matZ}) = \SL_2(\matZ/_{m\matZ})/_{\{ \pm I\}}.$$
The reduction modulo $m$ homomorphism $\matZ \to \matZ/_{m\matZ}$ induces the group homomorphisms $\SL_2(\matZ) \to \SL_2(\matZ/_{m\matZ})$ and 
$$\PSLZ \longrightarrow \matP\SL_2(\matZ/_{m\matZ}).$$
The kernel of this homomorphism is the \emph{principal congruence subgroup} $\Gamma(m)$ of $\PSLZ$. It is clearly discrete, since $\PSLZ$ is. It has finite index in $\PSLZ$ because $\matP\SL_2(\matZ/_{m\matZ})$ is finite. 
\index{principal congruence subgroup}
\begin{prop}
If $m\geqslant 4$ the group $\Gamma(m)$ acts freely on $H^2$.
\end{prop}
\begin{proof}
An element $A\in\Gamma(m)$ is a matrix $\matr abcd$ congruent to $\matr 1001$ modulo $m$. In particular $a+d$ is congruent to 2 modulo $m$, and hence is not $-1, 0, 1$. Therefore $A$ is never elliptic. 
\end{proof}

The quotient $\matH^2/_{\Gamma(m)}$ is a hyperbolic surface. We will construct many hyperbolic surfaces in Section \ref{surface:geometrisation:section} via some more geometric methods.

\subsection{Selberg's lemma}
The discrete group $\Gamma = \PSLZ$ does not act freely on $H^2$, but its finite-index normal subgroup $\Gamma(m)$ does as soon as $m\geqslant 4$. Is this the instance of a more general principle? Yes, it is. \index{Selberg lemma}

\begin{prop} \label{Selberg:prop}
Every finitely generated discrete group $\Gamma < \Iso(\matH^n)$ contains a finite-index normal subgroup $\Gamma'\triangleleft \Gamma$ that acts freely on $\matH^n$.
\end{prop}
\begin{proof}
The group $\Iso(\matH^n)$ is isomorphic to $O^+(n,1) < \GL(n+1,\matC)$, so Selberg's Lemma \ref{Selberg:lemma} applies to $\Gamma$, and it furnishes a finite-index torsion-free normal subgroup $\Gamma'\triangleleft \Gamma$. This subgroup acts freely by Proposition \ref{no:torsion:prop}.
\end{proof}

Every finitely generated discrete group $\Gamma < \Iso(\matH^n)$ contains at least one finite-index subgroup that acts freely. But how many such subgroups does it contain? Quite a lot, in fact.

\begin{prop} \label{FG:RF:prop}
Every finitely generated discrete group $\Gamma<\Iso(\matH^n)$ is residually finite.
\end{prop}
\begin{proof}
See Lemma \ref{RF:lemma}.
\end{proof}

The following corollary shows that there is an abundance of torsion-free subgroups. Although algebraic in nature, it has some remarkable geometric consequences, that will be revealed soon in Section \ref{finite:covers:subsection}.

\begin{cor}
Let $\Gamma < \Iso(\matH^n)$ be discrete and finitely generated. For every non-trivial $g \in \Gamma$ there is a finite-index normal subgroup $\Gamma'\triangleleft \Gamma$ that acts freely on $\matH^n$ and does not contain $g$.
\end{cor}

\section{Polyhedra} 
A polyhedron in $\matH^n$ is a natural geometric object, that may be used to visualise discrete groups in $\Iso(\matH^n)$ and hence hyperbolic manifolds. 
Polyhedra may sometimes be combined to form some tessellations of the space. 

\subsection{Polyhedra} \label{polyhedra:subsection}
A \emph{half-space} in $\matH^n$ is the closure of one of the two portions of space delimited by a hyperplane. We say that a set of half-spaces is locally finite if their boundary hyperplanes are. \index{polyhedron}

\begin{defn}
A $n$-dimensional \emph{polyhedron} $P$ in $\matH^n$ is the intersection of a locally finite set of half-spaces. We also assume that $P$ has non-empty interior.
\end{defn}

A subset $S\subset \overline{\matH^n}$ is \emph{convex} if $x,y\in S$ implies that the segment connecting $x,y$ is also contained in $S$ (such a segment is a half-line or a line if one or both points lie in $\partial \matH^n$). Every polyhedron $P$ is clearly convex because it is the intersection of convex sets. Its closure $\overline P$ in $\overline{\matH^n}$ is also convex.
\index{convex subset}

Let $H\subset \matH^n$ be a half-space containing the polyhedron $P$. If non-empty, the intersection $F=\partial H\cap P$ is called a \emph{face} of $P$. The \emph{supporting subspace} of $F$ is the smallest subspace of $\matH^n$ containing $F$; the \emph{dimension} of a face is the dimension of its supporting subspace. A face of dimension $0$, $1$, and $n-1$ is called a \emph{vertex}, an \emph{edge}, and a \emph{facet}.

\begin{ex} \label{intersection:faces:ex}
If non-empty, the intersection of faces of $P$ is a face.
\end{ex}

\begin{ex} Every $k$-dimensional face is a polyhedron in its supporting $k$-dimensional space.
\end{ex}

The \emph{convex hull} of a set $S\subset \overline{\matH^n}$ is the intersection of all the convex sets containing $S$. \index{convex hull}

\begin{ex} \label{hull:ex}
The convex hull of finitely many points in $\matH^n$ that are not contained in a hyperplane is a compact polyhedron. Conversely, every compact polyhedron has finitely many vertices and is the convex hull of them.
\end{ex}

Everything we said holds with no modifications for $\matR^n$. On $S^n$ some care should be taken: some definitions need to be modified slightly to take into account the annoying presence of antipodal points. We gloss over this technical point.

\subsection{Finite polyhedra}
We now enlarge slightly the class of compact polyhedra by admitting finitely many vertices at infinity. \index{polyhedron!finite polyhedron}

\begin{defn}
A \emph{finite polyhedron} is the convex hull of finitely many points $x_1,\ldots, x_k \in \overline{\matH^n}$ that are not contained in the closure of a hyperplane.
\end{defn}

The $x_i$'s that lie in $\partial \matH^n$ are called \emph{ideal vertices}, while the usual vertices of $P$ are the \emph{finite} or \emph{actual vertices}. The ideal vertices form the set $\overline P \setminus P$.

\begin{ex}
A finite polyhedron $P$ has finitely many faces and is the convex hull of its ideal and finite vertices.
\end{ex}

We want to estimate the volume of finite polyhedra, and to do this we need a lemma. Given a horosphere $O$ centred at $p\in\partial \matH^n$ and a domain $D\subset O$, the \emph{cone} $C$ of $D$ over $p$ is the union of all the half-lines exiting from $D$ and pointing towards $p$, see Figure \ref{volume:fig}.

\begin{lemma} \label{O:lemma}
Let $O$ be a horosphere centred at $p\in \partial \matH^n$, $D\subset O$ any domain and $C$ the cone over $D$. The following equality holds:
$$\Vol(C) = \frac {\Vol_O (D)}{n-1}$$
where $\Vol_O$ is the $(n-1)$-dimensional volume in the flat $(n-1)$-manifold $O$.
\end{lemma}
\begin{proof}
Let $O$ have some height $x_n=h$ as in Figure \ref{volume:fig}-(left). We obtain
$$\Vol(C) = \int_D dx \int_{h}^\infty \frac 1{t^{n}}dt = \frac 1{n-1}\int_D \frac{dx}{h^{n-1}} = \frac 1{n-1} \cdot \Vol_O(D).$$
The proof is complete.
\end{proof}

\begin{figure}
\begin{center}
\includegraphics[width = 12 cm] {\iftoggle{BW}{volume-BW}{volume}} 
\nota{The cone $C$ over a domain $D\subset O$ has volume proportional to the area of $D$ (left). If the domain is compact, the cone has finite volume: therefore a finite polyhedron has finite volume (right).}
\label{volume:fig}
\end{center}
\end{figure}

We now turn to finite polyhedra.

\begin{prop}
Every finite polyhedron has finite volume.
\end{prop}
\begin{proof}
For every ideal vertex of $P$, a small horoball centred at $p$ intersects $P$ into a cone that has finite volume. The polyhedron $P$ decomposes into finitely many cones and a bounded region, see Figure \ref{volume:fig}-(right).
\end{proof}

A finite polyhedron without finite vertices is called an \emph{ideal polyhedron}.
\index{polyhedron!ideal polyhedron}

\subsection{Polygons}
A polygon is just a polyhedron of dimension two. In contrast with Euclidean geometry, a strikingly simple formula relates the area of a finite polygon with its inner angles. We define the inner angle of an ideal vertex to be zero. 

\begin{figure}
\begin{center}
\includegraphics[width = 10 cm] {\iftoggle{BW}{triangle_area-BW}{triangle_area}}
\nota{A triangle with at least an ideal vertex (left). The area of a triangle with finite vertices can be derived as the area difference of triangles with one ideal vertex (right).}
\label{triangle_area:fig}
\end{center}
\end{figure}

\begin{prop} \label{polygon:area:prop}
A polygon $P$ with inner angles $\alpha_1,\ldots, \alpha_n$ has area
$$\Area (P) = (n-2)\pi - \sum_{i=1}^n \alpha_i.$$
\end{prop}
\begin{proof}
Every polygon decomposes into triangles, and it suffices to prove the formula on these. Consider first a triangle $T$ with at least one vertex at infinity. We use the half-plane model and send this vertex to $\infty$ as in Figure \ref{triangle_area:fig}-(left). We suppose that the \iftoggle{BW}{grey}{red} dot is the origin of $\matR^2$, hence 
$$T = \big\{(r\cos\theta, y)\ | \ \beta\leqslant \theta \leqslant \pi-\alpha, \ y\geqslant r\sin\theta\big\}$$
and we get
\begin{align*}
\Area(T) & = \int_T \frac 1{y^2}dxdy = \int_{\pi-\alpha}^{\beta} \int_{r\sin\theta}^\infty \frac{-r\sin \theta}{y^2}dyd\theta \\
 & = \int_{\pi-\alpha}^\beta -r\sin\theta \left[-\frac 1y\right]_{r\sin\theta}^\infty d\theta = 
\int_{\beta}^{\pi-\alpha} \frac{r\sin\theta}{r\sin\theta}d\theta \\
 & = \int_{\beta}^{\pi-\alpha} 1 = \pi-\alpha-\beta.
\end{align*} 
The area of a triangle with finite vertices $ABC$ is deduced as in Figure \ref{triangle_area:fig}-(right) using the formula
$$\Area (ABC) = \Area(AB\infty) + \Area(BC\infty) - \Area (AC\infty).$$
The proof is complete.
\end{proof}

The sum of the inner angles of a hyperbolic polygon is strictly smaller than that of a Euclidean polygon with the same number of sides, and the difference between these two numbers is precisely its area. 

\begin{cor} Every ideal triangle has area $\pi$.
\end{cor}

\subsection{Platonic solids} \label{platonici:subsection}
The theory of three-dimensional polyhedra in $\matH^3$ is very rich: we limit ourselves to the study of the platonic solids.

The five Euclidean platonic solids are the regular tetrahedron, the cube, the regular octahedron, icosahedron, and dodecahedron. We now see that each platonic solid $P\subset \matR^3$ generates a nice continuous family of solids in the three geometries $\matH^3, \matR^3$, and $S^3$.

To construct this family we fix any point $x$ in $\matH^3$ and represent $P$ centred at $x$ with varying size. To do this, we put $P$ inside the Euclidean tangent space $T_x\matH^3$ centred at the origin and with some radius $t>0$. Consider the image of its vertices by the exponential map and take their convex hull. We indicate by $P(-t)$ the resulting polyhedron in $\matH^3$. The polyhedron $P(-t)$ is combinatorially equivalent to $P$ and has the same symmetries of $P$.

We extend this family to the other geometries as follows. The polyhedron $P(0)$ is the Euclidean $P$ (unique up to dilations), and $P(t)$ with $t>0$ is the spherical $P$, constructed as above with $S^3$ instead of $\matH^3$: we fix $x\in S^3$, take a copy of $P$ inside $T_x S^3$ with radius $t$, project its vertices, and take the convex hull. We define the spherical $P(t)$ only for $t \in (0, \frac \pi 2]$: when $t= \frac\pi 2$ it degenerates to a hemisphere. We also define $P(-\infty)$ as the ideal hyperbolic platonic solid obtained by sending all the vertices at infinity.

\begin{figure}
\begin{center}
\includegraphics[width = 3.5 cm] {\iftoggle{BW}{Icosahedron-BW}{256px-Icosahedron}} \qquad
\includegraphics[width = 3.5 cm] {\iftoggle{BW}{Dodecahedron-BW}{245px-Dodecahedron}}
\nota{The regular icosahedron and dodecahedron.}
\label{icosaedro:fig}
\end{center}
\end{figure}

We have defined a polyhedron $P(t)$ for all $t\in [-\infty, \frac \pi 2]$, that lies in $\matH^3, \matR^3, S^3$ depending on whether $t$ is negative, null, or positive. In some sense the polyhedron $P(t)$ depends continuously on $t$ also when it crosses the value $t=0$: when $t\to 0$ the polyhedron in $\matH^3$ or $S^3$ shrinks, and every polyhedron tends to a Euclidean one when shrunk. 

In particular the dihedral angle $\theta(t)$ of $P(t)$ varies continuously with $t\in [-\infty, \frac \pi2]$. The function $\theta(t)$ is strictly monotone and we now determine its image. The \emph{vertex valence} of $P$ is the number of edges incident to each vertex.

\begin{prop} \label{dihedral:platonic:prop}
Let $n\in\{3,4,5\}$ be the vertex valence of $P$. It holds 
$$\theta \left(\left[-\infty, \frac \pi 2\right]\right) = \left[\frac{n-2}n \pi, \pi\right].$$
\end{prop}
\begin{proof}
Since $\theta$ is continuous and monotone increasing, it suffices to show that $\theta(-\infty) = \frac{n-2}n\pi$ and $\theta(\frac \pi 2) = \pi$.

By intersecting the ideal polyhedron $P(-\infty)$ with a small horosphere $O$ centred at an ideal vertex $v$ we get a regular $n$-gon in the Euclidean plane $O$, with interior angles $\frac{n-2}n\pi$. The dihedral angle at an edge $e$ incident to $v$ is measured by intersecting $P(-\infty)$ with any hypersurface orthogonal to $e$: since $O$ is orthogonal to $e$ we get $\theta(-\infty) =  \frac{n-2}n\pi$. 

The polyhedron $P(\frac \pi 2)$ is a hemisphere and hence $\theta(\frac \pi 2) = \pi$.
\end{proof}

For some values of $t$, the platonic solid $P(t)$ may have nice dihedral angles $\theta(t)$, for instance angles that divide $2\pi$. In the Euclidean world, the only platonic solid with such nice dihedral angles is the cube. In the hyperbolic and spherical world we find more.

The platonic solids with dihedral angles that divide $2\pi$ are listed in Table \ref{platonic:table}. The table is just a consequence of Proposition \ref{dihedral:platonic:prop}: it suffices to know the Euclidean dihedral angles for $P$, and all the angles bigger (smaller) than this value are spherical (hyperbolic).

\begin{table} 
\begin{center}
\begin{tabular}{c||cccc}
\phantom{\Big|} polyhedron & $\theta = \frac \pi 3$ & $\theta = \frac{2\pi}5$ & $\theta = \frac \pi 2$ &
$\theta = \frac{2\pi}3$ \\
\hline\hline
\rule{0pt}{3ex}
\phantom{\big|}  tetrahedron & ideal $\matH^3$ & $S^3$ & $S^3$ & $S^3$ \\
\phantom{\big|}  cube & ideal $\matH^3$ & $\matH^3$ & $\matR^3$ & $S^3$ \\
\phantom{\big|}  octahedron & & & ideal $\matH^3$ & $S^3$ \\
\phantom{\big|} icosahedron & & & & $\matH^3$ \\
\phantom{\big|} dodecahedron & ideal $\matH^3$ & $\matH^3$ & $\matH^3$ & $S^3$ \\
\end{tabular}
\end{center}
\nota{The platonic solids with dihedral angle $\theta$ that divide $2\pi$.}
\label{platonic:table}
\end{table}

In particular, there are four right-angled platonic solids: the spherical tetrahedron, the  Euclidean cube, the hyperbolic dodecahedron, and the ideal hyperbolic octahedron.

\section{Tessellations}
A tessellation is a nice paving of $\matH^n$ made of polyhedra. Not only tessellations are beautiful objects, but they are also tightly connected with discrete groups of $\Iso(\matH^n)$ and hence with hyperbolic manifolds.

\index{tessellation}

\begin{defn}
A \emph{tessellation} of $\matH^n$ (or $\matR^n$, $S^n$) is a locally finite set of polyhedra that cover the space and may intersect only in common faces.
\end{defn}

\begin{figure}
\begin{center}
\includegraphics[width = 5 cm] {\iftoggle{BW}{Uniform_tiling_532-t012-BW}{Uniform_tiling_532-t012}} \quad
\includegraphics[width = 5 cm] {\iftoggle{BW}{Truncated_octahedra-BW}{Truncated_octahedra}}
\nota{A tessellation of $S^2$ into squares, hexagons, and decagons, and a tessellation of $\matR^3$ into truncated octahedra.}
\label{tiling:fig}
\end{center}
\end{figure}

Some examples in spherical and Euclidean space are shown in Figure \ref{tiling:fig}. We now construct some families explicitly.

\subsection{Triangles} \label{triangular:subsection}
We want to define some nice tessellations of $\matH^2, \matR^2$, and $S^2$ into triangles. The following exercise says that every triple of acute angles is realised by some triangle in the appropriate geometry.

\begin{ex}
Given three real numbers $0<\alpha, \beta, \gamma \leqslant \frac \pi 2$ there is a triangle $\Delta$ with inner angles $\alpha, \beta, \gamma$ inside $\matH^2, \matR^2$, or $S^2$ depending on whether the sum $\alpha + \beta + \gamma$ is smaller, equal, or bigger than $\pi$. 
\end{ex}

Let $a,b,c \geqslant 2$ be three natural numbers and $\Delta$ be a triangle with inner angles $\frac \pi a, \frac \pi b, \frac \pi c$. The triangle $\Delta$ lies in $\matH^2$, $\matR^2$, or $S^2$ depending on whether $\frac 1a+\frac 1b + \frac 1c$ is smaller, equal, or bigger than 1. In all cases, by mirroring iteratively $\Delta$ along its edges we construct a tessellation of the space.

\begin{figure}
\begin{center}
\includegraphics[width = 3.5 cm] {\iftoggle{BW}{233_BW}{233}} \quad
\includegraphics[width = 3.5 cm] {\iftoggle{BW}{234_BW}{234}} \quad 
\includegraphics[width = 3.5 cm] {\iftoggle{BW}{235_BW}{235}} 
\nota{The tessellations $(2,3,3)$, $(2,3,4)$ and $(2,3,5)$ of the sphere.}
\label{triangular:fig}
\end{center}
\end{figure}

The triples realisable in $S^2$ are $(2,2,c)$, $(2,3,3)$, $(2,3,4)$, and $(2,3,5)$: the last three tessellations are shown in Figure \ref{triangular:fig} and are connected to the platonic solids. They consist of $24$, $48$, and $120$ triangles. 

\begin{figure}
\begin{center}
\includegraphics[width = 3.5 cm] {\iftoggle{BW}{Tile_V46b_BW}{200px-Tile_V46b}} \quad
\includegraphics[width = 3.5 cm] {\iftoggle{BW}{Tile_V488_bicolor_BW}{200px-Tile_V488_bicolor}} \quad 
\includegraphics[width = 3.5 cm] {\iftoggle{BW}{1024px-Tile_36_BW}{200px-Tile_36}} 
\nota{The tessellations $(2,3,6)$, $(2,4,4)$, and $(3,3,3)$ of the Euclidean plane.}
\label{triangular2:fig}
\end{center}
\end{figure}

The triples realisable in $\matR^2$ are $(2,3,6)$, $(2,4,4)$, and $(3,3,3)$: the tessellations are shown in Figure \ref{triangular2:fig}. There are infinitely many triples realisable in $\matH^2$, and some are shown in Figure \ref{triangular3:fig}. 

\begin{figure}
\begin{center}
\includegraphics[width = 3.7 cm] {\iftoggle{BW}{Order-3_heptakis_heptagonal_tiling-BW}{240px-Order-3_heptakis_heptagonal_tiling}} \quad
\includegraphics[width = 3.7 cm] {\iftoggle{BW}{Order-4_bisected_pentagonal_tiling-BW}{240px-Order-4_bisected_pentagonal_tiling}} \quad 
\includegraphics[width = 3.7 cm] {\iftoggle{BW}{Uniform_dual_tiling_433-t012-BW}{240px-Uniform_dual_tiling_433-t012}} 
\nota{The tessellations $(2,3,7)$, $(2,4,5)$, and $(3,3,4)$ of the hyperbolic plane.}
\label{triangular3:fig}
\end{center}
\end{figure}

In the hyperbolic plane we can also use triangles with vertices at infinity, which have inner angle zero by definition.
\begin{ex}
Given three real numbers $0\leqslant \alpha, \beta, \gamma \leqslant \frac \pi 2$ with sum smaller than $\pi$ there is a triangle $\Delta\subset \matH^2$ with inner angles $\alpha, \beta, \gamma$.
\end{ex}

For any triple $(a,b,c)$ of numbers in $\matN \cup \{\infty \}$ with $\frac 1a + \frac 1b + \frac 1c < 1$ we may take the triangle $\Delta \subset \matH^2$ with inner angles $\frac \pi a$, $\frac \pi b$, $\frac \pi c$ and reflect it iteratively to get a tessellation of $\matH^2$. 

\begin{figure}
\begin{center}
\includegraphics[width = 5 cm] {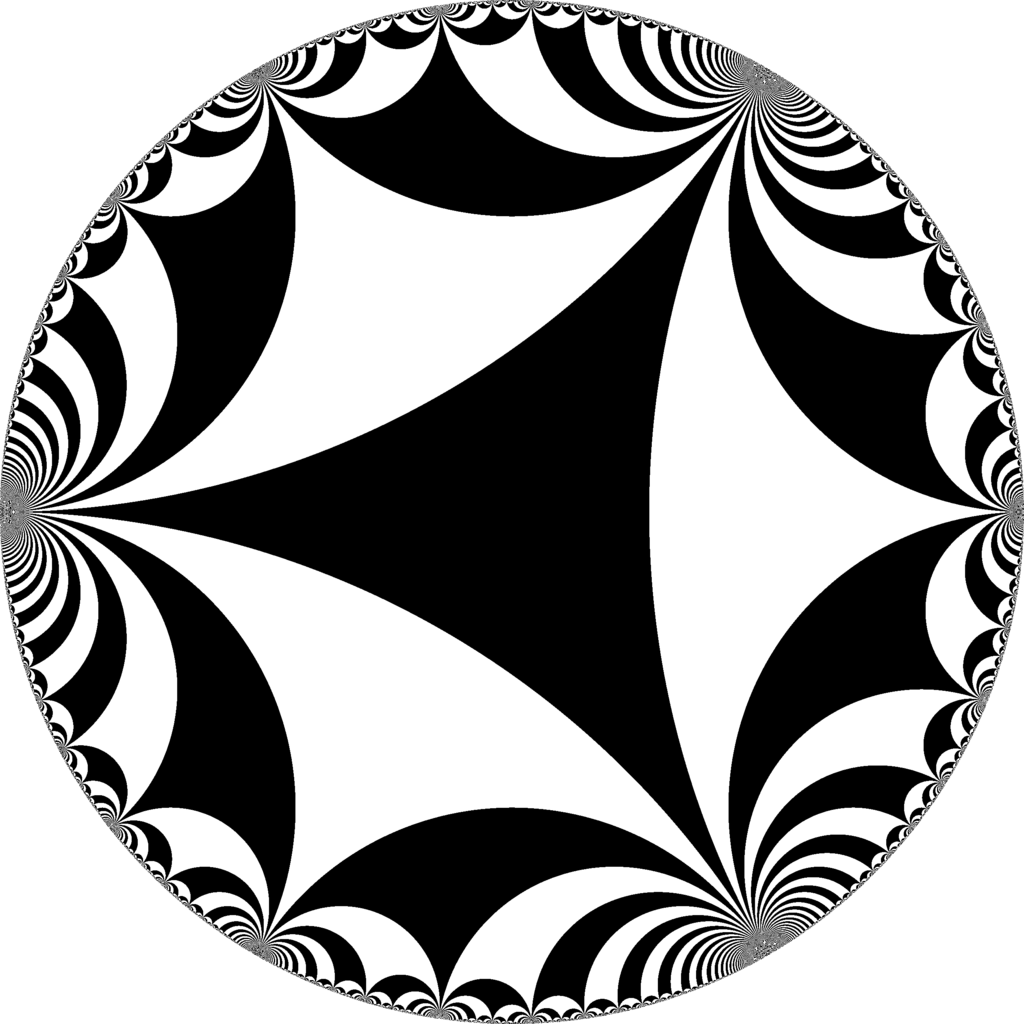}
\nota{The Farey tessellation $(\infty, \infty, \infty)$ of the hyperbolic plane.}
\label{tiling_Farey:fig}
\end{center}
\end{figure}

The triple $(\infty, \infty, \infty)$ gives a nice tessellation into ideal triangles called the \emph{Farey tessellation} and shown in Figure \ref{tiling_Farey:fig}. 
\index{tessellation!Farey tessellation}
 
\subsection{Platonic solids}
We now turn to tessellations of 3-dimensional spaces. Table \ref{platonic:table} displays a finite list of platonic solids $P$ with dihedral angles $\frac{2\pi}k$ contained in $\matH^3$, $\matR^3$, or $S^3$. For each solid in the list, by reflecting iteratively $P$ along its faces we get a tessellation of the space. 

The Figures \ref{honeycomb_dodecahedron_534:fig}, \ref{honeycomb_dodecahedron_536:fig}, and \ref{honeycomb_dodecahedron_533:fig} show the tessellations of $\matH^3, \matH^3,$ and $S^3$ by regular dodecahedra with dihedral angles $\frac \pi 2$, $\frac \pi 3$, and $\frac{2\pi}3$. The first two tessellations contain infinitely many polyhedra and appear in the figures as seen by an observer floating inside $\matH^3$. The reader is invited to observe the differences between the first two tessellations of the hyperbolic space. The third tessellation contains finitely many polyhedra and the figure shows its stereographic projection in $\matR^3$.

\begin{figure}
\begin{center}
\includegraphics[width = 9 cm] {\iftoggle{BW}{H3_534_CC_center-BW}{H3_534_CC_center}}
\nota{The tessellation of $\matH^3$ into right-angled regular dodecahedra, seen from inside $\matH^3$.}
\label{honeycomb_dodecahedron_534:fig}
\end{center}
\end{figure}

\begin{figure}
\begin{center}
\includegraphics[width = 9 cm] {\iftoggle{BW}{H3_536_CC_center-BW}{H3_536_CC_center}}
\nota{The tessellation of $\matH^3$ into ideal regular dodecahedra with dihedral angle $\frac \pi 3$, seen from inside $\matH^3$. Note that all vertices lie at infinity.}
\label{honeycomb_dodecahedron_536:fig}
\end{center}
\end{figure}

\begin{figure}
\begin{center}
\includegraphics[width = 9 cm] {\iftoggle{BW}{Stereographic_polytope_120cell_faces-BW}{600px-Stereographic_polytope_120cell_faces}}
\nota{The tessellation of $S^3$ into 120 regular dodecahedra with dihedral angle $\frac {2\pi} 3$, transposed into $\matR^3$ via the stereographic projection, that transforms straight faces into round ones but preserves the angles.}
\label{honeycomb_dodecahedron_533:fig}
\end{center}
\end{figure}

\subsection{Regular tessellations} \label{regular:subsection}
Platonic solids are \emph{regular}, that is they have many symmetries. A notion of regularity may be defined in all dimensions for polyhedra and tessellations as follows.
\index{tessellation!regular tessellation} \index{polyhedron!regular polyhedron}

Every tessellation in $\matH^n$, $\matR^n$, or $S^n$ has a symmetry group, consisting of all the isometries of the ambient space that preserve it. A \emph{flag} is a sequence of faces $f_0 \subset \ldots \subset f_n$ of the tessellation with $\dim f_i = i$. A tessellation is \emph{regular} if the symmetry group acts transitively on flags. 

These definitions apply also to polyhedra. We note that a regular polyhedron in $\matR^n$, considered up to similarities, may be interpreted as a regular tessellation in $S^{n-1}$, and vice versa.\index{Schl\"afli notation}

The \emph{Schl\"afli notation} beautifully encodes various regular tessellations, in all dimensions and in all geometries. The symbol $\{n\}$ with $n\geqslant 3$ denotes the regular $n$-gon, that is a regular tessellation of $S^1$ by $n$ congruent arcs, and by extension $\{\infty \}$ denotes the regular tessellation of $\matR$ by infinitely many congruent segments.

The symbol $\{p,q\}$ denotes a regular tessellation in $\matH^2, \matR^2,$ or $S^2$ by $p$-gons where $q$ of them meet at every vertex. The five platonic solids are
$$\{3,3\}, \{3,4\}, \{3,5\}, \{4,3\}, \{5,3\},$$
the symbols $\{3,6\}$, $\{4,4\}$, and $\{6,3\}$ denote the two regular tessellations of $\matR^2$ into equilateral triangles, squares, and regular hexagons, and every other pair $\{p,q\}$ denotes a tessellation of regular $p$-gons in $\matH^2$ with angles $\frac {2\pi}q$, that meet at $q$ at every vertex. We can also interpret $\{p, \infty\}$ as a regular tessellation of $\matH^2$ into ideal regular $p$-gons (the Farey tessellation is $\{3,\infty\}$) and $\{\infty, q\}$ as a regular tessellation of $\infty$-gons that meet at $q$ at every vertex, whose edges form a $q$-regular tree (a tree where $q$ edges meet at every vertex). Finally $\{\infty, \infty\}$ is a regular tessellation of ideal $\infty$-gons. See Figure \ref{infty:tessellations:fig}

\begin{figure}
\begin{center}
\includegraphics[width = 3.8 cm] {\iftoggle{BW}{600px-H2chess_iiif}{600px-H2chess_iiif}} \quad
\includegraphics[width = 3.8 cm] {\iftoggle{BW}{600px-H2chess_24ic}{600px-H2chess_24ic}} \quad 
\includegraphics[width = 3.8 cm] {\iftoggle{BW}{600px-H2chess_2iib}{600px-H2chess_2iib}} 
\nota{The tessellations of $\matH^2$ with Schl\"afli symbols $\{4,\infty\}$, $\{\infty, 4\}$, and $\{\infty, \infty\}$.}
\label{infty:tessellations:fig}
\end{center}
\end{figure}

The symbol $\{p,q,r\}$ denotes a regular tessellation in $\matH^3, \matR^3,$ or $S^3$ by polyhedra $\{p,q\}$ where $r$ of them meet at every edge. We deduce from Table \ref{platonic:table} that the regular tessellations of $S^3$ are
$$\{3,3,3\}, \{3,3,4\}, \{3,3,5\}, \{3,4,3\}, \{4,3,3\}, \{5,3,3\}.$$
These are the six regular polyhedra (usually called \emph{polytopes}) in dimension four: we will encounter them again in Chapter \ref{regular:polytopes:subsection}. The symbol $\{4,3,4\}$ indicates the tessellation of $\matR^3$ into cubes, and
$$\{3,3,6\}, \{3,4,4\}, \{3,5,3\}, \{4,3,5\}, \{4,3,6\}, \{5,3,4\}, \{5,3,5\}, \{5,3,6\}$$
denote the tessellations of $\matH^3$ into platonic solids. It is also possible to interpret more triples $\{p,q,r\}$ in an appropriate way, sometimes by representing the vertices of the polyhedra as space-like vectors in the hyperboloid model.

Every regular tessellation has a \emph{dual} regular tessellation obtained by taking the barycenters of all the cells involved. The dual of $\{a,b, \ldots, z\}$ is $\{z, \ldots, b,a\}$. 
Duals of ideal tessellations involve infinite polyhedra: for instance $\{6,3,3\}$ is shown in Figure \ref{honeycomb_dodecahedron_633:fig}.

\begin{figure}
\begin{center}
\includegraphics[width = 9 cm] {\iftoggle{BW}{H3_633_FC_boundary-BW}{H3_633_FC_boundary}}
\nota{In the tessellation $\{6,3,3\}$ of $\matH^3$ every polyhedron has infinitely many hexagonal faces, whose vertices all lie in a single horosphere.}
\label{honeycomb_dodecahedron_633:fig}
\end{center}
\end{figure}

The story continues in four dimensions: we denote by $\{p,q,r,s\}$ a regular tessellation made of polytopes of type $\{p,q,r\}$ that meet in $s$ at every codimension-two face. With the same techniques of Section \ref{platonici:subsection} we can identify the regular polytopes in $\matH^4, \matR^4$ and $S^4$ with dihedral angles that divide $2\pi$, and the corresponding tessellations of the ambient space -- that we are unfortunately unable to see.

\subsection{Voronoi tessellations} \label{Voronoi:subsection}
In the previous section we have constructed some tessellations by exploiting the symmetries of regular polyhedra: we can probably extend these methods to less symmetric polyhedra, but how far can we go? How can we construct highly non-regular tessellations?

There is a strikingly simple procedure that transforms every discrete set $S$ of points in $\matH^n$ into a tessellation, called the \emph{Voronoi tessellation} of $S$. The construction goes as follows. \index{tessellation!Voronoi tessellation}

For every point $p \in S$ we define
$$D(p) = \big\{ q\in \matH^n\ \big|\ d(q,p)\leqslant d(q,p') \quad \forall p'\in S \big\}.$$

\begin{prop}
The set $D(p)$ is a polyhedron and the polyhedra $D(p)$ as $p\in S$ varies form a tessellation of $\matH^n$.
\end{prop}
\begin{proof}
It is an easy exercise to show that the points in $\matH^n$ having the same distance from two distinct fixed points form a hyperplane. For every $p'\in S$ different from $p$ we define the half-space
$$H_{p'} = \big\{q \in \matH^n\ \big| \ d(q,p) \leqslant d(q,p') \big\}.$$
The set $D(p)$ is the intersection of the half-paces $H_{p'}$ as $p'$ varies in $S\setminus \{p\}$.

Since $S$ is discrete, these half-spaces are locally finite (there are finitely many points in $S$ at bounded distance from $p$, hence finitely many hyperplanes). Therefore $D(p)$ is a polyhedron. Since $S$ is discrete, every point $q\in \matH^n$ has at least one nearest point $p\in S$: therefore the polyhedra $D(p)$ cover $\matH^n$ as $p\in S$ varies. 

It remains to prove that the polyhedra $D(p)$ intersect along common faces. We have $D(p) \cap D(p') = D(p) \cap \partial H_{p'}$ and hence $D(p) \cap D(p')$ is either empty or a face of $D(p)$. The case of multiple intersections follows from Exercise \ref{intersection:faces:ex}.
\end{proof}

\begin{figure}
\begin{center}
\includegraphics[width = 5 cm] {\iftoggle{BW}{Coloured_Voronoi_2D-BW}{Coloured_Voronoi_2D}} 
\nota{A Voronoi tessellation of the Euclidean plane.}
\label{Voronoi:fig}
\end{center}
\end{figure}

Voronoi tessellations of course make sense also in $\matR^n$ and $S^n$, see Figure \ref{Voronoi:fig}. We have just proved that tessellations are not exoteric, but quite ordinary objects, and we now show that they are useful to study hyperbolic manifolds.

\section{Fundamental domains}
We turn back to our discrete subgroups $\Gamma < \Iso(\matH^n)$. The geometry of a discrete $\Gamma$ is nicely controlled by some polyhedra, called \emph{fundamental domains}.
We introduce these objects and make some important examples.

\subsection{Fundamental and Dirichlet domains}
Let $\Gamma$ be a discrete group of isometries of $\matH^n$. The group $\Gamma$ may or may not act freely on $\matH^n$.  \index{fundamental domain}
\begin{defn}
A \emph{fundamental domain} for $\Gamma$ is a polyhedron $D \subset\matH^n$ whose translates $g(D)$ as $g\in \Gamma$ varies are distinct and form a tessellation of the space $\matH^n$. 
\end{defn}
If $D$ is a fundamental domain, the group $\Gamma$ acts freely and transitively on the tessellation $\big\{g(D)\big\}_{g \in \Gamma}$, so in particular the polyhedra $g(D)$ are all isometric. 

We describe a procedure that builds a fundamental domain for any discrete subgroup $\Gamma < \Iso(\matH^n)$. Pick  a point $p\in \matH^n$ 
with trivial stabiliser $\Gamma_p$, which exists by Proposition \ref{dense:prop}. The group $\Gamma$ acts freely and transitively on the orbit $\Gamma(p)$, which is discrete by Proposition \ref{discrete:prop}.

The discrete orbit $\Gamma(p)$ defines a $\Gamma$-invariant Voronoi tessellation of $\matH^n$, and every polyhedron of the tessellation is a fundamental domain. The polyhedron $D(p)$ of the tessellation containing $p$ is called the \emph{Dirichlet domain} for $\Gamma$ centred at $p$. By construction we have \index{Dirichlet domain}
$$D(g(p)) = g(D(p))$$
for all $g\in \Gamma$. We have proved in particular the following.

\begin{prop}
Every discrete group $\Gamma < \Iso(\matH^n)$ has a fundamental domain.
\end{prop}

Fundamental domains are far from being unique. The Dirichlet domain $D(p)$ depends on $p$ in a continuous fashion, and many fundamental domains are not Dirichlet domains.

\begin{figure}
\begin{center}
\includegraphics[width = 11 cm] {\iftoggle{BW}{ModularGroup-FundamentalDomain-01-BW}{ModularGroup-FundamentalDomain-01}}
\nota{The shadowed triangle (with one ideal vertex) is a fundamental domain for $\Gamma = \PSLZ$ acting on the half-space $H^2$. The translates of the fundamental domain form the tessellation shown.}
\label{modular_group:fig}
\end{center}
\end{figure}

\begin{ex}
Prove that the shadowed triangle in Figure \ref{modular_group:fig} is a fundamental domain for the modular group $\PSLZ$. 
\end{ex}

Everything we say also holds for $\matR^n$ and $S^n$.

\subsection{Fundamental domain of manifolds}
Let $\Gamma< \Iso(\matH^n)$ be a discrete subgroup that acts freely on $\matH^n$ and $D$ be a fundamental domain for $\Gamma$. We can get some information on the hyperbolic manifold $M=\matH^n/_\Gamma$ by looking at $D$. 

\begin{prop}
Let $M=\matH^n/_\Gamma$ be a hyperbolic manifold and $D$ a fundamental domain for $\Gamma$. The projection $\pi\colon \matH^n \to M$ restricts to a surjective map $D \to M$ that sends $\interior D$ isometrically onto an open dense subset of $M$. In particular we have
$$\Vol(D) = \Vol (M).$$
If $D$ is a Dirichlet domain, it is compact if and only if $M$ is.
\end{prop}
\begin{proof}
The translates $g(D)$ cover $\matH^n$, hence the projection $D \to M$ is surjective.
The translates $g(\interior D)$ are disjoint, hence $\interior D$ is sent isometrically inside $M$. The boundary $\partial D$ has measure zero, hence $\Vol(D) = \Vol(\interior D) = \Vol(M)$. 

If $D$ is compact then $M$ clearly is. If $D = D(p)$ is a Dirichlet domain and $M$ is compact, $M$ has finite diameter $\delta$ and hence every point in $\matH^n$ is at distance $\leqslant \delta$ of some point in the orbit of $p$. Hence $D(p)$ is contained in the closure of $B(p,\delta)$.
\end{proof}

A fundamental domain $D$ alone however does not determine $M$: we will soon see in Section \ref{piatte:subsection} that a square in $\matR^2$ is the fundamental domain of two non homeomorphic compact flat surfaces.

\subsection{Asteroids}
The facets of a fundamental domain $D$ are naturally partitioned into isometric pairs, as follows. Ever facet $f$ of $D$ is incident to $D$ and to another fundamental domain $g(D)$ of the tessellation. The isometry $g^{-1}$ sends $g(D)$ to $D$ and hence sends $f$ to another facet $f'$ of $D$. One checks immediately that $(f')' = f$, so $f$ and $f'$ are paired isometrically. 

One should think of $M$ as obtained from $D$ by identifying these facets in pairs: we can picture an observer -- say, a spaceship -- floating and traveling inside $M$ by visualising it in $D$, and jumping from $f$ to $f'$ every time it crosses the interior of a facet $f$, like in the 1979 video game \emph{Asteroids}.

\subsection{Spine and cut locus}
Another picture that may help to understand $M$ geometrically is the \emph{spine}\index{spine} $S \subset M$ defined as $S=\pi(\partial D)$. This is a $(n-1)$-dimensional object in $M$ whose complement is an open ball. When $D = D(p)$ is the Dirichlet domain of a point $p$, the spine $S$ is called the \emph{cut locus} of the point $q=\pi(p) \in M$. \index{cut locus}

\begin{ex}
The cut locus of $q \in M$ consists of all points $q'$ such that there are at least two geodesics of minimal length connecting $q$ to $q'$.
\end{ex}

We now exhibit some important examples of discrete groups of isometries in the three geometries $\matH^n$, $\matR^n$, and $S^n$, and study their fundamental domains. 

\subsection{Triangle groups} \label{triangle:groups:subsection}
Let $a,b,c \geqslant 2$ be natural numbers and $\Delta$ be the triangle in $\matH^2, S^2$, or $\matR^2$ with inner angles $\frac \pi a$, $\frac \pi b$, $\frac \pi c$. By reflecting iteratively $\Delta$ along its sides we get a tessellation $T$, see Section \ref{triangular:subsection}. \index{group!triangle group}

The \emph{triangle group} $\Gamma = \Gamma (a,b,c)$ is the group of isometries of $S^2$, $\matR^2$, or $\matH^2$ generated by the reflections $x,y,z$ along the three sides of $\Delta$ opposite to the vertices with inner angles $a,b,c$, respectively.

\begin{prop} \label{triangle:gamma:prop}
The triangle group $\Gamma (a,b,c)$ acts freely and transitively on the triangles of the tessellation $T$. Hence it is discrete and $\Delta$ is a fundamental domain for $\Gamma$. A presentation for the group is
$$\langle x, y,z\ |\ x^2, y^2, z^2, (xy)^c, (yz)^a, (zx)^b \rangle.$$
\end{prop}
\begin{proof}
We restrict for simplicity to the hyperbolic case, the others being analogous. It is convenient to construct the tessellation $T$ abstractly (this also furnishes a rigorous proof that by mirroring $\Delta$ along its sides we get a tessellation of $\matH^2$). 

We denote by $x,y,z$ both the sides of $\Delta$ and the reflections along them. For every $g\in \Gamma$ we define an abstract copy $\Delta_g$ of $\Delta$, and then we glue all these abstract copies altogether by pairing their sides  as follows: for every side of $\Delta$, say $x$, and every $g\in \Gamma$, we identify the two copies of the side $x$ in $\Delta_g$ and $\Delta_{gx}$ with the obvious identity map.

Via these identifications we get an abstract space $T$ tessellated into the triangles $\Delta_g$. The group $\Gamma$ acts freely and transitively on the tessellation (the element $g'\in \Gamma$ sends $\Delta_g$ to $\Delta_{g'g}$). We now prove that $T$ has a natural structure of a hyperbolic surface: since $\Gamma$ acts transitively, it suffices to check this for the points $p$ lying in $\Delta_e$. If $p$ lies in the interior of $\Delta_e$ or of a side, this is clear by construction. If $p$ is a vertex, say with inner angle $a$, by construction a cycle of $2a$ triangles $\Delta_e, \Delta_y, \Delta_{yz}, \Delta_{yzy}, \Delta_{yzyz}, \ldots $ is attached around $p$ because $(yz)^a = e$, and hence $T$ is naturally a hyperbolic surface also near $p$.

The hyperbolic surface $T$ is easily seen to be connected and complete. There is a natural developing map $\varphi\colon T\to \matH^2$ that sends $\Delta_g$ to $g(\Delta)$. The map $\varphi$ is a local isometry, hence a covering by Proposition \ref{local:isometry:prop}, hence an isometry since $\matH^2$ is simply connected. We identify $T$ with $\matH^2$ through $\varphi$.

It remains to prove that $\Gamma$ may be presented as stated. Let $G$ be the group presented as
$$\langle x, y,z\ |\ x^2, y^2, z^2, (xy)^c, (yz)^a, (zx)^b \rangle.$$
There is a natural surjection $G\to \Gamma$. We may repeat all the arguments above using $G$ instead of $\Gamma$ and get another $T'$ isometric to $\matH^2$ that covers $T$. We deduce that $T'=T$ and $G=\Gamma$.
\end{proof}

The triangle group $\Gamma(2,2,c) \isom \matZ/_{2\matZ}\times D_{2c}$ has order $4c$ and is the symmetry group of a prism. The triangle groups $\Gamma(2,3,3)$, $\Gamma(2,3,4)$, and $\Gamma(2,3,5)$ have order $24$, $48$, and $120$ and are the symmetry groups of the regular tetrahedron, octahedron, and icosahedron: the reader is invited to check all these facts visually by looking at Figure \ref{triangular:fig}.

The triangle groups $\Gamma(2,3,6)$, $\Gamma(2,4,4)$, and $\Gamma(3,3,3)$ are discrete groups of isometries of $\matR^2$ with compact quotient $\matR^2/_\Gamma$. 
If $\frac 1a + \frac 1b + \frac 1c <1$ the group $\Gamma(a,b,c)$ is a discrete subgroup of $\Iso(\matH^2)$. It contains infinitely many elliptic elements, such as reflections along lines and finite-order rotations around the vertices of the triangles. By Selberg's Lemma, there is a torsion-free subgroup $\Gamma' < \Gamma$ of some finite index $h$. The quotient $\matH^2/_{\Gamma'}$ is a closed hyperbolic surface and is tessellated into $h$ triangles isometric to $\Delta$.

\subsection{Coxeter polyhedra}
Triangle groups may be generalised to all dimensions as follows. A polyhedron $P$ in $\matH^n$ (or $\matR^n$, $S^n$) is a \emph{Coxeter polyhedron} if the dihedral angle of every codimension-two face divides $\pi$.\index{polyhedron!Coxeter polyhedron} For instance, the regular ideal tetrahedron and octahedron are Coxeter polyhedra.

The following theorem generalises Proposition \ref{triangle:gamma:prop}. Let $P$ be a finite Coxeter polyhedron: it is the convex hull of finitely many vertices in $\overline{\matH^n}$ (or $\matR^n$, $S^n$) and has $k$ facets, that we number as $1, \ldots, k$; we denote by $r_i$ the reflection along the $i$-th facet, and by $\frac {\pi}{a_{ij}}$ the dihedral angle formed by the $i$-th and $j$-th facets, if they are incident. Let $\Gamma$ be the isometry group generated by the reflections along the facets of $P$.

\begin{teo}
By mirroring iteratively a finite Coxeter polyhedron $P$ along its facets we get a tessellation of $\matH^n$ (or $\matR^n, S^n$).
The group $\Gamma$ acts freely and transitively on the tessellation: hence it is discrete and $P$ is a fundamental domain for $\Gamma$. A presentation for $\Gamma$ is
$$\langle r_1, \ldots, r_k \ |\  r_i^2, (r_ir_j)^{a_{ij}} \rangle$$
where $i$ varies in $1,\ldots, k$ and the pair $i,j$ varies among the incident facets.
\end{teo}
\begin{proof}
Same proof as Proposition \ref{triangle:gamma:prop}, with in addition an induction on the dimension $n$. Here are the details. We consider for simplicity only the hyperbolic case.

For every $g\in \Gamma$ we define an abstract copy $P_g$ of $P$, and we identify the $i$-th facet of $P_g$ and $P_{gr_i}$ for all $g\in \Gamma$ and all $i$. To prove that the resulting space $T$ is naturally a hyperbolic manifold, we use the induction hypothesis as follows. Every $p\in P_e$ lies in the interior of some $h$-dimensional face $f$, and let $\Gamma_f<\Gamma$ be the subgroup generated by the reflections along all the facets of $P$ containing $f$. The point $p$ is adjacent in $T$ to the polyhedra $P_g$ such that $g\in \Gamma_f$. By intersecting each such $P_g$ with a small codimension $h+1$ sphere centred in $p$ and contained in the codimension $h$ subspace orthogonal to $f$, we get some spherical polytope $Q_g$ of dimension $n-h-1<n$. (The spherical polytope $Q_t$ is usually called the \emph{spherical link} of $f$.)

By construction $Q_e \subset S^{n-h-1}$ is a Coxeter polytope and the subgroup $\Gamma_f < \Iso(S^{n-h-1})$ is generated by the reflections along its facets. By applying the induction hypothesis on $Q_e$ we deduce that the polytopes $Q_g$ with $g\in \Gamma_f$ form a tessellation of $S^{n-h-1}$, hence the polyhedra $P_g$ incident to $p$ form naturally a hyperbolic ball locally near $p$. Therefore $T$ is a hyperbolic manifold.

Proving that $T$ is complete requires a bit of care if $P_e$ has some ideal vertex $v$. In that case $v$ has a Euclidean link (obtained by intersecting $P_e$ with a small horosphere) and we conclude by induction as above that the Euclidean links of the $P_g$ incident to $v$ in $T$ glue to form a Euclidean $\matR^{n-1}$, that is a horosphere, so by intersecting each $P_g$ with small horoballs we get a horoball near $v$, which is complete.

The hyperbolic space $T$ is complete and connected, and we conclude as in the proof of Proposition \ref{triangle:gamma:prop}. 
\end{proof}

A group generated by some reflections along hyperplanes is called a \emph{reflection group}.\index{group!reflection group}
The following proposition shows that Coxeter polyhedra generate all the interesting reflection groups.

\begin{prop}
Every discrete reflection group $\Gamma$ is generated by the reflections along the facets of some Coxeter polyhedron.
\end{prop}
\begin{proof}
Consider the mirror hyperplanes of all the reflections in $\Gamma$. Since $\Gamma$ is discrete, these form a locally finite set and hence define a tessellation of $\matH^n$ onto which $\Gamma$ acts transitively. Pick one polyhedron $P$ of the tessellation. The reflections along the facets of $P$ generate $\Gamma$ (exercise).
\end{proof}

Coxeter polyhedra are beautiful objects that can be used to construct hyperbolic manifolds: every finite Coxeter polyhedron $P$ generates a reflection group $\Gamma$ which contains, by Selberg's Lemma, a torsion-free subgroup $\Gamma'$ of some finite index $h$. The quotient $M=\matH^n/_{\Gamma'}$ is a hyperbolic manifold and is tessellated into $h$ copies of $P$, so that $\Vol(M) = h\Vol(P)$. By residual finiteness, there are plenty of such manifolds.

\subsection{Coxeter graphs}
If a Coxeter polyhedron $P$ has a reasonable number of facets (for instance, if $P$ is a simplex), then one can study many of its properties by looking at its \emph{Coxeter graph},\index{Coxeter graph} which is constructed as follows: draw one node for each facet of $P$, and for every pair of facets intersecting with dihedral angle $\frac \pi a$ connect the corresponding nodes with an edge labeled with $a$. One also uses thickened or dashed edges to denote non-incident faces that do or do not intersect at some ideal vertex, respectively (these cases cannot occur on simplexes). Given the abundance of right angles, one usually omit the edges with label 2. 

\begin{ex}
Let $P$ be a regular polytope or tessellation with Schl\"afli symbol $\{a,b, \ldots, z\}$. By quotienting $P$ via all its isometries we get a Coxeter simplex, called its \emph{characteristic simplex}, whose Coxeter graph is\index{simplex!characteristic simplex}
\begin{center}
\includegraphics[width = 3 cm]{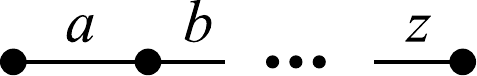}
\end{center}
\end{ex}

\subsection{Flat tori} \label{piatte:subsection}
The isometry group $\Iso(\matR^n)$ of $\matR^n$ contains the translations subgroup that we identify with $\matR^n$. The integer translations $\Gamma = \matZ^n$ form a discrete subgroup that acts freely on $\matR^n$, therefore the quotient $\matR^n /_\Gamma$ is a flat manifold. It is naturally diffeomorphic to the $n$-dimensional torus: 
$$\matR^n/_{\matZ^n} = \left(\matR/_{\matZ}\right)^n \isom \underbrace{S^1  \times \ldots \times S^1}_n.$$

\begin{ex} For every $p\in \matR^n$ the Dirichlet domain $D(p)$ is a $n$-dimensional unit cube centred at $p$.
\end{ex}

The flat $n$-torus may be seen as the unit $n$-cube with its opposite facets identified by translations. The two-dimensional case is shown in Figure \ref{domini:fig}-(left): by identifying the opposite sides of a square we get a torus. \index{flat torus}

\begin{figure}
\begin{center}
\includegraphics[width = 11.5 cm] {\iftoggle{BW}{domini-BW}{domini}}
\nota{A fundamental domain in $\matR^2$ for the torus (left) and the Klein bottle (centre): opposite sides should be identified as indicated by the arrows. A fundamental domain in $S^2$ for $\matRP^2$ (right).}
\label{domini:fig}
\end{center}
\end{figure}

The $n$-torus possesses a continuous family of non-isometric flat metrics. A \emph{lattice} $\Gamma<\matR^n$ is a discrete subgroup isomorphic to $\matZ^n$ which spans $\matR^n$ as a vector space. We see $\Gamma$ as a group of translations. \index{lattice}

\begin{ex} For every lattice $\Gamma$, the flat manifold $\matR^n/_\Gamma$ is diffeomorphic to the $n$-torus. A fundamental domain is the parallelotope spanned by $n$ generators of $\Gamma$.
\end{ex}

A Dirichlet domain for $\Gamma$ is almost never a parallelotope! In the \emph{hexagonal torus} $\matR^2/_\Gamma$ the group $\Gamma$ is the equilateral lattice generated by the translations\index{torus!hexagonal torus}
$$(1,0), \quad \big(\tfrac 12, \tfrac{\sqrt 3}2\big)$$ 
and the Dirichlet domain of any point is a regular hexagon. Orientable flat manifolds are easily classified (up to diffeomorphism) in dimension two. \index{hexagonal torus}

\begin{prop} \label{flat:is:torus:prop}
Every closed flat orientable surface is a torus.
\end{prop}
\begin{proof}
A closed flat surface $S$ is isometric to $\matR^2/_\Gamma$ for some discrete group $\Gamma$ of orientation-preserving isometries acting freely. Every fixed-point-free orientation-preserving isometry of $\matR^2$ is a translation (exercise), hence $\Gamma < \matR^2$. If $\Gamma$ has rank one then $S$ is not compact, hence $\Gamma$ is a lattice and $S$ is a torus.
\end{proof}

Flat manifolds more complicated than tori can be constructed by considering also non-translational isometries of $\matR^n$, which exist for $n\geqslant 3$ and also for $n=2$ in the non-orientable setting. For instance, in dimension two we find the \emph{Klein bottle} by taking $\Gamma$ as the group generated by the isometries:
$$\tau\colon(x,y) \mapsto (x+1, y), \qquad \eta\colon(x,y) \mapsto (1-x, y+1).$$
A fundamental domain for the Klein bottle is shown in Figure \ref{domini:fig}-(centre).
\index{Klein bottle}

Every subgroup $\Gamma' < \Gamma$ furnishes another flat manifold $\matR^2/_{\Gamma'}$ that covers the Klein bottle. For instance, the subgroups $\langle \tau \rangle$ and $\langle \eta \rangle$ generated respectively by $\tau$ and $\eta$ are both isomorphic to $\matZ$ but provide two different coverings: the manifold $\matR^2/_{\langle\tau\rangle}$ is an infinite cylinder while $\matR^2/_{\langle\eta\rangle}$ is an infinite M\"obius strip. 

The subgroup $\Gamma'$ generated by the translations $\tau$ and $\eta^2$ is isomorphic to $\matZ^2$ and has index 2 in $\Gamma$. The Klein bottle is doubly covered by the flat torus $\matH^2/_{\Gamma'}$. A fundamental domain for $\Gamma'$ is a rectangle with vertices $(0,0)$, $(1,0)$, $(0,2)$, $(1,2)$.

We have seen that the Klein bottle is covered by a torus. In fact, we will see in Section \ref{Bieberbach:section} that \emph{every} closed flat $n$-manifold is covered by a flat $n$-torus.

\begin{ex}
A square $Q\subset \matR^2$ is a Coxeter polygon and determines a reflection group $\Gamma$. What is the minimum index of a torsion-free subgroup $\Gamma' < \Gamma$? Which topological surfaces $\matR^2/_{\Gamma'}$ do we get?
\end{ex}

\subsection{Real projective spaces} 
We now exhibit some elliptic manifolds. Every elliptic manifold is covered by $S^n$ and is hence compact and has finite fundamental group (because coverings between compact manifolds have finite degree). \index{$\matRP^n$}

An important example is the real projective space $\matRP^n = S^n/_\Gamma$ where $\Gamma$ is the order-two group generated by the antipodal map $\iota (x) = -x$. 

\begin{ex}
For every $p\in S^n$, the Dirichlet domain $D(p)$ of $\Gamma$ is the hemisphere centred at $p$.
\end{ex}
The two-dimensional case is shown in Figure \ref{domini:fig}-(right). The only elliptic surfaces are $S^2$ and $\matRP^2$ in virtue of the following.
\begin{prop}
In even dimension $n$ the only elliptic manifolds are $S^n$ and $\matRP^n$.
\end{prop}
\begin{proof}
Let $M=S^n/_\Gamma$ be an elliptic manifold. Every matrix in $\SO(n+1)$ has an eigenvalue $+1$ because $n$ is even (exercise) and hence acts in $S^n$ with a fixed point. Therefore $\Gamma$ contains no non-trivial orientation-preserving isometries, and since these form a group of index at most two, either $\Gamma$ is trivial or is generated by an orientation-reversing fixed-point free involution in $O(n+1)$, and $-I$ is the only such element (exercise).
\end{proof}

\subsection{Lens spaces} \label{lens:subsection}
In dimension two $\matRP^2$ and $S^2$ are the only elliptic manifolds, but in dimension three there are more. Let $p\geqslant 1$ and $q\geqslant 1$ be coprime integers and set $\omega = e^{\frac{2\pi i}p}$. We identify $\matR^4$ with $\matC^2$ and see $S^3$ as
$$S^3 = \big\{(z,w)\in \matC^2\ \big|\ |z|^2+|w|^2 = 1 \big\}.$$
The map
$$f(z,w) = (\omega z, \omega^qw)$$
is an isometry of $\matR^4$ because it consists of two simultaneous rotations on the coordinate planes $w=0$ and $z=0$. The map $f$ hence induces an isometry of $S^3$. It has order $p$ and none of its iterates $f, f^2, \ldots, f^{p-1}$ has a fixed point. Therefore the group $\Gamma = \langle f \rangle$ generated by $f$ acts freely on $S^3$, and is discrete because it is finite. \index{lens space} \index{$L(p,q)$}

We have constructed an elliptic manifold $S^3/_\Gamma$, called a \emph{lens space} and indicated with the symbol $L(p,q)$. Its fundamental group is isomorphic to $\Gamma = \matZ/_{p\matZ}$. Note that the manifold depends on both $p$ and $q$.

\section{Geodesic boundary, non-complete, and cone manifolds}
We now introduce three important variations on the complete hyperbolic manifolds theme.
The first consists of admitting a totally geodesic boundary, the second is a brief overlook of some phenomena that may occur when the completeness hypothesis is dropped, and in the third we allow some conical singularity on a codimension-two geodesic stratum.

\subsection{Hyperbolic manifolds with geodesic boundary}
We reformulate a definition of hyperbolic (elliptic, flat) manifolds that allows the presence of some geodesic boundary. These manifolds are useful because they can be glued along their boundaries to produce new hyperbolic (elliptic, flat) manifolds.
\index{hyperbolic manifold!with geodesic boundary}
\index{flat manifold!with geodesic boundary}
\index{elliptic manifold!with geodesic boundary}
\begin{defn} A \emph{hyperbolic} (\emph{elliptic, flat}) \emph{manifold} $M$ \emph{with geodesic boundary}
is a Riemannian manifold with boundary where every point has an open neighbourhood isometric to an open set in a half-space in $\matH^n$ ($S^n$, $\matR^n$).
\end{defn}

The boundary $\partial M$ of a hyperbolic (elliptic, flat) $n$-manifold with geodesic boundary is a hyperbolic (elliptic, flat) $(n-1)$-manifold without boundary. Theorem \ref{simply:teo} extends appropriately to this context.

\begin{teo}
Every complete simply connected hyperbolic $n$-manifold $M$ with geodesic boundary is isometric to the intersection of some half-spaces in $\matH^n$ with disjoint boundaries.
\end{teo}
\begin{proof}
We construct a developing map $D\colon M\to\matH^n$ as in Theorem \ref{simply:teo}. The map $D$ is injective, because every two points $p,q \in M$ are connected by a geodesic (exercise: use completeness and geodesic boundary), which is sent to a geodesic via $D$, so $p,q$ are sent to distinct points $D(p), D(q)$. 

Therefore $D(M)\subset \matH^n$ is isometric to $M$. Since it is complete and with geodesic boundary, its boundary consists of disjoint hyperplanes, and hence $D(M)$ is the intersection of half-spaces bounded by them.  
\end{proof}

An example is sketched in Figure \ref{universale:fig}. Every complete hyperbolic manifold with geodesic boundary can be enlarged to a hyperbolic manifold without boundary in a canonical way.

\begin{figure}
\begin{center}
\includegraphics[width = 4cm] {\iftoggle{BW}{universale-BW}{universale}}
\nota{An intersection of (possibly infinitely many!) half-planes. The universal cover of a hyperbolic surface with boundary is isometric to such an object.}
\label{universale:fig}
\end{center}
\end{figure}

\begin{cor} Every complete hyperbolic $n$-manifold $M$ with geodesic boundary is contained in a complete hyperbolic $n$-manifold $N$ without boundary, such that $N\setminus \interior M$ is diffeomorphic to $\partial M \times [0,+\infty)$.
\end{cor}
\begin{proof}
The proof of Proposition \ref{corto:prop} applies and shows that $M = \tilde M/_\Gamma$ where $\tilde M\subset \tilde \matH^n$ is an intersection of half-spaces and $\Gamma< \Iso(\tilde M)$ acts freely and properly discontinuously on $\tilde M$.

Every local isometry in $\matH^n$ extends to a global isometry, therefore $\Gamma < \Iso(\tilde M) < \Iso(\matH^n)$. The group $\Gamma$ acts freely on $\tilde M$ and hence also on $\matH^n$: if it had a fixed point $x\in \matH^n$, it would fix also the unique point $y \in \tilde M$ that is the closest to $x$.
The manifold $N=\matH^n/_\Gamma$ contains naturally $M$.

For every $p \in \partial M$, let $\gamma_p(t)$ be the unit speed geodesic starting at $p$ orthogonal to $\partial M$ and directed outward. By looking at the universal cover one sees easily that the map $\partial M \times [0,+\infty) \to N \setminus \interior M$, $(p,t)\mapsto \gamma_p(t)$ is a diffeomorphism (we use here that $\tilde M\subset \matH^n$ is convex).
\end{proof}

\subsection{Cut and paste}
Hyperbolic manifolds with geodesic boundary are useful because they can be glued to produce new hyperbolic manifolds.

Let $M$ and $N$ be hyperbolic manifolds with geodesic boundary and $\psi\colon \partial M \to \partial N$ be an isometry. Let $M\cup_\psi N$ be the topological space obtained by quotienting the disjoint union $M\sqcup N$ by the equivalence relation that identifies $p$ to $\psi(p)$ for all $p\in \partial M$.

\begin{prop} \label{geodesic:glue:prop}
The space $M\cup_\psi N$ has a natural structure of hyperbolic manifold.
\end{prop}
\begin{proof}
The interiors of $M$ and $N$ inherit their hyperbolic metrics. When we glue $p$ to $\psi(p)$ we get a point $q$ that has two half-disc neighbourhoods on its sides, which glue to a honest hyperbolic disc, inducing a hyperbolic metric near $q$.
\end{proof}

Conversely, if an orientable hyperbolic manifold $M$ contains an orientable geodesic hypersurface $N$, we can cut $M$ along $N$ to get a hyperbolic manifold with geodesic boundary. The boundary will consist of two copies of $N$.

\subsection{Non-complete hyperbolic manifolds} \label{non-complete:subsection}
There is no classification of simply connected non-complete hyperbolic manifolds: for instance, we may get plenty of uninteresting examples by removing complicated closed sets from $\matH^n$. However, the first part of the proof of Theorem \ref{simply:teo} still applies and provides the following:

\begin{prop}
Let $M$ be a non-complete simply connected hyperbolic $n$-manifold. There is a local isometry
$$D\colon M \to \matH^n$$
which is unique up to post-composing with $\Iso(\matH^n)$.
\end{prop}
\begin{proof}
Construct $D$ as in the proof of Theorem \ref{simply:teo}: the completeness of $M$ is used there only in the last paragraph to show that $D$ is a covering.
As a local isometry, the map $D$ is determined by its first-order behaviour at any point $p\in M$, and is hence unique up to post-composing with an isometry of $\matH^n$.
\end{proof}

The map $D$ is called a \emph{developing map} and is neither injective nor surjective in general. We can push-forward isometries along $D$ as follows. \index{developing map}

\begin{prop}
Let $M$ be a hyperbolic $n$-manifold and $D\colon M \to \matH^n$ be a local isometry. For every $g\in \Iso(M)$ there is a unique $\rho(g)\in \Iso(\matH^n)$ such that
$$\rho(g) \circ D = D\circ g.$$
The resulting map $\rho \colon \Iso(M) \to \Iso(\matH^n)$ is a homomorphism.
\end{prop}
\begin{proof}
Pick a point $p\in M$ and define $\rho(g)$ as the unique isometry of $\matH^n$ such that $\rho(g) \circ D = D\circ g$ on a small neighbourhood of $p$. This equality between local isometries holds locally and hence globally on $M$. The map $\rho$ is easily checked to be a homomorphism.
\end{proof}

The homomorphism $\rho$ is the \emph{holonomy} associated to $D$. If $M$ is a non-complete hyperbolic $n$-manifold, its universal covering $\pi\colon \tilde M \to M$ inherits a non-complete hyperbolic metric. Therefore we get a developing map 
$$D\colon \tilde M \to \matH^n$$ 
together with a holonomy \index{holonomy}
$$\rho \colon \Aut(\pi) \longrightarrow \Iso(\matH^n)$$
which is the restriction of $\rho$ to the subgroup $\Aut(\pi)<\Iso(\tilde M)$. By fixing a point in $\tilde M$ we identify $\Aut(\pi)$ with $\pi_1(M)$ and get a holonomy
$$\rho\colon \pi_1(M) \longrightarrow \Iso(\matH^n).$$

As every metric space, a non-complete Riemannian manifold $M$ has a unique \emph{completion} $\overline M$, a complete metric space that contains $M$ as an open dense set. The completion $\overline M$ is however not necessarily a manifold, except in some lucky cases. These lucky cases are extremely important in dimension three, as we will see in Chapter \ref{three:chapter}. In that chapter we will understand the completion $\overline M$ of a hyperbolic three-manifold by studying its developing map and holonomy.

\subsection{Singularities} \label{singularities:subsection}
A particular class of non-complete hyperbolic manifolds deserves our attention, because their completions are some nice and natural objects: manifolds with \emph{cone angles} along some codimension-two singular geodesic stratum.

Let $S\subset \matH^n$ be a codimension-two subspace. The incomplete hyperbolic manifold $\matH^n \setminus S$ has fundamental group $\matZ$ and we denote by
$$\pi \colon X \longrightarrow \matH^n \setminus S$$
its universal covering: here $X$ is an interesting non-complete simply connected hyperbolic manifold. Note that $\pi$ may be interpreted as a developing map $\pi\colon X\to \matH^n$ which is neither injective nor surjective.

\begin{ex}
The metric completion of $X$ is $\overline X = X \sqcup \tilde S$ where $\tilde S$ is an identical copy of $S$. The covering $\pi$ extends to a surjective map
$$\pi\colon \overline X \to \matH^n$$
that sends $\tilde S$ to $S$. If $p\in \tilde S$ and $q\in \overline X$ then $d(p,q) = d(\pi(p), \pi(q)).$
\end{ex}

The space $\tilde S$ is the \emph{singularity} of $\overline X$ and should be interpreted as lying in $\overline X$ with infinite cone angle: the space $\overline X$ looks like a book with infinitely many pages (it is not locally compact at $\tilde S$). Singularities with finite cone angles will be defined in the next section as quotients of this model. The dimensions that will be interesting for us are of course $n=2$ and $3$, where $\tilde S$ is a point and a line, respectively.

\subsection{Cone manifolds}
For every $\theta\in \matR$ there is a well-defined rotation $\Rot_\theta\colon \overline X \to \overline X$ of angle $\theta$ around $\tilde S$, which projects to a rotation in $\matH^n$ of angle $\theta$ around $S$; note that $\Rot_\theta \neq \id$ for all $\theta \neq 0$, including $\theta = 2k\pi$.

We fix $\theta \neq 0$, so that $\Rot_\theta$ generates a free cyclic group $\Gamma < \Iso(\overline X)$.

\begin{ex}
The metric space $\overline X/_\Gamma$ is the completion of $X/_\Gamma$ and is homeomorphic to $\matR^n$.
\end{ex}

The image of $\tilde S$ along the quotient $\pi\colon \overline X \to \overline X /_\Gamma$ is a copy of $\tilde S$ and should now be interpreted as a singularity with cone angle $\theta$. When $\theta = 2k\pi$ for some integer $k$ the map $\pi\colon \overline X \to \matH^n$ factors as
$$\overline X \overset {\pi_1} {\longrightarrow} \overline X/_\Gamma \overset {\pi_2} \longrightarrow \matH^n.$$
When $k=\pm 1$ the map $\pi_2$ is an isometry. We introduce the following definition.\index{hyperbolic manifold!with cone angles}

\begin{defn}
A \emph{hyperbolic manifold $M$ with cone angles} is a metric space with charts in some $\overline X/_\Gamma$ and transition maps that are isometries.
\end{defn}

The term ``isometries'' should be interpreted in the strongest sense: the transition maps send singular points to singular points, and outside of the singular points these are isometries of Riemannian manifolds.

The points in $M$ that are mapped along charts to some singular set form a hyperbolic codimension-two submanifold in $M$ with some cone angle $\theta$, and different components may have different cone angles.
Points with cone angle $2\pi$ may be considered as ordinary points, while points with cone angle different from $2\pi$ are \emph{singular} and form the \emph{singular locus} of $M$.

The singular locus of $M$ is a geodesic codimension-two manifold, and its complement is an ordinary hyperbolic manifold, whose metric completion is $M$.
Flat and spherical manifolds with cone angles are defined in the same way. \index{elliptic manifold!with cone angles} \index{flat manifold!with cone angles}

\subsection{Examples}
The following simple construction is already a source of many two-dimensional examples: let $P$ be any polygon in $\matH^2, \matR^2,$ or $S^2$, with some inner angles $\alpha_1,\ldots, \alpha_k$. By doubling $P$ along its edges we get a topological sphere with $k$ cone points of angles $2\alpha_1, \ldots, 2\alpha_k < 2\pi$. This is a hyperbolic, flat, or spherical surface with cone angles.

Another simple but maybe more intriguing construction is the following: pick any polyhedron $P$ in $\matH^3, \matR^3$, or $\matS^3$. The faces of $P$ are (hyperbolic, flat, or elliptic) polygons that glue isometrically along the edges: hence $\partial P$ is a topological sphere with a natural structure of (hyperbolic, flat, or elliptic) manifold with cone angles at the vertices.

To construct examples in dimension three we need a bit more work because polyhedra have more strata. Non-compact examples are easier to build: let $P\subset \matH^3$ be any ideal hyperbolic polyhedron, with some $k$ edges with dihedral angles $\alpha_1,\ldots, \alpha_k$; the double of $P$ is naturally a hyperbolic cone manifold where the edges form the singular locus with cone angles $2\alpha_1, \ldots, 2\alpha_k$.

\section{Orbifolds} 
An orbifold is an object locally modelled on finite quotients of $\matR^n$. It naturally arises when we quotient a Riemannian manifold by a discrete group of isometries that may not act freely. Orbifolds behave like manifolds on many aspects.

\subsection{Definition} Let $\Gamma < \On(n)$ be a finite group of linear isometries and $V\subset \matR^n$ be a $\Gamma$-invariant open set. The resulting map
$$\varphi \colon V \longrightarrow V/_\Gamma$$
is called a \emph{local orbifold model}. For instance, we may pick $V=B(0,r)$. \index{orbifold!local model of an orbifold}

We now generalise manifolds by allowing local orbifold models in the atlas.\index{orbifold}
\begin{defn} Let $O$ be a Hausdorff paracompact topological space. An \emph{orbifold atlas} on $O$ is an open covering $\{U_i\}_{i\in I}$ of $O$, closed by finite intersections and equipped with local orbifold models 
$$\varphi_i\colon V_i \longrightarrow V_i/_{\Gamma_i} = U_i.$$
The local models are connected by some appropriate transition functions: for every inclusion $U_i\subset U_j$ there is an injective homomorphism 
$$f_{ij}\colon \Gamma_i \hookrightarrow \Gamma_j$$ 
and a $\Gamma_i$-equivariant smooth embedding $\psi_{ij} \colon V_i \hookrightarrow V_j$ compatible with the local models, that is 
$$\varphi_j \circ \psi_{ij} = \varphi_i.$$
\end{defn}

Two such atlases are equivalent if the are contained in some bigger atlas. The space $O$ equipped with an orbifold atlas is an \emph{orbifold} of dimension $n$.
See an example in Figure \ref{orbifold_charts:fig}. 

\begin{figure}
\begin{center}
\includegraphics[width = 8 cm] {\iftoggle{BW}{orbifold_charts-BW}{orbifold_charts}}
\nota{An orbifold $O$ with three local models $O=U_1 \supset U_2 \supset U_3$. The group $\Gamma_1 = \matZ_2\times \matZ_2$ acts on $V_1$ by reflections on the two coordinate axis, the group $\Gamma_2 = \matZ_2$ acts on $V_2$ by reflection on the horizontal axis, the group $\Gamma_3$ is trivial.}
\label{orbifold_charts:fig}
\end{center}
\end{figure}

\begin{oss}
We think at the maps $\psi_{ij}$ and $f_{ij}$ as defined only up to the action of $\Gamma_j$ (which acts on $\psi_{ij}$ by composition and on $f_{ij}$ by conjugation).  In particular, if $U_i\subset U_j \subset U_k$ then we can verify that the equalities $\psi_{ik} = \psi_{jk}\circ \psi_{ij}$ and $f_{ik} = f_{jk}\circ f_{ij}$ hold only up to this ambiguity. See Figure \ref{orbifold_charts:fig}.
\end{oss}
The \emph{isotropy group} $\Gamma_p$ of a point $p\in O$ is the stabiliser of any lift of $p$ in any local model $V$ with respect to the action of $\Gamma$. By definition $\Gamma_p$ is a finite subgroup of $\On(n)$. A point $p$ is \emph{regular} if its isotropy group is trivial, and \emph{singular} otherwise. \index{orbifold!isotropy group of an orbifold}

\begin{example} 
The isotropy group $O_p$ in the orbifold $O$ from Figure \ref{orbifold_charts:fig} is 
$\matZ_2 \times \matZ_2$ at the origin, $\matZ_2$ at the axis, and trivial elsewhere. 
\end{example}

\begin{prop}
The regular points in an orbifold form a dense subset. 
\end{prop}
\begin{proof}
On a local model $V \to V/_\Gamma$, a singular point is the image of the fixed point locus of some element in $\Gamma$, which is in turn a proper subspace.
\end{proof}

An orbifold is \emph{locally oriented} if all the finite groups $\Gamma$ lie in $\SO(n)$, and it is \emph{oriented} if, in addition, all the transition maps $\psi_{ij}$ are orientation-preserving. In a locally oriented orbifold reflections are not admitted and hence the singular locus has codimension $\geqslant 2$.
An open subset of an orbifold is naturally an orbifold.

\subsection{Examples} Orbifolds are natural objects and there are plenty of nice examples around.

\begin{example} A differentiable manifold is an orbifold whose points are all regular. A differentiable manifold with boundary may be interpreted as an orbifold whose boundary points have the local structure of type $\matR^n/_\Gamma$ where $\Gamma = \matZ_2$ is generated by a reflection along a hyperplane. The boundary should be interpreted as a \emph{mirror}.
\end{example}

\begin{ex}
Construct an orbifold structure on the triangle as suggested by Figure \ref{orbifold_triangle:fig}.
\end{ex}

\begin{example}
Let $\Gamma < \On(n)$ be finite and $V\subset \matR^n$ a $\Gamma$-invariant open set. The quotient $V/_\Gamma$ has an orbifold structure, defined by the unique local model $V \to V/_\Gamma$. 
\end{example}

\begin{figure}
\begin{center}
\includegraphics[width = 2.2 cm] {\iftoggle{BW}{orbifold_triangle-BW}{orbifold_triangle}}
\nota{An orbifold structure on the triangle constructed using the local models from Figure \ref{orbifold_charts:fig}. The isotropy group is $\matZ_2 \times \matZ_2$ at the vertices, $\matZ_2$ at the edges, and trivial in the interior.}
\label{orbifold_triangle:fig}
\end{center}
\end{figure}

We generalise the last example.

\begin{prop}
If $M$ a Riemannian manifold and $\Gamma<\Iso(M)$ is a discrete subgroup, the quotient $M/_\Gamma$ has a natural orbifold structure.
\end{prop}
\begin{proof}
Take a point $p\in M/_\Gamma$ and $\tilde p\in M$ a lift. Since $\Gamma$ is discrete, the stabiliser $\Gamma_{\tilde p}$ of $\tilde p$ is finite and there is a $r>0$ such that $\exp_{\tilde p}(B_r(0)) = B_r(\tilde p)$ and $g(B_r(\tilde p))$ intersects $B_r(\tilde p)$ for some $g\in \Gamma$ if and only if $g\in\Gamma_{\tilde p}$. 

The ball $B_r(\tilde p)$ is clearly $\Gamma_{\tilde p}$-invariant. The group $\Gamma_{\tilde p}$ acts linearly and orthogonally on $B_r(0)$ and we get an orbifold local model
$$B_r(0) \longrightarrow B_r(0) /_{\Gamma_{\tilde p}} \isom B_r(\tilde p)/_{\Gamma_{\tilde p}} = U_p$$
where the diffeomorphism is induced by the exponential map. We have constructed an orbifold local model $U_p$ around $p$. We extend the covering $\{U_p\}$ thus obtained by adding all the non-empty intersections. 
\end{proof}

This is the richest source of nice examples. The quotient map $M \to M/_\Gamma$ is a covering, in an appropriate sense that we now explain.

\subsection{Coverings}
Let $\Gamma' < \Gamma < \On(n)$ be finite groups and $V\subset \matR^n$ a $\Gamma$-invariant open set. The natural map
$$\varphi\colon V/_{\Gamma'} \longmapsto V/_{\Gamma}$$
between the two orbifolds is a \emph{local covering}. \index{orbifold!orbifold covering}

\begin{defn}
A continuous map $p\colon \tilde O \to O$ between orbifolds is a \emph{covering} if every point $p\in O$ has a neighbourhood $U$ with $p^{-1}(U) = \sqcup_{i \in I} U_i$ and every restriction $p|_{U_i}\colon U_i \to U$ is a local covering.
\end{defn}

The following is the main source of examples.

\begin{example}
Let $M$ be a Riemannian manifold and $\Gamma' < \Gamma < \Iso(M)$ be discrete groups. The natural map $M/_{\Gamma'} \to M/_\Gamma$ is an orbifold covering.
\end{example}

\begin{defn}[The good, the bad, and the very good]
An orbifold is \emph{good} if it is covered by a manifold, and it is \emph{bad} otherwise. It is \emph{very good} if it is finitely covered by a manifold. \index{orbifold!the bad, the good, and the very good orbifold}
\end{defn}

\subsection{Hyperbolic, flat, and elliptic orbifolds} \label{more:examples:subsection}
We define a \emph{hyperbolic, flat}, or \emph{elliptic} orbifold to be an orbifold whose local models are isometric quotients of balls in $\matH^n$, $\matR^n$, or $S^n$, and whose transition functions are also isometries. \index{orbifold!hyperbolic, flat, elliptic orbifold}

The quotient of a hyperbolic, flat, or spherical manifold by a discrete group of isometries is naturally a hyperbolic, flat, or elliptic orbifold.

\begin{example}
The triangle group $\Gamma(a,b,c)$ introduced in Section \ref{triangle:groups:subsection} is a discrete group of isometries of $S^2, \matR^2$, or $\matH^2$. It defines a triangular orbifold $\Delta$, which is hyperbolic, flat, or elliptic, according to whether $\frac 1a+\frac 1b + \frac 1c$ is smaller, equal, or bigger than 1. The isotropy groups are the dihedral $D_{2a}, D_{2b}$, and $D_{2c}$ at the vertices, $\matZ_2$ at the edges, and trivial elsewhere.

The index-two subgroup $\Gamma^{\rm or}(a,b,c) \triangleleft \Gamma(a,b,c)$ consisting of orientation-preserving isometries is sometimes called a \emph{von Dyck group}. It defines an orientable orbifold $O$ that double-covers $\Delta$. The orbifold $O$ is a topological sphere with three singular points with rotation isotopy groups $\matZ/_{a\matZ}, \matZ/_{b\matZ}, \matZ/_{c\matZ}$. See Figure \ref{orbifold_triangle_abc:fig}.
\index{group!von Dyck group}
\end{example}

\begin{ex} \label{VD:ex}
Show that the following is a presentation for $\Gamma^{\rm or}(a,b,c)$:
$$\langle\ r, s, t\ |\ r^{a}, s^{b}, t^{c}, rst\ \rangle.$$
The triple $r,s,t$ of generators is intrinsically determined in the group up to simultaneous conjugation or inversion. Different unordered triples $(a,b,c)$ produce non-isomorphic Von Dyck groups $\Gamma^{\rm or}(a,b,c)$.
\end{ex}
\begin{proof}[Hints]
The group $\Gamma^{\rm}(a,b,c)$ preserves the tessellation $T$ of $\matH^2$, $\matR^2$, or $S^2$ into triangles with angles $\frac \pi a, \frac \pi b, \frac \pi c$, so its finite-order elements are rotations along some vertices of $T$. Three rotations $r,s,t$ that satisfy $rst=1$ are of some very special kind.
\end{proof}

\begin{example}
More generally, every hyperbolic (or flat, elliptic) Coxeter polyhedron is a hyperbolic (or flat, elliptic) orbifold.
\end{example}

\begin{example}
Quotient the flat torus $T = \matR^2/_{\matZ^2}$ by the \emph{elliptic involution} $(x,y)\mapsto (-x,-y)$. This isometry of $T$ has four fixed points and the quotient flat orbifold is a sphere with four singular points, see Figure \ref{elliptic_involution:fig}.
\end{example}

\begin{figure}
\begin{center}
\includegraphics[width = 8 cm] {\iftoggle{BW}{orbifold_triangle_abc-BW}{orbifold_triangle_abc}}
\nota{The index-two orientation-preserving $\Gamma^{\rm or}(a,b,c) \triangleleft \Gamma(a,b,c)$ defines an index-two orbifold covering $O \to \Delta$, where $O$ is a 2-sphere with three singular points and $\Delta$ a triangle.}
\label{orbifold_triangle_abc:fig}
\end{center}
\end{figure}

\begin{figure}
\begin{center}
\includegraphics[width = 9 cm] {\iftoggle{BW}{elliptic_involution-BW}{elliptic_involution}}
\nota{The elliptic involution quotients the torus to a sphere with four singular points. The fixed points and their images are drawn in \iftoggle{BW}{grey}{red}. We show both a planar (left) and spacial (right) picture.}
\label{elliptic_involution:fig}
\end{center}
\end{figure}

\begin{figure}
\begin{center}
\includegraphics[width = 10 cm] {\iftoggle{BW}{orbifold_Farey-BW}{orbifold_Farey}}
\nota{The hyperbolic orbifold $\matH^2/_{\PSLZ}$ is obtained by mirroring the sides of a triangle with inner angles $\frac \pi 3$, $\frac \pi 2$, and zero (left). The orbifold is topologically a punctured sphere with two singular points with rotational isotropy $\matZ_2$ and $\matZ_3$ (right).}
\label{orbifold_Farey:fig}
\end{center}
\end{figure}

There are also many interesting non-compact orbifolds.

\begin{example}
A fundamental domain for $\PSLZ$ was shown in Figure \ref{modular_group:fig}. By gluing its sides according to the action we see quite easily that the hyperbolic orbifold $H^2/_{\PSLZ}$ is non-compact and contains two rotational points of order 2 and 3, see Figure \ref{orbifold_Farey:fig}.
\end{example}

\begin{ex} The 1-dimensional compact orbifolds are $S^1$ and the segment $[0,1]$ with mirrored endpoints. Show that there are coverings $S^1 \to [0,1]$ of degree 2, and also coverings $[0,1] \to [0,1]$ of any positive degree.
\end{ex}

When $O=X/_\Gamma$ and $X$ is a simply connected manifold, we say that $\Gamma$ is the \emph{fundamental group} of $O$. It is also possible to define the fundamental group for a generic orbifold using the appropriate orbi-notions of paths and homotopies. 

\subsection{Cone manifolds}
Hyperbolic (flat, elliptic) cone manifolds and orbifolds are different objects, but they have a wide common intersection.

A hyperbolic (flat, elliptic) orbifold $O$ whose isotropy groups $O_p$ are either trivial or generated by a $\frac{2\pi }p$-rotation along a codimension-two subspace is naturally a hyperbolic (flat, elliptic) cone manifold with cone angles $\frac{2 \pi}p$. The singular set may consist of various connected components and the natural number $p$ depends on the component. 

Conversely, every hyperbolic (flat, elliptic) manifold with cone angles that divide $2\pi$ can be given a natural hyperbolic (flat, elliptic) orbifold structure whose singular set consists of rotational points only. Moreover, we can easily prove that it is a good one:

\begin{prop} \label{from:cone:to:orbifold:prop}
Every hyperbolic (flat, elliptic) manifold $M$ with cone angles that divide $2\pi$ is a good hyperbolic (flat, elliptic) orbifold. If $M$ is complete, then $M=\matH^n/_\Gamma$ ($\matR^n/_\Gamma$, $S^n/_\Gamma$) for some discrete group $\Gamma$ of isometries.
\end{prop}
\begin{proof}
We consider for simplicity only the hyperbolic case. We must prove that $M$ is orbifold-covered by a hyperbolic manifold.

The singular set $S\subset M$ decomposes into connected components $S = \sqcup_i S_i$, each with a cone angle $\frac{2\pi}{p_i}$ for some integer $p_i \geqslant 2$. The complement $M' = M \setminus S$ is a non-complete hyperbolic manifold and as such it is equipped with a developing map $D\colon \tilde M' \to \matH^n$ and a holonomy $\rho\colon \pi_1(M') \to \Iso(\matH^n)$. 

We pick a loop $\mu_i$ around each component $S_i$ of $S$ and note that by definition $\rho(\mu_i)$ is an elliptic isometry that rotates $\matH^n$ around a codimension-two subspace by the angle $\frac{2\pi}{p_i}$. In particular $\rho(\mu_i)$ has order $p_i$. 

We consider the regular covering $N \to M'$ corresponding to the subgroup $\ker \rho$ and give $N$ the hyperbolic structure induced by $M'$. The completion $\overline N$ of $N$ is a hyperbolic manifold that orbifold-covers $M$: since $\rho(\mu_i)$ has order precisely $p_i$, the added singular points in $\overline N$ have cone angle $2\pi$ and are hence ordinary. 

If $M$ is complete, then $\overline N$ also is and hence $\overline N = \matH^n/_{\Gamma'}$. We define $\Gamma<\Iso(\matH^n)$ to be the group generated by $\Gamma'$ and any lifts of the deck transformations of $\overline N\to M$, and we get $M=\matH^n/_\Gamma$.
\end{proof}

\begin{cor} \label{very:good:cor}
If $M$ is complete and has finitely generated fundamental group, then it is very good. 
\end{cor}
\begin{proof}
We know that $M=\matH^n/_\Gamma$ and we apply Selberg's Lemma.
\end{proof}

Summing up, complete hyperbolic cone manifolds with angles that divide $2\pi$ can be promoted to good orbifolds and hence to quotients $\matH^n/_\Gamma$. We will encounter many examples in dimensions two and three later on.

\subsection{References}
All the material presented here is standard and introduced in various books, like Benedetti -- Petronio \cite{BP} and Ratcliffe \cite{R}. A beautiful introduction to regular polytopes is Coxeter's 1963 book \cite{Cox}, while the reader interested in Coxeter polytopes may consult Vinberg's survey \cite{Vin}.

The theory of orbifolds and cone manifolds presented here is quite limited: to get more background, the reader may consult Thurston's notes \cite[Chapter 13]{Th} for the orbifolds, and McMullen's paper \cite{McM-cone} for the cone manifolds.

%% file: Thick-thin.tex
\chapter{Thick-thin decomposition} \label{thick-thin:chapter}
A peculiar aspect of hyperbolic geometry is the existence of complete hyperbolic manifolds that have finite volume but are not compact. These manifolds behave like the compact ones in many aspects, but are sometimes easier to construct. They arise naturally when we study hyperbolic surfaces, and are a fundamental constituent of Thurston's geometrisation of three-manifolds.

We prove here a structure theorem for all such manifolds. The theorem says that every finite-volume complete hyperbolic manifold can be decomposed into two domains: a \emph{thick part} which is compact and has injectivity radius bounded from below, and a \emph{thin part} that consists of \emph{cusps}. A cusp is a (truncated) hyperbolic manifold of type $N\times [0,+\infty)$, where every section $N\times t$ has a flat metric that shrinks exponentially with $t$. 

The core of this theorem is a lemma about discrete subgroups of Lie groups called the \emph{Margulis Lemma}. We will apply this lemma also to the other geometries and prove Bierbach's Theorem, that states that every compact flat manifold is covered by a torus.

Throughout the discussion we will also study some general aspects of finite-volume complete hyperbolic manifolds, concerning in particular closed geodesics, isometry groups, and finite covers.

\section{Tubes and cusps} \label{tubes:cusps:section}
We introduce here two very simple kinds of hyperbolic manifolds, so simple that their fundamental groups will be called \emph{elementary} in the sequel: the tubes and the cusps. Before introducing them we show that the injectivity radius of a hyperbolic manifold behaves nicely.

We also study the closed geodesics in a hyperbolic manifold $M$, and prove that there is precisely one in every free homotopy class of closed curves of hyperbolic type.

\subsection{Injectivity radius}
Like many other geometric properties, the injectivity radius of a complete hyperbolic manifold may be nicely observed by looking at its universal cover. If $S\subset \matH^n$ is a discrete set, we define $d(S)$ as the infimum of $ d(x_1,x_2)$ among all pairs $x_1, x_2$ of distinct points in $S$. 

Recall that every complete hyperbolic manifold is a quotient $M=\matH^n/_\Gamma$ by some group $\Gamma$ of isometries acting freely on $\matH^n$. 

\begin{prop} \label{d:prop}
Let $M=\matH^n/_\Gamma$ be a complete hyperbolic manifold and $\pi\colon \matH^n \to M$ the projection. For every $x\in M$ we have
$$\inj_xM = \frac 12 \cdot d(\pi^{-1}(x)).$$
\end{prop}
\begin{proof}
The number $\inj_x M$ is the supremum of all $r>0$ such that $B(x,r)$ is isometric to a ball of radius $r$ in $\matH^n$. The open set $B(x,r)$ is a ball of radius $r$ if and only if its counterimage via $\pi$ consists of disjoint balls of radius $r$, and this holds $\Leftrightarrow$ two distinct points in $\pi^{-1}(x)$ stay at distance at least $2r$.
\end{proof}

Recall that we denote by $d(\gamma)$ the minimum displacement of $\gamma \in \Iso(\matH^n)$.

\begin{cor} \label{inj:cor}
Let $M=\matH^n/_\Gamma$ be a complete hyperbolic manifold. Then
$$\inj M = \frac 12 \cdot \inf \big\{d(\gamma)\ \big|\ \gamma \in \Gamma, \gamma \neq \id\big\}.$$
\end{cor}
\begin{proof}
We have
$$\inj_xM = \frac 12 \cdot d(\pi^{-1}(x)) = \frac 12 \cdot \inf \big\{ d(\tilde x, \gamma(\tilde x)) \ \big|\ \gamma\in \Gamma, \gamma\neq \id, \tilde x \in \pi^{-1}(x)\big\}. $$
Therefore
$$\inj M = \inf_{x\in M} \inj_x M = 
\frac 12 \cdot \inf \big\{ d(\gamma)\ \big| \ \gamma \in \Gamma, \gamma \neq \id \big\}.$$
The proof is complete.
\end{proof}

\begin{cor}
If $M=\matH^n/_\Gamma$ is a compact hyperbolic manifold then every non-trivial element in $\Gamma$ is hyperbolic.
\end{cor}
\begin{proof}
Every nontrivial element in $\Gamma$ is either hyperbolic or parabolic.
If $M$ is compact then $\inj M >0$. If $\Gamma$ contains a parabolic $\gamma$ then $d(\gamma)=0$ and hence $\inj M =0$.
\end{proof}

\subsection{Tubes} \label{tubes:subsection}
Consider the infinite cyclic group $\Gamma = \langle \varphi \rangle$ generated by a hyperbolic transformation $\varphi$ on $\matH^n$ with axis $l$ and minimum displacement $d>0$. The iterates $\varphi^k$ are again hyperbolic transformations with axis $l$ and displacement $kd$. Therefore $\Gamma$ acts freely on $\matH^n$.
The quotient manifold $M=\matH^n/_\Gamma$ is called an \emph{infinite tube}. \index{tube}

\begin{ex}
Fix $q\in l$. Let $q_1,q_2$ be the two points in $l$ at distance $\frac d2$ from $q$ and $\pi_1,\pi_2$ the two hyperplanes orthogonal to $l$ in $q_1,q_2$.
The Dirichlet domain $D(q)$ of $\Gamma$ is the space comprised between $\pi_1$ and $\pi_2$.
\end{ex}

The infinite tube $M = \matH^n/_\Gamma$ is obtained from $D(q)$ by identifying $\pi_1$ and $\pi_2$ along $\varphi$. Its fundamental group is isomorphic to $\Gamma \Isom \matZ$. The axis $l$ projects in $M$ onto a closed geodesic $\gamma$ of length $d$. We have $\inj M = \frac d2 $ by Corollary \ref{inj:cor} and the points in $\gamma$ are precisely those with minimum injectivity radius. The closed geodesic $\gamma$ is sometimes called the \emph{core geodesic} of the tube.

\begin{prop}
Every infinite tube is diffeomorphic to $S^1\times \matR^{n-1}$ or $S^1 \timtil \matR^{n-1}$ according to whether $\varphi$ is orientation-preserving or not.
\end{prop}
\begin{proof}
By projecting $\matH^n$ orthogonally onto $l$, we give $\matH^n$ the structure of a $\matR^{n-1}$-bundle over $l$ which is preserved by $\varphi$ and hence descends to a structure of $\matR^{n-1}$-bundle over $\gamma$. The conclusion follows from the classification of vector bundles over $S^1$, see Proposition \ref{S1:bundle:prop}.
\end{proof}

A \emph{tube of radius $R$}, or a \emph{$R$-tube}, is the quotient $N_R(l)/_\Gamma$ of the closed $R$-neighbourhood $N_R(l)$ of $l$, shown in Figure \ref{tube:fig}. It is diffeomorphic to $D^{n-1} \times S^1$ or $D^{n-1}\timtil S^1$, and in particular it is compact. Note that the boundary of a tube is not totally geodesic, see Figure \ref{tube:fig}.

\begin{figure}
\begin{center}
\includegraphics[width = 4 cm] {\iftoggle{BW}{tube-BW}{tube}}
\nota{The $R$-neighbourhood $N_R(l)$ of a vertical line in the half-space model is a Euclidean cone. The boundary of the cone is not totally geodesic: the cone is however convex.}
\label{tube:fig}
\end{center}
\end{figure}

\subsection{Cusps}
We now introduce another simple type of hyperbolic manifolds. In the previous example the discrete group $\Gamma$ consisted of hyperbolic transformations fixing the same line $l$, now $\Gamma$ will consist of parabolic transformations fixing the same point at infinity.

Let $\Gamma < \Iso(\matR^{n-1})$ be a non-trivial discrete group of Euclidean isometries acting freely on $\matR^{n-1}$: the quotient $M=\matR^{n-1}/_\Gamma$ is a flat $(n-1)$-manifold. If we use the half-space model for $\matH^n$ with coordinates $(x,t)$, every element $\varphi \in \Gamma$ acts as a parabolic transformation on $\matH^n$ by sending $(x,t)$ to $(\varphi (x),t)$. The whole group $\Gamma$ is a discrete group of parabolic transformations of $\matH^n$ fixing the point $\infty$. 

The quotient $\matH^n /_\Gamma$ is naturally diffeomorphic to $M \times \matR_{>0}$. The metric tensor at the point $(x,t)$ is 
$$g_{(x,t)} = \frac{g^M_x \oplus 1}{t^2}$$
where $g^M$ is the metric tensor of the flat $M$. The manifold $\matH^n/_\Gamma$ is called a \emph{cusp}. Since $\Gamma$ contains parabolics we have $\inj M = 0$. \index{cusp}

\begin{oss}
The vertical coordinate $t$ may be parametrized more intrinsically using arc-length. As we have seen in Proposition \ref{verticale:prop}, a vertical geodesic with unit speed is parametrized as $t=e^u$. Using $u$ instead of $t$ the cusp is isometric to $M\times \matR$ with metric tensor 
$$g_{(x,u)} = (e^{-2u}g_x^M) \oplus 1.$$ 
The lengths in the flat slice $M\times u$ are shrunk or dilated by the factor $e^{-u}$.
\end{oss}

\begin{figure}
\begin{center}
\includegraphics[width = 5 cm] {\iftoggle{BW}{Pseudosphere-BW2}{Pseudosphere}}
\nota{The \emph{pseudosphere} is a surface in $\matR^3$ isometric to the union of two truncated cusps, each with constant gaussian curvature $-1$.}
\label{pseudosphere:fig}
\end{center}
\end{figure}

A \emph{truncated cusp} is a portion $N=M \times [a,+\infty)$, bounded by the Euclidean manifold $M\times a$: note that the boundary $\partial N$ is Euclidean but not totally geodesic.
The volume of a truncated cusp is strikingly simple. \index{cusp!truncated cusp}
\begin{prop} \label{truncated:volume:prop}
Let $N$ be a truncated cusp. We have
$$\Vol \big(N\big) = \frac{\Vol(\partial N)}{n-1}. $$
\end{prop}
\begin{proof}
It follows from Lemma \ref{O:lemma}.
\end{proof}

This shows in particular that a (non truncated) cusp has infinite volume.

\begin{example}
In dimension $n=2$ there is only one cusp up to isometry. The group $\Gamma<\Iso(\matR)$ is the infinite cyclic group generated by a translation $x\mapsto x+b$ and up to conjugating in $\Iso(\matH^2)$ we may take $b=1$. The cusp is diffeomorphic to $S^1\times \matR$, and the circle $S^1\times \{u\}$ has length $e^{-u}$. Some truncated cusp (but not the whole cusp!) embeds in $\matR^3$ as in Figure \ref{pseudosphere:fig}.\index{pseudosphere}
\end{example}
\begin{rem}
Pick $p\in \matH^2$. Note that a cusp and $\matH^2\setminus \{p\}$ are both hyperbolic and diffeomorphic to an open annulus $S^1\times \matR$. However, they are not isometric because the cusp is complete while $\matH^2\setminus \{p\}$ is not.
\end{rem}

We sometimes employ the word \emph{cusp} to indicate a truncated cusp, for simplicity.

\subsection{Closed geodesics}
A \emph{closed curve} in a manifold $M$ is a smooth map $\gamma\colon S^1 \to M$. A (possibly closed) curve is \emph{simple} if it is an embedding.\index{geodesic!closed geodesic}

We consider $S^1$ as the unit circle in $\matC$. A \emph{closed geodesic} in a Riemannian manifold $M$ is a smooth map $\gamma\colon S^1 \to M$ whose lift $\gamma\circ \pi\colon \matR \to M$ along the universal covering $\pi(t) = e^{i t}$ is a non-constant geodesic. 
Two closed geodesics $\gamma_1$, $\gamma_2$ that differ only by a rotation, \emph{i.e.}~such that $\gamma_1(z) = \gamma_2(ze^{it})$ for some fixed $t\in\matR$, are implicitly considered equivalent. The closed geodesics $\gamma(z)$ and $\overline\gamma(z) = \gamma(\bar z)$ are however distinct (they have opposite orientations).

\begin{prop} \label{closed:geodesic:prop}
Let $\gamma$ be a closed geodesic in a Riemannian manifold $M$. Exactly one of the following holds:
\begin{enumerate}
\item the curve $\gamma$ is simple,
\item the curve $\gamma$ self-intersects transversely in finitely many points,
\item the curve $\gamma$ wraps $k\geqslant 2$ times along a curve of type (1) or (2).
\end{enumerate}
\end{prop}
\begin{proof}
If the geodesic is not simple, it self-intersects. If it self-intersects only with distinct tangents, then (2) holds. Otherwise (3) holds.
\end{proof}

The natural number $k$ in (3) is the \emph{multiplicity} of the closed geodesic. A closed geodesic $\gamma$ of multiplicity $k$ is of type $\gamma (e^{it}) = \eta(e^{kit})$ for some geodesic $\eta$ of type (1) or (2). 

\begin{oss}
A closed geodesic on a Riemannian manifold $M$ is determined by its support, its orientation, and its multiplicity.
\end{oss}

\subsection{Closed geodesics in a hyperbolic manifold}
Closed geodesics in hyperbolic manifolds have a particularly nice behaviour. 

Let $X,Y$ be topological spaces: as usual we indicate by $[X,Y]$ the homotopy classes of continuous maps $X \to Y$. Let $X$ be path-connected. There is a natural map $\pi_1(X,x_0) \to [S^1, X]$, and the following is a standard exercise in topology.\index{$[X,Y]$} 

\begin{ex}
The map induces a bijection between the conjugacy classes in $\pi_1(X,x_0)$ and $[S^1,X]$.
\end{ex}
A simple closed curve in $X$ is \emph{homotopically trivial} if it is homotopic to a constant. As a corollary, a simple closed curve $\gamma$ is homotopically trivial if and only if it represents the trivial element in $\pi_1(X,\gamma(1))$.

On a complete hyperbolic manifold $M=\matH^n/_\Gamma$ the conjugacy classes of $\pi_1(M)$ correspond to those in $\Gamma$ and we get a correspondence
$$\big\{ {\rm conjugacy\ classes\ in\ }\Gamma\big\} \longleftrightarrow [S^1, M].$$
This correspondence works as follows: given $\varphi \in \Gamma$, pick any point $x\in\matH^n$, connect $x$ to $\gamma(x)$ with any arc, and project the arc to a closed curve in $M$.

Two conjugate elements in $\Gamma$ are of the same type (trivial, parabolic, or hyperbolic) and have the same minimum displacement. Therefore every element in $[S^1, M]$ has a well-defined type and minimum displacement.

\begin{prop} \label{unique:closed:geodesic:prop}
Let $M$ be a complete hyperbolic manifold. Every hyperbolic element of $[S^1, M]$ is represented by a unique closed geodesic, of length $d$ equal to its minimum displacement. The trivial and parabolic elements are not represented by closed geodesics.
\end{prop}
\begin{proof}
Take $M=\matH^n/_\Gamma$. A hyperbolic isometry $\varphi \in \Gamma$ has a unique invariant geodesic in $\matH^n$, its axis, which projects to a closed geodesic of length $d$. Conjugate isometries determine the same closed geodesic in $M$.

On the other hand, a closed geodesic in $M$ lifts to a segment connecting two distinct points $x_0$ and $\varphi(x_0)$ for some $\varphi\in \Gamma$ which preserves the line passing through $x_0$ and $\varphi(x_0)$: since $\varphi$ preserves a line, it is hyperbolic.
\end{proof}
We get a bijection
$$\big\{{\rm hyperbolic\ conjugacy\ classes\ in\ }\Gamma\big\} \longleftrightarrow \big\{ {\rm closed\ geodesics\ in\ } M \big\}.$$

\begin{cor} \label{unique:closed:geodesic:cor}
Let $M$ be a closed hyperbolic manifold. Every non-trivial element in $[S^1, M]$ is represented by a unique closed geodesic.
\end{cor}
\begin{proof}
Since $M$ is compact there are no parabolics.
\end{proof}

\begin{cor} \label{minimum:cor}
Let $M$ be a complete hyperbolic manifold. Every closed geodesic has the minimum length among the closed curves in its homotopy class.
\end{cor}
\begin{proof}
If $\alpha$ is a closed geodesic, its length equals the minimum displacement $d$ of a corresponding hyperbolic transformation $\varphi$. Every closed curve $\beta$ homotopic to $\alpha$ lifts to an arc connecting two points $\tilde x$ and $\varphi(\tilde x)$ that have distance at least $d$; hence $\beta$ has length at least $d$.
\end{proof}

Not only a purely topological object like a hyperbolic homotopy class of closed curves has a unique geometric nice representative, but this representative is the shortest possible one. Informally, we may think that every closed curve may be shrunk until it becomes a closed geodesic, and the negative curvature forces this closed geodesic to be unique: we will soon see that the uniqueness is lost in the elliptic and flat geometries.

It is worth recalling that a closed geodesic may not be simple. When it is simple, we now show that the closed geodesic has some nice small neighbourhoods. We defined the $R$-tubes in Section \ref{tubes:subsection}.

\begin{prop} \label{tube:prop}
The $R$-neighbourhood of a simple closed geodesic $\gamma$ in a complete hyperbolic manifold is isometric to a $R$-tube, for some $R>0$.
\end{prop}
\begin{proof}
By compactness of $\gamma$ there is a sufficiently small $R>0$ such that the $R$-neighbourhood of $\gamma$ lifts to disjoint $R$-neighbourhoods of its geodesic lifts in $\matH^n$. Hence their quotient is a $R$-tube.
\end{proof}

\section{The Margulis Lemma}
We state and prove the \emph{Margulis Lemma}, that concerns arbitrary discrete groups in Lie groups and more specifically in $\Iso(\matH^n)$. The lemma implies that there is a constant $\varepsilon>0$, which depends only on the dimension $n$, such that the $\varepsilon$-thin part of any complete hyperbolic $n$-manifold consists of truncated cusps and tubes only. The $\varepsilon$-thin part is by definition the set of all the points with injectivity radius smaller than $\varepsilon$.

The proof of the lemma for general Lie groups is surprisingly simple and elegant; its application to the hyperbolic case needs however a more technical argument. We start by exposing some preliminary facts on the isometries of the hyperbolic space that are of independent interest.

\subsection{Isometries that commute or generate discrete groups} \label{that:commute:subsection}
Two isometries of $\matH^n$ that commute or generate a discrete group must be of a particular kind. We let $\Fix(\varphi)$ denote the fixed points in $\overline{\matH^n}$ of an isometry $\varphi$. 

\begin{lemma} \label{commute:lemma}
Let $\varphi_1,\varphi_2\in\Iso(\matH^n)$ be two hyperbolic or parabolic isometries. If they commute then $\Fix(\varphi_1) = \Fix(\varphi_2)$.
\end{lemma}
\begin{proof}
If they commute, the map $\varphi_1$ acts on $\Fix(\varphi_2)$ and vice versa. If $\varphi_2$ is hyperbolic, then $\Fix(\varphi_2)=\{p,q\}$ and $\varphi_1$ fixes the line with endpoints $p$ and $q$, hence is again hyperbolic with $\Fix(\varphi_1)=\{p,q\}$. If $\varphi_1$ and $\varphi_2$ are parabolic then they have the same fixed point $\Fix(\varphi_1) = \Fix(\varphi_2)$.
\end{proof}

\begin{lemma} \label{fissi:lemma}
Let $\varphi_1,\varphi_2\in\Gamma$ be two non-trivial isometries in a discrete group $\Gamma < \Iso (\matH^n)$ that acts freely on $\matH^n$. Either $\Fix(\varphi_1)\cap \Fix(\varphi_2) = \emptyset$, or one of the following holds:
\begin{itemize}
\item $\varphi_1$ and $\varphi_2$ are parabolics with the same fixed point,
\item $\varphi_1$ and $\varphi_2$ are powers of the same hyperbolic $\varphi \in \Gamma$.
\end{itemize}
\end{lemma}
\begin{proof}
Suppose that $\varphi_1$ and $\varphi_2$ have some common fixed points. We first consider the case where $\varphi_1$ is hyperbolic and $\varphi_2$ is parabolic. We pick the half-space model and suppose $\Fix(\varphi_1) = \{0,\infty\}$ and $\Fix(\varphi_2) = \{\infty\}$. 

Proposition \ref{tre:prop} says that 
$$\varphi_1(x,t) = \lambda(Ax,t), \quad \varphi_2(x,t) = (A'x+b, t)$$ 
with $A, A' \in \On(n-1)$ and $\lambda \neq 1$. Hence
\begin{align*}
\varphi_1^n\circ \varphi_2\circ  \varphi_1^{-n}(x,t) & = \varphi_1^n\big(A'(\lambda^{-n}A^{-n}x)+b,\lambda^{-n}t) \\
& = \big(A^nA'A^{-n}x + \lambda^n A^nb , t\big).
\end{align*}
Up to interchanging $\varphi_1$ and $\varphi_1^{-1}$ we may suppose $\lambda < 1$ and get
$$\lim_{n\to \infty} \varphi_1^n \circ \varphi_2 \circ \varphi_1^{-n}(0,t) = 
\lim_{n\to \infty} \left(\lambda^nA^nb, t\right) = (0,t).$$
A contradiction since $\Gamma$ is discrete.

We now suppose that both $\varphi_1$ and $\varphi_2$ are hyperbolic, with $\Fix(\varphi_1) = \{a,\infty\}$ and $\Fix (\varphi_2) = \{b,\infty\}$. The isometries $\varphi_1$ and $\varphi_2$ permute the horospheres centred at $\infty$ and
$$[\varphi_1,\varphi_2] = \varphi_1\circ \varphi_2 \circ \varphi_1^{-1} \circ \varphi_2^{-1} \in \Gamma$$
fixes every horizontal horosphere. Hence the commutator is parabolic or trivial: the first case is excluded by the previous discussion, in the second case we have $a=b$ by Lemma \ref{commute:lemma}. 

Both $\varphi_1$ and $\varphi_2$ have the same axis $l$, and since they generate a discrete group $\Gamma_0 < \Gamma$ they are both powers of some hyperbolic $\varphi\in\Gamma_0$ with that axis. To prove this, note that $\Gamma_0$ acts effectively on $l$ as a discrete group of translations, hence $\Gamma_0 \isom \matZ$. \end{proof}

Two isometries are contained in some discrete group if and only if they generate a discrete group, so the previous lemma is actually a fact on pairs of isometries that generate discrete groups.
The two lemmas have important geometric consequences.

\begin{cor}  \label{axis:cor}
Let $\matH^n/_\Gamma$ be a complete hyperbolic manifold. The axis in $\matH^n$ of two hyperbolic isometries in $\Gamma$ are either incident or ultra-parallel (not asymptotically parallel). 
\end{cor}

\begin{cor} \label{ZxZ:cor}
Let $\matH^n/_\Gamma$ be a complete hyperbolic manifold. Every subgroup of $\Gamma$ isomorphic to $\matZ\times\matZ$ consists of parabolic elements fixing the same point at infinity.
\end{cor}

\begin{cor} \label{no:ZxZ:cor}
The fundamental group $\pi_1(M)$ of a closed hyperbolic manifold $M$ does not contain subgroups isomorphic to $\matZ\times \matZ$.
\end{cor}

\begin{cor} \label{torus:no:hyperbolic:cor}
The torus has no hyperbolic structure.
\end{cor}

\subsection{The Margulis Lemma}
We now introduce the main character of this chapter, the Margulis Lemma. We first state a general version for Lie groups and then turn more specifically to $\Iso(\matH^n)$.

\begin{lemma} \label{G:lemma}
Let $G$ be a Lie group. There is a neighbourhood $U$ of $e\in G$ such that every discrete subgroup $\Gamma < G$ generated by some elements in $U$ is nilpotent.
\end{lemma}
Note that the discreteness of $\Gamma$ is essential here, since every connected Lie group $G$ is generated by arbitrarily small elements, and $G$ needs not to be nilpotent.

\begin{proof}
Consider the commutator map 
$$[\ ,\ ]\colon G\times G \to G$$ 
that sends $(g,h)$ to the commutator $[g,h]$. The map is smooth and sends $G\times \{e\}$ and $\{e\}\times G$ to the point $e$. Therefore its differential at $(e,e)$ vanishes. We identify a neighbourhood of $e$ in $G$ with $\matR^n$, so that the commutator map near $(e,e)$ can be read as 
$$[\ ,\ ]\colon U \times U \longrightarrow \matR^n$$
for some neighbourhood $U\subset \matR^n$ of the origin.
Since its differential at $(0,0)$ vanishes, up to restricting $U$ we may suppose that the map is $\frac 12$-Lipschitz. Therefore for all $x,y \in U$ we get
$$\|[x,y]\| < \|[x,0]\| + \frac 12\|y\| = \frac 12\|y\|, \quad \|[x,y]\| < \|[0,y]\| + \frac 12\|x\| = \frac 12\|x\|$$
and hence
$$\|[x,y]\| < \frac 12 \min \big\{\|x\|, \|y\|\big\}.$$
This implies that for every smaller neighbourhood $V\subset U$ of $0$ there is a $k>0$ such that
$$\underbrace{[U, [U, \ldots [U,U]]\cdots ]}_k \subset V.$$
We now turn back to $G$. Let $\Gamma$ be discrete and generated by some elements $S\subset U$. We choose a smaller neighbourhood $V$ of $e$ such that $V \cap \Gamma = \{e\}$, hence for any $a_1,\ldots, a_k,b\in S$ we get $[a_1,[a_2, \ldots[a_k,b]]\cdots ] \in V\cap\Gamma =\{e\}$. Proposition \ref{criterion:nilpotent:prop} says that $\Gamma$ is nilpotent.
\end{proof}

We now want to refine this lemma when $G=\Iso(\matH^n)$. Let $P$ be a property of groups, like being abelian, nilpotent, etc. A group is \emph{virtually $P$} if it has a finite-index subgroup which is $P$.\index{Margulis Lemma} 

\begin{lemma}[Margulis Lemma] \label{Margulis:lemma}
In every dimension $n\geqslant 2$ there is a constant $\varepsilon_n>0$ such that for all $x\in \matH^n$, every discrete group $\Gamma<\Iso(\matH^n)$ generated by elements that move $x$ at distance smaller than $\varepsilon_n$ is virtually nilpotent.
\end{lemma}
\begin{proof}
It suffices to prove the theorem for a fixed $x\in \matH^n$ since the isometries of $\matH^n$ act transitively on points. By Lemma \ref{G:lemma} there is a neighbourhood $U$ of $e$ in $G=\Iso(\matH^n)$ such that every discrete group generated by some elements in $U$ is nilpotent. 

\begin{figure}
\begin{center}
\includegraphics[width = 5 cm] {\iftoggle{BW}{Margulis-BW}{Margulis}}
\nota{The Margulis Lemma for Lie groups furnishes a nice neighbourhood $U$ of $e\in\Iso(\matH^n)$, but here we need a nice neighbourhood $W$ of the whole compact stabiliser $\On(n)$, drawn as a \iftoggle{BW}{dark grey}{blue} line.}
\label{Margulis:fig}
\end{center}
\end{figure}

Let $V_\varepsilon\subset\Iso(\matH^n)$ be the set of all the isometries that move $x$ at distance smaller than $\varepsilon>0$. It is relatively compact. As $\varepsilon>0$ varies, the sets $V_\varepsilon$ form a neighbourhood system for the compact stabiliser $G_x <\Iso(\matH^n)$ of $x$, isomorphic to $\On(n)$. Pick a $V_\varepsilon$. Being relatively compact, it can be covered by some $m-1$ translates of $U$. See Figure \ref{Margulis:fig}.

Let now $W=W^{-1}$ be a symmetric neighbourhood of $G_x$ satisfying $W^m \subset V_\varepsilon$. Let $\Gamma$ be a discrete group generated by some elements $S\subset W$. We want to prove that the nilpotent subgroup $\Gamma_U$ generated by $\Gamma \cap U$ has index at most $m$ in $\Gamma$. This concludes the proof, for it suffices to take a $\varepsilon_n$ such that $V_{\varepsilon_n} \subset W$.

Suppose by contradiction that $\Gamma/_{\Gamma_{U}}$ contains more than $m$ cosets. Let $r(i)$ be the number of cosets represented by elements in $\Gamma$ that are products of at most $i$ generators in $S$. If $r(i)=r(i+1)$ for some $i$, then it is easy to deduce that $r(i)$ stabilises forever. This easily implies that $r(m)\geqslant m$. Therefore there are $m$ elements in $W^m$ belonging to $m$ distinct cosets. By hypothesis two of them, say $g$ and $h$, belong to the same translate of $U$. Therefore $gh^{-1} \in \Gamma\cap U$, contradicting the fact that $g$ and $h$ lie in distinct cosets.
\end{proof}

For a discrete group $\Gamma < \Iso(\matH^n)$ and a point $x\in \matH^n$, we denote by $\Gamma_\varepsilon (x)< \Gamma$ the subgroup generated by all elements $g\in \Gamma$ that move $x$ of a distance smaller than $\varepsilon$. We have proved that $\Gamma_{\varepsilon_n} (x)$ is virtually nilpotent for every $x\in\matH^n$ and every discrete $\Gamma$.

A \emph{Margulis constant} in a fixed dimension $n$ is any number $\varepsilon_n>0$ for which Lemma \ref{Margulis:lemma} holds. Every sufficiently small number is a Margulis constant.\index{Margulis constant} 

\subsection{Elementary groups}
A non-trivial discrete group $\Gamma < \Iso(\matH^n)$ is \emph{elementary} if it preserves a finite set of points in $\overline{\matH^n}$.\index{group!elementary group}

\begin{prop} \label{elementary:prop}
An elementary $\Gamma$ acting freely on $\matH^n$ is:
\begin{itemize}
\item generated by a hyperbolic isometry, or
\item generated by parabolic isometries having the same fixed point at $\infty$.
\end{itemize}
\end{prop}
\begin{proof}
Since $\Gamma$ contains no elliptics, every non-trivial element fixes one or two points at infinity in $\overline{\matH^n}$, and no other finite set of points (because of Proposition \ref{iterate:prop}). We conclude using Lemma \ref{fissi:lemma}.
\end{proof}

\begin{prop} \label{elementary:index:prop}
Let $\Gamma<\Iso(\matH^n)$ be a discrete group acting freely on $\matH^n$. If $\Gamma'<\Gamma$ has finite index and is elementary, then $\Gamma$ also is.
\end{prop}
\begin{proof}
We know that every element in $\Gamma'$ is either hyperbolic with axis $l$, or parabolic with fixed point $p\in\partial\matH^n$. If $\varphi \in \Gamma$ then $\varphi^k\in \Gamma'$ for some $k$: hence $\varphi$ is also of that type.
\end{proof}

\begin{cor}
Let $\Gamma<\Iso(\matH^n)$ be discrete and acting freely. If $\Gamma$ is virtually nilpotent, it is either trivial or elementary.
\end{cor}
\begin{proof}
We know that $\Gamma$ contains a finite-index nilpotent subgroup $H$. If $H$ is trivial then $\Gamma$ is finite and hence trivial. If $H$ is non-trivial, it has a non-trivial centre by Proposition \ref{center:prop}. Lemma \ref{commute:lemma} then implies that all the elements in $H$ have the same fixed points and so $H$ is elementary. Therefore $\Gamma$ is elementary by Proposition \ref{elementary:index:prop}
\end{proof}

Let $\varepsilon_n$ be a Margulis constant. We can strengthen the Margulis Lemma.

\begin{cor}
Let $\Gamma<\Iso(\matH^n)$ be discrete and acting freely. For every point $x\in \matH^n$ the subgroup $\Gamma_{\varepsilon_n}(x)$ is either trivial or elementary.
\end{cor}

\subsection{Thick-thin decomposition}
A notable geometric consequence of the Margulis Lemma is that every complete hyperbolic manifold decomposes into a thick and a thin part, where the thin part consists only of particularly simple star-shaped pieces. This decomposition is known as the \emph{thick-thin decomposition} of hyperbolic manifolds.\index{thick-thin decomposition}

We define a \emph{star-shaped set} centred at $p\in\partial\matH^n$ to be any subset $U\subset\matH^n$ that intersects every half-line pointing to $p$ in a half-line. For instance, a horoball is star-shaped. A \emph{star-shaped neighbourhood} of a line $l\subset\matH^n$ is any neighbourhood $V$ of $l$ that intersects every line orthogonal to $l$ into a connected set. For instance, a $R$-neighbourhood of $l$ is star-shaped.

These definitions pass to quotients. A \emph{star-shaped cusp neighbourhood} is the quotient $U/_\Gamma$ of a $\Gamma$-invariant star-shaped set $U$ centred at $p\in\partial \matH^n$ via a discrete group $\Gamma$ of parabolic transformations fixing $p$ and acting freely. Analogously, a \emph{star-shaped simple closed geodesic neighbourhood} is the quotient $V/_\Gamma$ of a $\Gamma$-invariant star-shaped neighbourhood $V$ of $l$ via a discrete group $\Gamma\isom \matZ$ of hyperbolic transformations with axis $l$.

The truncated cusps and $R$-tubes studied in Section \ref{tubes:cusps:section} are particularly nice star-shaped cusp and geodesic neighbourhoods. 

Let $\varepsilon_n>0$ be a Margulis constant. We define $M_{[\varepsilon_n,\infty)}$ and $M_{(0,\varepsilon_n]}$ respectively as the set of all points $x \in M$ having $\inj_x M \geqslant \frac{\varepsilon_n}2$, and as the closure of the complementary set $M\setminus M_{[\varepsilon_n,\infty)}$. They form respectively the \emph{thick} and \emph{thin part} of $M$. 

\begin{oss} \label{exclude:thin:oss}
We do not define the thin part simply as the set of all points $x$ having $\inj_x M \leqslant \frac{\varepsilon_n}2$ because we want to discard the degenerate and unlucky case of a closed geodesic $\gamma$ having length precisely $\varepsilon_n$. The injectivity radius would be $\frac{\varepsilon_n}2$ at the points in $\gamma$ and strictly bigger than $\frac{\varepsilon_n}2$ near $\gamma$. With our definition, the geodesic $\gamma $ is contained in $M_{[\varepsilon_n, +\infty)}$ and not in $M_{(0,\varepsilon_n]}$.
\end{oss}

The following theorem is arguably the most important structural result on complete hyperbolic manifolds of any dimension $n$.

\begin{teo}[Thick-thin decomposition] Let $M$ be a complete hyperbolic $n$-manifold. The thin part $M_{(0,\varepsilon_n]}$ consists of a disjoint union of star-shaped neighbourhoods of cusps and of simple closed geodesics of length $< \varepsilon_n$.
\end{teo}
\begin{proof}
We have $M=\matH^n/_\Gamma$. For every isometry $\varphi\in \Gamma $ we define
$$S_\varphi(\varepsilon) = \big\{ x\in\matH^n\ \big|\ d(\varphi(x),x)\leqslant \varepsilon\big\}\subset \matH^n.$$
By Proposition \ref{d:prop} the thin part $M_{(0,\varepsilon_n]}$ is the image of the set 
\begin{align*}
 S & = \big\{x\in \matH^n\ \big|\ \exists \varphi\in \Gamma, \varphi\neq \id {\rm\ such\ that\ } d(\varphi(x),x)\leqslant \varepsilon_n\big\} \\
 & = \bigcup_{\varphi\in \Gamma, \varphi \neq \id} S_\varphi(\varepsilon_n).
 \end{align*}
More precisely, we should exclude the hyperbolic $\varphi \in \Gamma$ with $d(\varphi)=\varepsilon_n$, see Remark \ref{exclude:thin:oss}, but this is not an important point and we ignore it. It is easy to check that $S_\varphi(\varepsilon)$ is star-shaped, centred at a $p\in \partial \matH^n$ or at a line $l$ according to whether $\varphi$ is parabolic fixing $p$ or hyperbolic fixing $l$.

Suppose that $x\in S_\varphi(\varepsilon_n) \cap S_\psi(\varepsilon_n)$ for some non-trivial isometries $\varphi, \psi \in \Gamma$. By the Margulis Lemma both $\varphi$ and $\psi$ belong to the elementary group $\Gamma_{\varepsilon_n}(x)$ and hence by Proposition \ref{elementary:prop} both $\varphi$ and $\psi$ are either parabolic fixing the same point $p$ at infinity or hyperbolic fixing the same line $l$.

Therefore every connected component $S_0$ of $S$ is the union of all $S_\varphi(\varepsilon_n)$ where $\varphi$ varies in some maximal elementary subgroup $\Gamma_0<\Gamma$ of parabolics fixing the same point $p$ or hyperbolics fixing the same line $l$. The set $S_0$ is a union of star-shaped sets centred at $p$ or $l$ and is hence also star-shaped. 

The group $\Gamma$ preserves $S$ and the only isometries in $\Gamma$ that preserve $S_0$ are those in $\Gamma_0$, therefore the quotient $M_{(0,\varepsilon_n]} = S/_\Gamma$ consists of star-shaped neighbourhoods of cusps and of simple closed geodesics.
\end{proof}

Star-shaped neighbourhoods are particularly nice in low dimensions.

\begin{prop}
Let $M$ be a complete orientable hyperbolic manifold of dimension $n\leqslant 3$. The thin part $M_{(0,\varepsilon_n]}$ consists of truncated cusps and $R$-tubes.
\end{prop}
\begin{proof}
Pick a non-trivial $\varphi \in \Iso^+(\matH^n)$. We check that $S_\varphi(\varepsilon)$ is either empty, or a $R$-neighbourhood of a line $l$, or a horoball, for all $\varepsilon>0$. This proves the proposition.

Suppose that $\varphi$ is hyperbolic with axis $l$. The distance $d(x,\varphi(x))$ depends only on $d(x,l)$ since all the orientation-preserving hyperbolic transformations with axis $l$ commute with $\varphi$ and act transitively on the points at fixed distance from $l$. 
It is easy to check that $d(x,\varphi(x))$ increases with $d(x,l)$ and hence $S_\varphi(\varepsilon)$ is either empty or a $R$-neighbourhood of $l$.

If $\varphi$ is parabolic, it acts on each horosphere centred at $p$ like a
Euclidean fixed-point-free orientation-preserving isometry on $\matR^{n-1}$: this must be a translation when $n\leqslant 3$. Therefore $d(x,\varphi(x))$ depends only on the horosphere $O$ containing $x$ and decreases as $O$ moves towards $p$. Therefore $S_\varphi(\varepsilon)$ is a horoball.
\end{proof}

\begin{figure}
\begin{center}
\includegraphics[width = 11.5 cm] {\iftoggle{BW}{thick_thin-BW}{thick_thin}}
\nota{The thick-thin decomposition of a complete hyperbolic surface: the thin part (\iftoggle{BW}{light grey}{yellow}) consists of truncated cusps and neighbourhoods of short geodesics\iftoggle{BW}{}{ (blue)}.}
\label{thick_thin:fig}
\end{center}
\end{figure}

See the picture in Figure \ref{thick_thin:fig}. 

\begin{cor} \label{simple:small:geodesics:cor}
Let $M$ be a complete hyperbolic $n$-manifold. The closed geodesics in $M$ of length $< \varepsilon_n$ are simple and disjoint.
\end{cor}
\begin{proof}
These closed geodesics lie in the thin part. Star-shaped cusp neighbourhoods contain no closed geodesics, and each star-shaped geodesic neighbourhood contains only one closed geodesic, its core, which is simple.
\end{proof}

\subsection{Finite-volume hyperbolic manifolds} \label{finite-volume:subsection}
A specific feature of hyperbolic geometry is the existence of complete hyperbolic manifolds that have finite volume without being compact. The thin-thick decomposition furnishes a nice topological description of such manifolds.

\begin{prop}
A complete hyperbolic manifold $M$ has finite volume if and only if its thick part is compact.
\end{prop}
\begin{proof}
If the thick part is not compact, it contains an infinite number of points that stay pairwise at distance greater than $\varepsilon_n$. The open balls of radius $\frac{\varepsilon_n}2$ centred at these points are embedded and disjoint and all have the same volume: therefore their union has infinite volume.

If the thick part is compact, it has finite volume. Its boundary is also compact, and hence has finitely many connected components. Therefore the thin part
consists of finitely many star-shaped neighbourhoods of cusps and closed geodesics, each with compact boundary. Each such object has finite volume (because it is contained in a bigger abstract truncated cusp with compact base, or a $R$-tube, which has finite volume).
\end{proof}

\begin{cor} \label{tame:cor}
Every complete finite-volume hyperbolic manifold $M$ is diffeomorphic to the interior of a compact manifold $N$ with boundary. The boundary $\partial N$ consists of manifolds that admit some flat structure.
\end{cor}
\begin{proof}
We have seen in the previous proof that the thick part of $M$ is compact, and the thin part consists of finitely many star-shaped neighbourhoods of cusps and closed geodesics, each with compact boundary.

Every cusp neighbourhood has compact base and hence contains a smaller truncated cusp, diffeomorphic to $X\times [0,1)$ for some closed flat manifold $X$. The complement in $M$ of these truncated cusps is compact: it is obtained from the thick part by adding finitely many neighbourhoods of closed geodesics (that are compact) and compact portions of cusp neighbourhoods. Each truncated cusp is diffeomorphic to $X\times [0,1)$ and can be compactified by adding $X \times 1$. The resulting manifold $N$ is compact with boundary.
\end{proof}

Every boundary component $X$ of $N$ inherits a flat structure, uniquely determined up to rescaling (different truncations modify $X$ only by rescaling.) As already mentioned, a truncated cusp is often called simply a \emph{cusp}: for instance we say that the surface sketched in Figure \ref{thick_thin:fig} has two cusps.

We briefly discuss the effects of Corollary \ref{tame:cor} on low-dimensional manifolds. 
Let $M$ be a complete hyperbolic manifold of finite volume of dimension $n$. If $n=2$, the manifold $M$ is diffeomorphic to the interior of a compact surface $N$ with boundary. Every boundary component of $N$ is of course a circle.

If $n=3$, the manifold $M$ is diffeomorphic to the interior of a compact 3-manifold $N$ with boundary, and every boundary component $X$ of $N$ is a flat surface: if $M$ is orientable, then $\partial N$ also is, and by Proposition \ref{flat:is:torus:prop} it consists of tori.

\subsection{Geodesic boundary and cusps}
Most of the arguments of this chapter extend to hyperbolic manifolds with compact geodesic boundary.

\begin{prop} \label{geodesic:M:N:prop}
Every complete finite-volume hyperbolic manifold $M$ with compact geodesic boundary is diffeomorphic to a compact manifold $N$ with some boundary components removed. The removed components have a flat structure. The remaining components form the geodesic boundary of $M$ and hence have a hyperbolic structure.
\end{prop}
\begin{proof}
Double $M$ along the geodesic boundary to obtain a finite-volume hyperbolic manifold to which Corollary \ref{tame:cor} applies. 
\end{proof}

We briefly discuss the effects in low dimensions. In dimension $n=2$ the boundary $\partial N$ consists of circles: some are geodesic components of $\partial M$, while some others are removed and correspond to cusps. In dimension $n=3$, if $M$ is orientable the removed boundary of $\partial N$ consists of tori (the cusps), while the geodesic boundary of $M$ cannot contain any torus: a torus has no hyperbolic structure by Corollary \ref{torus:no:hyperbolic:cor}.

We add some information on fundamental groups. Let the manifolds $M$ and $N$ be as in Proposition \ref{geodesic:M:N:prop}.

\begin{prop} \label{pi1:M:prop}
For every boundary component $X$ of $N$ the homomorphism $\pi_1(X)\to\pi_1(N)$ induced by inclusion is injective. Two homotopically non-trivial closed curves in distinct boundary components of $N$ are not freely homotopic in $N$.
\end{prop}
\begin{proof}
We have $M=C/_\Gamma$ for some convex $C\subset \matH^3$ with boundary consisting of hyperplanes. If $X$ is a (cusp) flat component of $\partial N$, it is isometric to $\matR^{n-1}/_{\Gamma_X}$ and the map $\Gamma_X=\pi_1(X) \to \pi_1(N)=\Gamma$ is injective since it sends a Euclidean isometry to a corresponding parabolic isometry.
If $X$ is geodesic, it lifts to a hyperplane in $\partial C$, which is simply connected: hence $\pi_1(X) \to \pi_1(N)$ is injective.

A homotopically non-trivial closed curve $\gamma\subset X$ determines a conjugacy class of hyperbolic or parabolic transformations in $\Gamma$. Every element in this class preserves a unique boundary hyperplane or a point at infinity of $C$ that projects back to $X$, so the conjugacy class determines $X$. Therefore distinct boundary components contain different conjugacy classes. 
\end{proof}

We note that a finite-volume hyperbolic manifold $M$ may have non-compact geodesic boundary: an ideal polygon in $\matH^2$ is a simple example. 

\section{Geodesic spectrum, isometry groups, and finite covers} \label{geodesic:spectrum:section}
We now study the geodesic spectrum and the isometry group of hyperbolic manifolds. The geodesic spectrum of $M$ is the set of the lengths of all the closed geodesics in $M$. We prove that a finite-volume $M$ has a discrete geodesic spectrum (with finite multiplicities) and a finite isometry group.\index{geodesic spectrum}

Then we turn to finite covers and derive some geometric consequences from the residually finiteness of fundamental groups.

We are almost exclusively interested in finite-volume hyperbolic manifolds. The general strategy when proving something about finite-volume hyperbolic manifolds is the following: we first suppose that the manifold is compact for simplicity, and then we adapt the proof to the non-compact case by looking at what happens to the cusps.

\subsection{Geodesic spectrum} \label{geodesic:spectrum:subsection}
We start by proving the following.
\begin{prop} Let $M$ be a finite-volume complete hyperbolic manifold. For every $L>0$ there are finitely many closed geodesics in $M$ shorter than $L$.
\end{prop}
\begin{proof}
Suppose that there are infinitely many closed geodesics shorter than $L$. We know that $M$ decomposes into a compact part and a finite union of truncated cusps. The compact part has finite diameter $D$. 

Every closed geodesic intersects the compact part because a cusp contains no closed geodesic. Therefore we can fix a basepoint $x_0\in M$ and connect $x_0$ to these infinitely many geodesics of length $<L$ with arcs shorter than $D$. We use these arcs to freely homotope the geodesics into loops based at $x_0$ of length $<L+2D$, and lift the loops to arcs in $\matH^n$ starting from some basepoint $\tilde{x}_0\in \matH^n$.

If two such arcs end at the same point, the corresponding closed geodesics in $M$ are freely homotopic: this is excluded by Proposition \ref{unique:closed:geodesic:prop}, hence these endpoints are all distinct. The orbit of $\tilde{x}_0$ now contains infinitely many points in the ball $B(\tilde{x}_0, L+2D)$, a contradiction because the orbit is discrete.
\end{proof}

The lengths of the closed geodesics in $M$ form a discrete subset of $\matR$ called the \emph{geodesic spectrum} of $M$. Let $\ell_1(M)>0$ be the minimum of the spectrum of $M$: a \emph{shortest closed geodesic} in $M$ is a closed geodesic of shortest length $\ell_1(M)$. It is not necessarily unique (not even if considered up to orientation reversal), but there are finitely many of them.\index{shortest geodesic}

\begin{prop}
If $M$ is a closed hyperbolic manifold, then 
$$\inj M = \frac 12 \cdot \ell_1(M)$$
and every shortest closed geodesic is simple.
\end{prop}
\begin{proof}
The length of a closed geodesic is the minimum displacement of the corresponding hyperbolic transformation, so Corollary \ref{inj:cor} gives the equality.

Let now $\gamma$ be a shortest closed geodesic. Consider Proposition \ref{closed:geodesic:prop}:
if $\gamma$ wraps multiple times along a curve $\eta$, then $\eta$ is shorter than $\gamma$, a contradiction. If $\gamma$ self-intersects transversely at some point $p$, we can split $\gamma$ naturally as $\gamma_1*\gamma_2$ in $\pi_1(M,p)$, where both $\gamma_1$ and $\gamma_2$ are shorter than $\gamma$. Either $\gamma_1$ or $\gamma_2$ is non-trivial in $\pi_1(M,p)$ and is hence hyperbolic, but it has length smaller than $\ell_1(M)$: this contradicts Corollary \ref{minimum:cor}.
\end{proof}

\begin{rem}
A shortest geodesic may be non-simple when $M$ has cusps! By doubling an ideal triangle along its boundary we construct a hyperbolic surface called the \emph{thrice-punctured sphere} which has three cusps and contains various closed geodesics, none of which is simple. See Chapter \ref{surfaces:chapter}.
\end{rem}

\subsection{Isometry group}
We now study the isometry group $\Iso(M)$ of a hyperbolic manifold $M$.
Recall that the \emph{normaliser} $N(H)$ of a subgroup $H < G$ is the set of elements $g\in G$ such that $gH=Hg$. It is the biggest subgroup of $G$ containing $H$ as a normal subgroup.
The isometry group $\Iso(M)$ has an algebraic representation.\index{normaliser}

\begin{prop}
Let $M=\matH^n/_\Gamma$ be a complete hyperbolic manifold. There is a natural isomorphism 
$$\Iso(M) \Isom N(\Gamma)/_\Gamma.$$
\end{prop}
\begin{proof}
Every isometry $\varphi\colon M \to M$ lifts to an isometry $\tilde \varphi$ 
$$
\xymatrix{ 
\matH^n\ar@{.>}[r]^{\tilde \varphi} \ar[d]_\pi & \matH^n \ar[d]^\pi \\
M\ar[r]_\varphi & M
}
$$
such that $\tilde\varphi \Gamma = \Gamma \tilde\varphi$: hence $\tilde\varphi \in N(\Gamma)$. The lift is uniquely determined up to left- or right-multiplication by elements in $\Gamma$, hence we get a homomorphism
$$\Iso(M) \to N(\Gamma)/_\Gamma$$
which is clearly surjective (every element in $N(\Gamma)$ determines an isometry) and injective (if $\tilde\varphi\in\Gamma$ then $\varphi=\id$).
\end{proof}

Recall that the \emph{centraliser} of $H< G$ is the set of elements $g \in G$ such that $gh=hg$ for all $h$. It is a subgroup of $G$.\index{centraliser}

\begin{ex} \label{centralizer:ex}
Let $M=\matH^n/_\Gamma$ be a finite-volume hyperbolic manifold. The centralizer of $\Gamma$ in $\Iso(\matH^n)$ is trivial.
\end{ex}

\subsection{Outer automorphism group}
The \emph{automorphism group} of a group $G$ is the group $\Aut(G)$ of all the isomorphisms $G\to G$. The \emph{inner automorphisms} are those isomorphisms of type $g\mapsto hgh^{-1}$ for some $h\in G$, and they form a normal subgroup $\Int(G)\triangleleft \Aut(G)$. The quotient
$$\Out(G) = \Aut(G)/_{\Int(G)}$$
is called the \emph{outer automorphism group} of $G$.\index{group!automorphism group}\index{group!outer automorphism group}

If $x_0, x_1$ are two points in a path-connected topological space $X$ there is a non-canonical isomorphism
$\pi_1(X,x_0) \to \pi_1(X,x_1)$, unique only up to post-composing with an inner automorphism. Therefore there is a \emph{canonical} isomorphism $\Out(\pi_1(X,x_0))\to\Out(\pi_1(X,x_1))$. Hence $\Out(\pi_1(X))$ depends very mildly on the basepoint.

The group $\Omeo(X)$ of all homeomorphisms of $X$ does not act directly on $\pi_1(X)$ because of the inner-automorphism ambiguity, but we get a natural homomorphism 
$$\Omeo (X) \longrightarrow \Out(\pi_1(X))$$
which is neither injective nor surjective in general. We note that homotopic self-homeomorphisms give rise to the same element in $\Out(\pi_1(X))$.

\subsection{Finite isometry groups}
We turn back to hyperbolic manifolds. 

\begin{prop} \label{injective:prop}
If $M$ is a finite-volume complete hyperbolic manifold, the natural map 
$$\Iso(M) \longrightarrow \Out(\pi_1(M))$$
is injective.
\end{prop}
\begin{proof}
Set $M=\matH^n/_\Gamma$, identify $\Gamma$ with $\pi_1(M)$ and $\Iso(M)$ with $N(\Gamma)/_\Gamma$. With these identifications the map
$$N(\Gamma)/_\Gamma \longrightarrow \Out(\Gamma)$$
is just the conjugacy action that sends $h\in N(\Gamma)$ to the automorphism $g \mapsto h^{-1}gh$ of $\Gamma$. This is an inner automorphism if and only if there is a $f \in\Gamma$ such that $h^{-1}gh = f^{-1}gf$ for all $g\in\Gamma$, that is if $hf^{-1}$ commutes with $g$ for all $g\in\Gamma$. Exercise \ref{centralizer:ex} gives $h=f\in\Gamma$ and hence the map is injective.
\end{proof}

\begin{cor} \label{isometries:not:homotopic:cor}
Let $M$ be a finite-volume complete hyperbolic manifold. Distinct isometries of $M$ are not homotopic.
\end{cor}

This is a quite strong fact, from which we deduce the following.

\begin{cor} \label{finite:isometry:group:cor}
The isometry group of every finite-volume complete hyperbolic manifold $M$ is finite. 
\end{cor}
\begin{proof}
We know that $M$ is diffeomorphic to the interior of some compact $N$ having $k\geqslant 0$ boundary components. Let $C_1,\ldots, C_k \subset M$ be disjoint truncated cusps, all of the same small volume $V>0$. Every isometry of $M$ permutes these truncated cusps and fixes their complement $C = M\setminus \interior {C_1\cup \ldots \cup C_k}$.

The set $C$ is compact and preserved by the Lie group $\Iso(M)$, which is hence compact. To show that it is finite it suffices to prove that it is discrete. Suppose that a sequence of isometries $\varphi_i$ converges to the identity: we will deduce that $\varphi_i$ is homotopic to $\id$ for sufficiently large values of $i$, a contradiction.

For any $\varepsilon>0$ there is a $\varphi_i$ that moves all points of $C$ at distance $<\varepsilon$. 
If $\varepsilon>0$ is sufficiently small, the isometry $\varphi_i$ preserves each $C_i$, 
it restricts to an isometry of the flat torus $X_i = \partial C_i$ and its extension to $C_i$ is just a smaller rescaling of $\varphi_i|_{X_i}$ at every flat torus leaf in the cusp. 

Let also $\varepsilon>0$ be sufficiently small, so that $\varepsilon <\inj_x M$ for every $x\in C$. We get $d(x,\varphi_i(x))<\inj_xM$ for all $x\in C$, and since $\varphi_i$ acts on each $C_i$ just by rescaling its action on $X_i$, we get the same inequality for every $x\in M$.

For every $x\in M$ there is a unique geodesic $\gamma_x$ of length $d(x,\varphi_i(x))$ connecting $x$ and $\varphi_i(x)$. The geodesics $\gamma_x$ vary continuously with $x$ and we can use them to define a homotopy between $\varphi_i$ and $\id$: a contradiction.
\end{proof}

The finiteness of the isometry groups is peculiar to the hyperbolic world: we will see in Remark \ref{torus:isometries:oss} that the isometry group of a flat torus is infinite.

\subsection{Finite covers} \label{finite:covers:subsection}
Let $M=\matH^n/_\Gamma$ be a complete finite-volume hyperbolic manifold. How many finite covers $\tilde M \to M$ are there above $M$? Quite a lot, thanks to the following.

\begin{prop} \label{RF:hyperbolic:prop}
The fundamental group $\pi_1(M)$ of a complete finite-volume hyperbolic manifold $M$ is residually finite.
\end{prop}
\begin{proof}
The group $\pi_1(M)$ is finitely generated (actually, finitely presented) because $M$ is homeomorphic to the interior of a compact manifold with boundary, so Proposition \ref{FG:RF:prop} applies.
\end{proof}

This algebraic fact has some nice geometric consequences. Recall that $\ell_1(M)$ is the length of the shortest closed geodesic in $M$.

\begin{cor}
For every $L>0$ there is a finite cover $\tilde M \to M$ with $\ell_1(\tilde M)>L$.
\end{cor}
\begin{proof}
There are finitely many hyperbolic elements $a_1,\ldots, a_k \in \pi_1(M)$ such that every closed geodesic in $M$ of length $<L$ is freely homotopic to some of these.
By residually finiteness there is a finite-index normal $H \triangleleft \pi_1(M)$ that does not contain any of the elements $a_1,\ldots, a_k$, and hence none of their conjugates. 

The finite-index $H$ determines a finite cover $\pi\colon \tilde M \to M$ with $\pi_*(\pi_1(\tilde M)) = H$. A closed geodesic in $\tilde M$ cannot be shorter than $L$ because its image in $M$ would be freely homotopic to one of $a_1,\ldots, a_k$.
\end{proof}

\begin{cor}
If $M$ is a closed hyperbolic manifold, for every $R>0$ there is a finite cover $\tilde M \to M$ with $\inj\tilde M > R$.
\end{cor}

We can summarise this by saying that every closed hyperbolic manifold has arbitrarily fat finite covers.

\subsection{Subgroup separability}
The geometric consequences of residually finiteness are quite remarkable: we now introduce a stronger algebraic notion that will lead to more geometric applications.

Let $G$ be a group. A subgroup $H<G$ is \emph{separable} if the intersection of all the finite-index subgroups of $G$ containing $H$ is $H$ itself. In other words, for every non-trivial element $a\in G\setminus H$, there is a finite-index subgroup $G'< G$ that contains $H$ but not $a$ and hence ``separates'' $H$ from $a$.\index{separable subgroup} 

By definition, the trivial subgroup $\{e\}$ in $G$ is separable if and only if $G$ is residually finite. Note that, for general $H<G$, we cannot require the subgroup $G'$ to be normal in $G$, because $H$ itself may not be normal in $G$, as opposite to $\{e\}$.
The following proposition furnishes some interesting examples.

\begin{prop} \label{Long:prop}
Every maximal abelian subgroup $H$ in a residually finite group $G$ is separable.
\end{prop}
\begin{proof}
More generally, we suppose that $H$ is maximal with respect to some word relation $f(h_1,\ldots, h_n)=e$. (In our case, $f(h_1, h_2) = [h_1, h_2]$.) Since $G$ is residually finite, there is a sequence $N_i \triangleleft G $ of finite-index normal subgroups with $\cap N_i = \{e\}$. 

We show that $H = \cap_i HN_i$, and we conclude since $HN_i > N_i$ has finite index in $G$. We have $H \subset \cap_i HN_i$. We now show that the elements of $\cap_i HN_i$ satisfy $f=e$, and by maximality of $H$ we get $H=\cap_i HN_i$.

For every $N_i$ and $h_1,\ldots, h_n \in H$, we have $f(h_1N_i, \ldots, h_nN_i) \subset  N_i$ since its projection in $G/_{N_i}$ is trivial. Therefore $\cap_i f(HN_i, \ldots, HN_i) = \{e\}$ and the elements of $\cap_i HN_i$ satisfy $f=e$.
\end{proof}

\begin{cor} \label{separable:cor}
Let $M$ be a complete finite-volume hyperbolic $n$-manifold. The following subgroups are separable:
\begin{itemize}
\item the subgroup $\langle \varphi \rangle = \matZ$ generated by a primitive hyperbolic transformation $\varphi$;
\item the subgroup $\pi_1(T)$ generated by a $(n-1)$-torus cusp section $T$.
\end{itemize}
\end{cor}
\begin{proof}
Both subgroups are maximal abelian, see Section \ref{that:commute:subsection}.
\end{proof}

We deduce a couple of geometric consequences from this algebraic fact.

\begin{cor}
Let $M$ be a complete finite-volume hyperbolic $n$-manifold. Every closed geodesic $\gamma$ in $M$ lifts to a closed geodesic in some finite covers of $M$ of arbitrarily large degree.
\end{cor}
\begin{proof}
We have $M= \matH^n/_\Gamma$. The closed geodesic $\gamma$ is obtained by projecting the axis $l$ of a primitive hyperbolic element $\varphi \in \Gamma$.
Since $\langle \varphi \rangle$ is separable, there are subgroups $H<\pi_1(M)$ of arbitrarily large index that contain $\langle \varphi \rangle$.
\end{proof}

We have discovered in particular that there are many closed hyperbolic manifolds of bounded injectivity radius with arbitrarily large volume.
We can also use separability to promote some primitive closed geodesics to simple ones on some finite coverings.

\begin{cor}
Let $M$ be a complete finite-volume hyperbolic $n$-manifold. Every primitive closed geodesic $\gamma$ in $M$ lifts to a simple closed geodesic in some finite cover of $M$.
\end{cor}
\begin{proof}
We have $M= \matH^n/_\Gamma$. The closed geodesic $\gamma$ is obtained by projecting the axis $l$ of a primitive hyperbolic element $\varphi \in \Gamma$. If $\gamma$ is not simple, it self-intersects in $k>0$ points, which lift to $k$ transverse intersections between $l$ and some translates $\varphi_1(l), \ldots, \varphi_k(l)$, for some $\varphi_i \in \Gamma$. Therefore the elements $\psi \in \Gamma$ such that $\psi(l)$ intersects $l$ transversely are those of the form $\psi = \varphi^h\varphi_i\varphi^l$ for some $i=1,\ldots, k$ and $h,l \in \matZ$.

Since $\langle \varphi \rangle$ is separable, there is a finite-index $H<\Gamma$ with $\varphi_1, \ldots, \varphi_k \not\in H$ and $\langle \varphi \rangle \subset H$. No $\varphi^h\varphi_i\varphi^l$ lies in $H$, hence for every $\psi\in H$ the lines $l$ and $\psi(l)$ coincide or are disjoint. The line $l$ projects to a simple geodesic in $\matH^n/_H$.
\end{proof}

Concerning cusps, we deduce analogously the following.

\begin{cor}
Let $M$ be a complete finite-volume hyperbolic $n$-manifold and $T\subset M$ be a $(n-1)$-torus cusp section. The cusp section $T$ lifts to a cusp section in some finite covers of $M$ with arbitrarily large degree.
\end{cor}
\begin{proof}
The subgroup $\pi_1(T)$ is separable, so there are subgroups $H<\pi_1(M)$ of arbitrarily large finite degree containing it.
\end{proof}

\section{The Bieberbach Theorem} \label{Bieberbach:section}
The Margulis Lemma for Lie groups is fairly general and has important applications also in the elliptic and flat geometries. The most important one is the Bieberbach Theorem.

\subsection{Elliptic manifolds}
Recall that every complete elliptic manifold is isometric to $S^n/_\Gamma$ for some finite group $\Gamma < \On(n) = \Iso(S^{n-1})$ acting freely. We now improve the Margulis Lemma \ref{G:lemma} for $G=\On(n)$ by promoting ``nilpotent'' to ``abelian''. To this purpose we need the following.

\begin{lemma} \label{backwards:lemma}
There is a neighbourhood $U$ of $I\in \On(n)$ such that for every $A\in \On(n)$ and every $B\in U$ we have
$$[A,[A,B]] = I\ \Longrightarrow\ [A,B] = I.$$
\end{lemma}
\begin{proof}
We prove the lemma for the bigger Lie group $U(n)$ consisting of all unitary matrices. An immediate computation shows that
$$[A,[A,B]]=I \ \Longleftrightarrow\ [A, BA^{-1}B^{-1}] = I.$$
Let $\matC^n = \oplus V_i$ be the decomposition into orthogonal eigenspaces for $A$ and $A^{-1}$ (with respect to the standard hermitian product of $\matC^n$). The decomposition for the conjugate $BA^{-1}B^{-1}$ is just $\matC^n = \oplus B(V_i)$. 

Let $U$ be the neighbourhood of $I$ consisting of all matrices that move every vector of an angle $< \frac \pi 2$. Pick $B\in U$. Since the $V_i$ are orthogonal we get $B(V_i) \cap V_j = \{0\}$ for all $i\neq j$. On the other hand, the endomorphisms $A$ and $BA^{-1}B^{-1}$ commute and hence have a basis of common eigenvectors: hence the only possibility is that $B(V_i) = V_i$ for all $i$. 

The restriction $A|_{V_i}$ is just $\lambda_iI$ and hence commutes with $B|_{V_i}$, for all $i$. Therefore $A$ and $B$ commute everywhere.
\end{proof}

\begin{cor} \label{U:O(n):cor}
There is a neighbourhood $U$ of $I\in \On(n)$ such that every finite subgroup $\Gamma < \On(n)$ generated by some elements in $U$ is abelian.
\end{cor}
\begin{proof}
We know from Margulis Lemma \ref{G:lemma} that there is a $U$ where every such $\Gamma$ is nilpotent. The previous lemma promotes $\Gamma$ to an abelian group, because $[A,[A, \ldots [A,B]\cdots]] = I$ implies after finitely many steps that $[A,B] = I$ for every generators $A,B\in \Gamma \cap U$.
\end{proof}

\begin{prop} \label{U:finite:index:prop}
Let $G$ be a compact Lie group and $U$ a neighbourhood of $e\in G$. There is a $N>0$ such that for every group $\Gamma < G$, the subgroup $\Gamma_U< \Gamma$ generated by $\Gamma \cap U$ has index at most $N$ in $\Gamma$.
\end{prop}
\begin{proof}
Let $W\subset U$ be a smaller neighbourhood such that $W^{-1}=W$ and $W^2\subset U$ and set $N = \Vol (G) / \Vol (W)$ using the Haar measure for $G$.

We conclude by showing that if $g\Gamma_U$ and $g'\Gamma_U$ are distinct cosets of $\Gamma_U$ in $\Gamma$, then $gW \cap g'W = \emptyset$. This implies that $\Gamma_U$ has index at most $N$ in $\Gamma$.

Indeed, if $gW \cap g'W \neq \emptyset$ there are $w,w'\in W$ such that $gw = g'w'$ which implies that $g^{-1}g' = w(w')^{-1}\in W^2\subset U$ and hence $g^{-1}g' \in \Gamma_U$.
\end{proof}

We obtain a fairly interesting corollary about finite subgroups of $\On(n)$.

\begin{cor} \label{n:N:cor}
For every $n$ there is a $N>0$ such that every finite subgroup of $\On(n)$ contains an abelian subgroup of index at most $N$.
\end{cor}

\subsection{Isometries of Euclidean space}
We now turn to flat manifolds. Recall that a complete flat manifold is isometric to $\matR^n/_\Gamma$ for some discrete group $\Gamma < \Iso(\matR^n)$ acting freely. Every isometry $g$ of $\matR^n$ can be written uniquely as
$$g\colon x \longmapsto Ax +b$$
for some $A\in \On(n)$ and $b\in \matR^n$. The \emph{rotational} and \emph{translational part} of $g$ are $A$ and $b$ respectively. Let $\Fix(A)$ be the fixed points of $x\mapsto Ax$.

\begin{prop} \label{basic:isometries:Rn:prop}
Let $Ax+b$ be an isometry of $\matR^n$. Then:
\begin{itemize}
\item its inverse is $A^{-1}x-A^{-1}b$,
\item if it acts freely then $\Fix(A)\neq \{0\}$,
\item $\exists$ a translation of $\matR^n$ that conjugates it to $Ax+b'$ with $b'\in \Fix(A)$.
\end{itemize}
\end{prop}
\begin{proof}
If it acts freely then $Ax+b=x$ has no solution $x\in\matR^n$. By rewriting the equation as $(A-I)x = -b$ we see that $A-I$ is not surjective, hence not injective, hence $\Fix(A)$ is non-trivial.

A translation $x+d$ conjugates $Ax+b$ into 
$$\big(A(x+d)+b\big)-d = Ax + (A-I)d + b.$$
Since $A$ is orthogonal, we get $\Img (A-I) = \ker (A-I)^\perp = \Fix(A)^\perp$ and hence there is a $d$ such that $b'=(A-I)d+b\in \Fix(A)$.
\end{proof}

\begin{ex} \label{Euclidean:commutators:ex}
The commutator of two isometries is
$$[Ax+b, Cx+d] = [A,C]x + A(I-C)A^{-1}b + AC(I-A^{-1})C^{-1}d.$$
In particular we get
$$[Ax+b,x+d] = x + (A-I)d.$$
\end{ex}

\subsection{Discrete groups}
The homomorphism $r\colon \Iso(\matR^n)\to \On(n)$ that sends every isometry to its rotational part induces an exact sequence
$$0 \longrightarrow \matR^n \longrightarrow \Iso(\matR^n) \stackrel{r}\longrightarrow \On(n) \longrightarrow 0$$
where we indicate by $\matR^n$ the group of translations of $\matR^n$. 

Let $\Gamma < \Iso(\matR^n)$ be a discrete group. We get an exact sequence
$$0 \longrightarrow H \longrightarrow \Gamma \longrightarrow r(\Gamma) \longrightarrow 0$$
where $H \triangleleft \Gamma$ is the \emph{translation subgroup} of $\Gamma$. The subgroup $r(\Gamma)<\On(n)$ is not necessarily discrete, as the following example shows.

\begin{example} \label{rototranslation:example}
A \emph{rototranslation} in $\matR^3$ is a rotation of some angle $\theta$ along an axis $r$ composed with a translation of some distance $t>0$ in the direction of $r$. A rototranslation generates a discrete group $\Gamma = \matZ$ acting freely on $\matR^3$, whose quotient $\matR^3/_\Gamma$ is diffeomorphic to $\matR^2\times S^1$. If $\theta$ is not commensurable with $\pi$ the image $r(\Gamma)$ is not discrete and forms a dense subset of the circle in $\On(3)$ of all rotations along $r$.
\end{example}

We now extend Corollary \ref{n:N:cor} to this context.

\begin{teo} \label{discrete:virtual:abelian:teo}
For every $n$ there is a $N>0$ such that every discrete subgroup of $\Iso(\matR^n)$ contains an abelian subgroup of index at most $N$.
\end{teo}
\begin{proof}
We pick a small neighbourhood $U \subset \On(n)$ of $I$, such that
the commutator is contracting in $U$ (see the proof of the Margulis Lemma \ref{G:lemma}), it is symmetric ($U=U^{-1}$) and $A-I$ contracts vectors uniformly for all $A\in U$, more precisely $\|(A-I)v\| < \frac 14\|v\|$ for all $v\neq 0$.

Let $\Gamma < \Iso(\matR^n)$ be discrete and $r(\Gamma)_U$ be the group generated by $r(\Gamma)\cap U$. By Proposition \ref{U:finite:index:prop} the index of $r(\Gamma)_U$ in $r(\Gamma)$ is bounded by some $N$ depending only on $U$. Therefore its counterimage in $\Gamma$
$$\Gamma^* = r^{-1}(r(\Gamma)_U) \cap \Gamma$$ 
has index in $\Gamma$ bounded by $N$. The group $\Gamma^*$ is the subgroup of $\Gamma$ consisting of all $Ax+b$ with $A\in r(\Gamma)_U$. 
It remains to prove that $\Gamma^*$ is abelian.

Let $A_1x+b_1$ and $A_2x +b_2$ be two elements in $\Gamma^*$. Define
\begin{equation*} 
A_{i+1}x+b_{i+1} = [A_1x+b_1, A_ix+b_i]
\end{equation*}
for all $i\geqslant 2$. Exercise \ref{Euclidean:commutators:ex} gives
\begin{equation} \label{commutator:eqn}
A_{i+1}x+b_{i+1} = [A_1,A_i]x + A_1(I-A_i)A_1^{-1}b_1 + A_1A_i(I-A_1^{-1})A_i^{-1}b_i
\end{equation}
and hence
\begin{align*}
A_{i+1} & = [A_1,A_i], \\
 b_{i+1} & = A_1(I-A_i)A_1^{-1}b_1 + A_1A_i(I-A_1^{-1})A_i^{-1}b_i.
\end{align*}
Suppose $A_1, A_2 \in U$. Since the commutator is contracting in $U$ and $A_i-I$ contract vectors uniformly, we get $A_{i} \to I$ and $b_i\to 0$ as $i\to \infty$. Therefore $A_ix+b_i$ tends to the identity, and since $\Gamma^*$ is discrete it \emph{is} the identity for all $i$ bigger or equal than some $i_0$. In particular $A_{i_0} = I$ and Lemma \ref{backwards:lemma} used backwards gives $A_3=[A_1,A_2]=I$: hence $A_1$ and $A_2$ commute. 

Since $r(\Gamma^*)$ is generated by $r(\Gamma^*)\cap U$,
we deduce that $r(\Gamma^*)$ is abelian. Therefore (\ref{commutator:eqn}) may be restated as
\begin{equation} \label{simplified:commutator:eqn}
A_{i+1}x+b_{i+1} = x + (I-A_i)b_1 + (A_1 - I)b_i.
\end{equation}
We now consider the case $A_1 \in U$ and $A_2=I$. We get
$$A_{i+1}x+b_{i+1} = x + (A_1-I)^{i-1} b_2.$$
Since $A_1\in U$ we get $(A_1-I)^i b_2 \to 0$ and hence $(A_1-I)^ib_2=0$ for some $i$, which gives $(A_1-I)b_2=0$ since $A_1$ is diagonalisable. Thus $b_2 \in \Fix (A_1)$.

We have proved that if $\Gamma^*$ contains a translation $x+b$, then $b\in \Fix(A)$ for all $A\in r(\Gamma^*)\cap U$. Since these $A$ generate $r(\Gamma^*)$, the vector $b$ belongs to 
$$W = \Fix(r(\Gamma^*)) = \big\{x\ \big|\ Ax = x\ \forall A \in r(\Gamma^*)\big\}.$$
Pick now two arbitrary elements $Ax+b$ and $Cx+d$ in $\Gamma^*$.
By (\ref{simplified:commutator:eqn}) we get
$$[Ax+b, Cx+d] = x+(I-C)b+(A-I)d.$$
By what just said $(I-C)b+(A-I)d \in W$. On the other hand 
$$\Img (I-C) = \ker(I-C)^\perp = \Fix(C)^\perp \subset W^\perp$$ 
and hence $(I-C)b\in W^\perp$, and analogously $(A-I)d \in W^\perp$. We deduce that $(I-C)b+(A-I)d \in W \cap W^\perp$ is trivial and all elements in $\Gamma^*$ commute.
\end{proof}

\subsection{Crystallographic groups} \label{crystallographic:subsection}
A \emph{crystallographic group} is a discrete subgroup $\Gamma <\Iso(\matR^n)$ with compact quotient $\matR^n/_\Gamma$.\index{group!crystallographic group}

\begin{prop} \label{finite:image:r:prop}
The image $r(\Gamma)$ of a crystallographic group $\Gamma<\Iso(\matR^n)$ is finite. 
\end{prop}
\begin{proof}
By Theorem \ref{discrete:virtual:abelian:teo} we may suppose that $\Gamma$ is abelian. We now prove that $\Gamma$ abelian implies that $r(\Gamma)$ is trivial, \emph{i.e.}~all elements in $\Gamma$ are translations.

Suppose that $\Gamma$ contains a non-translation $Ax +b$. We conjugate $\Gamma$ by a translation as in Proposition \ref{basic:isometries:Rn:prop} to get $b\in \Fix(A)$. Pick another isometry $Cx+d$ in $\Gamma$. The commutator
$$[Ax +b, Cx +d] = x+(A-I)d-(C-I)b$$
is trivial. Since $A$ and $C$ commute, we get $(C-I)b \in \Fix(A)$, hence $(A-I)
d\in \Fix (A)$ and finally $d\in \Fix(A)$.

We have proved that $d\in \Fix(A)\subsetneq \matR^n$ for all elements $Cx+d$ in $\Gamma$. Therefore the $\Gamma$-orbit of $0\in\matR^n$ is contained in $\Fix(A)$. But the compactness of $\matR^n/_\Gamma$ implies that there is a compact fundamental domain, and hence a $R>0$ such that every point in $\matR^n$ is $R$-close to any fixed orbit: a contradiction since $\Fix(A)$ is a proper vector subspace of $\matR^n$.
\end{proof}

\begin{cor} \label{crystal:cor}
Every crystallographic group has a finite-index translation subgroup isomorphic to $\matZ^n$.
\end{cor}

Recall from Section \ref{piatte:subsection} that a \emph{flat torus} is a $n$-torus $\matR^n/_\Gamma$ where $\Gamma$ is a lattice, \emph{i.e.}~a discrete group isomorphic to $\matZ^n$ that spans $\matR^n$ as a vector space.\index{flat torus}\index{Bieberbach Theorem}

\begin{cor}[Bieberbach's Theorem] \label{Bieberbach:cor}
Every closed flat $n$-manifold is finitely covered by a flat torus.
\end{cor}

We conclude the discussion by noting that there are no ``cusps'' in Euclidean geometry.

\begin{prop} \label{finite:volume:flat:prop}
Every finite-volume complete flat manifold is closed.
\end{prop}
\begin{proof}
Let $M=\matR^n/_\Gamma$ be a finite-volume flat manifold. Up to taking finite indexes, the group $\Gamma$ is abelian by Theorem \ref{discrete:virtual:abelian:teo}. If all the elements are translations, we are done: finite-volume easily implies that $\Gamma = \matZ^n$, thus $M$ is closed. 

Suppose some element $Ax + b$ is not a translation: every element of $\Gamma$ commutes with it and hence preserves the space $\Fix (A)$ that has dimension $k<n$; hence $M$ is isometric to $\Fix (A)/_\Gamma \times \matR^{n-k}$ and has infinite volume.
\end{proof}

\subsection{References}
Most of the material introduced in this chapter is standard and can be found in Benedetti -- Petronio \cite{BP}, Ratcliffe \cite{R}, and of course in Thurston's notes \cite{Th}. We have also consulted Thurston's book \cite{Th_book} for the part on chrystallographic groups. Proposition \ref{Long:prop} was taken from a paper of Long \cite{L}.

%% file: Infinity.tex
\chapter{The sphere at infinity} \label{infinity:chapter}
We have discovered that every complete hyperbolic manifold is a quotient $M=\matH^n/_\Gamma$ for some discrete group $\Gamma < \Iso(\matH^n)$ acting freely, and we now raise the following question: how does $\Gamma$ act on the boundary $\partial \matH^n$ at infinity? Does the action of $\Gamma$ on $\partial \matH^n$ furnish some information on the geometry of $M$?

We show in this chapter that $\partial \matH^n$ subdivides naturally into two $\Gamma$-invariant subsets: an open set $\Omega(\Gamma)$ called the \emph{domain of discontinuity} where the action of $\Gamma$ is properly discontinuous (like in $\matH^n$) and a complementary closed set $\Lambda(\Gamma)$ called the \emph{limit set} where the action of $\Gamma$ is more chaotic.

We then use the limit set to define a canonical decomposition of cusped hyperbolic manifolds into ideal polyhedra, called the \emph{Epstein--Penner decomposition}.

We also devote some time to prove Theorem \ref{estensione:teo}, which states that every smooth homotopy equivalence $M\to N$ between two closed hyperbolic manifolds lifts to a map $\matH^n \to \matH^n$ that extends nicely to the compactifications $\overline{\matH^n}\to\overline{\matH^n}$. This fact will have important applications in the study of hyperbolic manifolds in the subsequent chapters.

\section{Limit set}
How does a discrete group $\Gamma $ of isometries of $\matH^n$ act on the boundary at infinity $\partial \matH^n$? We now prove that $\partial \matH^n$ divides canonically into two $\Gamma$-invariant subsets: an open zone $\Omega(\Gamma)$ where $\Gamma$ acts properly discontinuously, and a closed one $\Lambda(\Gamma)$ where it does not. These regions are called respectively the \emph{domain of discontinuity} and the \emph{limit set} of $\Gamma$. 

\subsection{The limit set}
Throughout this section, we let $\Gamma$ be a non-trivial discrete group of isometries of $\matH^n$. Fix a point $x\in\matH^n$. We now that the orbit $\Gamma(x)$ is a discrete subset of $\matH^n$. The \emph{limit set} $\Lambda (\Gamma)\subset \partial \matH^n$ of $\Gamma$ is the set of all the accumulation points of the orbit $\Gamma(x)$ in $\overline{\matH^n}$.\index{limit set}\index{$\Lambda(\Gamma)$}

\begin{ex}
The limit set does not depend on $x$.
\end{ex}

The limit set is clearly a closed $\Gamma$-invariant subset of $\partial \matH^n$. If $\Gamma' < \Gamma$ we obviously get $\Lambda(\Gamma')\subset \Lambda(\Gamma)$.

\begin{ex} \label{index:lambda:ex}
If $\Gamma'$ has finite index in $\Gamma$ then $\Lambda(\Gamma') = \Lambda(\Gamma)$.
\end{ex}

Recall that $\Gamma$ is elementary if it fixes a finite set of points in $\overline{\matH^n}$.

\begin{example} If $\Gamma$ is elementary and acts freely, that is $M = \matH^n/_\Gamma$ is wither a cusp or a tube (see Proposition \ref{elementary:prop}), then $\Lambda(\Gamma) = \Fix(\Gamma)$ consists of one or two points. 
\end{example}

We now characterise the elementary subgroups of $\Iso(\matH^n)$. Recall that $\Gamma$ is \emph{virtually P} if it has a finite-index subgroup that is P (where P is some property).\index{group!elementary group}

\begin{prop} \label{elementary:4:prop}
The following are equivalent:
\begin{enumerate}
\item $\Gamma$ is elementary,
\item $\Gamma$ fixes either a point $x\in\matH^n$, or a line, or a point $x\in\partial\matH^n$ and all the horospheres centred at $x$,
\item $\Gamma$ is virtually abelian,
\item $\Lambda(\Gamma)$ consists of 0, 1, or 2 points.
\end{enumerate}
\end{prop}
\begin{proof}
(1) $\Rightarrow$ (2). We expand on the proof of Proposition \ref{elementary:prop} by taking elliptic elements into account.

By hypothesis $\Gamma$ fixes a finite set of points in $\overline{\matH^n}$. If some of them lie in $\matH^n$, then $\Gamma$ fixes their barycenter and we are done. If $\Gamma$ fixes more than two points in $\partial \matH^n$, a barycenter is also defined: their convex hull is an ideal polyhedron of dimension at least 2, we can truncate (via horospheres) all the vertices by the same small volume and take the barycenter of the vertices of the resulting compact combinatorial polyhedron. So we are done also in this case. If $\Gamma$ fixes two points, it fixes a line.

We are left with the case where $\Gamma$ fixes a point $x\in \partial \matH^n$, and no other finite set of points. The proof of Lemma \ref{fissi:lemma} extends as is when $\varphi_2$ is elliptic, and shows that if $\Gamma$ contains a hyperbolic element $\varphi_1$ then every other non-trivial element of $\Gamma$ is either hyperbolic or elliptic and fixes the same axis of $\varphi_1$, but this is excluded: so $\Gamma$ contains only parabolics and elliptics, and these must fix all the horospheres centred at $x$, as required.

(2) $\Rightarrow$ (3). The group $\Gamma$ is virtually isomorphic to a discrete subgroup of $\Iso(\matR^m)$ for some $m$ and we apply Theorem \ref{discrete:virtual:abelian:teo}.

(3) $\Rightarrow$ (4). By Exercise \ref{index:lambda:ex} we may suppose that $\Gamma$ is abelian. By adapting the proof of Lemma \ref{commute:lemma} to elliptics we see (exercise) that $\Gamma$ satisfies (2) and hence we get (4).

(4) $\Rightarrow (1)$. If $\Lambda(\Gamma)$ is empty then $\Gamma$ is finite and fixes a point by Proposition \ref{finite:fixed:prop}.
If it consists of one or two points, these are preserved by $\Gamma$.
\end{proof}

\subsection{Minimality}
Let $\Gamma$ be a non elementary non-trivial discrete group of isometries of $\matH^n$.
We are interested in the action of $\Gamma$ on $\Lambda(\Gamma)$. We want to prove that the action is \emph{minimal}, that is $\Lambda(\Gamma)$ has no invariant non-empty proper closed subset.

Recall from Section \ref{polyhedra:subsection} that every closed subset $S\subset \overline{\matH^n}$ has a well-defined convex hull $C(S)\subset \overline{\matH^n}$. The convex hull $C(S)$ is closed and is the intersection of all the closed half-spaces containing $S$. We have $C(S) \cap \partial \matH^n = S \cap \partial \matH^n$.

\begin{prop}
If $\Gamma$ is not elementary, its action on $\Lambda(\Gamma)$ is minimal.
\end{prop}
\begin{proof}
Let $S\subset \Lambda(\Gamma)$ be a non-empty closed $\Gamma$-invariant subset. Since $S$ is $\Gamma$-invariant, its convex hull $C(S)$ also is. Since $\Gamma$ is not elementary, $S$ contains at least two points and hence the intersection $C(S)\cap \matH^n$ is not empty. If $x\in C(S)\cap \matH^n$, its $\Gamma$-orbit is confined in $C(S)$ and hence also its accumulation points are: therefore $\Lambda(\Gamma)\subset C(S)\cap \partial \matH^n = S$ and hence $\Lambda(\Gamma) = S$.
\end{proof}

\begin{cor} \label{normal:limit:cor}
If $\Gamma$ is not elementary and $\Gamma' \triangleleft \Gamma$ is an infinite normal subgroup, then $\Lambda(\Gamma') = \Lambda (\Gamma)$.
\end{cor}
\begin{proof}
For every $\gamma \in \Gamma$ we have $\gamma^{-1}\Gamma' \gamma = \Gamma'$, so $\gamma$ sends the orbit $\Gamma'(x)$ to the orbit $\Gamma'(\gamma(x))$ and hence preserves $\Lambda(\Gamma')$. Therefore $\Lambda(\Gamma')$ is $\Gamma$-invariant.

The limit set $\Lambda(\Gamma')$ is not empty because $\Gamma'$ is infinite: the minimality of $\Lambda(\Gamma)$ implies that $\Lambda(\Gamma') = \Lambda(\Gamma)$.
\end{proof}

\subsection{The convex core}
If $C\subset \overline{\matH^n}$ is a closed convex set, we define the \emph{nearest point retraction}\index{nearest point retraction} 
$$r\colon \overline{\matH^n} \to C$$ 
as the map that sends $x$ to the point $r(x)\in C$ that is closer to $x$. If $x\in\partial \matH^n$, we interpret $r(x)$ as the first point of $C$ that is contained in some horosphere centred at $x$. We have $r(x)=x$ if and only if $x\in C$. The map $r$ is continuous. Using the geodesic from $x$ to $r(x)$ we can construct a natural deformation retraction of $\overline{\matH^n}$ onto the closed convex set $C$.

Let $\Gamma < \Iso(\matH^n)$ be a non elementary non-trivial discrete group. The convex hull of $\Lambda(\Gamma)$ is $\Gamma$-invariant, so it makes sense to define the following.\index{convex core} 

\begin{defn} The \emph{convex core} of the orbifold $O=\matH^n/_\Gamma$ is the quotient $C(\Lambda(\Gamma))/_\Gamma \subset O$ of the convex hull $C(\Lambda(\Gamma))$ of the limit set $\Lambda(\Gamma)$.
\end{defn}

We are of course mostly interested in the case where $\Gamma$ acts freely and hence $M=\matH^n/_\Gamma$ is a manifold.
The deformation retraction defined above is $\Gamma$-invariant, therefore every complete hyperbolic manifold $M$ deformation retracts onto its convex core. In particular, $M$ is homotopically equivalent to its convex core.

\subsection{The domain of discontinuity}
Let $\Gamma < \Iso(\matH^n)$ be a non elementary non-trivial discrete group. The \emph{domain of discontinuity} of $\Gamma$ is the open set\index{domain of discontinuity}\index{$\Omega(\Gamma)$}
$$\Omega(\Gamma) = \partial\matH^n \setminus \Lambda(\Gamma).$$
The following proposition explains the terminology.
\begin{prop}
The action of $\Gamma$ on $\matH^n \cup \Omega(\Gamma)$ is properly discontinuous.
\end{prop}
\begin{proof}
The nearest point retraction $r$ sends $\matH^n \cup \Omega(\Gamma)$ to $C(\Lambda(\Gamma)) \setminus \Lambda(\Gamma)$ and commutes with $\Gamma$. The action of $\Gamma$ on the latter set is properly discontinuous (since it is contained in $\matH^n$), so the action on the former also is.
\end{proof}

\begin{prop}
If $\Vol(\matH^n/_\Gamma)<+ \infty $ then $\Lambda(\Gamma) = \partial \matH^n$.
\end{prop}
\begin{proof}
If $\Omega(\Gamma)\neq \emptyset$, pick $x\in \Omega(\Gamma)$. The point $x$ has a neighbourhood system consisting of half-spaces. Since $\Gamma$ acts properly discontinuously on $\matH^n \cup \Omega(\Gamma)$, there is some half-space $H$ which intersects only finitely many $\Gamma$-translates, contradicting the finite volume hypothesis.
\end{proof}

In this book we are mostly interested in finite-volume complete hyperbolic manifolds or orbifolds $\matH^n/_\Gamma$, and for these the limit set is just the whole boundary, regardless of $\Gamma$. This apparently disappointing piece of information has some interesting algebraic consequences. 

\begin{cor} \label{normal:abelian:cor}
If $\Vol(\matH^n/_\Gamma) <+\infty$ then
$\Gamma$ does not contain any non-trivial virtually abelian normal subgroup.
\end{cor}
\begin{proof}
Let $H\triangleleft \Gamma$ be virtually abelian. If $H$ is infinite, then $\Lambda(H) = \Lambda (\Gamma)=\partial \matH^n$ by Corollary \ref{normal:limit:cor}, contradicting Proposition \ref{elementary:4:prop}. If $H$ is finite, then $\Fix(H)$ is a non-empty proper subspace of $\matH^n$ and $\Gamma$ acts on $\Fix(H)$ because $H$ is normal. If $\Fix(H)$ is a point then $\Gamma$ is elementary, otherwise we get $\Lambda(\Gamma) \subset \partial \Fix(H)$, a contradiction in both cases.
\end{proof}

\begin{cor}
The fundamental group of a finite-volume complete hyperbolic manifold is never solvable.
\end{cor}
\begin{proof}
Solvable groups have non-trivial normal abelian subgroups, as proved in Proposition \ref{solvable:normal:prop}.
\end{proof}

\subsection{Schottky groups}
In all the examples encountered up to now the limit set $\Lambda(\Gamma)$ consists of either few points or the whole of $\partial \matH^n$. There are many cases where $\Lambda(\Gamma)$ is a more interesting (and often beautiful) set. The simplest examples are probably the following.

Choose a number $k\geqslant 2$ and $2k$ half-spaces $H_1,\ldots,H_{2k} \subset \matH^n$ with disjoint closures in $\overline{\matH^n}$. 
The closed complement $C = \matH^n \setminus \interior{H_1 \cup \ldots \cup H_{2k}}$ is a convex subset bounded by $2k$ disjoint hyperplanes. If we pair isometrically these hyperplanes we get a complete (exercise) hyperbolic manifold $M$, which is hence $M=\matH^n/_\Gamma$ for some $\Gamma$. The hyperbolic manifold $M$ has infinite volume and is the interior of a compact manifold that decomposes into one 0-handle and $k$ 1-handles. The group $\Gamma$ is free with $k$ generators.\index{group!Schottky group}

\begin{ex} If $n=2$, the limit set $\Lambda(\Gamma)\subset \partial \matH^2 = S^1$ is a Cantor set.
\end{ex}

\subsection{Parabolic and hyperbolic points}
Let $\Gamma < \Iso(\matH^n)$ be a non elementary non-trivial discrete group. A \emph{parabolic point} $x\in \partial \matH^n$ for $\Gamma$ is a point that is fixed by some parabolic element $\gamma \in \Gamma$. Similarly, a \emph{hyperbolic point} is a point that is fixed by some hyperbolic element.\index{parabolic point}\index{hyperbolic point}

\begin{prop}
The set of all parabolic (hyperbolic) points is either empty or a dense subset of $\Lambda(\Gamma)$.
\end{prop}
\begin{proof}
They clearly form a $\Gamma$-invariant subset, so its closure also is. We conclude because $\Lambda(\Gamma)$ is minimal.
\end{proof}

\begin{cor} \label{parabolic:hyperbolic:dense:cor}
If $\Vol(\matH^n/_\Gamma)<+\infty $, the set of all parabolic (hyperbolic) points is either empty or dense in $\partial \matH^n$.
\end{cor}

\subsection{Horoballs in the hyperboloid model}
In the conformal models, horoballs are just Euclidean balls tangent to the boundary at infinity. In the hyperboloid model $I^n$ there is a more algebraic description. We denote by
$$L = \{x\in \matR^{n+1}\ \big|\ \langle x,x \rangle = 0, \ x_{n+1}>0\}$$ 
the positive light cone and recall from Section \ref{projective:compactification:subsection} that the boundary at infinity $\partial I^n$ may be interpreted as the set of rays in $L$. More than that, every vector $x\in L$ determines a horoball $O$ centred at $[x]\in\partial I^n$, via the equation
$$O = \big\{y\in I^n \ |\ -1 \leqslant \langle x,y \rangle < 0\big\}.$$
The horoball gets smaller as $x$ goes to infinity.
\begin{ex}
This is really a horoball centred at $[x]$.
\end{ex}
\begin{proof}[Hint]
Prove that the boundary horosphere $\langle x,y \rangle =-1$ is orthogonal to all the lines pointing to $x$. 
\end{proof} 

There is also a simple (and maybe surprising) geometric relation between the vector $x$ and the corresponding horoball in the Poincar\'e disc model $D^n$. Represent $D^n$ inside the hyperplane $x_{n+1} = 0$, and let $\pi \colon I^n \to D^n$ the isometry obtained by projecting towards $P=(0,\ldots,0,-1)$, considered in Section \ref{Poincare:subsection}. By projecting towards $P$ we also get a homeomorphism $\pi\colon L \to D^n$.

\begin{ex} \label{Euclidean-centre:ex}
The Euclidean centre of the horoball $\pi(O)$ is $\pi(x)$.
\end{ex}

\subsection{The Epstein--Penner decomposition} \label{EP:subsection}
We now show that every cusped finite-volume hyperbolic manifold decomposes canonically into some ideal polyhedra. This decomposition is known as the \emph{canonical} or \emph{Epstein--Penner decomposition}.\index{Epstein--Penner decomposition}

Let $M = \matH^n/_\Gamma$ be a non-compact finite-volume complete hyperbolic manifold. The manifold $M$ is diffeomorphic to the interior of a compact $N$ with $c\geqslant 1$ boundary components. We fix $c$ disjoint truncated cusps $H_1, \ldots, H_c \subset M$. Their lifts in $\matH^n$ are $c$ disjoint $\Gamma$-orbits $B_1^i, \ldots, B_c^i$ of disjoint horoballs.

\begin{figure}
\begin{center}
\includegraphics[width = 12 cm] {\iftoggle{BW}{horoballs-BW}{horoballs}} 
\nota{By glueing two identical ideal triangles along their sides we get a thrice-punctured sphere $S$, and we choose three disjoint truncated cusps in it having the same area (left). In the universal cover, the two ideal triangles of $S$ lift to the Farey tessellation already shown in Figure \ref{tiling_Farey:fig}, and each of the three truncated cusps lifts to infinitely many disjoint horoballs (right). The Euclidean centre of the horoballs is shown here, in relation to Exercise \ref{Euclidean-centre:ex}.}
\label{horoballs:fig}
\end{center}
\end{figure}

\begin{example}
Consider the thrice-punctured sphere $S$ obtained by gluing two identical copies of an ideal triangle along their sides. This is a complete finite-volume hyperbolic surface with three cusps. We may fix three disjoint truncated cusps of the same area, as sketched in Figure \ref{horoballs:fig}-(left). These lift to three families of disjoint horoballs as in Figure \ref{horoballs:fig}-(right).
\end{example}

We use the hyperboloid model $I^n$ and interpret a horoball as a point in the positive light cone $L$ as explained in the previous section. The $\Gamma$-orbits $B_1^i, \ldots, B_c^i$ of horoballs form a discrete subset in $L$. If we modify the initial truncated cusp $H_j$, all the points $B_j^i\in L$ are rescaled by the same constant.

We make a crucial observation: the points $B_j^i$ are discrete in $L$, but the rays that contain them form a countable dense subset of $L$, because they correspond to the parabolic points of $\Gamma$, that are dense by Corollary \ref{parabolic:hyperbolic:dense:cor}.

\begin{figure}
\begin{center}
\includegraphics[width = 13 cm] {\iftoggle{BW}{convex-hull-BW}{convex-hull}} 
\nota{A portion of the convex hull $C$ for the thrice-punctured sphere of Figure \ref{horoballs:fig}. As suggested by Exercise \ref{Euclidean-centre:ex} we determine the horoballs in $L$ by projecting the Euclidean centers of their representations in the disc model. The black half-line above each horoball in $L$ is contained in $C$, see Proposition \ref{Bji:prop}.}
\label{convex-hull:fig}
\end{center}
\end{figure}

We now define $C\subset \matR^{n+1}$ as the convex hull of the points $B_j^i$ in $\matR^{n+1}$. We note that $C$ is $\Gamma$-invariant. An example is shown in Figure \ref{convex-hull:fig}.
The convex hull will define a decomposition of $M$ into ideal polyhedra. 

The convex hull $C$ is contained in the convex hull of $L$, that consists of all the positive timelike and lightlike vectors. The following proposition says that $C$ intersects $L$ into countably many half-lines based at the points $B_j^i$.

\begin{prop} \label{Bji:prop}
The set $C\cap L$ is the set of points $\alpha B_j^i$ for some $\alpha \geqslant 1$.
\end{prop}
\begin{proof}
Pick a point $x\in L$ not of this type: the segment $0x$ does not intersect $\{B_j^i\}$. Since $\{B_j^i\}\subset L$ is discrete, we can perturb the hyperplane $H$ tangent in $x$ to $L$ so that it intersects $L$ into a small ellipsoid around $0x$, and so that $\{B_j^i\}$ lies on the opposite side of $H$ of $x$. Therefore $x \not \in C$.

Conversely, pick a point $x = \alpha B_{j_0}^{i_0}$. The rays containing the points $B_j^i$ form a dense set in $L$, so we can find a sequence of $B_j^i\neq B_{j_0}^{i_0}$ such that $[B_j^i]$ converges to $[B_{j_0}^{i_0}]$. Since $\{B^i_j\}$ is discrete, the points $B_j^i$ go to infinity in the sequence and hence the segment with endpoints $B_j^i$ and $B_{j_0}^{i_0}$ approaches $x$ in the limit. Therefore $x\in C$.
\end{proof}

\begin{prop}
Every timelike ray intersects $\partial C$ exactly once.
\end{prop}
\begin{proof}
Every timelike ray $r$ enters in the interior of a polyhedron spanned by some $B_j^i$, so it intersects $C$ and hence $\partial C$. Consider a point $x \in \partial C \cap r$. Every supporting hyperplane for $x$ must be spacelike because the rays containing the $B_j^i$ are dense. Therefore $C \cap r$ is a half-line and $x$ is its endpoint.
\end{proof}

In particular $\partial C = (\partial_{\rm l}C) \sqcup (\partial_{\rm t}C)$ consists of the lightlike vectors $\partial_{\rm l} C = C\cap L$ and some timelike vectors $\partial_{\rm t} C$. The previous proposition furnishes a natural 1-1 correspondence $\partial_{\rm t} C \longleftrightarrow I^n$ with the hyperboloid model by projecting along timelike lines.

\begin{prop}
The timelike boundary $\partial_{\rm t} C$ is tessellated by countably many Euclidean $n$-dimensional polyhedra with vertices in $\{B_j^i\}$.
\end{prop}
\begin{proof}
Consider a point $x\in\partial_{\rm t} C$ and a supporting hyperplane $H$. As already mentioned, the hyperplane $H$ is spacelike and the intersection $H\cap \partial_{\rm t} C$ is some $k$-dimensional polyhedron, convex hull of finitely many points in $\{B_j^i\}$. If $k<n$ we can rotate $H$ until it meets one more point $B_j^i$ and after finitely many rotations we get a $n$-dimensional polyhedron.

We have proved that $\partial_{\rm t} C$ is paved by $n$-dimensional polyhedra with vertices in $\{B_j^i\}$ and intersecting in common faces. The polyhedra are locally finite because $\{B_j^i\}$ is discrete, hence they form a tessellation.
\end{proof}

\begin{figure}
\begin{center}
\includegraphics[width = 13 cm] {\iftoggle{BW}{hyperboloid-BW}{hyperboloid}} 
\nota{The $\Gamma$-invariant decomposition of the hyperbolic plane into ideal polygons obtained by projecting the Euclidean polygons of Figure \ref{convex-hull:fig} from $\partial C$ to $I^2$. This is in fact the Farey tessellation of Figure \ref{horoballs:fig}.}
\label{hyperboloid:fig}
\end{center}
\end{figure}

The tessellation of $\partial_{\rm t}C$ is $\Gamma$-invariant and projects to a $\Gamma$-invariant tessellation of the hyperbolic space $I^n$ into ideal polyhedra, which projects in turn to a tessellation of $M$ into finitely many ideal polyhedra, called the \emph{canonical} or \emph{Epstein-Penner decomposition}. An example is shown in Figure \ref{hyperboloid:fig}. We have discovered the following.

\begin{teo} \label{EP:teo}
Every finite-volume non-compact complete hyperbolic manifold can be tessellated into finitely many ideal polyhedra.
\end{teo}

\begin{figure}
\begin{center}
\includegraphics[width = 6 cm] {\iftoggle{BW}{TBD}{3ps-1}} 
\includegraphics[width = 6 cm] {\iftoggle{BW}{TBD}{3ps-2}} 
\nota{Different choices of truncated cusps may lead to distinct decompositions into ideal polyhedra.}
\label{3ps:fig}
\end{center}
\end{figure}

The canonical tessellation depends only on the chosen initial truncated cusps $H_1,\ldots, H_c$. Different choices may give different canonical tessellations, an example is sketched in Figure \ref{3ps:fig}. We can parametrize these choices by recording the volumes $V_1,\ldots, V_c$>0 of the truncated cusps (or the areas of $\partial H_i$, that are $V_i/(n-1)$ by Proposition \ref{truncated:volume:prop}).

\begin{prop}
If we multiply all volumes $V_1,\ldots, V_c$ by the same constant $\lambda>0$ we get the same canonical tessellation.
\end{prop}
\begin{proof}
The set $B_j^i$ changes by a global rescaling, and hence also the tessellation does. Its projection on $I^n$ is unaffected.
\end{proof}

A way to get a truly canonical decomposition of $M$, which depends on no choice, consists of taking $V_1=\ldots = V_c>0$. Note that when $c=1$ any choice leads to the same canonical decomposition.

\section{Extensions of homotopies}
This section is entirely devoted to the proof of the following theorem. Recall that every map $f\colon M\to N$ between connected manifolds lifts to a (non-unique) map $\tilde f\colon \tilde M \to \tilde N$ between their universal covers.
\begin{teo} \label{estensione:teo}
Let $f\colon M \to N$ be a smooth homotopy equivalence between closed hyperbolic $n$-manifolds. Every lift $\tilde f\colon \matH^n \to \matH^n$ extends to a continuous map $\tilde f\colon \overline{\matH^n} \to \overline{\matH^n}$ whose restriction $\tilde f|_{\partial \matH^n}\colon \partial \matH^n \to \partial \matH^n$ is a homeomorphism.
\end{teo}

Every smooth homotopy equivalence extends to a homeomorphism of the boundaries of the universal coverings. This theorem will have important consequences in the Chapters \ref{automorphisms:chapter} and \ref{Mostow:chapter}.
To prove it, we need to define a weaker notion of isometry.

\subsection{Quasi-isometries}

We introduce the following.\index{quasi-isometry}

\begin{defn} A map $F\colon X \to Y$ between metric spaces is a \emph{quasi-isometry} if there are two constants $C_1>0$, $C_2\geqslant 0$ such that
$$\frac 1 {C_1} d(x_1,x_2) - C_2 \leqslant d(F(x_1), F(x_2)) \leqslant C_1d(x_1,x_2) + C_2$$
for all $x_1,x_2 \in X$ and $d(F(X),y)\leqslant C_2$ for all $y\in Y$.
\end{defn}
A quasi-isometry is an isometry up to some error: note that $F$ might neither be continuous nor injective. Two metric spaces are \emph{quasi-isometric} if there is a quasi-isometry $F\colon X \to Y$ (which implies the existence of a quasi-isometry $G\colon Y\to X$, exercise) and quasi-isometry is an equivalence relation between metric spaces. Intuitively, looking at a space up to quasi-isometries is like watching it from some distance. Compact metric spaces are obviously quasi-isometric to a point. 

This notion is an important ingredient in \emph{geometric group theory}: one may for instance give every finitely-presented group $G$ a canonical metric (through its \emph{Cayley graph}), uniquely determined up to quasi-isometries.\index{geometric group theory}\index{Cayley graph} 

\subsection{Pseudo-isometries}
Let $f\colon M \to N$ be a homotopic equivalence between closed hyperbolic $n$-manifolds. Every continuous function is homotopic to a smooth one, hence we suppose that $f$ is smooth. The map lifts to a map $\tilde f\colon \matH^n \to \matH^n$. We will prove that $\tilde f$ is a quasi-isometry. Actually, the map $\tilde f$ is also continuous and (as we soon see) Lipschitz: it will be useful for us to retain this information on $\tilde f$ to simplify some arguments, so we introduce a different (more restrictive but less natural) version of a quasi-isometry.\index{pseudo-isometry}

\begin{defn} A map $F\colon X \to Y$ between metric spaces is a \emph{pseudo-isometry} if there are two positive constants $C_1, C_2>0$ such that
$$\frac 1 {C_1} d(x_1,x_2) - C_2 \leqslant d(F(x_1), F(x_2)) \leqslant C_1d(x_1,x_2)$$
for all $x_1,x_2 \in X$.
\end{defn}
In particular a pseudo-isometry is $C_1$-Lipschitz and hence continuous. Let $f\colon M \to N$ be a smooth map between Riemannian $n$-manifolds; the \emph{maximum dilatation} of $f$ at a point $x\in M$ is the maximum ratio $\frac{\|df_x(v)\|}{\|v\|}$ where $v$ varies among all the unitary vectors in $T_x$. The maximum dilatation of $f$ is the supremum of all maximum dilatations as $x \in M$ varies.\index{maximum dilatation} 

\begin{ex}
If $f\colon M\to N$ has maximum dilatation $C$ the map $f$ is $C$-Lipschitz.
\end{ex}

\begin{prop} \label{pseudo-isometria:prop}
Let $f\colon M\to N$ be a smooth homotopy equivalence of closed hyperbolic $n$-manifolds. The lift $\tilde f\colon \matH^n \to \matH^n$ is a pseudo-isometry.
\end{prop}
\begin{proof}
Since $M$ is compact, the map $f$ has some finite maximum dilatation $C$. Since $\tilde f$ is locally like $f$, it also has maximum dilatation $C$ and is hence $C$-Lipschitz. The same holds for the homotopic inverse $g\colon N\to M$. Therefore there is a $C_1>0$ such that
\begin{align*}
d\big(\tilde f(x_1), \tilde f(x_2)\big) & \leqslant C_1 \cdot d(x_1,x_2) \quad \forall x_1,x_2 \in \matH^n, \\
d\big(\tilde g(y_1), \tilde g(y_2)\big) & \leqslant C_1 \cdot d(y_1,y_2) \quad \forall y_1,y_2 \in \matH^n.
\end{align*}
We have $M = \matH^n/_\Gamma$.
Being a composition of lifts, the map $\tilde g \circ \tilde f$ commutes with $\Gamma$ and has maximum displacement bounded by some $K>0$ equal to the maximum displacement of the points in a (compact) Dirichlet domain. Hence
$$d(x_1,x_2) - 2K \leqslant d\big( \tilde g (\tilde f(x_1)), \tilde g (\tilde f(x_2)) \big) \leqslant C_1 \cdot d\big(\tilde f(x_1), \tilde f(x_2)\big) $$
for all $x_1,x_2 \in \matH^n$. Therefore $\tilde f$ is a pseudo-isometry with $C_2 = 2K/C_1$.
\end{proof}

\subsection{Boundary extension of a pseudo-isometry}
We now prove that not only pseudo-isometries $\matH^n \to \matH^n$ extend at infinity, but their extensions at $\partial \matH^n$ eliminate the ``errors'' and behave nicely.
\begin{teo} Every pseudo-isometry $F\colon \matH^n \to \matH^n$ extends to a continuous map $F\colon \overline{\matH^n} \to \overline{\matH^n}$ that injects $\partial\matH^n$ into itself.
\end{teo}
We separate the proof in some lemmas.

\begin{figure}
\begin{center}
\includegraphics[width = 12 cm] {\iftoggle{BW}{cosh-BW}{cosh}} 
\nota{We use the half-space model. The hyperbolic cosine of the distance between $x$ and $\pi(x)$ is the inverse of the cosine of $\theta$ (left). To determine the maximum dilatation we decompose the tangent space $T_x$ orthogonally as $U\oplus V$ (right).}
\label{cosh:fig}
\end{center}
\end{figure}

\begin{lemma} \label{cosh:cos:lemma}
Consider the picture in Figure \ref{cosh:fig}. We have
$$\cosh d(x,\pi(x)) = \frac 1{\cos \theta}.$$
\end{lemma}
\begin{proof} We can work with the half-plane model $H^2\subset \matC$ and up to translations and dilations suppose that $\pi(x)=i$. The geodesic $r$ is parametrized as $ie^t$. The M\"obius transformation $z \mapsto \frac{z+1}{-z+1}$ sends $r$ to $\gamma$ and fixes $i$, hence $\gamma(t) = \frac{ie^{t}+1}{-ie^{t}+1}.$
Set $s=d(x, \pi(x))$. We get $ x = \frac{ie^{s}+1}{-ie^{s}+1}$ and
$$\cos\theta = \Im x = \Im \frac{(ie^s+1)^2}{e^{2s}+1} = \frac{2e^s}{e^{2s}+1} = \frac 2{e^s+e^{-s}} = \frac 1 {\cosh s}.$$
The proof is complete.
\end{proof}

\begin{lemma} \label{dilatation:lemma}
Let $r\subset \matH^n$ be a line and $\pi\colon \matH^n \to r$ be the orthogonal projection to $r$. The maximum dilatation of $\pi$ at $x\in\matH^n$ is
$$d =\frac 1 {\cosh s}$$
where $s=d(x,r)$.
\end{lemma}
\begin{proof}
We use the half-space model with $r$ and $x$ as in Figure \ref{cosh:fig}-(left): we know that $\cosh s = \frac 1 {\cos \theta}$. We have $T_x=U\oplus V$ as in Figure \ref{cosh:fig}-(right) with $V= \ker d\pi_x$. A generator $u$ of $U$ is just rotated by $d\pi_x$ with respect to the Euclidean metric; with respect to the hyperbolic metric  we have
$$\frac{\|d\pi_x(u)\|}{\|u\|} = \frac {x_n}{\pi(x)_n} = \cos\theta = \frac 1{\cosh s}.$$
The proof is complete.
\end{proof}

We write $\overline{pq}$ for the segment with endpoints $p$ and $q$ and we denote by $N_r(A)$ the $r$-neighbourhood of $A$. 

\begin{lemma}
Let $F\colon\matH^n\to\matH^n$ be a pseudo-isometry. There is a $R>0$ such that
$$F(\overline{pq}) \subset N_R\big(\overline{F(p)F(q)}\big)$$ 
for all distinct points $p,q\in\matH^n$. 
\end{lemma}
\begin{proof}
Let $C_1, C_2$ be the pseudo-isometry constants of $F$. Fix a sufficiently big $R$ so that $\cosh R > 2C_1^2$. Let $l$ be the line containing $F(p)$ and $F(q)$.  We show that $F(\overline{pq})$ can exit from $N_R(l)$ only for a limited amount of time. Let $\overline {rs} \subset \overline {pq}$ be a maximal segment where $F(\overline{rs})$ is disjoint from the interior of $N_R(l)$, as the \iftoggle{BW}{dark grey}{blue} arc in Figure \ref{Mostow:fig}-(left). We have

\begin{figure}
\begin{center}
\includegraphics[width = 12 cm] {\iftoggle{BW}{Mostow0-BW}{Mostow0}} 
\nota{We use the \iftoggle{BW}{grey}{red} paths to estimate the distance between $F(r)$ and $F(s)$. On the left: since $F$ is $C_1$-Lipschitz, the \iftoggle{BW}{dark grey}{blue} path has length at most $C_1 d(r,s)$. Its projection onto $l$ has dilatation at most $1/{\cosh R}$ by Lemma \ref{dilatation:lemma}, hence the \iftoggle{BW}{grey}{red} path in $l$ has length at most $C_1d(r,s)/\cosh R$. Therefore $d(F(r), F(s)) \leqslant C_1 \frac {d(r,s)}{\cosh R} + 2R$. On the right we get $d(F(r), F(s)) \leqslant 2R$.}
\label{Mostow:fig}
\end{center}
\end{figure}

\begin{equation*}
\frac 1{C_1} d(r,s) - C_2 \leqslant d(F(r), F(s)) \leqslant C_1d(r,s). 
\end{equation*}
We can improve the right inequality as shown in Figure \ref{Mostow:fig}-(left) and write
$$\frac{1}{C_1}d(r,s) - C_2 \leqslant d(F(r), F(s)) \leqslant C_1\frac{d(r,s)}{\cosh R} + 2R.$$
Therefore
$$\left( \frac 1{C_1} - \frac{C_1}{\cosh R}\right)d(r,s) \leqslant 2R + C_2.$$
Since $\cosh R > 2C_1^2$ we deduce that $d(r,s) < M$ for some constant $M$ that depends only on $C_1$ and $C_2$. 

We have proved that $F(\overline {pq})$ may exit from $N_R(l)$ only on subsegments of $\overline{pq}$ with length $<M$. Since $F$ is $C_1$-Lipchitz the curve $F(\overline{pq})$ lies entirely in $N_{R+C_1M}(l)$, and we replace $R$ with $R+C_1M$.

It remains to prove that $F(\overline{pq})$ lies entirely (up to taking a bigger $R$) in the bounded set $N_R(\overline{F(p)F(q)})$: the proof is analogous and easier, since Figure \ref{Mostow:fig}-(right) shows that $d(F(r), F(s))\leqslant 2R$. 
\end{proof}

In the previous and following lemmas, the constant $R$ depends only on the pseudo-isometry constants $C_1$ and $C_2$.

\begin{figure}
\begin{center}
\includegraphics[width = 10 cm] {\iftoggle{BW}{Mostow2-BW}{Mostow2}} 
\nota{For every $0<u<t$, the point $F(l(u))$ is contained in the (\iftoggle{BW}{light grey}{yellow}) $R$-neighbourhood of $\overline{F(p) F(l(t))}$. If $u$ is big, the \iftoggle{BW}{}{blue }segment $\overline{F(p)F(l(u))}$ is long, while the \iftoggle{BW}{grey}{red} one is bounded by $R$: hence the angle $\alpha_{tu}$ between $v_t$ and $v_u$ is small. Therefore $v_t$ is a Cauchy sequence.}
\label{Mostow2:fig}
\end{center}
\end{figure}

\begin{lemma} \label{half:lemma}
Let $F\colon \matH^n \to \matH^n$ be a pseudo-isometry. There is a $R>0$ such that for all $p\in\matH^n$ and every half-line $l$ starting from $p$ there is a unique half-line $l'$ starting from $F(p)$ such that
$$F(l) \subset N_R(l').$$
\end{lemma}
\begin{proof}
We parametrize $l$ as a geodesic $l\colon [0,+\infty) \to \matH^n$ with unit speed. We have $l(0)=p$. Since $F$ is a pseudo-isometry we get 
$$\lim_{t\to\infty} d\big(F(p), F(l(t))\big) = \infty.$$
Let $v_t\in T_{F(p)}$ be the unitary tangent vector pointing towards $F(l(t))$:
Figure \ref{Mostow2:fig} shows that $\{v_t\}_{t\in\matN}$ is a Cauchy sequence, hence it converges to a unitary vector $v\in T_{F(p)}$. Let $l'$ be the half-line starting from $F(p)$ with direction $v$. It is easy to check that $F(l)\subset N_R(l')$ and $l'$ is the unique half-line from $p$ with this property. 
\end{proof}
The previous lemma gives a recipe to transform every half-line $l$ into a half-line $l'$ that approximates $F(l)$. Since $\partial \matH^n$ is an equivalence relation of half-lines, we define the extension $F\colon\partial\matH^n \to \partial \matH^n$ by sending $l$ to $l'$.

\begin{lemma}
The boundary extension $F\colon \partial \matH^n \to \partial\matH^n$ is well-defined and injective.
\end{lemma}
\begin{proof}
Let $l_1, l_2$ be two half-lines at bounded distance $d(l_1(t), l_2(t))<M$ for all $t$. If $d(l_1'(t), l_2'(t)) \to \infty$ we get $d\big(F(l_1(t)), F(l_2(t))\big) \to +\infty$, a contradiction since $F$ is Lipschitz. Therefore $l_1',l_2'$ are at bounded distance and $F$ is well-defined.

Injectivity is proved analogously: if $l_1$ and $l_2$ are divergent then $l_1'$ and $l_2'$ also are because $F$ is a pseudo-isometry.
\end{proof}

It remains to prove that the extension $F\colon \overline{\matH^n} \to \overline{\matH^n}$ is continuous. We start by extending Lemma \ref{half:lemma} from half-lines to lines.

\begin{lemma}
Let $F\colon \matH^n \to \matH^n$ be a pseudo-isometry. There is a $R>0$ such that for every line $l$ there is a unique line $l'$ with $F(l) \subset N_R(l')$.
\end{lemma}
\begin{proof}
Parametrize $l$ as $l\colon (-\infty,+\infty) \to \matH^n$ with unit speed. By cutting $l$ into two half-lines we know that $F(l(t))$ is a curve that tends to two distinct points $x_\pm\in\partial\matH^n$ as $t\to\pm\infty$. Let $l'$ be the line with endpoints $x_\pm$. For any $t>0$ we have
$$F\big(l([-t,t])\big) \subset N_R\left(\overline{F(l(-t))F(l(t))}\right)$$
and by sending $t\to +\infty$ we deduce that $F(l)\subset N_R(l')$.
\end{proof}

The next lemma says that a pseudo-isometry does not distort much lines and orthogonal hyperplanes. We will need it to prove continuity.

\begin{figure}
\begin{center}
\includegraphics[width = 12 cm] {\iftoggle{BW}{Mostow3-BW}{Mostow3}} 
\nota{Let $l$ and $H$ be a line and an orthogonal hyperplane. The orthogonal projection of $H$ onto $l$ is obviously a point $l\cap H$; the pseudo-isometry $F$ mildly distorts this picture: the image $F(H)$ projects to a bounded segment in $l'$.}
\label{Mostow3:fig}
\end{center}
\end{figure}

\begin{lemma}
Let $F\colon\matH^n \to \matH^n$ be a pseudo-isometry. There is a $R>0$ such that for any line $l$ and hyperplane $H$ orthogonal to $l$, the image $F(H)$ projects orthogonally to $l'$ onto a bounded segment length smaller than $R$.
\end{lemma}
\begin{proof}
See Figure \ref{Mostow3:fig}. Consider a generic line $s\subset H$ passing through $p=l\cap H$. By the previous lemmas $F(s)\subset N_R(s')$ with $s'\neq l'$, and the orthogonal projection on $l'$ sends any other line $s'$ onto a segment, bounded by the images of the endpoints of $s'$. 

\begin{figure}
\begin{center}
\includegraphics[width = 12 cm] {\iftoggle{BW}{Mostow4-BW}{Mostow4}} 
\nota{The lines $s_1$ and $s_2$ have a distance $d$ from $p$ which depends on nothing (in fact, $\cosh d = \sqrt 2$). The lines $l'$, $s_1'$, and $s_2'$ approximate up to an error $R$ the images of $l$, $s_1$, and $s_2$ along $F$. The projection $q$ of $F(p)$ on $l'$ is hence $R$-close to $F(p)$, which is in turn $(C_1d)$-close to the lines $s_i'$. Therefore $q$ is $(C_1d+2R)$-close to both $s_1'$ and $s_2'$. This easily implies that $f$ is $(C_1d+2R)$-close to $q$.
}
\label{Mostow4:fig}
\end{center}
\end{figure}

Consider as in Figure \ref{Mostow4:fig} the line $s$, with one endpoint $s^\infty$ and the corresponding endpoint $F(s^\infty)$ of $s'$. The figure shows that the projection $f$ of $F(s^\infty)$ to $l'$ is at bounded distance from a point $q$ which does not depend on $s$.
\end{proof}

Finally, we prove that the extended $F$ is continuous.

\begin{lemma}
The extension $F\colon \overline{\matH^n} \to \overline{\matH^n}$ is continuous.
\end{lemma}
\begin{proof}
Consider $x\in \partial \matH^n$ and its image $F(x) \in \partial \matH^n$. Let $l$ be a half-line pointing to $x$: hence $l'$ points to $F(x)$. The half-spaces orthogonal to $l'$ determine a neighbourhood system for $F(x)$: consider one such half-space $S$.

Let $R>0$ be as in the previous lemmas. The image $F(l)$ is $R$-close to $l'$, hence for sufficiently big $t$ the point $F(l(t))$ and its projection into $l'$ lie in $S$ at distance $>R$ from $\partial S$. By the previous lemma the image $F(H(t))$ of the  hyperplane $H(t)$ orthogonal to $l$ in $l(t)$ is also contained in $S$. Hence the entire half-space bounded by one such $H(t)$ goes inside $S$ through $F$.
This shows that $F$ is continuous at every point $x\in\partial \matH^n$.
\end{proof}

With some effort, we have proved that every pseudo-isometry of $\matH^n$ extends continuously to the boundary. Theorem \ref{estensione:teo} now follows easily.

\begin{proof}[Proof of Theorem \ref{estensione:teo}]
We know that $\tilde f$ is a pseudo-isometry and hence extends to a map $\tilde f \colon \overline {\matH^n} \to \overline{\matH^n}$ that sends injectively $\partial \matH^n$ to itself. It remains to prove that $\tilde f|_{\partial \matH^n}$ is a homeomorphism.

Pick a smooth homotopic inverse $g$ for $f$. The homotopy $\id_M \sim g \circ f$ lifts to a homotopy $\id_{\matH^n}\sim\tilde g\circ \tilde f$ for some lift $\tilde g$. Since $M$ is compact, the latter homotopy moves every point at uniformly bounded distance and hence $\tilde g\circ \tilde f$ extends continuously to the identity on $\partial \matH^n$, and the same holds for $\tilde f\circ \tilde g$. Therefore $\tilde g|_{\partial \matH^n}$ is the inverse of $\tilde f|_{\partial \matH^n}$ and they are both homeomorphisms.
\end{proof}

\subsection{References}
As in the previous chapters, we have mostly consulted Benedetti -- Petronio \cite{BP}, Ratcliffe \cite{R}, and Thurston's notes \cite{Th}. The Epstein -- Penner decomposition is taken from their 1988 paper \cite{EP}.

%% file: Surfaces.tex
\chapter{Surfaces} \label{superfici:chapter} \label{surfaces:chapter}

A surface is a differentiable manifold of dimension $n=2$, possibly with boundary. The closed orientable surfaces are classified topologically by their Euler characteristic, a complete invariant that also determines their possible geometries: a closed surface has a hyperbolic, flat, or spherical structure if and only if its Euler characteristic is negative, null, or positive.

We devote some time here to expose the topological classification and the geometrisation of closed surfaces. Then we describe some of the beautiful features of geometrisation: the geometry of a surface (typically, hyperbolic geometry) can be used to prove in an elegant way various non-trivial topological facts. 

We end this chapter by defining and studying the \emph{mapping class group}, a group that encodes the topological symmetries of a surface.

\section{Topological classification}
A surface can be topologically quite complicated: think for instance of $\matR^2$ with a Cantor set removed. We decide to restrict our investigation to the surfaces \emph{of finite type}, i.e.~obtained from a closed one by removing points and/or open discs: these include all the compact surfaces with or without boundary.
We introduce here these surfaces and classify them up to diffeomorphism.

We manipulate surfaces using various cut-and-paste tools: boundary gluings, removal of discs or points, handle decompositions, connected sums. See Section \ref{differential:topology:section} to refresh these notions.

\subsection{Gluing surfaces}
A simple way to construct a surface is by gluing simpler surfaces along their boundaries. We show here that the glued surface depends only on the orientation classes of the gluing maps.

Two self-diffeomorphisms of $S^1$ are \emph{co-oriented} if they both preserve or both invert the orientation of $S^1$.
\begin{lemma} \label{diffeo:isotopi:lemma}
Two co-oriented self-diffeomorphisms of $S^1$ are isotopic.
\end{lemma}
\begin{proof}
Let $f_0, f_1\colon S^1 \to S^1$ be two co-oriented self-diffeomorphisms. Their lifts $\tilde f_0, \tilde f_1\colon \matR\to\matR$ are periodic and monotone, hence $\tilde f_t = (1-t)\tilde f_0 + t\tilde f_1$ also is and descends to an isotopy $f_t$ connecting $f_0$ and $f_1$.
\end{proof}

\begin{cor} \label{not:depend:cor}
If we glue two oriented surfaces along their boundaries via orientation-reversing diffeomorphisms, the resulting oriented surface does not depend on the diffeomorphisms chosen.
\end{cor}
\begin{proof}
All the orientation-reversing gluing maps are isotopic by Lemma \ref{diffeo:isotopi:lemma}, and isotopic gluing maps produce diffeomorphic manifolds, as stated in Proposition \ref{depend:diffeo:prop}.
\end{proof}

In a non-oriented context, we have two possible maps for every glued boundary component.

\begin{cor} \label{sphere:cor}
If we glue two discs we get a sphere.
\end{cor}
\begin{proof}
The two discs are copies of $D^2\subset \matC$ and are glued along a diffeomorphism $\varphi\colon S^1\to S^1$. Up to mirroring one disc we may suppose that $\varphi$ is orientation-reversing, and by Corollary \ref{not:depend:cor} we may suppose that $\varphi(z) = \bar z$. The resulting surface is diffeomorphic to a sphere (exercise).
\end{proof}

\begin{warning}
In dimension $n\geqslant 7$, by gluing two discs we may get a manifold homeomorphic but \emph{not} diffeomorphic to a sphere! See Section \ref{self:disc:subsection}.
\end{warning}

\subsection{Classification of surfaces}
In Section \ref{sum:subsection} we introduced the connected sum, a two-steps operation which consists of first removing balls and then gluing the new sphere boundaries. Let $S_g$ be the connected sum 

\begin{figure}
\begin{center}
\includegraphics[width = 3 cm] {\iftoggle{BW}{Sphere_with_three_handles-BW}{500px-Sphere_with_three_handles}}
\quad \includegraphics[width = 3 cm] {\iftoggle{BW}{Triple_torus_array-BW}{500px-Triple_torus_array}}
\quad \includegraphics[width = 3 cm] {\iftoggle{BW}{Triple_torus_illustration-BW}{500px-Triple_torus_illustration}}
\nota{A surface of genus $3$ may be represented in various ways.}
\label{handles_surface:fig}
\end{center}
\end{figure}

$$S_g = \underbrace{T \# \ldots \# T}_g$$ 
of $g$ tori $T=S^1\times S^1$. By convention $S_0 = S^2$ is the sphere and $S_1=T$ is the torus. The number $g$ is the \emph{genus} of the closed surface $S_g$. The surface $S_g$ may be represented in $\matR^3$ in various ways, see Figure \ref{handles_surface:fig}.\index{surface!genus of a surface}

\begin{figure}
\begin{center}
\includegraphics[width = 11 cm] {\iftoggle{BW}{classification-BW}{classification}}
\nota{The 0-handle and $k$ 1-handles form a subsurface $S'\subset S$ with connected boundary, to which the 2-handle is attached (left). A 0-handle and two linked 1-handles (centre). If we attach a 2-handle to the centre figure we get a handle decomposition of the torus (right).}
\label{classification:fig}
\end{center}
\end{figure}

\begin{prop}
We have $\chi(S_g) = 2-2g$.
\end{prop}
\begin{proof}
Let $S,S'$ be surfaces and $D\subset \interior S, D'\subset \interior{S'}$ be discs. Then
$$\chi(S\# S') = \chi(S\setminus D) + \chi(S' \setminus D') - \chi(S^1) = \chi (S) + \chi(S') - 2.$$
Therefore $\chi(S_g) =2-2g$ by induction on $g$.
\end{proof}

\begin{teo}[Classification of surfaces] Every closed, connected, orientable surface is diffeomorphic to $S_g$ for some $g\geqslant 0$.\index{surface!classification of surfaces}
\end{teo}
\begin{proof}
Being a closed manifold, every closed orientable surface $S$ has a handle decomposition. By Proposition \ref{0:handle:prop} it has one with one 0-handle, a certain number $k$ of 1-handles, and one 2-handle. We get $\chi(S) = 2-k$. We prove by induction on $k$ that $k=2g$ is even and $S$ is diffeomorphic to $S_g$.

If $k=0$ then $S$ is obtained by gluing two discs (the 0- and 2-handle) and is hence a sphere by Corollary \ref{sphere:cor}.

Suppose $k>0$. The 0-handle is a disc and the 1-handles are rectangles attached to its boundary as in Figure \ref{classification:fig}-(left). Note that since $S$ is orientable every rectangle is attached without a twist, otherwise it would create a M\"obius strip. The 0- and 1-handles altogether form a compact surface $S'\subset S$ with only one boundary component, to which the 2-handle is attached.

Since $\partial S'$ is connected, every rectangle is linked to some other rectangle as in Figure \ref{classification:fig}-(centre). A pair of linked rectangles form a sub-subsurface $S''\subset S'\subset S$ as in Figure \ref{classification:fig}-(centre) with connected boundary. If we cut $S$ along the curve $\partial S''$ and then cap off with two discs we perform the inverse of a connected sum. 

Therefore $S = S_1\# S_2$, where $S_1$ is $S''$ with a disc attached, \emph{i.e.}~a torus  as Figure \ref{classification:fig}-(right) shows. The surface $S_2$ decomposes into a 0-handle, $k-2$ 1-handles, and one 2-handle. We conclude by induction on $k$.
\end{proof}

\subsection{Homology} \label{homolgy:surface:subsection}
The homology of $S_g$ is easily calculated.
\begin{prop}
We have
$$H_0(S_g,\matZ) = \matZ, \quad H_1(S_g,\matZ) = \matZ^{2g}, \quad H_2(S_g,\matZ) = \matZ.$$
\end{prop}
\begin{proof}
Since $S_g$ is closed, connected, and orientable, we have $H_0(S_g) = H_2(S_g) = \matZ$. Since $\chi(S_g) = 2-2g$, the group $H_1(S_g)$ has rank $2g$. By Poincar\'e duality $H_1(S_g) = H^1(S_g) = \Hom(\pi_1(S_g),\matZ)$ has no torsion.
\end{proof}

\begin{figure}
\begin{center}
\includegraphics[width = 10cm] {\iftoggle{BW}{surface_homology-BW}{surface_homology}}
\nota{These 6 oriented curves represent a symplectic basis for $H_1(S_3,\matZ) = \matZ^6$.}
\label{surface_homology:fig}
\end{center}
\end{figure}

Recall from Section \ref{intersection:subsection} that by fixing an orientation for $S_g$ we get a symplectic intersection form $\omega$ on $H_1(S_g,\matZ)$. 
A basis for $H_1(S_g, \matZ)$ is \emph{symplectic} if $\omega = \matr 0 {I_g} {-I_g} 0$ with respect to this basis: an example is in Figure \ref{surface_homology:fig}.\index{symplectic basis} 

The \emph{algebraic intersection} of two closed curves is the intersection form of their classes in $H_1(S_g, \matZ)$; if the curves are transverse, this is just the algebraic sum of their intersections, where each intersection counts as $\pm 1$ according to the local orientations.\index{algebraic intersection!algebraic intersection of curves}

\subsection{Surfaces of finite type}
We extend our investigation to a larger interesting class of surfaces.\index{surface!surface of finite type} 

\begin{figure}
\begin{center}
\includegraphics[width = 6cm] {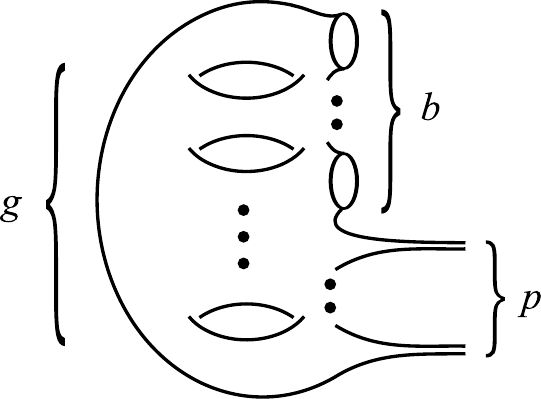}
\nota{The finite type surface $S_{g,b,p}$.}
\label{tipo_finito:fig}
\end{center}
\end{figure}

\begin{defn}
Let $g,b,p\geqslant 0$ be three natural numbers. The \emph{surface of finite type} $S_{g,b,p}$ is the surface obtained from $S_g$ by removing the interior of $b$ disjoint discs and $p$ points. 
\end{defn}
See Figure \ref{tipo_finito:fig}. We say that $S_{g,b,p}$ has genus $g$, has $b$ boundary components, and $p$ punctures. Its Euler characteristic is
$$\chi(S_{g,b,p}) = 2-2g-b-p.$$

We also use the notation $S_{g,b}$ to indicate $S_{g,b,0}$.

\begin{prop}
Every compact connected orientable surface with boundary is diffeomorphic to $S_{g,b}$ for some $g,b\geqslant 0$.
\end{prop}
\begin{proof}
Let $S$ be a compact orientable surface with some $b$ boundary components.
If we glue $b$ discs to $\partial S$ we get a closed orientable surface, hence diffeomorphic to $S_g$ for some $g\geqslant 0$. The original $S$ is obtained from $S_g$ by removing the interiors of $b$ disjoint open discs.
\end{proof}

The compact connected orientable surfaces with $\chi>0$ are $S^2=S_0$ and $D^2=S_{0,1}$, those with $\chi=0$ are the annulus $A=S_{0,2}$ and the torus $T=S_1$.

\begin{cor}[Smooth Jordan curve Theorem] \label{Jordan:cor}
Every smooth simple closed curve $\gamma \subset \matR^2$ bounds a disc.\index{smooth Jordan curve theorem}
\end{cor}
\begin{proof}
Consider the curve $\gamma$ inside $S^2 = \matR^2\cup\{\infty\}$. By cutting $S^2$ along $\gamma$ we get one or two compact orientable surfaces with non-empty boundary and with $\chi=2$ in total. The only possibility is that we get two discs.
\end{proof}

\subsection{Triangulations} \label{surface:triangulations:subsection}
Like every honest smooth compact manifold, compact surfaces admit smooth triangulations, see Section \ref{triangulations:subsection}. Conversely, we may use simplicial complexes to construct smooth surfaces combinatorially.

Let $X$ be a two-dimensional pure simplicial complex where every edge is incident to two faces, and the link of every vertex is a circle.

\begin{prop}
The complex $X$ is the smooth triangulation of a closed surface $S$, unique up to diffeomorphism.
\end{prop}
\begin{proof}
By dualising $X$ we get a handle decomposition: triangles, edges, and vertices transform into 0-, 1-, and 2-handles, and we get a smooth surface $S$ triangulated by $X$. The way the handles are attached is determined up to isotopy, therefore $S$ is determined up to diffeomorphism.
\end{proof}

\begin{warning}
It is worth noting that this procedure (getting a unique smooth structure from a simplicial complex) does not work in all dimensions (here we used implicitly Lemma \ref{diffeo:isotopi:lemma}).
\end{warning}

\section{Geometrisation} \label{surface:geometrisation:section}
We now prove that every surface $S_{g,b,p}$ of finite type can be \emph{geometrised}, that is it may be equipped with a hyperbolic, flat, or elliptic metric. The metric type is prescribed by the sign of the Euler characteristic of the surface.

\subsection{Hyperbolic pair-of-pants} \label{pantaloni:subsection}
The sphere has of course an elliptic structure, and the torus has many flat structures, see Section \ref{piatte:subsection}. We now construct hyperbolic structures on all the surfaces $S_g$ of genus $g\geqslant 2$, and more generally on all the surfaces $S_{g,b,p}$ of negative Euler characteristic.
We start with a simple block, the \emph{pair-of-pants} $S_{0,3}$, with Euler characteristic $-1$.\index{pair-of-pants}

\begin{prop} \label{pantaloni:prop}
Given three real numbers $a,b,c\geqslant 0$ there is (up to isometries) a unique complete finite-volume hyperbolic pair-of-pants with geodesic boundary, with boundary curves of length $a$, $b$, and $c$.
\end{prop}

\begin{figure}
\begin{center}
\includegraphics[width = 12.5 cm] {\iftoggle{BW}{pantalone_degenere-BW}{pantalone_degenere}}
\nota{A pair-of-pants, an annulus with one puncture, a disc with two punctures, and a thrice-punctured sphere. The last three surfaces may be considered as some degenerated hyperbolic pairs-of-pants where one or more boundary lengths $a,b$, or $c$ are zero, and we get cusps instead of geodesic boundary components there.}
\label{pantalone_degenere:fig}
\end{center}
\end{figure}

When some length in $a,b,c$ is zero, we mean that the geodesic boundary is actually a cusp (recall Section \ref{finite-volume:subsection}) and hence the surface is topologically a punctured annulus $S_{0,2,1}$, a twice punctured disc $S_{0,1,2}$, or a thrice-punctured sphere $S_{0,0,3}$: see Figure \ref{pantalone_degenere:fig}. 

To prove this proposition we construct some right-angled hexagons in $\matH^2$ as in Figure \ref{esagono:fig}-(left). Three \emph{alternate sides} on a hexagon are three pairwise non-incident sides, like the $a,b,c$ shown in the figure. A \emph{degenerate} hexagon is one where the length of some alternate sides is zero as in Figure \ref{esagoni_generalizzati:fig}.

\begin{figure}
\begin{center}
\includegraphics[width = 11cm] {\iftoggle{BW}{esagono-BW}{esagono}}
\nota{A right-angled hexagon with alternate sides of length $a$, $b$ and $c$ (left) and its construction (right), which goes as follows: take a line $l$ with two arbitrary points $A$ and $B$ in it  (bottom black). Draw the perpendiculars from $A$ and $B$ (\iftoggle{BW}{grey}{red}). At distances $a$ and $b$ we find two points $A'$ and $B'$ and we draw again two perpendiculars (black) $r$ and $s$, with some points at infinity $P$ and $Q$. 
Draw the (unique) perpendiculars to $l$ pointing to $P$ and $Q$ (\iftoggle{BW}{dark grey}{blue}): they intersect $l$ in two points $T$ and $U$. Note that $AT$ and $UB$ have some fixed length depending only on $a$ and $b$. We can vary the parameter $x=TU$: if $x>0$ the lines $r$ and $s$ are ultra-parallel and there is a unique segment orthogonal to both of some length $f(x)$.} 
\label{esagono:fig}
\end{center}
\end{figure}

\begin{figure}
\begin{center}
\includegraphics[width = 12.5 cm] {\iftoggle{BW}{esagoni_generalizzati-BW}{esagoni_generalizzati}}
\nota{A right-angled hexagon with parameters $a,b,c\geqslant 0$ degenerates to a pentagon, quadrilateral, or triangle with ideal vertices if one, two, or three parameters are zero.}
\label{esagoni_generalizzati:fig}
\end{center}
\end{figure}

\begin{lemma} \label{esagono:lemma}
Given three real numbers $a,b,c\geqslant0$ there exists (up to isometries) a unique (possibly degenerate) hyperbolic right-angled hexagon with three alternate sides of length $a$, $b$, and $c$.
\end{lemma}
\begin{proof}
We first suppose $a,b>0$. The construction of the hexagon is depicted in Figure \ref{esagono:fig}-(right). If $x=0$ the \iftoggle{BW}{dark grey}{blue} lines coincide, hence $P=Q$ and $f(0)=0$. The function $f\colon [0,+\infty) \to [0,+\infty)$ is continuous, strictly monotonic, and with $\lim_{x\to \infty} f(x) = \infty$: therefore there is precisely one $x$ such that $f(x)=c$. 

If exactly two parameters are zero, say $a=b=0$, a simpler construction works: take a segment of length $c$ as in Figure \ref{esagoni_generalizzati:fig}-(centre), draw the perpendiculars at their endpoints, and a line connecting the endpoints of these. 

If $a=b=c=0$, pick any ideal triangle. Ideal triangles are indeed unique up to isometry: use the half-space model and recall that $\PSLR$ acts transitively on the ordered triples of points in $\partial \matH^2$.
\end{proof}
The most degenerate case is so important that we single it out.
\begin{cor}
All the ideal triangles in $\matH^2$ are isometric.
\end{cor}

\begin{figure}
\begin{center}
\includegraphics[width = 10cm] {\iftoggle{BW}{pantalone-BW}{pantalone}}
\nota{By gluing two identical right-angled hexagons along their black sides we get a hyperbolic pair-of-pants with geodesic boundary.}
\label{pantalone:fig}
\end{center}
\end{figure}

By gluing two identical (possibly degenerate) hexagons along alternate sides as in Figure \ref{pantalone:fig} we construct a (possibly degenerate) hyperbolic pair-of-pants whose geodesic boundary consists of three simple closed geodesics of length $2a$, $2b$, and $2c$. 

\begin{proof}[Proof of Proposition \ref{pantaloni:prop}]
We have proved the existence of a geodesic pairs-of-pants with any parameters $a,b,c\geqslant 0$, and we now turn to its uniqueness. Let $P$ be a pair-of-pants with geodesic boundary curves $C_1, C_2, C_3$ of length $a, b, c>0$.

Since $P, C_1$, and $C_2$ are compact, there are points $x_1\in C_1$ and $x_2 \in C_2$ at minimum distance $d=d(x_1,x_2)$, connected by some curve $\gamma_3$ of length $d$. The curve is a simple geodesic orthogonal to both $C_1$ and $C_2$: if not, some other curve connecting $x_1$ and $x_2$ would be shorter. We construct analogously two orthogeodesics $\gamma_1$ and $\gamma_2$ connecting $C_2$ to $C_3$ and $C_3$ to $C_1$ having minimal length.

The fact that $\gamma_1,\gamma_2,\gamma_3$ have minimal length easily implies that they are disjoint (if they intersect, we find shorter curves). The three orthogeodesics subdivide $P$ into two hexagons, with alternate sides of length $L(\gamma_1)$, $L(\gamma_2)$, and $L(\gamma_3)$: by Lemma \ref{esagono:lemma} the two hexagons are isometric, and hence the three other alternating sides also have the same length $\frac a2, \frac b2$, and $\frac c2$. Hexagons are unique up to isometry and hence the original pair-of-pants also are.

We can extend the argument to the more general case $2a, 2b, 2c \geqslant 0$ as follows. If $a=0$, a neighbourhood of the puncture is a cusp, and we truncate it at some small horocycle $C_1$ (we do the same if $b=0$ or $c=0$). After these truncations we get a compact pair-of-pants (whose boundary is \emph{not} geodesic at the horocycles) and decompose it into two hexagons as above. The resulting curve $\gamma_3$ is orthogonal to the horocycle $C_1$, hence it extends to a half-line pointing towards the puncture. The curves $\gamma_1, \gamma_2, \gamma_3$ decompose the surface into degenerate hexagons.
\end{proof}

\subsection{Hyperbolic surfaces} \label{hyperbolic:surfaces:subsection}
The pairs-of-pants can be used as building blocks to construct topologically all finite type surfaces with $\chi<0$.

\begin{figure}
\begin{center}
\includegraphics[width = 8cm] {\iftoggle{BW}{surface-BW}{surface}}
\nota{Every surface of finite type with $\chi<0$ decomposes into pair-of-pants. We show here a decomposition of $S_3$.}
\label{surface:fig}
\end{center}
\end{figure}

\begin{prop} \label{decompone:pantaloni:prop}
If $\chi (S_{g,b,p}) <0$ then $S_{g,b,p}$ decomposes topologically into $-\chi(S_{g,b,p})$ (possibly degenerate) pairs-of-pants.
\end{prop}
\begin{proof}
If $b+p=0$ then $g\geqslant 2$ and the surface decomposes easily in many ways, see for instance Figure \ref{surface:fig}. If $b+p>0$ and $\chi < -1$, a decomposition for $S_{g,b,p}$ may be obtained from one of $S_{g,b-1,p}$ or $S_{g,b,p-1}$ by inserting one more (possibly degenerate) pair-of-pants. If $\chi=-1$ the surface is either a pair-of-pants, or a torus with a puncture or boundary component, which is in turn obtained by glueing two boundary components of a pair-of-pants.
\end{proof}

We can use this building block to construct hyperbolic structures.

\begin{cor} \label{ammette:cor}
If $\chi(S_{g,b,p})<0$ then $S_{g,b,p}$ admits a complete hyperbolic metric with $b$ geodesic boundary components of arbitrary length.
\end{cor}
\begin{proof}
Decompose $S_{g,b,p}$ in pair-of-pants, assign an arbitrary length to all the closed curves of the decomposition (the 6 \iftoggle{BW}{grey}{red} curves shown in Figure \ref{surface:fig}) and give each pair-of-pants the hyperbolic metric determined by the three assigned boundary lengths. Everything glues to a complete hyperbolic metric for $S_{g,b,p}$ thanks to Proposition \ref{geodesic:glue:prop}.
\end{proof}

Can we geometrise the few orientable surfaces with $\chi \geqslant 0$? Yes, but since there are no cusps in the elliptic and flat geometries we do not consider surfaces with punctures. The compact orientable surfaces with $\chi >0$ are the sphere and the disc, and they all have an elliptic metric with geodesic boundary (represent the disc as a hemisphere). Those with $\chi = 0$ are the torus and the annulus, and they admit flat metrics with geodesic boundary.

\subsection{Non-orientable surfaces}
The classification of all the finite-type non-orientable surfaces is also simple. Let $S^{\rm no}_g$ be the connected sum
$$S_g^{\rm no} = \underbrace{\matRP^2 \# \ldots \# \matRP^2}_g$$ 
of $g\geqslant 1$ copies of the projective plane $\matRP^2$.

\begin{figure}
\begin{center}
\includegraphics[width = 6 cm] {\iftoggle{BW}{Klein_bottle-BW}{Klein_bottle}}
\nota{The Klein bottle immersed in $\matR^3$.}
\label{Klein_bottle:fig}
\end{center}
\end{figure}

\begin{prop}
We have $\chi(S_g^{\rm no}) = 2-g$. The surface $S_2^{\rm no}$ is diffeomorphic to the Klein bottle $K$ shown in Figure \ref{Klein_bottle:fig}. We have
$$S_g \#\matRP^2 \isom S^{\rm no}_{2g+1}.$$
\end{prop}
\begin{proof}
The formula for $\chi$ follows from $\chi(\matRP^2)=1$. The Klein bottle $K$ may be cut along a closed curve into two M\"obius strips, and $\matRP^2$ minus an open disc is a M\"obius strip too, hence $\matRP^2\#\matRP^2 = K$. The latter equality is a consequence of $T\#\matRP^2 \isom K \# \matRP^2$ which is left as an exercise.
\end{proof}

\begin{prop}
Every closed, connected, non-orientable surface is diffeomorphic to $S_g^{\rm no}$ for some $g\geqslant 1$.
\end{prop}
\begin{proof}
Pick a handle decomposition of the surface $S$. Since it is non-orientable, at least one 1-handle is twisted and forms a M\"bius strip. We have proved that $S$ contains a M\"obius strip, and we now remove it and substitute it with a disc to get a new surface $S'$. We have $S=S'\#\matRP^2$ and we conclude by induction on $-\chi(S')$.
\end{proof}

We may also denote by $S^{\rm no}_{g,b,p}$ the surface obtained from $S^{\rm no}_g$ by removing the interiors of $b$ discs and $p$ points.

\begin{ex}
If $\chi(S^{\rm no}_{g,b,p})<0$ then $S^{\rm no}_{g,b,p}$ decomposes into pairs-of-pants and hence admits a complete hyperbolic metric with geodesic boundaries of arbitrary length.
\end{ex}

As in the orientable case, compact surfaces with $\chi \geqslant 0$ can be geometrised: these are $\matRP^2$, the M\"obius strip, and the Klein bottle. The first has an elliptic structure while the other two have flat structures with geodesic boundary.

\subsection{Orbifolds} \label{orbifold:surface:subsection}
We can push the classification and geometrisation further to compact two-dimensional orbifolds. 

The finite subgroups of $\On(2)$ are: the cyclic $\matZ_2$ generated by a reflection, the cyclic $\matZ_n$ generated by a rotation, and the dihedral $D_{2n}$ containing reflections and rotations. Therefore every singular point $x$ on a two-dimensional orbifold $O$ is locally of one of these types:
$$V/_{\matZ_2}, \quad V/_{\matZ_n}, \quad V/_{D_{2n}}$$
with $V$ the unit ball in $\matR^2$. The point $x$ is called respectively a \emph{mirror}, a \emph{cone}, and a \emph{corner reflector} point of $O$, see Figure \ref{orbifold_types:fig}.

\begin{figure}
\begin{center}
\includegraphics[width = 8 cm] {\iftoggle{BW}{orbifold_types-BW}{orbifold_types}}
\nota{A mirror, a cone point, and a corner reflector on a two-dimensional orbifold.}
\label{orbifold_types:fig}
\end{center}
\end{figure}

We consider for simplicity only \emph{locally orientable} orbifolds, \emph{i.e.}~orbifolds with isotropy groups in $\SO(2)$. In other words, we exclude mirrors and corner reflectors. One such orbifold $O$ is easily encoded as
$$(S, p_1, \ldots, p_k)$$
where $S$ is a surface, and the orbifold $O$ is $S$ with $k$ cone points with rotational isotropy groups of order $p_1, \ldots, p_k >1$. We do not require $O$ to be globally orientable, so $S$ can be a non-orientable surface like $\matRP^2$. We define the \emph{Euler characteristic} of that orbifold as\index{Euler characteristic!of an orbifold} 
$$\chi(O) = \chi(S) - \sum \left( 1 - \frac 1{p_i} \right).$$
This definition is designed to behave well under coverings. The \emph{degree} of an orbifold covering $O\to O'$ is the cardinality of the fibre of any non-singular point in $O'$.
\begin{prop}
If $O \to O'$ is a degree-$d$ orbifold covering then 
$$\chi(O) = d \cdot \chi(O').$$
\end{prop}
\begin{proof}
Pick a cellularisation of $O'$ whose vertices contain the cone points: it lifts to a similar cellularisation of $O$. We consider all the vertices of these cellularisations as cone points, possibly with order 1, and we get
$$\chi(O') = V - E + F - \sum \left(1-\frac 1{p_i'}\right) = - E + F + \sum \frac 1 {p_i'}$$
where $V, E, F$ is the number of vertices, edges, and faces in the cellularisation of $O'$. The same formula holds for $O$: 
$$\chi(O) = - dE + dF + \sum \frac 1 {p_j}.$$
A cone point of order $p_i'$ in $O'$ lifts to some cone points with orders $p_{j_1}, \ldots, p_{j_l}$ of $O$, such that
$$d = \sum_{a=1}^l \frac {p_i'}{p_{j_a}} \quad \Longrightarrow \quad 
d \cdot \frac 1 {p_i'} = \sum_{a=1}^l \frac 1{p_{j_a}}.$$
This implies that $\chi(O) = d\cdot \chi(O')$.
\end{proof}

\subsection{Geometrisation of orbifolds}
Recall that an orbifold is good if it is covered by a manifold.
With a couple of exceptions, all the closed two-orbifolds are good and geometric.\index{orbifold!geometrisation of orbifolds}

\begin{table}
\begin{center}
\begin{tabular}{c||c}
\phantom{\Big|} \!\! orbifold & $\Gamma$ \\
\hline \hline
\phantom{\Big|} \!\! $S^2$ & \{e\} \\
\phantom{\big|} \!\! $(S^2, n, n)$ & $\matZ_n$ \\ 
\phantom{\Big|} \!\! $\matRP^2$ & $\matZ_2$ \\ 
\phantom{\big|} \!\! $(\matRP^2, n)$ & $D_{2n}$
\end{tabular}
\quad
\begin{tabular}{c||c}
\phantom{\Big|} \!\! orbifold & $\Gamma$ \\
\hline \hline
\phantom{\Big|} \!\! $(S^2, 2, 2, n)$ & $D_{2n}$ \\
\phantom{\big|} \!\! $(S^2, 2, 3, 3)$ & $T_{12}$ \\
\phantom{\Big|} \!\! $(S^2, 2, 3, 4)$ & $O_{24}$ \\
\phantom{\big|} \!\! $(S^2, 2, 3, 5)$ & $I_{60}$ 
\end{tabular}
\vspace{.2 cm}
\nota{Except the two bad cases, every locally orientable closed orbifold $O$ with $\chi(O)>0$ can be geometrised as $O=S^2/_{\Gamma}$ for some discrete $\Gamma <\On(3)$. In all cases $\Gamma < \SO(3)$ and hence $O$ is orientable, except $\matRP^2$ and $(\matRP^2,n)$ whose $\Gamma$ contains the antipodal map. 
The groups on the right table are the spherical Von Dyck groups, see Section \ref{more:examples:subsection}. 
The groups $T_{12}, O_{24}, I_{60}$ are the orientation-preserving isometry groups of the tetrahedron, octahedron, and icosahedron, and are isomorphic to $A_4$, $S_4$, $A_5$ respectively.}
\label{orbifold:elliptic:table}
\end{center}
\end{table}

\begin{teo} \label{O:geometrisation:teo}
Every closed locally orientable 2-orbifold $O$ is good except
$$(S^2, p), \quad (S^2,p_1,p_2)$$
with $p_1 \neq p_2$. A good orbifold $O$ has an elliptic, flat, hyperbolic structure $\Longleftrightarrow$ $\chi(O)$ is positive, zero, negative. The elliptic orbifolds are listed in Table \ref{orbifold:elliptic:table}.
\end{teo}
\begin{proof}
The orbifolds $(S^2, p)$ and $(S^2,p_1,p_2)$ are bad: subdivide $S^2$ into two discs, each containing at most one cone point; each disc has a unique surface (disc) covering, but the two coverings do not match if $p\neq 1$ or $p_1\neq p_2$.
The remaining orbifolds with $\chi(O)>0$ are 
$$S^2,\ (S^2, n, n),\ \matRP^2,\ (\matRP^2, n),$$ 
$$(S^2, 2, 2, n),\ (S^2, 2, 3, 3),\ (S^2, 2, 3, 4),\ (S^2, 2, 3, 5)$$
and they can all be realised as $S^2/_\Gamma$ for an appropriate finite $\Gamma < \On(3)$, see Table \ref{orbifold:elliptic:table}. The orbifolds with $\chi(O)=0$ are
$$(S^2, 2,3,6), \ (S^2,2,4,4),\ (S^2, 3,3,3),\ (S^2, 2,2,2,2),$$
$$(\matRP^2, 2,2),\ K,\ T$$
where $K$ and $T$ are the Klein bottle and the torus. The first four orbifolds were obtained as $\matR^2/_\Gamma$ in Section \ref{more:examples:subsection}. The orbifold $(\matRP^2, 2,2)$ is $\matR^2/_\Gamma$ where $\Gamma$ is generated by two glide reflections with orthogonal axis:
$$(x,y) \mapsto (x+1, -y), \quad (x,y) \mapsto (-x, y+1).$$
To prove this, consider the fundamental domain $[-\frac 12, \frac 12]\times [-\frac 12, \frac 12]$.

If $\chi(O)<0$ we construct a hyperbolic metric as we did in the surface case. If $O = (S^2, p_1, p_2, p_3)$ we get $O = \matH^2/_{\Gamma^{\rm or}(p_1,p_2,p_3)}$ using Von Dyck groups, see Section \ref{more:examples:subsection}. Otherwise, the orbifold $O$ easily decomposes along disjoint simple closed curves (that avoid the cone points) into some basic pieces with $\chi<0$, which are of the following kind: 
$$P,\ (A, p),\ (D, p_1, p_2)$$
where $P,A,D$ are the pair-of-pants, the annulus, and the disc, and $(p_1,p_2)\neq (2,2)$. We can give a hyperbolic structure to $P$ with any fixed length at the boundaries. 

Analogously, we can give a cone manifold structure with geodesic boundary to $(A,p)$ (respectively, $(D,p_1,p_2)$) with one (two) cone point of angle $\frac{2\pi}p$ ($\frac{2\pi}{p_1}$ and $\frac{2\pi}{p_2}$), for any fixed lengths at the boundaries. We already know this for $P$, and the proof for $(A,p)$ and $(D,p_1,p_2)$ is similar: instead of constructing right-angled hexagons, we construct pentagons with angles $\frac \pi 2, \frac \pi 2, \frac \pi 2, \frac \pi 2, \frac \pi p$ and quadrilaterals with angles $\frac \pi 2, \frac \pi 2, \frac \pi {p_1}, \frac \pi {p_2}$, and we double them, see Figure \ref{orbifold_pieces:fig}.

\begin{figure}
\begin{center}
\includegraphics[width = 12.5 cm] {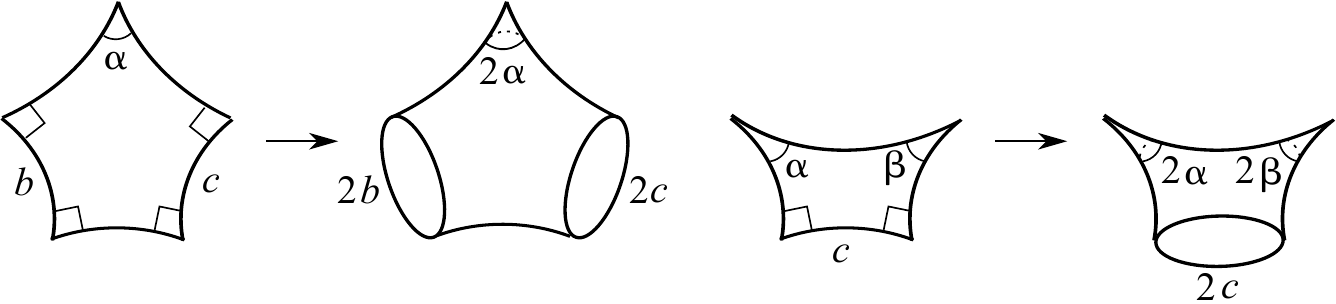}
\nota{For every angle $\alpha\leqslant \frac \pi 2$ and lengths $b,c>0$ there is a pentagon as shown (left); by doubling it along the unmarked edges we get a hyperbolic annulus with a cone point of angle $2\alpha$. For every angles $\alpha, \beta\leqslant \frac \pi 2$ and every length $c>0$ there is a quadrilateral as shown (right). By doubling it along the unmarked edges we get a hyperbolic disc with cone points of angle $2\alpha$ and $2\beta$.}
\label{orbifold_pieces:fig}
\end{center}
\end{figure}

By gluing all the pieces we obtain a hyperbolic cone structure on $O$ with cone angles $\frac{2\pi}{p_i}$ and we apply Proposition \ref{from:cone:to:orbifold:prop} to translate it to a hyperbolic orbifold structure.
\end{proof}

Recall that an orbifold is very good if it is finitely covered by a manifold. In the following corollary we use the geometrisation to prove shortly a non-trivial topological fact.

\begin{cor} \label{very:good:surface:cor}
Every good closed locally orientable 2-orbifold is very good.
\end{cor}
\begin{proof}
Every such orbifold is geometric, so Corollary \ref{very:good:cor} applies.
\end{proof}

The punctured case is also interesting when $\chi(O)<0$. We define a \emph{punctured orbifold} to be one of type
$$O = (S_{g,0,b}, p_1, \ldots, p_k).$$
\begin{teo}
Every punctured orbifold $O$ with $\chi(O)<0$ has a hyperbolic structure.
\end{teo}
\begin{proof}
Same proof as above, where we use triangular groups with possibly ideal vertices and we allow the pieces $P$, $(A, p)$, $(D, p_1, p_2)$ to have boundary components and/or cone points that degenerate to cusps. To construct them we build hexagons, pentagons, quadrilaterals with some lengths or angles that may be zero.
\end{proof}

These punctured orbifolds are also very good for the same reason above. It is also possible to consider orbifolds based on surfaces $S_{g,b,p}$ with boundary: every boundary component can be considered in two natural but distinct ways, either as a mirror (hence the orbifold is locally non-orientable), or as a boundary component: to do this we need to define an appropriate notion of orbifold with boundary.

\subsection{Gauss--Bonnet formula}
How can we compute the area of a hyperbolic surface? We can easily answer this question in the closed case.\index{Gauss--Bonnet formula} 

\begin{prop}[Gauss--Bonnet] \label{GB:prop}
Let $S$ be a closed hyperbolic surface. We have
$$\Area (S) = -2\pi \chi(S).$$
\end{prop}
\begin{proof}
Pick a Dirichlet domain $D$ for $S$. It is a $(2n)$-gon with inner angles $\alpha_1, \ldots, \alpha_{2n}$. By Proposition \ref{polygon:area:prop} we have
$$\Area(S) = \Area(D) = (2n-2)\pi - \sum_{i=1}^{2n} \alpha_i.$$
The surface $S$ is obtained from $D$ by gluing isometrically some pairs of edges. We get a cellularisation of $S$ with some $v$ vertices, $n$ edges, and one face $D$. Therefore 
$$\chi(S) = v - n + 1.$$
Each of the $v$ vertices is obtained by identifying some of the vertices of $D$, whose inner angles must sum to $2 \pi$. Therefore we get
$$\sum_{i=1}^{2n} \alpha_i = 2\pi v$$
and hence 
$$\Area(S) = (2n-2)\pi - 2\pi v = -2\pi(v-n+1) = -2\pi \chi(S).$$
The proof is complete.
\end{proof}
We can extend the formula to orbifolds and to the other geometries. 
\begin{cor}
Let $O$ be a closed hyperbolic, flat, or elliptic locally orientable 2-orbifold. We have
$$K\cdot\Area(O) = 2\pi \chi(O)$$
where $K=-1$, $0$, or $+1$ respectively.
\end{cor}
\begin{proof}
If $O$ is an orientable surface, then the equality holds: in the hyperbolic case this is Proposition \ref{GB:prop}, in the flat case we get $0=0$, and in the elliptic case $S=S^2$ and $\Area(S^2) = 4\pi = 2\pi \chi(S^2)$.

In general there is a degree-$d$ covering $S\to O$ of an orientable geometric surface $S$ by Corollary \ref{very:good:surface:cor}. The equality holds for $S$ and hence also for $O$ since 
$$\chi(S) = d \cdot \chi(O), \qquad \Area (S) = d \cdot \Area (O).$$
The proof is complete.
\end{proof}

\subsection{Lie subgroups of $\SO(3)$}
The following proposition is of general interest and is used in many different contexts.
It shows, among other things, that orientable elliptic orbifolds are \emph{rigid}: their geometry is fully determined by their fundamental group. (Flat and hyperbolic surfaces are not rigid, as we will soon see.)\index{$\SO(3)$}

\begin{prop} \label{subgroups:SO3:prop}
Every non-trivial proper Lie subgroup of $\SO(3)$ is conjugate to one of the following:
$$C_n, \ D_{2n},\ T_{12},\ O_{24},\ I_{60},\ \SO(2),\ \SO(2) \rtimes C_2.$$
These are the orientation-preserving isometry groups of: a regular $n$-pyramid, $n$-prism, tetrahedron, octahedron, icosahedron, cone, and cylinder.
\end{prop}
\begin{proof}
The Lie algebra $\frakg = \frakso(3,\matR)$ is generated by the matrices
$$A_x = \begin{pmatrix} 0 & 0 & 0 \\ 0 & 0 & -1 \\ 0 & 1 & 0 \end{pmatrix}, \quad
A_y = \begin{pmatrix} 0 & 0 & 1 \\ 0 & 0 & 0 \\ -1 & 0 & 0 \end{pmatrix}, \quad
A_z = \begin{pmatrix} 0 & -1 & 0 \\ 1 & 0 & 0 \\ 0 & 0 & 0 \end{pmatrix} $$
with the relations
$$[A_x, A_y] = A_z, \quad [A_y, A_z] = A_x, \quad [A_z, A_x] = A_y.$$
Therefore $\frakg$ is isomorphic to $\matR^3$ equipped with the standard vector product $\times$. We deduce that the only non-trivial proper sub-algebras of $\frakg$ have dimension 1, and hence every non-trivial proper Lie subgroup $G<\SO(3)$ has dimension $0$ or $1$.

If a Lie subgroup $G<\SO(3)$ has dimension zero, it is finite. Every non-trivial element in $G$ is a rotation along some axis, and hence it acts on $S^2$ with two (antipodal) fixed points.
Let $P\subset S^2$ be the set of the fixed points of all the non-trivial elements of $G$. 

The group $G$ clearly acts on $P$, and a double-counting argument gives
$$2\big(|G| - 1\big) = \sum_{p \in P}\big(|G_p| - 1\big)$$
where $G_p < G$ is the stabiliser of $p$. The fundamental theorem on group actions says that 
$$|G_p| = \frac{|G|}{|O(p)|}$$
where $O(p)$ is the orbit of $p$. Therefore
$$2\big(|G| - 1\big) = \sum_{p \in P}\left(\frac{|G|}{|O(p)|} - 1\right) = 
\sum_{O}\big({|G|} - |O|\big)
$$
where in the latter we sum on orbits $O$. We divide by $|G|$ to get
\begin{equation} \label{G:orbit:eqn}
2-\frac 2{|G|}  = \sum_{i=1}^r \left(1- \frac 1{a_i}\right)
\end{equation}
where $a_1,\ldots, a_r$ are the orders of the stabilisers of the $r$ orbits.

If $r=1$ we get $a_1 = |G| = 1$ and hence $G$ is trivial, a contradiction. If $r=2$ we get $a_1 = a_2 = |G|$ and hence every orbit is a single point, with stabiliser $G$. Therefore the points are antipodal and $G=C_n$.

If $r=3$ then $(a_1,a_2, a_3)$ is a triple with $\frac 1{a_1} + \frac 1{a_2} + \frac 1{a_3} > 1$, that is one of
$$(2,2,n), \quad (2,3,3), \quad (2,3,4), \quad (2,3,5).$$
In each case one calculates $|G|$ using (\ref{G:orbit:eqn}) and deduces the cardinalities of the three orbits. For instance in the $(2,2,n)$ case we have $|G|=2n$ and the orbits have cardinality $(n,n,2)$. One deduces that $G$ contains an index-two cyclic group $C_n$ and is the symmetry group $D_{2n}$ of the prism.

In the $(2,3,3)$ case we get $|G|=12$ and hence the orbits have order $(6,4,4)$. The stabiliser of a point in the order-4 orbit $O$ is a $\frac{2\pi}3$-rotation that rotates the other three points: one deduces that the points in $O$ are the vertices of a regular tetrahedron and $G$ is its symmetry group. The other cases are treated analogously.

If $G$ has dimension one, it contains the group $\SO(2)$ of all rotations along some axis $r$. Either $G=\SO(2)$, or there is a $g \in G \setminus \SO(2)$. If $g(r)\neq r$, then $G$ contains another $\SO(2)$ of rotations along $g(r)$, and these two copies of $\SO(2)$ generate a Lie subgroup of dimension $\geqslant 2$ in $G$, a contradiction. Therefore every element in $G\setminus \SO(2)$ preserves $r$ and is hence a $\pi$-rotation along an axis perpendicular to $r$. Therefore $G=\SO(2) \rtimes C_2$.
\end{proof}

\section{Curves on surfaces} \label{curves:surfaces:section}
It is natural to study a manifold by examining the lower-dimensional sub-manifolds that it contains: we now look at closed curves in surfaces and unveil an unexpectedly rich world. We will discover that the study of curves in surfaces is tightly connected to hyperbolic geometry: 
each of the two topics seems designed to give us a better understanding of the other.

We first prove various facts on curves, and then we use this knowledge to deduce some topological non-trivial results on surfaces. We will end up by showing that two self-diffeomorphisms of a closed surface $S_g$ are homotopic if and only if they are isotopic, a theorem that will be important also in our topological study of three-manifolds in the subsequent chapters.

We concentrate for simplicity on the genus-$g$ closed orientable surfaces $S_g$, although much of the discussion could be easily extended to finite-type surfaces $S_{g,b,p}$ with the appropriate modifications. 

\subsection{Definitions}
We start by recalling and fixing some basic definitions.
A \emph{curve} on a differentiable manifold $M$ is a smooth map $\gamma \colon I\to M$ defined on some interval $I$, while a \emph{closed curve} is a smooth map $\gamma\colon S^1\to M$.\index{closed curve}\index{regular curve} 

A (possibly closed) curve $\gamma$ is \emph{regular} if $\gamma'(t)\neq 0$ for all $t$. All the curves will be tacitly assumed to be regular; moreover, with a little abuse we will sometimes indicate by $\gamma$ the \emph{support} of the curve $\gamma$, that is its image.

A curve is \emph{simple} if it is an embedding, and the support of a simple closed curve is a one-dimensional submanifold of $M$ diffeomorphic to $S^1$. Recall that all isotopies are smooth by assumption, see Section \ref{isotopy:subsection}. Two simple closed curves are isotopic if and only if they are ambiently isotopic. Two simple closed curves with the same support and orientation are isotopic by Lemma \ref{diffeo:isotopi:lemma}. 

\subsection{Simple closed curves on the sphere}
We would like to classify the simple closed curves on a given closed surface up to isotopy, but we will soon realise that this task is harder than one could expect. A closed surface may contain many complicated closed curves, difficult to draw and visualise.

For the moment we content ourselves with a very simple case. A simple closed curve on a surface $S$ is \emph{trivial} if it bounds a disc.

\begin{prop}
All simple closed curves in $S^2$ are trivial and isotopic.
\end{prop}
\begin{proof}
Every simple closed curve in $S^2$ is trivial by Jordan's Theorem (see Corollary \ref{Jordan:cor}). 

Theorem \ref{dischi:teo} and Lemma~\ref{diffeo:isotopi:lemma} together imply that there are at most two isotopy classes of trivial simple closed curves on any connected surface $S$, depending on their orientations. When $S=S^2$ the situation is quite special and the two classes reduce to one: a $\pi$-rotation along a horizontal axis transforms an eastward-run equator into a westward-run one.
\end{proof}

The simple closed curves in a torus are much more interesting.

\subsection{Simple closed curves in the torus}
We now classify the simple closed curves in the torus $T=S^1\times S^1$ up to isotopy. The fundamental group $\pi_1(T) = \matZ\times \matZ$ is abelian, hence a closed curve $\gamma$ is determined up to free homotopy by its class $(m,n)\in\matZ\times\matZ$.

\begin{prop} \label{mn:prop}
The class $(m,n)\neq (0,0)$ is represented by a simple closed curve if and only if 
$m$ and $n$ are coprime. In that case, the simple closed curve is unique up to isotopy.
\end{prop}
\begin{proof}
We visualise the torus as a quotient $T = \matR^2/_{\matZ^2}$. If $m$ and $n$ are coprime, the vector line generated by $(m,n)\in \matR^2$ projects to a simple closed curve in $T$ representing the class $(m,n)$, see Figure \ref{toro32:fig}. 

Conversely, let $\gamma\subset T$ be a simple closed curve. We cut $T$ along $\gamma$ and get a surface $S$ with $\chi(S) =0$. If $\gamma$ separates, then by the classification of surfaces the only possibility is that $S$ consists of a one-holed torus and a disc, so $\gamma$ is trivial and $(m,n)=(0,0)$, which is excluded. If $\gamma$ does not separate, there is another curve $\eta$ intersecting it transversely in one point. This implies that $\eta$ and $\gamma$ have algebraic intersection $\pm 1$ and hence $(m,n)$ are coprime: if $(m,n) = k(p,q)$ with $k\geqslant 2$ the algebraic intersection would be divided by $k$.

Finally, we pick two simple closed curves $\gamma, \eta$ of the same type $(m,n)\neq (0,0)$ and prove that they are isotopic. We put $\gamma, \eta$ in transverse position, and we cut $T$ along $\gamma$ to get an annulus $A$. If $\eta$ does not intersect $\gamma$, by cutting $A$ along $\eta$ we get two annuli: hence $\eta$ and $\gamma$ cobound an annulus and we can use this annulus to build an isotopy between them. If $\eta$ intersects $\gamma$, the curve $\eta$ decomposes into arcs in $A$ as in Figure \ref{torus_simple_curve:fig}-(left). Since the algebraic intersection of $\eta$ and $\gamma$ is zero, there is an arc with both endpoints in the same boundary component of $A$. This arc forms a bigon, which we can slide by isotopy as in Figure \ref{torus_simple_curve:fig}-(right) to decrease the intersection points in $\gamma \cap \eta$ and conclude by induction.
\end{proof}

\begin{figure}
\begin{center}
\includegraphics[width = 9 cm] {\iftoggle{BW}{toro32-BW}{toro32}}
\nota{If $(m,n)$ are coprime the line generated by $(m,n)\in\matR^2$ projects to a simple closed curve in $T$. Here $(m,n) = (3,2)$.}
\label{toro32:fig}
\end{center}
\end{figure}

\begin{figure}
\begin{center}
\includegraphics[width = 12.5 cm] {\iftoggle{BW}{torus_simple_curve-BW}{torus_simple_curve}}
\nota{Two homotopic simple closed curves in a torus are isotopic: put them in transverse position, cut the first to get an annulus, then remove bigons by isotopy to destroy intersections.}
\label{torus_simple_curve:fig}
\end{center}
\end{figure}

Simple closed curves in surfaces of higher genus are more complicated to classify, but are still very important. A couple of techniques used in the previous proof (cutting along curves and simplifying bigons) will be employed again in the higher-genus context.

\subsection{Preliminaries on simple curves}
Let $S_g$ be a closed orientable surface of some genus $g\geqslant 1$. We now prove some topological facts on simple closed curves in $S_g$, sometimes employing hyperbolic geometry.
\begin{prop} \label{finitely:prop}
There are finitely many simple closed curves in $S_{g}$ up to diffeomorphisms of $S_g$.
\end{prop}
\begin{proof}
By cutting $S_g$ along a simple closed curve $\gamma$ we get a surface $S'$ with the same Euler characteristic as $S_g$, with one or two components, and with the boundary oriented as $\gamma$: there are only finitely many diffeomorphism types for $S'$. 

Suppose that by cutting $S_g$ along $\gamma_1$ and $\gamma_2$ we get two surfaces $S'_1$ and $S'_2$ of the same type. By hypothesis there is a diffeomorphism $\varphi\colon S'_1 \to S'_2$ that preserves the boundary orientations. By Lemma \ref{diffeo:isotopi:lemma} we may modify $\varphi$ near the boundary so that it extends to a diffeomorphism $\varphi\colon S_g \to S_g$ sending $\gamma_1$ to $\gamma_2$. This concludes the proof.
\end{proof}

\begin{ex} \label{one:non-sep:ex}
There is precisely \emph{one} non-separating simple closed curve in $S_g$ up to diffeomorphism! 
\end{ex}

\begin{figure}
\begin{center}
\includegraphics[width = 7 cm] {\iftoggle{BW}{surface_non_sep-BW}{surface_non_sep}}
\nota{There are infinitely many non-separating simple closed curves in $S_g$ up to isotopy when $g\geqslant 1$, but there is only one up to self-diffeomorphism of $S_g$. For instance, some diffeomorphism of $S_2$ sends the \iftoggle{BW}{left}{blue} simple closed curve to the \iftoggle{BW}{right}{red} one (a self-diffeomorphism may be quite hard to visualise directly...).}
\label{surface_non_sep:fig}
\end{center}
\end{figure}

Every non-separating simple closed curve in $S_g$ may be transformed into your favourite one by some diffeomorphism of $S_g$, see Figure \ref{surface_non_sep:fig}. This quite useful fact, analogous to the possibility of changing a basis in a vector space, was called the \emph{change of coordinates principle} by Farb and Margalit.\index{change of coordinates principle} 

Recall that a non-trivial element $g\in G$ in a group is \emph{primitive} if it cannot be written as $g=h^n$ for some $n\geqslant 2$ and some $h\in G$. This condition is conjugacy-invariant, so the following makes sense. 

\begin{prop} \label{disco:prop}
Let $\gamma$ be a simple closed curve in $S_g$. It holds:
\begin{itemize}
\item if $\gamma$ is homotopically trivial, it bounds a disc;
\item if $\gamma$ is not homotopically trivial, it is primitive in $\pi_1(S_g)$.
\end{itemize}
\end{prop}
\begin{proof}
Let $S'$ be the surface obtained by cutting $S=S_g$ along $\gamma$. The surface $S'$ may have one or two components and has the same Euler characteristic of $S$. If one component of $S'$ is a disc, we are done. If $S'$ is an annulus, then $S$ is a torus and we are done by Proposition \ref{mn:prop}.

In all the other cases there is a hyperbolic metric on $S$ where $\gamma$ is a geodesic: each component of $S'$ has negative Euler characteristic and hence can be given a hyperbolic structure with geodesic boundary curves of length $1$; by gluing them we get the hyperbolic metric on $S$. 

We now apply Proposition \ref{unique:closed:geodesic:prop} multiple times. Since $\gamma$ is a closed geodesic, it is not homotopically trivial. If $\gamma = \eta^k$ is not primitive in $\pi_1(S_g)$, then $\eta$ is homotopic to a closed geodesic $\bar\eta$ and hence $\gamma$ is homotopic to $\bar\eta$ run $k$ times: a closed curve cannot be homotopic to two distinct closed geodesics, a contradiction.
\end{proof}

We have employed hyperbolic geometry to prove a topological fact on surfaces: this will be a refrain in this chapter. Let the \emph{inverse} $\gamma^*$ of a closed curve $\gamma$ be $\gamma$ run with opposite orientation.

\begin{prop} \label{no:inverso:prop}
A non-trivial simple closed curve in $S_g$ is never freely homotopic to its inverse. 
\end{prop}
\begin{proof}
If $g=1$ the curves $\gamma$ and $\gamma^*$ represent distinct elements (and hence conjugacy classes) in $\pi_1(S_1) =  \matZ \times \matZ$. If $g\geqslant 2$, give $S_g$ a hyperbolic metric. The curve $\gamma$ is homotopic to a closed geodesic $\bar\gamma$ and hence $\gamma^*$ is homotopic to its inverse $\bar\gamma^*$, which is distinct from $\bar\gamma$ as a closed geodesic by definition. Distinct closed geodesics are not homotopic.
\end{proof}

\subsection{Simple closed geodesics}
Let now $S_g$ have genus $g\geqslant 2$ and be equipped with a hyperbolic metric.
We know from Corollary \ref{unique:closed:geodesic:cor}
that every homotopically non-trivial closed curve in $S_g$ has a unique geodesic representative. We now prove that, if the original curve is simple, then the geodesic representative also is. 

Recall that the \emph{$R$-neighbourhood} of an object in a metric space is the set of all points of distance at most $R$ from that object. The $R$-neighbourhoods of lines in $\matH^2$ are particularly simple. 

\begin{figure}
\begin{center}
\includegraphics[width = 11.5 cm] {\iftoggle{BW}{intorni-BW}{intorni}}
\nota{The cones and bananas $R$-neighbourhoods of a geodesic $l$ in the half-plane (left) and disc (right) models.}
\label{intorni:fig}
\end{center}
\end{figure}

\begin{prop} \label{banana:prop}
The $R$-neighbourhood of a line $l\subset\matH^2$ in a conformal model is bounded by two Euclidean lines or circle arcs having the same endpoints as $l$ as in Figure \ref{intorni:fig}.
\end{prop}
\begin{proof}
Put $l$ in the half-space model with endpoints at $0$ and $\infty$. A $R$-neighbourhood is invariant by the isometry $x\mapsto \lambda x$ and is hence a cone as in the figure. The other cases follow because isometries and inversions send lines and circles to lines and circles. 
\end{proof}

We use $R$-neighbourhoods to prove the following.

\begin{prop} \label{semplice:prop}
Let $S_g$ be equipped with a hyperbolic metric. Every non-trivial simple closed curve in $S_g$ is homotopic to a simple closed geodesic.
\end{prop}
\begin{proof}
Every non-trivial simple closed curve $\gamma$ in $S_g=\matH^2/_\Gamma$ is homotopic to a closed geodesic $\bar\gamma$ by Corollary \ref{unique:closed:geodesic:cor}, and we now prove that $\bar\gamma$ is simple. The counterimage of $\gamma$ in $\matH^2$ consists of disjoint simple \emph{arcs}, while the counterimage of $\bar{\gamma}$ consists of straight \emph{lines}: we prove that these lines are also disjoint. 

The homotopy between $\gamma$ and $\bar{\gamma}$ lifts to a homotopy between the arcs and the lines. The homotopy between $\gamma$ and $\bar{\gamma}$ has compact support, hence there is a $R>0$ such that every point in the arcs is moved in the lifted homotopy to some distance smaller than $R$. Therefore every arc is contained in the $R$-neighbourhood of a line as in Figure \ref{geodetica:fig}-(left). 

\begin{figure}
\begin{center}
\includegraphics[width = 11.5 cm] {\iftoggle{BW}{geodetica-BW}{geodetica}}
\nota{The lifts of $\gamma$ (black arcs) and of its geodesic representative (\iftoggle{BW}{grey}{red} lines) have the same endpoints in $\partial \matH^2$ (left). If two lines intersect, the corresponding arcs do (right).}
\label{geodetica:fig}
\end{center}
\end{figure}

This shows that lines and arcs have the same endpoints in $\partial \matH^2$. If two lines intersects, their endpoints are linked in the circle $\partial \matH^2$ and hence also the corresponding arcs intersect, see Figure \ref{geodetica:fig}-(right): a contradiction.

Since the lifts of $\bar\gamma$ do not intersect, the geodesic $\bar\gamma$ does not self-intersect transversely. Then $\bar\gamma$ is either simple or wraps multiple times a simple geodesic, but the second possibility is excluded by Proposition \ref{disco:prop}.
\end{proof}

We will soon promote ``homotopic'' to ``isotopic''.

\subsection{Geometric intersection}
Isotopy classes of simple closed curves on $S_g$ form a complicate and interesting set, and a way to study it consists of looking at the way these curves intersect each other. The \emph{algebraic} intersection is too weak a tool, because it detects only their homology classes. We now introduce the much finer \emph{geometric} intersection.

Let $\gamma_1$ and $\gamma_2$ be two simple closed curves in an orientable surface $S$. Let the \emph{geometric intersection} $i(\gamma_1,\gamma_2)$ be the minimum number of intersections of two transverse simple closed curves $\gamma_1', \gamma_2'$ homotopic to $\gamma_1$ and $\gamma_2$. The geometric intersection depends only on the homotopy classes of $\gamma_1$ and $\gamma_2$.\index{geometric intersection of curves}

\begin{prop} We have $i(\gamma,\gamma)=0$ for every simple closed curve $\gamma$.
\end{prop}
\begin{proof}
A tubular neighbourhood of $\gamma$ is diffeomorphic to $S^1\times [-1,1]$ because $S$ is orientable, hence $\gamma$ has two disjoint parallel representatives $S^1\times \left\{-\frac {1}2\right\}$ and $S^1 \times \left\{\frac 12\right\}$ which do not intersect.
\end{proof}

Recall that the \emph{algebraic} intersection of two curves counts the intersections with sign. Geometric and algebraic intersections behave much differently and are equal only modulo $2$.

The geometric intersection $i(\gamma_1,\gamma_2)$ remains unaffected if we substitute $\gamma_1$ with its inverse, that is if we reverse its orientation. 

\subsection{The bigon criterion}
Two simple closed curves $\gamma_1$ and $\gamma_2$ in an orientable surface $S$ are in \emph{minimal position} if they intersect transversely in exactly $i(\gamma_1,\gamma_2)$ points. How can we know if two transverse simple closed curves in $S_g$ are in minimal position? We now prove a nice and simple criterion: they are in minimal position if and only if they do not form \emph{bigons}.

The complement of two transverse simple closed curves is a finite disjoint union of open sets with polygonal boundaries; one such set is a \emph{bigon} if it is a disc with two sides as in Figure \ref{bigono:fig}. 

\begin{figure}
\begin{center}
\includegraphics[width = 4 cm] {\iftoggle{BW}{bigono-BW}{bigono}}
\nota{Two curves $\gamma_1$ and $\gamma_2$ are in minimal position if and only if they do not create bigons, like this one.}
\label{bigono:fig}
\end{center}
\end{figure}

\begin{teo}[Bigon criterion] \label{bigono:teo}
Two transverse simple closed curves $\gamma_1,\gamma_2$ in $S_g$ are in minimal position if and only if they do not form bigons.
\end{teo}
\begin{proof}
If $\gamma_1$ and $\gamma_2$ create a bigon, the homotopies in Figure \ref{bigono_proof:fig}-(left) transform $\gamma_1$ and $\gamma_2$ into two curves that intersect in a smaller number of points: hence $\gamma_1$ and $\gamma_2$ are not in minimal position.

Suppose now that $\gamma_1$ and $\gamma_2$ do not form bigons: we must show that they are in minimal position. If $\gamma_1$ is trivial, it bounds a disc $D$ as in Figure \ref{bigono_proof:fig}-(right). If $\gamma_2$ intersects $\gamma_1$, an \emph{innermost} argument shows that they form a bigon: the curve $\gamma_2$ intersects $D$ in arcs, each dividing $D$ into two parts; if one part contains no other arc it is a bigon, otherwise we iterate. Therefore $\gamma_1$ and $\gamma_2$ are disjoint and hence in minimal position.\index{innermost argument}

\begin{figure}
\begin{center}
\includegraphics[width = 9 cm] {\iftoggle{BW}{bigono_proof-BW}{bigono_proof}}
\nota{A bigon can be eliminated via homotopies (left). If $\gamma_1$ bounds a disc and $\gamma_2$ intersects $\gamma_1$, there is a bigon (right).}
\label{bigono_proof:fig}
\end{center}
\end{figure}

It remains to consider the case where both $\gamma_1$ and $\gamma_2$ are non-trivial. 
The torus case $g=1$ is obtained by readapting the proof of Proposition \ref{mn:prop} and is left as an exercise, so we suppose $g\geqslant 2$.

Fix an arbitrary hyperbolic metric $S_g=\matH^2/_\Gamma$ and let $\pi\colon \matH^2 \to \matH^2/_\Gamma$ be the projection. The closed curves $\gamma_1$ and $\gamma_2$ are now homotopic to two simple closed geodesics $\bar{\gamma}_1$ and $\bar{\gamma}_2$. The lifts
of $\gamma_i$ and $\bar{\gamma}_i$ in $\matH^2$ are \emph{arcs} and \emph{lines} and there is a $R>0$ such that every arc lies in the $R$-neighbourhood of a line, see the proof of Proposition \ref{semplice:prop}. Arcs and lines have the same endpoints at infinity as in Figure \ref{geodetiche0:fig}-(left). 

\begin{figure}
\begin{center}
\includegraphics[width = 10 cm] {\iftoggle{BW}{geodetiche0-BW}{geodetiche0}}
\nota{The lifts of $\gamma_i$ and $\bar{\gamma}_i$ have distance bounded by $R$ and hence have the same endpoints (left). Two curves that intersect in more than one point form a bigon (right).} 
\label{geodetiche0:fig}
\end{center}
\end{figure}

Two distinct arcs may intersect at most in one point: if they intersect more, an innermost argument shows that they form a bigon $D$ as in Figure \ref{geodetiche0:fig}-(right), which projects to a bigon $\pi(D)$ in $S_g$. (It is not immediate that $\pi(D)$ is a bigon! The two vertices of $\pi(D)$ might coincide, but this is easily excluded because $S_g$ is orientable.) 

We now show how to count the intersections between $\gamma_1$ and $\gamma_2$ directly in $\matH^2$. Let $C(\gamma_i)\subset \Gamma$ be the conjugacy class of all the hyperbolic transformations corresponding to $\gamma_i$. We know that the lifts of $\bar{\gamma}_i$ are precisely the axis of the hyperbolic transformations in $C(\gamma_i)$.

By Corollary \ref{axis:cor} any two such axis are either incident or ultra-parallel. Hence two lifts of $\gamma_1$ and $\gamma_2$ intersect (in a single point) if and only if the corresponding lifts of $\bar{\gamma}_1$ and $\bar{\gamma}_2$ intersect (in a single point), and this holds if and only if the endpoints are linked in $\partial \matH^2$. We have established two bijective correspondences
$$\pi^{-1}(\gamma_1) \cap \pi^{-1}(\gamma_2) \longleftrightarrow \pi^{-1}(\bar{\gamma}_1) \cap \pi^{-1}(\bar{\gamma}_2) \longleftrightarrow X$$
with
$$X = \big\{(\varphi_1,\varphi_2) \in C(\gamma_1) \times C(\gamma_2) \ \big| 
\ \Fix(\varphi_1) \ {\rm and}\ \Fix(\varphi_2) \ {\rm are\ linked} \big\}.$$
The bijective correspondences are $\Gamma$-equivariant. We quotient by $\Gamma$ and find
$$\gamma_1 \cap \gamma_2 \longleftrightarrow \bar{\gamma}_1 \cap \bar{\gamma}_2 \longleftrightarrow X/_\Gamma.$$
The cardinality $k$ of $X/_\Gamma$ depends only on the homotopy type of $\gamma_1$ and $\gamma_2$. Therefore any two curves homotopic to $\gamma_1$ and $\gamma_2$ will have at least these $k$ intersections. Hence $\gamma_1$ and $\gamma_2$ are in minimal position.
\end{proof}

The bigon criterion furnishes an efficient algorithm to calculate the geometric intersection of any pair $\gamma_1$ and $\gamma_2$ of simple closed curves in $S_g$: we put them in transverse position, and then we simplify bigons as much as we can. After finitely many steps we get two curves in minimal position.

Geodesic representatives are typically efficient: they minimise lengths (see Proposition  \ref{minimum:cor}), and they also minimise mutual intersections.

\begin{cor} \label{geodetiche:cor}
Let $g\geqslant 2$ and $S_g$ have a hyperbolic metric. Two simple closed geodesics with distinct supports are always in minimal position.
\end{cor}
\begin{proof}
Two geodesics do not form bigons: if they do, the bigon lifts to a bigon between two lines in $\matH^2$, but lines intersect at most once.
\end{proof}

Here is a very simple application of the bigon criterion.

\begin{cor} If two simple closed curves $\gamma, \eta$ intersect transversely in one point, we have $i(\gamma,\eta)=1$. In particular, they are both non-trivial.
\end{cor}
\begin{proof}
The curves $\gamma$ and $\eta$ cannot form bigons.
\end{proof}

The following exercise shows that the geometric intersection distinguishes a trivial curve from a non-trivial one.

\begin{ex} \label{eta:prop}
If a simple closed curve $\gamma$ is not trivial, there is another simple closed curve $\eta$ such that $i(\gamma, \eta)>0$.
\end{ex}
\begin{proof}[Hint] Use Proposition \ref{finitely:prop} to transform $\gamma$ into a comfortable curve and draw a $\eta$ which intersects $\gamma$ in at most $2$ points without bigons.
\end{proof}

On the contrary, the algebraic intersection does \emph{not} distinguish the trivial curve from any other separating curve.

\begin{ex} \label{i:torus:ex}
Let $\gamma$ and $\eta$ be non-trivial simple closed curves in the torus $T$ of type $(p,q)$ and $(r,s)$. We have
$$i(\gamma, \eta \big) = \left|\det \begin{pmatrix} p & r \\ q & s \end{pmatrix}\right|.$$
\end{ex}

\subsection{Homotopy and isotopy of curves}
We are now going to prove some ``homotopy implies isotopy'' theorems: we start with closed curves and we end in the next section with self-diffeomorphisms of surfaces. These slightly technical theorems have the remarkable and pleasant effect of making life easier in dimensions two and three (their impact to three-dimensional topology will be clear in the next chapter): topologists use these theorems everyday when they manipulate two- and three-manifolds.

We start by showing that two non-trivial simple closed curves are homotopic if and only if they are isotopic. We first consider a particular case.

\begin{lemma} \label{parallele:lemma}
Let $\gamma_1$ and $\gamma_2$ be two non-trivial simple closed curves in $S_g$. If they are disjoint and homotopic, they are parallel. 
\end{lemma}
\begin{proof}
Cut $S_g$ along $\gamma_1\cup \gamma_2$. We do not obtain discs because the curves are non-trivial and if we obtain an annulus the two curves are parallel. In all the other cases we obtain surfaces of negative Euler characteristic, and we may assign some hyperbolic metrics to them that glue to a hyperbolic metric on $S_g$ where both $\gamma_1$ and $\gamma_2$ are distinct geodesics: a contradiction.
\end{proof}

We now turn to the general case. 

\begin{prop}[Homotopy implies isotopy] \label{ambiente:prop}
Two non-trivial simple closed curves in $S_g$ are homotopic if and only if they are isotopic. 
\end{prop}
\begin{proof}
Let $\gamma_1$ and $\gamma_2$ be two non-trivial simple closed curves. With a small isotopy we put them in transverse position. 
Since $i(\gamma_1,\gamma_2) = i(\gamma_1,\gamma_1)=0$ the two curves are either disjoint or form some bigon. If they form a bigon, we can eliminate it via isotopies as in Figure \ref{bigono_proof:fig}-(left) and after finitely many steps we get two disjoint curves.

The curves $\gamma_1$ and $\gamma_2$ are parallel by Lemma \ref{parallele:lemma}, and we use the annulus they cobound to move $\gamma_2$ isotopically over $\gamma_1$. The two curves now have the same support and the same orientation by Proposition \ref{no:inverso:prop}: by Lemma \ref{diffeo:isotopi:lemma} they are isotopic. 
\end{proof}

\begin{warning}
Two homotopic simple closed curves in a \emph{three}-manifold are not necessarily isotopic, because they may be knotted differently: the \emph{knot theory} studies precisely this phenomenon.
\end{warning}

\begin{cor}
Let $g\geqslant 2$ and $S_g$ have a hyperbolic metric. Every non-trivial simple closed curve in $S_g$ is isotopic to a simple closed geodesic.
\end{cor}

\subsection{Multicurves}
Our next goal is to prove a ``homotopy implies isotopy'' theorem for self-diffeomorphisms of surfaces. The core of the proof is contained in the next section, and we will need there the notion of \emph{multicurve}: we now introduce this simple concept and extend some of the previous results from curves to multicurves.\index{multicurve} 

\begin{defn} A \emph{multicurve} $\mu$ in $S_g$ is a finite set of disjoint non-trivial simple closed curves.
\end{defn}
See an example in Figure \ref{multicurva:fig}, and note that every component of a multicurve is oriented. A multicurve is \emph{essential} if it has no parallel components. By cutting $S_g$ along an essential multicurve $\mu$ we get finitely many surfaces of negative Euler characteristic: if each such surface is a pair-of-pants, then $\mu$ is a \emph{pants decomposition}, already considered in Section \ref{hyperbolic:surfaces:subsection}.\index{pants decomposition}

\begin{prop} \label{parallele:prop}
An essential multicurve $\mu$ in $S_g$ with $g\geqslant 2$ has at most $3g-3$ components, and it has $3g-3$ if and only if it is a pants decomposition.
\end{prop}
\begin{proof}
By cutting $S_g$ along $\mu$ we get some surfaces $S^1,\ldots, S^k$ of negative Euler characteristic such that $\chi(S_g) = \chi(S^1)+\ldots +\chi(S^k)$. If each $S^i$ is a pair-of-pants then $\chi(S^i)=-1$ and $k=-\chi(S_g) = 2g-2$; the curves are $\frac 32(2g-2) = 3g-3$ because each curve separates two pants. If some $S^i$ is not a pair-of-pants it can be further subdivided into pair-of-pants.
\end{proof}

\begin{figure}
\begin{center}
\includegraphics[width = 7 cm] {\iftoggle{BW}{multicurva-BW}{multicurva}}
\nota{A multicurve in a surface of genus two.} 
\label{multicurva:fig}
\end{center}
\end{figure}

We define the \emph{geometric intersection} $i(\mu_1,\mu_2)$ of two multicurves $\mu_1$ and $\mu_2$ as the minimum number of intersections of two transverse multicurves $\mu_1'$, $\mu_2'$ isotopic to $\mu_1, \mu_2$. This definition extends the geometric intersection of simple closed curves (the original definition for simple closed curves is with ``homotopic'' instead of ``isotopic'', but these are equivalent by Proposition \ref{ambiente:prop}).

Two transverse multicurves $\mu_1$ and $\mu_2$ are in \emph{minimal position} if they intersect exactly in $i(\mu_1, \mu_2)$ points: the bigon criterion easily extends to this context. 

\begin{prop} \label{multicurve:sum:prop}
Let $\mu_1, \mu_2 \subset S_g$ be two transverse multicurves. The following equality holds:
$$i(\mu_1, \mu_2) = 
\sum_{\tiny{\begin{array}{c} \gamma_1 \subset \mu_1 \\ \gamma_2 \subset \mu_2 \end{array}} }  
i(\gamma_1, \gamma_2)$$
where the sum is taken on all components $\gamma_1$, $\gamma_2$ of $\mu_1$, $\mu_2$. The multicurves $\mu_1$ and $\mu_2$ are in minimal position if and only if they do not form bigons.
\end{prop}
\begin{proof}
If $\mu_1$ and $\mu_2$ form no bigons, then any two components $\gamma_1$ and $\gamma_2$ form no bigons too (exercise) and are therefore in minimal position. This proves the equality and that $\mu_1$ and $\mu_2$ are in minimal position.
\end{proof}

The formula says that if we consider a multicurve as a ``sum'' of disjoint simple closed curves, then the geometric intersection $i$ behaves like a bilinear form: this viewpoint will be further explored in Chapter \ref{automorfismi:chapter}.

Note again that $i(\mu,\mu)=0$. We extend Proposition \ref{ambiente:prop} to essential multicurves.

\begin{prop}[Homotopy implies isotopy] \label{ambiente:multicurve:prop}
Let
$$\mu_1 = \{\gamma_{1,1}, \ldots, \gamma_{1,n}\}, \quad \mu_2 = \{\gamma_{2,1}, \ldots, \gamma_{2,n}\}$$
be two essential multicurves in $S_g$.
If $\gamma_{1,i}$ is homotopic to $\gamma_{2,i}$ for all $i$ then there is an isotopy moving $\mu_1$ to $\mu_2$. 
\end{prop}
\begin{proof}
We adapt the proof of Proposition \ref{ambiente:prop}. Since
$i(\gamma_{1,i}, \gamma_{2,j}) = i(\gamma_{1,i},\gamma_{1,j}) = 0$ we get $i(\mu_1,\mu_2)=0$ by Proposition \ref{multicurve:sum:prop} and after an isotopy that destroys the bigons we get $\mu_1\cap\mu_2 = \emptyset$. Then Lemma \ref{parallele:lemma} implies that $\gamma_{1,i}$ and $\gamma_{2,i}$ are parallel and we are done.
\end{proof}
 
\begin{cor} \label{isotopia:geodetiche:cor}
Let $g\geqslant 2$ and $S_g$ have a hyperbolic metric. Every essential multicurve can be isotoped to a (unique) geodesic essential multicurve. 
\end{cor}

In particular every pants decomposition straightens to a (unique) geodesic one: this fact will be used in Chapter \ref{Teichmuller:chapter} to parametrize all the hyperbolic metrics on a given $S_g$.

\subsection{Minimal position}
We end this discussion on multicurves by showing that the minimal position is unique up  to isotopy.

\begin{figure}
\begin{center}
\includegraphics[width = 8 cm] {\iftoggle{BW}{minimale_unica-BW}{minimale_unica}}
\nota{A bigon between $\eta$ and $\eta'$ intersects $\mu$ in vertical arcs and can be removed via an ambient isotopy that preserves the support of $\mu$.}
\label{minimale_unica:fig}
\end{center}
\end{figure}

\begin{prop} \label{minimale:unica:prop}
Let $\mu, \eta$ be two multicurves in minimal position in $S_g$. The union $\mu\cup \eta$ of their supports depends up to ambient isotopy only on the isotopy classes of $\mu$ and $\eta$.
\end{prop}
\begin{proof}
Let $\mu', \eta'$ be multicurves in minimal position, individually isotopic to $\mu, \eta$. We need to prove that the supports $\mu'\cup \eta'$ and $\mu\cup \eta$ are ambiently isotopic.

By hypothesis there is an isotopy carrying $\mu'$ to $\mu$, which is ambient since multicurves are compact: hence we can suppose $\mu' = \mu$. We now construct an isotopy that fixes $\mu'=\mu$ as a set (not pointwise!) and carries $\eta'$ to $\eta$. 

Up to a small ambient isotopy fixing the set $\mu$ we may suppose that $\eta$ and $\eta'$ intersect transversely. If $\eta\cap \eta' \neq \emptyset$ then $\eta$ and $\eta'$ form a bigon as in Figure \ref{minimale_unica:fig}-(left): the multicurve $\mu$ intersects the bigon in arcs that join distinct edges as in the figure (because $\mu' = \mu$ forms no bigons with $\eta$ and $\eta'$). We can eliminate the bigon by an ambient isotopy that fixes $\mu$
as a set, as shown in Figure \ref{minimale_unica:fig}. 

After finitely many steps we get $\eta\cap \eta' = \emptyset$. Since $\eta$ and $\eta'$ are disjoint and isotopic, every component of $\eta$ is parallel to a component of $\eta'$. A maximal set of parallel curves in $\eta \cup \eta'$ is contained in a bigger annulus $[-1,1] \times S^1$ which intersects $\mu$ either into arcs $[-1,1] \times \{{\rm pt}\}$ or into circles $\{{\rm pt}\} \times S^1$. In both cases we easily see that $\mu\cup\eta$ is ambiently isotopic to $\mu\cup \eta'$.
\end{proof}

As an example, consider two essential multicurves $\eta$ and $\mu$ in $S_g$ with $g\geqslant 2$ (for instance, two non-trivial simple closed curves). A hyperbolic metric on $S_g$ produces two geodesic representatives $\bar{\eta}$ and $\bar{\mu}$ for them, and the following holds:
\begin{cor} \label{unica:configurazione:cor}
The support $\bar{\eta}\cup \bar{\mu}$ in $S_g$ does not depend (up to ambient isotopy) on the chosen hyperbolic metric.
\end{cor}
\begin{proof}
With the exception of the common components of $\bar{\eta}$ and $\bar{\mu}$, the rest intersects in minimal position. 
\end{proof}

This corollary is one of the many instances of the following nice phenomenon: a hyperbolic metric on $S_g$ may be used as an auxiliary structure to define objects or to prove statements which turn out \emph{a posteriori} not to depend on the chosen metric. We will employ this strategy many times in Chapter \ref{automorfismi:chapter}.

\section{Homotopy and isotopy}
We have proved that ``homotopy implies isotopy'' for non-trivial simple closed curves in $S_g$, and now we want to prove an analogous result for diffeomorphisms of $S_g$.
A crucial ingredient is the following theorem proved by Smale in the 1950s, which considers the self-diffeomorphisms of the disc $D^2$.

\begin{teo} \label{Smale:teo}
Two diffeomorphisms $\varphi, \psi\colon D^2 \to D^2$ that coincide on $\partial D^2$ are isotopic, via an isotopy that fixes $\partial D^2$ pointwise.
\end{teo}

There is a proof of Theorem \ref{Smale:teo} in the topological category that is surprisingly easy and works in any dimension. 

\subsection{The Alexander trick} \label{Alexander:subsection}
The proof of the following fact is so immediate, that it has a name: it is usually called the \emph{Alexander trick}.\index{Alexander trick}

\begin{prop}[Alexander trick] Two homeomorphisms $\varphi$, $\psi\colon D^n \to D^n$ that coincide on $\partial D^n$ are continuously isotopic, via an isotopy that fixes $\partial D^n$ pointwise.
\end{prop}
\begin{proof}
We consider $f=\varphi\circ \psi^{-1}$ and $\id_{D^n}$ and define an isotopy that transforms $f$ into $\id_{D^n}$ fixing $\partial D^n$. The following function does the job:
$$F(x,t) = \left\{ \begin{matrix} x & {\rm if \ } \|x\| \geqslant t, \\ t f\left(\frac xt \right) & {\rm if\  } \|x\| \leqslant t. \end{matrix} \right.$$
The proof is complete.
\end{proof}

Unfortunately this proof does not work in the smooth category, because the function $F$ is not smooth (and not smoothable) at $(0,0)$. Of course, as every continuous map, our $F$ may be approximated by smooth functions, but the injectivity of the slices $F(\cdot, t)$ can be lost in the approximation.  

Moreover, the statement in the smooth category is just false in dimension $n\geqslant 7$: there are self-diffeomorphisms of $D^n$ that are \emph{not} smoothly isotopic to $\id_{D^n}$ in high dimension, and this is connected to the existence of \emph{exotic spheres}, differentiable manifolds that are homeomorphic but not diffeomorphic to $S^n$.\index{exotic sphere} 

\subsection{Self-diffeomorphisms of the disc} \label{self:disc:subsection}
After this very brief excursion in the topological world, we turn back to the smooth category: every map considered so far is smooth by assumption like in the rest of the book. 

We prove here Theorem \ref{Smale:teo}. We need a fact on planar fields. 

\begin{prop} \label{PB:prop}
On a nowhere-vanishing vector field on $\matR^2$, a maximal integral curve is never trapped in a compact set.
\end{prop}
\begin{proof}
Suppose by contradiction that a maximal integral curve $\alpha\colon I \to \matR^2$ lies in a compact set $K$ for all $t>t_0$. Its \emph{$\omega$-limit} $\omega(\alpha)\subset K$ is the set of all points $x\in K$ to which $\gamma$ accumulates, that is the intersection of the closures of $\alpha(t_0,\infty)$ for all $t_0 \in I$. The $\omega$-limit is the intersection of a filtration of compact sets and is hence non-empty. Pick a point $p\in\omega(\alpha)$.

Up to a local diffeomorphism, the vector field near $p$ is constant vertical as in Figure \ref{PB:fig}-(1). The integral curves near $p$ are vertical. Since $p\in\omega(\alpha)$, the curve $\alpha$ contains infinitely many of them that tend to the one containing $p$. Two subsequent ones determine a closed curve as in Figure \ref{PB:fig}-(2), that bounds a disc by the Jordan curve theorem. 

We prove that such a disc cannot exist. If it existed, we could reverse all arrows and rotate it to obtain another disc as in Figure \ref{PB:fig}-(3), and by gluing the two discs we would construct as in Figure \ref{PB:fig}-(4) a nowhere-vanishing vector field on the sphere. Such a field of course does not exist since the sphere has non-zero Euler characteristic.
\end{proof}

\begin{figure}
\begin{center}
\includegraphics[width = 11 cm] {\iftoggle{BW}{PB-BW}{PB}}
\nota{At a point $p$, after a local diffeomorphism we may suppose that the vector field is constantly vertical (1). If $p$ is a limit point of an integral curve $\alpha$ we may find a closed curve that bounds a disc (2). If such a disc existed, we could rotate it and inverse all arrows (3) and glue the two portions to form a nowhere-vanishing vector field on the sphere (4). }
\label{PB:fig}
\end{center}
\end{figure}

We can now prove Theorem \ref{Smale:teo}.

\begin{figure}
\begin{center}
\includegraphics[width = 12.5 cm] {\iftoggle{BW}{trajectories-BW}{trajectories}}
\nota{A self-diffeomorphism $\varphi$ of $D^2$ fixing $\partial D^2$ pointwise (left). To get an isotopy fixing $\partial D^2$ pointwise we need to rescale via a diffeomorphism that maps each arc as in the right figure to itself, and varies smoothly from $\id$ to $f_t^{-1}$ from bottom to top (right).}
\label{trajectories:fig}
\end{center}
\end{figure}

\begin{proof}[Proof of Theorem \ref{Smale:teo}]
After composing with $\psi^{-1}$ we can suppose that $\psi=\id$. Now $\varphi|_{S^1}$ is the identity. Every self-diffeomorphism of any compact manifold that fixes the boundary is isotopic to one which is the identity on a collar of the boundary: hence we can suppose that $\varphi=\id$ on a collar of $S^1$.

We sketch $\varphi$ in Figure \ref{trajectories:fig}-(left): vertical lines are transformed into paths with the same endpoints. We want to construct an isotopy between $\varphi$ and $\id_{D^2}$ that ``straightens'' smoothly these lines, and we start by straightening their tangent vector fields.

Consider the constant unitary vertical vector field $X_0=(0,1)$ on $D^2$ and let $X_1$ be its image along the diffeomorphism $\varphi$. It is easy to construct a  homotopy $X_t$ between $X_0$ and $X_1$ through non-vanishing vector fields: we see both $X_0$ and $X_1$ as maps $D^2 \to \matC^*$, we lift them to maps $D^2\to \matC$ via the universal cover $\exp\colon \matC \to \matC^*$, we make a convex combination and we project it back to $\matC^*$.
By construction the vectors $X_t$ are constantly $(0,1)$ on a collar of $S^1$.

We now carefully integrate the homotopy $X_t$ of vector fields to an isotopy of diffeomorphisms. Let $S^1_+$ and $S^1_-$ be the upper and lower hemisphere of $S^1$. The vectors $X_t$ at $S^1_-$ point inside and those at $S^1_+$ point outside $D^2$. By Proposition \ref{PB:prop} the integral curve $\gamma_x^t$ of $X_t$ starting from a point $(x,-\sqrt{1-x^2}) \in S^1_-$ exits at some point in $S^1_+$ after some time $T_x^t$ that varies smoothly in $x$ and $t$.

The point $\gamma_x^t(u)$ is defined for $u\in [0, T_x^t]$ and varies smoothly on $x,t,u$. We obviously have $T_x^0= 2\sqrt{1-x^2}$ and also $T_x^1 = T_x^0$ since $\gamma_x^1 = \varphi\circ\gamma_x^0$. However $T_x^t$ may not be constant in $t$, so we rescale the integral curves as
$$\tilde\gamma_x^t(u) = \gamma_x^t\left( u \cdot \frac{T_x^t}{2\sqrt{1-x^2}}\right).$$
The new curve $\tilde\gamma_x^t$ is defined in the interval $\big[0, 2\sqrt{1-x^2}\big]$ not depending on $t$. We can now define
$$ \varphi_t \big(x, -\sqrt{1-x^2} + u \big) = \tilde\gamma_x^t\left(u\right).
$$
The diffeomorphism $\varphi_t\colon D^2 \to D^2$ furnishes an isotopy between $\id = \varphi_0$ and $\varphi=\varphi_1$. However, the diffeomorphism $\varphi_t$ does not fix the points in $S^1_+$ since the endpoint of $\tilde\gamma_x^t$ has a first coordinate $f_t(x)$ which might be distinct from $x$ when $0<t<1$ (but varies smoothly on $t$ and $x$). 

To fix that it suffices to compose each $\varphi_t$ with a diffeomorphism of $D^2$ that maps each arc as in Figure \ref{trajectories:fig}-(right) to itself and varies smoothly from $\id$ to $f_t^{-1}$ from bottom to top.
\end{proof}

\begin{oss} \label{Smale:collar:rem}
If $\varphi\colon D^2 \to D^2$ is the identity on a collar of $S^1$, we can easily transform the isotopy $\varphi_t$ so that it is the identity on this collar for all $t$. This will be useful to glue isotopies.
\end{oss}

\subsection{Homotopy and isotopy of diffeomorphisms}
We are now ready to promote homotopies to isotopies of diffeomorphisms. We start with the easier sphere case, which is a simple consequence of Theorem \ref{Smale:teo}.

\begin{teo} \label{S2:omotopi:teo}
Two diffeomorphisms $\varphi, \psi\colon S^2 \to S^2$ are isotopic $\Longleftrightarrow$ they are homotopic $\Longleftrightarrow$ they are co-oriented.
\end{teo}
\begin{proof}
We suppose that $\varphi$ and $\psi$ are co-oriented and we must prove that they are isotopic. Pick any disc $D\subset S^2$. Both $\varphi$ and $\psi$ send $D$ to some disc in $S^2$. Theorem \ref{dischi:teo} furnishes an ambient isotopy relating these maps, so we may suppose that $\varphi$ and $\psi$ coincide on $D$. The closed complement is another disc and we conclude by Theorem \ref{Smale:teo}.

More precisely, we pick a smaller closed disc $D' \subset \interior D$ and consider the closed complement $D'' = S^2\setminus \interior {D'}$. The diffeomorphisms $\varphi$ and $\psi$ of $D''$ coincide on a collar of $\partial D''$ and Theorem \ref{Smale:teo} together with Remark \ref{Smale:collar:rem} furnishes an isotopy that is the identity on that collar. We extend this isotopy as the identity on $D\setminus D''$ and we are done.
\end{proof}

The theorem has a very important three-dimensional application.\index{Smale theorem}

\begin{cor} [Smale's Theorem]
Every self-diffeomorphism of $S^2$ extends to a self-diffeomorphism of $D^3$.
\end{cor}
\begin{proof}
Every self-diffeomorphism of $S^2$ is isotopic either to the identity or to a reflection by Theorem \ref{S2:omotopi:teo}. We map a collar of $S^2$ to itself using this isotopy, then extend to the rest of $D^3$ using the identity or the reflection.
\end{proof}

We can finally promote homotopies to isotopies on a general genus-$g$ surface $S_g$. The proof uses many of the techniques introduced in the last sections: we employ multicurves to cut the surface into discs, and then use Theorem \ref{Smale:teo} individually on each disc.

\begin{figure}
\begin{center}
\includegraphics[width = 10cm] {\iftoggle{BW}{decomposizione2-BW}{decomposizione}}
\nota{Two essential multicurves (\iftoggle{BW}{dark and light grey}{blue and red}) in minimal position, which subdivide the surface into hexagons.}
\label{multicurva_decomposizione:fig}
\end{center}
\end{figure}

\begin{teo} \label{omotopi:isotopi:teo}
Two diffeomorphisms $\varphi$, $\psi\colon S_g \to S_g$ are isotopic if and only if they are homotopic.
\end{teo}
\begin{proof}
We have seen the case $g=0$, so we suppose $g\geqslant 1$. By composing with $\psi^{-1}$ we can suppose that $\psi = \id$.

Fix two essential multicurves $\mu$ and $\eta$ as in Figure \ref{multicurva_decomposizione:fig}. The figure easily generalises to any genus $g$, and the reader may check that $\mu$ and $\eta$ subdivide the surface into hexagons (or a single square when $g=1$). The important point is that $S_g \setminus (\mu \cup \eta)$ consists of polygons, none of which is a bigon.

By hypothesis $\varphi$ is homotopic to the identity, so the image multicurves $\varphi(\mu)$ and $\varphi(\eta)$ are curve by curve homotopic to $\mu$ and $\eta$. By Proposition \ref{ambiente:multicurve:prop} the multicurve $\varphi(\mu)$ is isotopic to $\mu$, and $\varphi(\eta)$ is isotopic to $\eta$.

The curves $\mu$ and $\eta$ are in minimal position because there are no bigons. The images $\varphi(\mu)$ and $\varphi(\eta)$ are also in minimal position because $\varphi$ is a diffeomorphism. 
By Proposition \ref{minimale:unica:prop} the supports $\mu \cup \eta$ and $\varphi(\mu\cup \eta)$ are ambiently isotopic, so we may suppose that they coincide.

The graph $\mu\cup \eta$ is made of vertices (the intersections $\mu\cap \eta$) and edges. The components of $\mu$ and $\eta$ are pairwise non-homotopic simple closed curves, hence $\varphi$ sends necessarily every component to itself, and it does so orientation-preservingly by Proposition \ref{no:inverso:prop}. This implies easily that vertices and edges are sent to themselves by $\varphi$. Hence $\varphi=\id$ on vertices and after an isotopy we may suppose that $\varphi=\id$ on edges too. 

After an isotopy we may also suppose that $\varphi = \id$ on a regular neighbourhood $U$ of $\mu\cup\eta$, obtained by thickening the cellularisation of $\mu\cup\eta$ to a handle decomposition with 0- and 1-handles. The complement of $U$ consists of discs (the hexagons). Consider one such disc $D$, enlarged a bit so that $\partial D \subset \interior U$. The map $\varphi$ sends $D$ to itself and is the identity on a collar of $\partial D$. By Theorem \ref{Smale:teo} and Remark \ref{Smale:collar:rem} there is an isotopy connecting $\varphi$ to $\id$ on every such disc $D$ that fixes pointwise this collar, so we can extend it constantly on the rest of $U$ and get a global isotopy on $S_g$ connecting $\varphi$ and $\id$.
\end{proof}

This theorem has important consequences in dimensions two and three.

\section{Mapping class group} \label{MCG:section}
We have just proved that two self-diffeomorphisms of $S_g$ are homotopic if and only if they are isotopic. We have also seen that when $g=0$ there are only two self-diffeomorphisms up to isotopy: the identity and a reflection. Do we get a more complicated picture when $g\geqslant 1$? Yes, we get an interesting group, called the \emph{mapping class group} of $S_g$. The group is naturally defined on all finite-type surfaces $S_{g,b,p}$.\index{group!mapping class group}

\subsection{Definition}
Recall that the finite-type orientable surface $S_{g,b,p}$ has genus $g$, it has $b$ boundary components, and $p$ punctures. 
\begin{defn}
The \emph{mapping class group} of $S_{g,b,p}$ is the group
$$\MCG(S_{g,b,p}) = \Diffeo^+(S_{g,b,p})/_\sim$$
where $\Diffeo^+(S_{g,b,p})$ indicates the group of all orientation-preserving diffeomorphisms $S_{g,b,p}\to S_{g,b,p}$ that fix pointwise the boundary and $\varphi\sim\psi$ if $\varphi$ and $\psi$ are connected by an isotopy that fixes the boundary pointwise at every level.
\end{defn}

\begin{example}
The groups $\MCG(S^2)$ and $\MCG(D^2)$ are trivial by Theorems \ref{Smale:teo} and \ref{S2:omotopi:teo}.
\end{example}

The group $\MCG(S_{g,b,p})$ acts on $H_1(S_{g,b,p}, \matZ)$ since homotopic functions induce the same maps in homology. We get a group homomorphism 
$$\MCG(S_{g,b,p}) \longrightarrow \Aut^+\big(H_1(S_{g,b,p}, \matZ)\big) = \Aut^+ \big(\matZ^{n}\big) = \SL_{n}(\matZ)$$
with $n=2g+ \max\{b+p-1,0\}$. This homomorphism is neither injective nor surjective in general. Its kernel is called the \emph{Torelli group} of $S_{g,b,p}$.\index{group!Torelli group} 

The mapping class group of a sphere is trivial, so we look at the torus.

\subsection{The torus}
The mapping class group of the torus is a familiar group of $2\times 2$ matrices.

\begin{prop} \label{Torelli:toro:prop}
The Torelli group of the torus $T$ is trivial and
$$\MCG(T) \Isom \Aut^+(H_1(T)) = \SLZ.$$
\end{prop}
\begin{proof}
Fix a meridian $m$ and longitude $l$ of $T$. A diffeomorphism $\varphi$ of $T$ that acts trivially on $H_1(T) = \pi_1(T) = \matZ^2$ sends $m$ and $l$ to two simple closed curves $\varphi(m)$ and $\varphi(l)$ homotopic and hence isotopic to $m$ and $l$: the proof of Theorem \ref{omotopi:isotopi:teo} shows that $\varphi$ is isotopic to the identity. Therefore the Torelli group is trivial.

A matrix $A\in\SLZ$ acts linearly on $\matR^2$ preserving the orientation and the lattice $\matZ^2$ and hence descends to a self-diffeomorphism of $T=\matR^2/_{\matZ^2}$. This shows that the map $\MCG(T) \to \Aut^+(H_1(T)) = \SL_2(\matZ)$ is surjective.
\end{proof}

The mapping class group of a surface $S_g$ of genus $g\geqslant 2$ is not isomorphic to a familiar group of matrices, at least as far as we know: indeed (except some low-genus cases) it is still unknown whether $\MCG(S_{g,b,p})$ is \emph{linear}, \emph{i.e.}~isomorphic to a subgroup of $\GL(n,\matC)$ for some integer $n$.

For the moment we content ourselves with a concrete description of some particularly simple elements of $\MCG(S_{g,b,p})$ called \emph{Dehn twists}.

\subsection{Dehn twists} \label{Dehn:twist:subsection}
Let $\gamma$ be a non-trivial simple closed curve in the interior of $S_{g,b,p}$. The \emph{Dehn twist} along $\gamma$ is the element $T_\gamma \in \MCG(S_{g,b,p})$ defined as follows.\index{Dehn twist}

\begin{figure}
\begin{center}
\includegraphics[width = 12.5 cm] {\iftoggle{BW}{Dehn_twist-BW}{Dehn_twist}}
\nota{A Dehn twist along a curve $\gamma$ maps a transverse arc $\mu$ onto an arc which makes a complete left turn.}
\label{Dehn_twist:fig}
\end{center}
\end{figure}

Pick a tubular neighbourhood of $\gamma$ orientation-preservingly diffeomorphic to $S^1\times [-1,1]$ where $\gamma$ lies as $S^1\times\{0\}$. Let $f\colon[-1, 1] \to \matR$ be a smooth function which is zero in $\left[-1, -\frac 1 2\right]$ and $2\pi$ on $\left[\frac 1 2,1\right]$. Let 
$$T_\gamma\colon S_{g,b,p} \longrightarrow S_{g,b,p}$$ 
be the diffeomorphism that acts on the tubular neighbourhood as $T_\gamma(e^{i\alpha},t) = (e^{i(\alpha+f(t))},t)$ and on its complementary set in $S_{g,b,p}$ as the identity. We may visualise $T_\gamma$ by noting that it gives a complete left turn to any arc $\mu$ that intersects $\gamma$ as in Figure \ref{Dehn_twist:fig}. 

\begin{prop} The element $T_\gamma \in \MCG(S_{g,b,p})$ is well-defined and depends only on the isotopy class of $\gamma$.
\end{prop}
\begin{proof}
In the definition of $T_\gamma$ we have chosen a tubular neighbourhood for $\gamma$ and a smooth function $f$. Tubular neighbourhoods are ambiently isotopic, and functions with fixed extremes are isotopic too: these facts imply easily that the isotopy class of $T_\gamma$ is well-defined and depends only on the isotopy class of $\gamma$.
\end{proof}

\begin{oss} To define $T_\gamma$ we needed the orientation of $S_{g,b,p}$, but not an orientation for $\gamma$. A bit surprisingly, if we change the orientation of $\gamma$ the element $T_\gamma$ remains unaffected. 

The inverse $T_\gamma^{-1}$ transforms every curve $\mu$ crossing $\gamma$ via a complete \emph{right}-turn and is sometimes called a \emph{negative Dehn twist}.
\end{oss}

We construct some examples in the torus $T$.
Let $m$ and $l$ be some fixed meridian and longitude on the oriented $T$ forming a positive basis, so that $m\cdot l = +1$ and we get an identification $\MCG(T) = \SLZ$.

\begin{prop}
The Dehn twists $T_m$ and $T_l$ are
$$T_m=\begin{pmatrix} 1 & {-1} \\ 0 & 1 \end{pmatrix}, \qquad T_l = \begin{pmatrix} 1 & 0 \\ 1 & 1 \end{pmatrix}.$$
\end{prop}
\begin{proof}
In homology we find
$$T_m(m) = m, \quad T_m(l) = l-m, \qquad T_l(l) = l, \quad T_l(m) =m+l.$$
The proof is complete.
\end{proof}

By Exercise \ref{one:non-sep:ex} there is a unique non-separating simple closed curve $\gamma$ in $S_g$ up to self-diffeomorphism for every $g\geq 1$. Therefore the Dehn twists $T_\gamma$ along non-separating curves in $S_g$ are all conjugate in $\MCG(S_g)$. In the torus case these can be easily identified algebraically.

\begin{cor}
An element $A\in \SL_2(\matZ) = \MCG(T)$ is a (positive or negative) Dehn twist $\Longleftrightarrow$ it is primitive and $\tr A = 2$.
\end{cor}
\begin{proof}
Every such matrix in $\SL_2(\matZ)$ is conjugate to $\matr 1 {\pm 1} 0 1$. 
\end{proof}

\begin{ex} The Dehn twists $T_m$ and $T_l$ generate $\MCG(T)$.
\end{ex}

We have seen that Dehn twists generate $\MCG(T)$, while many elements $A\in \MCG(T)$ are neither Dehn twists nor powers of Dehn twists (those with $\tr A \neq 2$). These two facts extend to higher-genus surfaces $S_g$; we start by generalising the first.

\subsection{Dehn twists generate}
Dehn twists are basic elements in the mapping class group, and we now prove that they generate the whole group. We restrict ourselves for simplicity to compact surfaces $S_{g,b}$.

\begin{teo} \label{twists:generate:teo}
The group $\MCG(S_{g,b})$ is generated by Dehn twists.
\end{teo}

To prove the theorem we will need some preliminary facts. We say that two non-separating simple closed curves in the interior of $S_{g,b}$ are \emph{related} if there is a combination of isotopies and Dehn twists transforming the first into the second.

\begin{figure}
\begin{center}
\includegraphics[width = 9 cm] {\iftoggle{BW}{related-BW}{related}}
\nota{Two curves $\alpha$ \iftoggle{BW}{}{(red) }and $\beta$ \iftoggle{BW}{}{(green) }intersecting in one point are related: we get $\alpha = T_\beta(T_\alpha(\beta))$.}
\label{related:fig}
\end{center}
\end{figure}

\begin{figure}
\begin{center}
\includegraphics[width = 12.5 cm] {\iftoggle{BW}{related2-BW}{related2}}
\nota{Pick two consecutive points in $\alpha$ \iftoggle{BW}{}{(green) }that intersect $\beta$\iftoggle{BW}{}{ (red)}. We can find a third non-separating curve $\gamma$ \iftoggle{BW}{}{(blue) }intersecting $\alpha$ and $\beta$ in $<k$ points. There are two cases (1) and (2): in (1) the curve $\gamma$ intersects $\beta$ in one point and is hence non-separating, in (2) we have two possibilities $\gamma_1, \gamma_2$, and one is certainly non-separating since $\beta$ is.}
\label{related2:fig}
\end{center}
\end{figure}

\begin{lemma}
The non-separating curves in $S_{g,b}$ are all related.
\end{lemma}
\begin{proof}
Let $\alpha$ and $\beta$ be two non-separating curves. 
Up to isotopy they intersect transversely into some $k$ points.
If $k=1$ they are related by a couple of Dehn twists as shown in Figure \ref{related:fig}. If $k=0$, since they are both non-separating, one sees easily that there is another curve $\gamma$ with $i(\alpha,\gamma) = i(\beta,\gamma)=1$: hence $\alpha$ and $\beta$ are both related to $\gamma$, so they are related themselves. If $k\geqslant 2$ we can find a curve $\gamma$ intersecting $\alpha$ and $\beta$ transversely into $<k$ points as shown in Figure \ref{related2:fig} and we proceed by induction on $k$.
\end{proof}

\begin{figure}
\begin{center}
\includegraphics[width = 12.5 cm] {\iftoggle{BW}{related3-BW}{related3}}
\nota{The arcs $\alpha$ \iftoggle{BW}{}{(red) }and $\beta$ \iftoggle{BW}{}{(green) }intersect at their endpoints and maybe at some $k$ interior points. If $k=0$ and they are oriented as in (1), a Dehn twist along $\gamma$ transforms $\beta$ into $\alpha$. If they are oriented as in (2), a negative Dehn twist along $\gamma$ transforms this configuration into (1). If $k>0$ we look at the first intersection point in $\alpha$. If the orientations are coherent as in (3) a Dehn twist along $\gamma$ decreases $k$. It the orientations are not coherent we change them as in (2).}
\label{related3:fig}
\end{center}
\end{figure}

We now consider $S_{g,b}$ with $b\geqslant 2$ and fix two points $p,q$ in two distinct boundary components of $S_{g,b}$. We consider all the properly embedded arcs in $S_{g,b}$ with endpoints at $p$ and $q$. As above, we say that two such arcs are \emph{related} if there is a combination of Dehn twists and isotopies transforming the first into the second.

\begin{lemma}
The arcs in $S_{g,b}$ with endpoints at $p$ and $q$ are all related.
\end{lemma}
\begin{proof}
Pick two arcs $\alpha$ and $\beta$ and put them into transverse position: they intersect at their endpoints $p$ and $q$ and maybe transversely at some $k$ other points. Figure \ref{related3:fig} shows that after isotopies and Dehn twists we get $\alpha = \beta$.
\end{proof}

We can finally prove Theorem \ref{twists:generate:teo}.
\begin{proof}[Proof of Theorem \ref{twists:generate:teo}]
Let $\varphi$ be a self-diffeomorphism of $S_{g,b}$ fixing pointwise the boundary. We prove that $\varphi$ is generated by isotopies and Dehn twists.

We first prove the case $g=0$ by induction on $b$. We know that $\MCG(S_{0,1})$ is trivial, so we suppose $b\geqslant 2$. Let $p,q$ be points in distinct boundary components of $S_{0,b}$, and $\alpha$ be an arc connecting them. All the arcs with endpoints in $p$ and $q$ are related, and hence $\alpha$ and $\varphi(\alpha)$ are. Therefore up to composing with Dehn twists and isotopies we may suppose that $\varphi$ is the identity of $\alpha$ and hence also on a tubular neighbourhood of $\alpha$. We cut $S_{0,b}$ along $\alpha$ and get $S_{0,b-1}$, with $\varphi$ transformed into a self-diffeomorphism of $S_{0,b-1}$. By induction on $b$ the new $\varphi$ is generated by Dehn twists and isotopies, so the original $\varphi$ also is.

We prove the case $g>0$ by induction on $g$. Let $\alpha$ be a non-separating simple closed curve. Since these are all related, up to isotopies and Dehn twists we may suppose that $\varphi$ is the identity on $\alpha$; as above we can cut $S_{g,b}$ along $\alpha$, get $S_{g-1,b+2}$ and conclude by induction on $g$. 
\end{proof}

\subsection{Action on simple closed curves}
We now show that every element $\varphi\in\MCG(S_g)$ is determined by the way it permutes the (isotopy classes of) simple closed curves in $S_g$.

\begin{prop} \label{action:curves:prop}
The action of $\varphi\in\MCG(S_g)$ on the isotopy classes of simple closed curves is faithful.
\end{prop}
\begin{proof}
The proof of Theorem \ref{omotopi:isotopi:teo} shows that if $\varphi$ fixes the isotopy classes of two essential multicurves $\mu$ and $\eta$ as in Figure \ref{multicurva_decomposizione:fig} then it is isotopic to the identity.
\end{proof}

Recall that every simple closed curve is oriented by assumption: if we considered \emph{unoriented} simple closed curves the proposition would be false, because in low genus there are some $\varphi$ that send every curve to its inverse: this happens for instance with the map $-I \in \SLZ = \MCG(T)$. 

\subsection{References}
The main source for this chapter is the book of Farb -- Margalit \cite{FM}, which contains a lot more information on the mapping class group of surfaces. We have also consulted Benedetti -- Petronio \cite{BP} and Thurson's notes \cite{Th} and book \cite{Th_book}, in particular for Smale's Theorem, that was originally proved in \cite{Sma}.

%% file: Teichmuller.tex
\chapter{Teichm\"uller space} \label{Teichmuller:chapter}
We have discovered in Chapter \ref{surfaces:chapter} that every surface $S_g$ of genus $g\geqslant 2$ can be equipped with a hyperbolic metric, and we have already noticed that this metric is not unique: this chapter is entirely devoted to studying this non-uniqueness phenomenon.

It turns out that $S_g$ admits a continuous family of non-equivalent hyperbolic metrics, that form altogether a nice topological space called the \emph{Teichm\"uller space} of $S_g$. We prove in his chapter that the Teichm\"uller space of $S_g$ is homeomorphic to an open ball of dimension $6g-6$. To prove this fact we will introduce and study concepts like \emph{length functions}, \emph{earthquakes}, and \emph{Fenchel--Nielsen coordinates}. The simple closed curves and their geodesic representatives will play a fundamental role in all the discussion.

\section{Introduction}
Let $S_g$ be as usual a closed orientable surface of genus $g$. We know that $S_g$ admits an elliptic, flat, or hyperbolic metric if and only if $g=0$, $g=1$, or $g\geqslant 2$ respectively. The elliptic metric on the two-sphere is unique up to isometries, but the flat and hyperbolic metrics on the other surfaces are not.

We want to define the space of all flat or hyperbolic metrics on $S_g$ when $g\geqslant 1$. There are two natural ways do to this:\index{Teichm\"uller space}\index{moduli space}

\begin{defn}
The \emph{moduli space} of $S_{g}$ is the set of all the flat or hyperbolic metrics on $S_{g}$ considered up to orientation-preserving isometries and rescaling.

The \emph{Teichm\"uller space} $\Teich(S_g)$ is the set of all the flat or hyperbolic metrics on $S_{g}$ considered up to isometries isotopic to the identity and rescaling.
\end{defn}

The \emph{rescaling} of the metric is a simple operation that takes place only on the flat metrics on the torus $T$. On the torus $T$ a flat metric $g$ can be rescaled by any constant $\lambda >0$ to give another flat metric $\lambda g$: the rescaling changes the lengths by a factor $\sqrt\lambda$ and the area by a factor $\lambda$. Up to rescaling, we may for instance require that $T$ has unit area.

At a first sight, the moduli space seems a more natural object to study. It turns out however that the Teichm\"uller space is homeomorphic (for some natural topology) to an open ball, while the moduli space is topologically more involved: it is then more comfortable to define and study the Teichm\"uller space first, and then consider the moduli space as a quotient of Teichm\"uller space.

\subsection{The Teichm\"uller space of the torus}
The flat metrics on the torus $T$ are classified quite easily. We have seen in Proposition \ref{flat:is:torus:prop} that every flat torus $T$ is of type $\matC/_\Gamma$ for some lattice $\Gamma <\matC$ isomorphic to $\matZ^2$. 

Fix two generators $m, l$ for $\pi_1(T)$. These are identified to two generators of $\Gamma$. Up to rescaling the metric, rotating $\matC$ around the origin, and reflecting along the real axis, we may suppose that these generators are 1 and some complex number $z$ lying in the upper half-plane $H^2$, so that $\Gamma = \langle 1,z\rangle$.

\begin{prop} \label{Teich:T:prop}
By sending the flat metric on $T$ to $z$ we get a bijection:
$$ \Teich(T) \longrightarrow H^2.$$
\end{prop}
\begin{proof}
The map is well-defined: two metrics related by an isometry isotopic to the identity produce the same $z$. The inverse $H^2 \to \Teich(T)$ is constructed by identifying $T$ with $\matC/_{\langle 1,z\rangle}$ sending $(m,l)$ to $(1,z)$.
\end{proof}

\begin{figure}
\begin{center}
\includegraphics[width = 8cm] {\iftoggle{BW}{toro_dominio-BW}{toro_dominio}}
\nota{The flat metric on the torus $T$ determined by $z\in H^2$ may be constructed by identifying the opposite sides of the parallelogram with vertices $0,1,z,z+1$. The lattice $\Gamma$ is generated by $1$ and $z$ and the parallelogram is a fundamental domain.}
\label{toro_dominio:fig}
\end{center}
\end{figure}

The flat metric that corresponds to $z \in H^2$ may be constructed by identifying the opposite sides of a parallelogram as in Figure \ref{toro_dominio:fig}. We will often tacitly identify $\Teich(T)$ with the projective plane $H^2$ via this bijective correspondence.

\begin{oss} \label{torus:isometries:oss}
Let $T=\matC/_\Gamma$ be a flat torus. Every translation $z \mapsto z +w$ in $\matC$ commutes with $\Gamma$ and hence descends to an isometry of $T$. Therefore the isometry group $\Iso^+(T)$ is not discrete: every flat torus is \emph{homogeneous}, \emph{i.e.}~for every pair of points $x,y \in T$ there is an isometry sending $x$ to $y$. 
\end{oss}

\begin{oss}
Let $T=\matC/_\Gamma$ be a flat torus. Every non-trivial element $\gamma\in\pi_1(T)$ is represented by a closed geodesic, unique up to translations. Its counterimage in $\matC$ consists of parallel lines whose slope depends only on $\gamma$.  The geodesic is simple if and only if $\gamma$ is primitive. 
\end{oss}

\subsection{Action of the mapping class group}
Recall that the mapping class group $\MCG(S_g)$ of $S_g$ is the group of all the orientation-preserving self-diffeomorphisms of $S_g$ considered up to isotopy (or equivalently, homotopy). We now show that the mapping class group of $S_g$ acts on its Teichm\"uller space.

A diffeomorphism $\varphi\colon S_g \to S_g$ transports a metric $m$ on $S_g$ into a new metric $\varphi_*m$ by pushing it forward as follows:
$$(\varphi_*m)_{\varphi(x)}\big(d\varphi_x(v), d\varphi_x(w)\big) = m_x(v,w).$$
If $m$ varies through an isotopy, the metric $\varphi_*m$ varies through a corresponding isotopy: therefore $\varphi$ acts on $\Teich(S_g)$ as follows
\begin{align*}
\Teich(S_g) & \longrightarrow \Teich(S_g) \\
[m] & \longmapsto [\varphi_*m]
\end{align*} 
If we vary $\varphi$ by an isotopy the action is unaffected. 
Therefore the mapping class group $\MCG(S_g)$ acts on $\Teich(S_g)$, and by definition the quotient
$$\Teich(S_g)/_{\MCG(S_g)}$$
is the moduli space of $S_g$. 

Everything can be written explicitly for the torus. Recall that we have identified $\Teich(T)$ with $H^2$ and $\MCG(T)$ with $\SLZ$, see Proposition \ref{Torelli:toro:prop}.

\begin{prop}  \label{toro:azione:prop}
The action of $\MCG(T)$ on $\Teich(T)$ is the following action of $\SLZ$ on $H^2$ as M\"obius transformations: 
$$\begin{pmatrix} a & b \\ c & d \end{pmatrix} \colon z \longmapsto \frac{az-b}{-cz+d}.$$
\end{prop}
\begin{proof}
The metric $z$ assigns to $T$ the structure $\matR^2/_\Gamma$ with $\Gamma = \langle 1,z \rangle$ and $(m,l)$ sent to $(1,z)$. 

Pick $\varphi =  \left(\begin{smallmatrix} a & b \\ c & d \end{smallmatrix}\right) \in \SLZ =\MCG(T)$. Since $\varphi^{-1} = \left(\begin{smallmatrix} d & -b \\ -c & a \end{smallmatrix}\right)$, in the new metric $\varphi_*$ we assign to $(m,l)$ the translations $(d-cz, -b+az)$, which transform via rotations and dilations into $(1,  \frac{az-b}{-cz+d}).$
\end{proof}

We note in particular that $\MCG(T)$ acts via isometries on the hyperbolic plane $H^2$. The kernel of the action is $\{\pm I\}$: two matrices $A$ and $-A$ act in the same way on $\Teich(T)$. 

\begin{cor}
The moduli space of $T$ is the orbifold $H^2/_{\PSLZ}$.
\end{cor}

The orbifold $H^2/_{\PSLZ}$ is described in Figure \ref{orbifold_Farey:fig}. It has two singular points of order 2 and 3: these represent the square and hexagonal torus, see Section \ref{piatte:subsection} and Figure \ref{modular_group:fig}.

We have constructed a bijection between $\Teich(T)$ and $H^2$, and we now want to construct for $g\geqslant 2$ an analogous identification between $\Teich(S_g)$ and some open set of $\matR^N$ for some $N$ depending on $g$. To this purpose we need to introduce some concepts.

\section{Earthquakes and length functions}
Simple closed geodesics are a formidable tool to study $\Teich(S_g)$: we can use a simple closed geodesic to \emph{twist} a metric (the operation is called an \emph{earthquake}), and by simply looking at the lengths of the other closed geodesics we can measure how the metric varies along this transformation. We introduce these operations here; later on, we will use them to parametrize $\Teich(S_g)$.

\subsection{Earthquakes} \label{earthquakes:subsection}
The hyperbolic, flat, and spherical metrics on surfaces may be twisted along simple closed geodesics: this operation is called an \emph{earthquake}.\index{earthquake} 

Let $m$ be a complete hyperbolic, flat, or elliptic metric on an oriented surface $S$ and $\gamma$ be a simple closed geodesic in $S$. Recall that $\gamma$ is a map $\gamma \colon S^1 \to S$. Fix an angle $\theta \in \matR$. Informally, a new complete hyperbolic, flat, or elliptic metric $m_\theta$ on $S$ is constructed by cutting $S$ along $\gamma$ and regluing with a counterclockwise twist of angle $\theta$. Formally, the new metric is defined as follows.

\begin{figure}
\begin{center}
\includegraphics[width = 10cm] {\iftoggle{BW}{prodotto-BW}{prodotto}}
\nota{A $R$-annulus around a simple closed geodesic $\gamma$ on a hyperbolic surface is the quotient of a $R$-neighbourhood of a line $l$ by a hyperbolic transformation. The orthogonal (\iftoggle{BW}{light grey}{green}) geodesic segments are parametrized by arc length as $[-R,R]$, hence the $R$-annulus is naturally parametrized as $S^1 \times [-R,R]$.}
\label{prodotto:fig}
\end{center}
\end{figure}

In the hyperbolic case, recall from Proposition \ref{tube:prop} that $\gamma$ has a $R$-neighbourhood isometric to a $R$-tube for some $R>0$. A $R$-tube in dimension two is a $R$-annulus as in Figure \ref{prodotto:fig}, defined by quotienting a $R$-neighbourhood of a line $l$ by a hyperbolic transformation. The $R$-annulus is naturally parametrized as $S^1\times [-R,R]$, where $\{e^{it}\}\times [-R,R]$ is the geodesic segment orthogonal to $\gamma$ in $\gamma(e^{it})$ parametrized by arc length.
The flat and elliptic cases are analogous.

\begin{figure}
\begin{center}
\includegraphics[width = 10cm] {\iftoggle{BW}{twist-BW}{twist}}
\nota{To define the earthquake we pick a diffeomorphism of the $R$-annulus that modifies the orthogonal segments as shown here.}
\label{twist:fig}
\end{center}
\end{figure}

We choose a diffeomorphism $\varphi$ of $S^1\times [-R, R]$ that curves the segments left-wise with step $\theta$ as in Figure \ref{twist:fig}-(right).  More precisely, let $f\colon[-R, R] \to \matR$ be a smooth function which is zero on $\left[-R, -\frac R 2\right]$ and is constantly $\theta$ on $\left[\frac R 2,R\right]$. We set $\varphi(e^{it},s) = (e^{i(t+f(s))},s)$.

We define a new metric $m_\theta$ on $S_g$ as follows: the metric tensor $m_\theta$ coincides with $\varphi_* m$ on the $R$-annulus and coincides with $m$ on the complement of the $\frac{R}2$-annulus $S^1\times \left[-\frac R2, \frac R2\right]$.

\begin{prop} The metric tensor $m_\theta$ is well-defined and gives a complete hyperbolic, flat, or elliptic metric to $S_g$. 
\end{prop}
\begin{proof}
It is well-defined because $m$ and $m_\theta$ coincide on $S^1\times \left[\frac R 2, R\right]$, since $(e^{it},s)\mapsto (e^{i(t+\theta)},s)$ is an isometry of the $R$-annulus. It is hyperbolic, flat, or elliptic because both patches $m$ and $\varphi_*m$ are.
\end{proof}

\begin{oss} \label{R:oss}
In the new metric $m_\theta$ the curve $\gamma$ is still a simple closed geodesic of the same length as before, and its $R$-neighbourhood is also unchanged. 
\end{oss}

Of course by deforming an elliptic metric in this way we get nothing new, because all the elliptic metrics on $S^2$ are isometric. Earthquakes are interesting only in the flat and hyperbolic geometries.

\subsection{The earthquake map}
We consider the surface $S_g$ with $g\geqslant 1$ and show that earthquakes define nice actions on the Teichm\"uller space $\Teich(S_g)$. 

If $m$ is any flat or hyperbolic metric and $\gamma$ is a non-trivial simple closed curve in $S_g$, we define $m_\theta^\gamma$ to be the flat or hyperbolic metric obtained from $m$ via an earthquake of angle $\theta$ performed along the unique (up to translations if $g=1$) simple closed geodesic homotopic to $\gamma$ in the metric $m$.\index{earthquake map}

\begin{prop}
The \emph{earthquake map}
\begin{align*}
E_\gamma\colon \matR \times \Teich(S_g) & \longrightarrow \Teich(S_g) \\
 (\theta, m) & \longmapsto m_\theta^{\gamma} 
\end{align*}
is a well-defined action of $\matR$ on $\Teich(S_g)$. The map $E_\gamma$ depends only on the homotopy class of $\gamma$.
\end{prop}
\begin{proof}
The only ambiguity in the definition of $m_\theta^\gamma$ is the choice of the function $f\colon [-R,R] \to \matR$. If we use another function $f'$ the resulting metric $(m_\theta^\gamma)'$ changes only by an isotopy: the diffeomorphism of $S_g$ which is the identity outside the $R$-annulus and sends $(e^{it},s)$ to $(e^{i(t+f(s)-f'(s))},s)$ is an isometry between $m_\theta^\gamma$ and $(m_\theta^\gamma)'$, and is clearly isotopic to the identity.

To prove that $E_\gamma$ is an action we need to check that 
$$m_{\theta+\theta'}^{\gamma} = \left(m_{\theta'}^\gamma\right)_\theta^\gamma.$$
By Remark \ref{R:oss} we can take the same $R$-annulus to compose two earthquakes and the equality follows.
\end{proof}

Like the Dehn twists $T_\gamma$ defined in Section \ref{Dehn:twist:subsection}, the action $E_\gamma$ depends on the orientation of $S_g$ but not on the orientation of $\gamma$. There is indeed a strong relation between earthquakes and Dehn twists on $\gamma$: as objects acting on Theichm\"uller space, the first generalise the second.

\begin{prop} \label{root:prop}
We have $T_\gamma(m) = E_\gamma(2\pi,m)$.
\end{prop}
\begin{proof}
It follows directly from the definitions.
\end{proof}

We now want to study the Teichm\"uller space and the action of $E_\gamma$ on it. In mathematics a space is often beautifully described by some natural functions defined on it: this role is played here by the length functions of closed curves.

\subsection{Length functions}
A homotopically nontrivial (possibly non simple) closed curve $\gamma$ in $S_g$ defines a \emph{length function}\index{length function}
$$\ell^\gamma\colon \Teich(S_g) \to \matR_{>0}$$
which assigns to a metric $m\in\Teich(S_g)$ the length $\ell^\gamma(m)$ of the unique closed geodesic homotopic to $\gamma$. 

When $g=1$ we must actually specify a couple of things in the definition: the closed geodesic $\gamma$ is unique only up to translations, which do not affect its length; 
on the other hand rescaling \emph{does} affect lengths, so to get a well-defined length function we always rescale the metric $m$ to have unit area. 

We want to study these length functions, and as usual we start by analysing the simpler torus flat world, where everything can be described explicitly.

\subsection{Length functions on the torus} \label{length:torus:subsection}
We denote every (isotopy class of) non-trivial simple closed curve on $T$ with a coprime pair $(p,q)$ of integers and we identify the Teichm\"uller space $\Teich(T)$ with $H^2\subset\matC$, see Propositions \ref{mn:prop} and \ref{Teich:T:prop}. The length functions may be written explicitly.

\begin{prop} \label{lunghezza:toro:prop}
The formula holds:
$$\ell^{(p,q)}(z) = \frac{|p+qz|}{\sqrt{\Im z }}$$
for every simple closed curve $(p,q)$ and every metric $z\in H^2$.
\end{prop}
\begin{proof}
Up to rescaling we have $T=\matR^2/_{\Gamma}$ with $\Gamma = \langle 1,z \rangle $. The translation in $\Gamma$ corresponding to $(p,q)$ is $p\cdot 1 + q\cdot z$ and the closed geodesic it produces has length $|p+qz|$. The area of the torus $T$ is $\Im z$ (see a fundamental domain in Figure \ref{toro_dominio:fig}) and hence we must rescale it by $1/\sqrt{\Im z}$.
\end{proof}

\begin{figure}
\begin{center}
\includegraphics[width = 9cm] {\iftoggle{BW}{toro_terremoto-BW}{toro_terremoto}}
\nota{A torus with metric $z$ (left) twisted along the horizontal curve $\gamma$ (right). The curve $\alpha$ is a closed geodesic in the new metric.}
\label{toro_terremoto:fig}
\end{center}
\end{figure}

We can also write the earthquake action on the meridian $m=(1,0)$ of $T$.

\begin{ex} \label{terremoto:toro:ex}
We have:
$$E_m(\theta, z) = z+\frac {\theta}{2\pi}.$$
\end{ex}
\begin{proof}[Hint]
Look at the closed geodesic $\alpha$ in Figure \ref{toro_terremoto:fig}.
\end{proof}

We have discovered in particular that $E_m$ is a parabolic isometry of $H^2$, and we can now consider the earthquake action along a generic curve $(p,q)$.

\begin{cor} 
The earthquake action $E_{(p,q)}$ is the 1-parameter family of parabolic transformations with fixed point $-\frac pq \in \partial H^2$. 
\end{cor}
\begin{proof}
We know this when $(p,q)=(1,0)=m$. In general, we send $(1,0)$ to $(p,q)$ via some element of the mapping class group: this element acts by isometries of $H^2$ by Proposition \ref{toro:azione:prop} and sends $\infty$ to $-\frac pq$. 
\end{proof}

The orbits of $E_{(p,q)}$ are the horospheres centred at $-\frac pq$. 
We have discovered that, quite unexpectedly, the hyperbolic geometry of the plane $H^2$ is well designed to model the Teichm\"uller space of the flat torus $T$. 

It is now natural to identify the \emph{unoriented} simple closed curve $\pm (p,q)$ with the rational point $-\frac pq$ in $\matR\cup\{\infty\} \subset \partial H^2$. We get the following convexity property.

\begin{cor} \label{convex:cor}
Let $\gamma$ and $\eta$ be two simple closed curves on $T$.
If $i(\gamma, \eta)>0$ then the length function $\ell^\eta$ is strictly convex on the orbits of $E_{\gamma}$.
\end{cor}
\begin{proof}
We may suppose $\gamma=m=(1,0)$ and note that the condition $i(\gamma,\eta)>0$ translates into $\eta = (p,q) \neq (\pm 1,0)$. The function in Proposition \ref{lunghezza:toro:prop} is strictly convex on the horospheres $\Im z = k$ when $(p,q)\neq (\pm 1,0)$.
\end{proof}

Summing up, the torus picture is the following: 
\begin{itemize}
\item the mapping class group acts on the Teichm\"uller space roughly like $\PSLZ$ acts isometrically on the hyperbolic plane $H^2$, 
\item the unoriented simple closed curves are in 1-1 correspondence with the rational points in $\partial H^2 = \matR\cup \infty$, 
\item earthquakes act like parabolic transformations centred at these rational points, \item the length function $\ell^\gamma$ is constant at the horospheres centred at $\gamma$ but strictly convex at the horospheres centred at all the other curves.
\end{itemize}

The main goal of Chapters \ref{Teichmuller:chapter} and \ref{automorfismi:chapter} is to draw a similar picture for surfaces $S_g$ of higher genus $g\geqslant 2$. Everything is more difficult in the hyperbolic world, because there are no nice explicit formulas describing the length functions, the action of the mapping class group, the simple closed curves, and the earthquakes. We now start by generalising the last point: the strict convexity of length functions.

\subsection{Convexity of the length functions} \label{convexity:lengths:subsection}
We consider $S_g$ with $g\geqslant 2$. Our aim now is to generalise Corollary \ref{convex:cor} to the hyperbolic setting and then later use this convexity property to parametrize $\Teich(S_g)$. We need a preliminary result.\index{convexity of the length functions}

\begin{ex} \label{convessa:ex}
Let $f\colon \matR^m \times \matR^n \to \matR_{\geqslant 0}$ be strictly convex and proper. The function
\begin{align*}
F\colon \matR^n & \longrightarrow \matR_{\geqslant 0} \\
y & \longmapsto \min \left\{f(x,y)\ \big| \ x \in \matR^m \right\}
\end{align*}
is well-defined, strictly convex, and proper. 
\end{ex}

The following generalisation of Corollary \ref{convex:cor} says that length functions are convex on the earthquakes orbits, and very often they are strictly convex and proper.

\begin{prop} \label{strettamente:convessa:prop}
Let $\eta$ and $\gamma$ be two homotopically non-trivial simple closed curves in $S_g$ and $m$ be a hyperbolic metric on $S_g$. The function 
\begin{align*}
\matR & \longrightarrow \matR_{\geqslant 0} \\
\theta & \longmapsto \ell^\eta(m_\theta^\gamma)
\end{align*}
is
\begin{itemize}
\item constant if $i(\eta,\gamma)=0$,
\item strictly convex and proper if $i(\eta, \gamma)>0$.
\end{itemize}
\end{prop}
\begin{proof}
The metric $m_\theta^\gamma$ is obtained by twisting $m$ of an angle $\theta$ along the simple closed geodesic $\gamma$. If $i(\eta,\gamma)=0$ the curves $\eta$ and $\gamma$ are disjoint geodesics in $m$ and $\eta$ is not affected by the earthquakes that we perform near $\gamma$, hence the function $\ell^\eta(m_\theta^\gamma)$ is constant.

Consider the case $n=i(\eta,\gamma)>0$. Denote by $\bar{\eta}^\theta$ the geodesic representative of $\eta$ in the twisted metric $m_\theta^\gamma$: it intersects $\gamma$ transversely in $n$ points. 

\begin{figure}
\begin{center}
\includegraphics[width = 10cm] {\iftoggle{BW}{spezza_new-BW}{spezza_new}}
\nota{A geodesic in $m_\theta^\gamma$ can be seen on the original metric $m$ as follows: it is a geodesic outside the annulus and deviates on the left by an angle $\theta$ every time it crosses it (centre). We may simplify the picture by describing it as a broken geodesic line that makes a left $\theta$-jump each time that it crosses $\gamma$ (right).}
\label{spezza:fig}
\end{center}
\end{figure}

Fix a sufficiently small $R$-annulus around $\gamma$ and note that the geodesics in $m_\theta^\gamma$ can be seen in the original metric $m$ as follows: these are curves that are geodesic outside the $R$-annulus and deviate smoothly on the left each time they cross it as in Figure \ref{spezza:fig}-(centre).
We may substitute each smooth deviation with a broken jump as shown in Figure \ref{spezza:fig}-(right) and get a bijection
$$\left\{\begin{array}{c} {\rm closed\ geodesics} \\ {\rm with\ respect\ to\ } m_\theta^\gamma \end{array} \right\}
\longleftrightarrow \left\{\begin{array}{c} {\rm broken\ geodesics} \\ {\rm with\ respect\ to\ } m \end{array}\right\}$$
where a \emph{broken geodesic} is a geodesic that at every crossing of $\gamma$ jumps to the left at distance $\frac{\theta L(\gamma)}{2\pi}$ and then keeps going on, leaving $\gamma$ with the same incidence angle (here $L(\gamma)$ is the length of $\gamma$). This correspondence is useful because it preserves the lengths: the length of the closed geodesic for $m_\theta^\gamma$ is equal to the length of the corresponding broken geodesic (which is the sum of the lengths of its components), because the segments in Figure \ref{spezza:fig}-(left) and (right) are isometric.

\begin{figure}
\begin{center}
\includegraphics[width = 6cm] {\iftoggle{BW}{spezza_su-BW}{spezza_su}}
\nota{The \iftoggle{BW}{}{(blue) }line $l$ is a lift of $\bar \eta$. Consider $n$ consecutive intersections with lifts $r_1,\ldots, r_{n+1}$ of $\gamma$\iftoggle{BW}{}{ (in red)}, with $\tau(r_1) = r_{n+1}$, and parametrize each $r_i$ with $\matR$ via arc-length. We set $\theta' = \theta L(\gamma)/2pi$.}
\label{spezza_su:fig}
\end{center}
\end{figure}

We lift this description to the universal cover $\matH^2$. We fix a lift $l$ of $\bar \eta = \bar{\eta}^0$ and pick $n+1$ consecutive intersections $r_1,\ldots, r_{n+1}$ of $l$ with the lifts of $\gamma$ as in Figure \ref{spezza_su:fig}. The hyperbolic transformation $\tau$ with axis $l$ corresponding to $\eta$ sends $r_1$ to $r_{n+1}$. We parametrise each $r_i$ with $\matR$ via arc length.

The closed geodesic $\bar{\eta}^\theta$, represented as a broken geodesic, lifts to a broken geodesic which starts at some point $x_1-\theta L(\gamma)/2\pi\in r_1$ and arrives at some other point $x_2\in r_2$, then jumps on the left at distance $\theta L(\gamma)/2\pi$ and starts again from $x_2-\theta L(\gamma)/2\pi$, and so on until it reaches the point $\tau(x_1)\in \tau(r_1) = r_{n+1}$. If we make the points $x_1\in r_1, \ldots, x_n \in r_n$ vary we get various broken paths in this way, but only one arrives and exits from each line $r_i$ with the same incidence angles and thus represents $\bar{\eta}^\theta$. The other broken paths represent piecewise-geodesic curves homotopic to $\bar{\eta}^\theta$ and are therefore longer than $\bar{\eta}^\theta$. Hence
$$\ell^\eta(m_\theta^\gamma) = \min \left\{\sum_{i=1}^n d\bigg(x_i-\frac{\theta L(\gamma)}{2\pi}, x_{i+1}\bigg) \ \Bigg|\ (x_1,\ldots, x_n) \in \matR^n \right\}$$
where $x_{n+1} = \tau(x_1)$.
We can now prove that the function $\theta \mapsto \ell^\eta(m_\theta^\gamma)$ is proper and strictly convex. The function
\begin{align*}
\psi\colon \matR^{2n} & \longrightarrow \matR \\
(x_1,y_1, \ldots, x_n,y_n) & \longmapsto \sum_{i=1}^n d(y_i,x_{i+1})
\end{align*}
where $x_{n+1}=\tau(x_1)$ is strictly convex and proper by Proposition \ref{distanza:convessa:prop}. The auxiliary function
\begin{align*}
\phi\colon \matR^{2n}\times \matR & \longrightarrow \matR \\
(x,\theta) & \longmapsto \psi(x)
\end{align*}
is only convex, but its restriction to the subspace
$$H = \bigg\{y_i = x_i - \frac{\theta l(\gamma)}{2\pi} \bigg\}$$
is strictly convex and proper, because the subspace $H$ is not parallel to the direction $(0,\ldots,0,1)$. The coordinates $x_i$ and $\theta$ identify $H$ with $\matR^n \times \matR$. The restriction $f = \phi|_H$ is hence a function $f\colon \matR^n \times \matR \to \matR$ and we obtain
$$\ell^\eta(m_\theta^\gamma) = \min \left\{f(x,\theta) \ \big| \ x \in \matR^n \right\}.$$
By Exercise \ref{convessa:ex} the function $\theta \mapsto \ell^\eta(m_\theta^\gamma)$ is strictly convex and proper. 
\end{proof}

We now employ this convexity property to prove some facts on the earthquakes and the Teichm\"uller space.

\subsection{Earthquakes on essential multicurves}
We use the convexity of the length functions to prove the following. We suppose again that $g\geqslant 2$.

\begin{cor} \label{distinti:cor}
For every simple closed curve $\gamma$, the earthquake action $E_\gamma$ on $\Teich(S_g)$ is free.
\end{cor}
\begin{proof}
Suppose by contradiction that $m = m_{\theta_0}^\gamma$ for some $\theta_0>0$. Then $m= m_{n\theta_0}^\gamma$ for every $n\in\matZ$. Let $\eta$ be a simple closed curve with $i(\eta,\gamma)>0$, which exists by Exercise \ref{eta:prop}; the function $\theta \mapsto \ell^\eta(m_\theta^\gamma)$ is strictly convex and constant on $\{n\theta_0, n \in \matZ\}$, a contradiction.
\end{proof}

The earthquake action is defined more generally for essential multicurves. An essential multicurve $\mu = \gamma_1\sqcup \cdots \sqcup \gamma_k$ of $S_g$ determines an action
\begin{align*}
E_\mu\colon \matR^{k} \times \Teich(S_g) & \longrightarrow \Teich(S_g) \\
 (\theta, m) & \longmapsto m_\theta^{\mu}
\end{align*}
where $\theta = (\theta_1,\ldots, \theta_{k})$ and $m_\theta^\mu = m_{\theta_1}^{\gamma_1} \circ \cdots \circ m_{\theta_{k}}^{\gamma_{k}}$. Note that the actions on disjoint curves commute.

\begin{figure}
\begin{center}
\includegraphics[width = 9cm] {\iftoggle{BW}{dual_curves-BW}{dual_curves}}
\nota{Choose for each component $\gamma_i$ of a pants decomposition $\mu$ a curve $\gamma_i'$ that intersects $\gamma_i$ in one or two points and is disjoint from the other components of $\mu$. There are two cases to consider, depending on whether the two pants adjacent to $\gamma_i$ are distinct (left) or not (right). }
\label{dual_curves:fig}
\end{center}
\end{figure}

\begin{cor}
For every essential multicurve $\mu$, the earthquake action $E_\mu$ on $\Teich(S_g)$ is free. 
\end{cor}
\begin{proof}
We may complete $\mu$ to a pants-decomposition $\mu=\gamma_1\sqcup \ldots \sqcup \gamma_{3g-3}$. Pick for every $i=1,\ldots, 3g-3$ a curve $\gamma_i'$ as in Figure \ref{dual_curves:fig} such that $i(\gamma_i,\gamma_i') > 0$ for all $i$ and $i(\gamma_i, \gamma_j')=0$ for all $i\neq j$. 

Suppose by contradiction that $m=m_\theta^\mu$ for some $\theta\neq 0$: hence $m_{n\theta}^\mu = m$ for all $n\in\matZ$. There is an $i$ such that $\theta_i\neq 0$. The length function $\ell^{\gamma_i'}(m_\theta^\mu)$ depends only on $\theta_i$ and not on the other coordinates of $\theta$: therefore it equals $\ell^{\gamma_i'}(m_{\theta_i}^{\gamma_i})$ which is strictly convex, a contradiction.
\end{proof}

\section{Fenchel--Nielsen coordinates}
It is now time to fix a global set of coordinates for the Teichm\"uller space when $g\geqslant 2$. These are the \emph{Fenchel--Nielsen coordinates} and they identify $\Teich(S_g)$ with $\matR^{6g-6}$, more precisely with $\matR_{>0}^{3g-3} \times \matR^{3g-3}$. 

\subsection{The coordinates}
We want to construct a parametrisation for $\Teich(S_g)$ when $g\geqslant 2$. To identify a finite-dimensional vector space with $\matR^n$ one needs to fix a basis; likewise, here the parametrisation depends on the choice of a  frame.\index{Fenchel--Nielsen coordinates} 

Let $S_g$ be oriented. A \emph{frame} for $S_g$ consists of two essential multicurves $\mu$ and $\nu$ in minimal position, such that:

\begin{enumerate}
\item the multicurve $\mu$ is a pants decomposition,
\item the multicurve $\nu$ decomposes every pair-of-pants in two hexagons.
\end{enumerate}

\begin{figure}
\begin{center}
\includegraphics[width = 10cm] {\iftoggle{BW}{decomposizione2-BW}{decomposizione}}
\nota{A frame for the Fenchel-Nielsen coordinates consists of a (\iftoggle{BW}{light grey}{red}) pants decomposition $\mu$ and a (\iftoggle{BW}{dark grey}{blue}) transverse multicurve $\nu$ that cuts each pair of pants into two hexagons. The number of components of $\mu$ is $3g-3$, that of $\nu$ can vary.}
\label{decomposizione:fig}
\end{center}
\end{figure}

An example that generalises easily to any genus $g\geq 2$ is shown in Figure \ref{decomposizione:fig}. The pants decomposition $\mu = \gamma_1\sqcup \ldots \sqcup \gamma_{3g-3}$ consists of $3g-3$ curves, while the number of curves in $\nu$ is not fixed a priori and depends on our choice of $\nu$.
We now show that a frame induces a \emph{Fenchel--Nielsen} map
\begin{align*}
\FN\colon \Teich(S_g) & \longrightarrow \matR_{>0}^{3g-3} \times \matR^{3g-3} \\
m & \longmapsto  (l_1,\ldots, l_{3g-3}, \theta_1, \ldots, \theta_{3g-3}).
\end{align*}
The map $\FN$ is defined as follows. Let $m\in \Teich(S_g)$ be a hyperbolic metric. The $3g-3$ \emph{length parameters} $l_i = \ell^{\gamma_i}(m)$ are defined using the length functions: the multicurve $\mu$ has a unique geodesic representative 
$$\bar \mu = \bar{\gamma}_1\sqcup \ldots \sqcup \bar{\gamma}_{3g-3}$$ 
in the metric $m$ by Corollary \ref{isotopia:geodetiche:cor}, and $l_i$ is the length of $\bar{\gamma}_i$. Note that these parameters depend only on $\mu$ and not on $\nu$.

The \emph{torsion angles} $\theta_i$ are more subtle to define: the geodesic multicurve $\bar \mu$ decomposes $S_g$ into geodesic pairs-of-pants, and the angle $\theta_i$ measures somehow the way the two geodesic pairs-of-pants are glued along the closed geodesic $\bar \gamma_i$. The precise definition of $\theta_i$ needs the auxiliary multicurve $\nu$. 

\begin{figure}
\begin{center}
\includegraphics[width = 11cm] {\iftoggle{BW}{Teichmuller0-BW}{Teichmuller0}}
\nota{A closed geodesic $\bar \gamma_1$ and the two adjacent pairs-of-pants. The torsion parameter $\theta_1$ measures the distance (in the universal covering) between two orthogeodesics (\iftoggle{BW}{in light grey}{coloured in green}) via the formula $\theta_i = \frac{2\pi s_i}{l_i}$.}
\label{Teichmuller:fig}
\end{center}
\end{figure}

We fix $i=1$ for simplicity and define $\theta_1$.
Figure \ref{Teichmuller:fig}-(left) shows the two geodesic pants adjacent to $\bar{\gamma}_1$ (they might coincide). The second multicurve $\nu$ intersects these pants in four \iftoggle{BW}{dark grey}{blue} arcs, two of which $\lambda$, $\lambda'$ intersect $\bar{\gamma}_1$: we pick one, say $\lambda$. We fix a lift $\tilde P\in\matH^2$ of $P=\bar{\gamma}_1\cap\lambda$ and we lift all the curves incident to $P$: the geodesic $\bar{\gamma}_1$ lifts to a line $\tilde{\gamma}_1$ and $\lambda$ lifts to a (non-geodesic) curve $\tilde\lambda$ that connects two lifts $\tilde{\gamma}_2$ and $\tilde{\gamma}_3$ of the closed geodesics $\bar{\gamma}_2$ and $\bar{\gamma}_3$. See Figure \ref{Teichmuller:fig}-(right).

We draw as in the figure the unique orthogeodesics connecting $\tilde{\gamma}_1$ to $\tilde{\gamma}_2$ and $\tilde{\gamma}_3$ and we denote by $s_1$ the signed length of the segment in $\tilde{\gamma}_1$ comprised between these two orthogeodesics, with positive sign if (as in the figure) an observer walking on a orthogeodesic towards $\tilde{\gamma}_1$ sees the other orthogeodesic on its left (here we use the orientation of $S_g$). 

By repeating this construction for each $\bar{\gamma}_i$ we find some real numbers $s_i$. Finally, the torsion parameter $\theta_i$ is 
$$\theta_i = \frac {2\pi s_i}{l_i}.$$

\begin{figure}
\begin{center}
\includegraphics[width = 5cm] {\iftoggle{BW}{Teichmuller2-BW}{Teichmuller2}}
\nota{If we pick $\lambda'$ instead of $\lambda$ we find a segment of the same length $s_1$. This holds because the two geodesic pairs-of-pants incident to $\gamma_1$ decompose into two hexagons isometric to $A$ and $B$ as shown, and the sides of $A$ and $B$ contained in $\tilde{\gamma}_1$ have the same length $\frac{l(\gamma_1)}2$, that is half the length of $\gamma_1$.}
\label{Teichmuller2:fig}
\end{center}
\end{figure}

\begin{teo}[Fenchel-Nielsen coordinates]
The map $\FN$ is well-defined and is a bijection.
\end{teo}
\begin{proof}
We first note that in the definition of the torsion parameters we could have chosen $\lambda'$ instead of $\lambda$. We would have obtained the same length $s_i$ as shown in Figure \ref{Teichmuller2:fig}. Moreover a hyperbolic metric $m'$ isometric to $m$ through a diffeomorphism $\varphi$ isotopic to the identity has the same parameters $l_i$ and $\theta_j$ since they depend only on the isotopy classes of $\mu$ and $\nu$. Therefore $\FN$ is well-defined.

We prove that $\FN$ is surjective. For every vector $(l_1,\ldots,l_{3g-3})\in\matR_{>0}^{3g-3}$ we may use Proposition \ref{pantaloni:prop} and construct a metric on $S_g$ by assigning to each pair-of-pants of the pants decomposition $\mu$ the (unique) hyperbolic metric with boundary lengths $l_i$. We get a metric with some arbitrary torsion angles $\theta$, which can be changed arbitrarily by an earthquake along $\mu$: it is easy to check that an earthquake with angles $\theta'$ changes the torsion angles from $\theta$ to $\theta+\theta'$, hence any torsion parameter can be realised and $\FN$ is surjective.

We prove that $\FN$ is injective. If $\FN(m) = \FN(m')$, up to acting via earthquakes we may suppose that $\FN(m)=\FN(m') = (l_1,\ldots,l_{3g-3},0,\ldots,0)$. Since the torsion parameter is zero, the orthogeodesics in Figure \ref{Teichmuller:fig}-(right) match and project in $S_g$ to a geodesic multicurve $\bar \nu$ isotopic to $\nu$ and orthogonal to $\bar\mu$. Therefore $S_g\setminus (\bar\mu \cup \bar\nu)$ is a tessellation of $S_g$ into right-angled hexagons, determined by the lengths $l_i$. Both metrics $m$ and $m'$ have the same tessellation and are hence isometric, via an isometry which is isotopic to the identity.
\end{proof}

\begin{oss}
As shown in the proof, the torsion parameters for $m$ are zero if and only if the geodesic representatives $\bar \nu$ and $\bar \mu$ of $\nu$ and $\mu$ are everywhere orthogonal. 
\end{oss}

\subsection{Length functions of $9g-9$ curves}
It is now natural to ask whether the length functions determine every point in $\Teich(S_g)$. The answer is positive; as usual, to warm up we start by examining the torus.

\begin{figure}
\begin{center}
\includegraphics[width = 7cm] {\iftoggle{BW}{torus_3_curves-BW}{torus_3_curves}}
\nota{The curves $\gamma$\iftoggle{BW}{}{ (red)}, $\gamma'$\iftoggle{BW}{}{ (blue)}, and $\gamma'' = T_{\gamma}(\gamma')$\iftoggle{BW}{}{ (green)} on the torus.}
\label{torus_3_curves:fig}
\end{center}
\end{figure}

Let $\gamma,\gamma'$ be two simple closed curves in the torus $T$ with $i(\gamma, \gamma')=1$ and let $\gamma'' = T_{\gamma}(\gamma')$ be obtained by Dehn twisting $\gamma'$ along $\gamma$, see Figure \ref{torus_3_curves:fig}.

\begin{prop} \label{3:prop}
The map
\begin{align*}
L\colon \Teich(T) & \longrightarrow \matR_{>0}^{3} \\
m & \longmapsto \big(\ell^{\gamma}(m), \ell^{\gamma'}(m), \ell^{\gamma''}(m)\big)
\end{align*}
is injective.
\end{prop}
\begin{proof}
After fixing $\gamma$ and $\gamma'$ as a homology basis, we have $\gamma = (1,0)$, $\gamma' = (0,1)$, and $\gamma'' = (1,-1)$. Proposition \ref{lunghezza:toro:prop} gives
$$L(z) = \left(\frac{1}{\sqrt{\Im z}}, \frac{|z|}{\sqrt{\Im z}}, \frac{|z-1|}{\sqrt{\Im z}} \right)$$
which is easily seen to be injective on $\Teich(T) = H^2$.
\end{proof}

A similar set of $9g-9$ curves does the job on $S_g$ when $g\geqslant 2$.
Let $\mu = \gamma_1\sqcup \ldots \sqcup \gamma_{3g-3}$ be a pants decomposition for $S_g$. For each $\gamma_i$ we choose a curve $\gamma_i'$ as in Figure \ref{dual_curves:fig}, and we indicate by $\gamma_i'' = T_{\gamma_i}(\gamma_i')$ the curve obtained by Dehn-twisting $\gamma_i'$ along $\gamma_i$, see an example in Figure \ref{dual_curves2:fig}.

\begin{figure}
\begin{center}
\includegraphics[width = 4cm] {\iftoggle{BW}{dual_curves2-BW}{dual_curves2}}
\nota{The curves $\gamma_i$\iftoggle{BW}{}{ (red)}, $\gamma_i'$\iftoggle{BW}{}{ (blue)}, and $\gamma_i'' = T_{\gamma_i}(\gamma_i')$\iftoggle{BW}{}{ (green)} when $\gamma_i$ is adjacent twice to the same pair of pants. It is an instructive exercise to draw $\gamma_i'$ and $\gamma_i''$ when $\gamma_i$ is incident to distinct pair of pants.}
\label{dual_curves2:fig}
\end{center}
\end{figure}

\begin{prop} \label{9g-9:prop}
The map
\begin{align*}
L\colon \Teich(S_g) & \longrightarrow \matR_{>0}^{9g-9} \\
m & \longmapsto \big(\ell^{\gamma_i}(m), \ell^{\gamma_i'}(m), \ell^{\gamma_i''}(m)\big)
\end{align*}
is injective.
\end{prop}
\begin{proof}
We compose $L$ with $\FN^{-1}$ and obtain a map
\begin{align*}
L\circ \FN^{-1}\colon \matR^{3g-3}_{>0} \times \matR^{3g-3} & \longrightarrow \matR_{>0}^{9g-9} \\
(l_i, \theta_i) \quad \ \ & \longmapsto (l_i, l_i', l_i'') 
\end{align*}
We prove that it is injective: it suffices to consider the case where the values $l_i$ are fixed and $\theta_i$ vary. Note that $\gamma_i'$ and $\gamma_i''$ intersect $\gamma_j$ if and only if $i=j$: hence $l_i'$ and $l_i''$ depend only on $\theta_i$ and not on the other torsion parameters $\theta_j$. Proposition \ref{strettamente:convessa:prop} says that $l_i' = f(\theta_i)$ is strictly convex and Proposition \ref{root:prop} gives $l_i'' = f(\theta_i+2\pi)$. 
A strictly convex proper function $f\colon\matR \to \matR$ is at most 2 to 1, hence the function 
\begin{align*}
\matR & \longrightarrow \matR \times \matR \\
\theta_i & \longmapsto \big(f(\theta_i), f(\theta_i + 2 \pi)\big)
\end{align*}
is injective. Therefore $L$ is injective. 
\end{proof}

\subsection{Collar lemma}
The thick-thin decomposition theorem implies that the closed geodesics of length smaller than a Margulis constant $\varepsilon_2$ on a complete hyperbolic surface are simple and have disjoint $R$-neighbourhoods (see Corollary \ref{simple:small:geodesics:cor}). When the curves are very short, one may choose $R$ to be very large: this fact is called the \emph{collar lemma} and we prove it directly using elementary tools.\index{collar lemma}

\begin{figure}
\begin{center}
\includegraphics[width = 11cm] {\iftoggle{BW}{collar-BW}{collar}}
\nota{Pick a segment of length $l$ and draw two perpendiculars at the endpoints: this determines a quadrilateral with two ideal vertices; let $f(l)$ be the distance between the opposite sides $l$ and $r$ (left). A geodesic pair-of-pants is the union of two isometric hexagons, which form in $\matH^2$ a right-angled octagon as drawn. The picture shows that the closed boundary geodesics $a$ and $b$ have disjoint $f(a)$ and $f(b)$-neighbourhoods, coloured here in \iftoggle{BW}{light grey}{yellow} (right).} 
\label{collar:fig}
\end{center}
\end{figure}

For any number $l>0$, draw the quadrilateral as in Figure \ref{collar:fig}-(left) and define $f(l)$ to be the distance between its opposite sides $l$ and $r$.

\begin{ex} The function $f\colon \matR_{>0}\to \matR_{>0}$ is strictly decreasing and a homeomorphism. In particular we get
$\lim_{l \to 0} f(l) = \infty$. 
Explicitly, we have 
$$\sinh f(l) = \frac{1}{\sinh \frac l2}.$$ 
\end{ex}
\begin{proof}[Hint]
Put $l$ in vertical position in $H^2$ and use Lemma \ref{cosh:cos:lemma}.
\end{proof}

The function $f$ is simple to define, and is particularly useful.

\begin{prop} \label{collar:pants:prop}
Let $P$ be a geodesic pair-of-pants with boundary lengths $a,b$, and $c$. The $f(a)$, $f(b)$, and $f(c)$-neighbourhoods of the boundary components form three disjoint collars.
\end{prop}
\begin{proof}
Consider two boundary components $a$ and $b$. 
The geodesic pair-of-pants $P$ is divided into two isometric hexagons, and we lift them to $\matH^2$ where they form a right-angled octagon as in Figure \ref{collar:fig}-(right). The picture shows that the $f(a)$ and $f(b)$-neighbourhoods of $a$ and $b$ are disjoint.
\end{proof}

Here is the collar lemma.

\begin{lemma}[Collar lemma] \label{collar:lemma}
Let $g\geqslant 2$ and $S_g$ have a hyperbolic metric. Disjoint simple closed geodesics $\gamma_1, \ldots, \gamma_k$ of length $l_1,\ldots, l_k$ have disjoint tubular $f(l_i)$-neighbourhoods. 
\end{lemma}
\begin{proof}
We may suppose (by adding more simple closed geodesics if necessary) that the closed geodesics form a pants decomposition, and it suffices to consider two curves that cobound the same pair-of-pants. Proposition \ref{collar:pants:prop} applies.
\end{proof}

We recall from Corollary \ref{simple:small:geodesics:cor} that every closed geodesic shorter than $\varepsilon_2$ is simple. Therefore very short closed geodesics are simple and have large disjoint tubular neighbourhoods: the more we shrink the curves, the larger are their neighbourhoods, and hence the larger is the diameter of the surface (the diameter of a metric space is the supremum of the distance of its points). In particular, if on a sequence of closed hyperbolic surfaces the injectivity radius tends to zero, their diameters must tend to infinity.

Among the many consequences of the collar lemma, we focus on a simple inequality which relates the geometric intersection of simple closed geodesics to their lengths. We denote by $L(\gamma)$ the length of $\gamma$.

\begin{cor}  \label{strizza:cor}
Let $g\geqslant 2$ and $S_g$ have a hyperbolic metric. Let $\gamma$ and $\eta$ be two simple closed geodesics in $S_g$. The following inequality holds:
$$L(\eta)\geqslant 2i(\eta, \gamma) \cdot f(L(\gamma)).$$
\end{cor}
\begin{proof}
The geodesic $\gamma$ has a tubular $f(L(\gamma))$-neighbourhood. The geodesic $\eta$ intersects $\gamma$ in $i(\eta,\gamma)$ points and hence crosses the tubular neighbourhood at least $i(\eta,\gamma)$ times, each with a segment of length $\geqslant 2f(L(\gamma))$.
\end{proof}

\subsection{A topology for the Teichm\"uller space} \label{Topology:Teich:subsection}
There are various equivalent ways to assign a topology to the Teichm\"uller space. 
On the torus $T$, we have seen that $\Teich(T)$ can be identified with the hyperbolic plane $\matH^2$, and the mapping class group acts as isometries on it: we could not hope for a better picture of $\Teich(T)$ and we are fully satisfied.

When $g\geqslant 2$ we could similarly use the Fenchel-Nielsen coordinates and give $\Teich(S_g)$ the topology of $\matR^{6g-6}$, but then to be honest we should also check that the topology does not depend on the frame... we prefer to equip the Teichm\"uller space with an intrinsic topology and then prove that the Fenchel-Nielsen coordinates are homeomorphisms.\index{$\calS$} 

We indicate by $\calS= \calS(S_g)$ the set of all the non-trivial simple closed curves in $S_g$, considered up to isotopy and orientation reversal (we say that the curves are \emph{unoriented}). Each element $\gamma\in\calS$ induces a length function
$$\ell^\gamma\colon \Teich(S_g) \longrightarrow \matR_{>0}.$$
We indicate as usual with $\matR^\calS$ the set of all functions $\calS \to \matR$ and give it the usual product topology (the weakest one such that all the projections are continuous). The natural map
\begin{align*}
\Teich(S_g) & \longrightarrow \matR^\calS \\
m & \longmapsto \big( \gamma \longmapsto \ell^\gamma(m) \big)
\end{align*}
is injective by Propositions \ref{3:prop} and \ref{9g-9:prop}. We may hence consider $\Teich(S_g)$ as a subspace of $\matR^\calS$ and assign it the subspace topology. This topology on $\Teich(S_g)$ is the weakest one where the length functions $\ell^\gamma$ are continuous.
Recall that a topological space is \emph{second-countable} if it has a countable base.

\begin{prop}
The space $\matR^\calS$ is Hausdorff and second-countable.
\end{prop}
\begin{proof}
Every product of Hausdorff spaces is Hausdorff, and every countable product of second-countable spaces is second-countable.
\end{proof}

We recall the following topological fact.

\begin{prop} \label{propria:chiusa:prop}
Let $f\colon X \to Y$ be a continuous and proper map between topological spaces. If $Y$ is Hausdorff and second-countable then $f$ is closed.
\end{prop}

Proper maps onto reasonable spaces are closed. If they are also injective, we can obtain more.

\begin{cor} \label{propria:chiusa:cor}
Let $f\colon X \to Y$ be a continuous, proper, and injective map between topological spaces. If $Y$ is Hausdorff and second-countable then $f$ is a homeomorphism onto its image. 
\end{cor}

We will use this corollary in a moment.
Recall that every isometry $\varphi\in \Iso^+(H^2) = \PSLR$ has a trace $\tr \varphi$ defined only up to sign, whose modulus is $>2$ precisely when $\varphi$ is hyperbolic, see Proposition \ref{PSLR:prop}.

\begin{prop} \label{trace:prop}
Let $S = \matH^2/_\Gamma$ be an orientable hyperbolic surface. Every hyperbolic transformation $\varphi\in\Gamma$ produces a closed geodesic $\gamma$ in $S$
with
$$|\tr\varphi| = 2\cosh \frac {L(\gamma)} 2.$$
\end{prop}
\begin{proof}
Up to conjugacy we have $\varphi(z) = e^{L(\gamma)} z$. The matrix is 
$$\varphi = \left(\begin{matrix} e^{\frac{L(\gamma)}2} & 0 \\ 0 & e^{-\frac{L(\gamma)}2} \end{matrix}\right)$$
hence $|\tr\varphi| = 2\cosh \frac{L(\gamma)}2$.
\end{proof}

In particular, the length of $\gamma$ depends continuously on the transformation $\varphi$. We will use this fact to prove the following. We suppose $g\geqslant 2$.

\begin{prop} \label{FN:omeom:prop}
The Fenchel-Nielsen map 
$$\FN\colon \Teich(S_g) \longrightarrow \matR_{>0}^{3g-3}\times\matR^{3g-3}$$
is a homeomorphism.
\end{prop}
\begin{proof}
We consider $\Teich(S_g)$ inside $\matR^\calS$ and examine the inverse map
$$\FN^{-1}\colon\matR_{>0}^{3g-3}\times\matR^{3g-3} \longrightarrow \matR^\calS.$$

We prove that $\FN^{-1}$ is continuous. The map $\FN^{-1}$ assigns to the parameters $(l_i, \theta_i)$ a metric on $S_g$ constructed by attaching right-angled hexagons. Both the hexagons and the attaching maps depend continuously on the parameters $(l_i, \theta_i)$ and lift to a tessellation of $\matH^2$ into hexagons. Since the decomposition into hexagons varies continuously, its deck transformations vary continuously in $\PSLR$ and hence the length functions too by Proposition \ref{trace:prop}. Therefore $\FN^{-1}$ is continuous.

We prove that $\FN^{-1}$ is proper. Take a diverging sequence of parameters $(l_i,\theta_i)$ (that is, without converging subsequences) in $\matR_{>0}^{3g-3}\times \matR^{3g-3}$: we need to show that its image is also a diverging subsequence. The thesis is equivalent to show that the length function of some curve goes to infinity. If $l_i\to +\infty$ for some $i$ we are done. If $l_i\to 0$, the length of any curve intersecting essentially the shrinking curve $\gamma_i$ goes to infinity by Corollary \ref{strizza:cor}. It remains to consider the case where the length parameters $l_i$ converge to some non-zero value, but some twist parameter $\theta_j$ goes to infinity: in that case the length of any curve intersecting the twisted curve $\gamma_j$ goes to infinity by Proposition \ref{strettamente:convessa:prop}.

Finally, the map $\FN^{-1}$ is a homeomorphism onto its image by Corollary \ref{propria:chiusa:cor}. The proof is complete.
\end{proof}

During the proof we have also discovered the following.

\begin{prop} \label{diverges:infinity:prop}
If a sequence $m_i\in\Teich(S_g)$ diverges, there is a $\gamma \in \calS$ such that $\ell^\gamma(m_i) \to \infty$ on a subsequence.
\end{prop}

Recall that the action of a topological group $G$ on a topological space $X$ is \emph{continuous} if the action map $G\times X \to X$ is continuous. This implies that $G$ acts on $X$ by homeomorphisms. We give the mapping class group $\MCG(S_g)$ the discrete topology.

\begin{prop}
The earthquakes and mapping class group actions on the Teichm\"uller space are continuous.
\end{prop}
\begin{proof}
The mapping class group acts on $\calS$ by permutations, hence its action on $\matR^\calS$ is continuous. On Fenchel-Nielsen coordinates the earthquake action sends $\theta$ to $\theta+\theta'$ and is hence continuous.
\end{proof}

The immersion in $\matR^{9g-9}$ is also a topological embedding.

\begin{prop}
The injective representation $\Teich(S_g) \hookrightarrow \matR^{9g-9}$
furnished by Proposition \ref{9g-9:prop} is a homeomorphism onto its image.
\end{prop}
\begin{proof}
Using Fenchel-Nielsen coordinates the map is clearly continuous. The proof that it is proper is as in Proposition \ref{FN:omeom:prop}.
\end{proof}

\subsection{Surfaces of finite type} \label{Teich:finite:type:subsection}
We have considered only closed surfaces $S_g$ for simplicity, but most of the arguments exposed in this chapter extend easily to all surfaces $S_{g,b,p}$ of finite type with negative Euler characteristic. 

The Teichm\"uller space $\Teich(S_{g,b,p})$ is the set of all the complete hyperbolic metrics with geodesic boundary, considered up to isometries that are isotopic to the identity. Fenchel-Nielsen coordinates are defined analogously: the surface decomposes into $-\chi(S_{g,b,p}) = 2g+b+p-2$ pairs-of-pants, and the interior curves of the decomposition are 
$$\frac 12 \big(3(2g+b+p-2) - b - p \big) = 3g+b+p-3.$$
The Fenchel-Nielsen coordinates are
$$\big(l_1, \ldots, l_{3g+b+p-3}, l_1^\partial, \ldots, l_b^\partial, \theta_1,\ldots, \theta_{3g+b+p-3}\big)$$
where the $l_i$ and $l_j^\partial$ are the length parameters of the $3g+b+p-3$ interior and $b$ boundary curves, and the $\theta_i$ are the torsion angles of the interior curves. We get a bijection
$$\FN\colon \Teich (S_{g,b,p}) \longrightarrow \matR_{>0}^{3g+2b+p-3} \times \matR^{3g+b+p-3}$$
which is a homeomorphism with respect to the natural topology on Teichm\"uller space as a subset of $\matR^\calS$. Therefore the Teichm\"uller space is homeomorphic to a ball of dimension $-3\chi(S_{g,b,p})-p$.

For instance, the Teichm\"uller space of a pair-of-pants is $\matR_{>0}^3$ parametrized by the lengths of the boundary geodesics, while that of a thrice-punctured sphere is a point (there is a unique metric).

An alternative description of the Teichm\"uller space for punctured surfaces, with ideal triangles playing the role of pairs-of-pants, is described in the next section.

\section{Shear coordinates}
The Teichm\"uller space of a punctured surface may also be parametrized using ideal triangles instead of pairs-of-pants: this viewpoint is maybe a bit simpler, and generalises successfully to dimension three (via ideal tetrahedra). The coordinates that it produces are called \emph{shear coordinates}.\index{shear coordinates}

\subsection{Ideal triangulations}
In dimension two and three it is customary to relax the definition of \emph{triangulation}, originally restricted to simplicial complexes, see Section \ref{triangulations:subsection}. We prefer to define a triangulation in a looser sense, as a finite set of triangles glued together by pairing their edges. 

Let $\Delta_1,\ldots, \Delta_{2k}$ be an even number of identical copies of the standard oriented 2-simplex. A \emph{triangulation} $\calT$ is a partition of the $6k$ edges of the triangles into $3k$ pairs, and for each pair a simplicial isometry between the two edges. The triangulation is \emph{oriented} if the simplicial isometries are orientation-reversing. If we glue the triangles along the isometries we get a compact surface $S$: we always suppose that $S$ is connected and $\calT$ is oriented, hence $S = S_g$ for some $g\geqslant 0$.

The surface $S$ is triangulated with $2k$ triangles, $3k$ edges, and some $p$ vertices. Vertices, edges, and triangles form a cellularisation of $S$, but not necessarily a simplicial complex: for instance two or three vertices of the same triangle $\Delta_i$ can be identified to a single one along the process, as the following exercise shows.

\begin{ex}
Construct a triangulation of the torus with one vertex, three edges, and two triangles.
\end{ex}

Let $\Sigma$ be the non-compact surface obtained by removing the $p$ vertices of the triangulation $\calT$ from $S$: we say that $\calT$ is an \emph{ideal triangulation} for $\Sigma$. The surface $\Sigma$ is a \emph{punctured surface}, \emph{i.e.}~$\Sigma = S_{g,0,p}$ with $p\geqslant 1$.\index{ideal triangulation}

\begin{figure}
\begin{center}
\includegraphics[width = 9 cm] {\iftoggle{BW}{triangulated_surface-BW}{triangulated_surface}}
\nota{The standard representation of $S_g$ by identifying the edges of a $4g$-gon has only one vertex; if we triangulate the $4g$-gon we get an ideal triangulation for $S_{g,0,1}$ (left). A move that modifies a triangulation by adding a new vertex (right).}
\label{triangulated_surface:fig}
\end{center}
\end{figure}

\begin{prop} We have $\chi(\Sigma) = -k <0$. Every punctured surface of negative Euler characteristic has an ideal triangulation.
\end{prop}
\begin{proof}
We have $\chi(\Sigma) = \chi(S) -p = p-3k+2k-p = -k$. An ideal triangulation for $S_{0,0,3}$ is constructed by attaching $\Delta_1$ and $\Delta_2$ via the obvious maps. An ideal triangulation for $S_{g,0,1}$ is constructed in Figure \ref{triangulated_surface:fig}-(left) for all $g\geqslant 1$. An ideal triangulation for $S_{g,0,p}$ with $\chi(S_{g,0,p})<0$ and $p\geqslant 2$ is obtained from one of these by increasing $p$ as in Figure \ref{triangulated_surface:fig}-(right).
\end{proof}

\subsection{Hyperbolic ideal triangulations}
Let $\calT$ be an oriented triangulation with triangles
$\Delta_1,\ldots, \Delta_{2k}$. We substitute every $\Delta_i$ with an ideal hyperbolic triangle (recall that it is unique up to isometry) and pair the edges with some orientation-reversing isometries. The resulting oriented punctured surface $\Sigma$ now inherits a hyperbolic structure of area $k\pi$, since an ideal triangle has area $\pi$. We call $\calT$ a \emph{hyperbolic ideal triangulation} for $\Sigma$.

\begin{figure}
\begin{center}
\includegraphics[width = 9 cm] {\iftoggle{BW}{triangoli_ideali-BW}{triangoli_ideali}}
\nota{On an ideal triangle every edge has a unique perpendicular which ends at the opposite vertex. The base of this perpendicular is the \emph{midpoint} of the edge and the three perpendiculars intersect in a point called \emph{barycenter} (left). There is a unique configuration of three pairwise tangent horocycles centred at the ideal vertices, and their tangency points are the midpoints of the edges (right).}
\label{triangoli_ideali:fig}
\end{center}
\end{figure}

\begin{figure}
\begin{center}
\includegraphics[width = 11 cm] {\iftoggle{BW}{shear-BW}{shear}}
\nota{The \emph{shear} is the signed distance $d$ between two midpoints after the gluing. By convention the sign is positive $d>0$ if the points are as shown here, that is an observer that arrives at the midpoint of an edge sees the other midpoint on its left: we need the orientation of $\calT$ here (left). We develop a triangulation around a puncture $v$: here $h=3$ (right).}
\label{shear:fig}
\end{center}
\end{figure}

The edges are infinite lines and the gluing isometry is not unique: indeed there is a 1-parameter family of isometries to choose from at every pair of edges, and the hyperbolic structure depends on that -- although the topology of $\Sigma$ does not. To encode this dependence, we note that every edge of an ideal triangle has a preferred \emph{midpoint} defined in Figure \ref{triangoli_ideali:fig}. The orientation-reversing isometry between two edges of two triangles is determined by the signed distance $d$ (called \emph{shear}) of their two midpoints as shown in Figure \ref{shear:fig}-(left). The hyperbolic structure on $\Sigma$ is determined by the \emph{shear coordinates} $d=(d_1, \ldots, d_{3k})$ of the $3k$ edges of the triangulation.

\subsection{Complete solutions}
The hyperbolic structure on $\Sigma$ may be incomplete! Recall that a hyperbolic surface is complete if and only if its universal cover is $\matH^2$, so we now look at the universal cover of $\Sigma$. 

\begin{figure}
\begin{center}
\includegraphics[width = 4cm] {\iftoggle{BW}{universale-BW}{universale}}
\nota{The universal covering of $\Sigma$ is isometric to the interior of a convex set delimited by some lines in $\matH^2$, which is the intersection of (possibly infinitely many) half-planes bounded by disjoint lines.}
\label{universal:fig}
\end{center}
\end{figure}

\begin{prop}
The universal cover $\tilde \Sigma$ of $\Sigma$ is isometric to
the interior of the intersection of some half-planes in $\matH^2$ with disjoint boundary lines (see Figure \ref{universal:fig}). The surface $\Sigma$ is complete if and only if $\tilde \Sigma \isom\matH^2$.
\end{prop}

\begin{proof}
Recall from Section \ref{non-complete:subsection} that there is a developing map $D\colon \tilde {\Sigma} \to \matH^2$. We prove that $D$ is injective. 

We can construct $D$ as follows: the hyperbolic ideal triangulation $\calT$ of $\Sigma$ lifts to a hyperbolic ideal triangulation $\tilde \calT$ of $\tilde{\Sigma}$ with infinitely many triangles. Send a triangle of $\tilde \calT$ to an arbitrary ideal triangle of $\matH^2$ and then develop the map $D$ by attaching subsequently all the triangles of $\tilde \calT$.  At each step the image of $D$ is an ideal polygon and we attach a new triangle to the exterior of some side of the polygon. Therefore at each step $D$ is injective and hence $D$ is globally injective.

Since $D$ is injective we identify $\tilde{\Sigma}$ with its open image in $\matH^2$. Pick a point $x\in\partial\tilde{\Sigma}$. Every neighbourhood of $x$ must intersect infinitely many triangles of $\tilde \calT$, hence there is a sequence of edges $e_i$ of $\tilde \calT$ such that $d(e_i,x)$ is monotone decreasing and tends to zero. Since the $e_i$ are disjoint, their endpoints in $\partial \matH^2$ are monotone and tend to two points, which determine a line $l$. We must have $d(l,x)=0$ and hence $x\in l$ and $l\subset \partial \tilde \Sigma$. 

We have discovered that $\partial \tilde \Sigma$ consists of disjoint lines. This proves that the closure of $\tilde\Sigma$ is the intersection of half-planes. A surface $\Sigma$ is complete if and only if $\tilde \Sigma$ is complete, and hence isometric to $\matH^2$.
\end{proof}

The hyperbolic structure on $\Sigma$ depends on the shearing coordinates $d\in \matR^{3k}$, and we now want to understand which values of $d$ produce a complete metric. An ideal vertex $v$ of the triangulation $\calT$ is adjacent to some $h$ triangles, which we denote by $\Delta_{1},\ldots, \Delta_{h}$ for simplicity although there might be repetitions, and to $h$ edges with some shearing coordinates that we also indicate by $d_1,\ldots, d_h$ for simplicity. The vertex $v$ is a puncture of $\Sigma$ and we let $N(v)$ be a small topological punctured closed disc around $v$.

\begin{prop}
The punctured disc $N(v)$ is complete if and only if $d_1 + \ldots + d_h = 0$.
\end{prop}
\begin{proof}
We construct a portion of the developing map $D\colon \tilde \Sigma \to \matH^2$ as follows: we use the half-plane model and send $\Delta_1$ to $\matH^2$ with $v=\infty$ as in Figure \ref{shear:fig}-(right), then we develop the triangulation horizontally. In the picture we have $h=3$ so we reproduce $\Delta_1, \Delta_2, \Delta_3$ and then another copy $\Delta_1'$ of $\Delta_1$. Let $\varphi\in\PSLR$ be the holonomy isometry that sends $\Delta_1$ to $\Delta_1'$.

If $d_1+\ldots + d_h=0$ the midpoints of $\Delta_1$ and $\Delta_1'$ are at the same height as in Figure \ref{shear:fig}-(right) and therefore $\varphi$ is a parabolic map $\varphi\colon z \mapsto z+b$ for some $b\in \matR$. Therefore $N(v)$ is a truncation of the cusp $\matH^2/_{\langle \varphi\rangle }$, which is complete.

\begin{figure}
\begin{center}
\includegraphics[width = 11 cm] {\iftoggle{BW}{shear2-BW}{shear2}}
\nota{When $d_1+\ldots + d_h \neq 0$, the isometry $\varphi$ sending $\Delta_1$ to $\Delta_1'$ is hyperbolic, with some axis $l$ which (up to translations) has endpoints $0$ and $\infty$. Hence $\varphi(z)=e^d z$ and the triangulation develops towards $l$ as shown here.}
\label{shear2:fig}
\end{center}
\end{figure}

If $d_1+\ldots + d_h \neq 0$ the midpoints of $\Delta_1$ and $\Delta_1'$ are at different heights and therefore $\varphi$ is not parabolic. Since $\infty \in \Fix(\varphi)$, the isometry $\varphi$ is a hyperbolic transformation having some other fixed point in $\matR$. Up to translating everything we may suppose that $\Fix(\varphi)=\{0,\infty\}$, hence the axis $l$ of $\varphi$ is the vertical coordinate axis, and $\varphi(z) =  e^d z$ with $d = d_1+\ldots + d_h$. The map $D$ develops towards $l$ as shown in Figure \ref{shear2:fig}: we get $l\subset \partial \tilde\Sigma$, so $\tilde\Sigma \neq\matH^2$ is not complete. 
\end{proof}
 
Summing up, we have a \emph{completeness equation} of type $d_{i_1}+\ldots + d_{i_h}=0$ for each of the $p$ ideal vertices of the triangulation.

\begin{cor}
The hyperbolic structure on $\Sigma$ is complete if and only if $d$ satisfies the $p$ completeness equations.
\end{cor}
\begin{proof}
Let $N(v)$ be a small punctured closed disc at $v$ for every ideal vertex $v$. The closure of $\Sigma \setminus \cup_v N(v)$ is compact. Therefore $\Sigma$ is complete if and only if each $N(v)$ is.
\end{proof} 

\subsection{Shear coordinates for Teichm\"uller space}
The solution space of the completeness equations is some linear subspace $V \subset \matR^{3k}$ of dimension at least $3k-p = -3\chi(\Sigma)-p$. Note that this is precisely the dimension of $\Teich(\Sigma)$, see Section \ref{Teich:finite:type:subsection}. We have constructed a map
$${\rm sh}\colon V \longrightarrow \Teich(\Sigma).$$
We will prove below that $\sh$ is a homeomorphism and hence $\dim V = 3k-p$. We will use the following.

\begin{prop} \label{straightens:prop}
Let $\Sigma$ be a complete hyperbolic punctured surface.
Every ideal triangulation is isotopic to a unique hyperbolic ideal triangulation. 
\end{prop}
\begin{proof}
We describe a procedure that \emph{straightens} every ideal triangulation $\calT$ of $\Sigma = \matH^2/_\Gamma$ to a hyperbolic one. 

Truncate the cusps of $\Sigma$ to get a compact sub-surface. The ideal triangulation $\calT$ is now truncated to a decomposition into hexagons, each with three boundary and three interior edges. We lift the truncated cusps and the hexagons to $\matH^2$. Truncated cusps lift to horoballs, bounded by horocycles, and hexagons lift to hexagons (the hexagons are not geodesic, only smooth). 

An interior edge of a hexagon in $\matH^2$ cannot have both its endpoints in the same horocycle, since it would form a bigon with it, and a bigon does not decompose into hexagons (by a Euler characteristic argument). Therefore every interior edge connects two distinct horocycles, centred at two distinct points of $\partial \matH^2$. Replace this edge with the geodesic line connecting these two ideal points.

If we do this at every interior edge we get a $\Gamma$-invariant hyperbolic ideal triangulation of $\matH^2$, which projects to a hyperbolic ideal triangulation for $\Sigma$ isotopic to $\calT$.  
\end{proof}

We can now parametrize $\Teich(\Sigma)$.
 
\begin{cor} The map $\sh$ is a bijection.
\end{cor}
\begin{proof}
Proposition \ref{straightens:prop} shows that in every metric the initial triangulation $\calT$ straightens to an ideal hyperbolic one and is hence realised by some $d$, so $\sh$ is surjective. Moreover the straightened triangulation is unique and the parameters $d_i$ are obtained intrinsically from it, so $\sh$ is also injective.
\end{proof}

\subsection{Incomplete metrics} \label{incomplete:metrics:subsection}
The shear coordinates $d$ may define an incomplete metric on $\Sigma$, and we now determine its metric completion $\overline \Sigma$.

Recall that every ideal vertex $v$ has a \emph{total shear} $d_v = d_{i_1}+\ldots +d_{i_h}$ which is zero precisely when $\Sigma$ is complete near $v$. 

\begin{figure}
\begin{center}
\includegraphics[width = 11 cm] {\iftoggle{BW}{shear3-BW}{shear3}}
\nota{Let $s$ be a Euclidean line at hyperbolic distance $R$ from $l$. If $R$ is sufficiently small, the line $s$ does not intersect the lower edges of $\Delta_1, \Delta_2, \Delta_3$. Let $S$ be the region lying above $s$.}
\label{shear3:fig}
\end{center}
\end{figure}

\begin{prop}  \label{completion:surface:prop}
The metric completion $\overline\Sigma$ is a hyperbolic surface with geodesic boundary, whose boundary lengths are the absolute values of the total shears of the vertices. 
\end{prop}
\begin{proof}
Let $N(v)$ be a small closed punctured disc around $v$.
If $d_v=0$ we already know that $N(v)$ is a truncated cusp. If $d_v\neq 0$ we develop $N(v)$ in $\matH^2$ as in Figure \ref{shear3:fig}. The hyperbolic transformation is $\varphi(z) = e^{d_v} z$ and the points at distance $R$ from $l$ form a Euclidean line $s$ starting from $0$. Pick $R$ small so that $s$ does not intersect the lower edges of $\Delta_1,\ldots,\Delta_h$. Since $s$ is $\varphi$-invariant, it is contained in the union of all the developed triangles of Figure \ref{shear3:fig} and hence in $\tilde\Sigma$.

The line $s$ projects to a (non-geodesic) circle in $\Sigma$ around $v$, and we can suppose that $\partial N(v)$ is that circle. We have $N(v) = S/_\gamma$ where $S\subset H^2$ is the region lying above $s$. The completion $\overline S $ equals $S\cup l$ and therefore $\overline{N(v)} = \overline S/_\gamma$ is an annulus with two boundary components: the original non-geodesic $\partial N(v)$ and a new geodesic $l/_\gamma$, a curve of length $|d_v|$.
\end{proof}

\begin{figure}
\begin{center}
\includegraphics[width = 10 cm] {\iftoggle{BW}{shear_completion2-BW}{shear_completion}}
\nota{Near an ideal vertex $v$, horocycles (in \iftoggle{BW}{black}{blue}) and edges of the triangulation (in \iftoggle{BW}{grey}{red}) are mutually orthogonal. They behave much differently in presence of a cusp (left) or of a geodesic boundary (right): in a cusp, horocycles form parallel circles and edges point toward the cusp; near a geodesic boundary, the horocycle point orthogonally toward the boundary and the edges spin and wind around it indefinitely.}
\label{shear_completion:fig}
\end{center}
\end{figure}

If $d_v\neq 0$, the edges of the triangulation pointing towards $v$ spin around the geodesic boundary as shown in Figure \ref{shear_completion:fig}-(right). The spinning direction induces an orientation on each component of $\partial \overline \Sigma$, which depends on the sign of $d_v$. 

Summing up, every shear coordinate $d = (d_1,\ldots, d_{3k})\in \matR^{3k}$ determines a hyperbolic surface $\overline{\Sigma}$ with some $p'$ cusps and $b$ oriented geodesic boundary, where $b+p'=p$ is the number of vertices of the ideal triangulation. The surface $\overline \Sigma$ is diffeomorphic to $S_{g,b,p'}$. Let $\Teich^{\rm or}(S_{g,b,p'})$ denote the Teichm\"uller space of the surface $S_{g,b,p'}$ equipped with an orientation of the $b$ boundary components: of course this is just $2^b$ copies of $\Teich(S_{g,b,p'})$.

\begin{prop} \label{straighten:incomplete:prop}
The shear coordinates induce a bijection
$$\matR^{3k} \longleftrightarrow \bigsqcup_{b+p' = p} \Teich^{\rm or}(S_{g,b,p'}).$$
\end{prop}
\begin{proof}
The map sends $d\in\matR^{3k}$ to the complete hyperbolic $\overline \Sigma$, and we now exhibit an inverse. Let $\overline\Sigma$ be a surface of finite type with a hyperbolic structure with $p'$ cusps and $b$ oriented geodesic boundary components, and $\Sigma = \interior{\overline\Sigma}$. We may straighten any ideal triangulation $\calT$ of $\Sigma$ similarly as we did in Proposition \ref{straightens:prop}, with the only difference that horocycles are replaced by \emph{oriented} geodesic lines as lifts of oriented geodesic boundary components; in the construction, we use the final endpoint of the oriented line instead of the ideal point of the horocycle. We get an ideal triangulation of $\Sigma$ whose completion is $\overline \Sigma$.
\end{proof}

For instance, $\matR^3$ parametrizes altogether all the hyperbolic metrics on the pair of pants, where each boundary component becomes either an oriented geodesic or ``degenerates'' to a cusp.

\subsection{References}
Most of the material presented here is standard and can be found in many books, starting from Thurston's notes \cite{Th}. The main reference is Farb -- Margalit \cite{FM}, which also contains the proof of Proposition \ref{strettamente:convessa:prop}, that appeared originally in a paper of Bestvina, Bromberg, Fujiwara, and Souto \cite{BBFS}. The proof of the Collar Lemma is taken from Hubbard \cite{Hu}.

%% file: Automorphisms.tex
\chapter{Surface diffeomorphisms} \label{automorfismi:chapter} \label{automorphisms:chapter}

We describe in this chapter an analogy between the hyperbolic space $\matH^n$ and the Teichm\"uller space $\Teich(S_g)$ of a closed orientable surface $S_g$ of genus $g\geqslant 1$. The theory, originated from Thurston in the late 1970s, provides a beautiful and powerful framework for the analysis of the geometric and dynamical properties of the surfaces $S_g$ and of their mapping class group $\MCG(S_g)$. 

Here is a quick sketch of this analogy. We have already seen that $\matH^n$ compactifies to a closed disc $\overline{\matH^n}$, that $\Iso(\matH^n)$ acts on it, that by Brouwer's theorem every non-trivial isometry $\varphi$ has a fixed point in $\overline{\matH^n}$, and we have called $\varphi$ \emph{elliptic}, \emph{parabolic}, or \emph{hyperbolic} according to the position of its fixed points.

We construct in this chapter a similar compactification of the open ball $\Teich(S_g)$ to a closed disc. The action of the mapping class group $\MCG(S_g)$ extends to this closed disc, and by Brouwer's theorem every non-trivial element $\varphi$ of $\MCG(S_g)$ has a fixed point there. According to the position of the fixed points of $\varphi$, we say that $\varphi$ is \emph{finite order}, \emph{reducible}, or \emph{pseudo-Anosov}. 

If $\varphi$ is an isometry for some hyperbolic structure of $S_g$, it belongs to the first type. Dehn twists belong to the second. The pseudo-Anosov maps are both the most mysterious and the most important, and for their study we need to introduce a wealth of beautiful new technology: \emph{geodesic currents}, \emph{laminations}, and \emph{train tracks}.

\section{Thurston's compactification}
Let $S_g$ be a surface of genus $g\geqslant 2$.
Recall that $\Teich(S_g)$ is the Teichm\"uller space of $S_g$ and $\calS = \calS(S_g)$ is the set of all non-trivial simple closed curves in $S_g$, considered up to isotopy and orientation reversal: all the closed curves considered in this chapter are \emph{unoriented}. 

Although $\Teich(S_g)$ and $\calS$ are very different in nature, we want to compactify both spaces by embedding them in a single bigger space. The model that we have in mind, and that we would like to extend to higher genus surfaces, is the flat torus picture that was painted in Section \ref{length:torus:subsection}. 

We briefly summarise it. The Teichm\"uller space $\Teich(T)$ of the torus $T$ is $\matH^2$ and it compactifies to $\overline {\matH^2}$; the mapping class group $\MCG(T)$ acts on $\overline {\matH^2}$ as a discrete group of isometries; as such, the non-trivial elements of $\MCG(T)$ are divided into three classes (hyperbolic, parabolic, elliptic); the simple closed curves form a dense countable set of rational points in $\partial \matH^2 = \matR\cup \{\infty\}$.

We now may wonder whether there is an analogous identification between $\Teich(S_g)$ and $\matH^{6g-6}$ that transforms $\Teich(S_g)$ into a discrete subgroup of $\Iso(\matH^n)$. This is unfortunately not the case, and the compactification of $\Teich(S_g)$ must be constructed from scratch via different methods: we compactify $\Teich(S_g)$ by embedding it in a bigger infinite-dimensional space.

\subsection{Projective immersion} \label{projective:immersion:subsection}
In Chapter  \ref{Teichmuller:chapter} we have used the length functions to construct an embedding
$$i\colon \Teich(S_g) \hookrightarrow \matR^\calS.$$
We know that $\Teich(S_g)$ is homeomorphic to an open ball of dimension $6g-6$, and we want to compactify it in a geometrically meaningful way. A first tentative could be to take its closure in $\matR^\calS$, but this does not work:

\begin{prop}
The subspace $i(\Teich(S_g))$ is closed in $\matR^\calS$.
\end{prop}
\begin{proof}
The inclusion is proper, hence closed by Proposition \ref{propria:chiusa:prop}.
\end{proof}

We are apparently stuck, so we turn back to our model, hyperbolic space, to get some inspiration. We recall that $\matH^n$ is properly embedded in the lorentzian $\matR^{n+1}$ as a hyperboloid $I^n$. To compactify $\matH^n$, we may consider its image $K^n$ in $\matRP^n$ (the Klein model) and take the closure $\overline{K^n}$ there.

We try to mimic this construction, by considering the infinite-dimensional projective space $\matP(\matR^\calS)$ with the projection 
$$\pi\colon \matR^\calS\setminus \{0\} \longrightarrow \matP(\matR^\calS).$$
We first need to check that $\Teich(S_g)$ embeds there.
\begin{prop}
The composition
$$\pi\circ i \colon \Teich(S_g) \longrightarrow \matP(\matR^\calS)$$
is injective.
\end{prop}
\begin{proof}
Suppose by contradiction that there are two distinct points $m, m' \in \Teich(S_g)$ with $\pi(i(m)) = \pi(i(m'))$; this implies that there is a constant $k>1$ such that $\ell^\gamma(m) = k\cdot\ell^\gamma(m')$ for all $\gamma\in\calS$. That sounds very unlikely, and we now prove that it easily leads to a contradiction.

Let $\gamma_1, \gamma_2 \in \calS$ be two curves with $i(\gamma_1,\gamma_2)=1$. We take $x_0 = \gamma_1 \cap \gamma_2$ as a basepoint for $\pi_1(S_g, x_0)$ and note that the elements $\gamma_2*\gamma_1$ and $\gamma_2*\gamma_1^{-1}$ are represented by two more non-trivial simple closed curves in $S_g$. It is easily checked that the formula
$$\tr(A)\cdot \tr(B) = \tr(AB) + \tr(A^{-1}B)$$
holds for any $A,B \in \SLR$. Proposition \ref{trace:prop} implies that
$$2\cosh \left(\frac{L(\gamma_1)}2 \right) \cdot \cosh\left(\frac{L(\gamma_2)}2 \right) = 
\cosh\left(\frac{L(\gamma_2*\gamma_1)}2 \right) + \cosh\left(\frac{L(\gamma_2*\gamma_1^{-1})}2 \right).$$
We have obtained a relation between the lengths of $\gamma_1, \gamma_2$, $\gamma_2* \gamma_1$, and $\gamma_2*\gamma_1^{-1}$ that holds for any hyperbolic metric on $S_g$. It may be rewritten as: 
$$ \cosh\left(\tfrac{L(\gamma_1) + L(\gamma_2)}2 \right) + \cosh\left(\tfrac{L(\gamma_1) - L(\gamma_2)}2 \right) = \cosh\left(\tfrac{L(\gamma_2*\gamma_1)}2 \right) + \cosh\left(\tfrac{L(\gamma_2*\gamma_1^{-1})}2 \right).
$$
By contradiction every $m'$-length is $k$ times a $m$-length: this equation is hence valid after multiplying every argument by $k$. It is easy to check that 
$$\cosh a + \cosh b = \cosh c + \cosh d, \quad \cosh ka + \cosh kb = \cosh kc + \cosh kd$$ 
if and only if $\{a,b\} = \{c,d\}$. This leads to a contradiction: the number $L(\gamma_1)+L(\gamma_2)$ is strictly bigger than $L(\gamma_2*\gamma_1)$ or $L(\gamma_2*\gamma_1^{-1})$, since $\gamma_2*\gamma_1$ and $\gamma_2*\gamma_1^{-1}$ have a non-geodesic representative of length $L(\gamma_1)+L(\gamma_2)$.
\end{proof}

We have embedded $\Teich(S_g)$ in $\matP(\matR^\calS)$, and we now turn to the set $\calS$.

\subsection{Thurston's compactification}
We now embed $\calS$ in $\matP(\matR^\calS)$. A simple closed curve $\gamma \in \calS$ defines a functional $i(\gamma) \in \matR^\calS$ by setting:
$$i(\gamma)(\eta) = i(\gamma, \eta).$$
We have constructed a map $i\colon \calS \to \matR^\calS$.
\begin{prop}
The composition
$$\pi \circ i\colon\calS \longrightarrow \matP(\matR^\calS)$$
is injective.
\end{prop}
\begin{proof}
Let $\gamma_1, \gamma_2\in \calS$ be distinct. There is always a curve $\eta\in\calS$ with $i(\gamma_1,\eta)\neq 0$ and $i(\gamma_2,\eta)=0$. (If $i(\gamma_1,\gamma_2)>0$, simply take $\eta=\gamma_2$. Otherwise, it is an easy exercise.)
\end{proof}

We will now tacitly consider both $\Teich(S_g)$ and $\calS$ as subsets of $\matR^\calS$.
\begin{prop}
The subsets $\Teich(S_g)$ and $\calS$ are disjoint in $\matP(\matR^\calS)$.
\end{prop}
\begin{proof}
For each $\gamma\in\calS$ we have $i(\gamma,\gamma)=0$, while every curve has positive length on any hyperbolic metric.
\end{proof}

We can now state Thurston's compactification theorem.\index{Thurston's compactification theorem} 

\begin{teo} \label{compattificazione:teo}
The closure $\overline{\Teich(S_g)}$ of $\Teich(S_g)$ in $\matP(\matR^\calS)$ is homeomorphic to the closed disc $D^{6g-6}$. Its interior is $\Teich(S_g)$ and its boundary sphere contains $\calS$ as a dense subset.
\end{teo}
In particular, the closure of $\calS$ is homeomorphic to a sphere $S^{6g-7}$. 
The proof of Theorem \ref{compattificazione:teo} occupies most of this chapter and will be completed in Section \ref{MCG:conclusion:subsection}. We will introduce in the process various new geometric objects that play an important role in the topology of manifolds in dimension two and three: \emph{geodesic currents}, \emph{laminations}, and \emph{train tracks}.

For the moment we content ourselves with checking that this projective embedding strategy works at least on the much simpler flat torus case.

\subsection{The torus}
The Teichm\"uller space of the torus $T$ is described in Section \ref{length:torus:subsection} and needs no further comment; nevertheless, we prove here that the projective embedding strategy works for $T$, just as a sanity check before approaching the more complicated higher genus surfaces.

\begin{prop} \label{T:compactification:prop}
The space $\Teich(T)$ embeds in $\matP(\matR^\calS)$ and its closure there is homeomorphic to a closed disc $D^2$. The interior of this disc is $\Teich(T)$ and its boundary contains $\calS$ as a dense subset.
\end{prop}
\begin{proof}
Everything can be written explicitly by identifying $\Teich(T)$ with the half-plane $H^2$ and $\calS$ with the set $\matQ\cup \{\infty\} \subset \partial H^2$. Exercise \ref{i:torus:ex} and Proposition \ref{lunghezza:toro:prop} give
\begin{align*}
i\left(\frac pq, \frac rs \right) & = \left|\det \begin{pmatrix} p & r \\ q & s \end{pmatrix}\right| = |ps-qr| = |s|\cdot \left|p-q\frac rs \right|, \\ 
\ell^{\frac pq}(z) & = \frac {|p+qz|}{\sqrt{\Im z}}.
\end{align*}
Therefore the images of $\frac rs\in\calS$ and $z\in H^2 = \Teich(T)$ in $\matP(\matR^\calS)$ are respectively the functionals
\begin{align*}
\frac pq & \longmapsto \left|p - q\frac rs\right| \ {\rm if\ } s \neq 0 {\rm \ and \ } |q|\ {\rm if } \ s=0, \\
\frac pq & \longmapsto \left|p+qz\right|.
\end{align*}
We could remove the constants $|s|$ and $\sqrt{\Im z}$ because we are considering functionals in $\matP(\matR^\calS)$ rather than in $\matR^\calS$. 
We define for all $z \in  \overline{H^2}$ the functional
\begin{align*}
f_z & \colon \frac pq \longmapsto |p+qz| \quad {\rm if\ } z\neq \infty,\\
f_\infty & \colon \frac pq \longmapsto |q|.
\end{align*}
and we get a continuous immersion $z \mapsto f_z$ of $\overline{H^2}$ into $\matP(\matR^\calS)$. The immersion is closed because it sends a compact space to a Hausdorff space, hence it is a homeomorphism onto its image. The image is the closure of $\Teich(T)$ and its boundary contains $\calS$ as a dense set.
\end{proof}

\section{Geodesic currents}
We will prove Thurston's compactification Theorem \ref{compattificazione:teo} using a slightly different perspective. We are guided by the analogy with the hyperbolic space: when we embed (actually, define) the space $\matH^n$ in $\matR^{n+1}$, we make an essential use of the lorentzian form $\langle,\rangle$ in $\matR^{n+1}$, so that $\matH^n$ consists of some points $v$ with $\langle v,v \rangle = -1$, and $\partial \matH^n$ may be identified with the light cone rays, that is the rays spanned by vectors $v$ with $\langle v , v \rangle = 0$.

It would be nice to have a similar nice bilinear form $\langle,\rangle$ on $\matR^\calS$, with the property that $\Teich(S_g)$ and $\calS$ consist of points with $\langle v,v\rangle = k$ and $\langle v, v \rangle = 0$ respectively, for some fixed $k\neq 0$. Unfortunately we are not able to define such a bilinear form in $\matR^\calS$, and in order to get one we now substitute $\matR^\calS$ with a similar (but more structured) infinite-dimensional space, the space of \emph{geodesic currents}.

In geometric measure theory, a \emph{current} is a measure on some space, which generalises the notion of $m$-dimensional submanifold in a $n$-manifold. A \emph{geodesic current} is a similar tool introduced by Bonahon in 1988 specifically designed for geodesics on hyperbolic surfaces.

We will see that the geodesic currents are indeed equipped with a bilinear form, which generalises beautifully both the intersection number $i(\alpha, \beta)$ of curves and the length functions $\ell^\gamma$ on the Teichm\"uller space.

\subsection{Geodesics}
From now on, and through all the rest of this chapter, we will consider geodesics only as subsets, neglecting their parametrisation.

More precisely, let $M$ a complete hyperbolic manifold. We indicate by $\calG(M)$ the set of the supports of all the complete non-trivial geodesics $\matR \to M$. With a little language abuse, we call an element of $\calG(M)$ a \emph{geodesic}. 

A geodesic in $\calG(M)$ is \emph{closed} if it is the support of a closed geodesic $S^1 \to M$, that is if it is compact; otherwise, it is \emph{open}. We say that a geodesic is \emph{simple} if it has a simple geodesic parametrisation, either as an open geodesic $\matR \to M$ or as a closed one $S^1 \to M$. 
  
We are particularly interested in the set $\calG = \calG(\matH^2)$ of lines in $\matH^2$. A line is determined by its extremes, hence there is a natural bijection 
$$\calG \longleftrightarrow \left(\partial \matH^2 \times \partial \matH^2\setminus \Delta\right)/_\sim$$
where $\Delta = \{(a,a)\ |\ a\in\partial\matH^2\}$ is the diagonal and $(a,b) \sim (b,a)$. We assign to $\calG$ the topology of $\left(\partial \matH^2 \times \partial \matH^2\setminus \Delta\right)/_\sim$. 
\begin{ex} \label{Mobius:ex}
The space $\calG$ is homeomorphic to an open M\"obius strip. The lines intersecting a compact set $K\subset \matH^2$ form a compact subset of $\calG$.
\end{ex}
The isometries of $\matH^2$ act naturally on $\calG$ by homeomorphisms.
\begin{prop}
If $S=\matH^2/_\Gamma$ is a complete hyperbolic surface there is a natural bijection
$$\calG(S) \longleftrightarrow \calG/_\Gamma.$$
\end{prop}
\begin{proof}
Every geodesic in $S$ lifts to a $\Gamma$-orbit of lines in $\matH^2$.
\end{proof}

It is typically more comfortable to lift objects from $S$ to the universal cover $\matH^2$, and to study them there: we will often see a geodesic $\gamma\in\calG(S)$ as a $\Gamma$-orbit of lines in $\matH^2$. 

We will since now consider only closed surfaces $S_g$ of genus $g\geqslant 2$.
We now prove the crucial fact that $\calG(S_g)$ depends only mildly on the chosen hyperbolic metric for $S_g$. 
An \emph{isomorphism} between two pairs $(\calG_1, \Gamma_1)$ and $(\calG_2, \Gamma_2)$ of groups $\Gamma_i$ acting on topological spaces $\calG_i$ is an isomorphism $\psi\colon\Gamma_1\to \Gamma_2$ together with a $\psi$-equivariant homeomorphism $\calG_1\to \calG_2$.

\begin{prop}
Let $g\geqslant 2$ and $S_g = \matH^2/_\Gamma$ have a hyperbolic metric. The pair $(\calG, \Gamma)$ does not depend (up to canonical isomorphisms) on the chosen hyperbolic metric.
\end{prop}
\begin{proof}
Let $m,m'$ be two hyperbolic structures on $S_g$, inducing two different coverings $\pi, \pi' \colon \matH^2 \to S_g$. The identity map $S_g\to S_g$ lifts to a map $\matH^2 \to \matH^2$ that extends to an equivariant homeomorphism $\partial \matH^2 \to \partial \matH^2$ by Theorem \ref{estensione:teo}. This induces an equivariant homeomorphism $\calG \to \calG$.
\end{proof}

Here and in the next sections, the hyperbolic metric on $S_g$ has only an auxiliary role: we need it to identify $S_g$ with $\matH^2/_\Gamma$ and to define and study some geometric objects like $\calG(S_g) = \calG/_\Gamma$, but most of our discoveries will be independent \emph{a posteriori} of the auxiliary hyperbolic metric. 

Let $S_g=\matH^2/_\Gamma$ be a hyperbolic surface. We note that $\Gamma$ does \emph{not} act properly discontinuously on $\calG$.

\begin{prop} \label{orbit:discrete:prop}
The $\Gamma$-orbit of a line $l\in\calG$ is discrete if and only if $l$ projects to a closed geodesic in $S_g$.
\end{prop}
\begin{proof}
Let $\pi(l)\subset S_g$ be the projection of $l$ in $S_g$. Since $S_g$ is compact, the projection $\pi(l)$ is \emph{not} a closed geodesic $\Leftrightarrow$ there is a small disc $D \subset S_g$ intersecting $\pi(l)$ into infinitely many distinct segments $\Leftrightarrow$ there is a small disc $D\subset \matH^2$ intersecting infinitely many lines of the $\Gamma$-orbit of $l$ $\Leftrightarrow$ the $\Gamma$-orbit is not discrete.
\end{proof}

The lines $l\subset \matH^2$ that project to closed geodesics in $S_g$ are precisely the axis of the hyperbolic isometries $\varphi$ in $\Gamma$, and the $\Gamma$-orbit of one such axis $l$ consists of the axis of all the isometries in $\Gamma$ conjugate to $\varphi$. Two such axis are either incident or ultraparallel by Corollary \ref{axis:cor}, and the closed geodesic in $S_g$ is simple if and only if all the distinct axis in the $\Gamma$-orbit are ultraparallel.

Proposition \ref{orbit:discrete:prop} implies in particular that the $\Gamma$-orbit of $l$ is discrete only for countably many lines $l\in\calG$. 

Recall from Corollary \ref{unique:closed:geodesic:cor} that every homotopically non-trivial closed curve in $S_g$ is homotopic to a unique closed geodesic. 

\subsection{Geodesic currents} \label{correnti:subsection}
We now introduce a measure-theoretical notion which is at first sight unrelated with everything we have seen up to now; we will later show that, on the contrary, it generalises many of the geometric objects that we have encountered in the last pages. The basic notions of measure theory that we will need are summarised in Section \ref{measure:theory:section}.

\begin{defn}
Let $S_g=\matH^2/_\Gamma$ be a hyperbolic surface. 
A \emph{geodesic current} on $S_g$ is a locally finite $\Gamma$-invariant Borel measure on $\calG = \calG(\matH^2)$.\index{geodesic current}
\end{defn}

We denote by $\calC = \calC (S_g)$ the set of all the geodesic currents in $S_g$. It is a subset of the space $\calM(\calG)$ of all the locally finite Borel measures on $\calG$, closed with respect to sums and products with non-negative scalars, and it inherits its topology, see Section \ref{topology:measure:subsection}.

The currents space $\calC$ is independent of the auxiliary hyperbolic metric on $S_g$ up to canonical isomorphisms, since the pair $(\calG, \Gamma)$ is.

\subsection{Closed geodesics}
We now introduce a fundamental example of geodesic current.

\begin{example}[Closed geodesics]
\label{scc:example}
A closed geodesic $\gamma$ on $S_g=\matH^2/_\Gamma$ lifts by Proposition \ref{orbit:discrete:prop} to a discrete $\Gamma$-orbit of lines in $\matH^2$. The Dirac measure on this discrete set is locally finite and $\Gamma$-invariant, hence it is a geodesic current.
\end{example}

We can therefore interpret every closed geodesic in $S_g$ as a particular geodesic current with discrete support. In particular we get an embedding
$$\calS \hookrightarrow \calC$$
of the set $\calS$ of all (isotopy classes of) unoriented non-trivial simple closed curves in $S_g$ into $\calC$. The embedding is defined by taking the geodesic representative of each curve. Recall that the hyperbolic metric plays only an auxiliary role.

The following proposition implies that -- conversely -- every current supported on a discrete set is a linear combination of closed geodesics. 

\begin{prop} \label{atomic:prop}
If $l\in \calG$ is an atomic point for a geodesic current $\mu$, that is if $\mu(\{l\})>0$, then $l$ projects to a closed geodesic in $S_g$.
\end{prop}
\begin{proof}
Since $l$ is atomic and $\mu$ is $\Gamma$-invariant and locally finite, the $\Gamma$-orbit of $l$ in $\calG$ is discrete. We conclude thanks to Proposition \ref{orbit:discrete:prop}.
\end{proof}

It is an important aspect of the theory that atomic points may be only of a very specific type. We remind that many closed geodesics are not simple.

\subsection{Pencils} \label{pencils:subsection}
We have determined the 0-dimensional subsets of $\calG$ that may have positive mass, and we now look at some natural 1-dimensional ones. Let a \emph{pencil} $p\subset\calG$ of lines centred at $a\in\partial \matH^2$ be some Borel set of lines all having one endpoint at $a$, as in Figure \ref{pencil-box:fig}-(left).\index{pencil}

\begin{figure}
\begin{center}
\includegraphics[width = 9 cm] {\iftoggle{BW}{pencil-box-BW}{pencil-box}}
\nota{A pencil and a box (in the Klein model).}
\label{pencil-box:fig}
\end{center}
\end{figure}

We note that a pencil may contain at most one axis of some hyperbolic transformation in $\Gamma$, because two axis are never asymptotically parallel by Corollary \ref{axis:cor}.

\begin{prop} \label{pencils:prop}
Let $\mu$ be a geodesic current. The mass of a pencil $p$ is zero, unless it contains the axis of some non-trivial element in $\Gamma$.
\end{prop}
\begin{proof}
It suffices to consider the case where $p$ contains no axis, and consists of all lines with one endpoint in $a$ and another in some segment $[b,c] \subset \partial \matH^2$ not containing $a$.

As the lines in $p$ point to $a$, they project in $S_g$ to a narrow beam of lines that go on running forever and exponentially-narrowing in $S_g$ but never close up. There is a disc $D\subset S_g$ intersecting them infinitely many times, hence a lift $D\subset \matH^2$ intersects infinitely many disjoint $\Gamma$-translates of $p$. The lines intersecting $D$ form a compact subset of $\calG$, hence $\mu(p)=0$.
\end{proof}

\subsection{Boxes}
We now investigate some natural 2-dimensional subsets of $\calG$. Let $a,b,c,d \in \partial \matH^2$ be four distinct counterclockwise-ordered points as in Figure \ref{pencil-box:fig}-(right): these determine two disjoint arcs $[a,b], [c,d] \subset \partial \matH^2$ and hence a compact set $B=[a,b] \times [c,d] \subset \calG$ consisting of all lines with endpoints in $[a,b]$ and $[c,d]$, see Figure \ref{pencil-box:fig}-(right). We call this compact set $B$ a \emph{box}.\index{box}

The topological boundary $\partial B$ of a box $B$ has four sides: the four pencils of lines in $B$ with one endpoint in $a,b,c$, or $d$. We say that the box $B$ is \emph{generic} if no line $l\in\partial B$ projects to a closed geodesic, \emph{i.e.}~no $l\in\partial B$ is the axis of some hyperbolic transformation in $\Gamma$. This is indeed a generic condition: if neither of $a,b,c,d$ is the endpoint of some axis, then $B$ is certainly generic, and recall that there are only countably many axis overall.

We equip the currents space $\calC$ with the weak-* topology introduced in Section \ref{topology:measure:subsection}. We now prove that a converging sequence of currents behave well on generic boxes.

\begin{prop}
If $\mu_i \rightharpoonup \mu$ is a converging sequence of geodesic currents and $B$ is a generic box, then $\mu_i(B) \to \mu(B)$.
\end{prop}
\begin{proof}
By Proposition \ref{delta:zero:prop} it suffices to prove that $\mu(\partial B) = 0$, and this follows from Proposition \ref{pencils:prop}.
\end{proof}

\subsection{Rigidity of atoms}
We now use boxes to prove that the atomic points in $\calC$ are somehow rigid. 

\begin{prop} \label{closed:discrete:currents:prop}
The closed geodesics form a discrete subset in $\calC$.
\end{prop}
\begin{proof}
Let $\mu \in \calC$ be a closed geodesic, \emph{i.e}~it is the Dirac measure on the $\Gamma$-orbit of an axis $l\subset \matH^2$ of some hyperbolic transformation $\varphi\in\Gamma$. 

\begin{figure}
\begin{center}
\includegraphics[width = 9 cm] {\iftoggle{BW}{boxes-BW}{boxes}}
\nota{The box $B$ containing $l$ and the sub-boxes $B'$ and $B''$ (left). A line $l\in U_\gamma$ parametrised by the pair $(t,\theta)$ (right).}
\label{boxes:fig}
\end{center}
\end{figure}

Let $B=[a,b]\times [c,d]$ be a small generic box containing $l$ and no other $\Gamma$-translate of $l$, so that $\mu(B) = 1$. Up to substituting $\varphi$ with $\varphi^{-1}$ we may suppose that $\varphi([a,b]) = [a',b'] \subset [a,b]$. Consider the sub-boxes $B' = [a,a']\times [c,d]$ and $B'' = [b',b]\times [c,d]$ shown in Figure \ref{boxes:fig}-(left). Since $\varphi$ contracts $[a,b]$ and expands $[c,d]$, we have 
\begin{equation} \label{Bl:eqn}
B\setminus p \subset \bigcup_{i=1}^\infty \varphi^i(B' \sqcup B'').
\end{equation}
where $p$ is the pencil containing $l$ and other lines with one endpoint varying in $[b,c]$. Note that $p$ contains no axis of $\Gamma$ other than $l$.
If necessary, we move $a'$ and $b'$ slightly farther from $a$ and $b$ to ensure that $B'$ and $B''$ (and hence all their $\varphi$-translates) are generic while (\ref{Bl:eqn}) still holds.

We have $\mu(B' \sqcup B'')=0$ and we define an open neighbourhood $U$ of $\mu$ in the current space $\calC$ as follows: 
$$U = \big\{ \eta\in \calC\ |\ \eta(B' \sqcup B'') < \varepsilon,\ \eta(B)>1-\varepsilon\big\}$$ 
for some fixed $0<\varepsilon < 1$. This is indeed an open set because $B, B',$ and $B''$ are generic. We now prove that $U$ contains no closed geodesic except $\mu$.

By contradiction, suppose $U$ contains a closed geodesic $\eta$. Since $\eta(B)>1-\varepsilon>0$, at least one atom $l'$ of $\eta$ is contained in $B$. If $l'=l$ then $\eta=\mu$ and we are done. Otherwise, the line $l'$ is contained in $B\setminus p$ and hence in $\varphi^i(B' \sqcup B'')$ for some $i$. This is impossible since this set has measure $\varepsilon < 1$.
\end{proof}

We have classified the elements in $\calC$ supported on discrete sets, and it is now time to  introduce some geodesic currents supported on the whole of $\calG$.

\subsection{The Liouville measure}
There is a continuous measure on $\calG = \calG(\matH^2)$ that is invariant under the action of the whole isometry group $\Iso(\matH^2)$. This measure is called the \emph{Liouville measure}, and is defined as follows.\index{Liouville measure}

Let $\gamma\colon \matR \to \matH^2$ be a geodesic parametrised by arc length, and $U_\gamma\subset \calG$ be the open set consisting of all the lines intersecting $\gamma$, except $\gamma$ itself. We can parametrise $U_\gamma$ via the homeomorphism 
$$\matR \times (0,\pi) \to U_\gamma$$ 
that sends $(t, \theta)$ to the line $l$ that intersects $\gamma$ at the point $\gamma(t)$ with angle $\theta$, see Figure \ref{boxes:fig}-(right). We define a volume 2-form on $U_\gamma$ by setting: 
$$L_\gamma = \frac 12 \sin\theta\, dt\wedge d\theta.$$
Since $\calG$ is non-orientable, we cannot hope to define a global area form on $\calG$. However all these local forms match up to sign and hence give rise to a measure.
\begin{prop}
The charts $U_\gamma$ form a differentiable atlas for $\calG$. The 2-forms $L_\gamma$ match up to sign and hence define a measure  $L$ on $\calG$.
\end{prop}
\begin{proof}
Every line in $\matH^2$ intersects some other line, hence the charts cover $\calG$. 
We consider a line $r\in U_\gamma\cap U_{\gamma'}$. The charts $U_\gamma$ and $U_{\gamma'}$ have parametrisations $(t,\theta)$ and $(t', \theta')$ and 2-forms 
$$L_\gamma = \frac 12 \sin\theta dt \wedge d\theta, \quad 
L_{\gamma'} = \frac 12 \sin\theta' dt' \wedge d\theta'.$$ 
We consider the jacobian $J = \frac{\partial (t', \theta')}{\partial (t,\theta)}$ 
and recall that 
$$dt' \wedge d\theta' = \det J \cdot dt \wedge d\theta$$
so we need to show that
\begin{equation} \label{J:eqn}
\det J = \pm\frac {\sin \theta}{\sin \theta'}.
\end{equation}
We can reparametrise $\gamma$ at our please, since this may change $L_\gamma$ only by a sign. We may suppose that $\gamma$ and $\gamma'$ are asymptotically parallel, because any two lines $\gamma, \gamma'$ intersecting $r$ are connected by a path $\gamma = \gamma_1, \ldots, \gamma_k = \gamma'$ of lines intersecting $r$, such that any two subsequent lines $\gamma_i$, $\gamma_{i+1}$ are asymptotically parallel (actually $k= 3$ suffices). 

\begin{figure}
\begin{center}
\includegraphics[width = 5 cm] {\iftoggle{BW}{Liouville-BW}{Liouville}}
\nota{The two asymptotic lines $\gamma, \gamma'$ in the half-space model $H^2$ and the line $r$ intersecting both.}
\label{Liouville:fig}
\end{center}
\end{figure}

We represent the lines $\gamma, \gamma',$ and $r$ in the half-plane model $H^2\subset \matR^2$ as in Figure \ref{Liouville:fig}. Up to isometries and reparametrisations we have 
$$\gamma(t) = (0, e^t), \qquad \gamma'(t') = (1,e^{t'}).$$
The line $r$ is a Euclidean half-circle with center at some point $C\in \matR$ and with some radius $R>0$. It intersects $\gamma$ and $\gamma'$ at some heights $e^t$ and $e^{t'}$, with some angles $\theta$ and $\theta'$. We now calculate the determinant of $J = \frac{\partial (t', \theta')}{\partial (t,\theta)}$. Figure \ref{Liouville:fig} shows that
$$ R = \frac{e^t}{\sin \theta} = \frac{e^{t'}}{\sin \theta'}$$
$$ C = -e^t \cot \theta, \qquad 1 - C = e^{t'} \cot \theta'.$$
Now
$$\frac{\partial (C,R)}{\partial(t,\theta)} = 
\begin{pmatrix} 
-e^t\cos\theta\sin^{-1}\theta & e^t\sin^{-2}\theta \\
e^t\sin^{-1}\theta & -e^t \cos\theta\sin^{-2} \theta
\end{pmatrix}
$$
whose determinant is $-e^{2t} \sin^{-1}\theta = -R\sin\theta$, and we get the same formula for $\frac{\partial (C,R)}{\partial(t',\theta')}$. Therefore
$$\det \frac{\partial (t', \theta')}{\partial (t,\theta)} =
\det \frac{\partial (t', \theta')}{\partial (C, R)}
\cdot 
\det \frac{\partial (C,R)}{\partial(t,\theta)} 
=
(-R\sin\theta')^{-1} \cdot (R \sin\theta) = \frac{\sin\theta}{\sin\theta'}$$
and we finally obtain (\ref{J:eqn}), as required.
\end{proof}
The measure $L$ on $\calG$ is called the \emph{Liouville measure} and is clearly invariant by the action of $\Iso(\matH^2)$. The renormalising factor $\frac 12$ in the definition was chosen to get the following property, which nicely characterises $L$.

\begin{prop} \label{l:prop}
Let $s\subset \matH^2$ be a geodesic segment of length $l$. The lines in $\matH^2$ intersecting $s$ form a set of measure $l$.
\end{prop}
\begin{proof}
The set has measure
$$\int_0^\pi\int_0^l \frac 12 \sin\theta \, dt\, d\theta = l \int_0^\pi \frac 12 \sin \theta \, d\theta = l.$$
The proof is complete.
\end{proof}

The Liouville measure $L$ is supported on the whole of $\calG$.

\begin{ex}
The Liouville measure of a box $B=[a,b]\times [c,d]$ is 
$$L(B) = \big| \log \beta(a,b,c,d) \big|$$
where $\beta$ is the cross-ratio of the four points.
\end{ex}

\subsection{The Liouville current}
Let now $S_g=\matH^2/_\Gamma$ be equipped with a hyperbolic metric. The Liouville measure on $\calG$ is $\Iso(\matH^2)$-invariant: in particular it is $\Gamma$-invariant and hence defines a current $L \in \calC(S_g)$, called the \emph{Liouville current}.\index{Liouville current} 

It is important to note that the space of currents $\calC = \calC(S_g)$ does not depend on the given hyperbolic metric for $S_g$, but the Liouville current does! Every metric $m\in \Teich(S_g)$ induces a Liouville current $L_m\in \calC$, and in this way we get a \emph{Liouville map}\index{Liouville map}
$$\Teich(S_g) \longrightarrow \calC$$
that sends $m$ to $L_m$.
We will soon see that this map is injective: as promised, we have mapped both $\calS$ and $\Teich(S_g)$ inside the currents space $\calC$. 

\begin{rem}
We constructed the Liouville current $L_m$ by assigning to $\calG$ a differentiable atlas -- hence a smooth structure -- and a 2-form defined up to sign. We note that the smooth structure on $\calG$ also depends crucially on $m$, because the boundary extension of Theorem \ref{estensione:teo} is not guaranteed to be a diffeomorphism.
\end{rem}

The crucial feature that makes $\calC$ preferable to $\matP(\matR^\calS)$ as a comfortable ambient space for both $\calS$ and $\Teich(S_g)$ is the existence of a nice bilinear form in $\calC$ that extends both the length and the geometric intersection of closed geodesics. We now introduce this bilinear form.

\subsection{The projective frame bundle}
We denote by $\calI \subset \calG \times \calG$ the open subset consisting of all the pairs of \emph{incident} distinct lines in $\matH^2$. We give $\calI$ the topology induced by $\calG \times \calG$, hence $\calI$ is an open topological 4-manifold. Since two distinct incident lines intersect in a single point $p\in\matH^2$, the set $\calI$ can be interpreted as the set of triples $(p,l_1,l_2)$ with $p\in \matH^2$ and $l_1,l_2$ two distinct vector lines in the tangent plane $T_p\matH^2$.

Let now $S_g=\matH^2/_\Gamma$ be equipped with a hyperbolic metric. The action of $\Gamma$ on $\calG$ is not properly discontinuous, but the diagonal action of $\Gamma$ on $\calI$ is, as the following shows.
\begin{prop}
The map $\calI \to \calI/_\Gamma$ is a topological covering.
\end{prop}
\begin{proof}
The group $\Gamma$ acts freely and properly discontinuously on $\matH^2$, and hence it does so \emph{a fortiori} on the triples $(p,l_1,l_2)$.
\end{proof}

In particular $\calI/_\Gamma$ is a four-manifold, and it is naturally a bundle over $S_g$, whose fibre above $p\in S_g$ consists of the ordered pairs of distinct vector lines $l_1,l_2$ in $T_pS_g$. Recall that the \emph{frame bundle} on a manifold $M$ is a bundle whose fibre above $p\in M$ is the set of all \emph{frames} (i.e.~basis) in $T_pM$. The space $\calI/_\Gamma$ can be seen as some projective quotient of the frame bundle on $S_g$.

\subsection{Intersection form}
Two geodesic currents $\alpha, \beta \in \calC = \calC(S_g)$ induce a $\Gamma$-invariant product measure $\alpha\times \beta$ on $\calG\times \calG$ and hence on $\calI$. This measure descends via the covering $\calI \to \calI/_\Gamma$ to a measure on $\calI/_\Gamma$ which we still indicate by $\alpha\times \beta$, see Section \ref{measure:descends:subsection}.\index{intersection form on currents}

\begin{defn}
The \emph{intersection} $i(\alpha,\beta)$ of two geodesic currents is the total volume of $\calI/_\Gamma$ in the measure $\alpha\times\beta$. 
\end{defn}

The intersection $i(\alpha, \beta)$ does not depend on the auxiliary hyperbolic metric for $S_g$. The finiteness of $i(\alpha,\beta)$ is not immediate since $\calI/_\Gamma$ is not compact, and it deserves a proof.

\begin{prop}
The intersection $i(\alpha,\beta)$ is finite.
\end{prop}
\begin{proof}
Let $D$ be a compact fundamental domain for $S_g = \matH^2/_\Gamma$. Since $D$ is compact, the lines intersecting $D$ form a compact subset $X\subset \calG$ which has finite measures $\alpha(X)$ and $\beta(X)$. The projection sends $(X \times X) \cap \calI$ surjectively onto $\calI/_\Gamma$. 
Therefore $i(\alpha,\beta) = (\alpha\times \beta) (\calI/_\Gamma) < \alpha(X) \cdot \beta(X)$.
\end{proof}

The form $i$ is clearly bilinear and symmetric. 

\subsection{Geometric intersection and length of closed geodesics}
We now show that this abstract-looking intersection form generalises both the geometric intersection and the length of simple closed geodesics.
We consider as usual the set $\calS$ of simple closed curves as a subset of $\calC$. The following proposition explains why we employed the notation $i$ for the intersection form of two geodesic currents.

\begin{prop}
If $\alpha, \beta \in \calS$, the number $i(\alpha,\beta)$ is the geometric intersection of the simple closed curves $\alpha$ and $\beta$.
\end{prop}
\begin{proof}
Represent $\alpha$ and $\beta$ as simple closed geodesics. The measure $\alpha\times \beta$ is the Dirac measure with support the pairs of incident lines in $\matH^2$ that cover respectively $\alpha$ and $\beta$. The $\Gamma$-orbits of these pairs are in natural bijection with the transverse intersection points in $\alpha\cap \beta$. Hence the volume of $\calI/_\Gamma$ is the cardinality of $\alpha\cap\beta$, except when $\alpha = \beta$ and in this case we get zero.
\end{proof}

We note in particular that $i(\alpha,\alpha)=0$ for every $\alpha\in\calS$, hence the set $\calS$ is contained in the ``light cone'' consisting of all the currents $\alpha$ with $i(\alpha,\alpha)=0$. (Of course $i(\alpha, \alpha)\geqslant 0$ for every geodesic current $\alpha$.) We now consider the Liouville current $L_m\in\calC$ determined by some hyperbolic metric $m$ on $S_g$.

\begin{prop} \label{i:metric:curve:prop}
If $\alpha\in \calS$ then $i(L_m, \alpha) = \ell^\alpha(m)$ is the length of the geodesic representative of $\alpha$ in the metric $m$.
\end{prop}
\begin{proof}
Represent $\alpha$ as a simple closed geodesic in the metric $m$, of length $\ell^\alpha(m)$. The measure $L_m \times \alpha$ has its support on the pairs $(l,l')$ of incident lines where $l$ is arbitrary and $l'$ is a lift of $\alpha$. 

Any segment $s'\subset l'$ of length $\ell^\alpha(m)$ in a fixed lift $l'$ of $\alpha$ is a fundamental domain for the action of $\Gamma$ on the lifts of $\alpha$. Therefore $i(L_m, \alpha)$ is the volume of the pairs $(l,l')$ where $l'$ is fixed, and $l$ is arbitrary and intersects $s'$. By Proposition \ref{l:prop} these pairs have volume $\ell^\alpha(m)$.
\end{proof}

The intersection form $i$ on $\calC$ generalises both the geometric intersection of curves and the length functions on Teichm\"uller space, two apparently unrelated objects! As promised, we can easily deduce that $\Teich(S_g)$ embeds in $\calC$.

\begin{cor}
The Liouville map $\Teich(S_g) \to \calC$ is injective.
\end{cor}
\begin{proof}
If $m\neq m'$ there is a curve $\gamma\in\calS$ with $\ell^\gamma(m) \neq \ell^\gamma(m')$ by Proposition \ref{9g-9:prop} and hence $i(L_m,\gamma) \neq i(L_{m'},\gamma)$. Therefore $L_m \neq L_{m'}$.
\end{proof}

We will since now consider both $\Teich(S_g)$ and $\calS$ as subsets of $\calC$.
We know the geometric meaning of the intersection of two curves, and of a curve and a hyperbolic metric. What is the intersection of two hyperbolic metrics? When they coincide, we get a positive constant.

\begin{prop}
If $m \in \Teich(S_g)$ we have $i(L_m, L_m) = -\pi^2 \chi(S_g)$.
\end{prop}
\begin{proof}
The metric $m$ produces a smooth structure on $\calG$ and hence on $\calI$, and $L_m$ is induced by a $\Gamma$-invariant 2-form on $\calG$ defined up to sign, that we also denote as $L_m$. We need to integrate the 4-form $L_m \times L_m$ on $\calI/_\Gamma$.

Recall that $\calI/_\Gamma$ may be interpreted as a bundle over $S_g$.
We make a first-order computation on a very small region of $S_g$, which we may suppose to be Euclidean: a very thin rectangle $R$ in $S_g$ of sides $l$ and $a$ with $a \ll l$. We compute the volume of the portion of bundle lying above $R$.

We integrate $L_m \times L_m$ on all the pairs of segments intersecting in some point in $R$, and we consider only segments with both endpoints on the big $l$ sides, neglecting the small $a$ sides. A segment meeting the $l$ sides at angle $\theta$ has length $L(\theta) = a/\sin \theta$. By Proposition \ref{l:prop} the $L_m \times L_m$ volume of the portion of $\calI$ above $R$ is (at a first order) equal to
$$\int_0^\pi \int_0^l \frac 12 \sin \theta \cdot L(\theta) \, dt\, d\theta \sim
\frac 12 l \int_0^\pi \sin \theta \cdot \frac a{\sin \theta} = \frac{\pi}2 la.$$

The first-order contribution of a small region is hence $\frac{\pi}2$ times its area. Since the volume is induced by a differentiable form, by taking the limit we find that the contribution of any region is precisely $\frac{\pi}2$ its area and we conclude by the Gauss-Bonnet theorem.
\end{proof}

We have already noted that $i$ is symmetric and bilinear. 
As we mentioned, there is an evident analogy between the embeddings $\Teich(S_g)\hookrightarrow \calC$ and $\matH^n\hookrightarrow \matR^{n+1}$, since there is a bilinear form on both spaces $\calC$ and $\matR^{n+1}$, and $\Teich(S_g)$ is also contained in a ``hyperboloid'' consisting of vectors $v$ with $i(v,v) = k$ for some fixed $k \neq 0$. On the other hand $\calS$ is contained in the ``light cone'' formed by all vectors $v$ with $i(v,v)=0$. 

We may compactify $\matH^n$ by projecting the hyperboloid to the Klein model $K^n\subset \matRP^{n}$. Now $K^n$ is an open disc, whose boundary is the image of the light cone. We would like to apply the same compactification strategy to $\Teich(S_g)$. Before doing that, we need to prove a slightly technical fact: that the bilinear form $i$ is continuous. 

\subsection{Continuity of the intersection form}
Let $S_g=\matH^2/_\Gamma$ be a closed hyperbolic surface.
The following theorem is not obvious because $\calI$ is not compact: the proof relies on the fact that closed geodesics form a discrete rigid set in $\calC$ and therefore atomic points cannot ``escape to infinity'' on a converging sequence of currents.

\begin{teo}
The form $i\colon \calC \times \calC \to \matR$ is continuous.
\end{teo}
\begin{proof}
Let $\beta_1, \beta_2$ be two currents converging to $\alpha_1, \alpha_2$. We need to prove that $i(\beta_1,\beta_2)$ converges to $i(\alpha_1, \alpha_2)$. The bundle $\calI/_\Gamma$ is not compact and we need to control that no mass escapes to infinity: the infinity here is the diagonal $\Delta \subset \calG \times \calG$, so our aim is to cover $\Delta$ with small boxes and prove that they contribute very little to $i(\beta_1,\beta_2)$. For any box $B$, we define
$$\Psi(B) = \big((B\times B)\cap \calI \big)/_\Gamma.$$
We show that for every $\varepsilon >0$ there are finitely many small generic boxes $B_i$ whose $\Gamma$-translates cover $\calG$, such that 
$$(\beta_1 \times \beta_2) (\cup_i\Psi(B_i)) < \varepsilon$$
as soon as $\beta_1,\beta_2$ are sufficiently close to $\alpha_1,\alpha_2$. Since the $\Gamma$-translates of $\cup_i (B_i\times B_i)$ form an open neighbourhood of $\Delta$, the set $(\calI/_\Gamma )\setminus (\cup_i \Psi(B_i))$ is compact with zero-measure boundary (because the boxes $B_i$ are generic) so its contribution to $i(\beta_1, \beta_2)$ tends to that to $i(\alpha_1,\alpha_2)$ if $\beta_1,\beta_2$ are sufficiently close, and the theorem is proved.

It remains to construct the boxes $B_i$. We first fix finitely many generic boxes $B_i$ whose $\Gamma$-translates cover $\calG$. Let $K>0$ be bigger than the total $\alpha_1$ and $\alpha_2$-mass of the boxes. 

We now pick a small number $0<\varepsilon' < \varepsilon / 2K$ and subdivide each box $B_i$ into finitely many generic sub-boxes (which we still call $B_i$) such that the following holds for each $i$ and each $j=1,2$:
\begin{enumerate}
\item either $\alpha_j(B_i) < \varepsilon'$, or
\item there is an atomic point $l \in B_i$ with $\alpha_j(l) > \varepsilon'$ and $\alpha_j(B_i \setminus l) < \varepsilon'$.
\end{enumerate}
If the second case holds for both $j=1,2$, we also require $l\in B_i$ to be the same line for $j=1,2$.
We note that the total $\alpha_1$ and $\alpha_2$-mass of the boxes is still smaller than $K$ after the subdivision.
We have
$$(\beta_1 \times \beta_2)(\Psi(B)) \leqslant \beta_1(B) \cdot \beta_2(B)$$
so if $\beta_1,\beta_2$ are sufficiently close to $\alpha_1, \alpha_2$ the boxes $B_i$ that are of the first kind for at least one value of $j=1,2$ contribute to $i(\beta_1, \beta_2)$ less than $2\varepsilon' K<\varepsilon$.

We are left to consider the boxes $B_i$ containing a line $l$ whose $\alpha_1$ and $\alpha_2$-mass are both bigger than $\varepsilon'$. This is a potentially bad situation since the point $(l,l) \in \calG\times \calG$ is atomic and contained in the frontier of $\calI$, and could enter inside $\calI$ abruptly for some arbitrarily small perturbations of $\alpha_1$ and $\alpha_2$. We show that this cannot happen because atomic points are rigid, as noted in Proposition \ref{closed:discrete:currents:prop}.

As in the proof of Proposition \ref{closed:discrete:currents:prop}, we set $B_i = B = [a,b] \times [c,d]$ and construct $B' = [a,a'] \times [c,d], B'' = [b',b] \times [c,d]$ whose $\varphi$-translates cover $B\setminus p$ where $p$ is a pencil containing $l$. We note that
$$(\beta_1 \times \beta_2) (\Psi(B)) \leqslant \beta_1(B' \sqcup B'') \cdot \beta_2(B) + \beta_2(B'\sqcup B'') \cdot \beta_1(B)$$
is smaller than $\varepsilon' (\beta_1(B) + \beta_2(B))$ and hence by summing on the $B_i$ we get again a contribution smaller than $2\varepsilon' K < \varepsilon$.
\end{proof}

\subsection{Filling geodesic currents}
We now define an interesting class of currents, whose importance will be evident in the next section.

We say that a geodesic current $\alpha\in\calC(S_g)$ \emph{fills} the surface $S_g=\matH^2/_\Gamma$ if every line in $\matH^2$ intersects transversely at least one line in the support of $\alpha$. For instance, a Liouville measure fills $S_g$ since its support is the whole of $\calG$. We say that $k$ closed geodesics $\gamma_1,\ldots, \gamma_k$ \emph{fill $S_g$} if the geodesic current $\gamma_1+\ldots + \gamma_k$ does.

\begin{prop}
Let $\gamma_1,\ldots, \gamma_k$ be closed geodesics. If $S_g \setminus (\gamma_1\cup \cdots\cup \gamma_k)$ consists of polygons, the curves fill $S_g$.
\end{prop}
\begin{proof}
Every geodesic in $S_g$ intersects these curves.
\end{proof}

As an example, we can pick two multicurves $\eta$ and $\mu$ in $S_g$ that intersect transversely forming only polygons but no bigons. By the bigon criterion $\eta$ and $\mu$ are in minimal position and by Proposition \ref{minimale:unica:prop} their geodesic representatives have the same configuration as $\eta\cup\mu$, so they fill $S_g$.

It is a bit more difficult to construct a single (non-simple) closed geodesic that fills $S_g$. The following remark is straightforward.

\begin{oss} \label{correnti:oss}
Let $\alpha$ and $\beta$ be currents. We have $i(\alpha,\beta)>0$ if and only if there are two lines in the supports of $\alpha$ and $\beta$ that intersect in a point.
\end{oss}

\begin{cor} \label{riempie:cor}
If $\alpha$ fills $S_g$ then $i(\alpha,\beta)>0$ for every non-trivial $\beta\in\calC$.
\end{cor}

\subsection{A compactness criterion}
We now state a simple and general compactness criterion for subsets of $\calC$ that has various nice (and apparently unrelated) geometric consequences.

\begin{prop}[Compactness criterion] \label{compattezza:C:prop}
If $\alpha\in\calC$ fills $S_g$, the set of all $\beta\in\calC$ with $i(\alpha,\beta)\leqslant M$ is compact, for all $M>0$.
\end{prop}
\begin{proof}
Let $C\subset\calC$ be the set of all $\beta$ with $i(\alpha,\beta)\leqslant M$. It is closed because $i$ is continuous.

Let $l$ be a line in $\matH^2$. By hypothesis there is another line $l'$ in the support of $\alpha$ which intersects $l$ in a point. Let $B'$, $B$ be two sufficiently small boxes neighbourhoods of $l',l$ in $\calG$, so that $B' \times B \subset \calI$ and
$B' \times B$ is mapped injectively into $\calI/_\Gamma$.
If $\beta\in C$ we have
$$\alpha(B')\beta(B) = (\alpha\times\beta)(B' \times B) \leqslant (\alpha\times\beta)\left(\calI/_\Gamma\right) = i(\alpha,\beta) \leqslant M.$$
Therefore every line $l$ in $\matH^2$ has a box neighbourhood $B$ such that
$$\beta(B) \leqslant K_l \quad \forall \beta \in C$$
for some constant $K_l = M/\alpha(B')$ that depends only on $l$. The set $C$ is relatively compact by Theorem \ref{compact:measure:teo}, and hence compact since it is closed. The proof is complete.
\end{proof}

We now state some corollaries. The first concerns the immersion of $\Teich(S_g)$ in $\calC$. 

\begin{cor}
The Liouville map $\Teich(S_g)\hookrightarrow\calC$ is proper and a homeomorphism onto its image.
\end{cor}
\begin{proof}
We prove that the map is proper: if $m_j\in \Teich(S_g)$ is a diverging sequence of metrics, we know from Proposition \ref{diverges:infinity:prop} that there is a closed curve $\gamma\in\calS$ such that $\ell^\gamma(m_j) = i(m_j, \gamma) \to \infty$ on a subsequence. Since $i$ is continuous, the sequence $m_j$ diverges also in $\calC$. 

We denote the Liouville map by $L$. We should now prove that $L$ is continuous, but we prefer to consider the inverse map $L^{-1}\colon L(\Teich(S_g)) \to \Teich(S_g)$. The map $L^{-1}$ is continuous because $i$ is and $\Teich(S_g)$ has the weakest topology where the length functions are continuous. We show that $L^{-1}$ is proper. Let $\gamma_1, \ldots, \gamma_k$ be simple closed curves that fill $S_g$. If $L(m_i)$ is a diverging sequence, by Proposition \ref{compattezza:C:prop} we have $i(L(m_i),\sum_t \gamma_t)\to \infty$ and hence $i(L(m_i),\gamma_t)\to \infty$ for some $t$. Therefore $m_i$ is divergent also in $\Teich(S_g)$.

Now $L^{-1}$ is continuous and proper and hence a homeomorphism by
Corollary \ref{propria:chiusa:cor}.
\end{proof}

A second immediate corollary is a general compactness criterion for the Teichm\"uller space.

\begin{cor} \label{bounded:metrics:compact:cor}
Let $\gamma_1,\ldots, \gamma_k$ be some closed geodesics that fill $S_g$. The metrics $m\in\Teich(S_g)$ with $\ell^{\gamma_i}(m) \leqslant M$ form a compact subset of $\Teich(S_g)$, for all $M>0$.
\end{cor}

We can similarly deduce the following.
\begin{cor}
Let $\gamma_1,\ldots, \gamma_k$ be some simple closed curves that fill $S_g$. For every $M$ there are only finitely many $\alpha \in \calS$ with $i(\gamma_i, \alpha) < M$ for all $i$.
\end{cor}

\begin{ex} Use the compactness criterion to re-prove that on a hyperbolic closed surface there are only finitely many closed geodesics of bounded length.
\end{ex}

We are finally ready to construct Thurston's compactification of the Teichm\"uller space, using Bonahon's geodesic currents. 

\subsection{Projective currents}
As usual, we pick a hyperbolic closed surface $S_g=\matH^2/_\Gamma$.
The currents space $\calC$ is equipped with a multiplication by positive scalars, hence we can define its projectivisation
$$\pi\colon \calC\setminus 0 \longrightarrow \matP\calC$$
where $\matP \calC = (\calC \setminus 0)/_\sim$ with $\alpha \sim \lambda \alpha$ for all $\lambda>0$. We give $\matP\calC$ the quotient topology.

\begin{prop} The space $\matP\calC$ is compact.
\end{prop}
\begin{proof}
Pick an $\alpha\in\calC$ that fills $S_g$. By the compactness criterion the set $C = \{\beta\in\calC\ |\ i(\alpha,\beta)=1\}$ is compact. By Corollary \ref{riempie:cor} we have $i(\alpha,\beta)>0$ for all $\beta$, hence $\lambda\beta \in C$ for some $\lambda>0$. Therefore $\pi(C) = \matP\calC$ and $\matP\calC$ is compact.
\end{proof}

We now want to embed both $\Teich(S_g)$ and $\calS$ in $\matP\calC$.

\begin{prop}
The composition $\calS \to \calC\setminus 0 \to \matP\calC$ is injective.
\end{prop}
\begin{proof}
Let $\gamma_1,\gamma_2 \in\calS$ be distinct. There is always a closed curve $\eta\in\calS$ with $i(\gamma_1,\eta)\neq 0$ and $i(\gamma_2,\eta)=0$.
\end{proof}

\begin{prop} \label{embed:Teich:prop}
The composition $\Teich(S_g) \to \calC \setminus 0\to \matP\calC$ is injective and a homeomorphism onto its image.
\end{prop}
\begin{proof}
We see $\Teich(S_g)$ already properly embedded in $\calC$.
Since $i(m,m) = -\pi^2\chi(S_g)$ is constant on $\Teich(S_g)$, the composition is injective. The restriction $\pi\colon \Teich(S_g) \to \pi(\Teich(S_g))$ is also continuous, and we now prove that it is proper (note that we need to restrict the codomain to get this).

Consider a diverging sequence $m_j\in \Teich(S_g)$. By compactness of $\matP\calC$ the sequence $[m_j]\in\matP\calC$ converges on a subsequence to some $[\alpha]\in\matP\calC$. For each $j$ there is a $\lambda_j>0$ such that $\lambda_jm_j \to \alpha$ in $\calC$. Since $m_j$ diverges in $\Teich(S_g)$ and hence in $\calC$ we get $\lambda_j\to 0$. We get
$$i(\alpha,\alpha) = \lim_{j\to \infty} i(\lambda_jm_j, \lambda_jm_j) =  -\pi^2\chi(S_g) \lim_{j\to \infty} \lambda_j^2 = 0.$$
In particular $[\alpha]\not\in \pi(\Teich(S_g))$. Therefore the restriction $\pi\colon \Teich(S_g) \to \pi(\Teich(S_g))$ is proper and hence a homeomorphism. 
\end{proof}

We will since now consider both $\calS$ and $\Teich(S_g)$ embedded in $\matP\calC$. 

\subsection{Thurston's compactification}
We consider the closure $\overline{\Teich(S_g)}$ of $\Teich(S_g)$ inside $\matP\calC$. This closure is compact since $\matP\calC$ is, and it is called the \emph{Thurston compactification} of Teichm\"uller space. We define its \emph{boundary} simply as
$$\partial\Teich(S_g) = \overline{\Teich(S_g)} \setminus \Teich(S_g).$$ 
Our aim is now to identify the topology of the Thurston boundary.
The proof of Proposition \ref{embed:Teich:prop} already shows the following.

\begin{prop} \label{alpha:zero:prop}
The Thurston boundary consists of projective currents $[\alpha]$ with $i(\alpha,\alpha)=0$.
\end{prop}

We will later prove that the Thurston boundary actually consists of \emph{all} the projective currents $[\alpha]$ with $i(\alpha,\alpha)=0$. For the moment we content ourselves with some examples. 

\subsection{Weighted pants decompositions}
We construct some geodesic currents based on geodesic pants decompositions. We will need the following.

\begin{ex} \label{pants:simple:line:ex}
Let $P$ be a hyperbolic geodesic pair-of-pants. Every simple geodesic $\gamma$ in the interior of $P$ is open and has both ends that converge to $\partial P$, winding around it.
\end{ex}
\begin{proof}[Hint]
Decompose $P$ into two right-angled hexagons.
\end{proof}

Let a \emph{weighted pants decomposition} for $S_g$ be a pants decomposition $\mu = \gamma_1 \sqcup \ldots \sqcup \gamma_{3g-3}$ equipped with some real numbers $\lambda_1,\ldots,\lambda_{3g-3}\geqslant 0$. A weighted pants decomposition defines a geodesic current $\mu = \lambda_1\gamma_1+\ldots + \lambda_{3g-3}\gamma_{3g-3}$ such that $i(\mu,\mu)=0$.\index{pants decomposition!weighted pants decomposition} 

We will soon prove that $[\mu]$ lies in the Thurston boundary if $\mu\neq 0$. We will need the following lemma that characterises the positively weighted pants decompositions among all currents.

\begin{lemma}
Let $\mu = \sum \lambda_i\gamma_i$ be a weighted pants decomposition with all positive weights $\lambda_i>0$. A current $\alpha \in \calC$ has $i(\alpha,\alpha) = i(\alpha,\mu)=0$ if and only if $\alpha = \sum \lambda_i'\gamma_i$ for some weights $\lambda_i'\geqslant 0$.
\end{lemma}
\begin{proof}
If $\alpha = \sum \lambda_i'\gamma_i$, then clearly $i(\alpha,\alpha) = i(\alpha,\mu)=0$. Conversely, let $\alpha$ be a current with $i(\alpha,\alpha) = i(\alpha, \mu)=0$. Since $\lambda_i>0$, the support of $\mu$ is the preimage in $\matH^2$ of the geodesic pants decomposition. Since $i(\alpha, \mu)=0$, the support of $\alpha$ consists of lines that are either contained or disjoint from that of $\mu$. 

By Exercise \ref{pants:simple:line:ex} a line $l$ of the latter type has both its endpoints at the endpoints of some line of $\mu$, so there are only countably many of them overall; they are not atomic since their projection in $S_g$ is not closed and hence have zero mass. So such lines $l$ do not occur. The support of $\alpha$ is hence contained in that of $\mu$ and the proof is complete.
\end{proof}

We now characterise simple closed curves.

\begin{lemma} \label{characterize:S:lemma}
Let $\gamma\in\calS$ be a simple closed curve and $\mu\in\calC$ be a non-trivial current. If $i(\gamma,\alpha)=0$ implies $i(\mu,\alpha)=0$ for all $\alpha\in\calS$, then $[\mu] = [\gamma]$.
\end{lemma}
\begin{proof}
Extend $\gamma=\gamma_1$ to a geodesic pants decomposition $\beta = \sum_i\gamma_i$. The previous lemma implies that $\mu = \sum_i\lambda_i\gamma_i$ and by letting the pants decomposition vary we get $\mu = \lambda_1\gamma_1$.
\end{proof}

\subsection{Pinching and twisting} \label{pinching:twisting:subsection}
How can we construct sequences of hyperbolic metrics on $S_g$ that converge to some point at infinity? We can do this quite easily by \emph{pinching} a simple geodesic, or by \emph{twisting} along it. 

To define the former, we need to fix some Fenchel--Nielsen coordinates $(l_i, \theta_i)$ for $\Teich(S_g)$ based on a pants decomposition $\mu = \gamma_1\sqcup \ldots \sqcup \gamma_{3g-3}$.
Fix a metric $(l_i,\theta_i) \in \Teich(S_g)$. A \emph{pinching} along $\gamma=\gamma_1$ is any sequence of metrics where the first length coordinate $l_1$ tends to zero and all the other coordinates are kept bounded (from above and below).

\begin{prop}
The limit of a pinching is $[\gamma]$.
\end{prop}
\begin{proof}
By the Collar Lemma \ref{collar:lemma}, the length of a closed geodesic $\alpha$ in the pinched metrics tends to infinity if $i(\gamma,\alpha)>0$ and stays bounded if $i(\gamma,\alpha)=0$. Therefore the pinched metrics tend in $\matP\calC$ to a class $[\mu]$ such that $i(\gamma,\alpha) = 0 \Rightarrow i(\mu,\alpha)=0$ for all $\alpha \in \calS$. Lemma 
\ref{characterize:S:lemma} gives $[\mu] = [\gamma]$.
\end{proof}

Let $m\in \Teich(S_g)$ be any metric and $\gamma$ be a closed geodesic.  We denoted by $m_\theta^\gamma$ the metric obtained from $m$ via an earthquake of angle $\theta$ along $\gamma$, see Section \ref{earthquakes:subsection}.

\begin{prop}
The limit of $m_\theta^\gamma$ as $\theta \to \pm \infty$ is $[\gamma]$.
\end{prop}
\begin{proof}
By Proposition \ref{strettamente:convessa:prop} the length of a closed geodesic $\alpha$ tends to infinity if $i(\gamma,\alpha)>0$ and is constant if $i(\gamma,\alpha)=0$. We conclude as above that the metrics tend to $[\gamma]$.
\end{proof}


To appreciate the qualitative difference between pinching and twisting along $\gamma$, recall the explicit torus case from Section \ref{length:torus:subsection}: if we pinch along $\gamma$, the sequence of metrics converge to $[\gamma]\in \partial \matH^2$ roughly like a geodesic pointing to $[\gamma]$, while if we twist along $\gamma$ it does so along a horosphere centred at $[\gamma]$.

We have discovered in particular that the set $\calS$ is entirely contained in Thurston's boundary. To fully identify the Thurston boundary, we now analyse the geodesic currents $\alpha$ lying in the ``light cone'', that is those with $i(\alpha,\alpha)=0$. We already know that the weighted pants decompositions are there, but there is more: we will soon see that the light cone contains many new fascinating objects called \emph{laminations}. 

\section{Laminations}

In this section we characterize geometrically the currents $\alpha\in\calC$ contained in the ``light cone,'' \emph{i.e.}~those with $i(\alpha,\alpha)=0$. Every such current may be represented as a \emph{measured geodesic lamination}, a closed subset of $S_g$ foliated by geodesics and equipped with a transverse measure. These unexpected objects play an important part in the topology of manifolds of dimension two and three.

\subsection{Geodesic laminations.}
Let $S = \matH^2/_\Gamma$ be a hyperbolic surface. Recall that in this chapter we consider geodesics only as subsets of $S$, neglecting their parametrisation. 

A \emph{geodesic lamination} $\lambda$ is a set of disjoint simple complete geodesics in $S$, whose union is a closed subset of $S$. Each geodesic may be closed or open and is called a \emph{leaf}; their union is the \emph{support} of $\lambda$. We will often confuse $\lambda$ with its support for simplicity, since the support determines the set of geodesics in all the interesting cases (see below).\index{geodesic lamination}

The following examples of geodesic laminations are fundamental:
\begin{itemize}
\item a geodesic multicurve (that is a finite set of disjoint simple closed geodesics) in $S$;
\item a set of disjoint lines in $\matH^2$ whose union is closed. 
\end{itemize}

A lamination in $\matH^2$ may be particularly complicated, see Figure \ref{lamination:fig}. Recall that $\calG = \calG(\matH^2)$ is the set of all lines in $\matH^2$, with its natural topology.

\begin{figure}
\begin{center}
\includegraphics[width = 5 cm] {\iftoggle{BW}{600px-lamination}{600px-lamination}}
\nota{A geodesic lamination in $\matH^2$.}
\label{lamination:fig}
\end{center}
\end{figure}

\begin{ex} 
A set $\lambda$ of disjoint lines in $\matH^2$ forms a closed set in $\matH^2$ if and only if $\lambda$ is closed when considered as a subset of $\calG$.
\end{ex}
If a set of disjoint lines in $\matH^2$ is not closed, it suffices to take its closure to get a lamination. Every lamination in $S=\matH^2/_\Gamma$ lifts to a $\Gamma$-invariant lamination  in $\matH^2$, hence the laminations in $S$ are in natural bijection with the $\Gamma$-invariant laminations in $\matH^2$. 

\subsection{Local behaviour}
A geodesic lamination is often too complicated to be determined with full precision. In fact, its topology is already quite involved, and we now start by examining it locally. Let $\lambda$ be a geodesic lamination on a hyperbolic surface $S=\matH^2/_\Gamma$. The following exercise may be proved by passing to the universal cover $\matH^2$.

\begin{ex} 
Every point $p\in \lambda$ has an open neighbourhood $U$ and a chart $U \to (-1,1) \times (-1,1)$ that sends $p$ to $(0,0)$ and $\lambda \cap U$ to $(-1,1) \times J$ for some closed subset $J\subset (0,1)$.
\end{ex}

What kind of $J\subset (0,1)$ may arise?
It is possible to construct some laminations in $\matH^2$ or in some simple hyperbolic surfaces like cusps and tubes, where their support is the whole manifold. In this case $J$ is the full interval $(0,1)$. These cases are however of no interest for us, because of the following.

\begin{prop} \label{empty:J:prop}
If $S$ has finite volume, then $J$ has empty interior.
\end{prop}
\begin{proof}
Suppose by contradiction that $J$ contains some open interval. This open interval determines a set of leaves of $\lambda$, whose lift in $\matH^2$ forms an open subset $U\subset \matH^2$ foliated into lines. Its limit in $\partial \matH^2$ has non-empty interior $S$. By Corollary \ref{parabolic:hyperbolic:dense:cor} there is a hyperbolic transformation $\gamma \in \Gamma$ with an attracting limit point in $S$. This implies easily that there are $\gamma$-translates of some lines in $U$ that intersect some lines of $U$ transversely, a contradiction. 
\end{proof}

\begin{cor}
If $S$ has finite volume, every geodesic lamination $\lambda$ in $S$ has empty interior and is determined by its support.
\end{cor}

\subsection{Complementary regions}
Let $\lambda \subset S$ be a geodesic lamination in a hyperbolic surface $S=\matH^2/_\Gamma$.
A \emph{complementary region} (shortly, a \emph{region}) is a connected component of the open complement $S\setminus \lambda$. The abstract completion of a complementary region is a hyperbolic surface with non-empty geodesic boundary consisting of lines and/or circles. 

In general, a complementary region may have genus and/or infinitely many boundary components. Note that in any case there are at most countably many regions, and hence only countably many leaves of $\lambda$ may be incident to a region: if $\lambda$ has uncountably many leaves (this will be the typical case), most leaves are not incident to any region (compare this with the Cantor subset in $[0,1]$, which contains uncountably many points, but only countably many of them are adjacent to some complementary open segment).

\begin{prop} \label{finite:regions:prop}
A geodesic lamination $\lambda\subset S_g$ in a closed hyperbolic surface $S_g=\matH^2/_\Gamma$ has at most $4g-4$ complementary regions. The boundary of each region has finitely many components.
\end{prop}
\begin{proof}
Every complementary region has area at least $\pi$ (the area of an ideal triangle) and a region with infinitely many boundary components has infinite area. By Gauss-Bonnet we get $\Area(S_g) = -2\pi\chi(S_g) = (4g-4)\pi$.
\end{proof}

We say that a lamination $\lambda \subset S_g$ is \emph{full} if every region is an ideal polygon.\index{geodesic lamination!full}



\subsection{Transverse measures}
Let $\lambda\subset S$ be a geodesic lamination in a hyperbolic surface $S$. A \emph{transverse arc} to $\lambda$ is the support of a simple regular curve $\alpha\colon [a,b] \to S$ transverse to each leaf of $\lambda$, whose endpoints $\alpha(a)$ and $\alpha(b)$ are not contained in $\lambda$.\index{geodesic lamination!measured geodesic lamination} 

\begin{defn} \label{measure:axioms:defn}
A \emph{transverse measure} for a lamination $\lambda\subset S$ is a locally finite Borel measure $L_\alpha$ on each transverse arc $\alpha$ such that:
\begin{enumerate}
\item if $\alpha' \subset \alpha$ is a sub-arc of $\alpha$, the measure $L_{\alpha'}$ is the restriction of $L_\alpha$;
\item the support of $L_\alpha$ is $\alpha\cap \lambda$;
\item the measure is invariant through isotopies of transverse arcs.
\end{enumerate}
\end{defn}
In particular every arc $\alpha$ transverse to $\lambda$ has a finite \emph{length}, defined as the total measure of the arc. The arc has length zero if and only if $\alpha\cap\lambda = \emptyset$. A \emph{measured geodesic lamination} is a geodesic lamination equipped with a transverse measure.

\begin{example} \label{weight:example}
A geodesic multicurve $\lambda \subset S$ has a natural transverse measure: for any transverse arc $\alpha$, the measure $L_\alpha$ on $\alpha$ is just the Dirac measure supported on the finite set $\alpha\cap\lambda$.

More generally, we may assign a positive weight $a_i>0$ at each component $\lambda_i$ of $\lambda$ and define a measured geodesic lamination by giving the weight $a_i$ at each intersection $\alpha \cap \gamma_i$. By varying weights we get distinct measured geodesic laminations with the same support.
\end{example}

\subsection{Currents and measured geodesic laminations}
Let $S_g=\matH^2/_\Gamma$ be a closed hyperbolic surface. We now construct a natural bijection between the measured geodesic laminations on $S_g$ and the geodesic currents $\lambda$ with $i(\lambda,\lambda)=0$, \emph{i.e.}~those lying in the ``light cone.'' We see a measured geodesic lamination in $S_g$ as a $\Gamma$-invariant measured geodesic lamination in $\matH^2$.

Let $\lambda \in \calC$ be a geodesic current with $i(\lambda,\lambda)=0$. By Remark \ref{correnti:oss} the support of $\lambda$ is a closed $\Gamma$-invariant subset of $\calG$ formed by disjoint lines, hence a $\Gamma$-invariant lamination in $\matH^2$. The geodesic current induces also a transverse measure as follows. Let $\alpha$ be an arc transverse to the lamination $\lambda$. Up to cutting $\alpha$ in finitely many arcs we may suppose that it intersects each leaf of $\lambda$ in at most one point. We define the measure of a Borel set $U\subset \alpha$ as the $\lambda$-measure in $\calG$ of the lines in $\lambda$ that it intersects.
\begin{ex}
This transverse measure satisfies the axioms of Definition \ref{measure:axioms:defn} and gives $\lambda$ the structure of a $\Gamma$-invariant measured geodesic lamination.
\end{ex}
\begin{prop}
We have just defined a bijection
$$\big\{ {\rm currents\ } \lambda {\rm \ with\ } i(\lambda,\lambda)=0 \big\}
\longleftrightarrow 
\big\{ {\rm measured\ geodesic\ laminations\ in\ } S_g\big\}
$$
\end{prop}
\begin{proof}
We define the inverse map by transforming a $\Gamma$-invariant measured geodesic lamination $\lambda$ into a geodesic current with support $\lambda \subset \calG$ as follows. For every leaf $l$ of $\lambda$, we pick a small transverse arc $\alpha$ that intersects all the leaves in a neighbourhood $U_l$ of $l$ in $\lambda$, each once. The measure on $\alpha$ translates into a measure on $U_l$, and all the measures on these small sets $U_l$ match to yield a measure on $\lambda$ thanks to Proposition  \ref{ricoprimento2:prop}. 
\end{proof}

We denote by $\calM\calL\subset \calC$ the set of all the measured geodesic laminations on $S_g$, henceforth identified with the currents lying in the ``light cone''. With this identification $\calML$ is independent on the auxiliary hyperbolic metric on $S_g$.\index{$\calML$}

Recall that a multicurve on $S_g$ is a finite collection of disjoint non-trivial simple closed curves. A multicurve determines a measured geodesic lamination: consider $n$ parallel components as a single one with weight $n$, take the geodesic representative of the resulting essential multicurve (keeping the same weights), and use Example \ref{weight:example}. We obtain the inclusions
$$\calS \subset \calM \subset \ML \subset \calC$$ 
where $\calS$ and $\calM$ are the sets of simple closed curves and multicurves in $S_g$.\index{$\calM$} 

\subsection{Some properties}
We can now use our knowledge of the geodesic currents to get some information about the measured geodesic laminations. 

\begin{prop} \label{at:most:two:prop}
Let $\lambda\subset S_g$ be a measured geodesic lamination. For every $p\in\partial \matH^2$, there are at most two lines in the preimage $\tilde\lambda\subset\matH^2$ incident to $p$, and if there are two, none of them projects to a closed geodesic in $\lambda$. 
\end{prop}
\begin{proof}
If there are at least three lines $l_1, l_2, l_3$ incident to $p$, with $l_2$ lying strictly between $l_1$ and $l_3$, we may find a small arc $\alpha$ transverse to $\tilde \lambda$, whose support lies between $l_1$ and $l_3$, that intersects $l_2$ and possibly other lines incident to $p$, and no other line of $\tilde \lambda$ (use Proposition \ref{empty:J:prop} near $l_2$ here). The lines intersected by $\alpha$ form a pencil with positive mass. By Proposition \ref{pencils:prop} the pencil consists of a single line, the axis of a hyperbolic element $\varphi\in \Gamma$.

On the other hand, if $p$ is the endpoint of the axis of some $\varphi$, then no other leaf $l\subset \tilde\lambda$ can be incident to $p$, otherwise infinitely many lines $\varphi^i(l)$ would also be, and the same argument above would give a contradiction.
\end{proof}

\begin{ex}
There are some geodesic laminations in $S_g$ that do not support any transverse measure.
\end{ex}
\begin{proof}[Hint]
Construct a geodesic lamination $\lambda \subset S_g$ consisting of some disjoint simple closed curves and finitely many open leaves that wind around them (see Section \ref{incomplete:metrics:subsection}). The pre-image $\tilde\lambda\subset\matH^2$ contradicts Proposition \ref{at:most:two:prop}.
\end{proof}


\subsection{Intersection form}
What is the geometric meaning of the intersection form $i$ of two geodesic currents, when applied to a measured geodesic lamination? Let $\lambda$ be a measured geodesic lamination in a hyperbolic closed surface $S_g=\matH^2/_\Gamma$. Three cases are relevant for us:

\begin{enumerate}
\item If $\gamma$ is a simple closed geodesic, then either $\gamma$ is contained in $\lambda$, and we get $i(\gamma, \lambda)=0$, or $\gamma$ is disjoint from $\lambda$, and we still get $i(\gamma,\lambda)=0$, or $\gamma$ is transverse to $\lambda$, and in that case $i(\lambda, \gamma)$ equals the full measure $L(\gamma)$ of $\gamma$ as a transverse curve to $\lambda$.
\item If $\lambda'$ is another measured geodesic lamination, then $i(\lambda, \lambda')$ ``counts'' the (possibly infinitely many) transverse intersections of $\lambda$ and $\lambda'$ with respect to the product of the two transverse measures.
\item If $m$ is a hyperbolic metric, then $i(\lambda, m)$ measures the ``length'' of $\lambda$ with respect to $m$, that is the integral over $\lambda$ of the usual length weighted with the transverse measure.
\end{enumerate}

We now consider a non-geodesic non-trivial simple closed curve $\gamma$ transverse to $\lambda$, and we denote by $L(\gamma)$ its full transverse $\lambda$-measure. The quantities $L(\gamma)$ and $i(\gamma, \lambda)$ need \emph{not} to be equal in this case (because $\gamma$ is not geodesic), but this phenomenon is easily controlled by a suitable version of the bigon criterion.

\begin{prop} \label{m:i:prop}
We have $L(\gamma) \geqslant i(\gamma, \lambda)$ and the equality holds $\Longleftrightarrow$ $\gamma$ forms no bigon with any leaf of $\lambda$.
\end{prop}
\begin{proof}
We follow the proof of Theorem \ref{bigono:teo}. Let $\tilde \lambda$ be the counterimage of $\lambda$ in $\matH^2$. A lift $\tilde\gamma$ of $\gamma$ has distinct endpoints in $\partial \matH^2$ and $\tilde \gamma$ must intersect the leaves of $\tilde\lambda$ whose endpoints are linked with them. The curve $\tilde \gamma$ intersects each such leaf only once and no other leaf $\Leftrightarrow$ $\tilde \gamma$ forms no bigon with any leaf of $\tilde \lambda$ $\Leftrightarrow$ $\gamma$ forms no bigon with any leaf of $\lambda$.
\end{proof}

\begin{prop} \label{sublamination:prop}
If $\lambda \in \calML$ is full, we have $i(\lambda, \alpha)>0$ for every current $\alpha$ that is not a measured sublamination of $\lambda$.
\end{prop}
\begin{proof}
If $i(\lambda, \alpha)=0$ then every line $l$ in the support of $\alpha$ is either a leaf of $\lambda$ or an infinite diagonal in some complementary ideal polygon; a diagonal would be an atomic point in $\alpha$ that would cover a closed geodesic and could be added to $\lambda$ contradicting 
Proposition \ref{at:most:two:prop}.
\end{proof}

\subsection{Train tracks}
How can one construct a measured geodesic lamination, concretely? There are some nice combinatorial tools designed to this purpose, called train tracks.\index{train track} 

\begin{figure}
\begin{center}
\includegraphics[width = 11 cm] {\iftoggle{BW}{switch-BW}{switch}}
\nota{A switch contains $a\geqslant 1$ and $b\geqslant 1$ branches on each side: here $a=3$ and $b=2$ (1). The complementary regions of a train track contain no discs, annuli, monogons, or bigons (2).}
\label{switch:fig}
\end{center}
\end{figure}

\begin{figure}
\begin{center}
\includegraphics[width = 8 cm] {\iftoggle{BW}{Triple-torus-train-track-BW}{Triple-torus-train-track}}
\nota{A train track on the genus-three surface $S_3$. }
\label{train_track:fig}
\end{center}
\end{figure}

A \emph{train track} in a closed surface $S_g$ is a closed subset $\tau\subset S_g$ built by taking a finite set of points (called \emph{vertices} or \emph{switches}) and joining them with disjoint arcs called \emph{branches}. We require that every switch looks locally like Figure \ref{switch:fig}-(1): all branches have the same tangent line, and there is at least one branch at each side. 

The \emph{complementary regions} of $\tau$ are the abstract closures of the connected components of $S_g\setminus \tau$. Each complementary region is naturally a compact surface with finitely many singular vertices in its boundary having ``zero interior angle''. In the definition of a train track, we also require that no complementary region be a disc, an annulus, a monogon, or a bigon, as in Figure \ref{switch:fig}. An example is shown in Figure \ref{train_track:fig}. The exclusion of these complementary regions is essential to get the following.

\begin{ex}
A train track contains at most $-6\chi(S_g)$ switches and at most $-2\chi(S_g)$ regions; it contains $-6\chi(S_g)$ switches if and only if every switch is trivalent and every complementary region is a triangle. The sphere and the torus contain no train track. 
\end{ex}
\begin{proof}[Hint] Make the appropriate Euler characteristic computation.
\end{proof}

\begin{ex}
Construct a train track with $-2\chi(S_g)$ triangular complementary regions in $S_g$ for every $g\geqslant 2$. 
\end{ex}
\begin{proof}[Hint] Take a pants decomposition and add three arcs (tangent to the boundary) inside each pants to decompose it into two triangular regions.
\end{proof}

We obtain a finiteness property.

\begin{cor}
There are only finitely many train tracks $\tau$ in $S_g$ up to diffeomorphism (but not up to isotopy!).
\end{cor}
\begin{proof}
There are only finitely many combinatorial types for $\tau$, its complementary regions, and the way they are adjacent.
\end{proof}

\subsection{Weight systems}

A \emph{weight system} on a train track $\tau\subset S_g$ is the assignment of a non-negative real number, called \emph{weight}, to each branch of $\tau$, such that at every vertex the \emph{switch condition} holds: the sum of the weights on the left branches should be equal to the sum of the weights on the right branches. 

Let $w$ be an \emph{integral} weight system, that is one whose weights are all integers. The integral weight system $w$ determines a closed 1-manifold in $S_g$ as follows: replace every branch with weight $a$ with $a$ parallel copies of it, and pair all the branches at every switch in the obvious way (this is possible thanks to the switch condition).

\begin{prop}
The resulting closed 1-manifold is a multicurve.
\end{prop}
\begin{proof}
Since no complementary region of $\tau$ is a disc or monogon, the complementary regions of the closed 1-manifold have non-positive Euler characteristic (exercise).
\end{proof}

Our aim now is to use the train tracks to parametrise all the multicurves in $S_g$. Later on, we will extend everything from multicurves to measured geodesic laminations.

\subsection{A parametrisation for $\calM$} \label{parametrisation:M:subsection}

\begin{figure}
\begin{center}
\includegraphics[width = 11 cm] {\iftoggle{BW}{pantaloni_anelli_ex-BW}{pantaloni_anelli_ex}}
\nota{Decompose $S$ into pairs of pants and annuli, so that each simple closed curve separates a pair of pants from an annulus. We mark a \iftoggle{BW}{}{blue }point in each curve.}
\label{pantaloni_anelli_ex:fig}
\end{center}
\end{figure}

\begin{figure}
\begin{center}
\includegraphics[width = 6 cm] {\iftoggle{BW}{pantaloni_anelli_new2-BW}{pantaloni_anelli_new2}}
\nota{Identify each pair of pants and annulus with one of these two fixed models.}
\label{pantaloni_anelli:fig}
\end{center}
\end{figure}

Let $S_g$ have genus $g\geqslant 2$. We now construct finitely many train tracks that parametrise all the multicurves in $S_g$, and to this purpose we fix a frame similar to the one used in the Fenchel-Nielsen parametrisation of the Teichm\"uller space. 

We decompose $S_g$ into pairs of pants and annuli as in Figure \ref{pantaloni_anelli_ex:fig}. We mark \iftoggle{BW}{}{(in blue) }an arbitrary point in each closed curve and we also fix once for all a diffeomorphism (preserving orientation and marked points) between each pair of pants and annulus with the corresponding model in Figure \ref{pantaloni_anelli:fig}. The resulting parametrisation of $\calM$ will depend also on these fixed diffeomorphisms.

Let a \emph{colouring} be the assignment of a triple $a_i,b_i,c_i$ of non-negative numbers to each annulus of the decomposition, such that one of the following equalities holds:
$$a_i = b_i+c_i, \quad b_i = c_i+a_i, \quad c_i = a_i +b_i.$$
We assign the colour $a_i$ to each of the two curves that bound the annulus.
Since there are $-\frac 32\chi(S_g) = 3g-3$ annuli, we get $-\frac 92\chi(S_g) = 9g-9$ colours overall. The colouring is \emph{integral} if the colours are integers and also the following holds: 
for each pair of pants, the sum $a_i+a_j+a_k$ of the colours of its boundary curves must be even.

\begin{figure}
\begin{center}
\includegraphics[width = 12.5 cm] {\iftoggle{BW}{pants_train_tracks-BW}{pants_train_tracks}}
\nota{The portion of weighted train track in a pair of pants $P$ is determined by the point $Q=[a_i,a_j,a_k]\in\matRP^2$. Having non-negative entries, the point $Q$ is contained in the triangle with vertices $[1,0,0], [0,1,0], [0,0,1]$. This triangle subdivides into four sub-triangles by adding the vertices $[0,1,1]$, $[1,0,1]$, $[1,1,0]$ and the shape of the train track depends on which sub-triangle contains $Q$. The three edges of the train track are given the weights indicated in the figure: these are the only weights that sum to $a_i, a_j, a_k$ at the boundaries. When $Q$ lies in the frontier of two or more triangles some branch has weight zero and we delete it: the different shapes prescribed by the adjacent triangles coincide after deleting this branch.}
\label{pants_train_tracks:fig}
\end{center}
\end{figure}

We now transform each colouring into a weighted train track $\tau$ as follows. The boundary components of each pair-of-pants $P$ are coloured by some triple $a_i, a_j, a_k$. We insert in $P$ a portion of train track as in Figure \ref{pants_train_tracks:fig}: its shape depends on the position of the point $[a_i,a_j,a_k]$ in $\matRP^2$ and its weights depend linearly on $a_i, a_j, a_k$. Note that there are finitely many possible shapes overall.

\begin{figure}
\begin{center}
\includegraphics[width = 8 cm] {\iftoggle{BW}{anelli_train_tracks-BW}{anelli_train_tracks}}
\nota{The portion of weighted train track in each annulus $A$ is determined by the point $Q=[a,b,c] = [a_i,b_i,c_i]\in \matRP^2$. By hypothesis $Q$ is contained in the boundary of the triangle with vertices $[1,0,0], [0,1,0], [0,0,1]$. The boundary subdivides into three segments and the shape of the train track depends on which segment contains $Q$. When $Q$ is a vertex some branch has weight zero and we delete it.}
\label{anelli_train_tracks:fig}
\end{center}
\end{figure}

We extend the train track inside each annulus $A$ coloured with $a_i, b_i, c_i$ as shown in Figure \ref{anelli_train_tracks:fig}. Again, the shape of the portion depends on the position of $[a_i,b_i,c_i]$ in $\matRP^2$. As a result we get a train track $\tau$ with integral weights and hence a multicurve.

\begin{prop} \label{parametrize:multicurves:prop}
The construction produces a bijection
$$\big\{ {\rm integral\ colourings}\big\} \longleftrightarrow \calM.$$
\end{prop}
\begin{proof}
We start by proving surjectivity. Given a multicurve $\mu\in \calM$, we minimise its intersections with our decomposition into pants and annuli: now $\mu$ intersects every pair of pants $P$ into non-trivial simple closed curves and arcs (an arc is trivial if it forms a bigon with $\partial P$). 

By decomposing $P$ into two hexagons, and minimising the intersections of the curves and arcs with them, one deduces easily that $P$ contains only three isotopy classes of non-trivial closed curves (one parallel to each component of $\partial P$) and six isotopy classes of non-trivial arcs (one connecting every pair of components of $\partial P$). Note that isotopies are \emph{not} required to fix $\partial P$ pointwise.

We isotope all the closed curves away from $P$ and place them inside the adjacent annuli, so there are only arcs. As above, we check easily that the isotopy class of the whole set of arcs in $P$ is determined by the intersection numbers $a_i, a_j, a_k$ with the components of $\partial P$ and is hence represented by the corresponding train track portion from Figure \ref{pants_train_tracks:fig}. 

We now turn to annuli. On each annulus $A$, there is only one isotopy class of non-trivial arcs and one isotopy class of non-trivial closed curves. However, an arc may wind many times around $A$ and we are not allowed to use isotopies that move $\partial A$ anymore because the configuration in $P$ is already fixed: one deduces easily that the triples $(a,b,c)$ from Figure \ref{anelli_train_tracks:fig} describe precisely all the possible configurations. 

The injectivity is surprisingly simple to prove: the configuration of $\mu$ that minimises its intersections with the decomposition into pairs of pants and annuli is unique thanks to Proposition \ref{minimale:unica:prop} and hence the numbers $a_i,b_i,c_i$ are easily determined by $\mu$. The proof is complete.
\end{proof}

We have found a complete combinatorial parametrisation for the set $\calM$, and we now want to extend it to $\calML$.

\subsection{Topological laminations} \label{topological:laminations:subsection}

\begin{figure}
\begin{center}
\includegraphics[width = 11 cm] {\iftoggle{BW}{switch_weights-BW}{switch_weights}}
\nota{How to construct a lamination from a weighted train track. Here we have a switch with weights $a$, $b$, $c$ (1) and the switch condition $a=b+c$ allows all the Euclidean rectangles (2) to glue (3).}
\label{switch_weights:fig}
\end{center}
\end{figure}

A train track $\tau$ in $S_g$ with integer weights parametrises a multicurve. We now show that, more generally, a weighted train track parametrises a measured geodesic lamination. 

The construction goes as follows. First, we remove all the branches with zero weight, and we replace every branch of weight $a>0$ with a Euclidean rectangle of width $a$ and with arbitrary length as in Figure \ref{switch_weights:fig}-(2). Thanks to the switch conditions, these rectangles glue nicely at each switch as in Figure \ref{switch_weights:fig}-(3).

In Figure \ref{switch_weights:fig}-(3) we see that all the lines from left- and right-rectangles are matched in a 1-1 correspondence, with a finite number of exceptions. At each exception, some $m \geqslant 1$ left-lines are matched to some $n \geqslant 1$ right-lines at some point, and we have $m+n=3$ or 4. We call these lines and points \emph{singular}. 

After gluing the rectangles we get a closed subset $\lambda$ of $S_g$ foliated by lines, with finitely many singular points and lines. We now eliminate the singular points and lines by cutting $\lambda$ carefully along them: if $m+n=4$ this amounts simply to doubling the singular point; if $m+n=3$ we double the singular leaf by opening a small open corridor in the foliation, starting from the singular point and digging along the leaf. Since the singular leaf may be non-compact, for the process to converge in $S_g$ we need to shrink the width of the corridor sufficiently fast. 

After this cut we get a closed subset of $S_g$, which we still name $\lambda$, nicely partitioned into disjoint lines that may be either open or closed: we call it a \emph{topological lamination}. After the cut every rectangle of type $L\times [0,a]$ as in Figure \ref{switch_weights:fig} transforms into a set $L \times J$ where $J$ is obtained by cutting $[0,a]$ along (at most) countably many points. Note that $J$ may be a Cantor set.\index{topological lamination}

The set $J$ inherits from $[0,a]$ a Borel measure with total mass $a$. This measure gives a transverse measure to the leaves of the rectangle $L\times J$, and $\lambda$ inherits the structure of a measured topological lamination (the notion of a transverse measure is exactly the same as in the geodesic case).

It only remains to promote the topological lamination to a geodesic lamination: this is usually done by \emph{straightening} its leaves. We describe this procedure only for a particular class of train tracks: the ones introduced in the previous section to parametrise $\calM$.

\subsection{A parametrisation for $\ML$} \label{parametrizzazione:subsection}
We now extend the arguments of Section \ref{parametrisation:M:subsection}
from multicurves to laminations.

\begin{figure}
\begin{center}
\includegraphics[width = 11 cm] {\iftoggle{BW}{foglia_new-BW}{foglia_new}}
\nota{The decomposition into pants and annuli of $S_g$ is homotopic to a geodesic pants decomposition (with each curve counted twice) whose counterimage in $\matH^2$ consists of infinitely many disjoint ultraparallel lines as in the picture.
We show that the lift of a leaf $l$ of $\lambda$ in $\matH^2$ has two distinct endpoints in $\partial \matH^2$: if $l$ is a closed curve in an annulus $A$, then it is isotopic to a closed geodesic and we are done; if $l$ intersects annuli and pairs of pants into non-trivial arcs, its lift intersects the ultraparallel lines forming no bigons, hence it intersects each line at most once: therefore it has disjoint limits at $\pm \infty$ (left). The straightening replaces the lift with the unique line with these endpoints (right).}
\label{foglia_new:fig}
\end{center}
\end{figure}

We fix a decomposition of $S_g$ into 
pants and annuli. Every colouring $(a_i,b_i,c_i)$ produces a weighted train track, which in turn parametrises a measured topological lamination $\lambda$. We fix a hyperbolic structure $S_g = \matH^2/_\Gamma$. 

The lamination $\lambda$ can be \emph{straightened} to a measured geodesic lamination $\bar\lambda$ as follows. A leaf $l$ of $\lambda$ is either a closed curve in some annulus $A$, or it intersects every annulus $A$ and pair of pants $P$ in non-trivial arcs: in both cases, every lift of $l$ in $\matH^2$ is a curve with two distinct limit endpoints in $\partial \matH^2$, see Figure \ref{foglia_new:fig}. We replace every lift of $l$ with the unique line having these endpoints, and we do this for every leaf $l$ of $\lambda$. 

The result is a $\Gamma$-invariant closed set of disjoint lines in $\matH^2$ that project to a geodesic lamination $\bar\lambda$ in $S_g$. The transverse measure on $\lambda$ easily induces one on the straightened $\bar\lambda$: it suffices to consider transverse arcs contained in the decomposition into pairs of pants and annuli. Note that some parallel closed leaves of $\lambda$ may have collapsed to a single atom closed geodesic in $\bar\lambda$.

\begin{prop} \label{parametrize:laminations:prop}
The construction induces a bijection
$$\big\{ {\rm colourings}\big\} \longleftrightarrow \calML.$$
\end{prop}
\begin{proof}
We adapt the proof of Proposition \ref{parametrize:multicurves:prop}, starting with surjectivity: given a measured geodesic lamination $\lambda$, we determine a colour $(a_i,b_i,c_i)$ representing it.

The decomposition into pairs of pants and annuli is homotopic to a geodesic pants decomposition $\mu=\gamma_1\sqcup \ldots \sqcup \gamma_{3g-3}$ (with each curve counted twice). We set $a_i = i(\gamma_i,\lambda)$.

The measured geodesic lamination $\lambda \in \calML$ decomposes as $\lambda = \lambda_0 \sqcup \lambda_1$ where $\lambda_0 = \sum_ik_i\gamma_i$ is a weighted pants decomposition with support in $\mu$ and $\lambda_1$ is transverse to $\mu$. If $k_i\neq 0$ then necessarily $a_i=0$ and we set $b_i=c_i=k_i$.

At every pair of pants $P$ the three colours $a_i,a_j,a_k$ determine a portion of weighted train track as prescribed by Figure \ref{pants_train_tracks:fig}. The intersection $\lambda_1\cap P$ consists of geodesic and hence non-trivial arcs: therefore the weighted train track describes faithfully this portion of measured geodesic lamination, up to isotopy. The numbers $b_i$ and $c_i$ are then determined by the way these two portions wind and match along $\gamma_i$. The colours $(a_i,b_i,c_i)$ parametrise $\lambda$.

The choices of the colours $(a_i,b_i,c_i)$ were forced by the intersection of $\lambda$ with the two pants $P$ adjacent to $\gamma_i$ (exercise). This shows injectivity.
\end{proof}

We have parametrised $\calML$, and we now investigate its topology. The space of all colourings $(a_i,b_i,c_i)$ forms the subset $C\times \ldots \times C \subset \matR^3 \times \ldots \times \matR^3$ where $C \subset \matR^3$ is the cone based on the origin over the sides of the triangle with vertices $(1,1,0)$, $(0,1,1)$, and $(1,0,1)$, and is homeomorphic to $\matR^2$. The space of colourings is homeomorphic to $\matR^{2\cdot (3g-3)} = \matR^{6g-6}$, that is to the Teichm\"uller space itself!

The parametrisation identifies $\calML$ with  $C\times \ldots \times C$. We want to show that this identification is a homeomorphism, and to this purpose we study the intersection form.

\subsection{Intersection form}
We now show that the intersection form $i$ between curves and laminations has a surprisingly simple behaviour after that we parametrise the space $\calML$ as $C\times \ldots \times C$. 

We note that $C\times \ldots \times C$ is a piecewise-linear object in $\matR^{9g-9}$, that is it is the support of a simplicial complex. Recall that a continuous map between simplicial complexes is \emph{simplicial} if it sends each simplex onto a simplex in an affine linear way, and it is \emph{piecewise linear} if it restricts to a simplicial map on some subdivisions. Finally, the map is \emph{half-integral piecewise linear} if all the affine linear maps have half-integer (possibly integer) coefficients.

\begin{prop} \label{integral:PL:map:prop}
For every $\gamma\in\calS$ the map
\begin{align*}
C\times \ldots \times C & \to \matR \\
\lambda & \longmapsto i(\gamma, \lambda)
\end{align*}
is half-integral piecewise linear.
\end{prop}
\begin{proof}
A colour $(a_i,b_i,c_i)\in C\times \ldots \times C$ defines a weighted train track $\tau$ and hence a measured topological lamination $\lambda$, obtained by substituting every edge of $\tau$ with some product $L\times J$ of horizontal leaves. Thanks to Proposition \ref{m:i:prop} the intersection $i(\gamma, \lambda)$ is realized by some representative $\gamma$ which is either a leaf of $\lambda$ or is transverse to $\lambda$ and forms no bigon with any leaf of $\lambda$.

\begin{figure}
\begin{center}
\includegraphics[width = 12.5 cm] {\iftoggle{BW}{monotonic-BW}{monotonic}}
\nota{A monotonic curve $\gamma$ intersects the products $L \times J$ in monotonic paths (left) and forms no bigons outside the products (right).}
\label{monotonic:fig}
\end{center}
\end{figure}

We say that a closed curve transverse to $\lambda$ is \emph{monotonic} if 
\begin{enumerate}
\item it intersects every product $L\times J$ into arcs that are monotonic in both coordinates, and the same monotonicity is preserved when the curve goes from one product $L\times J$ to an adjacent one as in Figure \ref{monotonic:fig}-(left);
\item it makes no bigons outside the rectangles as in Figure \ref{monotonic:fig}-(right).
\end{enumerate}

Concerning (1), a horizontal path in $L\times J$ disjoint from $\lambda$ is allowed, but it must keep being horizontal on the adjacent products. 

\begin{figure}
\begin{center}
\includegraphics[width = 12.5 cm] {\iftoggle{BW}{monotonic2-BW}{monotonic2}}
\nota{If a monotonic curve $\gamma$ forms a bigon with a leaf $\alpha$, by cutting carefully $\tau$ along $\gamma$ we construct a portion of train track on this bigon containing $\alpha$ and with truncated edges exiting orthogonally from $\gamma$ (left). We then simplify the triangles adjacent to $\gamma$ (right).}
\label{monotonic2:fig}
\end{center}
\end{figure}

A simple closed curve that forms no bigons with any leaf of $\lambda$ can be easily isotoped to be monotonic. Conversely, we now prove that a monotonic curve $\gamma$ forms no bigons with any leaf of $\lambda$. If it did, by cutting $\tau$ along $\gamma$ we would get an abstract bigon with a portion of train track as in Figure \ref{monotonic2:fig}-(left). If $\gamma$ is adjacent to some triangles as in Figure \ref{monotonic2:fig}-(right) we close them as shown there. Finally, by doubling the bigon along $\gamma$ we get a train track on a disc, and by doubling again we build a train track in a sphere, which is absurd (one checks easily that no complementary region is a disc, annulus, monogon, or bigon: there are no bigons because we have closed the triangles as in Figure \ref{monotonic2:fig}-(right)).

The $\lambda$-transversal length $L(\gamma)$ of a monotonic $\gamma$ is the sum of the lengths of the sub-paths intersecting sequences of products $L\times J$ as in Figure \ref{monotonic:fig}-(left). One checks easily that the length of each sub-path is a half-integer linear combination of the colours $(a_i,b_i,c_i)$. Therefore $L(\gamma) = i(\gamma,\lambda)$ is a half-integer combination of the colours.

If we vary the colouring $(a_i,b_i,c_i)$ a little, the curve $\gamma$ keeps being monotonic in the same way as before, except when a portion of $\gamma$ is horizontal: in that case the new monotonicity depends on how the colouring varies. In all cases, the new $\gamma$ is still monotonic, hence bigonless, hence $L(\gamma)=i(\gamma, \lambda)$ again. 

If there are no horizontal portions in $\gamma$, the length $L(\gamma)$ varies linearly with the same integral formula found above. If there are horizontal portions, then the colour lies in a hyperplane of $\matR^{9g-9}$, and there are two different linear formulas joining there. Hence $i(\gamma,\lambda)$ is half-integral piecewise linear.
\end{proof}

The map $\lambda \mapsto i(\gamma, \lambda)$ is also obviously homogeneous, in the sense that $i(\gamma, t\lambda) = ti(\gamma, \lambda)$ for all $t\geqslant 0$. 

The proof of Proposition \ref{integral:PL:map:prop} contains a recipe for calculating $i(\gamma, \lambda)$. It suffices to put $\gamma$ in monotonic position and then add the contributions of each monotonic arc. For instance, the $i$-th curve $\gamma_i$ of the pants decomposition has an obvious monotone position giving $i(\gamma_i,\lambda) = a_i$.

The colours $b_i$ and $c_i$ do not have such an immediate description, but the following holds anyway.

\begin{ex} \label{abc:ex}
There exist finitely many simple closed curves $\gamma_j$ whose intersections $i(\lambda, \gamma_j)$ determine the colouring $(a_i,b_i,c_i)$ representing $\lambda$ in a half-integral piecewise-linear continuous fashion.
\end{ex}
\begin{proof}[Hint] Choose some additional curves as in Proposition \ref{9g-9:prop} 
and compute their intersections with $\lambda$.
\end{proof}

\subsection{The Thurston boundary}
We can finally determine the topology of the measured geodesic laminations space $\calML$ and of the Thurston boundary of the Teichm\"uller space.

\begin{prop}
The colouring parametrisation induces a homeomorphism $\calML \isom \matR^{6g-6}$.
\end{prop}
\begin{proof}
The map $\calML \to C\times \ldots \times C \isom \matR^{6g-6}$ is continuous by Exercise \ref{abc:ex}. It is proper: if the colours $(a_i,b_i,c_i)$ stay bounded, the intersections $i(\lambda, \gamma_i)$ with finitely many filling curves $\gamma_i$ stay bounded, and hence $\lambda$ moves in a compact set by Proposition \ref{compattezza:C:prop}. The continuous map is a homeomorphism by Proposition \ref{propria:chiusa:cor}.
\end{proof} 

Let now $\matP\calML \subset \matP\calC$ be the set of all \emph{projective measured laminations}, that is the image of $\calML\setminus 0$ in $\matP\calC$. The set $\matP\calML$ contains the projective simple closed curves $\matP\calS$ and the projective multicurves $\matP\calM$.\index{geodesic lamination!projective measured geodesic lamination}

\begin{teo}
The following homeomorphism holds
$$\partial \Teich(S_g) = \matP\calML \isom S^{6g-7}.$$ 
The set $\matP\calS$ is dense in $\matP\calML$.
\end{teo}
\begin{proof}
The homeomorphism $\calML \isom \matR^{6g-6}$ induces $\matP\calML \isom S^{6g-7}$. Concerning the Thuston boundary, we know from Section \ref{pinching:twisting:subsection} and Proposition \ref{alpha:zero:prop} that 
$$\matP\calS \subset \partial \Teich(S_g) \subset \matP\calML.$$
We now show that $\matP\calS$ is dense in $\matP\calML$, and this concludes the proof. 

First, we prove that $\matP\calM$ is dense in $\matP\calML$. Rational colours $(a_i,b_i,c_i)$ form a dense subset of $C\times \ldots \times C$ and project to a dense subset of $\matP\calML$. Every rational colour is a multiple of an integer colour, which represents a multicurve.

Second, we prove that $\matP\calS$ is dense in $\matP\calM$. Let $\mu$ be a multicurve in $S_g$. Up to acting via $\MCG(S_g)$ we may suppose that $\mu$ is supported on the pants decomposition $\gamma_1,\ldots, \gamma_{3g-3}$ used to define $C\times \ldots \times C$, hence $\mu$ is defined by some integral colour $(a_i,b_i,c_i) = (0,b_i,b_i)$. Now pick any simple closed curve $\gamma$ with $a_i' = i(\gamma, \gamma_i)>0$ for all $i$, represented by some colours $(a_i', b_i', c_i')$. By Dehn twisting $\gamma$ along the $\gamma_i$ we vary the pair $b_i', c_i'$ as we please, and by twisting at the correct rates we construct a sequence of simple closed curves that converge projectively to $[\mu]$.
\end{proof}

The mapping class group $\MCG(S_g)$ acts naturally on the whole set $\calC$ of currents and in particular it acts by homeomorphisms on $\partial \Teich(S_g)$. 

The identification of $\partial\Teich(S_g) = \matP\calML$ with the projectivisation of $C\times \ldots \times C$ is of course not canonical because it depends on a decomposition of $S_g$ into pairs of pants and annuli. However, Proposition \ref{integral:PL:map:prop} and Exercise \ref{abc:ex} imply that any two different identifications differ by some projective integral piecewise-linear homeomorphisms, hence the Thurston boundary has a natural \emph{projective integral piecewise linear} structure, called PIP for short, which is preserved by $\MCG(S_g)$. In particular $\matP\calML$ contains some natural rational points, and these are $\matP\calM$.\index{PIP structure}

\subsection{A projection}
The reader has probably noted that the coordinates that yield the homeomorphisms $\Teich(S_g) \cong \matR^{6g-6}$ and $\partial \Teich(S_g)\cong S^{6g-7}$ are quite similar: they both depend on a pants decomposition $\mu$ plus some additional marking, and every curve of the pants decomposition contributes roughly with two parameters, a ``length'' and a ``twist''. It would now be reasonable to expect that both these coordinates merge nicely to give a global homeomorphism $\overline{\Teich(S_g)} \cong D^{6g-6}$, but this is unfortunately not the case. 

The proof that $\overline{\Teich(S_g)}$ is homeomorphic to $D^{6g-6}$ is disappointingly indirect. We first prove that $\overline{\Teich(S_g)}$ is a topological manifold with boundary, by constructing some charts. To construct a chart, we first build a map
$$q\colon\Teich(S_g) \to \ML$$ 
that depends only on a fixed pants decomposition $\mu$. For a given hyperbolic metric $m\in\Teich(S_g)$, we construct a measured geodesic lamination $q(m)$ as follows. We straighten $\mu$ to its geodesic representative (with respect to $m$), and we consider separately each geodesic pair of pants $P$ of the decomposition. We first show that $P$ has a natural partial foliation, that depends only on its metric, whose leaves are not geodesics. We then glue these partial foliations to get a partial foliation on $S_g$, and then straighten it to a geodesic lamination $q(m)$.

\begin{figure}
\begin{center}
\includegraphics[width = 10 cm] {\iftoggle{BW}{pants_q-BW}{pants_q}}
\nota{Every geodesic pair of pants $P$ has a natural partial foliation: the shape of the foliation is determined by the sub-triangle containing $Q=[a_1,a_2,a_3]$ in Figure \ref{pants_train_tracks:fig}, where $a_1,a_2,a_3$ are the lengths of $\partial P$; the two pictures shown here correspond to two cases. In each case, there are three rectangles: each is the $R$-neighbourhood of the corresponding (unique) orthogeodesic\iftoggle{BW}{}{ (drawn in blue)}, and $R$ is the linear combination of $a_1,a_2,a_3$ prescribed by Figure \ref{pants_train_tracks:fig}, that is the unique combination that guarantees that the three rectangles match nicely as in the figure (1). We cut the geodesic pants $P$ along the three\iftoggle{BW}{}{ blue} orthogeodesics and get two right-angled hexagons. Since by hypothesis the three \iftoggle{BW}{small}{red} sides have length $\geqslant \varepsilon/2$, one sees easily that the central triangle has diameter bounded by some constant $C'$ depending only on $\varepsilon$, and hence every leaf in each of the three rectangles is shorter than some $C''$ that also depends only on $\varepsilon$ (2).}
\label{pants_q:fig}
\end{center}
\end{figure}

The natural partial foliation on $P$ is constructed as follows. Let $g_{ij}$ be the unique orthogeodesic connecting the $i$-th and $j$-th boundary component of $P$, for all $1 \leqslant i,j \leqslant 3$. Let $a_1, a_2, a_3$ be the lengths of the boundaries of $P$. We consider the point $Q=[a_1,a_2,a_3]$ in the triangle of Figure \ref{pants_train_tracks:fig}, pick the three orthogeodesics $g_{ij}$ isotopic to the three curves indicated in the sub-triangle of Figure \ref{pants_train_tracks:fig} containing $Q$, and thicken each $g_{ij}$ to a metric $R$-neighbourhood of $g_{ij}$, where $R$ is the linear combination of $a_1,a_2,a_3$ indicated in Figure \ref{pants_train_tracks:fig}. The $R$-neighbourhood is naturally foliated into (non-geodesic) arcs staying at fixed distance from $g_{ij}$ and the foliated neighbourhoods cover nicely much of the pair of pants $P$, as shown in Figure \ref{pants_q:fig}-(1). 
Every foliated rectangle is equipped with a natural transverse measure induced by the orthogonal distance between leaves.

The partial foliations of all the pairs of pants of $S_g$ glue to a singular measured foliation for $S_g$, which straightens as prescribed in Section \ref{parametrizzazione:subsection} to a measured geodesic lamination $q(m)$. 

Let $\mu = \gamma_1\sqcup \ldots \sqcup \gamma_{3g-3}$ be our original pants decomposition.

\begin{prop} \label{q:homeom:prop}
The map $q$ restricts to a homeomorphism
$$q\colon \Teich(S_g) \longrightarrow \big\{\lambda \in \ML \ \big|\ i(\lambda, \gamma_i)>0\ \forall i\ \big\}.$$
\end{prop}
\begin{proof}
Compare the coordinates $(l_i, \theta_i)$ for $\Teich(S_g)$ and $(a_i,b_i,c_i)$ for $\ML$. The map $q$ sends $(l_i,\theta_i)$ to $(l_i,b_i,c_i)$ for some $(b_i,c_i)$ that depends homeomorphically on $\theta_i\in\matR$ for each $i$.
Note that $a_i = i(\lambda, \gamma_i)$.
\end{proof}

\subsection{The fundamental lemma}
We now show that the projection $q$ distorts very little the lengths of the simple closed curves, as long as we put a lower bound on the lengths of the curves $\gamma_1,\ldots \gamma_{3g-3}$ of the fixed pants decomposition $\mu$. This technical fact is called the Fundamental Lemma \cite{FLP}.

For every $\varepsilon > 0$, we define $V(\varepsilon) \subset \Teich(S_g)$ to be the open subset consisting of all metrics $m$ such that $i(m,\gamma_i)>\varepsilon$ for all $i$.

\begin{lemma}
For every simple closed curve $\alpha\in \calS$ there exists a constant $C>0$ such that, for all $m\in V(\varepsilon)$, we have
$$i(q(m),\alpha) \leqslant i (m,\alpha) \leqslant i(q(m), \alpha) + C.$$
\end{lemma}
\begin{proof}
We represent $q(m)$ as a singular partial foliation, without straightening it. The transverse measure of $q(m)$ is just the length of orthogonal geodesics, hence the $q(m)$-measure of any piecewise-transverse closed curve is smaller or equal than its length: this proves that  $i(q(m),\alpha) \leqslant i (m,\alpha)$.

We prove the other inequality. If $\alpha = \gamma_i$ then $i(m,\alpha) = i(q(m),\alpha)$, so we suppose that $i(\alpha, \mu)>0$. Up to isotopy we may take $\alpha$ to be monotonically transverse to $q(m)$, recall the proof of Proposition \ref{integral:PL:map:prop}. The curve $\alpha$ intersects each pair of pants $P$ into some essential arcs $\beta$ that cross each foliated rectangle monotonically. 

We can easily homotope each of these arcs $\beta \subset \alpha\cap P$ with fixed endpoints into a (not necessarily simple!) piecewise smooth arc that decomposes into finitely many sub-arcs, that are alternatively contained either in a component of $\partial P$ or in a leaf of some foliated rectangle. Each sub-arc in a component of $\partial P$ may make many full turns (and hence may not be injective). Each of the three rectangles contains at most one leaf that is a sub-arc of $\beta$. We arrange the homotopy efficiently so that the new $\beta$ has the same transverse $q(m)$-measure as before.

The transverse $q(m)$-measure of $\beta$ on the sub-arcs in $\partial P$ is equal to its length, whereas on the leaf sub-arcs it is zero. Therefore the $m$-length of $\beta$ is equal to its measure plus the length of at most 3 leaves of some rectangles in $P$. Figure \ref{pants_q:fig}-(2) shows that each such leaf is shorter than some constant $C''$ that depends only on $\varepsilon$. 

Note that $\alpha$ is decomposed into $i(\alpha, \mu)$ arcs like $\beta$.
By homotoping each $\beta$ as above we find a homotopic representative for $\alpha$ whose length is at most the $q(m)$-measure of $\alpha$ plus $C=3i(\alpha,\mu)C''$. The length of any homotopic representative is greater or equal than the length $i(m,\alpha)$ of the geodesic one, and this proves the second inequality.
\end{proof}

\begin{cor} \label{same:limit:cor}
Let $m_i \in V(\varepsilon)$ be a diverging sequence in $\Teich(S_g)$. The sequence converges in $\overline{\Teich(S_g)}$ $\Longleftrightarrow$ the sequence $[q(m_i)]\in\PML$ does, and in this case they tend to the same limit.
\end{cor}
\begin{proof}
If $m_i$ converges to some $[\alpha]\in\PML$, we have $\lambda_im_i \to \alpha$ for some real numbers $\lambda_i \to 0$. By the fundamental lemma, for every simple closed curve $\gamma \in \calS$ we get
$$\big|i \big(\lambda_im_i, \gamma\big) - i\big(\lambda_iq(m_i), \gamma\big) \big| \to 0.$$
Therefore $\lambda_iq(m_i)$ converges (on a subsequence) to a $\beta\in\calML$ such that $i(\alpha, \gamma) = i(\beta, \gamma)$ for all $\gamma \in \calS$. Exercise \ref{abc:ex} gives $\beta =\alpha$.

The other implication is analogous.
\end{proof}

\subsection{A topological chart}
Pick an arbitrary $[\lambda]\in \PML = \partial \Teich(S_g)$. We now construct an explicit neighbourhood of $[\lambda]$ in $\overline{\Teich(S_g)}$. 

Let $\mu = \gamma_1\sqcup \ldots \sqcup \gamma_{3g-3}$ be a pants decomposition such that $i(\gamma_i, \lambda)>0$ for all $i$. We fix $\varepsilon>0$ and define as above $V(\varepsilon)\subset \Teich(S_g)$ as the set of all the metrics $m$ such that $i(\gamma_i, m)> \varepsilon$ for all $i$. We define similarly $W\subset \PML$ as the set of all $[\alpha]$ such that $i(\gamma_i, \alpha) > 0 $ for all $i$. We have $[\lambda]\in W$. 

\begin{prop} We have 
$$\pi(q(V(\varepsilon))) = \pi(q(\Teich(S_g))) = W.$$
The set $W\cup V(\varepsilon)$ is an open neighbourhood of $[\lambda]$ in $\overline{\Teich(S_g)}$.
\end{prop}
\begin{proof}
Proposition \ref{q:homeom:prop} implies the second equality, the first holds because every $[\alpha] \in W$ is represented by an $\alpha$ with $i(\gamma_i, \alpha)\geqslant \varepsilon$.

The sets $W$ and $V(\varepsilon)$ are open in $\PML$ and $\Teich(S_g)$ respectively.
If $W\cup V(\varepsilon)$ were not open in $\overline{\Teich(S_g)}$, there would be a sequence $m_i \not\in V(\varepsilon)$ of metrics converging to some $[\alpha]\in W$. On a subsequence we may suppose that $i(m_i,\gamma_j)<\varepsilon$ for some fixed $j$ and therefore $i(\alpha, \gamma_i)=0$, a contradiction.
\end{proof}

We have found an open neighbourhood of $[\lambda]$ and we now want to determine its topology. 

\begin{prop}
There is a homeomorphism 
$$\phi\colon W\cup V(\varepsilon) \longrightarrow U $$
onto an open subset $U\subset H$ of a half-space $H\subset \matR^{6g-6}$, with $\phi(W) = U\cap \partial H$.
\end{prop}
\begin{proof}
We complete the pants decomposition $\mu$ to a filling set $\mu'\in \calC$ of simple closed curves in $S_g$, considered as a current. We know from Proposition \ref{bounded:metrics:compact:cor} that the metrics $m$ with $i(m,\mu')\leqslant M$ form a compact subset in $\Teich(S_g)$. We define the map
\begin{align*}
\phi \colon W \cup V(\varepsilon) & \longrightarrow W\times [0,1] \\
x & \longmapsto 
\left\{
\begin{array}{ll} (x,0) & {\rm if\ } x\in W, \\
(\pi(q(x)), e^{-i(q(x),\mu')}) & {\rm if\ } x\in V(\varepsilon).
\end{array}
\right.
\end{align*}

The map $\phi$ is continuous: if a sequence $m_i\in V(\varepsilon)$ of metrics tends to $[\alpha]\in\PML$, then $i(m_i,\mu')\to \infty$, hence $i(q(m_i),\mu')\to \infty$ by the fundamental lemma. Moreover $\pi(q(m_i)) \to [\alpha]$ by Corollary \ref{same:limit:cor}. Summing up, we get that $\phi(m_i) \to ([\alpha],0) = \phi([\alpha])$.

The map $\phi$ is injective: given two metrics $m,m'$, if the first components of $\phi(m), \phi(m')$ are equal then $q(m) = \lambda q(m')$, and if the second components are equal we get $\lambda = 1$.

By similar methods one proves that the image of $\phi$ is open and the inverse there is continuous. The set $W\times [0,1]$ embeds in $H$ since $W$ embeds in $\matR^{6g-7}$.
\end{proof}

\begin{cor} The space $\overline{\Teich(S_g)}$ is a topological manifold with boundary $\partial \Teich(S_g)$.
\end{cor}

\subsection{Conclusion} \label{MCG:conclusion:subsection}
We have discovered that $\overline{\Teich(S_g)}$ is a compact topological manifold with boundary, and we now invoke Theorem \ref{topological:disc:teo} to deduce that it is a disc. We have finally completed the proof of the following theorem, which is the main achievement of the whole chapter.

\begin{teo} \label{current:compactification:teo}
The closure $\overline{\Teich(S_g)}$ of $\Teich(S_g)$ in $\matP\calC$ is homeomorphic to the closed disc $D^{6g-6}$. Its interior is $\Teich(S_g)$ and its boundary contains $\calS$ as a dense subset.
\end{teo}

Thurston's original compactification theorem embeds everything in $\matR^\calS$ instead of $\calC$. We can easily deduce it from Theorem \ref{current:compactification:teo}.

\begin{proof}[Proof of Theorem \ref{compattificazione:teo}]
The natural map $\calC \to \matR^\calS$ induced by the intersection form $i$ is continuous since $i$ is, and it induces a continuous map $\varphi\colon\matP\calC \to \matP(\matR^\calS)$ on their projective spaces. 

The Teichm\"uller space is embedded in both projective spaces and $\varphi$ restricts to a homeomorphism between the two embeddings, and to a continuous surjective map from their compactifications. The map on compactifications is actually injective since $\varphi|_{\PML}$ is (intersections with simple closed curves distinguish laminations) and is hence a homeomorphism.
\end{proof}

\section{Surface diffeomorphisms} \label{surface:diffeomorphisms:section}
We have discovered that Thurston's compactification $\overline{\Teich(S_g)}$ of the Teichm\"uller space is homeomorphic to a closed disc. The mapping class group $\MCG(S_g)$ acts on it naturally. We can now apply Brouwer's fixed point theorem to every element $\varphi \in \MCG(S_g)$, and characterise $\varphi$ according to the position of its fixed points, much similarly as we did for the isometries of $\matH^n$.

\subsection{The torus case}
As usual, the flat torus case is very instructive because everything can be written explicitly. 

We know from Proposition \ref{toro:azione:prop} that $\MCG(T)$ acts on $\Teich(T)$ like the M\"obius transformations $\PSLZ$ do on the hyperbolic half-plane $H^2$. The action of course extends to the compactification $\overline{H^2} = H^2\cup \matR\cup\{\infty\}$ and every non-trivial isometry $A \in\PSLZ$ is elliptic, parabolic, or hyperbolic according to the position of its fixed points.

As an integral matrix, the isometry $A$ has also more properties that are easy to check: if it is elliptic, it has finite order because $\PSLZ$ is discrete; if it is parabolic, it is conjugate in $\PSLZ$ to a matrix $\matr 1 n 0 1$ for some $n\neq 0$ and is thus the $n$-th power of a Dehn twist; if it is hyperbolic, it is conjugate in $\PSLR$ to a diagonal matrix $\matr \lambda 0 0 {\lambda^{-1}}$ for some $\lambda >1$, with some basis of eigenvectors $v,w$. The two foliations of $\matR^2$ into lines parallel to $v$ or $w$ are both preserved by $A$ and descend to foliations in $T = \matR^2/_{\matZ^2}$ that are preserved by $A$. The two foliations are both irrational (because $\lambda$ is), that is every leaf is dense in $T$, and $A$ stretches one foliation by $\lambda$ and contracts the other by $1/\lambda$. One such diffeomorphism of $T$ is called \emph{Anosov}.\index{Anosov diffeomorphism}

Summing up, the non-trivial elements in $\MCG(T)$ either have finite order, or preserve a simple closed curve, or are Anosov. We now define an analogous trichotomy for the elements in $\MCG(S_g)$ when $g\geqslant 2$, where the foliations are replaced by measured geodesic laminations, and the measure is there to encode stretchings and contractions.

\subsection{The trichotomy}
Let $S_g$ have genus $g\geqslant 2$. The mapping class group $\MCG(S_g)$ acts naturally on the whole space $\calC$ of currents and in particular on the compactification $\overline{\Teich(S_g)} \cong D^{6g-6}$ of the Teichm\"uller space.

Let $\varphi \in \MCG(S_g)$ be a non-trivial element. By Brouwer's fixed point theorem, $\varphi$ fixes at least one point in $\overline{\Teich(S_g)}$. We say that $\varphi$ is:
\begin{enumerate}
\item \emph{finite order} if it fixes a hyperbolic metric $m \in \Teich(S_g)$;
\item \emph{reducible} if it fixes a multicurve $\mu \in \calM$;
\item \emph{pseudo-Anosov} in all the other cases.
\end{enumerate}

We now analyse the three cases individually.\index{finite-order mapping class}\index{reducible mapping class}\index{pseudo-Anosov mapping class}

\subsection{Finite order elements}
We must first explain the terminology.

\begin{prop} \label{finite:finite:prop}
A non-trivial element $\varphi\in\MCG(S_g)$ is finite order if and only if it has indeed finite order in $\MCG(S_g)$.
\end{prop}
\begin{proof}
Suppose that $\varphi$ preserves the isotopy class $[m]\in\Teich(S_g)$ of a hyperbolic metric $m$ in $S_g$. We can choose a representative for $\varphi$ that fixes $m$. This representative is an isometry for $m$.
Since the isometry group of a closed hyperbolic manifold is finite (see Corollary \ref{finite:isometry:group:cor}) we have $\varphi^n=\id$ for some $n>1$ and $\varphi$ has indeed finite order in $\MCG(S_g)$. 

Conversely, let $\varphi$ be an element having finite order in $\MCG(S_g)$. The subgroup $\langle \varphi \rangle$ generated by $\varphi$ cannot act freely on $\Teich(S_g) \isom  \matR^{6g-6}$, otherwise it would quotient $\matR^{6g-6}$ to an aspherical manifold with finite fundamental group, contradicting Theorem \ref{no:torsion:aspherical:teo}. Therefore some non-trivial power of $\varphi$ has a fixed point in $\Teich(S_g)$.

If $\varphi$ has prime order we easily conclude that also $\varphi$ has a fixed point and we are done. However, if $\varphi$ has order $p_1\cdots p_s$ for some primes $p_i$, we need to do more work. By induction, $\varphi' = \varphi^{p_1\cdots p_{s-1}}$ has a fixed point $[m]\in\Teich(S_g)$ and is hence represented by an isometry for $S_g$ with metric $m$. The isometry $\varphi'$ quotients $S_g$ to a hyperbolic orbifold, and the fixed points $\Fix(\varphi')$ of $\varphi'$ in $\Teich(S_g)$ can be identified naturally with the (suitably defined) Teichm\"uller space of this orbifold, which is homeomorphic to $\matR^N$ for some $N>0$ like in the surface case (exercise).

Since $\varphi$ and $\varphi'$ commute, the first act as a mapping class on $\Fix(\varphi')$ with order $p_1\cdots p_{s-1}$: we conclude by induction on $s$ that $\varphi$ has a fixed point in $\Teich(S_g)$.
\end{proof}

We get in particular the following corollary, which is far from obvious because $\MCG(S_g)$ does \emph{not} act itself on $S_g$, and algebraic relations in $\MCG(S_g)$ do not translate into algebraic relations between representatives, except in some very lucky cases.

\begin{cor}
If $\varphi\in\MCG(S_g)$ has order $k$, it may be represented by a diffeomorphism $\varphi\colon S_g\to S_g$ such that $\varphi^k = \id$.
\end{cor}
\begin{proof}
The class $\varphi$ has a representative $\varphi\colon S_g \to S_g$ that is an isometry for some hyperbolic metric; the isometry $\varphi^k$ is isotopic to the identity and is hence the identity by Corollary \ref{isometries:not:homotopic:cor}.
\end{proof}

\subsection{Reducible elements}
We must explain the terminology also in this case.
If $\varphi$ fixes a multicurve $\mu$, one can cut $S_g$ along $\mu$ and look at the restriction of $\varphi$ to the resulting pieces: after extending all the theory to surfaces with boundary (that we have not done here for simplicity), we can hence study inductively each piece, and this explains the word \emph{reducible}.

The cases (1) and (2) are not exclusive: there are isometries of hyperbolic surfaces that preserve some multicurve. On the other hand, there are finite order elements that are not reducible (exercise) and reducible mapping classes that are not of finite order: for instance, Dehn twists. 

\subsection{The action on $\calML$}
The mapping class group $\MCG(S_g)$ of $S_g$ acts on the currents space $\calC$ and hence on the ``light cone'' subspace $\calML$ of all measured geodesic laminations, which contains the space $\calM$ of multicurves. 

The action of $\MCG(S_g)$ on $\calML$ can be seen concretely inside $S_g$. It suffices to consider $\calML$ as the space of all measured \emph{topological} laminations (defined in Section \ref{topological:laminations:subsection}) considered up to isotopy and collapsing of parallel closed leaves. 
Now $\varphi\in\MCG(S_g)$ acts on $\calML$ simply by sending the measured topological lamination $\mu$ to $\varphi(\mu)$. 

We now prove that if $\varphi$ fixes a non-trivial point in $\calML$ we fall back into one of the two cases already considered.

\begin{prop} \label{varphi:fixes:mu:prop}
If $\varphi(\mu) = \mu$ for some non-trivial $\mu\in\ML$, then $\varphi$ is either finite order or reducible.
\end{prop}
\begin{proof}
We fix a hyperbolic metric on $S_g$ and represent $\mu$ as a measured geodesic lamination there. 

By Proposition \ref{finite:regions:prop} there are finitely many complementary regions in $S_g\setminus \mu$. If $\mu$ is not full, some region is not an ideal polygon and hence deformation-retracts onto a subsurface $S'\subset S_g$, whose boundary components $\partial S'$ are non-trivial simple closed curves in $S_g$. The union of all such curves $\partial S'$ as $S'$ varies produces a multicurve preserved by $\varphi$, which is hence reducible.

If $\mu$ is full, all the complementary regions are ideal polygons, and after substituting $\varphi$ with a finite power we can suppose that each region is fixed (not pointwise) by $\varphi$, together with each of its boundary components. We consider the preimage $\tilde \mu \subset \matH^2$ and note that its complementary regions in $\matH^2$ are still ideal finite polygons, the lifts of the ones in $S_g$. 

The diffeomorphism $\varphi$ lifts to a homeomorphism $\tilde\varphi \colon \overline{\matH^2} \to \overline{\matH^2}$ that fixes $\tilde \mu$, and we may suppose that it fixes (not pointwise) a complementary polygon and each of its boundary lines. In particular $\varphi$ fixes orientation-preservingly a line $l\subset \tilde \mu$. 

Remember that $\calG$ is the set of all lines in $\matH^2$.
We now prove that, if $\tilde\varphi$ fixes (orientation-preservingly) a 
line $l\subset \tilde \mu \subset \calG$, then it fixes (still orientation-preservingly) all the lines in $\tilde \mu$ contained in some neighbourhood of $l\in \calG$. This is done as follows: since the boxes form a neighbourhood system for $\calG$, and $\tilde\varphi$ acts via homeomorphisms on $\calG$, there are two boxes $B\subset B'$ containing $l$ such that $\tilde\varphi(B) \subset B'$. Now $B\cap \tilde \mu \subset B'\cap \tilde \mu$ consists of some parallel lines, which look like $J\times \matR \subset J' \times \matR$ for some measured ordered set $J'$ and some subsegment $J\subset J'$: the map $\varphi$ sends  $J$ to $J'$ preserving both the ordering and the measure, and fixing the point corresponding to $l$: hence it is the identity.

On the other hand, every complementary region of $\tilde\mu$ is a finite ideal polygon: hence if $\tilde\varphi$ fixes (orientation-preservingly) one side of the polygon, it fixes all the others. 

The two properties just listed together easily imply that, since $l$ is fixed, the whole of $\tilde\mu$ is. Since the lines are fixed orientation-preservingly, their endpoints are fixed: the map $\tilde\varphi$ fixes all the endpoints of all lines in $\tilde\mu$, and since these endpoints form a dense subset of $\partial \matH^2$ the map $\tilde\varphi$ fixes $\partial \matH^2$ pointwise and hence $\varphi$ is trivial in $\MCG(S_g)$ by Proposition \ref{action:curves:prop}.
\end{proof}

The previous proposition shows in particular that pseudo-Anosov elements act freely on $\calML\setminus \{0\}$. We now investigate more closely these mysterious mapping classes.

\subsection{Pseudo-Anosov elements}
A pseudo-Anosov element $\varphi$ is by definition neither finite order nor reducible. We have just seen that $\varphi$ acts freely on $\calML\setminus \{0\}$, but this does not prevent it from having some fixed points in $\PML$; indeed we now show that there are two fixed points there, one attracting and the other repelling, so that $\varphi$ looks very much like a hyperbolic isometry on the hyperbolic space.

\begin{figure}
\begin{center}
\includegraphics[width = 5 cm] {\iftoggle{BW}{pAdynamics-BW}{pAdynamics}}
\nota{The appropriate lift $\tilde \varphi$ acts on $\partial \matH^2$ with $2k$ fixed points that are alternatively attractive and repelling. By joining the repelling points we find another lamination $\mu_u$ fixed by $\varphi$. Here $k=5$.}
\label{pAdynamics:fig}
\end{center}
\end{figure}

\begin{teo} \label{su:teo}
Let $\varphi \in \MCG(S_g)$ be a pseudo-Anosov element. There are two measured geodesic laminations $\mu_s, \mu_u \in \ML$ and a real number $\lambda >1$ such that 
$$\varphi(\mu_s) = \lambda \mu_s, \qquad \varphi(\mu_u) = \frac 1\lambda \mu_u.$$
The laminations $\mu_s$ and $\mu_u$ are full, and they altogether fill $S_g$. 
\end{teo}
\begin{proof}
By Brouwer's fixed point theorem, a pseudo-Anosov element $\varphi$ has a fixed point in $\overline{\Teich(S_g)}\isom D^{6g-6}$ which is (by definition) neither a metric nor a multicurve. Therefore $\varphi$ fixes a projective measured lamination $[\mu]$ which is full (otherwise $\varphi$ would be reducible: see the proof of Proposition \ref{varphi:fixes:mu:prop}).

Since $[\mu]$ is a projective class, we have $\varphi(\mu) = \lambda\mu$ for some real number $\lambda >0$. Up to replacing $\varphi$ with its inverse $\varphi^{-1}$ we may suppose that $\lambda \geqslant 1$, and Proposition \ref{varphi:fixes:mu:prop} shows that $\lambda>1$. We denote this $\mu$ by $\mu_s$.

We now construct $\mu_u$. As in the proof of Proposition \ref{varphi:fixes:mu:prop}, we consider the preimage $\tilde\mu_s\subset \matH^2$ of $\mu_s$, and after replacing $\varphi$ with a finite power we may choose a lift $\tilde \varphi$ of $\varphi$ that fixes a complementary polygonal region $R$ of $\tilde\mu_s$ and its sides, hence in particular the vertices of $R$, see Figure \ref{pAdynamics:fig}.

The $k$ vertices of $R$ divide $\partial \matH^2$ into $\tilde\varphi$-invariant arcs $I_1,\ldots, I_k$, corresponding to the sides $s_1,\ldots, s_k$ of $R$, see Figure \ref{pAdynamics:fig}. Since $\mu_s$ is full, each $s_i$ is the limit of a sequence of leaves in $\tilde\mu_s$ with both endpoints in $I_i$ converging to the endpoints of $s_i$ but distinct from them by Proposition \ref{at:most:two:prop}. Since $\lambda > 1$, the map $\tilde\varphi$ pushes these leaves towards $s_i$, so in particular the vertices of $R$ are local attractors for the action of $\tilde\varphi$ on $\partial \matH^2$, see Figure \ref{pAdynamics:fig}.

Since the endpoints $p$ and $q$ of $l_i$ are attractors, the map $\tilde\varphi$ fixes at least one point in the interior of $I_i$, and we show that it cannot fix two: if $\tilde \varphi$ fixes $r$ and $s$ in the interior of $I_i$, then it fixes the box $[p,r] \times [s,q]$ that has non-zero measure, a contradiction since $\lambda > 1$. There is a single fixed point in $l_i$, and it must be repulsive.

The dynamics of $\tilde \varphi$ on $\partial \matH^2$ is described in Figure \ref{pAdynamics:fig}. 
The closure of the projection of the $k$ \iftoggle{BW}{light grey}{green} lines is another invariant geodesic lamination $\mu_u$, which must be full 
because $\varphi$ is pseudo-Anosov.

We now prove that $\mu_u$ admits some (non unique) transverse measure and can hence be considered as an element of $\ML$. The dynamics shows that for some curve $\gamma\in\calS$ the supports of $\varphi^{-k}(\gamma)$ tend to that of $\mu_u$. The sequence $\varphi^{-k}([\gamma])\in \PML$ hence converges on a subsequence to a projective measured geodesic lamination with support $\mu_u$.

In principle, the measure supported by $\mu_u$ needs not to be unique, not even up to rescaling: hence the element $[\mu_u]$ may not be uniquely determined and we cannot conclude that it is fixed by $\varphi$, unfortunately. However, distinct measures on the same support form obviously a convex cone in the current space $\calC$ and hence a closed disc in $\matP\calC$, the class $\varphi$ acts on this disc and therefore has a fixed point there by Brouwer's fixed point theorem again. 

We have $\varphi(\mu_u) = \lambda' \mu_u$ for some $\lambda'$. 
It is clear that $\mu_u$ and $\mu_s$ fill $S_g$ altogether, and in particular $i(\mu_u, \mu_s) > 0$. Therefore
$$0 < i(\mu_u, \mu_s) = i(\varphi(\mu_u), \varphi(\mu_s)) = \lambda \lambda' i(\mu_u, \mu_s)$$
gives $\lambda' = \frac 1 \lambda$. The proof is complete.
\end{proof}

The laminations $\mu_s$ and $\mu_u$ are the \emph{stable} and the \emph{unstable} measured geodesic laminations fixed by $\varphi$, and $\lambda$ is the \emph{dilatation} of $\varphi$. We now prove a converse to Theorem \ref{su:teo}.\index{geodesic lamination!stable and unstable measured geodesic lamination}

\begin{prop} \label{pA:converse:prop}
If $\varphi \in \MCG(S_g)$ is such that $\varphi(\mu) = \lambda \mu$ for some full 
$\mu \in \calML$ and $\lambda > 1$, then $\varphi$ is pseudo-Anosov.
\end{prop}
\begin{proof}
We need to prove that $\varphi$ is neither finite order nor reducible, that is that $\varphi$ fixes no non-trivial current $\alpha \in \Teich(S_g) \cup \calM$. We have 
$$i\big(\varphi^k(\alpha), \mu\big) = i\big(\alpha, \varphi^{-k}(\mu)\big) =
\lambda^{-k}i\big(\alpha, \mu) \to 0$$
as $k\to +\infty$. Therefore $\varphi(\alpha) \neq \alpha$ unless $i(\alpha, \mu)=0$, which is excluded since $\mu$ is full.
\end{proof}

\subsection{Examples}
We now construct plenty of pseudo-Anosov diffeomorphisms. We need a bit of simple linear algebra.

We say that a matrix or vector is \emph{positive} if all its entries are.
A square matrix $M$ with non-negative integral entries is \emph{Perron-Frobenius} if $M^k$ is positive for some $k>0$.\index{Perron-Frobenius matrix} 

\begin{prop} \label{PF:prop}
Every Perron-Frobenius integral matrix $M$ has a positive eigenvector $v$ with some eigenvalue $\lambda>1$.
\end{prop}
\begin{proof}
The matrix $M$ has non-negative entries and hence preserves the standard simplex $\Delta = \{x_i\geqslant 0, \sum x_i = 1\}$, so by Brouwer's fixed point theorem it has a fixed point $v\in \Delta$ there. Since $M^k > 0$ and $v$ is an eigenvector for $M^k$, we get $v>0$. Since $M^k$ is integral we get $\lambda^k > 1$ and hence $\lambda >1$.
\end{proof}

The following examples were constructed by Penner in 1988. Recall that a multicurve $\alpha$ in $S_g$ is \emph{essential} if it contains no parallel components: this holds for instance if $\alpha$ is a simple closed curve or a pants decomposition.

\begin{figure}
\begin{center}
\includegraphics[width = 11 cm] {\iftoggle{BW}{Penner-BW}{Penner}}
\nota{We smoothen $\alpha \cup \beta$ to a bigon track that contains both $\alpha$ and $\beta$ (left). The bigon track may have some bigons as complementary regions (right).}
\label{Penner:fig}
\end{center}
\end{figure}

\begin{teo} \label{Penner:teo}
Let $\alpha$ and $\beta$ be two essential multicurves that altogether fill $S_g$. Let $\varphi\in\MCG(S_g)$ be any composition of Dehn twists $T_a^{+1}$ and $T_b^{-1}$ where $a$ and $b$ vary among the curves in $\alpha$ and $\beta$. If every component $a$, $b$ of $\alpha$, $\beta$ occurs at least once, then $\varphi$ is pseudo-Anosov.
\end{teo}
The Dehn twists may occur in any order, for instance $\varphi = T_{a}^2T_b^{-3}T_{b'}^{-1}T_{a'}$.
\begin{proof}
We put $\alpha$ and $\beta$ in minimal position and smoothen the transverse intersections as in Figure \ref{Penner:fig}-(left) to get a \emph{bigon track} $\tau$. A bigon track is like a train track, except that it may contain some complementary bigon as in Figure \ref{Penner:fig}-(right). The straightening procedure described in Section \ref{parametrizzazione:subsection} works also in this case (exercise), so every weight system on $\tau$ determines a measured geodesic lamination in $S_g$.\index{train track!bigon track} (Different weight systems may determine the same geodesic laminations because of the bigons.)

We have $\alpha = a_1 \sqcup \ldots \sqcup a_m$ and $\beta = b_1 \sqcup \ldots \sqcup b_n$. Each $a_i$ and $b_j$ may be represented by assigning weights 1 or 0 to the edges of $\tau$, and these weights form $m+n$ independent (exercise) vectors $v_1,\ldots, v_{m+n}$
in the weights space of $\tau$. Let $V$ be the $(m+n)$-dimensional sub-cone of the weights space generated by $v_1,\ldots, v_{m+n}$ via combinations with non-negative coefficients. 
Every vector in $V$ models a measured geodesic lamination in $S_g$.

Let $\Omega$ be the $m+n$ square matrix $\Omega_{ij} = i(v_i, v_j)$. We have $\Omega = \matr 0HK0$. The crucial point here is that both $T_{a_j}$ and $T_{b_k}^{-1}$ act on $V$ like the matrix 
$$Q_i = I + D_i\Omega$$
where $I$ is the identity and $D_i$ has 1 on the $i$-th entry of the diagonal and 0 everywhere else. The map $\varphi$ therefore acts on $V$ as a product $M = Q_{i_1} \cdots Q_{i_h}$ of non-negative matrices. Since each $i=1,\ldots,m+n$ occurs at least once as an index, we may deduce (exercise) that $M$ is Perron-Frobenius. Therefore $M$ has a positive eigenvector $v$ with eigenvalue $\lambda > 1$.

The positive eigenvector $v$ determines an element $\mu \in \ML$ such that $\varphi(\mu) = \lambda \mu$. We leave to the reader the proof that $\mu$ is full, 
and we conclude using Proposition \ref{pA:converse:prop}.
\end{proof}

\subsection{References}
The material contained in this chapter is well known to experts, but it is hard to find a source in the literature that contains everything in a fully self-contained way. The whole theory was presented by Thurston in a very nice and readable paper \cite{Th_diff} that however contained no proof. The most complete book on the subject is then Fathi -- Laudenbach -- Po\'enaru \cite{FLP}, and this is also the main source that we have used for writing this chapter. 
Another important source is Casson -- Bleiler \cite{CB}.

We have chosen to describe the whole theory using Bonahon's geodesic currents, that originated in the papers \cite{B1, B2}. To this purpose we have also consulted McMullen \cite{McM}, Aramayona -- Leininger \cite{AL}, and Calegari \cite{Ca}. We also borrowed some arguments on train tracks and measured geodesic laminations from a nice and self-contained paper of Hatcher \cite{H2}. The proof of Proposition \ref{finite:finite:prop} was taken from \cite[Theorem 6.1]{FM}. Theorem \ref{Penner:teo} was proved by Penner in \cite{Pe}.

%% file: Three_topology.tex
\chapter{Topology of three-manifolds} \label{Three:topology:chapter}

The three-manifolds world is topologically much richer than the surfaces realm, while yet not so crazy as the four-manifolds universe. In dimension two a simple homological invariant (the Euler characteristic) suffices to classify topologically all the closed orientable manifolds: at the complete opposite, the closed four-manifolds cannot be classified in any reasonable sense. The three-manifolds lie in the middle: we do not have yet a complete satisfactory picture, but we understand them a good deal.

The rest of this book is devoted to three-manifolds, more specifically to compact orientable three-manifolds, possibly with boundary (the orientability assumption is not essential, but it helps to simplify many arguments). We split their study into some parts. First, we state some universal topological facts, mostly concerning the way surfaces can be contained in three-manifolds. Then we construct some classes of three-manifolds, focusing mostly on Seifert manifolds. After classifying the Seifert manifolds topologically, we assign a geometry to each: there are eight interesting geometries in dimension three, and we introduce them with some care. Finally, we concentrate on the most interesting and beautiful of the eight: hyperbolic geometry.

In this chapter we start to study the topology of three-manifolds. We begin with some algebraic topology, then we show that the connected sum behaves like products of numbers: every closed three-manifold splits in a unique way in a unique list of prime factors. Finally we introduce the important notion of incompressible surface.\index{three-manifold}

\section{Algebraic topology}
The algebraic topology of compact three-manifolds is not very complicate. 

\subsection{Integral homology}
In this section all the homology groups are considered over $\matZ$. We first note that, for closed orientable 3-manifolds, the fundamental group determines everything.

\begin{prop}
The homology $H_*(M)$ of a closed orientable 3-manifold $M$ is determined by $\pi_1(M)$. 
\end{prop} 
\begin{proof}
As in every path connected space, the group $H_1(M)$ is the abelianisation of $\pi_1(M)$ and $H^1(M) = \Hom(H_1(M),\matZ)$, which is hence isomorphic to $H_1(M)$ modulo its torsion. By Poincar\'e duality $H_2(M) = H^1(M)$ and $H^2(M) = H_1(M)$. Finally $H_3(M) = \matZ$.
\end{proof}

\subsection{Homology spheres}
A \emph{homology $n$-sphere} is a closed orientable $n$-manifold $M$ whose homology is the same as that of $S^n$. That is, we have $H_0(M) = H_n(M) = \matZ$ and $H_i(M)$ vanishes for all $0<i<n$.\index{homology sphere} 

By Poincar\'e duality, a closed 3-manifold $M$ is a homology sphere if and only if $H_1(M)$ vanishes. Since $H_1(M)$ is the abelianization of $\pi_1(M)$, this happens precisely when $\pi_1(M)$ is a \emph{perfect group}, that is a group with trivial abelianization.\index{group!perfect group} 

In 1900 Poincar\'e conjectured that every homology 3-sphere should be homeomorphic to $S^3$. Four years later, he found himself a counterexample by constructing what is known today as \emph{Poincar\'e's homology sphere}, a closed 3-manifold with a perfect fundamental group of order 120. He then modified his original conjecture by asking whether every simply connected closed 3-manifold should be homeomorphic to $S^3$. This fact, widely known as \emph{Poincar\'e's Conjecture}, was proved only in 2002 by Perelman.

We will construct Poincar\'e's homology sphere (and many more homology spheres) in Chapter \ref{Seifert:chapter}, and we will discuss the Poincar\'e Conjecture in Section \ref{geometrisation:section}, as a part of Thurston's wider \emph{Geometrisation Conjecture}, also proved by Perelman in 2002.

\subsection{The boundary}
The Euler characteristic $\chi(M)$ of a closed odd-dimensional manifold vanishes, and on manifolds with boundary we have the following.

\begin{prop}
If $M$ is a compact 3-manifold with boundary, then 
$$\chi(M) = \frac {\chi(\partial M)}2.$$
\end{prop}
\begin{proof}
We only use that $M$ has odd dimension $n$. If $M$ is closed and orientable, we have $\chi(M) = \sum_{i=0}^n (-1)^i b_i$ and the Betti numbers $b_i$ and $b_{n-i}$ are equal by Poincar\'e duality, hence $\chi(M)=0$. If $M$ is non-orientable then it has an orientable double-cover $N$ and $\chi(N) = d\chi(M)$ on degree-$d$ covers, hence $\chi(N)=0$ implies $\chi(M)=0$.
If $M$ has boundary then 
$$0 = \chi(DM)=2\chi (M) - \chi(\partial M)$$
where $DM$ is the \emph{double} of $M$, constructed by taking two identical copies of $M$ and identifying their boundaries in the obvious way.
\end{proof}

The manifold $M$ has half the Euler characteristic of $\partial M$, and is also responsible for half of the first homology group of $\partial M$. More precisely, we are interested in the boundary map
$$\partial\colon H_2(M,\partial M) \longrightarrow H_1(\partial M)$$
and we want to prove that its image is a particular half-dimensional subgroup.
Recall from Proposition \ref{homology:prop} that every class in $H_2(M,\partial M)$ is represented by an oriented properly embedded surface $\Sigma$, and $\partial$ sends $[\Sigma]$ to $[\partial \Sigma]$.

Let $M$ be oriented. The boundary $\partial M$ may be disconnected and inherits an orientation. Recall that $H_1(\partial M)\isom \matZ^{2n}$ for some $n$ and $H_1(\partial M)$ is equipped with a symplectic intersection form $\omega$, see Section \ref{homolgy:surface:subsection}. A subgroup $L<H_1(\partial M)$ is \emph{lagrangian} if $\omega|_L\equiv 0$.\index{lagrangian subgroup} 

\begin{ex}
If $L$ is lagrangian then $\rk L \leqslant n$.
\end{ex}
When $\rk L = n$ we say that $L$ has maximal rank.

\begin{prop}  \label{lagrangian:prop}
Let $M$ be an oriented compact 3-manifold with boundary. The image of the map 
$$\partial \colon H_2(M,\partial M, \matZ) \longrightarrow H_1(\partial M,\matZ)$$
is a lagrangian subgroup of $H_1(\partial M, \matZ)$ of maximal rank.
\end{prop}
\begin{proof}
Consider the long exact sequence
$$\ldots \longrightarrow H_2(M,\partial M) \stackrel{\partial}\longrightarrow H_1(\partial M) \stackrel{i_*}\longrightarrow H_1(M) \longrightarrow \ldots$$
We have two pairings
\begin{align*}
\omega \colon H_1(\partial M) \times H_1(\partial M) & \longrightarrow  \matZ, \\
\eta \colon H_2(M,\partial M) \times H_1(M) & \longrightarrow \matZ.
\end{align*}
The latter is provided by Lefschetz duality and is non-degenerate after quotienting the torsion subgroups. We have
$$\omega (\partial \alpha, \beta) = \eta (\alpha, i_*\beta)$$
for all $\alpha \in H_2(M,\partial M)$ and $\beta \in H_1(M)$: this equality can be proved easily by representing the elements of $H_2(M,\partial M)$ as surfaces (which we can do thanks to Proposition \ref{homology:prop}). Now if $\beta = \partial\alpha'$ we get 
$$\omega(\partial \alpha, \partial \alpha') = \eta(\alpha, i_*\partial \alpha') = \eta(\alpha,0) = 0$$
and hence $\Img \partial$ is lagrangian. It has maximal rank since $H_2(M,\partial M)=H^1(M)$ and $H_1(M)$ have the same rank: if $\rk\, \Img \partial < \frac 12 b_1(\partial M)$, then we get $\rk \,\Img i_* > \frac 12 b_1(\partial M)$ and $\rk \ker\partial > b_1(M) - \frac 12 b_1(\partial M)$, a contradiction since $\Img i_*$ and $\ker\partial$ are $\eta$-orthogonal.
\end{proof}

\begin{cor} \label{rank:cor}
Let $M$ be an oriented compact 3-manifold. We have 
$$b_1(M)\geqslant \frac {b_1(\partial M)}2.$$
\end{cor}
\begin{proof}
The rank of $H_1(M)$ equals that of $H^1(M) = H_2(M,\partial M)$.
\end{proof}

\begin{cor}
The boundary of a simply connected compact 3-manifold consists of spheres.
\end{cor}

\subsection{Non-orientable surfaces}
Let $M$ be an orientable 3-manifold and $S\subset M$ be a connected surface. If $S$ is orientable, then by Proposition \ref{unique:line:bundle:prop} every tubular neighbourhood for $S$ is diffeomorphic to the product $S \times \matR$.

The orientable manifold $M$ may also contain a \emph{non}-orientable surface $S$: for instance the orientable projective space $\matRP^3$ contains the non-orientable projective plane $\matRP^2$. In that case a tubular neighbourhood of $S$ is diffeomorphic to the unique orientable interval bundle $S\timtil \matR$ over $S$, see Proposition \ref{unique:line:bundle:prop} again. A compact tubular neighbourhood of $S$ is an interval bundle $S\timtil I$, whose boundary is the orientable double cover of $S$. For instance, the boundary of $\matRP^2 \timtil I$ is a sphere.

A non-orientable properly embedded surface $S\subset M$ does not define a homology class in $\matZ$, but it defines one in $\matZ_2$, that is  $[S] \in H_2(M,\partial M;\matZ_2)$. As opposite to orientable surfaces, this class $[S]$ is always non-trivial.

\begin{prop} \label{non:ori:prop}
Let $M$ be orientable. Every non-orientable surface $S$ determines a non-trivial class $[S] \in H_2(M,\partial M;\matZ_2)$. The manifold $M$ cannot contain more than $\dim H_2(M,\partial M;\matZ_2)$ disjoint non-orientable surfaces.
\end{prop}
\begin{proof}
A tubular neighbourhood $S\timtil I$ of $S$ contains a simple closed loop $\alpha$ intersecting $S$ transversely in one point. The class $[S] \in  H_2(M,\partial M;\matZ_2) = H^1(M;\matZ_2) $ sends $\alpha$ to $1\in\matZ_2$ and is hence non-trivial. If $S=S_1\sqcup\ldots \sqcup S_k$ are all non-orientable, each $S_i$ has its own $\alpha_i$ and therefore the elements $[S_1], \ldots, [S_k] \in H_2(M,\partial M;\matZ_2)$ are linearly independent.
\end{proof}

\begin{cor} A simply-connected three-manifold does not contain any closed non-orientable surface.
\end{cor}

\section{Prime decomposition}
In this section we study spheres and discs in three-manifolds. We prove that every sphere in $\matR^3$ bounds a disc, and we call \emph{irreducible} a three-manifold that has this property. Then we study the connected sum of three-manifolds, and we show that it behaves like multiplication of natural numbers: every closed oriented three-manifold decomposes uniquely into some \emph{prime} factors.

\subsection{Balls}
In general dimension $n$, we have used the term \emph{disc} to denote the closed Euclidean disc $D^n$, and the term \emph{ball} for its open interior. Speaking about three-dimensional spaces, we will henceforth use a more intuitive terminology and call \emph{disc} and \emph{ball} respectively the closed discs $D^2$ and $D^3$, and we will use the symbols $D$ and $B$ for them.

Let $M$ be a connected 3-manifold. Theorem \ref{dischi:teo} says that all the closed balls $B\subset \interior M$ in $M$ are isotopic and hence the removal of the interior of $B$ from $M$ produces a new manifold $N$ with boundary, which does not depend on $B$. The sphere $\partial B$ is a new boundary component of $N$.

The inverse operation consists of \emph{capping off} a boundary component of $N$ diffeomorphic to $S^2$ by attaching a ball $B$ to it, via some diffeomorphism. This operation depends only on the boundary component that is capped off.

\begin{prop} \label{cap:prop}
The manifold $M$ obtained by capping off a spherical boundary component of $N$ does not depend on the diffeomorphism chosen.
\end{prop}
\begin{proof}
There are only two diffeomorphisms up to isotopy, see Theorem \ref{S2:omotopi:teo}, and they are related by a reflection of $B$, so they produce diffeomorphic manifolds.
\end{proof}

In dimension 3 we can therefore freely remove and add balls without affecting much the topology of the manifold. In particular, by removing the interior of a ball from $S^3$ we get another ball $B$, and by attaching a ball to $B$ we get $S^3$ back. We must thank Smale's Theorem \ref{Smale:teo} for that: the situation in higher dimensions is more complicated.

\subsection{Connected sums}
Connected sums exist in any dimension $n$, but when $n\leqslant 3$ they may be redefined in slightly simpler terms:\index{connected sum} 
\begin{defn}
The \emph{connected sum} $M=M_1\#M_2$ of two oriented connected $3$-manifolds $M_1$, $M_2$ is constructed by removing the interiors of two closed balls from $M_1$ and $M_2$, and then gluing the two resulting spheres via any orientation-reversing diffeomorphism. 
\end{defn}
In dimension $n\leqslant 3$ these diffeomorphisms are all isotopic, see Theorem \ref{S2:omotopi:teo},
 hence this is a good definition (in arbitrary dimension we require the diffeomorphism to extend to the removed discs). We usually work with oriented 3-manifolds to have a uniquely defined connected sum, but in some cases this is not necessary. An orientable manifold $M$ is \emph{mirrorable} if it admits an orientation-reversing self-diffeomorphism.\index{manifold!mirrorable manifold}

\begin{ex} If $M_1$ is mirrorable, the manifold $M = M_1\# M_2$ does not depend (up to diffeomorphisms) on the orientations of $M_1$ and $M_2$.
\end{ex}

If both $M_1$ and $M_2$ are oriented and not mirrorable, it may happen that $M_1\# M_2$ and $M_1 \# \overline{M_2}$ are not diffeomorphic (here $\overline{M_2}$ is $M_2$ with the reverse orientation).
There is also a boundary-version of connected sum:

\begin{defn}
The \emph{$\partial$-connected sum} $M=M_1\#_\partial M_2$ of two oriented $3$-manifolds with boundary is constructed by gluing two discs $D_1\subset \partial M_1$ and $D_2\subset \partial M_2$ via an orientation-reversing diffeomorphism.
\end{defn}

The operation depends only on the components of $\partial M_1, \partial M_2$ containing $D_1, D_2$, because of Theorem \ref{dischi:teo}. The following holds:
$$M\# S^3 = M, \qquad M\#_\partial B = M, \qquad \partial (M_1\#_\partial M_2) = (\partial M_1)\#(\partial M_2).$$
In the latter equality we suppose that $\partial M_1$ and $\partial M_2$ are connected.

\begin{ex} Let $M = M_1\# M_2$ or $M=M_1\#_\partial M_2$. We have
$$\pi_1(M) = \pi_1(M_1) * \pi_2(M_2).$$
\end{ex}
\begin{proof}[Hint]
Use Van Kampen. 
\end{proof}
\begin{cor}
We have $H_1(M,\matZ) = H_1(M_1,\matZ) \oplus H_1(M_2,\matZ)$.
\end{cor}

\subsection{Irreducible 3-manifolds}
We introduce an important definition. Let $M$ be a connected, oriented 3-manifold with (possibly empty) boundary.\index{three-manifold!irreducible three-manifold}
\begin{defn}
The manifold $M$ is \emph{irreducible} if every sphere $S\subset \interior M$ bounds a ball.
\end{defn}
It is important to recall here that the sphere must be a smooth surface. If we admitted also topological spheres, no 3-manifold would be irreducible: there are ``wild'' topological spheres inside every ball, as shown in Figure \ref{Alexander:fig}. 

\begin{figure}
\begin{center}
\includegraphics[width = 7 cm] {\iftoggle{BW}{Morse-BW}{Morse}}
\nota{Non-degenerate points of index 0, 1, and 2.}
\label{Morse:fig}
\end{center}
\end{figure}

\begin{figure}
\begin{center}
\includegraphics[width = 3.5 cm] {\iftoggle{BW}{3D-Leveltorus-BW}{240px-3D-Leveltorus}}
\nota{The height function on this torus is a Morse function. It has four critical points of index (from bottom to top) 0, 1, 1, 2. }
\label{Morse_torus:fig}
\end{center}
\end{figure}

\subsection{Alexander theorem} The first 3-manifold to look at is of course $\matR^3$. We prove here that $\matR^3$ is irreducible. 

We need some Morse theory. Let $S\subset \matR^3$ be a closed surface and $f(x,y,z)=z$ be the height function. The function $f$ is a \emph{Morse function} for $S$ if $f|_S$ has finitely many critical points, and at each critical point the Hessian of $f|_S$ is non-singular (the Hessian is read on a chart for $S$, but this definition is chart-independent). The critical point is a local minimum, a saddle, or a local maximum, according to the signature of the Hessian, see Figure \ref{Morse:fig}. These critical points have \emph{index} 0, 1, and 2, respectively.\index{Morse function}

\begin{ex}
A critical point $p$ for $f|_S$ is non-degenerate if and only if $p$ is a regular (\emph{i.e.}~non-critical) point for the Gauss map $S \to S^2$.
\end{ex}

\begin{lemma} Let $S\subset \matR^3$ be a closed surface. After rotating $S$ of an arbitrarily small angle, the height function $f$ is a Morse function for $S$. 
\end{lemma}
\begin{proof}
Consider the Gauss map $\psi\colon S\to S^2$. By Sard lemma there are regular values arbitrarily close to $v=(0,0,1)$. Rotate $S$ so that $v$ is a regular value. Now $\psi^{-1}(v)$ is the set of critical points for $f$ and they are all non-degenerate.
\end{proof}

We are now ready to prove the following.\index{Alexander theorem}

\begin{figure}
\begin{center}
\includegraphics[width = 12.5 cm] {\iftoggle{BW}{Morse_cap-BW}{Morse_cap}}
\nota{The plane $P$ intersects $S$ into circles. Starting from the innermost ones, we cut $S$ along the circles and cap them off by adding pairs of discs. The resulting surface does not intersect $P$ anymore.}
\label{Morse_cap:fig}
\end{center}
\end{figure}

\begin{teo}[Alexander's Theorem] \label{Alexander:teo}
The space $\matR^3$ is irreducible.
\end{teo}
\begin{proof}
Let $S\subset \matR^3$ be a 2-sphere. Up to a small rotation we suppose that the height function $f|_S$ is a Morse function, and after a further small rotation we may suppose that the $k$ critical points of $f|_S$ stay at distinct heights $z_1< \ldots <z_k$. Pick a regular value $u_i \in (z_i, z_{i+1})$ for every $i=1,\ldots, k-1$. The horizontal plane $P$ at height $u_i$ intersects $S$ transversely into circles. Starting from the innermost ones, we cut $S$ along these circles and cap them off by adding pairs of discs as in Figure \ref{Morse_cap:fig}. The resulting surface is disjoint from $P$.

\begin{figure}
\begin{center}
\includegraphics[width = 10 cm] {\iftoggle{BW}{Morse2-BW}{Morse2}}
\nota{After capping off at each $u_i$ we end up with many spheres of these basic types, which clearly bound balls.}
\label{Morse2:fig}
\end{center}
\end{figure}

At every cut a sphere decomposes into two spheres. If we do this for every $i=1,\ldots, k-1$ we end up with many spheres of the types shown in Figure \ref{Morse2:fig}, that clearly bound balls in $\matR^3$.

Now we reverse the process and undo all the cuts: we prove inductively that at each backward step we have a set of spheres bounding balls (note that the balls are not disjoint!). At each backward step we replace two spheres $S_1, S_2$ bounding balls $B_1, B_2$ with one sphere $S$. Isotope $S_1$ and $S_2$, so that they intersect in a disc $D$. If the interiors of $B_1$ and $B_2$ are disjoint, then $S$ bounds the ball $B_1\cup B_2$. If they are not disjoint, then one is contained in the other, say $B_1 \subset B_2$ and $S$ bounds the ball $B_2\setminus \interior{B_1}$.
\end{proof}

\begin{cor}
Every sphere contained in a sub-ball $B\subset M$ of a 3-manifold $M$ bounds a ball $B'\subset B$.
\end{cor}
\begin{cor} \label{both:cor}
Every sphere in $S^3$ bounds a ball on both sides.
\end{cor}
\begin{proof}
Choose two points $p, q \not \in S$ on opposite sides with respect to $S$. We have $S^3 \setminus p = S^3 \setminus q = \matR^3$ hence $S$ is contained in $\matR^3$ in two ways and bounds a ball in each.
\end{proof}

Alexander's Theorem generalises the smooth Jordan curve Theorem to dimension $3$. The situation in higher dimensions is much more problematic: it is still unknown whether every smooth 3-sphere in $\matR^4$ bounds a smooth 4-disc (this is usually called the \emph{Sch\"onflies problem}).\index{Sch\"onflies problem}

\subsection{Prime manifolds}
A connected sum $M_1\# M_2$ is \emph{trivial} if either $M_1$ or $M_2$ is a sphere. \index{three-manifold!prime three-manifold}

\begin{defn}
A connected, oriented 3-manifold $M$ is \emph{prime} if every connected sum $M=M_1\# M_2$ is trivial.
\end{defn}

Being prime is equivalent to be irreducible, with a single exception.

\begin{figure}
\begin{center}
\includegraphics[width = 8 cm] {\iftoggle{BW}{prime-BW}{prime}}
\nota{If $S\subset M$ is non-separating, there is a simple closed curve $\alpha$ intersecting $S$ transversely in one point: here $S$ is drawn as a disc whose boundary should be collapsed to a point, and $\alpha$ is a line whose endpoints should be identified (left). Pick two tubular neighbourhoods of $S$ and $\alpha$ and consider their union $N$ (right).}
\label{prime:fig}
\end{center}
\end{figure}

\begin{prop} \label{S2S1:prop}
Every oriented 3-manifold $M\neq S^2 \times S^1$ is prime if and only if it is irreducible.
\end{prop}
\begin{proof}
The inverse operation of a connected sum $M=M_1\# M_2$ consists of cutting along a separating sphere $S\subset M$ and then capping off the two resulting manifolds $N_1, N_2$ with balls. The capped $N_i$ is $S^3$ if and only if $N_i$ is a ball. Therefore the connected sum is trivial if and only if $S$ bounds a ball on one side. Therefore $M$ is prime if and only if every separating sphere $S\subset M$ bounds a ball.

If $M$ is irreducible, then it is clearly prime. If $M$ is prime and not irreducible, there is a non-separating sphere $S\subset M$. There is a simple closed curve $\alpha\subset M$ intersecting $S$ transversely in one point as in Figure \ref{prime:fig}-(left). Pick two tubular neighbourhoods of $S$ and $\alpha$ as in Figure \ref{prime:fig}-(right): their union is a manifold $N$ with a boundary sphere $\partial N = S'$. The sphere $S'$ is separating and $M$ is prime, hence $S'$ bounds a ball $B$ on the other side and $M=N\cup B$. 

We now prove that $M=S^2\times S^1$. We embed $S\cup \alpha$ naturally in $S^2\times S^1$ as $S = S^2\times y$ and $\alpha = x \times S^1$. Decompose $S^2= D\cup D'$ in two discs and $S^1 = I\cup I'$ in two intervals. The manifold
$N$ also embeds as $S^2\times I \cup D \times S^1$ and its complement $B=D'\times I'$ is a ball. Therefore $M = S^2\times S^1$.
\end{proof}

\subsection{Some irreducible manifolds}
The Alexander theorem generates many more examples of irreducible 3-manifolds.

\begin{prop} \label{connected:boundary:prop}
Every compact three-dimensional submanifold $M\subset S^3$ with connected boundary is irreducible.
\end{prop}
\begin{proof}
Every sphere $S \subset M$ bounds two balls in $S^3$. Since $\partial M$ is connected, it is contained in one of them, so the other is contained in $M$.
\end{proof}

We turn to coverings.

\begin{prop} \label{irreducible:cover:prop}
Let $p\colon  M \to N$ be a covering of 3-manifolds. If $ M$ is irreducible then $N$ also is.
\end{prop}
\begin{proof}
A sphere $S\subset N$ lifts to many spheres in $ M$, each bounding at least one ball. Pick an innermost such ball $B$. We prove that $p(B)$ is a ball with boundary $S$ and we are done.

To do this, note that $p|_{\partial B}$ is a diffeomorphism onto $S$ and $p(\interior B)$ is disjoint from $S$ since $B$ is innermost. This implies that $p|_{B}\colon B \to p(B)$ is a covering and since it has degree one on $S$ it is a diffeomorphism.
\end{proof}

\begin{cor}
Elliptic, flat, hyperbolic 3-manifolds are irreducible.
\end{cor}
\begin{proof}
Their universal covering is diffeomorphic to $S^3$ or $\matR^3$.
\end{proof}

Finally, we consider the exception $S^2\times S^1$. It is not irreducible since it contains a non-separating sphere, but it is prime.

\begin{prop} The manifold $S^2\times S^1$ is prime.
\end{prop}
\begin{proof}
Let $S\subset S^2\times S^1$ be a separating sphere: we must prove that it bounds a ball. It separates $S^2\times S^1$ into two manifolds $M$ and $N$, and on fundamental groups we get $\matZ = \pi_1(M)*\pi_1(N)$. This implies easily that either $\pi_1(M)$ or $\pi_1(N)$ must be trivial: suppose the first. 

Since $M$ is simply connected, a copy $M'$ of $M$ lifts to the universal cover $S^2\times \matR$ of $S^2\times S^1$. We identify $S^2\times\matR = \matR^3\setminus 0$. This copy $M'$ now lies in $\matR^3$, and $\partial M' = S^2$ implies that $M'$ is a ball by Alexander theorem.
\end{proof}

\subsection{Handlebodies and line bundles}
We now introduce some basic, but important, compact 3-manifolds with boundary. A \emph{handlebody} is a connected orientable 3-manifold that decomposes into 0- and 1- handles.\index{handlebody} 

\begin{figure}
\begin{center}
\includegraphics[width = 3.5 cm] {\iftoggle{BW}{Sphere_with_three_handles-BW}{500px-Sphere_with_three_handles}}
\nota{A handlebody of genus $3$.}
\label{handlebody:fig}
\end{center}
\end{figure}

\begin{prop}
The boundary of a handlebody is $S_g$ for some $g\geqslant 0$. Two handlebodies are diffeomorphic if and only if they have the same $g$.
\end{prop}
\begin{proof}
By simplifying handles we may decompose the handlebody into one 0-handle and some $g$ 1-handles, so that the boundary is a genus-$g$ surface, see Figure \ref{handlebody:fig}. Each 1-handle is attached along a pair of discs; the orientability assumption together with Theorem \ref{dischi:teo} easily imply that the result of attaching a 1-handle depends on nothing and the handlebody depends only on $g$ up to diffeomorphisms.
\end{proof}

The \emph{genus} of a handlebody $H_g$ is the genus $g$ of its boundary surface. 
Proposition \ref{connected:boundary:prop} implies the following.

\begin{cor} \label{handle:irr:cor}
Handlebodies are irreducible.
\end{cor}

Some other simple manifolds are irreducible.
\begin{prop} \label{product:irr:prop}
If $g\geqslant 1$ then $S_g\times [0,1]$ is irreducible.
\end{prop}
\begin{proof}
Its universal cover $\matR^2 \times [0,1]$ is irreducible, because its interior is diffeomorphic to $\matR^3$.
\end{proof}

\begin{ex}
If $b\geqslant 1$ the manifold $S_{g,b}\times [0,1]$ is homeomorphic to a handlebody of genus $2g+b-1 = -\chi(S_{g,b}) +1$.
\end{ex}

\subsection{Normal surfaces}
Let $M$ be a compact 3-manifold with (possibly empty) boundary. As every smooth compact manifold, $M$ has a triangulation $T$, made of a certain number $t$ of tetrahedra. We now show that $T$ can be used to treat combinatorially many interesting surfaces in $M$.

A properly embedded surface $S\subset M$ is \emph{transverse} to $T$ if it is transverse to all its simplexes. In particular $S$ does not intersect the vertices of $T$, and it intersects every edge, face, and tetrahedron respectively into a finite number of points, curves, and surfaces. Every properly embedded surface $S\subset M$ can be perturbed to be transverse to $T$.

\begin{figure}
\begin{center}
\includegraphics[width = 7 cm] {\iftoggle{BW}{Normal_surface-BW}{640px-Normal_surface}}
\nota{A normal surface intersects every tetrahedron in triangles or squares.}
\label{normal_surface:fig}
\end{center}
\end{figure}

\begin{defn}
A \emph{normal surface} is a properly embedded surface $S$ transverse to $T$ that intersects every tetrahedron into triangles or squares as in Figure \ref{normal_surface:fig}.\index{normal surface}
\end{defn}

\begin{example}
For every vertex $v$ of $T$ lying in $\interior M$ we may take a small sphere centred at $v$ that intersects every incident tetrahedron in a small triangle: we get a normal sphere. If $v$ lies in $\partial M$ we get similarly a normal disc. A surface of this type is called  \emph{vertex-linking}. Vertex-linking spheres are not very interesting since they bound balls.\index{normal surface!vertex-linking normal surface}  
\end{example}

\begin{ex}
A connected normal surface is vertex-linking if and only if it consists of triangles only (no squares).
\end{ex}

\begin{figure}
\begin{center}
\includegraphics[width = 8 cm] {\iftoggle{BW}{incompressible-BW}{incompressible}}
\nota{We can \emph{surger} a surface $S$ along a disc $D$ with $\partial D = D\cap S$. The operation consists of removing an annular tubular neighbourhood of $\partial D$ in $S$ and adding two parallel copies of $D$. We get a new surface $S'$.}
\label{incompressible:fig}
\end{center}
\end{figure}

\begin{figure}
\begin{center}
\includegraphics[width = 11 cm] {\iftoggle{BW}{incompressible2-BW}{incompressible2}}
\nota{We can also surger a surface $S$ along a disc $D$ touching the boundary in a segment. The result is a new properly embedded surface $S'$.}
\label{incompressible2:fig}
\end{center}
\end{figure}

Let $S\subset M$ be a properly embedded, possibly disconnected, compact surface.
An \emph{elementary transformation} on $S$ is one of the following moves: 
\begin{itemize}
\item the removal of a connected component of $S$ contained in some ball;
\item let $D\subset \interior M$ be a disc with $\partial D = D \cap S$: we \emph{surger} $S$ along $D$ as shown in Figure \ref{incompressible:fig};
\item let $D\subset M$ be a disc with $D\cap (S\cup \partial M) = \partial D = \alpha \cup \beta$, where $\alpha$ is an arc in $S$ and $\beta$ an arc in $\partial M$ as in Figure \ref{incompressible2:fig}-(left): we \emph{surger} $S$ as shown in the figure.
\end{itemize}
Every elementary transformation is local, \emph{i.e.}~it takes place in a ball. It transforms $S$ into a new surface $S'$. 
We fix a triangulation $T$ for $M$.

\begin{figure}
\begin{center}
\includegraphics[width = 12.5 cm] {\iftoggle{BW}{normal2-BW}{normal2}}
\nota{Starting from the innermost curves in $\partial \Delta \cap S$ we surger along discs so that $S\cap \Delta$ consists only of discs and closed components, which we then remove.}
\label{normal2:fig}
\end{center}
\end{figure}

\begin{figure}
\begin{center}
\includegraphics[width = 8 cm] {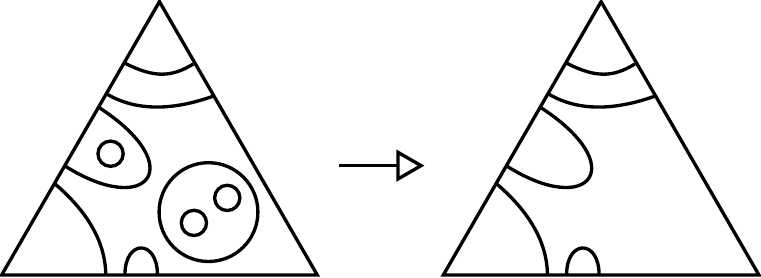}
\nota{The closed curves contained in a triangle are the intersections of small spherical components of $S$ with that triangle. We just remove them.}
\label{normal:fig}
\end{center}
\end{figure}

\begin{figure}
\begin{center}
\includegraphics[width = 9 cm] {\iftoggle{BW}{normal3-BW}{normal3}}
\nota{If the boundary of a (\iftoggle{BW}{light grey}{yellow}) disc $D$ intersects an edge $e$ of $\partial \Delta$ twice, we pick two innermost intersections, that must be in opposite directions. There must be a \iftoggle{BW}{grey}{red} disc as shown here (exercise: use Alexander's theorem), which can be used to push $S$ in the direction indicated by the arrow. The resulting surface is shown in the right. If $e\subset \partial M$, we cannot push $S$ outside $\partial M$: instead, we surger along the \iftoggle{BW}{grey}{red} disc.}
\label{normal3:fig}
\end{center}
\end{figure}

\begin{prop}
Every properly embedded surface $S\subset M$ becomes normal after finitely many isotopies and elementary transformations.
\end{prop}
\begin{proof}
Put $S$ transverse to $T$. It intersects every tetrahedron $\Delta$ into surfaces as in Figure \ref{normal2:fig}-(left). The intersection $S\cap \partial \Delta$ consists of closed curves: the innermost ones bound discs in $\partial \Delta$, which may be pushed a bit inside $\Delta$ and then used to surger $S$. We do this with all curves (starting with the innermost ones) to transform $S$ as in Figure \ref{normal2:fig}, so that at the end $S\cap\Delta$ consists only of discs and closed surfaces. Then we remove the closed surfaces: this is an elementary transformation since $\Delta$ is a ball. We do this for every tetrahedron $\Delta$.

Now $S$ intersects every tetrahedron $\Delta$ in discs. Consider the boundary curve of one such disc: if it is contained in a face as in Figure \ref{normal:fig}-(left), it belongs to one or two discs (depending on whether the face belongs to $\partial M$ or not), which form a component of $S$ contained in a ball that can be removed. 

Suppose the boundary curve crosses twice an edge of $\Delta$: if the edge lies in the interior of $M$, we isotope $S$ as shown in Figure \ref{normal3:fig} to decrease the number of intersections of $S$ with the 1-skeleton of $T$; if the edge lies in $\partial M$ we simplify analogously by surgerying along the \iftoggle{BW}{grey}{red} disc shown in the figure. In all cases then we go back to the first step of the algorithm. 

After finitely many steps we get a surface $S$ intersecting every tetrahedron $\Delta$ in discs, whose boundary curves intersect every edge of $\Delta$ at most once. One such disc is either a triangle or a square.
\end{proof}

\begin{figure}
\begin{center}
\includegraphics[width = 11 cm] {\iftoggle{BW}{normal4-BW}{normal4}}
\nota{If we cut a tetrahedron $\Delta$ along triangles and squares we get arbitrarily many prisms with triangular or quadrilateral basis, and at most $6$ other pieces (here we get $4$ pyramids and two esahedra).}
\label{normal4:fig}
\end{center}
\end{figure}

\subsection{Interesting surfaces cannot have too many components}
We know from Proposition \ref{non:ori:prop} that a compact orientable $M$ cannot contain too many disjoint non-orientable surfaces. It can however contain arbitrarily many orientable surfaces: for instance, small surfaces contained in disjoint balls. A crucial aspect of 3-manifolds theory is that $M$ cannot contain arbitrarily many ``interesting'' surfaces, as we now see.

Two disjoint connected diffeomorphic surfaces $\Sigma,\Sigma'\subset M$ are \emph{parallel} if they cobound a region diffeomorphic to $\Sigma\times [0,1]$ with $\Sigma= \Sigma\times 0$ and $\Sigma'=\Sigma\times 1$. Two parallel surfaces are obviously isotopic. 

Let $T$ be a triangulation of a compact $M$ with (possibly empty) boundary. Let $t$ be the number of tetrahedra in $T$ and set $b_2 = \dim H_2(M,\partial M,\matZ_2)$.

\begin{lemma} \label{non:parallel:lemma}
Let $S$ be an orientable normal surface. If $S$ has more than $10t+b_2$ components, then two components $\Sigma, \Sigma'$ of $S$ are parallel and cobound a $\Sigma \times [0,1]$ which is disjoint from the other components.
\end{lemma}
\begin{proof}
The complement $M\setminus S$ intersects every tetrahedron $\Delta$ into polyhedra: there are many prisms lying between parallel triangles or squares, and at most $6$ other pieces, see Figure \ref{normal4:fig}. These at most 6 pieces are adjacent to at most 1+1+1+1+3+3 = 10 triangles and squares.

This implies that, except at most $10t$ of them, the components of $S$ are only adjacent (on both sides) to prisms. These prisms glue to form $I$-bundles. Therefore at least $b_2+1$ components of $S$ are adjacent to $I$-bundles on both sides. The twisted $I$-bundles are at most $b_2$ by Proposition \ref{non:ori:prop}, and each is adjacent to one surface. Therefore at least one surface is adjacent to a product $I$-bundle $\Sigma \times [0,1]$. \end{proof}

We get a topological corollary. A \emph{ball with holes} is a 3-manifold obtained by removing some $k\geqslant 0$ disjoint small open balls from a ball.
A \emph{sphere system} for a 3-manifold $M$ is a surface $S\subset \interior{M}$ consisting of disjoint spheres, such that no component of $M\setminus S$ is a ball with holes disjoint from $\partial M$. (Balls with holes adjacent to spherical components of $\partial M$ are allowed.) \index{ball with holes}\index{sphere system}

\begin{cor} \label{K:cor}
Let $M$ be a compact orientable 3-manifold that does not contain any non-separating sphere. There is a $K>0$ such that every sphere system in $M$ contains less than $K$ spheres.
\end{cor}
\begin{proof}
Pick a triangulation $T$ of $M$ with some $t$ tetrahedra: we prove that $K=10t + b_2+1$ works, with $b_2 = \dim H_2(M,\partial M;\matZ_2)$. Suppose by contradiction that there is a sphere system $S$ with $\geqslant K$ spheres. Via isotopies and elementary transformations we transform $S$ into a normal surface $S'$.

We now examine the effect of elementary transformations in detail. No component of $S$ is contained in a ball, otherwise by Alexander theorem an innermost component of $M\setminus S$ would be a ball. Therefore an elementary transformation cannot cancel a component of $S$. 

A surgery along a disc splits a sphere $S_0$ into two spheres $S_1, S_2$. 
We now prove that by substituting $S_0$ with either $S_1$ or $S_2$ we still get a sphere system. If we push $S_1$ and $S_2$ away from $S_0$ the surfaces $S_0, S_1, S_2$ altogether bound a ball with two holes $B$. 

By our hypothesis on $M$ the spheres $S_0, S_1$, and $S_2$ are separating. Let $N_i$ be the component of $M\setminus (S_0\cup S_1\cup S_2)$ adjacent to $S_i$ distinct from $B$, for $i=0,1,2$. If both $N_1, N_2$ are balls with holes disjoint from $\partial M$, then $N_1\cup N_2 \cup B$ also is, which is excluded. Therefore one, say $N_1$, is not a ball with holes disjoint from $\partial M$. On the other side, if $N_2 \cup B \cup N_0$ were a ball with holes disjoint from $\partial M$, then $N_0$ would also be (by Alexander's theorem), which is excluded. Therefore by substituting $S_0$ with $S_1$ we still get a sphere system.

This proves that the final normal surface $S'$ contains a sphere system with the same number $\geqslant K$ of components as $S$. Lemma \ref{non:parallel:lemma} gives a contradiction.
\end{proof}

\subsection{Prime decomposition}
We now show that the connected sum operation on oriented three-manifolds behaves like the product of natural numbers: every object decomposes uniquely into prime factors.\index{prime decomposition of three-manifolds} 

\begin{figure}
\begin{center}
\includegraphics[width = 6 cm] {\iftoggle{BW}{holes-BW}{holes}}
\nota{A reducing set of spheres (\iftoggle{BW}{grey}{red}) for $M = M_1\# M_2 \# (S^2\times S^1)$. Here $B$ is a ball with four holes and $N_i$ is $M_i$ with one hole.}
\label{holes:fig}
\end{center}
\end{figure}

\begin{teo} \label{prime:teo}
Every compact oriented 3-manifold $M$ with (possibly empty) boundary decomposes into prime manifolds:
$$M = M_1\# \ldots \# M_k$$
This list of prime factors is unique up to permutations and adding/removing copies of $S^3$.
\end{teo}
\begin{proof}
We first show the existence of a decomposition. If $M$ contains a non-separating sphere, then the proof of Proposition \ref{S2S1:prop} shows that $M = M'\# (S^2\times S^1)$. Since $H_1(M) = H_1(M')\oplus \matZ$, up to factoring finitely many copies of $S^2\times S^1$ we may suppose that every sphere in $M$ is separating.

If $M$ is prime we are done. If not, it decomposes as $M=M_1\# M_2$. We keep decomposing each factor until all factors are prime: this process must end, because a decomposition $M=M_1\# \ldots \#M_k$ gives rise to a system of $k-1$ spheres, and $k$ cannot be arbitrarily big by Corollary \ref{K:cor}.

We turn to uniqueness. Let 
$$M = M_1\# \ldots \# M_k \#_h(S^2\times S^1), \qquad M = M_1'\# \ldots \# M_{k'}'\#_{h'}(S^2\times S^1)$$
be two prime decompositions with $M_i, M_j' \neq S^2\times S^1$, so $M_i, M_j'$ are irreducible for all $i,j$. We say that a set $S \subset M$ of disjoint spheres is a \emph{reducing set of spheres} for the decomposition $M=M_1\# \ldots \# M_k\#_h(S^2\times S^1)$ if $M\setminus S$ consists of precisely one $M_i$ with some holes for each $i$, and some balls with holes disjoint from $\partial M$. An example is drawn in Figure \ref{holes:fig}. In general, we may construct $S$ by taking the spheres of the prime decomposition, plus one non-separating sphere inside each $S^2\times S^1$ summand. Similarly, let $S'$ be a reducing set of spheres for the other decomposition.

The first observation we make is that if we add to $S$ any sphere $\Sigma$ disjoint from $S$, then we still get a reducing set of spheres for the same decomposition as before. This is because $\Sigma$ is contained in a holed $M_i$ or $S^3$, and since $M_i$ is irreducible $\Sigma$ bounds a ball $B$ there. Therefore by adding $\Sigma$ we still get the same holed $M_i$, plus a possibly holed (if $B\cap S \neq \emptyset$) ball $B$.

We assume $S$ and $S'$ intersect transversely in circles and pick an innermost circle in a component of $S$ bounding a disc $D\subset S$. We surger $S'$ along $D$, thus substituting a component $S_0'$ of $S'$ with two spheres $S_1' \sqcup S_2'$. We check that the result is another sphere system for the same decomposition. We isotope the spheres $S_0', S_1', S_2'$ so that they are disjoint and cobound a ball with two holes $B_2$: the system $S'\sqcup S'_1 \sqcup S'_2$ is still reducing by the observation above. The removal of $S_0'$ then adds $B_2$ to the outside of $S_0'$, and this is equivalent to making one more hole there. 

After finitely many surgeries we get $S\cap S'=\emptyset$. By the same observation above $S\cup S'$ is a reducing set of spheres for both decompositions: therefore the pieces $M_i$ and $M_j'$ of the decompositions are pairwise diffeomorphic.

Finally we must have $h=h'$ since $M=N\#_h(S^2\times S^1) = N\#_{h'}(S^2\times S^1)$ and $H_1(M) = H_1(N) \oplus \matZ^h = H_1(N) \oplus \matZ^{h'}$.
\end{proof}

This important result is known as the \emph{Prime decomposition Theorem} for 3-manifolds: the existence of a decomposition was proved by Kneser \cite{Kn} in 1929, and its uniqueness by Milnor \cite{Mi} in 1962. In light of this result, topologists have since long restricted their attention to prime 3-manifolds, or almost equivalently to irreducible 3-manifolds.

The strategy of cutting canonically a three-manifold along surfaces has proved successful with spheres, so we now try to do the same with other surfaces. We start by studying properly embedded discs, which of course occur only in manifolds with boundary. It is tempting to guess that discs should behave roughly like spheres, because by doubling a 3-manifold along its boundary we transform properly embedded discs into spheres. We now see that this is indeed the case. As for spheres, we need to distinguish the interesting ones, that we call \emph{essential}, from the others.

\subsection{Essential discs} \label{essential:discs:subsection}
Let $M$ be a compact 3-manifold with (possibly empty) boundary. A properly embedded surface $S\subset M$ is \emph{$\partial$-parallel} if it is obtained by slightly pushing inside $M$ the interior of a compact surface $S'\subset \partial M$, possibly with boundary.\index{surface!$\partial$-parallel surface} 

In what follows $D$ and $S$ are properly embedded. 
We say that:\index{sphere!essential sphere}\index{disc!essential disc}
\begin{itemize}
\item a sphere $S\subset M$ is \emph{essential} if it does not bound a ball, 
\item a disc $D\subset M$ is \emph{essential} if it is not $\partial$-parallel.
\end{itemize}
Now, the manifold $M$ is:\index{three-manifold!$\partial$-irreducible three-manifold}
\begin{itemize}
\item \emph{irreducible} if it does not contain essential spheres, and
\item \emph{$\partial$-irreducible} if it does not contain essential discs.
\end{itemize}

To warm up, we show the following. A \emph{solid torus} is a handlebody of genus 1, that is a three-manifold diffeomorphic to $D\times S^1$.\index{solid torus}

\begin{prop} \label{irreducible:solid:torus:prop}
Let $M$ be a compact, irreducible, orientable 3-manifold with boundary, and let $D\subset M$ be an essential disc. Let $\Sigma \subset \partial M$ be the boundary component containing $\partial D$. Then:
\begin{itemize}
\item the curve $\partial D$ is non-trivial in $\Sigma$;
\item if $\Sigma$ is a torus then $M$ is a solid torus.
\end{itemize}
\end{prop}
\begin{proof}
If $\partial D$ bounds a disc $D'\subset \Sigma$ then $D\cup D'$ is a sphere, which bounds a ball because $M$ is irreducible; this ball furnishes a parallelism between $D$ and $D'$ and hence $D$ is $\partial$-parallel, a contradiction. 

If $\Sigma$ is a torus, by surgerying $\Sigma$ along $D$ we get a sphere which must bound a ball $B$. Therefore $M$ is obtained by adding a one-handle to $B$, that is it is a solid torus.
\end{proof}

\subsection{Decomposition along discs}
We now show that essential discs behave roughly like essential spheres, in the sense that there is a kind of prime decomposition theorem also for discs. After stating and proving this theorem, we will essentially forget about essential spheres and discs and focus on the 3-manifolds that do not contain them.

Let $M$ be a compact 3-manifold with (possibly empty) boundary. A \emph{disc system} in $M$ is a set of pairwise disjoint and non-parallel essential discs. 
 
\begin{prop} \label{K2:prop}
There is a $K>0$ such that every disc system in $M$ cannot contain more than $K$ discs.
\end{prop}
\begin{proof}
The proof is analogous to that of Corollary \ref{K:cor}
\end{proof}

We now want to cut irreducible manifolds along essential discs.

\begin{rem}
The opposite operation of cutting a manifold along a properly embedded disc is a 1-handle addition.
\end{rem}

We now state the analogue of Theorem \ref{prime:teo} in this context. Note that $M$ is irreducible by hypothesis: we have already eliminated the essential spheres, and we now remove the essential discs. 

\begin{teo}
Every compact oriented irreducible 3-manifold $M$ is obtained by adding 1-handles to a finite list
$$M_1,\ldots, M_k$$
of connected irreducible and $\partial$-irreducible 3-manifolds. The list is unique
up to permutations and adding/removing balls.
\end{teo}
\begin{proof}
The proof is similar to Theorem \ref{prime:teo}. 
If we cut $M$ along a maximal disc system we get some connected manifolds $M_1,\ldots, M_k$. Each $M_i$ is $\partial$-irreducible because the disc system is maximal, and it is irreducible because every sphere in $M_i\subset M$ bounds a ball $B$ in $M$, and $B$ is necessarily contained in $M_i$.

We prove uniqueness. Pick two sets of discs $S, S' \subset M$ such that by cutting along each we get irreducible and $\partial$-irreducible components. We prove that these components are the same up to adding/removing balls. We assume that $S$ and $S'$ intersect transversely. 

\begin{figure}
\begin{center}
\includegraphics[width = 7 cm] {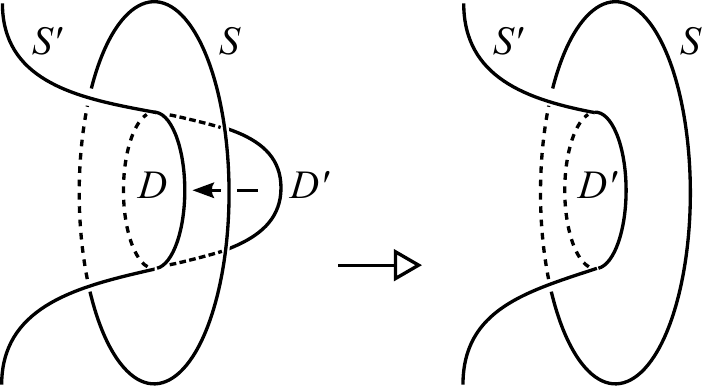}
\nota{The discs $D$ and $D'$ form a sphere which bounds a ball since $M$ is irreducible. We can use this ball to isotope $S'$ away from $S$.}
\label{irreducible_isotopy:fig}
\end{center}
\end{figure}

A circle $S\cap S'$ which is innermost with respect to $S$ is bounded by two sub-discs $D\subset S$ and $D'\subset S'$, which together form a sphere: this sphere in turn bounds a ball (because $M$ is irreducible) as in Figure \ref{irreducible_isotopy:fig}. We can use the ball to isotope $S'$ away from $S$ and reduce the intersections in $S\cap S'$. After finitely many isotopies there are no circles in $S\cap S'$.

We are left to consider the arcs in $S\cap S'$. We end precisely as in Theorem \ref{prime:teo}. We observe that by adding more discs to $S$ or $S'$ we get a set of discs that $\partial$-reduces $M$ into the same list of $M_i$'s, plus possibly some balls.

Thanks to an innermost argument, we may surger some disc $D_0\subset S$ along an arc in $D_0\cap S'$ and transform it into two discs $D_1 \sqcup D_2$, both disjoint from $S$ and with less intersections with $S'$. The set $S\sqcup D_1 \sqcup D_2$ is still a reducing set for the $M_i$'s. If we remove $D_0$ from this new system, in the complement we only remove a ball: hence we can replace $D_0$ with $D_1\sqcup D_2$. After finitely many steps we get $S\cap S' = \emptyset$ and therefore $S\cup S'$ is a reducing set producing the same $M_i$'s of $S$ and $S'$, up to balls.
\end{proof}

\subsection{Projective planes}
The previous discussions show how to deal with spheres and discs, and we now turn to the only remaining compact surface with $\chi>0$: the projective plane. There is not much to say about it.

\begin{prop}
The only compact irreducible orientable 3-manifold $M$ containing $\matRP^2$ is $\matRP^3$.
\end{prop}
\begin{proof}
The orientable $M$ contains the orientable tubular neighbourhood $N=\matRP^2 \timtil I$ of $\matRP^2$. The boundary $\partial N$ is the orientable cover $S^2$ or $\matRP^2$. Since $M$ is irreducible, the sphere $\partial N$ bounds a ball $B$. Therefore $M$ is uniquely determined as the manifold obtained by capping off the boundary sphere of $N$. 

The manifold $\matRP^3$ is irreducible (because its cover $S^3$ is) and contains $\matRP^2$, hence we must have $M=\matRP^3$.
\end{proof}

This proposition does not hold for non-orientable three-manifolds, which might contain $\matRP^2$ with a product neighbourhood $\matRP^2\times I$. These 
cases are difficult to deal with, so one typically excludes them \emph{a priori} by requiring the manifold to be \emph{$\matRP^2$-irreducible}, which means that it does not contain two-sided projective planes.\index{three-manifold!$\matRP^2$-irreducible three-manifold}

\section{Incompressible surfaces}
We have proved that every compact oriented 3-manifold decomposes along essential spheres and discs, into some canonical pieces that do not contain essential spheres or discs anymore. 

We would like to pursue this strategy with the next simplest surfaces, namely tori and annuli. To this purpose we define the important notion of \emph{incompressible} surface in a three-manifold, which applies to all surfaces of non-positive Euler characteristic.

\subsection{Incompressible surfaces}
Throughout all this section $M$ denotes a compact orientable 3-manifold with (possibly empty) boundary. Let $S\subset M$ be a properly embedded orientable surface. A \emph{compressing disc} for $S$ is a disc $D\subset M$ with $\partial D = D \cap S$, such that $\partial D$ does not bound a disc in $S$. With this hypothesis, the surgery in Figure \ref{incompressible:fig} is called a \emph{compression}: it transform $S$ into a new surface $S'\subset M$ which is simpler than $S$.\index{compressing disc}

\begin{prop}  \label{S':simpler:prop}
The surface $S'$ may have one or two components $S_i'$, and $\chi(S_i')>\chi(S)$ for each component. 
\end{prop}
\begin{proof}
We have $\chi(S') = \chi(S) + 2$. If $S'$ has one component we are done, so suppose $S' = S_1' \sqcup S_2'$. Since $\partial D$ did not bound a disc in $S$, no $S_i'$ is a sphere, hence $\chi(S_i')\leqslant 1$ that implies $\chi(S_i')> \chi(S)$ for $i=1,2$.
\end{proof}

A properly embedded connected orientable compact surface $S\subset M$ with $\chi(S)\leqslant 0$ is \emph{compressible} if it has a compressing disc, and \emph{incompressible} otherwise. See Figure \ref{incompressible3:fig}-(top).\index{surface!compressible and incompressible surface} 

\begin{figure}
\begin{center}
\includegraphics[width = 11 cm] {\iftoggle{BW}{incompressible3-BW}{incompressible3}}
\nota{A surface $S$ is incompressible ($\partial$-incompressible) if the existence of a disc $D$ as in the top-left (bottom-left) picture implies the existence of another disc $D'\subset S$ as in the top-right (bottom-right) picture. If in addition $M$ is irreducible ($\partial$-irreducible) the two discs $D$ and $D'$ form a sphere (disc) which bounds a ball (is $\partial$-parallel), and hence by substituting $D'$ with $D$ we get two isotopic surfaces.}
\label{incompressible3:fig}
\end{center}
\end{figure}

\begin{cor} Let $S\subset M$ be any properly embedded orientable surface. After compressing it a finite number of times it transforms into a disjoint union of spheres, discs, and incompressible surfaces.
\end{cor}
\begin{proof}
We compress $S$ as much as we can; after finitely many steps we must stop because of Proposition \ref{S':simpler:prop}.
\end{proof}

\begin{rem}
By definition, an orientable 3-manifold is $\partial$-irreducible if and only if its boundary consists of incompressible components. (Strictly speaking, a boundary component is not properly embedded and hence cannot be incompressible: we implicitly push it inside the 3-manifold.)
\end{rem}

A simple incompressibility criterion is the following.

\begin{prop} \label{pi1:prop}
Let $S\subset M$ be an orientable, connected, properly embedded surface with $\chi(S)\leqslant 0$. If the map $\pi_1(S) \to \pi_1(M)$ induced by inclusion is injective, then $S$ is incompressible.
\end{prop}
\begin{proof}
Suppose that a disc $D$ compresses $S$. Its boundary $\partial D$ is a non-trivial element in $\pi_1(S)$ because it does not bound a disc in $S$ by Proposition \ref{disco:prop}, but is clearly trivial in $\pi_1(M)$: a contradiction.
\end{proof}

The converse is also true, but its proof is much harder! We will complete it at the end of this chapter. For the moment we content ourselves with the following.

\begin{prop} \label{incompressible:boundary:prop}
If $S\subset M$ is incompressible, every component of $\partial S$ is non-trivial in $\partial M$.
\end{prop}
\begin{proof}
If a component of $\partial S$ is trivial in $\partial M$, it bounds a disc $D\subset \partial M$ there. By taking an innermost one we get $D\cap S = \partial D$, and by pushing $D$ inside $M$ we find a compressing disc for $S$.
\end{proof}

\subsection{Tori}
The first closed surface to look at is the torus.

\begin{figure}
\begin{center}
\includegraphics[width = 9 cm] {\iftoggle{BW}{torus-BW}{torus}}
\nota{A torus in an irreducible 3-manifold is either incompressible, or it bounds a solid torus (left), or is contained in a ball (right).}
\label{torus_irred:fig}
\end{center}
\end{figure}

\begin{prop} \label{T:M:prop}
Let $T\subset M$ be a torus in an irreducible 3-manifold. One of the following holds:
\begin{enumerate}
\item $T$ is incompressible,
\item $T$ bounds a solid torus,
\item $T$ is contained in a ball.
\end{enumerate}
\end{prop}
\begin{proof}
If $T$ is not incompressible, it compresses along a disc $D$. The result of the compression is necessarily a sphere $S\subset M$ which bounds a ball $B$ since $M$ is irreducible. If $B$ is disjoint from $T$, then $T$ bounds a solid torus as in Fig.\ref{torus_irred:fig}-(left). If $B$ contains $T$, then case (3) holds as shown in Figure \ref{torus_irred:fig}-(right).
\end{proof}

If $M=S^3$ the case (3) can be excluded: if $T$ is contained in a ball $B$, then it bounds a solid torus ``on the outside'' using the complementary ball $S^3\setminus B$. We will also exclude the case (1) when we show that $S^3$ contains no incompressible surfaces.

\subsection{$\partial$-incompressible surfaces}
There is of course also a $\partial$-version of incompressibility. Let $S\subset M$ be a properly embedded orientable surface in a 3-manifold $M$. A \emph{$\partial$-compressing disc} for $S$ is a disc $D$ with $\partial D = \alpha \cup \beta$, where $\alpha$ lies in $S$ and $\beta$ in $\partial M$ as in Figure \ref{incompressible2:fig}-(left); we also require that there is no sub-disc $D'\subset S$ with $\partial D' = \alpha\cup \beta' $ and $\beta'\subset\partial S$. 
The move in Figure \ref{incompressible2:fig} is a \emph{$\partial$-compression} and transforms $S$ into a surface $S'\subset M$ simpler than $S$:

\begin{prop}
The surface $S'$ may have one or two components $S_i'$, and $\chi(S_i')>\chi(S)$ for each component. 
\end{prop}
\begin{proof}
We have $\chi(S') = \chi(S) + 1$. If $S'$ has one component we are done, so suppose $S' = S_1' \sqcup S_2'$. Since $\alpha$ did not bound a disc in $S$, no $S_i'$ is a disc, hence $\chi(S_i')\leqslant 0$ that implies $\chi(S_i')> \chi(S)$ for $i=1,2$.
\end{proof}

A properly embedded connected orientable compact $S\subset M$ with $\chi(S)\leqslant 0$ is \emph{$\partial$-compressible} if it has a $\partial$-compressing disc, and \emph{$\partial$-incompressible} otherwise. See Figure \ref{incompressible3:fig}.\index{surface!$\partial$-compressible and $\partial$-incompressible surface}

\begin{cor} Let $S\subset M$ be any properly embedded orientable surface. After $\partial$-compressing it a finite number of times it transforms into a disjoint union of spheres, discs, and $\partial$-incompressible surfaces.
\end{cor}

\subsection{Annuli}
The first non-closed surface to look at is the annulus. Let a \emph{tube} be the tubular neighbourhood of a properly embedded arc.

\begin{figure}
\begin{center}
\includegraphics[width = 12.5 cm] {\iftoggle{BW}{annulus-BW}{annulus}}
\nota{An annulus $A$ in an irreducible and $\partial$-irreducible 3-manifold is either incompressible and $\partial$-incompressible, or is parallel to an annulus in $\partial M$ (left), or bounds a tube (centre), or is contained in a ball intersecting $\partial M$ in a disc (right).}
\label{annulus:fig}
\end{center}
\end{figure}

\begin{prop} \label{A:M:prop}
Let $A\subset M$ be a properly embedded annulus in an irreducible and $\partial$-irreducible 3-manifold. One of the following holds:
\begin{enumerate}
\item $A$ is incompressible and $\partial$-incompressible,
\item $A$ bounds a tube,
\item $A$ is parallel to an annulus in $\partial M$,
\item $A$ is contained in a ball $B$ intersecting $\partial M$ in a disc.
\end{enumerate}
\end{prop}
\begin{proof}
If $A$ compresses along a disc $D$, it transforms into two discs that are parallel to two discs $D_1, D_2 \subset \partial M$ since $M$ is $\partial$-irreducible. If $D_1\cap D_2 = \emptyset$ then $A$ bounds a tube as in Figure \ref{annulus:fig}-(right); if $D_1 \subset D_2$ then $A$ is contained in a ball $B$ intersecting $\partial M $ in $D_2$.

 If $A$ $\partial$-compresses along a disc $D$, it transforms into a disc which is again $\partial$-parallel and hence $A$ is as in Figure \ref{annulus:fig}-(left) or bounds a tube.
\end{proof}

As a corollary, we get a simple criterion for detecting incompressible and $\partial$-incompressible annuli:

\begin{cor} \label{A:non-parallel:cor}
Let $A\subset M$ be a properly embedded annulus in an irreducible and $\partial$-irreducible 3-manifold. If the components of $\partial A$ are non-trivial and non-parallel in $\partial M$, the annulus $A$ is incompressible and $\partial$-incompressible.
\end{cor}

\subsection{Handlebodies}
We study some examples. We investigate the incompressible and $\partial$-incompressible surfaces in the three-sphere and in the handlebodies. Recall that a ($\partial$)-incompressible surface is always compact, orientable, connected, properly embedded, and with non-positive Euler characteristic by hypothesis.

\begin{prop}
There are no incompressible surfaces in $\matR^3$.
\end{prop}
\begin{proof}
Let $S$ be a surface in $\matR^3$. By applying as is the proof of Alexander Theorem \ref{Alexander:teo} we find that $S$ transforms into spheres after surgerying along discs. Therefore $S$ compresses somewhere.
\end{proof}

A high genus closed surface may be embedded in $\matR^3$ in a rather complicated way (which may be hard to imagine) and the proposition says that no matter how intricate the surface is, there is always a compressing disc somewhere that one can use to simplify the picture. After finitely many compressions the complicated surface is atomised into some trivial spheres. 

\begin{cor} \label{no:inc:B:cor}
There are no incompressible surfaces in $S^3$.
\end{cor}
\begin{cor} \label{B:cor}
There are no incompressible surfaces in the ball $B$.
\end{cor}
\begin{proof}
Use Proposition \ref{incompressible:boundary:prop}.
\end{proof}

\begin{cor}
Every torus in $S^3$ bounds a solid torus and every properly embedded annulus in $B$ bounds a tube.
\end{cor}
\begin{proof}
Apply Propositions \ref{T:M:prop} and \ref{A:M:prop}. Since the manifold is $S^3$ or $B$, the cases (3) or (4) easily imply (2).
\end{proof}

The solid torus in $S^3$ and the tube in $B$ may of course be knotted! We now turn to handlebodies. Recall that a handlebody of genus $g\geq 1$ is irreducible by Corollary \ref{handle:irr:cor}, but it is clearly not $\partial$-irreducible.

\begin{prop}[Handlebodies] \label{handlebodies:prop}
The genus-$g$ handlebody $H_g$ contains no incompressible and $\partial$-incompressible surfaces.
\end{prop}
\begin{proof}
Suppose that $S\subset H_g$ is incompressible and $\partial$-incompressible. Pick disjoint essential discs $D_1,\ldots, D_g$ that cut $H_g$ into a ball $B$. Put $S$ in transverse position with respect to $\sqcup_i D_i$, so that the intersection of $S$ with $\sqcup_i D_i$ consists of circles and properly embedded arcs. Since $H_g$ is irreducible and $S$ is incompressible and $\partial$-incompressible, all these intersections can be removed by an isotopy of $S$, as explained below. Then Corollary \ref{B:cor} gives a contradiction.

The intersections removal goes as follows. Let $\alpha$ be one arc or circle in $S \cap D_j$. Since $S$ is incompressible and $\partial$-incompressible, there is a disc $D'\subset S$ bounded by $\alpha$ (if $\alpha$ is a circle) or by $\alpha \cup \beta$ with $\beta\subset \partial S$ (if $\alpha$ is an arc) as in Figure \ref{incompressible3:fig}. After substituting $\alpha$ with an innermost intersection of $D'$ with $\sqcup_i D_i$ we may suppose that $D'$ is entirely contained in $B$ and hence we may eliminate $\alpha$ with an isotopy as explained in Figure \ref{incompressible3:fig}. 
\end{proof}

We can be more specific on the solid torus $H_1$.

\begin{prop} \label{incompressible:solid:torus:prop}
Every incompressible surface in a solid torus is a $\partial$-parallel annulus.
\end{prop}
\begin{proof}
Every incompressible surface $S\subset H_1$ is $\partial$-compressible by Proposition \ref{handlebodies:prop}. If we $\partial$-compress it, we get either a disc or an incompressible surface again. Therefore $S$ is constructed iteratively from some discs by a sequence of moves opposite to the $\partial$-compression, that produce incompressible surfaces at each step. 
One easily sees that the only incompressible surface that this move can produce at the first step is a $\partial$-parallel annulus (from a single disc), and then one gets stuck.
\end{proof}

\begin{warning}
The handlebody $H_2$ contains many complicated properly embedded incompressible surfaces! However, these are all $\partial$-compressible. The main difference between $H_1$ and $H_2$ is that $\pi_1(H_2)=\matZ*\matZ$ is a free group of rank two and contains many free groups of arbitrarily high rank, so there is a lot of space for $\pi_1$-injective incompressible surfaces in $H_2$ with boundary (whose fundamental group is free).
\end{warning}

\subsection{Line bundles}
We now turn to product line bundles. Recall that $S_g$ is the closed orientable surface of genus $g$. 

\begin{prop}[Line bundles] \label{line:bundles:prop}
Fix $g\geqslant 1$. The product $M=S_g\times [-1,1]$ is irreducible and $\partial$-irreducible. The incompressible and $\partial$-incompressible surfaces in $M$ are precisely the following (up to isotopy):
\begin{itemize}
\item the \emph{horizontal} surface $S_g\times 0$,
\item a \emph{vertical} annulus $\gamma \times [-1,1]$ for each non-trivial simple closed curve $\gamma\subset S_g$.
\end{itemize}
\end{prop}
\begin{proof}
We know that $M$ is irreducible by Proposition \ref{product:irr:prop}. The horizontal surface and the vertical annuli are both incompressible and $\partial$-incompressible by 
Propositions \ref{pi1:prop} and \ref{A:non-parallel:cor}. Propositions \ref{pi1:prop} also shows that $M$ is $\partial$-irreducible.

\begin{figure}
\begin{center}
\includegraphics[width = 12.5 cm] {\iftoggle{BW}{prisms-BW}{prisms}}
\nota{A move that decreases the intersection number of $S$ with the vertical edges (left). The surface $S$ intersects every prism into horizontal triangles (centre) or vertical rectangles (right).}
\label{prisms:fig}
\end{center}
\end{figure}

We now prove that every incompressible and $\partial$-incompressible $S$ is either vertical or horizontal, up to isotopy. We use a version of normal surface theory, with prisms instead of tetrahedra. 

A triangulation $\Delta$ of $S_g$ determines a decomposition of $S_g\times [-1,1]$ into prisms. Every vertex $v$ of $\Delta$ gives rise to a vertical edge $e=v\times [-1,1]$. We suppose that $S$ is transverse to the prisms and has minimum intersection number with the vertical edges. The intersection of $S$ with a vertical rectangular face consists of arcs and circles. The simplification argument described below shows that after an isotopy we get only horizontal or vertical arcs (that is, arcs joining opposite sides of the rectangle).
By further isotopies we can also supppose that $S$ intersects every horizontal triangle of the prism into arcs joining distinct sides. Therefore $S$ intersects necessarily every prism into either horizontal triangles or vertical rectangles as in Figure \ref{prisms:fig}-(centre) and (right). These pieces glue up to give a horizontal or vertical surface.

Here is the simplification argument. There cannot be arcs with endpoints in two consecutive edges, because they could be removed by an isotopy as in Figure \ref{prisms:fig}-(left) contradicting minimality. Analogously there are no arcs with both endpoints in the same vertical edge. Circles and arcs with endpoints in the same horizontal edge are removed as in the proof of Proposition \ref{handlebodies:prop}.  
\end{proof}

\section{Haken manifolds}
There are two classes of irreducible closed three-manifolds: those that contain incompressible surfaces, and those that do not. Both classes are very important and contain a wealth of interesting manifolds.

The manifolds belonging to the first class are called \emph{Haken manifolds} and are somehow easier to study, because they can be cut into balls via a standard procedure called \emph{hierarchy}: we cut the manifold along an essential surface, then along another, and we iterate until we get balls. We study these manifolds here. We also prove the converse of Proposition \ref{pi1:prop}, that is that a closed surface is incompressible if and only if it is $\pi_1$-injective, see Theorem \ref{incompressible:teo}.

\subsection{Introduction}
If not otherwise mentioned, all the 3-manifolds $M$ we will consider in this section will be connected, compact, oriented and with (possibly empty) boundary.
We introduce a definition.

\begin{defn} A \emph{Haken manifold} is a compact, connected, oriented 3-manifold $M$ with (possibly empty) boundary, which is irreducible, $\partial$-irreducible, and contains an incompressible and $\partial$-incompressible surface.\index{three-manifold!Haken three-manifold}
\end{defn} 

The reader should not be frightened by the abundance of adjectives: this definition is really clever because it summarises various reasonable hypothesis in a unique word. The rest of this chapter is mainly devoted to the study of Haken manifolds. We start by looking at their boundaries.

\begin{prop}
Every boundary component $X$ of a Haken manifold $M$ has $\chi(X)\leqslant 0$ and is incompressible.
\end{prop}
\begin{proof}
No component $X$ of $\partial M$ is a sphere: if it were so, it would bound a ball $B$ and we would have $M=B$, contradicting Corollary \ref{B:cor}. Hence $\chi(X)\leqslant 0$ and $X$ is incompressible
because $M$ is $\partial$-irreducible.
\end{proof}

We now prove that there are plenty of Haken manifolds. We start with a general proposition.

\begin{prop}
Let $M$ be an oriented, compact, irreducible, and $\partial$-irreducible 3-manifold with (possibly empty) boundary.
Every non-trivial homology class $\alpha \in H_2(M,\partial M;\matZ)$ is represented by a disjoint union of incompressible and $\partial$-incompressible oriented surfaces.
\end{prop}
\begin{proof}
Every class $\alpha$ is represented by a properly embedded oriented surface $S$ by Proposition \ref{homology:prop}. A compression as in Figure \ref{incompressible:fig} and \ref{incompressible2:fig} does not alter the homology class of the surface: indeed in homology we have $S'-S = \partial B$ where $B=D\times [-1,1]$ is a tubular neighbourhood of the compressing disc $D$. Hence $[S'] = [S] = \alpha$.

We compress $S$ until its connected components are either incompressible and $\partial$-incompressible surfaces, discs, or spheres. Since $M$ is irreducible and $\partial$-irreducible, discs and spheres bound balls and are hence homologically trivial, so they can be removed.
\end{proof}

\begin{cor} 
Let $M$ be oriented, compact, irreducible, and $\partial$-irreducible. If $H_2(M,\partial M; \matZ)\neq \{e\}$ then $M$ is Haken. 
\end{cor}
\begin{cor} \label{boundary:implies:Haken:cor}
Let $M$ be oriented, compact, irreducible, and $\partial$-irreducible. 
If $\partial M \neq \emptyset$ and $M\neq B$, then $M$ is Haken.
\end{cor}
\begin{proof}
If $\partial M$ contains a sphere, it bounds a ball $B$ and hence $M=B$. Otherwise $H_1(\partial M)$ has positive rank, and hence $H_2(M,\partial M)=H^1(M)$ also has positive rank by Corollary \ref{rank:cor}.
\end{proof}

We recall that every compact orientable 3-manifold decomposes along spheres and disc into irreducible and $\partial$-irreducible pieces. If one such piece has non-empty boundary then either it is ball, or it is Haken.

The following lemma will be useful soon. It says that every Haken manifold contains an interesting ``spanning'' surface $S$, that touches all the boundary components (the spanning surface need not to be connected).

\begin{lemma} \label{chiappa:tutto:lemma}
Every Haken manifold $M$ contains an oriented surface $S$, whose components are incompressible and $\partial$-incompressible, such that $[\partial S\cap X]\in H_1(X,\matZ)$ is non-trivial for every boundary component $X$ of $M$.
\end{lemma}
\begin{proof}
We have $\partial M = X_1\sqcup \ldots \sqcup X_k$ with $\chi(X_i)\leqslant 0$ for all $i$. Proposition \ref{lagrangian:prop} says that the image of 
$$\partial \colon H_2(M,\partial M, \matZ) \longrightarrow H_1(\partial M,\matZ)$$
is a lagrangian subgroup $L$ of maximal rank. We have 
$$H_1(\partial M, \matZ) = \oplus_{i=1}^k H_1(X_i, \matZ).$$
The decomposition is orthogonal with respect to the symplectic intersection form $\omega$. There is an $\alpha\in L$ whose projection to each $H_1(X_i,\matZ)$ is non-trivial: if not, the space $L$ would be $\omega$-orthogonal to some $H_1(X_i,\matZ)$ and hence contained in a smaller symplectic subspace, a contradiction since $L$ has maximal rank.
Pick any incompressible and $\partial$-incompressible surface $S$ such that $\partial [S] = \alpha$.
\end{proof}

\subsection{Normal surfaces}
On Haken manifolds, incompressible surfaces are efficiently detected by normal surfaces.  

\begin{prop} \label{normal:incompressible:prop}
Let $M$ be Haken and $T$ be a triangulation for $M$. Every compact surface $S\subset M$ whose components are all incompressible and $\partial$-incompressible is isotopic to a normal surface.
\end{prop}
\begin{proof}
We know that $S$ becomes normal after surgeries as in Figure \ref{incompressible:fig} or Figure \ref{incompressible2:fig}. These surgeries are actually isotopies since $S$ is incompressible and $\partial$-incompressible, and $M$ is irreducible and $\partial$-irreducible (see Figure \ref{incompressible3:fig}).
\end{proof}

Note that in the proof of Proposition \ref{normal:incompressible:prop}, as in many other proofs, it is crucial that $M$ be irreducible and $\partial$-irreducible.
We can now apply Lemma \ref{non:parallel:lemma} to get the following.

\begin{cor} \label{incompressible:cor}
Let $M$ be a Haken manifold. There is a $K>0$ such that every set $S$ of pairwise disjoint and non-parallel incompressible and $\partial$-incompressible surfaces in $M$ consists of at most $K$ elements.
\end{cor}

Our aim is now to cut a Haken manifold iteratively along incompressible and $\partial$-incompressible surfaces. The two-dimensional analogue to keep in mind is the following: every surface $S_g$ of genus $\geqslant 2$ can be cut into pairs-of-pants; a pair-of-pants is a quite simple surface, but we are still not satisfied and we further cut it along three arcs into two discs (two hexagons). We have constructed a two-step decomposition of $S_g$ into discs: this is what we would like to extend from two to three dimensions.

\subsection{Cutting along surfaces}
When we cut a 3-manifold along an incompressible surface, some nice properties of the manifold are preserved.

\begin{prop} \label{cut:surface:prop}
Let $M$ be compact and irreducible, and $S\subset M$ be either an essential disc or an incompressible surface. Let $M'$ be obtained by cutting $M$ along $S$. The following holds:
\begin{itemize}
\item the manifold $M'$ is irreducible;
\item a closed $\Sigma\subset M'$ is incompressible in $M'$ $\Longleftrightarrow$ it is so in $M$.
\end{itemize}
\end{prop}
\begin{proof}
We prove that $M'$ is irreducible. Let $\Sigma\subset M'$ be a sphere. Since $M$ is irreducible, the sphere $\Sigma$ bounds a ball $B\subset M$. The ball $B$ cannot contain $S$ because all surfaces in a ball are compressible. Therefore $B\subset M'$ and $M'$ is irreducible.

To prove the second assertion, we show that $\Sigma$ has a compressing disc $D$ in $M$ if and only if it has one in $M'$. If $D$ lies in $M'$ then of course it lies also in $M$. Conversely, suppose $D$ lies in $M$. Put $D$ in transverse position with respect to $S$ and pick an innermost intersection circle in $D$, bounding a disc $D'\subset D$. Since $D'$ cannot compress $S$, and since $M$ is irreducible, the disc $D'$ can be isotoped away from $S$. This simplifies $D\cap S$ and after finitely many steps we get $D\cap S = \emptyset$ and hence $D\subset M'$.
\end{proof}

\begin{cor} \label{Haken:cor}
If we cut a Haken 3-manifold along a closed incompressible surface, we get a disjoint union of Haken 3-manifolds.
\end{cor}

The following consequence is also interesting.

\begin{cor} \label{handlebody:closed:cor}
Let $M$ be compact with non-empty boundary. If $M$ is irreducible, then either it is a handlebody or it contains a closed incompressible surface. 
\end{cor}
\begin{proof}
The manifold $M$ decomposes along essential discs into manifolds $M_1, \ldots, M_k$ that are irreducible and $\partial$-irreducible. If every $M_i$ is a ball then $M$ is a handlebody. If $M_i$ is not a ball, the closed surface $S=\partial M_i$ is not a sphere and is incompressible in $M_i$, and hence also in $M$ by Proposition \ref{cut:surface:prop}.
\end{proof}

\subsection{Hierarchies}
We want to use incompressible surfaces to cut every Haken manifold into simpler pieces. The procedure goes as follows.

A \emph{hierarchy} for a Haken 3-manifold $M$ is a sequence of 3-manifolds\index{hierarchy}
$$M=M_0 \stackrel{S_0}\rightsquigarrow M_1 \stackrel{S_1}\rightsquigarrow M_2 \stackrel{S_2}\rightsquigarrow \ldots \stackrel{S_{h-1}}\rightsquigarrow M_h $$
where each $M_{i+1}$ is obtained cutting $M_i$ along a properly embedded (possibly disconnected) surface $S_i\subset M_i$, such that the following holds:
\begin{itemize}
\item every component of $S_i$ is an incompressible and $\partial$-incompressible surface or an essential disc, for all $i$;
\item the final manifold $M_h$ consists of balls.
\end{itemize}
The number $h$ is the \emph{height} of the hierarchy. We now show that every Haken manifold can be ``atomised'' into balls in three steps.
\begin{teo} \label{hierarchy:teo}
Every Haken manifold has a hierarchy of height $3$.
\end{teo}

\begin{proof}
Let $S_0$ be a maximal family of pairwise disjoint and non-parallel closed incompressible surfaces in $M$, which exists by Corollary \ref{incompressible:cor}. We cut $M_0=M$ along $S_0$ and get $M_1$.

Every connected component $M_1^i$ of $M_1$ is Haken by Corollary \ref{Haken:cor}. 
By Lemma \ref{chiappa:tutto:lemma} for every $i$ there is a ``spanning'' surface $S_1^i\subset M_1^i$ made of incompressible and $\partial$-incompressible components that intersects every boundary component of $M_1^i$. We cut $M_1$ along the spanning $S_1 = \sqcup S_1^i$ and get a new manifold $M_2$. 

We now prove that $M_2$ contains no closed incompressible surface. Indeed, if $\Sigma\subset M_2$ were closed and incompressible, then it would be so also in $M$ by Proposition \ref{cut:surface:prop}. Since $S_0$ is maximal, the surface $\Sigma$ would be parallel to a component of $S_0$, that is it would cobound a $\Sigma \times [0,1]$ with it. Since the spanning surface $S_1$ intersects all the boundary components of $M_1$, a component of $S_1$ would be contained in $\Sigma \times [0,1]$ and would intersect the boundary only on the side of $S_0$: this is excluded by Proposition \ref{line:bundles:prop} (products do not contain incompressible and $\partial$-incompressible surfaces with boundary only on one side).

Every component of $M_2$ is a handlebody by Corollary \ref{handlebody:closed:cor}. We cut it along a set $S_2$ of essential discs to get balls.
\end{proof}

\begin{figure}
\begin{center}
\includegraphics[width = 8 cm] {\iftoggle{BW}{hierarchy-BW}{hierarchy}}
\nota{The strata $S_0$, $S_1$, and $S_2$ intersect in vertices (left). The boundary of each ball $B\subset M_2$ is tessellated into domains belonging to $S_0, S_1$, or $S_2$. Three of them intersect at a 3-valent vertex (right).}
\label{hierarchy:fig}
\end{center}
\end{figure}

Hierarchies may be used to prove theorems on Haken manifolds. We will now use them to prove that incompressible surfaces must be $\pi_1$-injective: we need a preliminary discussion and a lemma.

In our hierarchy of height $3$, the surfaces in $S_0$ are closed, the spanning surfaces in $S_1$ have boundary, and $S_2$ consists of discs. To simplify notations, we redefine $S_0$ as $\partial M \cup S_0$. 

It is convenient to consider all the surfaces $S_0$, $S_1$, and $S_2$ inside $M$, without cutting $M$ along them.  With this perspective $\partial S_1$ is glued to $S_0$ and $\partial S_2$ is glued to $S_0\cup S_1$, via transverse maps. Every intersection $S_0\cap S_1 \cap S_2$ is a \emph{vertex} as in Figure \ref{hierarchy:fig}-(left). The space $X=S_0\cup S_1 \cup S_2$ is a two-dimensional cell complex whose complement in $M$ is a union of open balls: such an object $X$ is usually called a \emph{spine} for $M$.

We say that our hierarchy of height $3$ is \emph{minimal} if the essential discs in $S_2$ are chosen to minimise the total number of vertices. We can of course suppose that the hierarchy is minimal.

\begin{figure}
\begin{center}
\includegraphics[width = 12 cm] {\iftoggle{BW}{hierarchy2-BW}{hierarchy2}}
\nota{If $\gamma$ intersects at most three domains, it bounds a disc $D$ intersecting only these domains, as shown here.}
\label{hierarchy2:fig}
\end{center}
\end{figure}

The final manifold $M_3$ is the abstract closure of $M\setminus (S_0\cup S_1 \cup S_2)$ and consists of balls.
The boundary of every such ball $B$ is tessellated into domains belonging to $S_0$, $S_1$, or $S_2$. Three domains intersect at 3-valent vertices as in Figure \ref{hierarchy:fig}-(right). 
We say that $B$ is \emph{essential} if every simple closed curve $\gamma\subset \partial B$ transverse to the tessellation and intersecting $\leqslant 3$ domains is the boundary of a disc $D\subset \partial B$ intersecting only these domains: see Figure \ref{hierarchy2:fig}. 

\begin{lemma} \label{minimal:essential:lemma}
If the hierarchy is minimal, every ball in $M_3$ is essential.
\end{lemma}
\begin{proof}
Let $\gamma\subset \partial B$ intersect $k\leqslant 3$ domains. The curve $\gamma$ obviously bounds a properly embedded disc $D'\subset B$. 
If $k=1$, the curve is entirely contained in $S_i$ for some $i$. Since $S_0$ and $S_1$ are incompressible and $S_2$ consists of discs, the curve $\gamma$ bounds a disc $D$ also in $S_i$. Since $S_1$ is incompressible and $S_2$ consists of essential discs, the boundary components of $S_j$ with $j>i$ are not attached in the interior of $D$, so $D$ is entirely contained in the domain containing $\gamma$ and we are done.

If $k=2$, the curve $\gamma$ is contained in $S_i\cup S_j$ for some $i<j$. If $(i,j) = (0,1)$ we use that $S_1$ is $\partial$-incompressible to get a disc $D$ as in Figure \ref{hierarchy2:fig}-(top-right). If $j=2$, then $\gamma$ cuts a disc in $S_2$ into two discs $D_1 \cup D_2$. 

We have three half-discs $D', D_1, D_2$ intersecting in an arc, all contained in a handlebody $H\subset M_2$. Recall that $H\setminus S_2$ consists of balls. 
If we replace $D_1\cup D_2$ by either $D_1 \cup D'$ or $D_2 \cup D'$ (say, the first) we still get a disc system that cuts $H$ into balls, and hence another hierarchy. By minimality $D_2=D$ is adjacent to no vertices, \emph{i.e.}~it is as in Figure \ref{hierarchy2:fig}-(top-right) and we are done. 

The case $k=3$ is analogous: the curve $\gamma$ cuts a disc in $S_2$ into two parts $D_1\cup D_2$ and we may replace it with $D_1\cup D'$. By minimality $D_2$ is incident to at most one vertex $v$: we isotope $\gamma$ through $v$ as in Figure \ref{hierarchy3:fig}, and then conclude using the $k=2$ case.
\end{proof}

\begin{figure}
\begin{center}
\includegraphics[width = 10 cm] {\iftoggle{BW}{hierarchy3-BW}{hierarchy3}}
\nota{By minimality the region in $S_2$ is incident to a single vertex and we can slide $\gamma$ through it. Then we resume to the $k=2$ case.}
\label{hierarchy3:fig}
\end{center}
\end{figure}

\subsection{Dehn's Lemma}
Proposition \ref{pi1:prop} says that an orientable connected properly embedded $\pi_1$-injective surface $S\subset M$ with $\chi(S)\leqslant 0$ is incompressible. We can finally prove the converse, at least for closed surfaces.

The following result is often proved as a corollary of a famous topological fact called the \emph{Dehn Lemma}. We do not state Dehn's Lemma here, and we prove directly the following using hierarchies.\index{Dehn lemma}

\begin{teo} \label{incompressible:teo}
Let $M$ be a compact oriented 3-manifold. A connected oriented closed $S\subset M$ with $\chi(S)\leqslant 0$ is incompressible if and only if the induced map $i_*\colon \pi_1(S) \to \pi_1(M)$
is injective.
\end{teo}
\begin{proof}
We know one implication from Proposition \ref{pi1:prop}; here we suppose that $S$ incompressible and prove that $i_*\colon \pi_1(S) \to \pi_1(M)$ is injective. 

The decomposition of $M$ into irreducible and $\partial$-irreducible factors is made by cutting $M$ along 
essential spheres and discs transverse to $S$. Since $S$ is incompressible, we may surger $S$ along these spheres and discs without altering $i_*$, so that $S$ is disjoint from them and hence contained in a single factor. Therefore we may suppose $M$ is irreducible and $\partial$-irreducible. 

Now $M$ contains the incompressible $S$ and is hence Haken. Theorem \ref{hierarchy:teo} furnishes a hierarchy of height $\leqslant 3$. We may suppose that $S_0$ is a maximal system of closed incompressible surfaces containing $S$ and that the hierarchy is minimal. The balls in $M_3$ are essential by Lemma \ref{minimal:essential:lemma}.

Suppose by contradiction that $i_*$ is not injective: there is a loop $\gamma\colon S^1 \to S$ which is homotopically trivial in $M$ but not in $S$. The triviality in $M$ furnishes a continuous map $f\colon D^2 \to M$ which extends $\gamma$. We homotope $f$ to a smooth map, transverse to all the strata of $S_0\cup S_1\cup S_2$.

\begin{figure}
\begin{center}
\includegraphics[width = 5 cm] {\iftoggle{BW}{Dehn-BW}{Dehn}}
\nota{The counterimage of $S_0\cup S_1 \cup S_2$ along $f$. An edge with label $i$ goes to $S_i$ (unlabelled edges go to $S_0$). Here there is a $0$-gon (\iftoggle{BW}{light grey}{yellow}) and three $2$-gons (\iftoggle{BW}{grey}{green}). The correct notion of $k$-gon should be clear from the picture.}
\label{Dehn:fig}
\end{center}
\end{figure}

By transversality, the counterimage $f^{-1}(S_0\cup S_1 \cup S_2)$ is a graph in $D^2$ as in Figure \ref{Dehn:fig}, which divides $D^2$ into regions. The graph is itself a hierarchy, with edges of type $0,1,2$ attached iteratively.
An easy Euler characteristic argument shows that at least one region $R$ is a $k$-gon with $k\leqslant 3$, see Figure \ref{Dehn:fig}. The region $R$ is mapped inside a ball $B\subset M_3$, and $\partial R$ is mapped to an immersed curve $\alpha\subset\partial B$ intersecting $k\leqslant 3$ domains of the tessellated $\partial B$. 

\begin{figure}
\begin{center}
\includegraphics[width = 11 cm] {\iftoggle{BW}{Dehn2-BW}{Dehn2}}
\nota{The immersed curve $\alpha$ intersects $k\leqslant 3$ regions and is hence of one of these types (because $B$ is essential).}
\label{Dehn2:fig}
\end{center}
\end{figure}

Since $B$ is essential, the curve $\alpha$ is of one of the types shown in Figure \ref{Dehn2:fig}. In all cases we may slide the disc $f(R)$ away from $B$ and decrease the number of regions in $D$ by destroying $R$. After finitely many homotopies of this kind we get $f(D)\cap (S_0 \cup S_1 \cup S_2) = \emptyset$ and hence $f(D)$ is entirely contained in a ball $B\subset M_3$. Therefore $\gamma$ is trivial in $S$, a contradiction.
\end{proof}

\begin{cor}
The fundamental group of a Haken manifold is infinite.
\end{cor}
\begin{proof}
It contains the fundamental group of a closed surface with $\chi\leqslant 0$, which is infinite.
\end{proof}
\begin{cor}
Elliptic 3-manifolds are not Haken.
\end{cor}

We will see in the subsequent chapters that every flat 3-manifold is Haken. Hyperbolic 3-manifolds may or may not be Haken.

\subsection{Essential surfaces} \label{essential:subsection}
Topologists sometimes use the term ``essential'' to summarise various reasonable notions in a single word. We already know what an essential disc or sphere is (see Section \ref{essential:discs:subsection}) and we now turn to surfaces with non-positive Euler characteristic. 

Let $M$ be a compact oriented three-manifold and $S\subset M$ be a properly embedded connected compact surface with $\chi(S)\leqslant 0$. We say that $S$ is \emph{essential} if it is incompressible, $\partial$-incompressible, and not $\partial$-parallel.\index{surface!essential surface}

\subsection{Simple manifolds} \label{simple:manifold:subsection}
Let $M$ be irreducible and $\partial$-irreducible. We introduce yet some more definitions. We say that
\begin{itemize}
\item $M$ is \emph{atoroidal} if it does not contain essential tori,
\item $M$ is \emph{acylindrical} if it does not contain essential annuli.
\end{itemize}

Finally, the manifold $M$ is \emph{simple} if it is atoroidal and acylindrical. We can summarise this definition as follows:\index{three-manifold!atoroidal and acylindrical three-manifold}\index{three-manifold!simple three-manifold}

\begin{defn} A compact oriented 3-manifold $M$ with (possibly empty) boundary is \emph{simple} if it does not contain any essential sphere, disc, torus, or annulus.
\end{defn}

Many examples come from elliptic and hyperbolic geometry:

\begin{prop}
Every closed elliptic or hyperbolic 3-manifold $M$ is simple.
\end{prop}
\begin{proof}
We know that $M$ is irreducible. The manifold $M$ does not contain incompressible tori because $\pi_1(M)$ does not contain $\matZ \times \matZ$: if $M$ is elliptic then $\pi_1(M)$ is finite, if it is hyperbolic we use Corollary \ref{no:ZxZ:cor}.
\end{proof}

The flat geometry is an exception: the three-torus $S^1\times S^1\times S^1$ contains many incompressible tori and is hence not simple. For instance, the two-torus $S^1\times S^1 \times p$ is incompressible (because it is $\pi_1$-injective).

Our next goal will be to decompose every irreducible and $\partial$-irreducible manifold $M$ along some canonical set of essential tori and annuli into some pieces. These pieces will be either simple, or belong to a particular class: the \emph{Seifert manifolds}. We introduce this class in the next chapter. 

\subsection{References}
The main sources that we have used for this chapter are an unfinished book of Hatcher \cite{H} and Fomenko--Matveev \cite{FoMa}. The proof of Theorem \ref{incompressible:teo} through hierarchies is due to Aitchison and Rubinstein \cite{AR}.

%% file: Seifert.tex
\chapter{Seifert manifolds} \label{Seifert:chapter}

In the previous chapter we have proved various general theorems on three-manifolds, and it is now time to construct examples. A rich and important source is a family of manifolds built by Seifert in the 1930s, which generalises circle bundles over surfaces by admitting some ``singular'' fibres. The three-manifolds that admit such kind of fibration are now called \emph{Seifert manifolds}. 

In this chapter we introduce and completely classify (up to diffeomorphisms) the Seifert manifolds. In Chapter \ref{eight:chapter} we will then show how to \emph{geometrise} them, by assigning a nice Riemannian metric to each. We will show, for instance, that all the elliptic and flat three-manifolds are in fact particular kinds of Seifert manifolds.

\section{Lens spaces}
We introduce some of the simplest 3-manifolds, the lens spaces. These manifolds (and many more) are easily described using an important three-dimensional construction, called \emph{Dehn filling}.\index{Dehn filling}

\subsection{Dehn filling} \label{Dehn:filling:subsection}
If a 3-manifold $M$ has a spherical boundary component, we can cap it off with a ball. If $M$ has a toric boundary component, there is no canonical way to cap it off: the simplest object that we can attach to it is a solid torus $D\times S^1$, but the resulting manifold depends on the gluing map. This operation is called a \emph{Dehn filling} and we now study it in detail.

Let $M$ be a 3-manifold and $T\subset \partial M$ be a boundary torus component. 
\begin{defn}
A \emph{Dehn filling} of $M$ along $T$ is the operation of gluing a solid torus $D\times S^1$ to $M$ via a diffeomorphism $\varphi\colon \partial D\times S^1 \to T$. 
\end{defn}

The closed curve $\partial D \times \{x\}$ is glued to some simple closed curve $\gamma\subset T$, see Figure \ref{DF:fig}. The result of this operation is a new manifold $M^{\rm fill}$, which has one boundary component less than $M$.

\begin{figure}
\begin{center}
\includegraphics[width = 9 cm] {\iftoggle{BW}{DF-BW}{DF}}
\nota{The Dehn filling $M^{\rm fill}$ of a 3-manifold $M$ is determined by the unoriented simple closed curve $\gamma \subset T$ to which a meridian $\partial D$ of the solid torus is attached.}
\label{DF:fig}
\end{center}
\end{figure}

\begin{lemma} \label{pq:lemma}
The manifold $M^{\rm fill}$ depends only on the isotopy class of the unoriented curve $\gamma$.
\end{lemma}
\begin{proof}
Decompose $S^1$ into two closed segments $S^1 = I\cup J$ with coinciding endpoints. The attaching of $D\times S^1$ may be seen as the attaching of a 2-handle $D\times I$ along $\partial D \times I$, followed by the attaching of a 3-handle $D\times J$ along its full boundary. 

If we change $\gamma$ by an isotopy, the attaching map of the 2-handle changes by an isotopy and hence gives the same manifold. The attaching map of the 3-handle is irrelevant by Proposition \ref{cap:prop}.
\end{proof}

We say that the Dehn filling \emph{kills} the curve $\gamma$, since this is what really happens on fundamental groups, as we now see. 

The \emph{normaliser} of an element $g\in G$ in a group $G$ is the smallest normal subgroup $N(g) \triangleleft G$ containing $g$. The normaliser depends only on the conjugacy class of $g^{\pm 1}$, hence  the subgroup $N(\gamma)\triangleleft \pi_1(M)$ makes sense without fixing a basepoint or an orientation for $\gamma$.\index{normaliser}

\begin{prop} 
We have
$$\pi_1(M^{\rm fill}) = \pi_1(M) /_{N(\gamma)}.$$
\end{prop}
\begin{proof}
The Dehn filling decomposes into the attachment of a 2-handle over $\gamma$ and of a 3-handle. By Van Kampen, the first operation kills $N(\gamma)$,
and the second leaves the fundamental group unaffected.
\end{proof}

Let a \emph{slope} on a torus $T$ be the isotopy class $\gamma$ of an unoriented homotopically non-trivial simple closed curve.\index{slope} The set of slopes on $T$ was indicated by $\calS$ in Chapter \ref{Teichmuller:chapter}. If we fix a basis $(m, l)$ for $H_1(T,\matZ)=\pi_1(T)$, every slope may be written as $\gamma = \pm (pm + ql)$ for some coprime pair $(p,q)$. Therefore we get a 1-1 correspondence
$$\calS \longleftrightarrow \matQ\cup \{\infty\}$$
by sending $\gamma$ to $\frac pq$. If $T$ is a boundary component of $M$, every number $\frac pq$ determines a Dehn filling of $M$ that kills the corresponding slope $\gamma$. 

Different values of $\frac pq$ typically produce non-diffeomorphic manifolds $M^{\rm fill}$: this is not always true - a notable exception is described in the next section - but it holds in ``generic'' cases.

\subsection{Lens spaces}
The simplest manifold that can be Dehn-filled is the solid torus $M=D\times S^1$ itself. The oriented \emph{meridian} $m=S^1 \times \{y\}$ and \emph{longitude} $l=\{x\}\times S^1$ form a basis for $H_1(\partial M,\matZ)$.\index{meridian and longitude}

\begin{defn} The \emph{lens space} $L(p,q)$ is the result of a Dehn filling of $M=D\times S^1$ that kills the slope $qm+pl$.\index{lens space}
\end{defn}

A lens space is a three-manifold that decomposes into two solid tori. We have already encountered lens spaces in the more geometric setting of Section \ref{lens:subsection}, and we will soon prove that the two definitions are coherent. Since $L(p,q)=L(-p,-q)$ we usually suppose $p\geqslant 0$.

\begin{ex}
We have $\pi_1\big(L(p,q)\big) = \matZ/_{p\matZ}$. 
\end{ex}

\begin{prop}
We have $L(0,1) = S^2\times S^1$ and $L(1,0) = S^3$.
\end{prop}
\begin{proof}
The lens space $L(0,1)$ is obtained by killing $m$, that is by mirroring $D\times S^1$ along its boundary. The lens space $L(1,0)$ is $S^3$ because the complement of a standard solid torus in $S^3$ is another solid torus, with the roles of $m$ and $l$ exchanged (exercise).
\end{proof}

\begin{ex} \label{fill:T:ex}
Every Dehn filling of one component of the product $T\times [0,1]$ is diffeomorphic to $D\times S^1$. Therefore by Dehn-filling both components of $T\times [0,1]$ we get a lens space.
\end{ex}

The solid torus $D\times S^1$ has a non-trivial self-diffeomorphism 
$$(x, e^{i\theta}) \longmapsto (xe^{i\theta}, e^{i\theta})$$
called a \emph{twist} along the disc $D\times \{y\}$. The solid torus can also be \emph{mirrored} via the map
$$(x,e^{i\theta}) \longmapsto (x,e^{-i\theta}).$$

\begin{ex} \label{pqpq:ex}
We have $L(p,q)\isom L(p,q')$ if $q' \equiv \pm q^{\pm 1}$ (mod $p$).
\end{ex}
\begin{proof}[Hint]
Twist, mirror, and exchange the two solid tori giving $L(p,q)$.
\end{proof}

\begin{rem}
The meridian $m$ of the solid torus $M=D\times S^1$ may be defined intrinsically as the unique slope in $\partial M$ that is homotopically trivial in $M$. The longitude $l$ is \emph{not} intrinsically determined: a twist sends $l$ to $m+l$. The solid torus contains infinitely many non-isotopic longitudes, and there is no intrinsic way to choose one of them.
\end{rem}

\subsection{Equivalence of the two definitions}
When $p>0$, we have defined the lens space $L(p,q)$ in two different ways: as the $(q,p)$-Dehn filling of the solid torus, and as an elliptic manifold in Section \ref{lens:subsection}. In the latter description we set
$$\omega = e^{\frac{2\pi i}p}, \qquad f(z,w) = (\omega z, \omega^qw)$$ 
and define $L(p,q)$ as $S^3/_\Gamma$ where $\Gamma = \langle f \rangle$ is generated by $f$. We now show that the two definitions produce the same manifolds.

\begin{prop}
The manifold $S^3/_{\langle f \rangle}$ is the $(q,p)$-Dehn filling of the solid torus.
\end{prop}
\begin{proof}
The isometry $f$ preserves the central torus 
$$T = \big\{(z,w)\ \big|\ |z|=|w| = \tfrac {\sqrt 2}2 \big\}$$
that divides $S^3$ into two solid tori
\begin{align*}
N^1 & = \big\{(z,w)\ \big|\ |z| \leqslant \tfrac {\sqrt 2}2,\ |w| = \sqrt{1-|z|^2}\big\}, \\
N^2 & = \big\{(z,w)\ \big|\ |w| \leqslant \tfrac {\sqrt 2}2,\ |z| = \sqrt{1-|w|^2}\big\}.
\end{align*}
Identify $T$ with $S^1\times S^1 = \matR^2/_{\matZ^2}$ in the obvious way, so that  $H_1(T) = \matZ \times \matZ$. The meridians of $N^1$ and $N^2$ are $(1,0)$ and $(0,1)$. The isometry $f$ act on $T$ as a translation of vector $v=\big(\frac 1p, \frac qp\big)$. The quotient $T/_{\langle f\rangle}$ is again a torus, with fundamental domain the parallelogram generated by $v$ and $w=(0,1)$.

The quotients $N^1/_{\langle f \rangle}$, and $N^2/_{\langle f \rangle}$ are again solid tori. Therefore $S^3/_{\langle f \rangle}$ is also a union of two solid tori. Their meridians are the projections of the horizontal and vertical lines in $\matR^2$ to $T/_{\langle f \rangle} = \matR^2/_{\langle v,w \rangle}$. In the basis $(v,w)$ these meridians are $pv -qw$ and $w$ respectively. Therefore $S^3/_{\langle f\rangle}$ is a $(-q,p)$-Dehn filling on the solid torus, which is diffeomorphic to the $(q,p)$-Dehn filling by mirroring the solid torus.
\end{proof}

\begin{cor}
We have $L(1,0) = S^3$ and $L(2,1) = \matRP^3$.
\end{cor}
\begin{proof}
We have $f = \id$ and $f=-\id$, correspondingly.
\end{proof}

\subsection{Classification of lens spaces}
Which lens spaces are diffeomorphic? It is not so easy to answer this question, because many lens spaces like $L(5,1)$ and $L(5,2)$ have the same homotopy and homology groups, while there is no evident diffeomorphism between them. 
A complete answer was given by Reidemeister in 1935, who could distinguish lens spaces using a new invariant, now known as the \emph{Reidemeister torsion}.\index{Reidemeister torsion} More topological proofs were discovered in th 1980s by Bonahon and Hodgson. We follow here Hatcher \cite{H}.

\begin{teo} \label{Lpq:teo}
The lens spaces $L(p,q)$ and $L(p',q')$ are diffeomorphic $\Longleftrightarrow$ $p=p'$ and $q' \equiv \pm q^{\pm 1}$ (mod $p$).
\end{teo}
\begin{proof}
We may suppose $p>1$.
Exercise \ref{pqpq:ex} furnishes one implication, so we start with a lens space $L(p,q)$ and we prove that there is (up to isotopy) only one torus $T$ dividing $L(p,q)$ into two solid tori: this fact then implies that one can recover $p$ and $q$ intrinsically from $L(p,q)$, up to the stated ambiguity for $q$ (which depends on the chosen longitudes and orientations, and changes by switching the solid tori). So we suppose that there is another torus $T'$.

Let $\Sigma$ be the core circle of one solid torus adjacent to $T'$. If we can isotope $\Sigma$ inside the torus $T$, we are done: in that case $T'$ is isotopic to the boundary of a small tubular neighbourhood of $\Sigma$, hence both $T$ and $T'$ are cut into two annuli $T = A_1 \cup A_2$ and $T' = A_1' \cup A_2'$ such that $\Sigma \subset A_1$, all four annuli share the same boundaries, and $A_2$ is contained in the large outside solid torus bounded by $T'$; the annulus $A_2$ is incompressible in this solid torus by Proposition \ref{pi1:prop} (otherwise $\Sigma$ would be homotopically trivial, contradicting $p>1$) and hence $\partial$-parallel by Proposition \ref{incompressible:solid:torus:prop}, so it is isotopic to either $A_1'$ or $A_2'$, suppose to $A_2'$; since $A_1$ is clearly isotopic to $A_1'$ we conclude that $T$ and $T'$ are isotopic.

Our aim is now to prove that $\Sigma$ can be isotoped inside $T$. To this purpose we construct two different objects from $T'$ and $T$, a \emph{spine} $\Delta$ and a \emph{foliation} $\calF$. The spine $\Delta\subset M$ is built by adding to $\Sigma$ the meridian $D$ of the other solid torus incident to $T'$, enlarged so that $\partial D$ is glued along $\Sigma$ like a degree-$p$ covering. Note that $M\setminus \Delta$ is an open ball.

We construct the foliation $\calF$ of $M$ by subdividing each solid torus bounded by $T$ into concentric tori, with a central singular circle in each. We represent $\calF$ as the level sets of a map $f\colon L(p,q) \to [0,1]$ where the singular circles are the extreme levels $f^{-1}(0)$ and $f^{-1}(1)$.

We now put $\Delta$ in some good position with respect to $\calF$. We first perturb $\Sigma$ so that it is disjoint from the singular circles of $\calF$ and $f|_\Sigma$ is a Morse function with singular points in distinct levels. At every local maximum (minimum) of $f|_\Sigma$ we isotope $\Delta$ so that the $p$ sheets of $\Delta$ lie above (below) $\Sigma$ as in Figure \ref{sheets_above:fig}. Finally, we require $\Delta$ to be transverse to the two singular circles of $\calF$ and $f|_{\Delta \setminus \Sigma}$ to be a Morse function, with critical points at distinct levels (also distinct from those of the critical points of $f|_\Sigma$).

\begin{figure}
\begin{center}
\includegraphics[width = 3 cm] {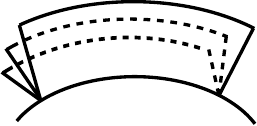}
\nota{We require the $p$ sheets of the spine $\Delta$ to lie above every local maximum for $f|_\Sigma$.}
\label{sheets_above:fig}
\end{center}
\end{figure}

\begin{figure}
\begin{center}
\includegraphics[width = 12.5 cm] {\iftoggle{BW}{lens_singular-BW}{lens_singular}}
\nota{The singular points in $D$ can be centres (top left), saddles (top right) and half-saddles (bottom).}
\label{lens_singular:fig}
\end{center}
\end{figure}

The foliation $\calF$ induces a singular foliation on $\Delta$ that pulls back to $D$ via the map $D\to \Delta$. There are three types of singular points in the foliation, shown in Figure \ref{lens_singular:fig}: centres, saddles, and half-saddles.

\begin{figure}
\begin{center}
\includegraphics[width = 12.5 cm] {\iftoggle{BW}{lens_leaves-BW}{lens_leaves}}
\nota{A singular leaf exiting from a singular point may be of one of these kinds. The dots in (e) indicate that there might be another half-saddle at the other endpoint: this is because there are half-saddles at the same level (precisely $p$ for every local maximum or minimum for $f|_\Sigma$).}
\label{lens_leaves:fig}
\end{center}
\end{figure}

A leaf in $D$ incident to a singular point is \emph{singular}. Since we minimised the critical points that may stay at the same level, the singular leaves exiting from a singular point may be only of the six possible kinds shown in Figure \ref{lens_leaves:fig}.  
The saddles and half-saddles of type (a) and (e) are called \emph{essential}, and the others \emph{inessential}.
In (a) the singular leaves divide $D$ into four \emph{quarter discs}, and in (e) they cut off two or more \emph{half discs}.

Among all possible good configurations of $\Delta$, we pick one that minimises first the number of critical points in $f|_\Sigma$, and second the number of essential saddles.

Let $D' \subset D$ be a quarter or half-disc not containing any smaller quarter or half-disc. Define $\alpha = D' \cap \partial D$. Suppose first that $D'$ is a half-disc and both endpoints of $\alpha$ are singular. The singular disc $D'$ can contain some singular point of type (b), and it certainly contains one singular point of type (f), see Figure \ref{lens_leaves1:fig}-(left). (It contains no other singular point except these.) The image of $D'$ in $\Delta$ is shown in Figure \ref{lens_leaves1:fig}-(right). In this case we can slide $\alpha$ along $D'$ to the \iftoggle{BW}{light grey}{green} arc shown in the figure which lies entirely in a level torus of $f$. If we could do this for all the arcs of $\partial D$ cut by singular points, we would happily isotope $\Sigma$ inside a torus level and hence into $T$, and we would be done.

There are however two other cases to consider, and both will be excluded by our minimality assumption on $\Delta$. One is that $D'$ may be a half-disc with only one singular endpoint for $\alpha$, as in Figure \ref{lens_leaves2:fig}. In that case we can isotope $\alpha$ to the \iftoggle{BW}{light grey}{green} curve through $D'$, dragging the spine $\Delta$ behind. In the new configuration $f|_\Sigma$ has strictly less critical points, a contradiction.

\begin{figure}
\begin{center}
\includegraphics[width = 8.5 cm] {\iftoggle{BW}{lens_leaves1-BW}{lens_leaves1}}
\nota{The disc $D'$ is a half-disc and $\alpha$ is incident to two singular points.}
\label{lens_leaves1:fig}
\end{center}
\end{figure}

\begin{figure}
\begin{center}
\includegraphics[width = 8.5 cm] {\iftoggle{BW}{lens_leaves2-BW}{lens_leaves2}}
\nota{The disc $D'$ is a half-disc and $\alpha$ is incident to one singular point.}
\label{lens_leaves2:fig}
\end{center}
\end{figure}

\begin{figure}
\begin{center}
\includegraphics[width = 9 cm] {\iftoggle{BW}{lens_leaves3-BW}{lens_leaves3}}
\nota{The disc $D'$ is a quarter disc.}
\label{lens_leaves3:fig}
\end{center}
\end{figure}

In the last case $D'$ is a quarter disc as in Figure \ref{lens_leaves3:fig} and we can isotope $\alpha$ to the \iftoggle{BW}{light grey}{green} curve dragging $\Sigma$ as above, but more carefully: we enlarge the other $p-1$ sheets of $\Delta$ by parallel copies of $D'$. In the new configuration the number of critical points in $f|_\Sigma$ is unchanged, but there are strictly fewer essential saddles: at least one is destroyed and no new one is created. (Note that many inessential saddles may be created.) This is also excluded.
\end{proof}

We have discovered, in particular, that there are closed three-manifolds like $L(5,1)$ and $L(5,2)$ that are not diffeomorphic, although their fundamental groups are both isomorphic to $\matZ/_{5\matZ}$; the two manifolds are both covered by $S^3$, so they also have isomorphic higher homotopy groups.

\section{Circle bundles}
We now introduce another simple class of 3-manifolds, the orientable circle bundles over some compact surface $S$. We will discover that there is essentially only one circle bundle if $S$ has boundary, and infinitely many if $S$ is closed, distinguished by an integer called the \emph{Euler number}.

\subsection{The trivial circle bundle}
Let $S$ be a compact connected surface. As every connected manifold, the surface $S$ has a unique orientable line bundle
$$S\times I \quad {\rm or} \quad S\timtil I$$ 
depending on whether $S$ is orientable or not. We denote by 
$$M=S\times S^1 \quad {\rm or} \quad S\timtil S^1$$ 
respectively the double of $S\times I$ and $S\timtil I$ along its boundary. If we do not know whether $S$ is orientable or not, we use the symbols $S\timforsetil I$ and $S\timforsetil S^1$ to denote these objects. The manifold $S\timforsetil S^1$ is an orientable circle bundle over $S$, called the \emph{trivial} one. 

\subsection{Circle bundles with boundary}
We start by exploring the case where the base surface $S$ has non-empty boundary: in this case every bundle $M$ over $S$ is a 3-manifold with boundary; the boundary consists of tori, one fibering above each circle in $\partial S$, because the torus is the unique orientable surface that fibres over $S^1$.

It turns out that there is essentially only one bundle over $S$, the trivial one:

\begin{lemma} \label{unique:lemma}
If $\partial S \neq \emptyset$, the orientable circle bundles on $S$ are all isomorphic.
\end{lemma}
\begin{proof}
Let $N\to S$ be an orientable circle bundle. Decompose $S$ as a disc $D$ with some pairs of disjoint segments $(I_i, J_i)$ in $\partial D$ to be glued. Since $D$ is contractible the restriction of $N$ to $D$ is a product $D\times S^1$ and $N$ is obtained from it by gluing the annuli $I_i\times S^1$ and $J_i\times S^1$ via orientation-reversing fibre-preserving maps. Two such maps are always isotopic (exercise) and hence $N$ is uniquely determined.
\end{proof}

\begin{figure}
\begin{center}
\includegraphics[width = 11 cm] {\iftoggle{BW}{section-BW}{section}}
\nota{A twist along a fibered annulus $A$ modifies a section $S$ as shown.}
\label{section:fig}
\end{center}
\end{figure}

We now want to study the \emph{sections} of the trivial bundle $M\to S$, because these will be useful in the study of bundles over closed surfaces. We now discover that, although the bundle is trivial, it contains many non isotopic sections, and we want to classify them.

Recall that a \emph{section} of the bundle $\pi\colon M\to S$ is a map $i\colon S\to M$ such that $\pi\circ i = \id$. Since the section $i$ is determined by its image $i(S)$, we simply consider the surface $i(S)$ as a section of $\pi$. By construction $M$ is the double of an interval bundle over $S$ and as such it contains the \emph{zero-section} $S$ there. However, this section is not unique in general, not even up to isotopy. 

To modify a section, pick a properly embedded arc in $S$. The arc determines a fibered annulus $A\subset M$ above it, which we may use to \emph{twist} the section as shown in Figure \ref{section:fig}. This operation modifies the curves $\partial S\subset \partial M$ via two Dehn twists (one positive and one negative) on the tori $\partial M$ along the two curves in $\partial A$.

By twisting along annuli we may construct all the sections of $M$:

\begin{lemma} \label{section:lemma}
Two sections of $S\timforsetil S^1$ are connected by a composition of twists along fibered annuli and fibre-preserving isotopies.
\end{lemma}
\begin{proof}
Let $i_1, i_2$ be two sections. Decompose $S$ as a 0-handle $D$ with some rectangular 1-handles attached to $\partial D$. The bundle is trivial on each handle. Since $D$ is contractible, every map $i\colon D\to S^1$ is homotopic to a constant and we may hence isotope $i_1$ and $i_2$ so that they coincide on $D$. 

See a 1-handle as $[-1,1]^2$ with $\{\pm 1\}\times [-1,1]$ glued to $D$. The sections $i_1, i_2\colon [-1,1]^2 \to S^1$ coincide on $\{\pm 1\}\times [-1,1]$, and up to reparametrising we may suppose that they are constantly $1$ there. The section $i_1 \colon [-1,1]^2 \to S^1$ defines an element $d\in \pi_1(S^1,1) = \matZ$, determined by the loop $t\mapsto i_1(t,u)$, whose homotopy class does not depend on the fixed $u\in [-1,1]$. The section is in fact determined up to isotopy by $d$, and by twisting along the annulus $0 \times [-1,1]$ we can vary this integer arbitrarily. Therefore after some twists and isotopies the sections $i_0$ and $i_1$ are both constantly 1 on 1-handles, and thus coincide everywhere.
\end{proof}

\begin{cor} \label{unique:cor}
If $S$ has only one boundary component, the boundary of a section of $M=S\timforsetil S^1$ is a slope in $\partial M$ that does not depend on the section.
\end{cor}
\begin{proof}
Distinct sections are connected by finitely many Dehn twists along annuli. One such twist acts on the torus $\partial M$ as a composition of two opposite Dehn twists, which cancel each other. Hence it does not affect the boundary slope of a section.
\end{proof}

\begin{ex} Prove this corollary using Proposition \ref{lagrangian:prop}.
\end{ex}

\subsection{Closed circle bundles} \label{closed:circle:subsection}
We now turn to closed circle bundles. In this section we prove that 
the oriented circle bundles over a closed surface are parameterised by an integer called the \emph{Euler number} of the fibering.\index{Euler number}

We prefer to see the bundles over closed surfaces as Dehn fillings of bundles over surfaces with boundary. Here are the details.

Let $S$ be a compact surface with non-empty boundary. Pick $M=S\timforsetil S^1$ and fix an orientation for $M$. Recall that we denote by $S$ the zero-section of $M$. Every boundary component $T$ of $M$ is an oriented torus, which contains two natural unoriented simple closed curves: the boundary $m=T\cap \partial S$ of the section $S$ and the fibre $l$ of the bundle. If oriented, the curves $m$ and $l$ form a basis $(m,l)$ for $H_1(T,\matZ)$. We choose orientations for $m$ and $l$ such that $(m,l)$ form a positively oriented basis: there is a unique choice up to reversing both $m$ and $l$.

A Dehn filling on $T$ is determined as usual by a pair $(p,q)$ of coprime integers that indicate the slope $\pm (pm+ql)$ to be killed. 

Suppose now that $S$ has only one boundary component and let $M^{\rm fill}$ be obtained by Dehn filling $M$ along the slope $(1,q)$. Let $\hat S$ be the closed surface obtained by capping $S$ with a disc.

\begin{prop} \label{classify:bundles:prop}
The circle bundle $M\to S$ extends to a circle bundle $M^{\rm fill}\to \hat S$. Every oriented circle bundle on $\hat S$ is obtained in this way, and distinct values of $q$ yield vector bundles that are not orientation-preservingly isomorphic.
\end{prop}
\begin{proof}
The meridian of the filling solid torus is $m'=m + ql$. The fibre $l$ has geometric intersection 1 with $m'$ and is hence a longitude for the filling solid torus. We may represent the filling solid torus as $D \times S^1$ with $m' = S^1 \times \{y\}$ and $l=\{x\}\times S^1$. The circle bundle $M\to S$ extends naturally to a circle bundle $M^{\rm fill} \to \hat S$ with $\hat S = S\cup D$.

Every closed circle bundle $N\to \hat S$ arises in this way: the bundle above a disc $D\subset \hat S$ is the trivial $D\times S^1$, and if we remove it we get $M \to S$ back. 
The number $q$ is intrinsically determined: the meridian $m$ does not depend on the section of $M\to S$ by Corollary \ref{unique:cor}, and the equality $m' = m+ql$ determines $q$. Therefore distinct values of $q$ yield non-isomorphic bundles.
\end{proof}

The integer $q$ is the \emph{Euler number} of the circle bundle and is usually denoted with the letter $e$. We summarise our discovery:
\begin{cor} \label{unique:e:cor}
For every $e\in \matZ$ and every closed surface $S$ there is a unique oriented circle bundle over $S$ with Euler number $e$.
\end{cor}

A change of orientation for $M$ transforms $e$ into $-e$.
Recall that every closed surface $S$ has a trivial circle bundle $S\timforsetil S^1$ constructed by doubling the unique oriented line bundle on $S$.

\begin{ex} \label{circle:bundle:ex}
An oriented circle bundle over a closed surface is trivial $\Longleftrightarrow$ $e=0$ $\Longleftrightarrow$ the bundle has a section. 
\end{ex}

We may see the Euler number of a bundle $M\to S$ over a closed $S$ as an obstruction for the existence of a section.

\begin{oss} \label{bundle:oss}
Every oriented $n$-dimensional vector bundle $E \to S$ over a closed oriented $n$-manifold $S$ has a \emph{Euler number} defined by taking two generic sections and counting their signed intersections. We briefly explain how this number is closely related to the one we defined here.\index{Euler number}

Each vector bundle $E \to S$ induces a sphere bundle $M\to S$: it suffices to fix a Riemannian metric on $E$ and take the sub-bundle consisting of unit tangent vectors. When $n=2$ we get a circle bundle $M\to S$ and the Euler number of $E\to S$ coincides with that of $M\to S$ that we defined above.

When $E$ is the tangent bundle of $S$, the Euler number is the Euler characteristic $\chi(S)$. For instance, the unit tangent bundle of $S^2$ has Euler number $e=\chi(S^2)=2$ and hence it is diffeomorphic to $L(2,1) = \matRP^3$. 
\end{oss}

\begin{ex} \label{bundle:ex}
Let $M$ be a circle bundle over the genus-$g$ surface $S_g$ with Euler number $e$. We have 
$H_1(M,\matZ) = \matZ^{2g} \times \matZ/_{e\matZ}$.
\end{ex}

\begin{cor} \label{distinct:bundle:cor}
Let $M \to S_g$ and $M' \to S_{g'}$ be circle bundles with Euler numbers $e$ and $e'$. The manifolds $M$ and $M'$ are diffeomorphic $\Longleftrightarrow$ $g=g'$ and $|e|=|e'|$.
\end{cor}

\begin{ex} \label{e:ex}
The circle bundle $M$ over $S^2$ with Euler number $e$ is diffeomorphic to the lens space $L(|e|,1)$.
\end{ex}
\begin{proof}[Hint] The base sphere $S^2$ decomposes into two discs, and the fibration over each disc is a solid torus. Therefore $M$ is the union of two solid tori.
\end{proof}

We end the discussion by defining explicitly a famous non-trivial circle bundle over $S^2$.

\subsection{The Hopf fibration} \label{Hopf:lens:subsection}
The quotient map 
\begin{align*}
(\matC^2)^* & \longrightarrow \matCP^1 \\
(z,w) & \longmapsto [z,w]
\end{align*}
restricts to a circle bundle $S^3 \to \matCP^1 = S^2$ called the \emph{Hopf fibration}.
The fibre over $[z,w]$ is the circle $(\omega z, \omega w)$ as $\omega \in S^1$ varies. See Figure \ref{Hopf:fig}.\index{Hopf fibration}

\begin{figure}
\begin{center}
\includegraphics[width = 9 cm] {\iftoggle{BW}{Hopf_Fibration-BW}{Hopf_Fibration}}
\nota{The fibre of every point in the Hopf fibration $S^3\to S^2$ is a circle, and the counterimage of a circle in $S^2$ is a torus in $S^3$. The picture shows the counterimage of three circle arcs: each is a portion of standard torus in $S^3 = \matR^3\cup \infty$, foliated by circles of type $(1,1)$.}
\label{Hopf:fig}
\end{center}
\end{figure}

Exercise \ref{e:ex} implies that the Euler number of the Hopf fibration $S^3 \to S^2$ is $\pm 1$, and hence it has no section by Exercise \ref{circle:bundle:ex}.

\section{Seifert manifolds}
We now enlarge the class of circle bundles over surfaces by admitting some kind of singular fibres. We introduce the \emph{Seifert fibrations}, whose total spaces are called \emph{Seifert manifolds}. These manifolds were introduced by Herbert Seifert in 1933 and fit perfectly in the much more recent geometrisation perspective: in the next chapter we will introduce the eight three-dimensional geometries, and prove that six of them are realised precisely by the Seifert manifolds.

\subsection{Definition}
We define the Seifert manifolds as Dehn fillings of trivial bundles over surfaces with boundary. Here are the details.

Let $M$ be the (unique) oriented bundle $S\timforsetil S^1$ over a compact connected (possibly non-orientable) surface $S$ with boundary. We denote by $S$ the zero-section.

Let $T_1,\ldots, T_k$ be the boundary tori of $M$. On each $T_i$ we choose an orientation for the meridian $m_i=T_i\cap \partial S$ and for the fibre $l_i$ of the bundle so that the basis $(m_i,l_i)$ for $H_1(T_i,\matZ)$ be positively oriented.
 
A $(p_i,q_i)$-Dehn filling on $T_i$ kills the slope $p_im_i + q_il_i$. We say that the Dehn filling is \emph{fibre-parallel} if $p_i=0$, \emph{i.e.}~if it kills a fibre.

\begin{defn}
A \emph{Seifert manifold} is any 3-manifold $N$ obtained from $M$ by Dehn filling some $h\leqslant k$ boundary tori in a non-fibre-parallel way, that is with $p_i\neq 0$ for all $i$.\index{Seifert manifold} 
\end{defn}

The Seifert manifold is closed if $h=k$, and has $k-h$ boundary tori otherwise. It is not important to know which $h$ tori are filled, in virtue of the following.

\begin{prop}
Every permutation of the boundary tori is realised by a self-diffeomorphism of $M$ that preserves the pairs $\pm (m_i,l_i)$.
\end{prop}
\begin{proof}
Every permutation of the boundary circles of $S$ is realised by a self-diffeomorphism of $S$, that extends orientation-preservingly to the orientable $I$-bundle and its double $M$.
\end{proof}

The pair $(p_i,q_i)$ is determined up to sign, so we can always suppose $p_i>0$ and we fully encode the Seifert manifold $N$ using the following notation: 
\begin{equation} \label{Seifert:notation:eqn}
N = \big(\hat S, (p_1,q_1), \ldots, (p_h, q_h)\big)
\end{equation}
where $\hat S$ is $S$ with $h$ boundary components capped. The reason for using $\hat S$ instead of $S$ is that $N$ has a particular fibration onto $\hat S$, as we will soon see. Before constructing this fibration we list some simple examples that should hopefully help the reader to familiarise with the notation (\ref{Seifert:notation:eqn}), that will be used extensively in the whole chapter.

\begin{example}
The Seifert manifold $\big(S_g,(1,e)\big)$ is the circle bundle over the orientable genus-$g$ surface $S_g$ with Euler number $e$, by construction. In particular $\big(S_g,(1,0)\big) = S_g \times S^1$.
\end{example}

\begin{example}
The Seifert manifold $\big( S^2, (p,q) \big)$ is diffeomorphic to the lens space $L(q,p)$.
\end{example}

The following facts follow from Exercise \ref{fill:T:ex}.

\begin{ex} \label{solid:ex}
The Seifert manifold $\big(D^2, (p,q)\big)$ is a solid torus.
\end{ex}

In the following exercise, it should be clear that the resulting manifold is a lens space (by Exercise \ref{fill:T:ex}) and the hard part is to check carefully that the proposed lens space parameters are correct.
\begin{ex} \label{two:lens:ex}
The Seifert manifold $\big(S^2, (p_1,q_1), (p_2, q_2)\big)$ is the lens space $L (p_1q_2 + q_1p_2, rq_2 + sp_2)$
where $r,s$ are such that $p_1s-q_1r=\pm 1$. In particular, this Seifert manifold is $S^3$ when $p_1q_2+q_1p_2 = \pm 1$.
\end{ex}

\subsection{Seifert fibrations}
As we anticipated, a Seifert manifold $N$ as in (\ref{Seifert:notation:eqn}) possesses some kind of singular fibration over the filled surface $\hat S$. We clarify this point here by defining the notion of \emph{Seifert fibration}. 

Let $(p,q)$ be two coprime integers with $p>0$. A \emph{standard fibered solid torus} with coefficients $(p,q)$ is the solid torus\index{solid torus!standard fibered solid torus} 
$$D \times [0,1]/_\psi$$
where $\psi\colon D \times 0 \to D \times 1$ is a rotation of angle $2\pi \frac qp$. The fibration into vertical segments $\{{\rm pt}\} \times [0,1]$ extends to a fibration into circles of the solid torus. The central fibre obtained by identifying the endpoints of $0 \times [0,1]$ is the core of the solid torus, and every non-central fibre winds $p$ times around the core of $M$: see Figure \ref{Seifert_fibration:fig}.

\begin{figure}
\begin{center}
\includegraphics[width = 2.5 cm] {\iftoggle{BW}{Seifert_fibration-BW}{Seifert_fibration}}
\nota{A standard fibered solid torus. We identify the top and bottom discs by a $2\pi \frac qp$ rotation, for some $q$ coprime with $p$. Here $p=5$. Every non-central fibre (\iftoggle{BW}{light grey}{green}) winds $p$ times along the central fibre (\iftoggle{BW}{dark grey}{red}).}
\label{Seifert_fibration:fig}
\end{center}
\end{figure}

The positive number $p$ is the \emph{multiplicity} of the central fibre. If $p=1$ the fibered solid torus is diffeomorphic to the usual product fibration $D\times S^1$ and the central fibre is \emph{regular}. If $p>1$ the central fibre is \emph{singular}.

\begin{defn}
A \emph{Seifert fibration} is a partition of a compact oriented 3-manifold $N$ with (possibly empty) boundary into circles, such that every circle has a fibered neighbourhood diffeomorphic to a standard fibered solid torus.\index{Seifert fibration}
\end{defn}

Let $S$ be the topological space obtained from $N$ by quotienting circles to points.

\begin{prop}
The space $S$ is a compact connected surface with (possibly empty) boundary.
\end{prop}
\begin{proof}
Every standard fibered solid torus quotients to a disc: a sector $z= \rho e^{i\theta}$ with $0\leqslant \rho \leqslant 1$ and $0\leqslant \theta \leqslant \frac{2\pi}p$ with two edges $\theta=0$ and $\theta = \frac{2\pi}p$ identified by a rotation. This disc actually looks like a ``cone'', and this motivates the following discussion.
\end{proof}

The map $N \to S$ is in fact what we call a Seifert fibration. The surface $S$ may have boundary and may be non-orientable, and its interior has a natural orbifold structure: if the preimage of $x\in S$ is a fibre of order $p$, we see $x$ as a cone point of order $p$, see Section \ref{orbifold:surface:subsection}. If $N$ has boundary, then $S$ also has, and we say that $S$ itself is an orbifold for simplicity although we actually mean only its interior. Morally, we should consider the fibration $N\to S$ as a circle bundle over the orbifold $S$.

A Seifert fibration without singular fibres is just an ordinary circle bundle. 
We now show that Seifert manifolds and Seifert fibrations are more or less the same thing. 

\begin{prop} 
The Seifert manifold
\begin{equation*} 
N = \big(S, (p_1,q_1), \ldots, (p_h, q_h)\big)
\end{equation*}
has a Seifert fibration $N \to S$ over the orbifold 
$$(S, p_1,\ldots, p_h).$$
Every Seifert fibration arises in this way.
\end{prop}
\begin{proof}
The Seifert manifold $N$ is obtained by filling $h$ components of a bundle $M=S'\timforsetil S^1$. The bundle $M \to S'$ extends to a Seifert fibration $N \to S$ where the orbifold $S$ is obtained from $S'$ by attaching $h$ discs with cone points $p_1,\ldots, p_h$.

More precisely, we fill each boundary torus $T_i$ of $M$ with a solid torus having meridian $\mu_i= p_im_i + q_il_i$. We fix a longitude $\lambda_i= r_im_i + s_il_i$ for this solid torus by choosing $r_i,s_i$ with $p_is_i-q_ir_i=1$. We get $l_i = p_i\lambda_i - r_i\mu_i$. By hypothesis the Dehn filling is not fibre-parallel, hence $p_i\neq 0$ and the fibration $M\to S'$ extends to a standard fibration of the solid torus with coefficients $(p_i,-r_i)$.

Every Seifert fibration $N\to S$ arises in this way: if we remove the singular fibres (or a regular one, if there are not) we get an ordinary circle bundle over a surface with boundary, which is trivial. Therefore $N$ is a Dehn filling of this trivial bundle, hence a Seifert manifold.
\end{proof}

We have seen that the notation
\begin{equation} \label{notation:eqn}
N = \big(S, (p_1,q_1), \ldots, (p_h, q_h)\big)
\end{equation}
defines a Seifert fibration $N \to S$ and a Seifert manifold $N$.

\begin{example} \label{bad:example}
If the orbifold $S$ is a disc with at most one singular point then $N$ is a standard fibered solid torus. If $S$ is $S^2$ with at most $2$ singular points then $N$ is a lens space: see Exercises \ref{solid:ex} and \ref{two:lens:ex}.
\end{example}

\subsection{Classification of Seifert fibrations}
We say that two Seifert fibrations $\pi_1\colon N_1\to S$, $\pi_2\colon N_2 \to S$ are \emph{isomorphic} if there is a diffeomorphism $\psi\colon N_1 \to N_2$ such that $\pi_1 = \pi_2\circ \psi$. Two different notations as in (\ref{notation:eqn}) may describe isomorphic fibrations, but this phenomenon is completely understood.

\begin{prop} \label{S:moves:prop}
Two notations as in (\ref{notation:eqn}) describe two orientation-preservingly isomorphic Seifert fibrations if and only if they are related by a finite sequence of the following moves and their inverses:
\begin{align}
(p_i,q_i), (p_{i+1},q_{i+1}) & \longmapsto (p_i, q_i+p_i), (p_{i+1},q_{i+1}-p_{i+1}), \label{Seifert:1:eqn} \\
(p_1,q_1), \ldots, (p_h, q_h) & \longmapsto (p_1,q_1), \ldots, (p_h, q_h), (1,0), \label{Seifert:2:eqn} \\
(p_i,q_i) & \longmapsto (p_i,q_i+p_i) \quad {\rm if}\ \partial N \neq\emptyset, \label{Seifert:3:eqn}
\end{align}
and permutations of the pairs $(p_i,q_i)$'s.
\end{prop}
\begin{proof}
Recall that $N$ is a Dehn filling of $M=S'\timforsetil S^1$. Move (\ref{Seifert:1:eqn}) is the result of twisting $M$ along a fibered annulus $A$ connecting the tori $T_i$ and $T_{i+1}$ in $\partial M$: this self-diffeomorphism of $M$ acts on $T_i$ and $T_{i+1}$ like two opposite Dehn twists and extends to an isomorphism of the two fibrations. In (\ref{Seifert:3:eqn}) we twist along an annulus connecting $T_i$ and $\partial N$. The move (\ref{Seifert:2:eqn}) corresponds to drilling a nonsingular fibered torus and refilling it back.

We now prove that these moves suffice to connect two isomorphic Seifert fibrations. 
Suppose that two distinct notations as in (\ref{notation:eqn}) describe isomorphic fibrations. 
We use the moves (\ref{Seifert:1:eqn}), (\ref{Seifert:2:eqn}), and (\ref{Seifert:3:eqn}) to eliminate the parameters $p_i=1$ as much as possible from both notations. If there are no singular fibres, we end up with a single parameter $(1,e)$ if $S$ is closed, and no parameters at all if $\partial S \neq \emptyset$. We conclude by Corollary \ref{unique:e:cor}. 

If the fibration has at least one singular fibre, we can eliminate all $p_i=1$ and the $(p_i,q_i)$ correspond to singular fibres. An isomorphism of Seifert fibrations sends singular fibres to singular fibres and hence induces an isomorphism of their complement $S'\timforsetil S^1$. The parameters $(p_i,q_i)$ are determined by the choice of a section in $S'\timforsetil S^1$. Different sections are related by Dehn twist along annuli and hence the parameters are related by the moves (\ref{Seifert:1:eqn}) and (\ref{Seifert:3:eqn}).
\end{proof}

Proposition \ref{S:moves:prop} classifies all the Seifert fibrations up to isomorphism. A classification of Seifert manifolds up to \emph{diffeomorphism} would also be desirable, but it is much harder to obtain because a three-manifold may admit many non-isomorphic Seifert fibrations. For instance, 
Exercise \ref{two:lens:ex} shows that the lens spaces may fibre in many different ways; a manifold as familiar as $S^3$ fibres over the orbifold $(S^2, p_1, p_2)$ if $p_1$ and $p_2$ are coprime and hence has infinitely many non-isomorphic fibrations. It is a stimulating exercise to try to visualise these Seifert fibrations of $S^3$.

We now start a long journey in Seifert manifolds theory, whose ultimate goal is to classify them completely up to diffeomorphism. We will see at the end that the only Seifert manifolds admitting non-isomorphic fibrations are the ``smallest'' ones, like $S^3$, the lens spaces, and few more that will be classified using some \emph{ad hoc} argument.

\begin{ex}
The number of non-isomorphic Seifert fibrations over a fixed orbifold $S$ is finite $\Longleftrightarrow$ $\partial S \neq \emptyset$.
\end{ex}

\subsection{Euler number}
We now extend the notion of Euler number from ordinary to Seifert fibrations.
We define the \emph{Euler number} of the fibration (\ref{notation:eqn}) to be the rational number\index{Euler number}
$$e(N) = \sum_{i=1}^h \frac{q_i}{p_i}.$$
The Euler number is only defined modulo $\matZ$ when $N$ has boundary. The good definition follows from Proposition \ref{S:moves:prop} (the moves do not affect $e$, except (\ref{Seifert:3:eqn}) that modifies $e$ into $e+1$ and applies only when $N$ has boundary) and is coherent with the circle bundle case. The Euler number depends on the fibration and not only on $N$, but we write $e(N)$ anyway for simplicity. Proposition \ref{S:moves:prop} easily implies the following.

\begin{cor} \label{classification:Seifert:cor}
Two Seifert fibrations
$$\big(S, (p_1,q_1), \ldots, (p_h, q_h)\big), \quad \big(S', (p_1',q_1'), \ldots, (p_{h'}', q_{h'}')\big)$$
with $p_i, p_i'\geqslant 2$ are orientation-preservingly isomorphic if and only if $S=S'$, $h=h'$, $e=e'$, and up to reordering $p_i=p_i'$ and $q_i \equiv q_i'$ (mod $p_i$) for all $i$.
\end{cor}
The numbers $e$ and $e'$ indicate the Euler numbers of the two fibrations, and recall that they are only defined modulo $\matZ$ when $\partial S \neq \emptyset$.

\begin{rem} \label{reverse:rem}
The move
$$\big(S,(p_1,q_1), \ldots, (p_h,q_h)\big) \longmapsto \big(S,(p_1,-q_1), \ldots, (p_h,-q_h)\big) $$
corresponds to a change of orientation for the three-manifold and transforms $e$ into $-e$.
\end{rem}

\subsection{Homology} \label{Seifert:homology:sphere:subsection}
We now characterise the Seifert manifolds that have finite homology groups, and in particular the homology spheres.

\begin{prop} \label{Seifert:homology:prop}
The homology group $H_1(M, \matZ)$ of 
$$M=\big(S,(p_1,q_1), \ldots, (p_h,q_h)\big)$$
is finite $\Longleftrightarrow$ one of the following holds:
\begin{itemize}
\item $S=S^2$ and $e\neq 0$, and we get $|H_1(M,\matZ)| = |e|p_1\cdots p_h$;
\item $S=\matRP^2$, and we get $|H_1(M,\matZ)| = 4p_1\cdots p_h$.
\end{itemize}
\end{prop}
\begin{proof}
If $M$ is not closed then $H_1(M,\matZ)$ is infinite because $H_1(\partial M,\matZ)$ is. If $S\neq S^2, \matRP^2$ then $H_1(S)$ is infinite. It is easy to check that the fibration $M\to S$ induces a surjection $H_1(M)\to H_1(S)$ and hence $H_1(M)$ is also infinite.

Suppose $S=S^2$. The manifold $M$ is a $(p_i,q_i)$-Dehn filling of $S_{0,h}\times S^1$ where $S_{0,h}$ is the sphere with $h$ holes. The homology of $S_{0,h}\times S^1$ is generated by $m_1,\ldots,m_h,l$ with the relation $m_1+\ldots + m_h=0$. The $i$-th Dehn filling adds the relation $p_im_i + q_il_i=0$. The relations form a square $(h+1)$-matrix
$$
\begin{pmatrix}
1 & p_1 & 0 & \cdots & 0 \\
1 & 0 & p_2 & \ddots & \vdots \\
\vdots & \vdots & \ddots & \ddots & 0 \\
\phantom{\Big|} \! 1 & 0 & 0 & \cdots & p_h \\
0 & q_1 & q_2 & \cdots & q_h
\end{pmatrix}
\sim
\begin{pmatrix}
1 & p_1 & 0 & \cdots & 0 \\
1 & 0 & p_2 & \ddots & \vdots \\
\vdots & \vdots & \ddots & \ddots & 0 \\
\phantom{\Big|} \! 1 & 0 & 0 & \cdots & p_h \\
-e & 0 & 0 & \cdots & 0
\end{pmatrix}.
$$
We have used Gauss moves to simplify the last row. The determinant of this matrix is $\pm ep_1\cdots p_h$. The order $|H_1(M,\matZ)|$ is the absolute value of the determinant if it is non-zero, and is infinite if it is zero.

Suppose $S=\matRP^2$ and let $S^{\rm no}_{1,h}$ be $\matRP^2$ with $h$ holes. The homology of $S^{\rm no}_{1,h} \timtil S^1$ is generated by $a,m_1,\ldots, m_h, l$, where $a$ is an orientation-reversing curve in $S^{\rm no}_{1,h}$, with the relations $2a+m_1+\ldots +m_h=0$ and $2l=0$. The relations now form a square $(h+2)$-matrix with determinant
$$
\det \begin{pmatrix}
2 & 0 & 0 & \ldots & 0 \\
1 & 0 & p_1 & \ddots & 0 \\
\vdots & \vdots & \vdots & \ddots & \vdots \\
\phantom{\Big|} \! 1 & 0 & 0 & \cdots & p_h \\
0 & 2& q_1 & \cdots & q_h
\end{pmatrix}
= \pm 2 \det
\begin{pmatrix}
2 & 0 & \cdots & 0 \\
1 & p_1 & \ddots & \vdots \\
\vdots & \vdots & \ddots & \vdots \\
1 & 0 & \cdots & p_h
\end{pmatrix} = \pm 4 p_1\cdots p_h.
$$
The proof is complete.
\end{proof}

We deduce an elegant description of all the Seifert homology spheres.

\begin{cor} \label{sigma:cor}
For every set $p_1, \ldots, p_h$ of pairwise coprime integers $p_i\geqslant 2$ there is a unique homology sphere $\Sigma(p_1,\ldots, p_h)$ Seifert-fibering over $(S^2,p_1,\ldots,p_h)$. Every homology sphere Seifert manifold arises in this way.
\end{cor}
\begin{proof}
By the previous proposition a Seifert manifold 
$$M=\big(S,(p_1,q_1), \ldots, (p_h,q_h)\big)$$
is a homology sphere if and only if $S=S^2$ and $|ep_1\cdots p_h|=1$. We have
$$ep_1\cdots p_h = \sum_{i=1}^h \frac{q_i}{p_i}p_1\cdots p_h = \sum_{i=1}^h q_ip_1\cdots \widehat {p_i}\cdots p_h = \sum_{i=1}^h q_ip_i'$$
where we set $p_i' = p_1\cdots \widehat {p_i}\cdots p_h$. The integers $p_1,\ldots, p_h$ are \emph{pairwise} coprime if and only if $p_1',\ldots, p_h'$ are \emph{globally} coprime (no prime number divides all of them). The equation $\sum_{i=1}q_ip_i' = 1$ is satisfied by some $q_1,\ldots, q_h$ $\Longleftrightarrow$ they are globally coprime. Different solutions $q_i$ are related by moves as in Proposition \ref{S:moves:prop} (exercise) and produce the same Seifert fibration.
\end{proof}

For instance, the homology spheres $\Sigma(p_1)$ and $\Sigma(p_1,p_2)$ are just $S^3$ by Example \ref{bad:example}. The simplest homology sphere with three singular fibres is $\Sigma(2,3,5)$. This manifold is called the \emph{Poincar\'e homology sphere} and we will soon see that its fundamental group has order 120.\index{Poincar\'e homology sphere} All the other Seifert homology spheres have infinite fundamental group: this is related to the fact that $\sum \frac 1{p_i}>1$ for all choices of $(p_1,\ldots, p_h)$ with $h\geqslant 3$ except $(2,3,5)$, as we will soon see. A Seifert homology sphere of type $\Sigma(p_1,p_2,p_3)$ is called a \emph{Brieskorn homology sphere}\index{Brieskorn homology sphere}.

\subsection{Coverings} 
We now start to investigate the coverings of Seifert fibered spaces. Our main goal will be to subdivide the Seifert manifolds into nine classes: this will be done in Section \ref{commensurability:Seifert:subsection}.

Like ordinary fibrations, Seifert fibrations behave well with respect to coverings.
Let $M \to S$ be a Seifert fibration and $\tilde M \to M$ a covering. The foliation into circles of $M$ lifts to a foliation into circles or lines in $\tilde M$, with some quotient space $\tilde S$. 

\begin{prop} \label{Seifert:lifts:prop}
The quotient $\tilde S$ is an orbifold covering of $S$. 
\begin{itemize}
\item If $\tilde M$ foliates in circles then $\tilde M \to \tilde S$ is a Seifert fibration, 
\item If $\tilde M$ foliates in lines then $\tilde M \to \tilde S$ is a line bundle.
\end{itemize}
In the second case $\tilde S$ has no singular points.
\end{prop}
\begin{proof}
This holds on all coverings of a standard fibered solid torus and hence holds everywhere.
\end{proof}

We now concentrate ourselves on the finite-degree case.

\subsection{Finite-degree coverings}
We define a \emph{finite-degree covering} of a Seifert fibration $M\to S$ to be a commutative diagram
$$
\xymatrix{ 
\tilde M\ar[r] \ar[d] & M \ar[d] \\
\tilde S\ar[r] & S
}
$$
where $\tilde M \to \tilde S$ is a Seifert fibration, $\tilde M \to M$ is a finite-degree covering, and $\tilde S \to S$ is an orbifold covering. Proposition \ref{Seifert:lifts:prop} implies the following.

\begin{cor}
Let $M\to S$ be a Seifert fibration. Every finite-degree covering $\tilde M \to M$ induces a finite-degree covering of Seifert fibrations:
$$
\xymatrix{ 
\tilde M\ar[r] \ar@{.>}[d] & M \ar[d] \\
\tilde S\ar@{.>}[r] & S
}
$$
\end{cor}
The dotted arrows indicate the maps that are induced. The degree $d$ of such a covering $\tilde M \to M$ splits into two parts: 
\begin{itemize}
\item the \emph{horizontal} degree $d_{\rm h}$ is the degree of the covering $\tilde S \to S$, 
\item the \emph{vertical} degree $d_{\rm v}$ is the degree with which a regular fibre of $\tilde M$ covers a regular fibre of $M$. 
\end{itemize}
The vertical degree is well-defined since regular fibres in $\tilde M$ form a connected set. The pre-image of a regular fibre in $M$ consists of $d_{\rm h}$ regular fibres in $\tilde M$, each fibering with degree $d_{\rm v}$. Therefore
$$d = d_{\rm h} \cdot d_{\rm v}.$$
A covering $\tilde M \to M$ is \emph{horizontal} or \emph{vertical} if respectively
$d_{\rm v} = 1$ or $d_{\rm h} = 1$.

\begin{prop}[Pull-back] \label{horizontal:prop}
Let $M\to S$ be a Seifert fibration. Every finite-degree orbifold covering $\tilde S \to S$ is induced by a unique horizontal covering of Seifert fibrations:
$$
\xymatrix{ 
\tilde M\ar@{.>}[r] \ar@{.>}[d] & M \ar[d] \\
\tilde S\ar[r] & S
}
$$
\end{prop}
\begin{proof}
Like for ordinary bundles, there is a unique way to define $\tilde M$ by pulling back the Seifert fibration on fibered solid tori.
\end{proof}

Recall that an orbifold is \emph{very good} when it is finitely covered by a surface. Every locally orientable 2-orbifold is very good except the bad orbifolds $S^2(p_1)$ and $S^2(p_1,p_2)$ with $p_1\neq p_2$, see Theorem \ref{O:geometrisation:teo} and Corollary \ref{very:good:surface:cor}.

\begin{cor} \label{good:bundle:cor}
If $S$ is good, every Seifert fibration $M\to S$ is finitely covered by a circle bundle over a surface.
\end{cor}
\begin{proof}
Pull-back the fibration along the surface cover $\tilde S \to S$.
\end{proof}

\begin{prop} \label{hv:prop}
Every finite-degree covering between Seifert fibrations is a composition of one vertical and one horizontal covering. 
\end{prop}
\begin{proof}
Let $p\colon \tilde M\to M$ be a covering of Seifert fibrations, with base spaces $\tilde S \to S$. If we pull-back $M$ to $\tilde S$ we get a horizontal covering $p_{\rm h}\colon M_{\rm h} \to M$. There is a natural vertical $p_{\rm v}\colon \tilde M \to M_{\rm h}$ such that $p = p_{\rm h} \circ p_{\rm v}$.
\end{proof}

\subsection{Circle bundle coverings}
The Euler number of the fibration and the Euler characteristic of the base orbifold behave well with coverings.

\begin{prop} \label{e:prop}
Let $p\colon \tilde M \to M$ be a finite covering of closed Seifert fibrations with base orbifolds $\tilde S$ and $S$, with degrees $(d_{\rm h}, d_{\rm v})$. We have
\begin{align*}
\chi(\tilde S) & = d_{\rm h} \cdot \chi (S), \\
e(\tilde M ) & = \frac{d_{\rm h}}{d_{\rm v}} \cdot e(M).
\end{align*}
\end{prop}
\begin{proof}
The first equality holds for every orbifold covering $\tilde S \to S$.
Concerning the second one, by Proposition \ref{hv:prop} we may suppose $p$ is either vertical or horizontal. Write
$$M = \big(S, (p_1,q_1), \ldots (p_h, q_h)\big)$$
and recall that $M$ is obtained from a circle bundle $N$ by $(p_i,q_i)$-filling the boundary torus $T_i$ for all $i=1,\ldots, h$. Set $\tilde T_i = p^{-1}(T_i)$ and $\tilde N = p^{-1}(N)$.

If $p$ is vertical, we fix a section of $\tilde  N$ and note that it projects to a section of $N$: these sections induce meridians on $\tilde T_i$ and $T_i$. Here $\tilde T_i$ is a single torus and $\tilde T_i \to T_i$ wraps the fibre with degree $d=d_{\rm v}$. Since the covering extends to the filled solid tori, we have $d|q_i$ and $\tilde T_i$ is filled with parameters $(p_i,q_i/d)$. Therefore $e(\tilde M) = e(M)/d$.

If $p$ is horizontal, we fix a section of $N$, and its counterimage is a section of $\tilde N$: these sections induce meridians on $T_i$ and $\tilde T_i$.
Now $\tilde T_i$ consists of some tori $\tilde T_i = \tilde T_{i,1} \sqcup \ldots \sqcup \tilde T_{i,k_i}$ and $\tilde T_{i,j} \to T_i$ wraps the meridians with some degree $d_{i,j}$. The total sum of these local degrees is $d_{i,1}+\ldots +d_{i,k_i}=d = d_{\rm h}$. Similarly as above we have $d_{i,j}|p_i$ and the filling solid torus at $\tilde T_{i,j}$ has parameters $(p_i/d_{i,j},q_i)$. Therefore 
$$e(\tilde M) = \sum_{i=1}^h \sum_{j=1}^{k_i} \frac{d_{i,j}q_i}{p_i} = \sum_{i=1}^h d\cdot \frac{q_i}{p_i} = d \cdot e(M).$$
This concludes the proof.
\end{proof}

An important consequence is that the signs of $\chi(S)$ and $e(M)$ are invariant under finite coverings. 

\subsection{Commensurability classes} \label{commensurability:Seifert:subsection}
We now would like to subdivide the Seifert manifolds into few classes, and to this purpose we introduce a general equivalence relation between manifolds.

\begin{defn}
Two manifolds $M$ and  $N$ are \emph{commensurable} if there is a manifold that covers both $M$ and $N$ with finite degrees.\index{commensurability of manifolds}
\end{defn}

\begin{prop}
Commensurability is an equivalence relation.
\end{prop}
\begin{proof}
If $M$ is commensurable with $N_1$ and $N_2$, it has finite-sheeted coverings $M_1, M_2$ that cover $N_1$ and $N_2$ corresponding to finite-index subgroups $\Gamma_1, \Gamma_2 < \pi_1(M)$. The subgroup $\Gamma = \Gamma_1\cap \Gamma_2$ has also finite index and determines a manifold that covers both $N_1$ and $N_2$.
\end{proof}

\begin{ex}
There are three commensurable classes of closed surfaces, determined by their Euler characteristic being positive, null, or negative.
\end{ex}

We now classify the commensurability classes of Seifert manifolds. The \emph{3-torus} is of course $T\times S^1 = S^1\times S^1\times S^1$.\index{3-torus}

\begin{prop} \label{cover:Seifert:prop}
A closed Seifert fibration $M\to S$ has:
\begin{itemize}
\item $\chi(S)>0$ and $e=0$ $\Longleftrightarrow M$ is covered by $S^2\times S^1$,
\item $\chi(S)>0$ and $e\neq 0$ $\Longleftrightarrow M$ is covered by $S^3$,
\item $\chi(S)=0$ and $e=0$ $\Longleftrightarrow M$ is covered by the 3-torus,
\item $\chi(S)=0$ and $e\neq 0$ $\Longleftrightarrow M$ is covered by a twisted bundle over $T$,
\item $\chi(S)<0$ and $e=0$ $\Longleftrightarrow M$ is covered by $S_g \times S^1$ for some $g>1$,
\item $\chi(S)<0$ and $e\neq 0$ $\Longleftrightarrow M$ is covered by a twisted bundle over $S_g$ for some $g>1$.
\end{itemize}
\end{prop}
\begin{proof}
If $S$ is a bad orbifold, then $S=S^2(p)$ or $S^2(p_1,p_2)$ with $p_1\neq p_2$ and hence we get $e\neq 0$ and $\chi(S)>0$. The manifold $M$ is a lens space by Example \ref{bad:example} and is hence covered by $S^3$.

If $S$ is good, the fibration $M\to S$ is covered by a circle bundle $\tilde M\to \tilde S$ over an orientable closed surface $\tilde S$ by Corollary \ref{good:bundle:cor}. By Proposition \ref{e:prop} the numbers $\chi(\tilde S)$ and $e(\tilde M)$ have the same signs of $\chi(S)$ and $e(M)$. Therefore $\tilde S = S^2, T$, or $S_g$ with $g> 1$, depending on whether $\chi(S)$ is positive, null, or negative. 
Exercise \ref{circle:bundle:ex} says that the circle bundle $\tilde M \to \tilde S$ is trivial $\Longleftrightarrow e(\tilde M) =0 \Longleftrightarrow e(M) =0$. Note that a non-trivial bundle over $S^2$ is a lens space $L(e,1)$ with $e\neq 0$ and is hence covered by $S^3$.

We have proved that $M$ is covered (according to the signs of $\chi$ and $e$) by a manifold belonging to one of the six types: $S^3$, $S^2\times S^1$, the 3-torus, a twisted bundle over $T$, $S_g\times S^1$, and a twisted bundle over $S_g$. It remains to prove that $M$ cannot be covered by two manifolds $M_1, M_2$ belonging to two different types: this holds because manifolds of distinct types are not commensurable. To prove that, note that the finite cover of a manifold of one of the six types is a manifold of the same type, and Corollary \ref{distinct:bundle:cor} implies that a manifold cannot belong to two different types. If two manifolds of distinct types were commensurable they would be covered by a manifold belonging to both types, yielding a contradiction.
\end{proof}

\begin{table}
\begin{center}
\begin{tabular}{c||ccc}
        & \phantom{\Big|} $\chi>0$ & $\chi = 0$ & $\chi<0$ \\
 \hline \hline
 $e=0$ & \phantom{\Big|} $S^2\times S^1$ & $S^1\times S^1\times S^1$ & $S_2 \times S^1$ \\
 $e\neq 0$ & \phantom{\big|} $S^3$ & $\big(T,(1,1)\big)$ & $\big(S_2,(1,1)\big)$ \\
 $\partial M\neq \emptyset$ & \phantom{\Big|} $D\times S^1$ & $T \times [0,1]$ & $P \times S^1$
\end{tabular}
\vspace{.2 cm}
\nota{There are 9 commensurability classes of Seifert manifolds: 6 closed and 3 with boundary. Every Seifert manifold $M$ is commensurable with one (and only one) of these 9 manifolds. The commensurability class of $M$ is easily detected by looking at its invariants $\chi$ and $e$ (the latter only in the closed case). The surfaces $T,S_2,D,P$ are the torus, the genus-2 surface, the disc, and the pair-of-pants. Note that $S^3 = \big(S^2,(1,1)\big)$.}
\label{commensurability:table}
\end{center}
\end{table}

We can now easily classify all closed Seifert manifolds up to commensurability. To conclude we just need to solve the following exercise.

\begin{ex} \label{e:construct:ex}
Pick $S=S_g$ and $e>0$. Construct:
\begin{itemize}
\item a degree-$e$ vertical covering $\big(S,(1,1)\big)\to \big(S,(1,e)\big)$,
\item a degree-$e$ horizontal covering $\big(\tilde S,(1,e)\big) \to \big(S,(1,1)\big)$ if $g\geqslant 1$.
\end{itemize}
\end{ex}
\begin{proof}[Hint] If $g\geqslant 1$ then $S$ has covers of any degree $e$. To construct them, pick surjective homomorphisms $\pi_1(S)\to H_1(S) \to \matZ/_{e\matZ}$.
\end{proof}

\begin{cor}
There are six commensurability classes of closed Seifert manifolds, depending on $\chi$ and $e$ as shown in Table \ref{commensurability:table}.
\end{cor}
\begin{proof}
By the previous exercise the non-trivial (or trivial) bundles over closed surfaces with $\chi < 0$ (or $\chi =0$, $\chi >0$) are all commensurable.
\end{proof}

\begin{cor}
A closed Seifert fibration $M\to S$ has $e=0$ $\Longleftrightarrow$ it is finitely covered by a trivial circle bundle.
\end{cor}

We recall that a Seifert manifold $M$ may have non-isomorphic Seifert fibrations $M\to S$ and $M\to S'$, and the invariants $\chi, \chi'$ and $e, e'$ of the two fibrations are not necessarily equal; however $\chi$ is positive, null, or negative if and only if $\chi'$ is, and $e$ vanishes if and only if $e'$ does (in the closed case). This holds because the commensurability class of $M$ does not depend on its fibration. For instance, Exercise \ref{two:lens:ex} shows that a lens space $M=L(p,q)$ with $p>0$ has many fibrations, but each with $\chi>0$ and $e\neq 0$.

We now consider the boundary case, which is simpler because circle bundles over orientable surfaces with boundary are always trivial. Let $A$ denote the annulus and $S_{g,b}$ be the surface of genus $g$ with $b$ open discs removed.

\begin{prop}
A Seifert fibration $M\to S$ with boundary has
\begin{itemize}
\item $\chi(S)>0 \Longleftrightarrow M = D\times S^1$,
\item $\chi(S)=0 \Longleftrightarrow M$ is covered by $A \times S^1 = T\times [0,1]$,
\item $\chi(S)<0 \Longleftrightarrow M$ is covered by $S_{g,b} \times S^1$ for some $g+b>2$.
\end{itemize}
\end{prop}
\begin{proof}
Every Seifert manifold with boundary is covered by a circle bundle over an orientable surface, and such a bundle is trivial here. If $\chi>0$ the base surface is a disc with at most one cone point and hence $M$ is a standard fibered solid torus. 
\end{proof}

\begin{cor} There are three commensurability classes of Seifert manifolds with boundary, depending on $\chi$ as shown in Table \ref{commensurability:table}.
\end{cor}
\begin{proof}
The surfaces $S_{g,b}$ with $g+b>2$ are commensurable (exercise).
\end{proof}

We now characterise some commensurability classes by looking at the fundamental groups.

\begin{table}
\begin{center}
\begin{tabular}{c||ccc}
        & \phantom{\Big|} $\chi>0$ & $\chi = 0$ & $\chi<0$ \\
 \hline \hline
 $e=0$ & \phantom{\Big|} $\matZ$ & $\matZ^3$ & \\
 $e\neq 0$ & \phantom{\big|} $\{e\}$ & &  \\
 $\partial M\neq \emptyset$ & \phantom{\Big|} $\matZ$ & $\matZ^2$ &
\end{tabular}
\vspace{.2 cm}
\nota{There are five Seifert manifolds $M$ with $\pi_1(M) = \matZ^h$ and they belong to distinct commensurability classes. The fundamental group $\pi_1(N)$ of every other manifold $N$ in these five classes is virtually abelian: it contains $\matZ^h$ as a finite-index subgroup.}
\label{abelian:group:table}
\end{center}
\end{table}

\subsection{Virtually abelian fundamental groups} \label{virtually:abelian:subsection}
A group $G$ is \emph{virtually abelian of rank $h$} if it contains $\matZ^h$ as a finite-index subgroup. 

This is a finite-index-independent property: if $G'<G$ has finite index, the group $G'$ is virtually abelian of rank $h$ if and only if $G$ is. Therefore a manifold $M$ has a virtually abelian fundamental group of rank $h$ $\Longleftrightarrow$ every manifold $N$ commensurable with $M$ also has $\Longleftrightarrow$ there is a manifold $N$ in the commensurability class with $\pi_1(N)=\matZ^h$.

\begin{prop} \label{virtually:abelian:prop}
There are five commensurability classes of Seifert manifolds with virtually abelian fundamental groups, shown in Table \ref{abelian:group:table}.
\end{prop}
\begin{proof}
These five classes contain $S^2\times S^1$, $S^1\times S^1\times S^1$, $S^3$, $D\times S^1$, and $A\times S^1$. The remaining four classes do not contain manifolds with abelian fundamental group: every manifold in these classes is covered by some $S \times S^1$ with $\chi(S)<0$ or $\big(S_g, (1,e)\big)$ with $g\geqslant 1$ and $e\geqslant 1$ (use Exercise \ref{e:construct:ex}). The fundamental group of such manifolds is not $\matZ^h$: the former because $\pi_1(S)$ is not abelian, and the latter because its abelianization contains torsion, see Exercise \ref{bundle:ex}.
\end{proof}

It could be reasonable to expect these five commensurability classes to contain manifolds of some simpler topological nature, which are easier to study and classify. We will see however in Section \ref{Seifert:finite:subsection} that the Seifert manifolds covered by $S^3$ are quite interesting and their topological classification is certainly not immediate: we already experienced that for lens spaces in Theorem \ref{Lpq:teo}, whose classification is quite involved. The other four classes are indeed of a simpler nature: we now prove that each contains finitely many manifolds, and we classify them completely.

\subsection{Finite commensurability classes} \label{finite:classes:prop}
We now completely identify the commensurability classes of Seifert manifolds that contain only finitely many manifolds.

We start by considering manifolds $M$ with boundary. We have already seen that the class $\chi>0$ contains only the solid torus. The class $\chi=0$ also contains few manifolds:

\begin{prop} \label{diffeo:fiberings:prop}
Every Seifert fibration with boundary and $\chi = 0$ is isomorphic to one of the following:
$$ A \times S^1, \qquad S \timtil S^1, \qquad \big(D, (2,1), (2,-1) \big).$$
The Seifert manifold is correspondingly diffeomorphic to the interval bundles
$ T \times I$, $K \timtil I$, and $K \timtil I$ again. Here $A$, $S$, and $K$ are the annulus, the M\"obius strip, and the Klein bottle.
\end{prop}
\begin{proof}
The orbifolds with boundary and $\chi=0$ are $A$, $S$, and $(D,2,2)$.
Every Seifert fibration over $(D,2,2)$ is isomorphic to $\big(D,(2,2n+1),(2,2m+1)\big)$ and does not depend on $n,m\in \matZ$ by Proposition \ref{S:moves:prop}. Hence there is only one Seifert fibration over the orbifold $(D,2,2)$, which we write as $\big(D,(2,1),(2,-1)\big)$ only for aesthetic reasons (to get $e=0$).

\begin{figure}
\begin{center}
\includegraphics[width = 12 cm] {\iftoggle{BW}{two_fibrations-BW}{two_fibrations}}
\nota{Draw the orientable line bundle over the Klein bottle as $S^1\times [-1,1]\times [-1,1]/_\psi$ where $\psi$ identifies the lower and upper (\iftoggle{BW}{light grey}{yellow}) annuli $S^1\times [-1,1]\times \{\pm 1\}$ by gluing $(e^{i\theta},t)$ to $(e^{-i\theta},-t)$. The Klein bottle is drawn in \iftoggle{BW}{grey}{green} (left). 
This manifold has two fibrations: a horizontal one with fibres $S^1\times x \times y$ and base surface the \iftoggle{BW}{grey}{green} M\"obius strip (centre), and a vertical one with fibres $x\times y \times [-1,1]$ and base surface the \iftoggle{BW}{grey}{green} disc with two cone points of order two (the \iftoggle{BW}{grey}{green} annulus becomes a disc after identifying the points $(e^{i\theta}, 0)$ and $(e^{-i\theta}, 0)$). The two singular fibres of order two lie above $(\pm 1,0)$ and are drawn in \iftoggle{BW}{grey}{red} (right).}
\label{two_fibrations:fig}
\end{center}
\end{figure}

We get three manifolds and we now prove that they are all diffeomorphic to some interval bundle. The diffeomorphism $A\times S^1 \isom T\times I$ is obvious. Draw $K\timtil I$ as $S^1\times [-1,1]\times [-1,1]/_\psi$ as in Figure \ref{two_fibrations:fig}. The manifold has two fibrations: a horizontal one by circles $S^1\times x\times y$, and a vertical one by segments $x\times y \times [-1,1]$. The horizontal one gives $S\timtil S^1$ and the vertical one closes to a Seifert fibration over the disc with two singular fibres of order two, and hence is $(D,(2,1),(2,-1))$.
\end{proof}

\begin{cor}
Every Seifert manifold with boundary and $\chi\geqslant 0$ is diffeomorphic to $D\times S^1$, $T\times I$, or $K\timtil I$.
\end{cor}

We now turn to closed Seifert manifolds. We have seen that $K\timtil I$ fibres in two non-isomorphic ways, and this has some consequences.

\begin{cor} \label{Seifert:diffeo:cor}
The following diffeomorphisms hold:
\begin{align*}
\big(S^2,(2,1),(2,-1),(p,q)\big) & \isom \big(\matRP^2, (q,p)\big), \\
\big(S^2,(2,1),(2,1),(2,-1),(2,-1)\big) & \isom K\timtil S^1.
\end{align*}
\end{cor}
\begin{proof}
Both equalities follow from $\big(D, (2,1), (2,-1)\big) \isom S\timtil S^1$. In the first we look at Figure \ref{two_fibrations:fig} to check that a $(p,q)$ curve in $\big(D,(2,1),(2,-1)\big)$ becomes a $(q,p)$ curve in $S\timtil S^1$.
The second equality is obtained by doubling the line bundle $K\timtil I$.
\end{proof}

We want to classify the closed Seifert manifolds with $\chi\geqslant 0$ and $e=0$. We start with the case $\chi >0$.

\begin{prop} \label{rp:prop}
Every closed Seifert fibration with $\chi>0$ and $e=0$ is isomorphic to one of the following: 
$$S^2 \times S^1, \qquad \matRP^2 \timtil S^1, \qquad \big(S^2, (p,q), (p,-q)\big).$$
The manifolds of the last type are all diffeomorphic to $S^2\times S^1$.
\end{prop}
\begin{proof}
If the base surface $S$ is a sphere with $\leqslant 2$ singular fibres, we use Exercise \ref{two:lens:ex}. Otherwise $S$ is one of the following orbifolds (see Table \ref{orbifold:elliptic:table}):
$$(S^2,2,2,p), \quad (S^2,2,3,3), \quad (S^2,2,3,4), \quad (S^2,2,3,5), \quad \matRP^2, \quad (\matRP^2, p) $$
with $p\geqslant 2$. In all cases except $\matRP^2$, we get $e\neq 0$ for any choice of Dehn filling parameters: for instance 
$$e\big(S^2,(2,q_1),(2,q_2),(p,q_3)\big) = \frac{q_1+q_2}2 + \frac{q_3}p \neq 0.$$ 
The other cases are analogous.
\end{proof}

The manifold $\matRP^2 \timtil S^1$ is not diffeomorphic to $S^2\times S^1$, because they have non-isomorphic fundamental groups. Moreover, $\matRP^2 \timtil S^1$ is not prime:

\begin{ex} The two manifolds $\matRP^2\timtil S^1$ and $\matRP^3\#\matRP^3$ are diffeomorphic.
\end{ex}

\begin{table}
\begin{center}
\begin{tabular}{c||c}
\phantom{\Big|} $M$ & $H_1(M,\matZ)$ \\
\hline \hline
\phantom{\Big|} $ T \times S^1$ & $\matZ^3$ \\ 
\phantom{\big|} $ K \timtil S^1$ & $\matZ \times \matZ/_{2\matZ} \times \matZ/_{2\matZ}$ \\
\phantom{\Big|} $ \big(S^2, (2,1), (2,1), (2,-1), (2,-1) \big) $ & $\matZ \times \matZ/_{2\matZ} \times \matZ/_{2\matZ}$ \\
\phantom{\big|} $ \big(S^2, (3,1), (3,1), (3,-2) \big) $ & $\matZ \times \matZ/_{3\matZ} $ \\
\phantom{\Big|} $ \big(S^2, (2,1), (4,-1), (4,-1) \big) $ & $\matZ \times \matZ/_{2\matZ} $ \\
\phantom{\big|} $ \big(S^2, (2,1), (3,1), (6,-5) \big) $ & $\matZ$ \\
\phantom{\Big|} $ \big(\matRP^2, (2,1), (2,-1)\big) $ & $\matZ/_{4\matZ} \times \matZ/_{4\matZ}$
\end{tabular}
\vspace{.2 cm}
\nota{The seven closed Seifert fibrations with $\chi=0$ and $e=0$. Two of these manifolds are actually diffeomorphic, so we get six closed Seifert manifolds up to diffeomorphism, distinguished by their homology.}
\label{seven:table}
\end{center}
\end{table}

We will soon see that $\matRP^2 \timtil S^1$ is the unique non-prime Seifert manifold. We now turn to the $\chi=0$ case.

\begin{prop} \label{chi:zero:prop}
Every closed Seifert fibration with $\chi = 0$ and $e=0$ is isomorphic to one of the seven listed in Table \ref{seven:table}. These seven manifolds are all pairwise non-diffeomorphic, except
$$\big(S^2,(2,1),(2,1),(2,-1),(2,-1)\big) \isom K\timtil S^1.$$
\end{prop}
\begin{proof}
The closed orbifolds with $\chi =0$ are 
$$T,\ K,\ (\matRP^2, 2,2),\ (S^2, 2,2,2,2),\ (S^2,2,3,6),\ (S^2,3,3,3),\ (S^2,2,4,4).$$ 
It is easy to show that by imposing $e=0$ we get the fibrations listed.
The homology calculation is an easy exercise and luckily distinguishes all the manifolds except (of course) the two diffeomorphic ones (see Corollary \ref{Seifert:diffeo:cor}).
\end{proof}

\begin{table}
\begin{center}
\begin{tabular}{c||ccc}
       & \phantom{\Big|} $\chi>0$ & $\chi = 0$ & $\chi<0$ \\
 \hline \hline
 $e=0$ & \phantom{\Big|} $2$ & $6$ & $\infty$ \\
 $e\neq 0$ & \phantom{\big|} $\infty$ & $\infty$ & $\infty$ \\
 $\partial M \neq \emptyset$  & \phantom{\Big|} $1$ & $2$ & $\infty$
\end{tabular}
\vspace{.2 cm}
\nota{The number of Seifert manifolds in each commensurability class.}
\label{finite:table}
\end{center}
\end{table}

Summing up, there are four commensurability classes containing finitely many manifolds, and their number is shown in Table \ref{finite:table}.

\subsection{Universal cover}
We now determine the universal cover of all the closed Seifert manifolds. Nothing strange happens: we either get $S^3$, $S^2\times \matR$, or $\matR^3$.

\begin{table}
\begin{center}
\begin{tabular}{c||ccc}
       & \phantom{\Big|} $\chi>0$ & $\chi = 0$ & $\chi<0$ \\
 \hline \hline
 $e=0$ & \phantom{\Big|} $S^2 \times \matR$ & $\matR^3$ & $\matR^3$ \\
 $e\neq 0$ & \phantom{\big|} $S^3$ & $\matR^3$ & $\matR^3$
\end{tabular}
\vspace{.2 cm}
\nota{The universal cover of a closed Seifert manifold depends on its invariants $e$ and $\chi$.}
\label{universal:table}
\end{center}
\end{table}

\begin{prop} \label{universal:Seifert:prop}
The universal cover of a closed Seifert manifold is shown in Table \ref{universal:table}.
\end{prop}
\begin{proof}
The universal cover of a circle bundle over a surface $S$ with $\chi(S)\leqslant 0$
is a line bundle over the universal cover $\matR^2$ of $S$. The line bundle is trivial since $\matR^2$ is contractible and we get $\matR^2\times \matR = \matR^3$.
\end{proof}

\begin{cor}
Every Seifert manifold $M$ is either irreducible or covered by $S^2\times \matR$. In the latter case $M$ is diffeomorphic to $S^2\times S^1$ or $\matRP^2\timtil S^1$.
\end{cor}
\begin{proof}
If the universal cover is $S^3$ or $\matR^3$, it is irreducible and hence also $M$ is. When $M$ has boundary, we apply the proof of Proposition \ref{universal:Seifert:prop} to the interior of $M$ and find that its universal cover is $\matR^3$.
\end{proof}

Therefore every Seifert manifold is prime except $\matRP^2 \timtil S^1 = \matRP^3 \# \matRP^3$.

\subsection{Fibre-Parallel Dehn filling}
We have investigated various topological properties of Seifert manifolds, and we are now curious: what happens if we perform a forbidden fibre-parallel Dehn filling? We start by looking at the basic block $P\times S^1$, where $P$ is a pair-of-pants.

\begin{prop}
A fibre-parallel Dehn filling on $P\times S^1$ produces the connected sum of two solid tori; the fibres of $P\times S^1$ become the meridians of the solid tori.
\end{prop}
\begin{proof}
Let $\gamma \subset \partial P$ be the component whose torus $\gamma \times S^1$ is filled. Pick an essential arc $\alpha$ in $P$ with both endpoints in $\gamma$: the fibered annulus $A = \alpha \times S^1$ closes up to a two-sphere $S$ in the Dehn filling. 

The sphere $S$ separates the filled manifold into two portions, each diffeomorphic to $T\times [0,1]$ with a two-handle attached along a non-trivial curve: this is a holed solid torus. Every fibre now bounds a disc there.
\end{proof}

We now turn to arbitrary bundles.

\begin{ex}
Let $S$ be a compact surface with $b\geqslant 1$ boundary components.
The fibre-parallel Dehn filling of $S \timtil S^1$ is diffeomorphic to the connected sum of $b-1$ solid tori and $-\chi(S)+2-b$ copies of $S^2\times S^1$.
\end{ex}

We now deduce a more general fact: a fibre-parallel Dehn filling of a Seifert manifold is a connected sum of lens spaces and solid tori.

\begin{cor} \label{fiber:parallel:cor}
The fiber-parallel Dehn filling of
$$\big(S, (p_1,q_1), \ldots, (p_h, q_h)\big)$$
is diffeomorphic to the connected sum
$$L(p_1,q_1) \# \ldots \# L(p_h,q_h) \#_k (S^2 \times S^1) \#_{b-1} (D^2 \times S^1)$$
where $S$ has $b$ boundary components and $k=-\chi(S) +2-b$.
\end{cor}
If $S$ is orientable then $k$ is twice the genus of $S$. 

\begin{cor} \label{DF:Seifert:cor}
Every Dehn filling of a Seifert manifold is a connected sum of Seifert manifolds.
\end{cor}

If we have a Seifert manifold $M$ with many boundary components and we Dehn fill some of them, we are guaranteed to produce a new Seifert manifold, unless one of the Dehn filling kills a fiber-parallel slope: in that case the manifold ``degenerates'' to a connected sum of lens spaces and solid tori. This is an interesting phenomenon to keep in mind, because it will reproduce in Chapter \ref{Hyperbolic:Dehn:chapter} in the hyperbolic world: we will show that by Dehn filling a cusped hyperbolic three-manifold we always get a new hyperbolic manifold if we avoid a finite number of ``exceptional slopes'' on every boundary torus; if some of the exceptional slopes is employed, the manifold is possibly not hyperbolic, and it typically ``degenerates'' and breaks into simpler pieces.

\section{Classification}
In the previous section we have classified the Seifert fibrations up to isomorphism and the Seifert manifolds up to commensurability. We now want to complete our study by classifying Seifert manifolds up to diffeomorphism: our final achievement will be Theorem \ref{Seifert:teo} that determines precisely the Seifert manifolds that admit non-isomorphic fibrations.

The proof of Theorem \ref{Seifert:teo} is not straightforward: we will apply different techniques to different classes of Seifert manifolds. We start by studying the fundamental group of Seifert manifolds, showing in particular that it fits into a nice short exact sequence. Then we study essential surfaces: we show that these can always be isotoped to be either in ``vertical'' or ``horizontal'' position with respect to the fibration.

In ``most'' cases, a Seifert manifold contains many incompressible vertical tori and these can be used to characterise the manifold. When the Seifert manifold is ``small'' it contains no such tori, and we must use different techniques: for instance we look at its fundamental group, which may be finite or infinite.

\subsection{Fundamental group}
We study the fundamental group of Seifert manifolds. Recall that every closed Seifert fibration has two fundamental invariants $\chi$ and $e$ that determine the commensurability class of the Seifert manifold. Finite fundamental groups are easily detected:

\begin{prop}
A closed Seifert manifold has finite fundamental group $\Longleftrightarrow$ it is covered by $S^3$ $\Longleftrightarrow$ it has $\chi>0$ and $e\neq 0$.
\end{prop}
\begin{proof}
A compact manifold has finite fundamental group $\Longleftrightarrow$ its universal cover is also compact. Table \ref{universal:table} applies.
\end{proof}

Ordinary fibrations generate exact sequences in homotopy, and also Seifert fibrations do. Recall that a good orbifold has a well-defined fundamental group, and that all locally orientable surface orbifolds are good except $(S^2, p)$ and $(S^2, p_1, p_2)$ with $p_1\neq p_2$. We will ignore bad orbifolds, since the Seifert manifolds fibering over them are lens spaces.

\begin{prop}
Let $M \to S$ be a Seifert fibration over a good orbifold $S$. There is an exact sequence
$$1 \longrightarrow K \longrightarrow \pi_1(M) \longrightarrow \pi_1(S)\longrightarrow 1$$
where $K$ is the normal cyclic subgroup of $\pi_1(M)$ generated by a regular fiber and $\pi_1(S)$ is the orbifold fundamental group of $S$. 
\end{prop}
\begin{proof}
The universal cover $\tilde M$ of $M$ fibres over the universal cover $\tilde S$ of $S$. 
The fiber is a circle or a line, depending on whether $\pi_1(M)$ is finite or infinite, see Proposition \ref{universal:Seifert:prop}.

The group $\pi_1(M)$ acts fibre-preservingly on $\tilde M$ and hence acts also on $\tilde S$ as a covering automorphism for $\tilde S \to S$. This induces a natural homomorphism $\pi_1(M) \to \pi_1(S)$. Its kernel consists of all deck transformations of $\tilde M$ that fix the fibres: these are precisely $K$. The homomorphism is surjective because every loop in $S$ lifts to a loop in $M$.
\end{proof}

We now look more closely at the normal cyclic subgroup $K$.

\begin{prop} \label{infinite:K:prop}
The group $K$ is infinite if and only if $\pi_1(M)$ is.
\end{prop}
\begin{proof}
The group $K$ acts freely and proper discontinuously on each fiber of the universal cover $\tilde M$, and it quotients it to a circle. The fiber in $\tilde M$ is compact if and only if $\pi_1(M)$ is finite.
\end{proof}

\begin{prop} \label{central:prop}
The group $K$ lies in the centre of $\pi_1(M)$ if and only if $S$ is orientable or $K$ has order two.
\end{prop}
\begin{proof}
If $S$ is orientable, the manifold $M$ is a Dehn filling of a product $S' \times S^1$ where the fiber is obviously in the centre, and it remains so after the Dehn filling. If $S$ is non-orientable, pick an orientation-reversing loop $\alpha$ in $S$: there is a Klein bottle fibering above $\alpha$ and we get $\alpha g \alpha^{-1} = g^{-1}$ for the fiber $g\in K$. Therefore $K$ is central if and only if $g^2 = e$.
\end{proof}

In many cases the subgroup $K$ may be characterised intrinsically.

\begin{prop}
Let $M$ be a Seifert manifold whose fundamental group is not virtually abelian. The subgroup $K\triangleleft \pi_1(M)$ is the unique maximal cyclic normal subgroup of $\pi_1(M)$.
\end{prop}
\begin{proof}
Suppose by contradiction that $\pi_1(M)$ contains a cyclic normal subgroup $K'$ which is not contained in $K$. Therefore its image in $\pi_1(S)$ is a non-trivial cyclic normal subgroup. 

If $\chi(S)<0$ then $\interior S = \matH^2/_\Gamma$ for a discrete group $\Gamma=\pi_1(S)$ of isometries of $\matH^2$. Corollary \ref{normal:abelian:cor} implies that $\Gamma$ does not contain non-trivial cyclic normal subgroups.

If $\chi(S) = 0$ then $\interior S=\matR^2/_\Gamma$ for a discrete group $\Gamma = \pi_1(S)$ of isometries of $\matR^2$, which in fact may contain cyclic normal subgroups: so we need a more careful analysis. Since $\pi_1(M)$ is not virtually abelian, the surface $S$ is closed and $e\neq 0$.

The image of $K'$ in $\pi_1(S)$ is non-trivial, normal, and cyclic, and is hence infinite (finite cyclic groups are generated by rotations and cannot be normal). Therefore $K'$ is infinite and it intersects non-trivially every finite-index subgroup of $\pi_1(M)$: up to substituting $M$ with a finite cover we may suppose that $S$ is a torus. 

The image of $K'$ is a non-trivial subgroup of $\pi_1(S) =\matZ \times \matZ$. 
Pick two elements $a' \in K', b \in \pi_1(M)$ whose images in $\pi_1(S)$ generate a finite-index subgroup of $\pi_1(S)$. Pick a generator $a\in K$. The three elements $a,a',b$ generate a finite-index subgroup of $\pi_1(M)$. The element $a$ is central by Proposition \ref{central:prop}. We have $b^{-1}a'b = (a')^{\pm 1}$ because $K' = \matZ $ is normal.  The elements $a,a',b^2$ commute and generate a finite-index abelian $\matZ^3<\pi_1(M)$: a contradiction.
\end{proof}

\begin{cor} \label{center:cor}
Let $M$ be a Seifert manifold whose fundamental group is not virtually abelian. The centre of $\pi_1(M)$ is $K$ if $S$ is orientable and trivial otherwise.
\end{cor}
\begin{proof}
Use Propositions \ref{central:prop} and \ref{infinite:K:prop}.
\end{proof}

We have collected some important information on the fundamental group of Seifert manifolds. We now move from algebra to topology and study the essential surfaces in Seifert manifolds: we will prove a theorem analogous to Proposition \ref{line:bundles:prop}, namely that every essential surface in a Seifert manifold is either horizontal or vertical.

\subsection{Horizontal and vertical surfaces}
We now want to study how Seifert manifolds may contain interesting surfaces.
Let $M \to S$ be a Seifert fibration. A properly embedded surface $\Sigma \subset M$ is
\begin{itemize}
\item \emph{vertical} if it is a union of some regular fibres,
\item \emph{horizontal} if it is transverse to all fibres.
\end{itemize}

If $\Sigma$ is vertical, it is either an annulus, a torus, or a Klein bottle, projecting respectively to an arc, an orientation-preserving, or an orientation-reversing simple closed curve that avoids the cone points. Vertical surfaces are in 1-1 correspondence with 1-dimensional objects in $S$ and are thus easily determined.

Horizontal surfaces are more subtle. We first note that if $\Sigma$ is horizontal, the natural projection $\Sigma \to S$ is an orbifold covering. When does $M$ contain a horizontal surface? It certainly does when $M$ has boundary.

\begin{figure}
\begin{center}
\includegraphics[width = 5 cm] {\iftoggle{BW}{Euler_zero-BW}{Euler_zero}}
\nota{We cut $S$ along $k$ (\iftoggle{BW}{grey}{blue}) arcs to get $S=S_0\cup D_1\cup \ldots \cup D_k$ with $D_i$ a disc containing the (\iftoggle{BW}{grey}{red}) cone point $p_i$: here $k=3$. }
\label{Euler_zero:fig}
\end{center}
\end{figure}

\begin{prop}
Every Seifert fibration $M\to S$ with boundary contains a horizontal surface.
\end{prop}
\begin{proof}
We have $M=\big(S,(p_1,q_1),\ldots,(p_k,q_k)\big)$. We cut $S$ along 
$k$ arcs as in Figure \ref{Euler_zero:fig} to get $S=S_0\cup D_1\cup\ldots \cup D_k$ with $D_i$ a disc containing the cone point $p_i$. Let $p$ be a common multiple of $p_1,\ldots, p_k$. The bundle over $S_0$ has no singular fibres and is hence $S_0\timforsetil S^1$. The bundle $S_0\timforsetil S^1$ contains a (possibly disconnected) horizontal surface $\Sigma$ intersecting the fibres in any fixed number $n$ of points (exercise) and we pick one $\Sigma$ with $n=p$.

For every $i=1,\ldots, k$ there is a standard fibered solid torus lying above $D_i$, attached to $S_0 \timforsetil S^1$ via a vertical annulus $A$. The meridian of this solid torus intersects $A$ into $p_i$ horizontal segments. Pick $\frac{p}{p_i}$ parallel meridians: both $\Sigma$ and these meridians intersect $A$ in $p$ horizontal segments and can hence be glued to form a horizontal surface for $M$.
\end{proof}

On closed Seifert manifolds the existence of a horizontal surface is fully detected by the Euler number.

\begin{prop}
A closed Seifert fibration $M \to S$ contains a horizontal surface if and only if $e=0$.
\end{prop}
\begin{proof}
Suppose $M\to S$ contains a horizontal surface $\Sigma$. We pull-back the fibration along the orbifold covering $\Sigma \to S$ to get a $\tilde M \to \Sigma$. The horizontal surface $\Sigma$ lifts to a section of $\tilde M \to \Sigma$. Since the fibration has a section, we get $e(\tilde M) = 0$ and hence $e(M) = 0 $.

Conversely, suppose $e(M)=0$. Drill one open fibered solid torus from $M$: we know by the previous proposition that the resulting fibration has a horizontal surface $\Sigma$. Its boundary consists of parallel curves of some type $(p,q)$ with $p>0$. If we Dehn-fill the manifold by killing these curves we get another Seifert manifold $M'$, to which the fibration and the section extend: therefore $e(M')=0$. At most one pair $(p,q)$ may produce a manifold with $e=0$, and hence $M=M'$.
\end{proof}

\subsection{Essential surfaces}
Recall the definition of essential surface from Sections \ref{essential:discs:subsection} and \ref{essential:subsection}. We now show that essential surfaces in irreducible Seifert manifolds are either vertical or horizontal. 

\begin{prop} \label{v:h:prop}
Let $M\to S$ be a Seifert fibration and $M$ be irreducible. Every essential surface $\Sigma$ is isotopic to a vertical or horizontal surface.
\end{prop}
\begin{proof}
See $M$ as a Dehn filling of $S\timforsetil S^1$. We suppose that $\Sigma$ intersects transversely the cores of the filling solid tori in the minimum number of points up to isotopy. After an isotopy $\Sigma$ intersects the filling solid tori into parallel horizontal discs.

Decompose $S$ into one 0-handle and some $g$ one-handles. We see the 0-handle as a $2g$-gon and each 1-handle as a rectangle. Above each edge of the $2g$-gon or of a rectangle there is a vertical annulus $A$. We put $\Sigma$ in transverse position with respect to these vertical annuli and up to isotopy we may suppose that $\Sigma\cap A$ consists either of vertical fibres or horizontal arcs for each $A$. Indeed we can easily eliminate trivial circles (because $\Sigma$ is incompressible and $M$ is irreducible) and arcs forming bigons with $\partial A$ (because $\Sigma$ is $\partial$-irreducible and has minimal intersection with the solid tori cores).

Above every polygonal handle of $S$ there is a prismatic solid torus $W$. The closed curves $\Sigma \cap \partial W$ are made of horizontal or vertical arcs and are hence essential, forming some parallel slopes in $\partial W$. We can suppose that $\Sigma\cap W$ consists of essential discs or incompressible surfaces (if there is a compressing disc $D$ inside $W$, then $\partial D\subset \Sigma$ bounds a disc in $\Sigma$ which we isotope to $D$ reducing the intersection with the vertical annuli). 
We can also suppose that $\Sigma\cap W$ is $\partial$-incompressible with respect to every vertical annulus $A\subset \partial W$, in the sense that there is no $\partial$-compressing disc $D$ with $\partial D\subset A \cup \Sigma$ (otherwise we isotope $\Sigma$ and reduce intersections).

By Proposition \ref{incompressible:solid:torus:prop} the surface $\Sigma\cap W$ consists either of essential discs or $\partial$-parallel annuli (not both). In the first case $\Sigma\cap \partial W$ consists of horizontal curves and $\Sigma\cap W$ consists of horizontal discs. In the second case, since the annuli are $\partial$-incompressible with respect to the vertical annuli in $\partial W$, their boundaries must be vertical circles lying in distinct vertical annuli. Up to isotopy $\Sigma\cap W$ consists of vertical annuli.

We have decomposed $\Sigma$ into horizontal discs or vertical annuli, and both cannot coexist. Hence $\Sigma$ is either horizontal or vertical.
\end{proof}

\begin{cor} \label{partial:cor}
Every Seifert manifold is irreducible and $\partial$-irreducible, except $S^2\times S^1$, $\matRP^2 \timtil S^1$, and $D\times S^1$.
\end{cor}
\begin{proof}
We already know it is irreducible with the first two exceptions. If it contains an essential disc $D$, it contains a horizontal one which covers the base surface $S$ of the fibration, hence $\chi (S) >0$ and we get $D\times S^1$.
\end{proof}

\begin{figure}
\begin{center}
\includegraphics[width = 10 cm] {\iftoggle{BW}{vertical_incompressible-BW}{vertical_incompressible}}
\nota{A vertical torus or annulus is essential, unless its projection is: the boundary of a disc containing zero (1) or one (2) singular cone points, an arc parallel to the boundary (3), or a $\partial$-parallel closed curve (4).}
\label{vertical_incompressible:fig}
\end{center}
\end{figure}

We now prove a converse to Proposition \ref{v:h:prop}. 

\begin{prop} \label{incompressible:vertical:prop}
Let $M\to S$ be a Seifert fibration and $M$ be irreducible. Let $\Sigma \subset M $ be an orientable connected surface. Suppose that
\begin{itemize}
\item $\Sigma$ is horizontal, or
\item $\Sigma$ is vertical and its projection is not as in Figure \ref{vertical_incompressible:fig}.
\end{itemize}
Then $\Sigma$ is essential.
\end{prop}
\begin{proof}
If $\Sigma$ is horizontal, it finitely covers $S$ and hence $\pi_1(\Sigma)$ injects in $\pi_1(S)$. Therefore it injects also in $\pi_1(M)$ and $\Sigma$ is incompressible. By doubling everything along $\partial M$ we get a horizontal $D\Sigma$ inside $DM$, which must also be incompressible: this implies that $\Sigma$ is $\partial$-incompressible. It is also clearly not $\partial$-parallel, so it is essential.

If $\Sigma$ is vertical, by cutting along it we get one or two Seifert fibrations. Since $\Sigma$ is not as in Figure \ref{Euler_zero:fig}-(right), the base orbifolds of these fibrations have $\chi\leqslant 0$, so their boundary is incompressible by Corollary \ref{partial:cor}. Therefore $\Sigma$ is incompressible and $\partial$-incompressible. Moreover it is not a $\partial$-parallel torus because none of the cut Seifert manifold is $A\times S^1$.
\end{proof}

\begin{figure}
\begin{center}
\includegraphics[width = 5 cm] {\iftoggle{BW}{vertical_examples-BW}{vertical_examples}}
\nota{Some essential vertical annuli and tori in a Seifert manifold fibering over a disc with 2 or 3 singular points.}
\label{vertical_examples:fig}
\end{center}
\end{figure}

\begin{example}
The curves in Figure \ref{vertical_examples:fig} determine vertical essential annuli and tori in Seifert manifolds fibering over the disc with 2 or 3 singular fibres.
\end{example}

We can now fully detect essential tori in almost all Seifert manifolds.

\begin{cor} \label{vertical:tori:cor}
Let $M\to S$ be a Seifert fibration and $M$ be irreducible. If $M$ is not covered by a 3-torus or $T\times [0,1]$, every essential torus or annulus is vertical.
\end{cor}
\begin{proof}
Horizontal tori or annuli may arise only when $\chi=0$ and $\partial M\neq 0$ or $e=0$.
\end{proof}

\subsection{Simple Seifert manifolds} \label{simple:Seifert:subsection}
We want to classify Seifert manifolds up to diffeomorphism. To do so, we group them into some classes, and use different techniques on each class. Recall that a 3-manifold is \emph{simple} if it contains no essential sphere, disc, annulus, or torus. 

\begin{prop} \label{simple:Seifert:prop}
A Seifert manifold is simple $\Longleftrightarrow$ it fibres over $S^2$ with at most three singular fibres and is not covered by $S^2\times S^1$ or the 3-torus. 
\end{prop}
\begin{proof}
Let $M\to S$ be a Seifert fibration. The manifold $M$ contains an essential vertical torus or annulus $\Longleftrightarrow$ the orbifold $S$ contains an arc or an orientation-preserving simple closed curve that is not as in Figure \ref{vertical_incompressible:fig}. The only orbifolds that do not contain such arcs or curves are: $D$ with at most one fiber, $S^2$ with at most three fibres, and $\matRP^2$ with at most one fiber. In the first case $M$ is a solid torus, which contains an essential disc, and in the third case $M$ also fibres over $S^2$ with at most three singular fibres, see Corollary \ref{Seifert:diffeo:cor}. 

We are left to consider horizontal closed essential surfaces with $\chi \geqslant 0$. These arise only when $e=0$ and $\chi(S)\geqslant 0$, \emph{i.e.}~when $M$ is covered by $S^2\times S^1$ or the 3-torus.
\end{proof}

We note in particular that all the simple Seifert manifolds are closed.

\begin{rem}
Among the nine commensurability classes of Seifert manifolds, three contain simple manifolds: those with empty boundary and $e\neq 0$. The Euler characteristic $\chi$ of a sphere with 3 singular fibres may in fact be positive, null, or negative. In particular there are simple Seifert manifolds with finite and with infinite fundamental group.
\end{rem}

\subsection{Seifert manifolds with boundary}
To classify Seifert manifolds up to diffeomorphism, we will show that (except a few explicit exceptions) a generic Seifert manifold has a unique fibration up to isomorphism. We start with the easier non-empty-boundary case.

\begin{prop} \label{except:non-empty:prop}
Every Seifert manifold $M$ with non-empty boundary admits only one fibration up to isomorphism, except in the following cases:
\begin{itemize}
\item $M = D \times S^1$ fibres as $\big( D, (p,q) \big)$,
\item $M = S \timtil S^1$ fibres as $\big(D, (2,1), (2,-1) \big)$.
\end{itemize}
Here $S$ is the M\"obius strip.
\end{prop}
\begin{proof}
If $M$ fibres over an orbifold $S$ with $\chi(S)\geqslant 0$, we have already proved this in Section \ref{finite:classes:prop}. Suppose that $M$ has two fibrations $M\to S$ and $M\to S'$ with $\chi(S), \chi(S')< 0$. 

Pick a minimal collection of properly embedded arcs in $S$ that avoid the cone points and decompose $S$ into discs, each containing at most one cone point. Each arc determines an essential vertical annulus in $M\to S$. The complement of these annuli consists of vertical fibered solid tori.

By Proposition \ref{incompressible:vertical:prop} these annuli are essential and by  Proposition \ref{v:h:prop} they are isotopic to vertical annuli with respect to the other fibration $M\to S'$ (they cannot be horizontal since $\chi(S')<0$). Therefore we may isotope the fibration $M\to S'$ so that the annuli are vertical in both fibrations, and by further isotoping we can in fact suppose that the two fibrations coincide on a neighbourhood of these annuli. 

The complement of this neighbourhood consists of solid tori. The two fibrations are obtained from the same fibration by Dehn-filling along the same slopes, and hence they are isomorphic.
\end{proof}

\subsection{Seifert manifolds with infinite fundamental group}
We now turn to closed Seifert manifolds. We start by examining the $\chi\leqslant 0$ case.

\begin{prop}
Every closed Seifert manifold $M$ not covered by $S^3$ and $S^2\times S^1$ admits only one fibration up to isomorphism, except:
$$\big(S^2,(2,1),(2,1),(2,-1),(2,-1)\big) \isom K\timtil S^1. $$
\end{prop}
\begin{proof}
If the Seifert manifold $M$ is covered by the 3-torus, we have already proved this in Proposition \ref{chi:zero:prop}. Suppose that $M$ is not covered by $S^3$, $S^2\times S^1$, or the 3-torus. Let $M$ have two fibrations $M\to S$ and $M\to S'$. 

We first suppose that $M$ contains an essential torus. We try to proceed as in the proof of Proposition \ref{except:non-empty:prop} using vertical tori instead of annuli. Essential tori are vertical to both fibrations by Corollary \ref{vertical:tori:cor}. In particular both $S$ and $S'$ are not spheres with at most three singular points.

Let $\dot S\subset S$ be $S$ without its singular points. Pick two multicurves $C_1,C_2 \subset \dot S$ without puncture-parallel components, that fill $\dot S$ and intersect minimally. By filling $\dot S$ we mean that $\dot S \setminus (C_1\cup C_2)$ consists of discs or once-punctured discs, and minimality implies that no such disc is an unpunctured bigon.

Let $T_1, T_2$ be the collections of tori fibering above $C_1, C_2$. The set $T_1$ consists of disjoint essential tori, so up to isotopy they are vertical also with respect to $M\to S'$. The two fibrations now coincide on a neighbourhood of $T_1$. 

We turn to the tori $T_2$, that are vertical with respect to $M\to S$. Up to isotopy, each simple closed curve in $T_1\cap T_2$ is either horizontal or vertical with respect to $M\to S'$. Each annulus in $T_2\setminus T_1$ is essential in $M\setminus T_1$ and hence it is correspondingly horizontal or vertical after an isotopy. Horizontal annuli glue to a horizontal torus in $T_2$, which is excluded: hence the tori in $T_2$ are also vertical with respect to $M\to S'$. The two fibrations now coincide on a neighbourhood of $T_1\cup T_2$ after an isotopy.

The complement $M\setminus (T_1\cup T_2)$ consists of vertical solid tori and we conclude as in Proposition \ref{except:non-empty:prop} that the two fibrations are isomorphic.

We are left to consider the case where $M$ contains no essential tori, and hence is a simple manifold. Both $S,S'$ are spheres with exactly three singular points by Proposition \ref{simple:Seifert:prop} (not less than three singular points since $M$ is not covered by $S^3$ or $S^2\times S^1$). 

By Corollary \ref{center:cor}, the fundamental groups $\pi_1(S)$ and $\pi_1(S')$ are both isomorphic to $\pi_1(M)$ quotiented by its centre. Exercise \ref{VD:ex} implies that the two orbifolds are isomorphic, so $S=S' = (S^2, p_1,p_2,p_3)$ for some $p_1, p_2, p_3 \geqslant 2$. 
Consider the exact sequence
$$1 \longrightarrow K \longrightarrow \pi_1(M) \longrightarrow \pi_1(S)\longrightarrow 1.$$
Fix a generator $l\in K$ and coherently an orientation for the fibres of both fibrations $M\to S$ and $M\to S'$. 
The Von Dyck group $\pi_1(S^2,p_1,p_2,p_3)$ has a presentation 
$$\langle\ r_1, r_2, r_3\ |\ r_1^{p_1}, r_2^{p_2}, r_3^{p_3}, r_1r_2r_3\ \rangle$$
and the three generators $r_1, r_2, r_3$ are intrinsically determined up to simultaneous conjugation or inversion, see Exercise \ref{VD:ex}. Fix three lifts $m_1, m_2, m_3$ of $r_1, r_2, r_3$ in $\pi_1(M)$ with $m_1m_2m_3 = 1$: these lifts determine sections for the fibrations $M\to S$ and $M\to S'$ with boundary meridians $m_1,m_2,m_3$. We use these sections to determine the parameters $q_i$ and $q_i'$ in both fibrations, for $i=1,2,3$.

The centre $K$ is infinite since $\pi_1(M)$ is, see Proposition \ref{infinite:K:prop}. Therefore $m_i^{p_i} \in K$ equals $l^{q_i''}$ for some unique $q_i''\in \matZ$. By construction we have $q_i'' = q_i$ and $q_i'' = q_i'$, hence $q_i=q_i'$. The two fibrations $M\to S$ and $M\to S'$ are isomorphic.
\end{proof}

\subsection{Seifert manifolds with finite fundamental group} \label{Seifert:finite:subsection}
This long journey through Seifert manifolds is almost finished: it remains to classify the Seifert manifolds with finite fundamental group, that is that are covered by the three-sphere.

\begin{table}
\begin{center}
\begin{tabular}{c||ccccc}
\phantom{\Big|} fibration        & $e(M)$ & $|\pi_1(S)|$ & $|\pi_1(M)|$ & $|H_1(M,\matZ)|$ \\
 \hline \hline
\phantom{\Big|} $\big(S^2,(2,1),(2,1),(p,q)\big)$ & $\frac {p+q}p$ & $2p$ & $4p|p+q|$ & $4|p+q|$ \\
\phantom{\Big|} $\big(S^2,(2,1),(3,1),(3,q)\big)$ & $\frac{5+2q}6$ & $12$ & $24|5+2q|$ & $3|5+2q|$ \\
\phantom{\Big|} $\big(S^2,(2,1),(3,1),(4,q)\big)$ & $\frac{10+3q}{12}$ & $24$ & $48|10+3q|$ & $2|10+3q|$ \\
\phantom{\Big|} $\big(S^2,(2,1),(3,1),(5,q)\big)$ & $\frac{25+6q}{30}$ & $60$ & $120|25+6q|$ & $|25+6q|$
\end{tabular}
\vspace{.2 cm}
\nota{The non-lens Seifert manifolds $M$ with finite fundamental group. For each fixed base orbifold $S$, the fibration is determined by the order of $H_1(M,\matZ)$: when two different parameters $q$ give the same order of $H_1(M,\matZ)$, the fibrations are actually the same (with opposite orientations and hence Euler numbers).}
\label{elliptic:table}
\end{center}
\end{table}

\begin{prop} \label{elliptic:center:prop}
Let $M$ be a Seifert manifold covered by $S^3$ that is not a lens space. It has a unique fibration over one of the orbifolds
$$S = (S^2,2,2,p), \quad (S^2,2,3,3), \quad (S^2,2,3,4), \quad (S^2,2,3,5)$$
for some $p\geqslant 2$. The centre of $M$ is the cyclic group $K$ in 
$$1 \longrightarrow K \longrightarrow \pi_1(M) \longrightarrow \pi_1(S)\longrightarrow 1.$$
The manifold $M$ is determined by $S$ and the order of $H_1(M)$, see Table \ref{elliptic:table}.
\end{prop}
\begin{proof}
We have $M\to S$ with $\chi(S)>0$ and hence $S$ is either a sphere with at most two singular points (so $M$ is a lens space), one of the orbifolds listed, or $\matRP^2$ with one singular point; in the latter case $M$ also fibres over $(S^2,2,2,p)$ by Corollary \ref{Seifert:diffeo:cor}.

The subgroup $K$ is central by Proposition \ref{central:prop}. Suppose that $\pi_1(M)$ contains a central element disjoint from $K$: its image in $\pi_1(S)$ is a non-trivial central element. However $S$ is spherical and $\pi_1(S)$ is a non-cyclic group of rotations of $S^2$: two rotations with different axis never commute.

The subgroup $K<\pi_1(M)$ is intrinsically determined as the centre, hence the quotient $\pi_1(S)$ is also determined, and $S$ also is by
Exercise \ref{VD:ex}.
Therefore $M$ cannot fiber on two different orbifolds in the list.

To conclude, we prove that the fibration $M\to S$ is determined by $S$ and the finite number $|H_1(M,\matZ)|$. We analyse each case separately. First, we transform the fibration as in Table \ref{elliptic:table}: using the move from Proposition \ref{S:moves:prop}-(\ref{Seifert:1:eqn}) we manage to transform each $(2,q')$ in $(2,1)$ and $(3,q')$ in either $(3,1)$ or $(3,2)$, and by reversing orientation (see Remark \ref{reverse:rem}) we transform $(3,2)$ into $(3,-2)$ and hence into $(3,1)$.

Using Proposition \ref{Seifert:homology:prop} we see that $|H_1(M)|$ is as in Table \ref{elliptic:table}. The only cases where different values of $q$ give the same $|H_1(M)|$ are the following:
\begin{itemize}
\item $S=(S^2,2,2,p)$ and $|p+q|=|p+q'|$, hence $q'=-q-2p$, 
\item $S=(S^2,2,3,3)$ and $|5+2q| = |5+2q'|$, hence $q' = -q-5$ which implies $q,q' \equiv 2 \ ({\rm mod}\ 3)$.
\end{itemize}
In these cases it is easy to verify that the fibrations with $q$ and $q'$ are isomorphic using Proposition \ref{S:moves:prop}-(\ref{Seifert:1:eqn}) and Remark \ref{reverse:rem}.

The order $|\pi_1(M)|$ in Table \ref{elliptic:table} is obtained using the formula
$$|\pi_1(M)| = \frac{4e(M)}{\chi(S)^2}$$
that we now prove. Since the orbifolds are good, the fibration is covered horizontally by a fibration over $S^2$, which is in turn covered vertically by $\big(S^2,(1,1)\big) = S^3$. The total degree is $d=d_{\rm h} \cdot d_{\rm v}$. Proposition \ref{e:prop} gives
$$d_{\rm h} = \frac{\chi(S^2)}{\chi(S)} = \frac 2{\chi(S)}, \quad 1=e(S^2,(1,1)) = \frac{d_{\rm h}}{d_{\rm v}} \cdot e(M).$$ 
Therefore $d_{\rm v} = d_{\rm h}\cdot e(M)$ and
$$d = d_{\rm h} \cdot d_{\rm v} = d_{\rm h}^2 \cdot e(M) = \frac{4e(M)}{\chi(S)^2}$$
as required.
\end{proof}

\subsection{Summary} \label{summary:Seifert:subsection}
We now summarise the topological classification of Seifert manifolds in a single statement. 

The Seifert fibrations $M\to S$ are fully classified by Proposition \ref{S:moves:prop} and Corollary \ref{classification:Seifert:cor}. The latter says that two Seifert fibrations
$$\big(S, (p_1,q_1), \ldots, (p_h, q_h)\big), \quad \big(S', (p_1',q_1'), \ldots, (p_{h'}', q_{h'}')\big)$$
with $p_i, p_i'\geqslant 2$ are orientation-preservingly isomorphic if and only if $S=S'$, $h=h'$, $e=e'$, and up to reordering $p_i=p_i'$ and $q_i \equiv q_i'$ (mod $p_i$) for all $i$.

We easily understand when two Seifert fibrations are isomorphic. To classify Seifert manifolds up to diffeomorphism it only remains to understand which Seifert manifolds can have non-isomorphic fibrations. A long discussion has shown the following. We write $S^2\times S^1$ as the lens space $L(0,1)$.

\begin{teo} \label{Seifert:teo}
Every Seifert manifold has a unique Seifert fibration up to isomorphism, except the following:
\begin{itemize}
\item $L(p,q)$ fibres over $S^2$ with $\leqslant 2$ singular points in many ways,
\item $D\times S^1$ fibres over $D$ with $\leqslant 1$ singular point in many ways,
\item $\big( D, (2,1), (2,-1)\big) \isom S\timtil S^1$,
\item $\big(S^2,(2,1),(2,-1),(p,q)\big) \isom \big(\matRP^2, (q,p)\big)$, 
\item $\big(S^2,(2,1),(2,1),(2,-1),(2,-1)\big) \isom K\timtil S^1. $
\end{itemize}
Here $S$ and $K$ are the M\"obius strip and the Klein bottle.
\end{teo}

\subsection{References}
The main sources that we have consulted for this long chapter are Hatcher \cite{H} and Scott \cite{S}. Some material can also be found in Fomenko -- Matveev \cite{FM}. Two classical references are Seifert's original paper \cite{Se} and Orlik \cite{O}.

%% file: Construction.tex
\chapter{Constructions of three-manifolds} \label{construction:chapter}

In the previous chapter we have introduced and fully classified an important family of three-manifolds called \emph{Seifert manifolds}, and we now address the following question: how can we construct more three-manifolds?

The most popular techniques employed for the construction and manipulation of three-manifolds are of cut-and-paste type: we build three-manifolds by gluing some blocks altogether, and we try to describe both the blocks and their gluing with some reasonable combinatorial formalism.

The choice of the right blocks is of course fundamental, and different choices lead to quite different environments.
The first reasonable option may be to use tetrahedra as blocks, and in that case we talk about \emph{triangulations} of three-manifolds: this construction has a strong combinatorial flavour and can be easily carried on by a computer. Other choices involve blocks without ``ridges'', that is manifolds with boundary: by using handlebodies we get \emph{Heegaard splittings}, with knot/link complements and solid tori we get \emph{Dehn surgery}, and with product manifolds $\Sigma \times [-1,1]$ we get \emph{surface bundles}. We introduce here all these topological constructions.

We end this chapter by showing that every prime three-manifold has a canonical decomposition along disjoint embedded tori, called the \emph{geometric decomposition} -- the reason for adopting this name will be evident in the next chapters.

\section{Heegaard splittings}
A \emph{Heegard splitting} is a decomposition of a closed three-manifold in two manifolds of the simple kinds, the handlebodies. 

\subsection{Definition}
The following proposition is quite surprising, because it shows that every closed orientable three-manifold decomposes into two pieces of a very simple type. 

\begin{prop}
Every closed orientable 3-manifold $M$ decomposes into two handlebodies of some genus $g$.
\end{prop}
\begin{proof}
The 3-manifold $M$ has a handle decomposition with 0-, 1-, 2-, and 3-handles. The 0- and 1-handles altogether form a handlebody. The 2- and 3-handles may be turned upside down to form a handle decomposition into 0- and 1-handles, so another handlebody.
\end{proof}

A decomposition of $M$ into two handlebodies is called a \emph{Heegaard splitting}.\index{Heegaard splitting} The two handlebodies have necessarily the same genus $g$, since their boundaries are glued together and are hence diffeomorphic surfaces.

\begin{defn}
The \emph{Heegaard genus} $g(M)$ of a closed orientable $M$ is the minimum genus of a Heegaard splitting for $M$.\index{Heegaard genus}
\end{defn}

\subsection{Examples}
The manifolds of genus zero and one are perfectly understood.

\begin{prop}
The 3-sphere has genus zero, lens spaces (except $S^3$) have genus one, and all the other closed orientable three-manifolds have genus at least two.
\end{prop}
\begin{proof}
By gluing two discs we get a 3-sphere, and by definition by gluing two solid tori we get a lens space.
\end{proof}

\begin{figure}
\begin{center}
\includegraphics[width = 8 cm] {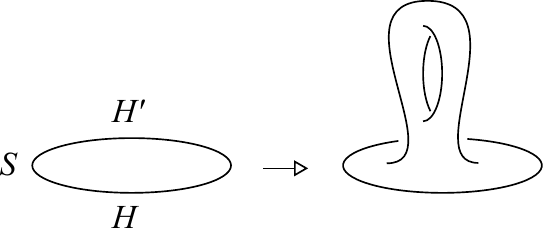}
\nota{A stabilisation of a Heegaard splitting.}
\label{stabilization:fig}
\end{center}
\end{figure}

Classifying closed manifolds of genus two is a much harder task. We limit ourselves to a class of examples.

\begin{ex}
Every Seifert manifold $\big(S^2,(p_1,q_1), (p_2,q_2), (p_3,q_3)\big)$ with $p_1,p_2,p_3\geqslant 2$ has genus two.
\end{ex}
\begin{proof}[Hint]
Pick two singular fibres and connect them with a horizontal arc. A regular neighbourhood of this graph is a genus-two handlebody: prove that its complement is also a handlebody.
\end{proof}

\subsection{Stabilisation}
The same manifold may have various non isotopic Heegaard splittings: for instance there is a simple move that modifies a Heegaard splitting by increasing its genus in a somehow trivial way.

A \emph{stabilisation} is a move as in Figure \ref{stabilization:fig} that transforms a Heegaard splitting $M=H\cup_S H'$ of genus $g$ into one of genus $g+1$ of the same manifold $M$. We add an unknotted handle to the surface $S=\partial H = \partial H'$, so that both $H$ and $H'$ are transformed into handlebodies of genus $g+1$. 

\begin{example} \label{S3:example}
The complement of a standardly embedded genus-$g$ handlebody in $S^3$ as in Figure \ref{handlebody:fig} is another handlebody, and together they form a genus-$g$ Heegaard splitting of $S^3$ obtained by stabilising $g$ times the genus-0 Heegaard splitting of $S^3$.
\end{example}

\subsection{Triangulations}
Smooth triangulations are somehow related to Heegaard splittings, at least in one direction. 

Every compact manifold admits a smooth triangulation (see Section \ref{triangulations:subsection}) and in dimension three we can luckily prove a converse statement, much as we did in dimension two (in Section \ref{surface:triangulations:subsection}).

Let $X$ be a three-dimensional pure simplicial complex, where: every face is incident to two tetrahedra, every edge is contained in a cycle of adjacent tetrahedra, and the link of every vertex is a sphere.

\begin{prop} \label{X:triangulation:prop}
The complex $X$ is the smooth triangulation of a closed three-manifold $M$, unique up to diffeomorphism.
\end{prop}
\begin{proof}
By dualising $X$ we get a handle decomposition: tetrahedra, triangles, edges, and vertices determine 0-, 1-, 2-, and 3-handles. This procedure constructs a smooth closed three-manifold $M$ triangulated by $X$. The way handles are attached is determined up to isotopy, and hence $S$ is determined up to diffeomorphism.
\end{proof}

As we did with surfaces, it is worth noting that this procedure (getting a unique smooth structure from a simplicial complex) does not work in all dimensions (here we used implicitly the non-obvious Proposition \ref{cap:prop}).

Topologists usually prefer to loosen the notion of ``triangulation'' by allowing self- and multiple adjacencies between tetrahedra, see Section \ref{ideal:subsection}. What is important to note here is that, no matter what the definition of ``triangulation'' is, a triangulated three-manifold always has a well-defined smooth structure determined only by the combinatorics of the triangulation. This shows that smooth three-manifolds can be treated combinatorially, for instance by a computer.

A relation between triangulations and Heegaard splittings is the following: a smooth triangulation of a closed manifold with $t$ tetrahedra gives rise to a dual handle decomposition and hence to a Heegaard splitting of genus $t$. A fundamental difference between the two notions is that there is a bounded number of closed three-manifolds triangulated by at most $t$ tetrahedra for every $t$, while there are infinitely many manifolds with Heegaard genus at most $g$ for every $g\geqslant 1$.

\section{Knots and links}
Knots and links are fundamental and beautiful objects in geometric topology. Knots and links in $S^3$ have a combinatorial and mildly two-dimensional nature because they can be treated as planar diagrams, but they should be considered as intrinsically three-dimensional objects.

\subsection{Definition}
A \emph{link} in a 3-manifold $M$ is a compact submanifold of dimension one. Being compact, it consists of finitely many circles, and a connected link is called a \emph{knot}. Links and knots are usually considered up to ambient isotopy.\index{link}\index{knot}

\begin{figure}
\begin{center}
\includegraphics [width = 4 cm] {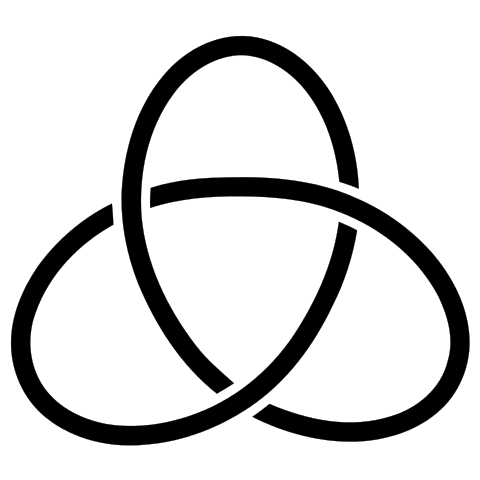}
\quad
\includegraphics [width = 4 cm] {\iftoggle{BW}{Blue_Trefoil_Knot-BW}{512px-Blue_Trefoil_Knot}}
\nota{Every knot may be described via a planar diagram with over/under crossings (left) and its tubular neighbourhood is a knotted solid torus (right).}
\label{trefoil_knots:fig}
\end{center}
\end{figure}

Every link in $S^3$ can be projected to a plane: if the link is generic with respect to the projection, its image is a \emph{diagram} as in Figure \ref{trefoil_knots:fig}-(left). Every link in $S^3$ can be described via a diagram and can be thus treated combinatorially, but it is important to note that the same link can be represented via infinitely many diagrams and it is hard in general to tell whether two given diagrams represent the same link. The \emph{trivial knot} (or \emph{unknot}) $K\subset S^3$ is the one that has a planar diagram without crossings.\index{unknot}

\subsection{Chirality}
The \emph{mirror image} of a knot or link is obtained by reflecting it with respect to any plane in $S^3$. On a diagram, this operation is realised simply by inverting all crossings simultaneously.

A knot is \emph{chiral} if it is not isotopic to its mirror image, and \emph{achiral} otherwise. For instance, the trivial knot is achiral, but the trefoil knot (shown in Figure \ref{trefoil_knots:fig}) is chiral.\index{knot!trefoil knot}
As in the rest of this book, we mostly ignore orientation issues and often consider implicitly two mirrored knots or links as equivalent.\index{knot!achiral and chiral knot}

\subsection{Link complement}
The tubular neighbourhood of a link $L$ in an orientable $M$ consists of solid tori by Proposition \ref{S1:bundle:prop}, and the \emph{link complement} of $L$ is the three-manifold with boundary obtained from $M$ by removing the interiors of these solid tori. If $M$ is compact, the link complement also is.\index{link complement}

\begin{prop}
Let $K\subset S^3$ be a knot. The complement of $K$ is an irreducible manifold. Moreover, the following facts are equivalent:
\begin{enumerate}
\item $K$ is trivial,
\item $K$ bounds a disc in $S^3$,
\item the complement of $K$ is a solid torus.
\end{enumerate}
\end{prop}
\begin{proof}
The complement of $K$ is irreducible by Proposition \ref{connected:boundary:prop}.
The implication (1)$\Rightarrow$(2) is the smooth Jordan curve theorem, and (2)$\Rightarrow$(1) holds because all discs in a connected three-manifold are isotopic. 

A trivial knot thickens to a standardly embedded solid torus, whose complement is a solid torus: hence (1)$\Rightarrow$(3). Conversely, if the complement of $K$ is a solid torus then $S^3$ decomposes into two solid tori as $S^3=L(1,q)$. The meridian of the complement solid torus extends to a disc with boundary $K$, hence (3)$\Rightarrow$(2).
\end{proof}

\begin{cor} \label{H:tori:cor}
All Heegaard tori for $S^3$ are isotopic
\end{cor}
\begin{proof}
A Heegaard torus $T\subset S^3$ is by definition a torus that decomposes $S^3$ into two solid tori. By the proposition the core of one of these solid tori is always isotopic to a trivial knot and hence $T$ is isotopic to a standard torus.
\end{proof}

A link complement in $S^3$ can be reducible: this holds precisely when there is a sphere $\Sigma \subset S^3$ disjoint from the link $L$ that cut $S^3$ into two balls, each containing some components of $L$. In that case we say that the link is  \emph{split}.\index{link!split link}

\begin{cor}
Let $L\subset S^3$ be a non-trivial knot or non-split link. The link complement is Haken.
\end{cor}
\begin{proof}
The link complement $M$ is irreducible by assumption; by Proposition \ref{irreducible:solid:torus:prop}, if $M$ contains an essential disc then it is a solid torus and hence $L$ is the trivial knot, which is excluded. Now $M$ is Haken by Corollary \ref{boundary:implies:Haken:cor}.
\end{proof}

\subsection{Prime knots}
There is an operation on knots in $S^3$ called \emph{connected sum}, similar to the one on manifolds and described in Figure \ref{sum_of_knots:fig}. 
A connected sum is \emph{trivial} if one of the two knots is trivial; a non-trivial knot is \emph{prime} if it is not the result of a non-trivial connected sum, and \emph{composite} otherwise.\index{knot!prime and composite knot}

\begin{figure}
\begin{center}
\includegraphics [width = 12.5 cm] {\iftoggle{BW}{Sum_of_knots-BW}{Sum_of_knots}}
\nota{Connected sum of knots: put two knots in disjoint balls and connect them with a band as shown.}
\label{sum_of_knots:fig}
\end{center}
\end{figure}

\begin{figure}
\begin{center}
\includegraphics [width = 10 cm] {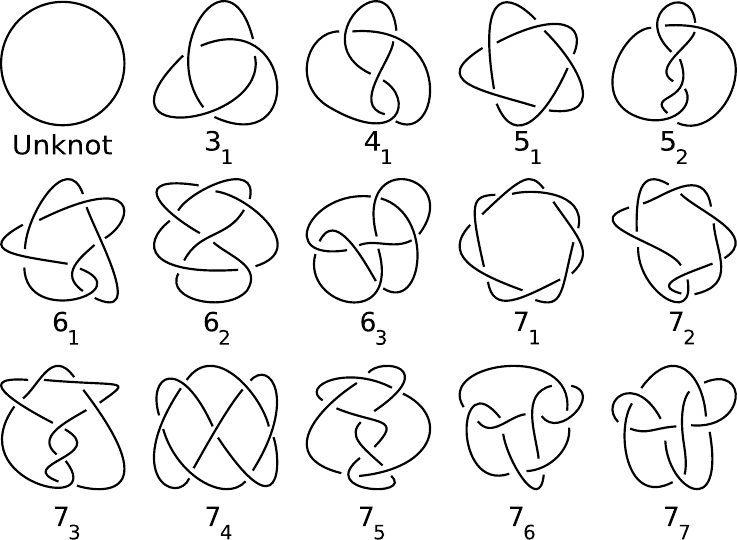}
\nota{Prime knots in $S^3$ with at most $7$ crossings. The knots $3_1$ and $4_1$ are the \emph{trefoil} and the \emph{figure-eight} knots respectively.}
\label{Knot_table:fig}
\end{center}
\end{figure}

The \emph{crossing number} of a knot is the minimum number of crossings in a diagram describing it.\index{knot!crossing number} Prime knots with small crossing number have been tabulated since the XIX century: at present all the prime knots with up to 16 crossing have been classified (by Hoste -- Thistlethwaite -- Weeks \cite{HTW} in 1998) and the first 14 numbers of them are listed in Table~\ref{prime:knots:table}. The prime knots with $\leqslant 7$ crossings are shown in Figure \ref{Knot_table:fig}: the first three in the list are the unknot, the \emph{trefoil}, and the \emph{figure-eight knot}.\index{knot!figure-eight knot} 

\begin{table}
\begin{center}
\begin{tabular}{c||ccccccccccccccc}
 $c$       & \phantom{\Big|}\!\! $0$ & $1$ & $2$ & $3$ & $4$ & $5$ & $6$ & $7$ & $8$
 & $9$ & $10$ & $11$ & $12$ & $13$ & $14$ \\
 \hline \hline
$n$ & 
\phantom{\Big|}\!\!
1 & 0 & 0 & 1 & 1 & 2 & 3 & 7 & 21 & 49 & 165 & 552 & 2176 & 9988 & 46972 
\end{tabular}
\vspace{.2 cm}
\nota{The number $n$ of prime knots with $c$ crossings, for all $c\leqslant 14$.}
\label{prime:knots:table}
\end{center}
\end{table}

\begin{figure}
\begin{center}
\includegraphics [width = 6 cm] {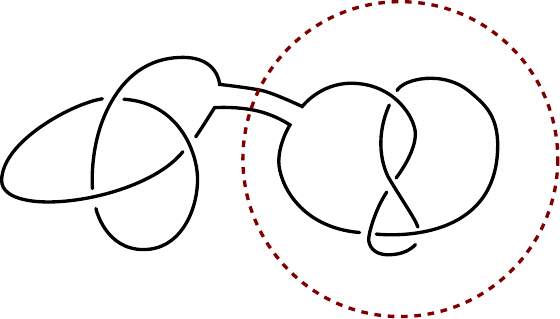}
\nota{A connected sum produces a sphere (dotted in the figure) intersecting the knot in two points.}
\label{sum_of_knots3:fig}
\end{center}
\end{figure}

A connected sum gives rise to a sphere which intersects the new knot in two points, see Figure \ref{sum_of_knots3:fig}. The sphere intersects the knot complement into an essential annulus. 

\subsection{Torus knots}
We now introduce a simple and natural family of prime knots. 
Consider the standardly embedded torus $T\subset S^3$ with meridian $m$ and longitude $l$ oriented as in Figure \ref{torus:ml:fig}. A pair $(p,q)$ of coprime integers determines a simple closed curve in $T$ that is homologically $qm + pl$, and such a curve is called a \emph{$(p,q)$-torus knot}. More generally, a pair of integers $(a,b)$ determines a multicurve in $T$ that is homologically $bm+al$, called an \emph{$(a,b)$-torus link}.\index{knot!torus knot}\index{link!torus link}

\begin{figure}
\begin{center}
\reflectbox{\includegraphics[width = 5 cm] {\iftoggle{BW}{800px-Toroidal_coord_resized-BW}{800px-Toroidal_coord_resized}}}
\nota{The meridian $m$ (\iftoggle{BW}{light grey}{red}) and longitude $l$ (\iftoggle{BW}{dark grey}{blue}) of a standard torus are oriented as shown here, like a right-hand screw.}
\label{torus:ml:fig}
\end{center}
\end{figure}

\begin{figure}
\begin{center}
\includegraphics[width = 5 cm] {\iftoggle{BW}{Blue_7_1_Knot-BW}{480px-Blue_7_1_Knot}}
\includegraphics[width = 5 cm] {\iftoggle{BW}{3D-Link-BW}{3D-Link}}
\nota{The torus knot $(2,-7)$ and the torus link $(2,-8)$.}
\label{torus_knot:fig}
\end{center}
\end{figure}

For instance, the trefoil knot is a $(2,-3)$-torus knot. More examples are in Figure \ref{torus_knot:fig}. The knots $3_1$, $5_1$, and $7_1$ from Figure \ref{Knot_table:fig} are torus knots.

\begin{ex} \label{trivial:torus:knot:ex}
If $|p|\leqslant 1$ or $|q|\leqslant 1$ the torus knot is the unknot. The torus knots with parameters $(p,q)$, $(-p,-q)$, and $(q,p)$ are isotopic, the torus knot with parameters $(p,-q)$ is the mirror image of $(p,q)$.
\end{ex}

The complement of a torus knot is a quite simple kind of Seifert manifold. In the following we suppose that $p,q \geqslant 1$ are coprime. 

\begin{prop} \label{torus:knot:complement:prop}
The complement of a $(p,q)$-torus knot is a Seifert manifold fibering over the orbifold $(D,p,q)$. More precisely, it is
$$\big(D, (p,r), (q,s)\big)$$
where $(r,s)$ is any pair such that $ps+qr = 1$.
\end{prop}
\begin{proof}
Let $K \subset T \subset S^3$ be the $(p,q)$-torus knot. The pair $(r,-s)$ determines another simple closed curve $\alpha \subset T$ that intersects $K$ in one point. The complement of $K$ in a tubular neighbourhood $T\times [-1,1]$ of $T$ is diffeomorphic to $P\times S^1$, where $P$ is a pair of pants. On the tori $T\times \{-1\}$ and $T \times 1$ the curves $\partial P \times \{{\rm pt}\}$ and $\{{\rm pt}\} \times S^1$ are isotopic to $\alpha$ and $K$.

The complement of $T\times [-1,1]$ in $S^3$ consists of two solid tori, with meridians $(1,0)$ and $(0,1)$. Read in the basis $(\alpha, K)$ the meridians are $(q,s)$ and $(p,r)$. The complement of $K$ in $S^3$ is obtained from $P\times S^1$ by filling these curves and hence we get $\big(D, (p,r), (q,s)\big)$.
\end{proof}

In particular the complement of the trefoil knot is $\big(D,(2,1), (3,1)\big)$. If we suppose $p,q\geqslant 2$, two torus knots with distinct (unordered) parameters have non-diffeomorphic complements and hence are not isotopic.
Our intuition says that a torus knot is prime, and we can prove this rigorously.

\begin{prop}
Every torus knot is prime. 
\end{prop}
\begin{proof}
If $K$ is obtained as a non-trivial connected sum, there is an essential annulus in the complement, whose boundary curves are meridians of $K$. There is only one essential annulus in $\big(D, (p,r), (q,s)\big)$ and its boundary curves are not meridians.
\end{proof}

The $(2,2)$-torus link is called the \emph{Hopf link} and is drawn in Figure \ref{Hopf_link:fig}.\index{link!Hopf link}

\begin{ex} The complement of the Hopf link is $T\times [0,1]$. Any two distinct fibres of the Hopf fibration form a Hopf link.
\end{ex}

\begin{figure}
\begin{center}
\includegraphics[width = 4 cm] {\iftoggle{BW}{Hopf_Link-BW}{617px-Hopf_Link}}
\nota{The Hopf link is the $(2,2)$-torus link.}
\label{Hopf_link:fig}
\end{center}
\end{figure}

\subsection{Satellite knots}
We have discovered that che composite knots contain essential annuli in their complement, and we now look for essential tori. It is convenient to introduce the following definition.

\begin{defn}
A knot $K\subset S^3$ is a \emph{satellite} if its complement contains an essential torus.\index{knot!satellite knot}
\end{defn}

We now formulate an equivalent and more inspiring definition of satellite knots which justifies the terminology. A knot in a solid torus $D^2 \times S^1$ is \emph{local} if it is contained in a ball, and a \emph{core} if it is isotopic to $\{x\} \times S^1$. An embedding $\varphi\colon D^2 \times S^1 \hookrightarrow S^3$ is \emph{trivial} (or \emph{unknotted}) if the image of a core is a trivial knot.

\begin{prop}
A knot $K\subset S^3$ is satellite $\Longleftrightarrow$ it is the image of a knot $K' \subset D^2 \times S^1$ which is neither local nor a core, along a non-trivial embedding $\varphi \colon D^2\times S^1 \hookrightarrow S^3$.
\end{prop}
\begin{proof}
A knot $K=\varphi(K')$ constructed in this way is satellite, because the torus $\varphi (S^1\times S^1)$ is essential in the complement of $K$: it is incompressible (on one side because $\varphi$ is non-trivial, and on the other because $K'$ is not local) and not $\partial$-parallel (because $K'$ is not a core).

Conversely, if $K$ is a satellite knot then its complement contains an essential torus $T\subset S^3$. As every torus in $S^3$, the torus $T$ bounds a solid torus. Since $T$ is essential, the knot $K$ is contained in this solid torus in a non-local and non-core way.
Moreover the solid torus is knotted, otherwise $T$ would be compressible on the other side.
\end{proof}

\begin{figure}
\begin{center}
\includegraphics[width = 5.5 cm] {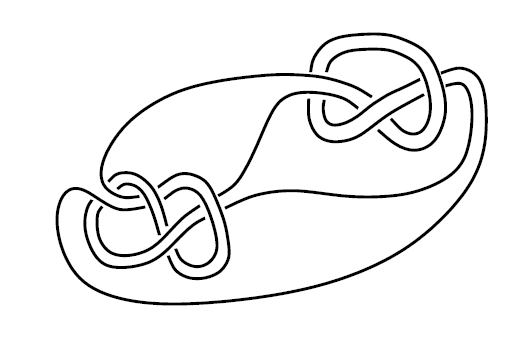}
\includegraphics[width = 5 cm] {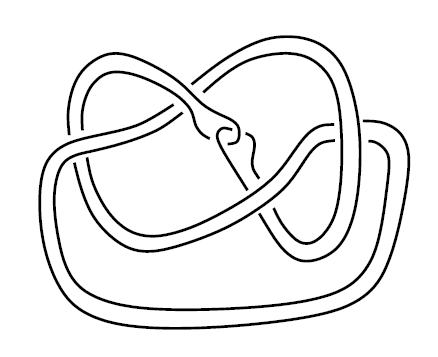}
\nota{Two satellite knots.}
\label{satellites:fig}
\end{center}
\end{figure}

The non-trivial embedding $\varphi$ sends the core curve of $D^2\times S^1$ to some non-trivial knot $H\subset S^3$ called the \emph{companion} of $K$: we should think of $K$ as orbiting as a ``satellite'' around its companion $H$, with orbit path $K'$. When $K'$ is contained in the boundary of the solid torus $D^2\times S^1$ the satellite knot $K$ is called a \emph{cable} knot. Some examples are shown in Figure \ref{satellites:fig}.\index{knot!cable knot}

\subsection{Simple complements}
Recall from Section \ref{simple:manifold:subsection} that a compact three-manifold is simple if it contains no essential sphere, disc, torus, and annulus. Which knots have a simple complement? To answer to this question, we need a general lemma on three-manifolds that identifies the few cases where there are essential annuli without essential tori.

\begin{lemma} \label{only:annuli:lemma}
Let $M$ be irreducible and $\partial$-irreducible, with boundary consisting of tori. The manifold $M$ contains no essential tori but contains some essential annuli $\Longleftrightarrow$ it is diffeomorphic to one of the following:
$$\big(D,(p_1,q_1), (p_2,q_2)\big), \quad \big( A, (p,q) \big), \quad P\times S^1$$
with $p_1,p_2 \geqslant 2$. Here $P$ is a pair of pants.
\end{lemma}
\begin{proof}
Propositions \ref{v:h:prop} and \ref{incompressible:vertical:prop} imply that the Seifert manifolds listed contain vertical essential annuli, but not essential tori. 

Conversely, let $M$ contain an essential annulus $A$. Suppose that $A$ connects two distinct boundary tori $T,T'$ of $M$. A regular neighbourhood of $T\cup T'\cup A$ is diffeomorphic to $P\times S^1$, and its boundary contains a third torus $T''\subset M$. Since $T''$ cannot be essential, it is either boundary parallel or bounds a solid torus in $M$, and $M$ is diffeomorphic respectively  to $P\times S^1$ or a Dehn filling of it. In the latter case, the Dehn filling is not fibre-parallel because $M$ is $\partial$-reducible (use Corollary \ref{fiber:parallel:cor}), hence we get a Seifert manifold of type $\big(A, (p,q) \big)$.

If $A$ connects one boundary component to itself we conclude similarly and may also get $M = \big(D, (p_1,q_1), (p_2,q_2)\big)$. 
\end{proof}

We can now state and prove an important trichotomy on knots. 

\begin{prop} \label{trichotomy:knots:prop}
Every knot $K\subset S^3$ is either a torus knot, a satellite knot, or has a simple complement.
\end{prop}
\begin{proof}
We only need to prove that if the complement $M$ contains an essential annulus and no essential tori then $K$ is a torus knot. The previous lemma gives $M=\big(D,(p_1,q_1),(p_2,q_2)\big)$, and to get $S^3$ back we must have a Dehn filling $S^3 = \big(S^2,(p_1,q_1), (p_2,q_2), (1,n)\big)$. In particular $K$ is isotopic to a fibre and hence contained in a vertical Heegaard torus for $S^3$. Heegaard tori are standard by Corollary \ref{H:tori:cor}, and hence $K$ is a torus knot.
\end{proof}

The relevance of this proposition will be magnified later, after introducing geometrisation. Note that in particular all the composite knots are satellite knots: where is the essential torus?

\section{Dehn surgery}
While a Dehn \emph{filling} consists of attaching a solid torus to a boundary component, a Dehn \emph{surgery} is a two-step operation that consists of drilling a tubular neighbourhood of a knot and then re-gluing it via a different map. 

Using this fundamental operation we can describe \emph{every} closed three-manifold via some reasonable and combinatorial drawing on the plane: a link diagram with some rational numbers attached to its components. Such a combinatorial description has a strong four-dimensional flavour and is called a \emph{Kirby diagram}.

In this section we introduce some basic knot theory concepts: longitudes, Seifert surfaces, Dehn surgery, and the Lickorish-Wallace Theorem which asserts that every closed orientable 3-manifold is obtained by surgerying some link in $S^3$.

\subsection{Canonical longitudes} \label{canonical:longitudes:subsection}
The tubular neighbourhood $N$ of a knot $K\subset S^3$ is a solid torus. As usual with solid tori, a \emph{meridian} is a simple closed curve $m\subset \partial N$ bounding a disc in $N$ and a \emph{longitude} is any other simple closed curve $l$ such that $m$ and $l$ generate $H_1(\partial N,\matZ)$. The meridian $m$ is unique up to sign, but the longitude $l$ is not: if $l$ is a longitude, then $l+km$ also is for any $k\in \matZ$. 

The purpose of this section is to define a canonical longitude.

\begin{prop} Let $L\subset S^3$ be a link with $k$ components and $M$ its complement. We have $H_1(M,\matZ) = \matZ^k$, generated by the $k$ meridians.
\end{prop}
\begin{proof}
Let $N = N_1\sqcup \ldots \sqcup N_k$ be the solid tori neighbourhoods of $L$ and $T_i= \partial N_i$. The Mayer--Vietoris sequence on $S^3 = M\cup N$ gives
$$0 \longrightarrow H_1(T_1\sqcup \ldots \sqcup T_k) \longrightarrow H_1(M) \oplus H_1(N_1\sqcup \ldots \sqcup N_k) \longrightarrow 0$$
since $H_2(S^3) = H_1(S^3) = 0$. The equalities $H_1(T_i) = \matZ\times \matZ$ and $H_1(N_i) = \matZ$ imply that $H_1(M) = \matZ^k$.
The group $H_1(T_i)$ is generated by $(m_i, l_i)$ and $m_i$ goes to zero in $H_1(N_i)$. Therefore the meridians $m_1,\ldots, m_k$ go to generators of $H_1(M)$.
\end{proof}

\begin{cor} \label{longitude:cor}
Let $K\subset S^3$ be a knot and $M$ be its complement. A unique (up to sign) longitude $l \subset \partial M$ vanishes in $H_1(M, \matZ)$. 
\end{cor}
\begin{proof}
In the map $\matZ\times \matZ = H_1(\partial M) \to H_1(M)  =\matZ $ the meridian $m$ goes to a generator, hence the kernel is generated by a longitude $l$.
\end{proof}

We call $l$ the \emph{canonical longitude} of $K$. The torus $T=\partial M$ is hence equipped with a canonical basis $(m,l)$ for $H_1(T,\matZ)$; we orient $m$ and $l$ as shown in Figure \ref{torus:ml:fig}, like a right-hand screw. The pair $(m,l)$ is well-defined up to reversing both their signs. 

\subsection{Seifert surfaces}
We now show that the canonical longitude has a concrete geometric interpretation.

\begin{figure}
\begin{center}
\includegraphics [width = 3 cm] {\iftoggle{BW}{Noeud_de_trefle_et_surface_de_seifert-BW}{500px-Noeud_de_trefle_et_surface_de_seifert}}
\nota{A Seifert surface for the trefoil knot: it is a punctured torus.}
\label{Seifert_surface:fig}
\end{center}
\end{figure}

\begin{defn} A \emph{Seifert surface} for a knot $K\subset S^3$ is any orientable connected compact surface $S\subset S^3$ with $\partial S = K$.\index{Seifert surface}
\end{defn}

See an example in Figure \ref{Seifert_surface:fig}. Every Seifert surface $S$ determines a longitude $l$ for $K$: pick a small tubular neighbourhood $N$ of $K$ and set $l=S\cap \partial N$. The same knot $K$ has plenty of non-isotopic Seifert surfaces, but luckily these all induce the same longitude:

\begin{prop}
Every knot $K$ has a Seifert surface $S$. Every Seifert surface for $K$ induces the canonical longitude $l$.
\end{prop}
\begin{proof}
Let $M$ be the complement of $K$. Let $S$ be a surface representing a generator of $H_2(M,\partial M) = H^1(M) = \matZ$. The long exact sequence
$$\ldots \longrightarrow H_2(M,\partial M) \longrightarrow H_1(\partial M) \longrightarrow H_1(M) \longrightarrow \ldots$$ 
implies that $[S]$ is mapped to a non-trivial primitive element $\alpha \in H_1(\partial M)=\matZ\times \matZ$ that is trivial in $H_1(M)$. Therefore $[\partial S]=\alpha = [l]$ and $\partial S$ consists of an odd number of parallel copies of $l$ and some homotopically trivial simple closed curves in $T = \partial M$. The homotopically trivial components may be eliminated by isotoping them inside $M$ and capping them with discs; if the parallel copies are more than 1, since their signed sum is 1, there must be two of them that are close and with opposite signs, that can be canceled by isotoping them inside $M$ and capping them with an annulus. At the end we get $\partial S = K$.
\end{proof}

The \emph{Seifert genus} of $K$ is the minimum genus of a Seifert surface for $K$.\index{knot!Seifert genus}

\begin{prop} The unknot is the only knot with Seifert genus zero.
\end{prop}
\begin{proof}
A knot has genus zero $\Longleftrightarrow$ it bounds a disc.
\end{proof}

Figure \ref{Seifert_surface:fig} shows that the trefoil knot has genus one.

\subsection{Dehn surgery}
Let $L\subset M$ be a link with some $k$ components in an orientable 3-manifold $M$. A \emph{Dehn surgery} on $L$ is a Dehn filling of the complement of $L$. In other words, it is a two-step operation that consists of:\index{Dehn surgery}
\begin{enumerate}
\item (\emph{drilling}) the removal of small open tubular neighbourhoods of $L$, that creates new boundary tori $T_1,\ldots, T_k$; 
\item (\emph{filling}) a Dehn filling of the new boundary tori $T_1,\ldots, T_k$.
\end{enumerate}
We remove the tubular neighbourhoods of $L$ and glue them back differently: this explains the use of the word \emph{surgery}.
The outcome is a new manifold $N$ that has the same boundary of $M$, but $N$ may be not diffeomorphic to $M$.

The surgered manifold $N$ is determined by the slopes in $T_1, \ldots, T_k$ that are killed by the Dehn filling, see Section \ref{Dehn:filling:subsection}. When $M=S^3$, every torus $T_i $ is equipped with a canonical basis $m_i, l_i$, the slope is of the form $\pm(p_im_i + q_il_i)$ and is hence determined by the rational number $\frac {p_i}{q_i}\in\matQ \cup \{\infty\}$. The full surgery on the link is comfortably encoded by assigning the number $\frac {p_i}{q_i}$ to the $i$-th component of $L$, for every $i=1,\ldots,k$.

The result of a Dehn surgery along $L\subset S^3$ is a closed orientable 3-manifold. The slope $\infty = \frac 10$ indicates the meridian $m_i$.

\begin{prop} An $\infty$-surgery on a knot $K$ has no effect.
\end{prop}
\begin{proof}
It consists of removing a solid torus neighbourhood of $K$ and regluing it back with the same map.
\end{proof}

\begin{figure}
\begin{center}
\includegraphics[width = 5 cm] {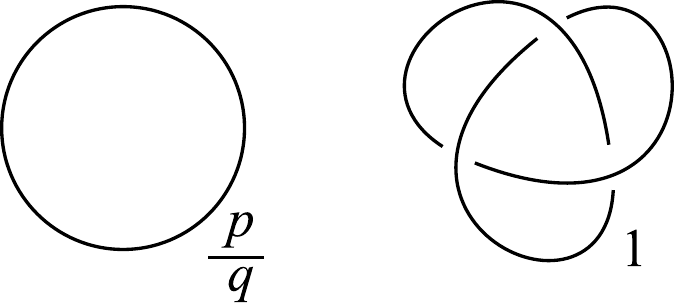}
\nota{Two Kirby diagrams describing the lens space $L(p,q)$ and the Poincar\'e homology sphere $\Sigma(2,3,5)$.}
\label{Kirby:fig}
\end{center}
\end{figure}

A \emph{Kirby diagram} is a link diagram on the plane with a rational number $\frac{p_i}{q_i}$ assigned to each component. Such a diagram defines a Dehn surgery and hence a closed orientable three-manifold.\index{Kirby diagram} Some examples are shown in Figure \ref{Kirby:fig}.

\begin{ex} 
The $\frac pq$-surgery on the unknot yields the lens space $L(p,q)$.
\end{ex}

The slope $0$ indicates the canonical longitude $l_i$. More generally, the integer $n_i \in \matZ$ encodes the longitude $n_im_i + l_i$. The surgery is \emph{integral} if the coefficients $\frac {p_i}{q_i}$ are all integers.

\begin{rem} The notion of \emph{integral} Dehn surgery exists for any link $L\subset M$ in any 3-manifold $M$: a Dehn surgery is integral if the killed slopes are longitudes of the previously removed solid tori. However on a generic $M$ there is no rule for choosing a canonical longitude and hence to transform slopes into numbers.
\end{rem}

Let $N$ be obtained by Dehn surgery on a knot $K\subset S^3$ with coefficient $\frac{p}{q}$.

\begin{prop}
The surgered manifold $N$ has $H_1(N, \matZ) = \matZ/_{p\matZ}$.
\end{prop}
\begin{proof}
Let $M$ be the complement of $K$. We know that the meridian $m$ generates $H_1(M) = \matZ$ while the longitude $l$ is zero there. The Dehn filling kills the element $pm + ql = pm$.
\end{proof}

Recall that a \emph{homology sphere} is a closed 3-manifold $M$ having the same integral homology as $S^3$, that is with trivial $H_1(M, \matZ)$.

\begin{cor}
If the coefficient is $\frac 1{q}$ the surgered manifold $N$ is a homology sphere.
\end{cor}

We have just found a simple method to construct plenty of homology spheres.

\subsection{Torus knots}
As a first example, we study the effect of a Dehn surgery on torus knots. We suppose $p,q \geqslant 1$.

\begin{prop}
The $\frac tu$-Dehn surgery on a $(p,q)$-torus knot gives
\begin{align} 
\big(S^2, (p,r), (q,s), (t-pqu,u)\big) & \quad {\rm if}\ \frac {t}u \neq pq, \label{Stu:eqn} \\
L(p,r) \# L(q,s) & \quad {\rm if}\ \frac tu = pq
\end{align}
where $(r,s)$ is any pair with $ps+qr = 1$.
\end{prop}
\begin{proof}
Recall from the proof of Proposition \ref{torus:knot:complement:prop} that the complement of $K$ in $T\times [-1,1]$ is diffeomorphic to $P\times S^1$. The meridian $m$ of $K$ is isotopic to $\partial P \times \{{\rm pt}\}$, so the canonical longitude $l$ of $K$ is isotopic to $l'+km$ for some $k\in \matZ$, where $l'$ is isotopic to $\{{\rm pt}\}\times S^1$, \emph{i.e.}~it is the framing induced by $T$.

By Proposition \ref{torus:knot:complement:prop} the Dehn filling is $\big(S^2, (p,r), (q,s), (t-ku,u)\big)$ and we need to determine $k$. Recall that if $\frac tu = \frac 01$ the manifold has infinite cyclic homology, hence it is a Seifert manifold with $e=\frac 1{pq}-\frac 1k = 0$ by Proposition \ref{Seifert:homology:prop}: we deduce that $k=pq$.

When $t-pqu = 0$ the filling is fiber-parallel and Corollary \ref{fiber:parallel:cor} shows that we actually get $L(p,r) \# L(q,s)$.
\end{proof}
Recall from Corollary \ref{sigma:cor} that for every triple $(p,q,n)$ of pairwise coprime natural numbers $p,q,n\geqslant 2$ there is a unique Seifert homology sphere $\Sigma(p,q,n)$ fibering over the orbifold $(S^2,p,q,n)$.

\begin{cor}
Fix two coprime $p,q \geqslant 2$. The Seifert homology sphere obtained by $\frac 1u$-surgery on the $(p,q)$-torus knot is $\Sigma(p,q,|pqu-1|)$.
\end{cor}
In particular the $1$-surgery on the trefoil knot (seen as a $(2,3)$-torus knot) depicted in Figure \ref{Kirby:fig} produces the Poincar\'e homology sphere $\Sigma(2,3,5)$.
 
\begin{ex} Every homology sphere of type $\Sigma(2,3,k)$ may be obtained via a Dehn surgery of the trefoil knot.
\end{ex}

Note that when $t-pqu = \pm 1$ the Seifert manifold (\ref{Stu:eqn}) has only two singular fibres and is hence a lens space.

\begin{cor}
Fix two coprime $p,q \geqslant 2$ and $n\geqslant 1$. When $\frac tu = pq \pm \frac 1n$ the $\frac tu$-Dehn surgery on a $(p,q)$-torus knot produces the lens space 
$$L(t,uq^2) = L(npq \pm 1, nq^2).$$
\end{cor}
\begin{proof}
The Seifert manifold (\ref{Stu:eqn}) is
\begin{align*}
\big(S^2, (p,r), (q,s), (1, \pm n) \big) & = \big(S^2, (p,r), (q,s \pm nq) \big) \\
& = L(p(s\pm nq) +rq, (-q)(s\pm nq) + sq) \\
& = L(1 \pm npq, \mp nq^2)
\end{align*}
using Exercise \ref{two:lens:ex}.
\end{proof}
We note that we never get $S^2\times S^1$ or $S^3$ via non-trivial surgeries on non-trivial torus knots. The lens space with smallest fundamental group that we may get is $L(5,1)$, which arises from a $5$-Dehn surgery on the trefoil knot.

\subsection{The Lickorish-Wallace theorem}
We now prove that the Dehn surgery construction is as general as possible.

\begin{teo}[Lickorish-Wallace theorem] \label{LW:teo}
Every closed orientable 3-manifold can be described via an integral Dehn surgery along a link $L\subset S^3$.
\end{teo}
\begin{proof}
Let $M$ be a closed orientable 3-manifold. Pick a Heegaard splitting $M=H_1 \cup_\psi H_2$ where $H_1$ and $H_2$ are genus-$g$ handlebodies and $\psi\colon \partial H_1\to \partial H_2$ is a diffeomorphism. We fix an identification of both $H_1$ and $H_2$ with a model handlebody $H$, so that $\psi$ can be interpreted as an element of the mapping class group $\MCG(S)$ of the genus-$g$ surface $S=\partial H$. 

\begin{figure}
\begin{center}
\includegraphics[width = 6 cm] {\iftoggle{BW}{LW-BW}{LW}}
\nota{We push $\gamma$ inside a collar for $\partial H_2$ and drill a solid torus neighbourhood, represented here as a \iftoggle{BW}{light grey}{yellow} parallelepiped.}
\label{LW:fig}
\end{center}
\end{figure}

Example \ref{S3:example} shows that the three-sphere also decomposes as $S^3 = H_1\cup_\varphi H_2$ for some $\varphi\in \MCG(S)$. Theorem \ref{twists:generate:teo} says that $\psi\circ \varphi^{-1}$ is a composition of Dehn twists
$$\psi\circ \varphi ^{-1} = T_{\gamma_k}^{\pm 1} \circ \ldots \circ T_{\gamma_1}^{\pm 1}$$
along some curves $\gamma_i \subset S$. Set $M_i = H_1\cup_{\psi_{i}} H_2$ with 
$$\psi_i = T_{\gamma_i}^{\pm 1} \circ \ldots \circ T_{\gamma_1}^{\pm 1}\circ \varphi.$$
We have $M_0=S^3$ and $M_k=M$. We prove that $M_i$ can be described via an integral Dehn surgery along a $i$-components link in $S^3$ by induction on $i$. To obtain that it suffices to check that $M_{i+1}$ can be obtained from $M_i$ via integral Dehn surgery along a knot. We have 
$$M_i = H_1\cup_{\psi_i} H_2, \qquad M_{i+1} = H_1\cup_{T_{\gamma_{i+1}}^{\pm 1}\circ \psi_i} H_2.$$
Now fix a collar of $S=\partial H_2$ in $H_2$ and push $\gamma_{i+1}$ inside the collar as in Figure \ref{LW:fig}. Drill from $H_2$ a solid torus around $\gamma_{i+1}$ (a \iftoggle{BW}{light grey}{yellow} parallelepiped in the figure) to get a submanifold $H_2^{\rm drill} \subset H_2$. We see $T_{\gamma_{i+1}^{\pm 1}}$ as a Dehn twist supported in the annulus $A\subset S$ drawn in the figure. 

The Dehn twist supported on $A$ extends product-wise to the solid torus $A\times [0,1]$ lying between $A$ and the drilled \iftoggle{BW}{light grey}{yellow} parallelepiped, and extends trivially to a self-diffeomorphism $T\colon H_2^{\rm drill} \to H_2^{\rm drill}$ such that $T|_S = T_{\gamma_{i+1}}^{\pm 1}$. We define
$$M_j^{\rm drill} = H_1 \cup_{\psi_j} H_2^{\rm drill}$$
for $j=i, i+1$. Since $T|_S = T_{\gamma_{i+1}}^{\pm 1}$ the map $T$ extends to a diffeomorphism
$$T\colon M_i^{\rm drill} \longrightarrow M_{i+1}^{\rm drill}.$$
Therefore $M_{i+1}$ is obtained from $M_i$ by Dehn surgery along $\gamma_{i+1}$. The surgery is integral since $T$ sends a meridian of the drilled \iftoggle{BW}{light grey}{yellow} solid torus to a longitude. The proof is complete.
\end{proof}

The theorem shows that every closed orientable three-manifold can be constructed using Kirby diagrams. As for knots and many other combinatorial descriptions of topological objects, it is important to note that many different Kirby diagrams may describe the same three-manifold, and in general it is hard to tell, given two diagrams, whether they define the same manifold or not.

\subsection{Four-manifolds}
The theorem of Lickorish and Wallace has a nice four-dimensional interpretation, which shows in particular that every closed orientable three-manifold is the boundary of a simply-connected four-manifold.

Recall from Section \ref{handles:subsection} that in dimension four a two-handle $H=D^2 \times D^2$ is attached to an orientable four-manifold $W$ via an embedding $\psi\colon S^1 \times D^2 \to \partial W$. After attaching the handle $H$ to $W$, the boundary of $W$ changes by substituting the solid torus $\psi(S^1 \times D^2)$ with another solid torus $D^2 \times S^1$: in other words, it changes by Dehn surgery. The Dehn surgery is integral because the meridians of the two solid tori intersect in one points. 

Conversely, every integral Dehn surgery on a knot $K\subset \partial W$ may be interpreted as the effect of attaching some 2-handle: it suffices to interpret the drilled and filled solid tori as horizontal and vertical boundaries of a four-dimensional two-handle.

With such an interpretation, we may consider $S^3$ as the boundary of $D^4$ and see the integral Dehn surgery on a link $L\subset S^3$ as the result of attaching some two-handles to $D^4$. The Lickorish-Wallace theorem implies the following.

\begin{cor}
Every closed orientable three-manifold is the boundary of a simply-connected compact four-manifold.
\end{cor}
\begin{proof}
Every closed orientable three-manifold is the result of an integral surgery on a link $L\subset S^3$ and is hence the boundary of a four-manifold obtained by attaching some two-handles to $D^4$. Every such manifold is simply-connected by Van Kampen's theorem.
\end{proof}

This theorem reveals in particular that it is fairly easy to construct plenty of simply connected four-manifolds with boundary. Note that the simply connected four-manifolds constructed in this way are not contractible: indeed the following proposition shows that they have a non-trivial second homology group.

\begin{prop}
A four-manifold with one 0-handle and k 2-handles is homotopy equivalent to a bouquet of $k$ two-spheres.
\end{prop}
\begin{proof}
We first deformation retract the two-handles over their core discs, and then we shrink the whole 0-handle to a point.
\end{proof}

To construct contractible four-manifolds, we need to employ 1-handles.

\subsection{Three-manifolds with boundary}
The Lickorish-Wallace Theorem extends to compact manifolds with boundary, in the appropriate way. 

\begin{prop} \label{LW:bordo:prop}
Every compact orientable 3-manifold $M$ with boundary is obtained from $S^3$ as follows: pick $L\sqcup C\subset S^3$ where $L$ is a link and $C$ is a 1-complex; remove an open regular neighbourhood of $C$ and perform an integral surgery on $L$.
\end{prop}
\begin{proof}
Cap off the boundary of $M$ by adding handlebodies $H_g$, to get a closed manifold $M'$. Choose a 1-complex $C\subset M'$ whose regular neighbourhood consists of these handlebodies, so that by removing $C$ we get $M$ back. 

Now $M'$ is obtained from $S^3$ via integral surgery along some link $L$, and by general position we may suppose that $C$ is disjoint from the cores of the surgered solid tori, so that we can see both $C$ and $L$ disjointly in $S^3$.
\end{proof}

If $\partial M$ consists of tori we may suppose that $C$ is also a link. 

\section{Surface bundles}
Seifert manifolds have base orbifolds of dimension 2 and smooth fibres of dimension 1. We now introduce a complementary construction where the base orbifolds have dimension 1 and the smooth fibres have dimension 2. These are called \emph{surface bundles}.

\subsection{Surface bundles}
A \emph{surface bundle over $S^1$} is a fibre bundle $M\to S^1$ of a compact orientable 3-manifold $M$ (possibly with boundary) over the circle, whose fibre $\Sigma$ is a connected compact orientable surface. If $M$ has boundary, then $\Sigma$ also has, and $\partial M$ consists of tori fibering over $S^1$.\index{surface bundle over $S^1$}

\begin{prop}
Every surface bundle over $S^1$ is constructed by taking $\Sigma\times [0,1]$ and glueing $\Sigma\times 0$ to $\Sigma\times 1$ via an orientation-preserving diffeomorphism $\psi$.
\end{prop}
\begin{proof}
One such glueing clearly gives rise to a surface bundle over $S^1$. Conversely, by cutting a surface bundle over $S^1$ along a fibre we get a surface bundle over the interval, which is a product $\Sigma\times [0,1]$.
\end{proof}

The diffeomorphism $\psi$ is the \emph{monodromy} of the surface bundle $M_\psi$. Since isotopic glueings produce diffeomorphic manifolds, the three-manifold $M_\psi$ depends only of the class of $\psi$ in the mapping class group $\MCG(\Sigma)$ on $\Sigma$. More than that, it actually depends only on its conjugacy class:

\begin{prop}
If $\psi$ and $\phi$ are conjugate in $\MCG(\Sigma)$, then $M_\psi$ and $M_\phi$ are diffeomorphic.
\end{prop}
\begin{proof}
The diffeomorphism $g\colon \Sigma\to \Sigma$ that conjugates them extends to $\Sigma\times [0,1]$ and gives a diffeomorphism $M_\psi \to M_\phi$.
\end{proof}

The manifolds $M_\psi$ and $M_{\psi^{-1}}$ are (orientation-reversingly) diffeomorphic.

\subsection{Properties}
We now start to investigate the topological properties of surface bundles.
Let $M\to S^1$ be a surface bundle with fibre $\Sigma$.

\begin{ex}
The maps $\Sigma\to M \to S^1$ induce an exact sequence
$$0 \longrightarrow \pi_1(\Sigma) \longrightarrow \pi_1(M) \longrightarrow \pi_1(S^1) \longrightarrow 0.$$
\end{ex}

This implies in particular that $\pi_1(M)$ surjects onto $\matZ$ and therefore:

\begin{cor}
We have $b_1(M)\geqslant 1$.
\end{cor}

In other terms, the surface fibre $\Sigma$ is non-separating and hence $[\Sigma]\in H_2(M,\partial M)$ is non-trivial (and has infinite order, since there is no torsion there).

Note that there is an obvious degree-$n$ regular covering $M_{\psi^n} \to M_{\psi}$ for every $n$ and an infinite regular covering $\Sigma\times \matR \to M_{\psi}$ induced by the normal subgroup $\pi_1(\Sigma) \triangleleft \pi_1(M)$. 

\begin{prop}
The fibre $\Sigma$ is an essential surface.
If $\chi(\Sigma)>0$ then $M$ is diffeomorphic to $D\times S^1$ or $S^2\times S^1$.
If $\chi(\Sigma)\leqslant 0$ the universal cover of $\interior M$ is $\matR^3$ and $M$ is Haken.
\end{prop}
\begin{proof}
If $\chi(\Sigma)>0$ then $\MCG(\Sigma)$ is trivial and we are done, so we suppose $\chi(\Sigma)\leqslant 0$. The fibre $\Sigma$ is incompressible because $\pi_1(\Sigma)$ injects, and is also $\partial$-incompressible by a doubling argument (the double $DM$ fibres to $S^1$ with incompressible fibre $D\Sigma$, hence $\Sigma$ is $\partial$-incompressible). The fibre is clearly not $\partial$-parallel, hence it is essential.

The manifold $M$ is covered by $\Sigma \times \matR$, whose interior is covered by $\matR^2 \times \matR = \matR^3$: hence $M$ is irreducible. It is also $\partial$-irreducible because its double also fibres and hence is irreducible. Therefore $M$ is Haken. 
\end{proof}

\subsection{Semi-bundles}
The Seifert fibrations are circle bundles over two-dimensional orbifolds, and likewise it is natural to consider surface bundles over 1-orbifolds. The compact 1-orbifolds are $S^1$ and the closed segment with mirrored endpoints. We have already considered the $S^1$ case and we now define some surface bundles over the closed segment, called \emph{semi-bundles}.

Let $\Sigma$ be a non-orientable surface, with orientable double cover $p\colon \tilde \Sigma \to \Sigma$ and deck transformation $\tau$ that gives $\Sigma=\tilde \Sigma/_\tau$. We have 
$$\Sigma\timtil (-1,1) = \big(\tilde \Sigma \times (-1,1)\big)/_{(\tau, \iota)}$$
with $\iota(x) = -x$. A \emph{local semi-bundle} is the map $\Sigma\timtil (-1,1) \to [0,1)$ that sends $(p,x)$ to $|x|$. The fibre over $0$ is $\Sigma$, that over $x\in (0,1)$ is $\tilde \Sigma$.

A \emph{semi-bundle} $M \to [-1,1]$ is a map which is a local semi-bundle when restricted to $[-1, 1-\varepsilon)$ and $(\varepsilon, 1]$. The fibre over $\pm 1$ is $\Sigma$ and the fibre of $x\in (-1,1)$ is $\tilde \Sigma$. We should think about this object as a surface bundle over the segment orbifold $[-1,1]$ with mirror points $\pm 1$. 

\begin{ex}
Every semi-bundle is constructed by gluing two copies of $\Sigma\timtil [-1,1]$ along their boundaries via some diffeomorphism.
\end{ex}

Set $I = [-1,1]$. Let $M\to I$ be a semi-bundle with fibres $\Sigma$ and $\tilde \Sigma$.
Note that $I=\matR/_\Gamma$ where $\Gamma<\Iso(\matR)$ is generated by the reflections at the points $\pm 1$. The orbifold fundamental group of $I$ is 
$$\pi_1(I) = \Gamma = \matZ/_{2\matZ} * \matZ/_{2\matZ}.$$

\begin{ex}
The maps $\tilde\Sigma\to M \to I$ induce an exact sequence
$$0 \longrightarrow \pi_1(\tilde \Sigma) \longrightarrow \pi_1(M) \longrightarrow \pi_1(I) \longrightarrow 0.$$
\end{ex}

\begin{prop}
If $\chi(\Sigma)>0$ then $M$ is diffeomorphic to $\matRP^2 \timtil S^1$. If $\chi(\Sigma)\leqslant 0$ the orientable fibre $\tilde \Sigma$ is essential, the universal cover of $\interior M$ is $\matR^3$, and $M$ is Haken.
\end{prop}
\begin{proof}
Same proof as in the standard bundle case.
\end{proof}

\begin{ex} \label{semi:iff:ex}
An orientable 3-manifold $M$ has a (semi-)bundle structure if and only if there is an orientable $\Sigma \subset M$ that cuts $M$ into interval bundles.
\end{ex}

Let $M\to I$ be a semi-bundle. We may pull-back the semi-bundle along the orbifold double cover $S^1 \to I$ and get an ordinary surface bundle $\tilde M \to S^1$, so that the following diagram commutes:
$$
\xymatrix{ 
\tilde M\ar@{.>}[r] \ar@{.>}[d] & M \ar[d] \\
S^1\ar[r] & I
}
$$
Every semi-bundle is thus covered by an ordinary bundle.

\subsection{Seifert manifolds} \label{Seifert:fibers:subsection}
We now classify the Seifert manifolds that have a surface (semi-)bundle structure.

\begin{prop} \label{Seifert:bundle:prop}
A Seifert manifold $M$ has a (semi-)bundle structure if and only if one of the following holds: 
\begin{itemize}
\item $\partial M \neq 0$, 
\item $e(M)=0$, 
\item $M = \big(T,(1,e)\big)$ or $\big(K,(1,e)\big)$,
\item $M = \big(S^2, (2,1),\! (2,1),\! (2,1),\! (2,2q+1)\big)$ or $\big(\matRP^2, (2,1),\! (2,2q+1)\big)$.
\end{itemize}
The manifolds in the last line occur only as semi-bundles.
\end{prop}
\begin{proof}
If $\partial M \neq \emptyset$ or $e(M)=0$ then $M\to S$ has an orientable section $\Sigma$, which cuts $M$ into interval bundles: hence $\Sigma$ is a fibre of a (semi-)bundle (see Exercise \ref{semi:iff:ex}).

Let $S$ be one of the orbifolds
$$T, \quad K, \quad (S^2,2,2,2,2), \quad (\matRP^2, 2,2).$$
It contains a circle that splits $S$ into one or two annuli, M\"obius strips, or $(D,2,2)$. 
Every circle fibering over these pieces is an interval bundle by Proposition \ref{diffeo:fiberings:prop}. Therefore every fibering $M \to S$ contains a vertical torus that splits $M$ into interval bundles and is hence a (semi-)bundle.

Conversely, if $M$ has a (semi-)bundle structure the orientable fibre is essential and is hence either horizontal or vertical: we get one of the types listed.
\end{proof}

\subsection{Torus bundles}
A \emph{torus bundle} is of course a surface bundle $M\to S^1$ with fibre a torus $T$. We fix a basis for $\pi_1(T)$, so that $\MCG(T) = \SLZ$. By what said above, every matrix $A\in \SLZ$ defines a torus bundle $M_A$\index{torus bundle} with monodromy $A$. We want to understand when $M_A$ is a Seifert manifold.

\begin{ex} \label{Me:ex}
For every $e\in\matZ$ there are diffeomorphisms 
$$M_{\matr 1e01} \isom \big(T,(1,e)\big), \quad M_{\matr {-1}e0{-1}} \isom \big(K,(1,e)\big). $$
\end{ex}

The classification of torus bundles reduces to linear algebra.

\begin{prop} \label{torus:bundle:classification:prop}
Two torus bundles $M_A$ and $M_{A'}$ are diffeomorphic if and only if $A'$ is conjugate to $A^{\pm 1}$ in $\SLZ$.
\end{prop}
\begin{proof}
Let $M_A \isom M_{A'}$ be a diffeomorphism. Consider the tori $T$ and $T'$ of the two fibrations both inside $M_A$ and minimise their transverse intersection. 

If $T\cap T' = \emptyset$, by Proposition \ref{line:bundles:prop} the tori are parallel and the two fibrations are isotopic: hence $A'$ is conjugate to $A^{\pm 1}$. Otherwise by the same proposition $T'$ decomposes into fibered annuli in the product line bundle $M_A\setminus T$. One such fibered annulus identifies a circle $\gamma\subset T$ preserved by $A$. Therefore $A$ is conjugate to $\matr 1 e 0 1 $ or $\matr {-1} e 0 {-1}$, and $A'$ is conjugate to $\matr 1 {e'} 0 1 $ or $\matr {-1} {e'} 0 {-1}$ for the same reason. Exercise \ref{Me:ex} implies that $A' = A^{\pm 1}$.
\end{proof}

\begin{ex} A matrix $A\in \SLZ$ has finite order if and only if $A=\pm I$ or $|\tr A| < 2$. Every finite-order $A$ is conjugate to one of
$$ 
\begin{pmatrix} 1 & 0 \\ 0 & 1 \end{pmatrix}, \quad
\begin{pmatrix} -1 & 0 \\ 0 & -1 \end{pmatrix}, \quad
\begin{pmatrix} -1 & 1 \\ -1 & 0 \end{pmatrix}, \quad
 \begin{pmatrix} 0 & -1 \\ 1 & 0 \end{pmatrix}, \quad
\begin{pmatrix} 1 & -1 \\ 1 & 0 \end{pmatrix}.
$$
These matrices have order 1, 2, 3, 4, 6.
\end{ex}
\begin{proof}[Hint]
Use Proposition \ref{PSLR:prop}.
\end{proof}

We can easily determine whether a torus bundle is Seifert by looking at its monodromy $A$.

\begin{prop} \label{torus:bundle:prop}
Let $M = M_A$ be a torus bundle with monodromy $A\neq \pm I$. The following holds:
\begin{itemize}
\item if $|\tr A| < 2$ then $M$ is a Seifert manifold with $e=0$ and $\chi = 0$,
\item if $|\tr A|=2$ then $M$ is a Seifert manifold with $e\neq 0$ and $\chi = 0$,
\item if $|\tr A|>2$ then $M$ is not a Seifert manifold.
\end{itemize}
\end{prop}
\begin{proof}
Consider $T \times [-1,1]$ foliated by lines $\{x\}\times [-1,1]$. The foliation extends to $M_A$. If $A$ has finite order, then $M_A$ is finitely covered by $M_I = T\times S^1$ and hence all fibers are compact. Therefore $M_A$ is Seifert fibered and covered by $T\times S^1$, and we get $e=\chi=0$. 

If $|\tr A|=2$ then $A$ is conjugate to $\matr {\pm 1}e0{\pm 1}$ and we use Exercise \ref{Me:ex}. Proposition \ref{Seifert:bundle:prop} easily shows that all the Seifert manifolds that are torus bundles are realised with $|\tr A| \leqslant 2$, hence if $|\tr A|>2$ the manifold $M$ is not Seifert by Proposition \ref{torus:bundle:classification:prop}.
\end{proof}

When $|\tr A| > 2$ we say that the monodromy $A$ is \emph{Anosov}.

\subsection{Bundles with $\chi(\Sigma) <0$}
Proposition \ref{torus:bundle:classification:prop} does not extend to surface bundles with $\chi(\Sigma)<0$; indeed it may happen that non-conjugate monodromies in $\MCG(\Sigma)$ give rise to diffeomorphic manifolds and understanding when this happens is a hard problem.

Proposition \ref{torus:bundle:prop} extends nevertheless and reflects the trichotomy of mapping classes. Let $\Sigma$ be a closed orientable surface with $\chi(\Sigma)<0$.
Recall from Section \ref{surface:diffeomorphisms:section} that every element $\psi\in\MCG(\Sigma)$ is either finite order, reducible, or pseudo-Anosov.

\begin{prop} \label{trichotomy:prop}
Let $M_\psi$ be a surface bundle with fibre $\Sigma$ and monodromy $\psi\in\MCG(\Sigma)$. The following holds:
\begin{itemize}
\item if $\psi$ has finite order, then $M_\psi$ is Seifert with $\chi <0$ and $e=0$,
\item if $\psi$ is reducible, then $M_\psi$ contains an essential torus,
\item if $\psi$ is pseudo-Anosov, then $M_\psi$ is simple and not Seifert.
\end{itemize}
\end{prop}
\begin{proof}
Same proof as Proposition \ref{torus:bundle:prop}. If $\psi$ has finite order, it is an isometry for some hyperbolic metric on $\Sigma$, and the line fibration of $\Sigma \times [-1,1]$ glues to a Seifert fibration for $M_\psi$. 

If $\psi$ is reducible there are disjoint essential simple closed curves $\gamma_1,\ldots, \gamma_k$ with $\psi(\gamma_i) = \gamma_{i+1}$ cyclically; by gluing the annuli $\gamma_i \times [-1,1]$ we get a torus $T\subset M_\psi$. It is essential because by cutting along it we still get a fibration over $S^1$ with fibers having $\chi \leqslant 0$.

If $\psi$ is pseudo-Anosov there are no essential tori $T\subset M_\psi$, for by minimising $T\cap \Sigma$ then $M_\psi \setminus T$ would consist of essential annuli of type $\gamma \times [-1,1]$ and hence $\psi$ would be reducible.
The manifold $M_\psi$ is not Seifert because the fibre $\Sigma$ would become a horizontal surface: then $M$ would be covered by $\Sigma \times S^1$ and hence $\psi$ would be of finite order.
\end{proof}

\section{JSJ decomposition}
In Chapter \ref{Three:topology:chapter} we cut every closed three-manifold along spheres, and it is now time to decompose it further along tori. 
This two-steps cutting operation is called the \emph{JSJ decomposition} of the manifold, after the names of Jaco, Shalen, and Johansson who discovered it in the mid 1970s.\index{JSJ decomposition}

The core of this decomposition is the existence of a canonical set of disjoint essential tori, unique up to isotopy. 

\subsection{Canonical torus decomposition}
Let $M$ be an orientable irreducible and $\partial$-irreducible compact 3-manifold with (possibly empty) boundary consisting of tori. How can we define a canonical set of disjoint essential tori in $M$? The answer is not obvious: for instance, if $M$ is a Seifert manifold, it may contain many vertical incompressible tori and there is no canonical way to choose among them.

We will soon see that the Seifert manifolds are in fact the only possible source of ambiguity. Let
$$S=T_1\sqcup \cdots \sqcup T_k$$ 
be a set of disjoint essential tori $T_i\subset \interior M$. 
We say that $S$ is a \emph{torus decomposition} of $M$ if it decomposes $M$ into blocks that are either:\index{torus decomposition} 
\begin{itemize}
\item torus (semi-)bundles,
\item Seifert manifolds, or
\item simple manifolds.
\end{itemize}

A torus decomposition is \emph{minimal} if no proper subset of $S$ is a torus decomposition. We prove here the following.

\begin{teo}[JSJ decomposition] \label{JSJ:teo}
Let  $M$ be an orientable irreducible and $\partial$-irreducible compact 3-manifold with (possibly empty) boundary consisting of tori. A minimal torus decomposition for $M$ exists and is unique up to isotopy.
\end{teo}

Such a minimal decomposition is called the \emph{canonical torus decomposition} or the \emph{JSJ decomposition} of $M$. The canonical torus decomposition may be empty: this holds precisely when $M$ is itself a torus (semi-)bundle, Seifert, or simple.

\begin{oss}
Torus (semi-)bundles are closed: therefore if $M$ is not itself a torus (semi-)bundle, the blocks of its canonical decomposition are either Seifert or simple.
\end{oss}

\subsection{Existence and uniqueness}
Let  $M$ be an orientable irreducible and $\partial$-irreducible compact 3-manifold with (possibly empty) boundary consisting of tori. We now prove Theorem \ref{JSJ:teo}. Let us start by showing existence.
\begin{prop}
The manifold $M$ has a torus decomposition.
\end{prop}
\begin{proof}
Let $T_1,\ldots, T_k$ be a maximal set of disjoint non-parallel essential tori in $M$, which exists by Corollary \ref{incompressible:cor}. We now prove that $S=T_1 \sqcup \cdots \sqcup T_k$ is a torus decomposition.

Suppose it is not: one block $N$ of the decomposition is neither a (semi-)bundle, nor Seifert, nor simple. The block $N$ is irreducible and $\partial$-irreducible since these properties are preserved after cutting along incompressible surfaces. Being not simple, it contains an essential annulus $A$ or an essential torus $T$. 

In the latter case we can add $T$ to the family $T_1,\ldots, T_k$ and get a contradiction since $S$ is maximal. In the former case Lemma \ref{only:annuli:lemma} applies and $N$ is Seifert.
\end{proof}

Since $M$ has a torus decomposition, it certainly has a minimal one. We now prove that it is unique.

\begin{prop} \label{minimal:torus:unique:prop}
The manifold $M$ has a unique minimal torus decomposition up to isotopy.
\end{prop}
\begin{proof}
Let $S = T_1\sqcup \cdots \sqcup T_k$ and $S' = T_1' \sqcup \cdots \sqcup T_{k'}'$ be two minimal torus decompositions for $M$. We minimise their transverse intersections, so that $S\cap S'$ consists of essential circles cutting some tori into annuli. 

Let $T_i'$ be decomposed into some annuli. 
Each such annulus is essential in $M\setminus S$, hence it is contained in some non-simple block, \emph{i.e.}~a Seifert one. It is contained there horizontally or vertically: in the former case, the block is $\big(D, (2,1), (2,1) \big)$, $S_*\timtil S^1$, or $A\times S^1$ with $S_*$ the M\"obius strip. The first two blocks are diffeomorphic, and by swapping the fibration the annulus becomes vertical. The third block $T\times I$ is excluded since $S$ is minimal.

Now all annuli in $T_i'$ are vertical. Two consequent vertical annuli are separated by some torus $T_j$; since the two annuli are fibered, the fibers of the two Seifert blocks incident to $T_j$ are isotopic: hence the two blocks glue to a bigger Seifert block and $T_j$ can be removed, a contradiction since $S$ is minimal.

We have shown that $S\cap S' = \emptyset$. If $T_i$ is parallel to $T_j'$ we superpose the two tori, cut $M$ along $T_i=T_j'$ and proceed by induction. Now we suppose by contradiction that there is no parallelism.

Every $T_i'$ is an essential vertical torus in a Seifert block of $M\setminus S$, and vice versa. This easily implies that all the blocks in $M\setminus S$, $M\setminus S'$, and all their intersections are Seifert! Pick one such intersection. It has a unique Seifert fibration, unless it is $K\timtil I$ which may fiber in two ways. Since $\partial (K\timtil I)$ is connected, one block is $K\timtil I$ itself and we change the fibration on this block if necessary. Now all intersections and all blocks have unique fibrations and they all glue to a Seifert fibration for $M$, a contradiction.
\end{proof}

The proof of Theorem \ref{JSJ:teo} is complete.

\begin{oss}
The sphere decomposition of Theorem \ref{prime:teo} and the torus decomposition of Theorem \ref{JSJ:teo} differ in two aspects: (i) the set of decomposing spheres is \emph{not} canonical up to isotopy, while the set of tori is; (ii) on the other hand, after cutting along the spheres and capping off we get a canonical set of prime manifolds, whereas if we cut along the tori we get some canonical manifolds with toric boundaries, but there is no canonical way to cap them off. 
\end{oss}

\subsection{Geometric decomposition} \label{geometric:decomposition:subsection}
The \emph{geometric decomposition} is a slight variation of the canonical torus decomposition that is more suited to the geometrisation perspective that we will encounter in the next chapters. It is constructed from the canonical torus decomposition $S$ for $M$ as follows. Whenever a block $N$ of the torus decomposition is diffeomorphic to $K\timtil I$,  
we substitute the torus $\partial N$ in $S$ with the core Klein bottle $K$ of $K\timtil I$. This substitution has the effect of deleting $N$ from the list of blocks of the decomposition.\index{geometric decomposition}

The geometric decomposition consists of incompressible tori and Klein bottles. One reason for preferring the geometric decomposition to the canonical torus one is that it contains no interval bundles and every Seifert block has a unique fibration up to isotopy (because we have eliminated the $K\timtil I$ blocks). In particular we get the following easy criterion, whose proof is straightforward. 
\begin{prop} \label{criterion:geometric:prop}
Let  $M$ be an orientable irreducible and $\partial$-irreducible compact 3-manifold with (possibly empty) boundary consisting of tori. 
A non-empty torus decomposition of $M$ is the geometric one if and only if:
\begin{itemize}
\item every block is simple or Seifert with $\chi<0$,
\item the fibrations of two adjacent Seifert blocks do not match.
\end{itemize}
\end{prop}

Another reason for preferring the geometric decomposition is that every Seifert block can be geometrised with finite volume, whereas $K\timtil I$ needs infinite volume (this will be shown in the next chapter). 

Geometric decompositions also behave well under finite coverings. Let $\tilde M \to M$ be a finite covering between two orientable irreducible and $\partial$-irreducible compact 3-manifolds with (possibly empty) boundary consisting of tori. 

\begin{prop} \label{JSJ:lifts:prop}
The geometric decomposition of $\tilde M$ is the counterimage of that of $M$.
\end{prop}
\begin{proof}
The criterion of Proposition \ref{criterion:geometric:prop} lifts from $M$ to $\tilde M$.
\end{proof}

\subsection{Graph manifolds}
Waldhausen introduced in the 1960s a simple but non-trivial class of three-manifolds using only two blocks. Here $D$ and $P$ are the disc and the pair-of-pants.

\begin{defn}
A \emph{graph manifold} is any orientable three-manifold that decomposes along disjoint tori into pieces diffeomorphic to $D\times S^1$ or $P\times S^1$.\index{graph manifold}
\end{defn}

One may describe any such manifold via a graph with vertices of valence 1 and 3 representing the blocks, and some appropriate $2\times 2$-matrices labelling the edges telling the way the two incident blocks are glued. As usual, different graphs may represent the same manifold. 

\begin{ex}
Let $M$ be an orientable three-manifold. The following are equivalent:
\begin{enumerate}
\item $M$ is a graph manifold,
\item $M=M_1\#\ldots \#M_h$ for some prime manifolds $M_i$ whose geometric decompositions consist of Seifert manifolds or torus (semi-)bundles.
\end{enumerate} 
\end{ex}
\begin{proof}[Hint]
Use Corollary~\ref{DF:Seifert:cor} and Proposition \ref{criterion:geometric:prop}. 
\end{proof}

\subsection{References}
Most of the arguments contained in this chapter are well-known to three-dimensional topologists and can be found in many books. A standard introduction to knots and links is Rolfsen \cite{Ro}, and much more on Dehn surgeries and their four-dimensional interpretations is contained in Gompf -- Stipsicz \cite{GS}. A proof of the JSJ decomposition can be found in Hatcher \cite{H}, the original papers of Jaco -- Shalen and Johannson are \cite{JS} and \cite{Jo}. Fomenko -- Matveev \cite{FM} contains a chapter devoted to graph manifolds; the original paper of Waldhausen is \cite{Wa}.

%% file: Eight.tex
\chapter{The eight geometries} \label{eight:chapter}
We have concluded the previous chapter by defining a \emph{geometric decomposition} of three-manifolds along sphere and tori. The reason for using this terminology is the famous \emph{geometrisation conjecture}, proposed by Thurston in 1982 and proved by Perelman in 2003, which states that each of the blocks of the decomposition should be \emph{geometric}, in the sense that it may be equipped with a nice Riemannian metric.

There are eight nice Riemannian metrics available in dimension three. Three of them are the constant curvature ones (hyperbolic, elliptic, and flat), while the other five are some kind of (sometimes twisted) products of low-dimensional geometries. All these metrics are \emph{homogeneous}: distinct points have isometric neighbourhoods.

The Seifert manifolds studied in the previous chapters occupy precisely six of these eight geometries, and we analyse them in detail here.

\section{Introduction}
A connected Riemannian manifold $M$ is \emph{homogeneous} if for every $p,q\in M$ there is an isometry of $M$ sending $p$ to $q$, and is \emph{isotropic} if at each point $p$ every isometry of $T_pM$ is realised by an isometry of $M$. It is easy to prove that a complete isotropic manifold is also homogeneous and has constant sectional curvature: the fundamental examples of isotropic spaces are $S^n$, $\matR^n$, and $\matH^n$.\index{manifold!Riemannian manifold!homogeneous Riemannian manifold}\index{manifold!Riemannian manifold!isotropic Riemannian manifold}

The homogeneous condition alone (without isotropy) is more relaxed and produces manifolds that may not have constant sectional curvature.
We introduce in this chapter eight important homogeneous simply-connected complete Riemannian 3-manifolds:
$$S^3,\ \matR^3,\ \matH^3,\ S^2\times \matR,\ \matH^2\times \matR,\ \Nil,\ \Sol,\ \widetilde{\SL_2}.$$
The first three manifolds are also isotropic and have constant sectional curvature, the other five are not. 

Let $M$ be one of these eight model manifolds. We say that a Riemannian 3-manifold $N$ has a \emph{geometric structure modelled on $M$} if $N$ is locally isometric to $M$, that is if every point $p\in N$ has an open neighbourhood isometric to some open set in $M$. This implies that $N$ is \emph{locally homogeneous}: every two points $p,q\in N$ have isometric neighbourhoods $U(p) \isom U(q)$, both isometric to an $\varepsilon$-ball at any point of $M$.

If $N$ is complete, the developing map construction of Theorem \ref{simply:teo} applies also in this context and shows that $N = M/_\Gamma$ for some discrete group $\Gamma <\Iso(M)$ acting freely. \index{geometric structure}

\begin{table}
\begin{center}
\begin{tabular}{c||ccc}
       & \phantom{\Big|} $\chi>0$ & $\chi = 0$ & $\chi<0$ \\
 \hline \hline
 $e=0$ & \phantom{\Big|} $S^2 \times \matR$ & $\matR^3$ & $\matH^2 \times \matR$ \\
 $e\neq 0$ & \phantom{\big|} $S^3$ & $\Nil$ & $\widetilde {\SL_2}$
\end{tabular}
\vspace{.2 cm}
\nota{The closed manifolds modelled on six geometries are precisely the six commensurable classes of Seifert manifolds, distinguished by the numbers $e$ and $\chi$.}
\label{Seifert:geometry:table}
\end{center}
\end{table}

We now start a long journey through these eight geometries. The final goal of this chapter is to prove the following.

\begin{teo} \label{Seifert:geometrisation:teo}
A closed orientable 3-manifold has a geometric structure modelled on one of the following six geometries:
$$S^3,\ \matR^3,\ S^2\times \matR,\ \matH^2\times \matR,\ \Nil,\ \widetilde{\SL_2}$$
if and only if it is a Seifert manifold of the appropriate commensurability class, as prescribed by Table \ref{Seifert:geometry:table}. It has a $\Sol$ geometric structure if and only if it is a torus (semi-)bundle of Anosov type.
\end{teo}

It is a surprising (and maybe disappointing) fact that, despite its elegance and generality, the only known proof of this theorem available today works by investigating each geometry separately and carefully, often employing quite different techniques. We start with the elliptic case.

\section{Elliptic three-manifolds}
We start our journey by investigating elliptic 3-manifolds, that is manifolds modelled on $S^3$. We want to prove the following.

\begin{teo} \label{S3:teo}
A closed 3-manifold $M$ admits an elliptic metric if and only if it is a Seifert manifold  with $e\neq 0$ and $\chi>0$.
\end{teo}

An important ingredient of the proof is the complete classification of elliptic three-manifolds: to achieve this goal we need to study the isometries of $S^3$, and these are described elegantly via quaternions. 

\subsection{Unit quaternions}
We write as usual a \emph{quaternion} as\index{quaternion}
$$q= a+bi+cj+dk$$
with $a,b,c,d \in \matR$ and $i^2 = j^2 = k^2 =ijk = -1$. Quaternions form a non-commutative algebra, whose group structure is identified with $\matR^4$ by sending $q$ to $(a,b,c,d)$ and with $\matC^2$ by sending $q$ to $(a+bi, c+di)$. The norm is
$$|q| = \sqrt{a^2 + b^2 + c^2 + d^2}$$
and we have $|qq'| = |q||q'|$. The conjugate of $q$ is
$$\bar q = a-bi-cj-dk$$
and we have $|q|^2 = q\bar q$.
Unit quaternions are identified with $S^3$ and are closed under multiplication: this gives $S^3$ a Lie group structure. 
\begin{ex}
The centre of $S^3$ is $\{\pm 1\}$.
\end{ex}
Left or right multiplication by a fixed element $q\in S^3$ is an orientation-preserving isometry of both $\matR^4$ and $S^3$. We consider the homomorphism
\begin{align*}
\Psi\colon S^3 \times S^3 & \longrightarrow \Iso^+(S^3) = \SO(4) \\
(q_1, q_2) & \longmapsto \left\{x \longmapsto q_1xq_2^{-1}\right\}
\end{align*}

\begin{prop}
The homomorphism $\Psi$ is a degree-2 covering with kernel $\{\pm (1,1)\}$. Therefore it induces an isomorphism
$$\SO(4) = S^3\times S^3 /_{\{\pm (1,1)\}}.$$
\end{prop}
\begin{proof}
If $(q_1,q_2)$ lies in the kernel, by setting $x=1$ we get $q_1q_2^{-1}=1$ and so $q_1=q_2$. The general $x$ implies that $q_1=q_2$ lies in the centre and hence $(q_1,q_2) = \pm(1,1)$. Since $\SO(4)$ has the same dimension $6$ of $S^3\times S^3$ and is connected, we get a covering by Proposition \ref{Lie:covering:prop}.
\end{proof}

We now specialise $\Psi$ to the case $q=q_1=q_2$. The isometry
$$x \longmapsto qxq^{-1}$$
fixes $1$ and hence preserves the orthogonal 3-space generated by $i,j,k$, which we identify with $\matR^3$. We get a homomorphism
\begin{align*}
\Phi\colon S^3 & \longrightarrow \Iso^+(\matR^3) = \SO(3) \\
q & \longmapsto \left\{x \longmapsto qxq^{-1}\right\}
\end{align*}
\begin{prop}
The homomorphism $\Phi$ is a degree-2 covering with kernel $\pm 1$. Therefore it induces an isomorphism
$$\SO(3) = S^3/_{\{\pm 1\}}.$$
\end{prop}
\begin{proof}
The centre of $S^3$ is $\{\pm 1\}$ and $\SO(3)$ is connected and has dimension 3 like $S^3$, hence $\Phi$ is a covering.
\end{proof}
The real part of $q=a+bi+cj+dk$ is of course $a$.
\begin{cor}
Two elements $q,q' \in S^3$ are conjugate if and only if they have the same real part.
\end{cor}
\begin{proof}
Conjugations are isometries that fix the real axis and hence preserve the real part; rotations in $\SO(3)$ connect any two elements with the same real part.
\end{proof}
\begin{cor} \label{unit:complex:cor}
Every unit quaternion is conjugate to a unit complex quaternion $q=a\pm bi$, unique up to complex conjugation.
\end{cor}
The conjugacy classes in $S^3$ are the poles $+1$ and $-1$, and the parallel two-spheres between them. Let the \emph{imaginary two-sphere} be the maximal two-sphere consisting of all elements with zero real part.
\begin{cor}
The only element in $S^3$ of order two is $-1$. The elements of order four form the imaginary two-sphere and are all conjugate.
\end{cor}
\begin{proof}
An element of order four is conjugate to a complex one $a+bi$, which must be $\pm i$ and is hence purely imaginary.
\end{proof}

\begin{ex}
By sending the unit quaternion $(w,z)\in\matC^2$ to the matrix
$$\begin{pmatrix} w & z \\ -\bar z & \bar w \end{pmatrix}$$
we get a Lie group isomorphism between $S^3$ and $\SU(2)$.
\end{ex}

Note that $\bar q = q^{-1}$ on unit quaternions.
The inversion $q \mapsto q^{-1}$ is an orientation-reversing isometry of $S^3$. 

\subsection{Finite groups of quaternions}
We classify the finite subgroups of the Lie group $S^3$. Recall from Proposition \ref{subgroups:SO3:prop} that the finite subgroups of $\SO(3)$ up to conjugation are:
$$C_n, \quad D_{2m}, \quad T_{12} \isom A_4, \quad O_{24} \isom S_4, \quad I_{60} \isom A_5$$
with $n\geqslant 1$ and $m\geqslant 2$. Here $C_n$ is cyclic generated by a $\frac{2\pi}n$-rotation, and $D_{2m}$, $T$, $O$, $I$ are the orientation-preserving isometry groups of the regular $m$-prism, tetrahedron, octahedron (or cube), and icosahedron (or dodecahedron). The group $D_{2m}$ is the dihedral group. The subscript always indicates the order of the group (except in the alternating $A_n$ and symmetric $S_n$ that have order $\frac{n!}2$ and $n!$).

\begin{table}
\begin{center}
\begin{tabular}{c||cc}
\phantom{\Big|} \!\! & name        & elements \\
 \hline \hline
\phantom{\Big|} \!\!  $C_n$ & cyclic & $\big\{e^{\frac{2a\pi i}n}\big\}_{a=1,\ldots, n}$ \\
\phantom{\Big|} \!\!  $D_{4m}^*$ & binary dihedral & $\big\{e^{\frac{a\pi i}m}$, $e^{\frac{a\pi i}m}j\big\}_{a=1,\ldots, 2m}$  \\
\phantom{\Big|} \!\!  $T_{24}^*$ & binary tetrahedral & $\big\{\pm 1 ,\ \pm\, i,\ \pm\, j,\ \pm\, k,\ \frac 12 \big(\!\pm 1 \pm i \pm j \pm k\big)\big\}$ \\
\phantom{\Big|} \!\!  $O_{48}^*$ & binary octahedral & 
$T_{24}^* \cup \big\{\frac {\pm 1 \pm i}{\sqrt 2},\ \frac {\pm 1 \pm j}{\sqrt 2}, \ldots, \frac {\pm i \pm k}{\sqrt 2},\ \frac {\pm j \pm k}{\sqrt 2} \big\}$ \\
\phantom{\Big|}  \!\!  $I_{120}^*$ & binary icosahedral &
$T_{24}^* \cup \big\{\pm \frac 1 2i \pm \frac{\sqrt 5 -1}4 j \pm \frac{\sqrt 5 + 1}4 k, \ldots \big\} $ \\
\end{tabular}
\vspace{.2 cm}
\nota{Every finite subgroup of $S^3$ is conjugate to one of these groups. 
The group $O_{48}^*$ consists of $T_{24}^*$ and the $6\times 4 = 24$ numbers obtained from $\frac {\pm 1 \pm i}{\sqrt 2}$ by permuting the elements $1,i,j,k$ and varying the signs. The group $I^*_{120}$ consists of $T_{24}^*$ and the $12 \times 8 = 96$ numbers obtained from $\pm \frac 12i \pm \frac{\sqrt 5 -1}4 j \pm \frac{\sqrt 5 + 1}4 k$ by permuting the elements $1,i,j,k$ with an \emph{even} permutation and varying the signs.}
\label{binary:table}
\end{center}
\end{table}

We defined above a degree-2 covering $\Phi\colon S^3 \to SO(3)$. Let 
$$D_{4m}^*, \quad T_{24}^*, \quad O_{48}^*, \quad I_{120}^*$$ 
be the counterimages of $D_{2m}$, $T_{12}$, $O_{24}$, $I_{60}$ along $\Phi$. They are called the \emph{binary} dihedral, tetrahedral, octahedral, and icosahedral group. 
We now classify the finite subgroups of $S^3$.\index{group!binary group}

\begin{prop}
Every finite subgroup of $S^3$ is either binary or cyclic. The finite subgroups up to conjugation are listed in Table \ref{binary:table}.
\end{prop}
\begin{proof}
If $G<S^3$ is finite then $\Phi(G) = C_n, D_{2m}, T, O$, or $I$ up to conjugation. If $\Phi^{-1}(\Phi(G)) = G$ we are done (the counterimage of $C_n$ is still cyclic: exercise). This holds precisely when $-1 \in G.$

If $G$ has even order, it contains an order-2 element which is necessarily $-1$ and we are done. If $G$ has odd order then $\Psi(G)$ has odd order and is cyclic, hence $G$ is contained in the cyclic $\Psi^{-1}(\Psi(G))$ and is hence cyclic.
\end{proof}

The smallest non-abelian group in the list is the \emph{quaternion group} $Q_8 = D_8^*$ consisting of the \emph{Lipschitz units}\index{group!quaternion group}\index{Lipschitz unit} 
$$Q_8 = \{ \pm 1, \pm i, \pm j, \pm k \}.$$
The quaternion group is contained in $T_{24}^*$ as a normal subgroup, giving an exact sequence
\begin{equation} \label{Q8:eqn}
0 \longrightarrow Q_8 \longrightarrow T_{24}^* \longrightarrow \matZ_3 \longrightarrow 0.
\end{equation}
The binary dihedral group $D_{4m}^*$ contains the index-two cyclic group $C_{2m}$, giving an exact sequence
\begin{equation} \label{D4n:eqn}
0 \longrightarrow C_{2m} \longrightarrow D_{4m}^* \longrightarrow \matZ_2 \longrightarrow 0.
\end{equation}

\begin{rem}
Every group $\Gamma$ in Table \ref{binary:table} acts freely and isometrically on $S^3$ by right multiplication, hence $S^3/_\Gamma$ is an elliptic manifold with fundamental group $\Gamma$. The only perfect group in the list (\emph{i.e.}~with trivial abelianisation) is $I^*_{120}$ and we will soon see that $S^3/_{I^*_{120}}$ is the ubiquitous Poincar\'e homology sphere $\Sigma(2,3,5)$, defined in Section \ref{Seifert:homology:sphere:subsection}.
\end{rem}

\subsection{Regular polytopes} \label{regular:polytopes:subsection}
There are six regular polytopes in dimension four, listed in Table \ref{polytopes:table}, and the finite groups of $S^3$ can be used to describe five of them. The convex hull of $Q_8$ in $\matR^4$ is the \emph{16-cell}, whose dual is the \emph{hypercube}. 
The group $T_{24}^*$ consists of the vertices of the 16-cell and of the dual hypercube altogether.\index{regular polytope} 

\begin{table}
\begin{center}
\begin{tabular}{c||ccccc}
\phantom{\Big|} \!\! name & vertices & edges & faces & facets & Schl\"afli \\
 \hline \hline
\phantom{\Big|} \!\!  simplex & 5 & 10 & 10 & 5 tetrahedra & $\{3,3,3\}$ \\
\phantom{\Big|} \!\!  hypercube & 16 & 32 & 24 & 8 cubes & $\{4,3,3\}$ \\
\phantom{\Big|} \!\!  $16$-cell & 8 & 24 & 32 & 16 tetrahedra & $\{3,3,4\}$ \\
\phantom{\Big|} \!\!  $24$-cell & 24 & 96 & 96 & 24 octahedra & $\{3,4,3\}$ \\
\phantom{\Big|} \!\!  $120$-cell & 600 & 1200 & 720 & 120 dodecahedra & $\{5,3,3\}$ \\
\phantom{\Big|} \!\!  $600$-cell & 120 & 720 & 1200 & 600 tetrahedra & $\{3,3,5\}$
\end{tabular}
\vspace{.2 cm}
\nota{The six regular polytopes in dimension four, with their Schl\"afli notation (see Section \ref{regular:subsection}). The groups $Q_8$, $T_{24}^*$, and $I_{120}^*$ consist of the vertices of the 16-cell, 24-cell, and 600-cell respectively.}
\label{polytopes:table}
\end{center}
\end{table}

The convex hull of $T_{24}^*$ is the \emph{24-cell}: this is the unique self-dual regular polytope in all dimensions $n\geqslant 3$ which is not a simplex! See Figure \ref{24cell:fig}. The group $O_{48}^*$ consists of the vertices of the 24-cell and of its dual. Finally, the convex hull of $I_{120}^*$ is the \emph{600-cell}, whose dual is the \emph{120-cell}. 

\begin{figure}
\begin{center}
\includegraphics[width = 7 cm] {\iftoggle{BW}{Stereographic_polytope_24cell_faces-BW}{768px-Stereographic_polytope_24cell_faces}}
\nota{A stereographic projection of the tessellation of $S^3$ induced by the 24-cell. Its $24$ vertices form the binary tetrahedral group $T_{24}^*$. Its facets are 24 regular octahedra.}
\label{24cell:fig}
\end{center}
\end{figure}

\subsection{Classification of elliptic 3-manifolds}
We now want to classify all the elliptic 3-manifolds. We start with a linear algebra exercise.
\begin{ex}
Let $A \in \On(4)$ be an orthogonal transformation. Then
\begin{itemize}
\item if $\det A = -1$, it fixes pointwise a line,
\item if $\det A = +1$, it preserves two orthogonal planes and acts as a rotation on each.
\end{itemize}
\end{ex}

Every elliptic 3-manifold is isometric to $S^3/_\Gamma$ for some finite subgroup $\Gamma < \On(4)$ acting freely. The exercise implies that $\Gamma < \SO(4)$, since orientation-reversing elements fix some point in $S^3$. Therefore we get:

\begin{cor}
Every elliptic 3-manifold is orientable.
\end{cor}

We have already encountered the lens space $L(p,q) = S^3/_\Gamma$ where the cyclic group $\Gamma = C_{p,q}$ is generated by the isometry 
$$(z,w) \longmapsto (\omega z, \omega^q w)\quad {\rm with}\  \omega = e^\frac{2\pi i}p.$$ 
Another consequence of the exercise is the following. 
\begin{cor} \label{abelian:SO4:cor}
Every finite abelian group $\Gamma < \SO(4)$ acting freely is conjugate to $C_{p,q}$ for some coprime $p,q$.
\end{cor}
\begin{proof}
An abelian $\Gamma$ fixes two orthogonal planes $U$ and $V$ and acts on each as rotations. The restriction map $\Gamma \to \Iso^+(U)$ is injective because $\Gamma$ acts freely on $S^3$. Therefore $\Gamma$ is cyclic, generated by an element that acts on both planes by rotations of order $p=|\Gamma|$, so $\Gamma$ is conjugate to $C_{p,q}$.
\end{proof}

We now want to construct non-abelian examples. It is convenient to look at $\SO(4)$ as the image of $S^3\times S^3$ via the surjective map $\Psi$ defined above. Which elements of $S^3\times S^3$ act freely on $S^3$? The answer is particularly simple.
\begin{prop} \label{conj:fixed:prop}
The isometry $\Psi(q_1,q_2)$ of $S^3$ has a fixed point if and only if $q_1$ and $q_2$ are conjugate in $S^3$.
\end{prop}
\begin{proof}
We have $q_1xq_2^{-1}=x$ for some $x\in S^3$ if and only if $q_1 = xq_2x^{-1}$, that is $q_1$ and $q_2$ are conjugate.
\end{proof}
The following corollary will help to rule out many cases.
\begin{cor} \label{no:four:cor}
If $q_1$ and $q_2$ have both order four, the isometry $\Psi(q_1,q_2)$ has a fixed point.
\end{cor}
We now construct more examples.
\begin{prop} If two finite subgroups $G,H < S^3$ have coprime orders, the image $\Psi(G \times H)$ acts freely on $S^3$.
\end{prop}
\begin{proof}
If $(g_1, g_2) \in G \times H$ is non-trivial, the elements $g_1$ and $g_2$ have coprime and hence distinct orders, thus they are not conjugate.
\end{proof}
\begin{cor}
If $G$ is a group from Table \ref{binary:table} and $n$ is coprime with the order of $G$, then $\Psi(G\times C_n)$ acts freely on $S^3$.
\end{cor}
This corollary produces most of the non-cyclic subgroups of $\SO(4)$ acting freely on $S^3$, but not all of them! Table \ref{finite:freely:table} lists all the finite subgroups of $\SO(4)$ acting freely, up to conjugation in $\On(4)$. The first line shows the cyclic groups $C_{p,q}$, and the other families of groups are all of type $\Psi(G\times C_n)$, except the ones in the third and fifth line which we now explain.

\begin{table}
\begin{center}
\begin{tabular}{c||cccc}
\phantom{\Big|}\!\! $\Gamma$      & conditions & $|\Gamma|$ & $Z$ & $\Gamma/_Z$ \\
 \hline \hline
\phantom{\Big|}\!\! $C_{p,q}$ & $p>0$, $(p,q)=1$ & $p$ & $C_p$ & $\{e\}$ \\
\phantom{\Big|}\!\! $\Psi\big(D_{4m}^* \times C_n\big)$ & $m>1 ,n>0$, $(4m,n)=1$ & $4mn$ & $C_{2n}$ & $D_{2m}$ \\
\phantom{\Big|}\!\! $\Gamma \underset 2< \Psi\big(D_{4m}^* \times C_{4n}\big)$ & $m>1,n>0$, $n$ even, $(m,n)=1$ & $4mn$ & $C_{2n}$ & $D_{2m}$ \\
\phantom{\Big|}\!\! $\Psi\big(T_{24}^* \times C_n\big)$ & $n>0$, $(24,n)=1$ & $24n$ & $C_{2n}$ & $T_{12}$ \\
\phantom{\Big|}\!\! $\Gamma \underset 3< \Psi\big(T_{24}^* \times C_{6n}\big)$ & $n>0$, $n$ odd, $3|n$ & $24n$ & $C_{2n}$ & $T_{12}$ \\
\phantom{\Big|}\!\! $\Psi\big(O_{48}^* \times C_n\big)$ & $n>0$, $(48,n)=1$ & $48n$ & $C_{2n}$ & $O_{24}$ \\
\phantom{\Big|}\!\! $\Psi\big(I_{120}^* \times C_n\big)$ & $n>0$, $(120,n)=1$ & $120n$ & $C_{2n}$ & $I_{60}$
\end{tabular}
\vspace{.2 cm}
\nota{The finite subgroups $\Gamma < \SO(4)$ that act freely on $S^3$, up to conjugation in $\On(4)$. For each we show its centre $Z$ and the quotient $\Gamma/_Z$. If $\Gamma$ is not cyclic, the two groups $Z$ and $\Gamma/_Z$ determine $\Gamma$ up to conjugation (note that $n$ satisfies some congruence equality which separates lines 2,3 and 4,5).
The symbol $\Gamma <_i G$ indicates that $\Gamma$ has index $i$ in $G$.}
\label{finite:freely:table}
\end{center}
\end{table}

\begin{rem} \label{image:product:rem}
For all groups in Table \ref{finite:freely:table} of type $G\times C_n$, we have $\Psi(G\times C_n) \isom G\times C_n$ because $n$ is odd and hence $(-1,-1)\not\in G\times C_n$.
\end{rem}

In the third line of Table \ref{finite:freely:table}, the group $\Gamma$ is an index-two subgroup of $\Psi(D_{4m}^* \times C_{4n})$. It is the image along $\Psi$ of the kernel of the map
\begin{align*}
D^*_{4m} \times C_{4n} & \longrightarrow \matZ_2 \\
(g_1, g_2) & \longmapsto f_1(g_1) + f_2(g_2)
\end{align*}
where $f_1$ is the map in (\ref{D4n:eqn}) and $f_2$ is the surjective homomorphism $C_{4n} \to \matZ_2$. Analogously, the group $\Gamma$ in the fifth line is an index-three subgroup of $\Psi(T_{24}^* \times C_{6n})$: it is the image along $\Psi$ of the kernel of the map
\begin{align*}
T^*_{24} \times C_{6n} & \longrightarrow \matZ_3 \\
(g_1, g_2) & \longmapsto f_1(g_1) + f_2(g_2)
\end{align*}
where $f_1$ is the map in (\ref{Q8:eqn}) and $f_2$ is any surjective homomorphism $C_{6n} \to \matZ_3$. Both $f_1$ and $f_2$ are well-defined only up to automorphisms of $\matZ_3$, but the kernel is well-defined up to conjugation in $S^3\times S^3$ anyway: there are four possibilities and they are related by conjugations via the elements $(1,j)$ and $\big(\frac{1+i}{\sqrt 2},1\big)$, as one can easily check (these conjugations permute the cosets of $C_{2n} \triangleleft C_{6n}$ and $Q_8 \triangleleft T_{24}^*$).

\begin{prop} \label{finite:SO4:prop}
The finite subgroups of $SO(4)$ shown in Table \ref{finite:freely:table} act freely on $S^3$.
\end{prop}
\begin{proof}
We only need to prove this for the groups $\Gamma$ that belong to the third and fifth family. Concerning the third, pick an element $(g_1,g_2) \in S^3\times S^3$ that projects to a non-trivial element in $\Gamma$. We want to prove that $g_1$ and $g_2$ have distinct orders and are hence non-conjugate.

By hypothesis $f_1(g_1) + f_2(g_2)=0$ in $\matZ_2$ and hence $f_1(g_1) = f_2(g_2)$ is either 0 or 1. In the first case we get $g_1 \in C_{2m}$ and $g_2 \in C_{2n}$ which have distinct orders since $(m,n)=1$ and $(g_1,g_2) \neq \pm (1,1)$. In the second case $g_1 \in D_{4m}^* \setminus C_{2m}$ has order four (check from Table \ref{binary:table}) and $g_2 \in C_{4n} \setminus C_{2n}$ does not have order four since $n$ is even. 

The fifth family is treated similarly. Pick $(g_1, g_2) \in \Gamma$ with non-trivial image in $SO(4)$. If $f_1(g_1) = f_2(g_2) = 0$ then $g_1 \in Q_8$ has order 1, 2, or 4 and $g_2 \in C_{2n}$ does not have order 4 since $n$ is odd. If $f_1(g_1) = - f_2(g_2) \neq 0$ then $g_1\in T_{24}^* \setminus Q_8$ has order 3 or 6 (check from Table \ref{binary:table}) while $g_2 \in C_{6n} \setminus C_{2n}$ cannot have order 3 or 6 since $3$ divides $n$.
\end{proof}

We now show that Table \ref{finite:freely:table} exhausts all possibilities.

\begin{prop}
Every finite subgroup of $SO(4)$ acting freely on $S^3$ is conjugate in $\On(4)$ to one in Table \ref{finite:freely:table}.
\end{prop}
\begin{proof}
Let a finite $\Gamma < SO(4)$ act freely on $S^3$. We consider its counterimage $G = \Psi^{-1}(\Gamma)$ in $S^3\times S^3$. 

We note that the orientation-reversing isometry $q \mapsto q^{-1}$ of $S^3$ conjugates $\Psi(q_1,q_2)$ to $\Psi(q_2,q_1)$, since
$$x \longmapsto (q_1x^{-1}q_2^{-1})^{-1} = q_2xq_1^{-1}.$$
Therefore via conjugation in $\On(4)$ we may swap the factors of $S^3\times S^3$.

If $G$ is a product $G=G_1 \times G_2$, then $G_1$ and $G_2$ are some groups from Table \ref{binary:table}. If they are both cyclic, then $\Gamma$ is abelian and we conclude by Corollary \ref{abelian:SO4:cor}.
If they are both non-cyclic, they both contain order-4 elements and we get a contradiction from Proposition \ref{no:four:cor}.

Therefore $G = G_1 \times C_k$ with $G_1$ equal to $D_{4m}^*$, $T_{24}^*$, $O_{48}^*$, or $I_{120}^*$, and $k$ not divisible by four. Since $(-1,-1) \in G$ we get $k=2n$, and $n$ odd implies that $\Psi (G_1 \times C_n) = \Psi(G_1 \times C_{2n})$. Moreover $n$ is coprime with the order of $G_1$ otherwise Proposition \ref{conj:fixed:prop} would easily give a fixed point. Therefore we get a product group as in Table \ref{finite:freely:table}.

We are left to consider the case $G$ is not a product. Let $G_1 \times G_2$ (resp.~$G_1' \times G_2'$) be the smallest (resp.~biggest) product subgroup of $S^3\times S^3$ such that 
$$G_1' \times G_2' < G < G_1 \times G_2.$$
The subgroup $G_1$ (resp.~$G_1'$) consists of all $g_1\in S^3$ such that $(g_1,g_2) \in G$ for some $g_2$ (resp.~such that $(g_1,1)\in G$).
It is easy to check that $G_i'\triangleleft G_i$ and 
$${G}/_{G_1' \times G_2'} \isom {G_1}/_{G_1'} \isom {G_2}/_{G_2'}.$$
Each $G_i$ and $G_i'$ is finite and hence conjugate to one in Table \ref{binary:table}. 
If both $G_1$ and $G_2$ are cyclic, the group $\Gamma$ is abelian and we are done. We henceforth suppose that $G_1$ and $G_2$ are not both cyclic.

By the same reasoning above, up to reordering $G_2'$ contains no order-4 elements and hence $G_2' = C_s$ with $s$ not divisible by four. 

We now prove that $G_1'$ contains all the order-4 elements of $G_1$.
To do so we pick an order-4 element $g_1 \in G_1$ and prove that $\pm g_1 \in G_1'$: this suffices since $g_1 = (-g_1)^3$. We have $(g_1,g_2) \in G$ for some $g_2$ of some order $2^tk$ with $k$ odd. We have $(g_1^k,g_2^k) = (\pm g_1, g_2^k)$, so up to replacing $g_2$ with $g_2^k$ we may suppose $k=1$ and $g_2$ has order $2^t$. 

We have $(1,g_2^4) = (g_1^4,g_2^4) \in G$ and hence $g_2^4 \in G_2'$. Suppose $s$ is odd. Therefore $g_2^4 \in G_2' = C_s$ is trivial and $g_2$ has order $2^t$; thus $g_2$ has order $1,2,$ or $4$. It cannot have order $4$ by Corollary \ref{no:four:cor}, therefore $g_2= \pm 1$ and hence $\pm g_1 \in G_1'$ as required.

If $s$ is even we get the same conclusion: now $(-1,g_2^2) = (g_1^2,g_2^2) \in G$ hence $(1,-g_2^2) \in G$ gives $-g_2^2 \in G_2'$ and thus $g_2^2 \in G_2'$ since $n$ is even. Again this implies that $g_2$ has order 1,2, or 4 and we conclude as above.

We have proved that $G_1'$ contains all order-4 elements of $G_1$. The groups $D_{4n}^*$, $O_{48}^*$, and $I_{120}^*$ are generated by their order-4 elements (exercise), hence $G_1' \neq G_1$ implies that $G_1$ is not one of them. Then it is either cyclic or $T_{24}^*$.

If $G_1= C_h$ is cyclic, then $G_1'$ and $G_1/_{G_1'} \isom G_2/_{G_2'}$ also are. By assumption $G_2$ is not cyclic, but it contains the cyclic $G_2'$ with non-trivial cyclic quotient: the only possibility from Table \ref{binary:table} is that $G_2 = D_{4m}^*$ and $G_2' = C_{2m}$. Therefore $G$ is an index-two subgroup of $G_1\times G_2 = C_h \times D_{4m}^* $ as in the third line of Table \ref{finite:freely:table}. 

We need to prove that $h=4n$ with $n$ even and $(n,m)=1$. Since $G_1' < G_1$ has index two, $h$ is even. Since $(1,-1) \in G_1' \times C_{2m}< G$, we also have $(-1,1)\in G$ and hence $-1\in G_1'$. Therefore $|G_1'|$ is even and $4|h$. Since $G_1'$ contains all order-4 elements of $G_1$, four divides $|G_1'|$ and hence $8|h$. Therefore $h=4n$ with $n$ even, as required.
Moreover $(m,n)=1$ since $G_1' \times G_2' = C_{2n} \times C_{2m} \subset G$.

If $G_1 = T_{24}^*$ the order-4 elements generate the index-three subgroup $Q_8$ and hence $G_1' = Q_8$ and $G_2/_{G_2'} \isom G_1/_{G_1'} \isom \matZ_3$. Recall that $G_2' = C_s$ with $s$ not divisible by four. As above, $-1 \in G_1'$ implies $-1 \in G_2'$ and hence $s = 2n$ with $n$ odd. Therefore $G_2$ has order $6n$ and must hence be cyclic. 

The group $G$ is an index-3 subgroup of $G_1\times G_2 = T_{24}^* \times C_{6n} $ as in the third line of Table \ref{finite:freely:table}. We must have $3|n$ otherwise $G$ contains an element $(g_1,g_2)$ with both $g_1,g_2$ having order 3 (see the proof of Proposition \ref{finite:SO4:prop}).
\end{proof}

We can finally summarise our discoveries:

\begin{cor} 
Table \ref{finite:freely:table} lists all the finite subgroups $\Gamma < SO(4)$ acting freely on $S^3$ up to conjugation in $\On(4)$, without repetitions. 
\end{cor}
\begin{proof}
There are no repetitions because the non-cyclic groups listed are all non-isomorphic: the isomorphism types of the centre $Z$ and $\Gamma/_Z$ suffice to determine $\Gamma$ in that list, see Table \ref{finite:freely:table} (note that $n$ satisfies some congruence equality which separates the lines 2,3 and 4,5).

To compute $|\Gamma|$, we use Remark \ref{image:product:rem} and note that $(-1,-1)$ belongs to $(D_{4m}^* \times C_{4n})$ and $(T_{24}^* \times C_{6n})$.

To compute $Z$ and $\Gamma/_Z$, we note that in all the non-cyclic cases $\Gamma$ is the image of a subgroup in some product $G^* \times C_{kn}$ which is ``diagonal'', in the sense that it maps surjectively to both factors $G^*$ and $C_{kn}$. The map onto $G^*$ pushes-forward to a surjection $\Gamma \to G$, whose kernel is easily detected as being the centre $Z$ and isomorphic to $C_{2n}$ in all cases.
\end{proof}

\subsection{Seifert fibrations of elliptic 3-manifolds}
We now turn back to Seifert manifolds and prove Theorem \ref{S3:teo}. Let $S^1\subset S^3$ be the unit complex numbers. Recall that $S^3\times S^3$ acts on $S^3$ via $\Psi$. The following proposition is crucial: it says that a big Lie subgroup of $S^3\times S^3$ preserves the Hopf fibration; later on we will discover that every finite subgroup $\Gamma$ of $S^3\times S^3$ acting freely may be conjugated into this big group and hence every quotient $S^3/_\Gamma$ inherits a Seifert structure from the Hopf fibration, hence every elliptic three-manifold is Seifert. 

\begin{prop}
The group $S^1\times S^3$ preserves the Hopf fibration.
\end{prop}
\begin{proof}
Represent quaternions as pairs $(z_1,z_2)$ of complex numbers. The Hopf fibration $S^3 \to \matCP^1$ is the map $(z_1,z_2) \to [z_1,z_2]$. 
One checks easily that quaternion multiplication acts as follows:
$$(z_1, z_2) (w_1,w_2) = (z_1w_1 - z_2 \bar{w}_2, z_2 \bar{w}_1 + z_1w_2).$$
Therefore right-multiplication by $(w_1,w_2)\in S^3$ acts $\matC$-linearly on $\matC^2$ and hence preserves the Hopf fibration.

Left multiplication by $(w_1,w_2)$ is \emph{not} $\matC$-linear in general, but it is so when $(w_1,w_2) = (w_1,0) \in S^1$, since $(w_1,0)(z_1,z_2) = (w_1z_1, w_1z_2).$
Therefore $S^1\times S^3$ preserves each fibre of the Hopf fibration.
\end{proof}

We can finally prove Theorem \ref{S3:teo}.
\begin{teo}
A closed 3-manifold $M$ admits an elliptic metric if and only if it is a Seifert manifold  with $e\neq 0$ and $\chi>0$.
\end{teo}
\begin{proof}
If $M$ is elliptic then $M=S^3/_\Gamma$ for some finite subgroup $\Gamma < SO(4)$ acting freely. If $\Gamma$ is conjugate to $C_{p,q}$ we get a lens space and we are done. Otherwise Table \ref{finite:freely:table} shows that $\Gamma$ is conjugate to the image of a subgroup of $S^1\times S^3$ (because $C_n < S^1$). Therefore up to conjugation $\Gamma$ preserves the Hopf fibration of $S^3$, which descends to a Seifert fibration. We have $\chi>0$ and $e\neq 0$ because the universal cover of $M$ is $S^3$.

Conversely, we now show that every Seifert manifold with $\chi >0$ and $e \neq 0$ arises as an elliptic manifold $S^3/_\Gamma$. The Seifert manifolds were listed in Table \ref{elliptic:table}, and the corresponding $\Gamma$ is shown in Table \ref{elliptic:Seifert:table}.
To verify the correspondence, note that the centre $Z$ of $\Gamma = \pi_1(M)$ and the quotient $\Gamma/_Z$ are shown in Table \ref{finite:freely:table} and they fully determine $\Gamma$. Proposition \ref{elliptic:center:prop} says that $Z<\pi_1(M)$ is generated by a regular fibre and $\Gamma/_Z$ is the orbifold fundamental group of the base.
\end{proof}

\begin{table}
\begin{center}
\begin{tabular}{c||c||c}
\phantom{\Big|} fibration        & condition & $\Gamma$ \\
 \hline \hline
\phantom{\Big|} $\big(S^2,(2,1),(2,1),(p,q)\big)$ & ${p+q}>0$ odd &
$\Psi\big(D_{4p}^* \times C_{p+q}\big)$  \\
\phantom{\Big|} $\big(S^2,(2,1),(2,1),(p,q)\big)$ & ${p+q}>0$ even &
$\Gamma \underset 2< \Psi\big(D_{4p}^* \times C_{4(p+q)}\big)$ \\
\phantom{\Big|} $\big(S^2,(2,1),(3,1),(3,q)\big)$ & $3 \nmid 5+2q$ &
$\Psi\big(T_{24}^* \times C_{|5+2q|}\big)$ \\
\phantom{\Big|} $\big(S^2,(2,1),(3,1),(3,q)\big)$ & $3|5+2q>0$ &
$\Gamma \underset 3< \Psi\big(T_{24}^* \times C_{6(5+2q)}\big)$ \\
\phantom{\Big|} $\big(S^2,(2,1),(3,1),(4,q)\big)$ & &
$\Psi\big(O_{48}^* \times C_{|10+3q|}\big)$ \\
\phantom{\Big|} $\big(S^2,(2,1),(3,1),(5,q)\big)$ & &
$\Psi\big(I_{120}^* \times C_{|25+6q|}\big)$
\end{tabular}
\vspace{.2 cm}
\nota{The non-lens elliptic manifolds, listed without repetitions. For each we show its Seifert fibration and fundamental group $\Gamma < \SO(4)$. Here $p>1$. The integer $q$ may be negative, if not forbidden explicitly. In the fourth line $3$ divides $5+2q$, which must be positive. }
\label{elliptic:Seifert:table}
\end{center}
\end{table}

In particular we have the following:
\begin{align*}
\pi_1\big(S^2,(2,1),(2,1),(2,-1)\big) & = Q_8, \\
\pi_1\big(S^2,(2,1),(2,1),(p,1-p)\big) & = D_{4p}^*, \\
\pi_1\big(S^2,(2,1),(3,1),(3,-2)\big) & = T_{24}^*, \\
\pi_1\big(S^2,(2,1),(3,1),(4,-3)\big) & = O_{48}^*, \\
\pi_1\big(S^2,(2,1),(3,1),(5,-4)\big) & = I_{120}^*. 
\end{align*}

The latter is Poincar\'e's homology sphere, the unique elliptic three-manifold with perfect fundamental group (see Section \ref{Seifert:homology:sphere:subsection}).\index{Poincar\'e homology sphere}

\section{Flat three-manifolds}
We now turn to flat three-manifolds. In dimension two, every orientation-preserving isometry of $\matR^2$ is a translation, and this easily implies that every flat orientable closed surface is a torus. In dimension three we also have \emph{rototranslations}, which produce more orientable manifolds. 

We now classify all the closed flat orientable three-manifolds up to diffeomorphism and prove the following.
\begin{teo} \label{main:flat:teo}
A closed orientable 3-manifold $M$ admits a flat metric if and only if it is a Seifert manifold with $e=\chi=0$.
\end{teo}
There are six such manifolds, listed in Table \ref{seven:table}. 

\subsection{Classification}
Every closed flat 3-manifold is isometric to $\matR^3/_\Gamma$ for some crystallographic group $\Gamma< \Iso(\matR^3)$ acting freely, see Section \ref{crystallographic:subsection}.

\begin{ex}
Every element in $\Gamma$ is either a translation or a rototranslation (defined in Example \ref{rototranslation:example}).
\end{ex}

We prove one half of Theorem \ref{main:flat:teo}.

\begin{prop} \label{flat:to:Seifert:prop}
Every closed orientable flat 3-manifold is a Seifert manifold with $e=\chi=0$.
\end{prop}
\begin{proof}
We have $M=\matR^3/_\Gamma$. Recall the exact sequence
$$0 \longrightarrow H \longrightarrow \Gamma \longrightarrow r(\Gamma) \longrightarrow 0$$
where $H \triangleleft \Gamma$ is the translation subgroup and $r(\Gamma)<SO(3)$ is finite by Proposition \ref{finite:image:r:prop}. We now prove that $\Gamma$ preserves a foliation of $\matR^3$ into parallel lines that projects to a Seifert structure on $M$.

If $r(\Gamma)$ is trivial, then $\Gamma = H$ consists of translations and preserves many foliations into parallel lines that project to a Seifert structure on the quotient 3-torus $M$. If $r(\Gamma)$ is non-trivial, it is isomorphic to $C_n$, $D_{2m}$, $T_{12}$, $O_{24}$, or $I_{60}$. If $r(\Gamma) = C_n$ or $D_{2m}$, it has a common fixed vector line $l\subset \matR^3$ and $\Gamma$ preserves the foliation of lines parallel to $l$.

If $r(\Gamma) = T_{12}, O_{24}$, or $I_{60}$ we obtain a contradiction as follows. In all cases we have $T_{12}\subset r(\Gamma)$. The group $T_{12}$ consists of the identity, the $\pi$-rotations along the three coordinate axis, and the $\pm\frac {2\pi}3$-rotations along the axis spanned by $(1,1,1), (1,-1,-1), (-1,1,-1)$, and $(1,-1,-1)$. 

Pick a rototranslation $h\in \Gamma$ with axis $l$ parallel to $(1,1,1)$, and up to conjugating $\Gamma$ by a translation we may suppose that $l$ contains the origin $0\in\matR^3$. We have $h(0) = (d,d,d)$ for some $d\neq 0$. Since $h^3$ is a translation we get $(3d,3d,3d) \in H$. 

The group $\Gamma$ and hence $T_{12}$ acts on $H$ via conjugation. Therefore $H$ is $T_{12}$-symmetric and $(3d,-3d,-3d)\in H$ using the $\pi$-rotation along the first axis. Hence $t=(6d,0,0) = (3d,3d,3d) + (3d,-3d,-3d) \in H$. The composition $t\circ h^2$ has a fixed point, because it sends $0$ to $(2d,2d,2d)-(6d,0,0) = (-4d,2d,2d)$ which is orthogonal to $l$: a contradiction.

In all cases $M$ has a Seifert structure. By Bieberbach's Theorem (stated as Corollary \ref{Bieberbach:cor}) the manifold $M$ is finitely covered by the 3-torus, and hence $\chi=e=0$ by Proposition \ref{cover:Seifert:prop}.
\end{proof}

We now prove the other half.

\begin{figure}
\begin{center}
\includegraphics[width = 10 cm] {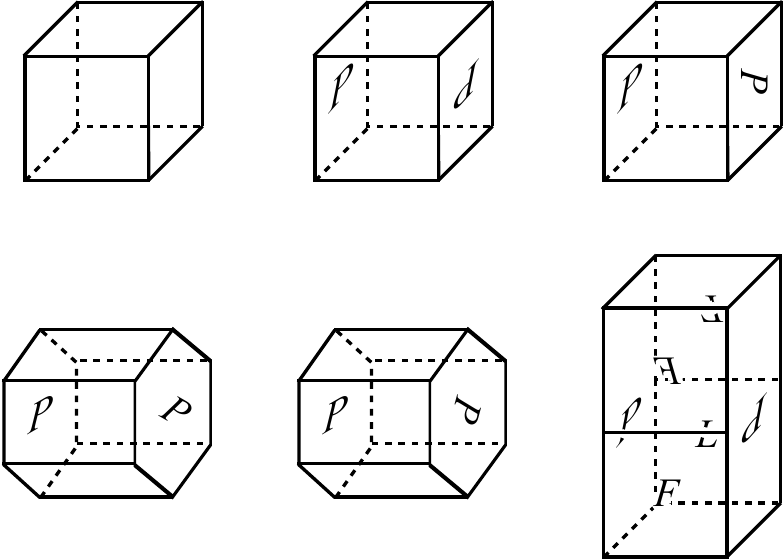}
\nota{The six closed orientable flat 3-manifolds, up to diffeomorphism. Each is constructed by pairing isometrically the faces of a polyhedron in $\matR^3$ according to the labels. When a face has no label, it is simply paired to its opposite by a translation. The polyhedra shown here are three cubes, two prisms with regular hexagonal basis, and one parallelepiped made of two cubes.}
\label{flat:fig}
\end{center}
\end{figure}

\begin{prop}
Every closed Seifert manifold $M$ with $e=\chi = 0$ admits a flat metric.
\end{prop}
\begin{proof}
There are six such Seifert manifolds up to diffeomorphism, listed in Table \ref{seven:table}. We build a flat metric for each in Figure \ref{flat:fig}.

The figure shows six flat manifolds, constructed by identifying isometrically the faces of a polyhedron in $\matR^3$. The reader is invited to check that each construction gives indeed a flat manifold, by verifying that the flat metric extends to the edges and to the vertices.

In all cases the foliation by parallel horizontal lines (orthogonal to the $P$ faces) descends to a Seifert fibration on the flat manifold. By looking at these lines one checks that the base surface of the fibration is respectively
$$T,\ (S^2,2,2,2,2),\ (S^2,2,4,4),\ (S^2,2,3,6),\ (S^2,3,3,3),\ (\matRP^2,2,2).$$
These orbifolds are obtained respectively from the figures by considering: the square torus, its quotient via a $\pi$-rotation, via a $\frac{\pi}4$-rotation, the quotient of a hexagon torus by a $\frac \pi 3$-rotation, by a $\frac{2\pi}3$-rotation, and the quotient of a Klein bottle via a $\pi$-rotation.

These flat Seifert manifolds have $e=0$ by Proposition \ref{flat:to:Seifert:prop}, and hence they are precisely those listed in Table \ref{seven:table}.
\end{proof}

We have proved Theorem \ref{main:flat:teo}.

\begin{figure}
\begin{center}
\includegraphics[width = 2.5 cm] {\iftoggle{BW}{HW-BW}{HW}}
\nota{The Hantsche-Wendt manifold $M=\matR^3/_\Gamma$. The group $\Gamma$ is generated by three rototranslations along the three \iftoggle{BW}{grey}{red} axis shown, each of angle $\pi$ and with unit displacement (this is a unit cube). A fundamental domain for $\Gamma$ is made of two unit cubes as in Figure \ref{flat:fig}-(bottom, right).}
\label{HW:fig}
\end{center}
\end{figure}

\begin{rem}
The rotational image $r(\Gamma) < SO(3)$ of $\Gamma$ for the six flat manifolds $M=\matR^3/_\Gamma$ constructed in Figure \ref{flat:fig} is respectively $\{e\}$, $C_2$, $C_4$, $C_6$, $C_3$, and the dihedral $D_4 = C_2\times C_2$. The sixth manifold is called the \emph{Hantzsche-Wendt manifold}: the group $\Gamma$ is generated by three rototranslations as in Figure \ref{HW:fig}.\index{three-manifold!Hantzsche-Wendt manifold}
\end{rem}

In contrast with the hyperbolic case, a finite-volume complete flat manifold is necessarily closed (see Proposition \ref{finite:volume:flat:prop}). The Seifert manifolds with boundary and $\chi=0$ are diffeomorphic to the bundles $T\times I$ and $K\timtil I$. Their interiors $T\times \matR$ and $K\timtil \matR$ may be given an \emph{infinite}-volume complete flat structure. The lack of a finite-volume complete flat structure for $K\timtil \matR$ is a reason for preferring the geometric decomposition to the canonical torus decomposition, see Section \ref{geometric:decomposition:subsection}.

\section{The product geometries}
The eight three-dimensional geometries include the three isotropic ones, plus five more. Among the five non-isotropic geometries, two are products of lower-dimensional geometries. We analyse them here.

\subsection{$S^2 \times \matR$ geometry.}
We equip $S^2\times \matR$ with the product metric. The product $S^2\times \matR$ is the poorest of the eight geometries, in the sense that there are very few manifolds modelled on $S^2\times \matR$. 

Recall that the sectional curvature is a number assigned to every plane in the tangent space of every point $p\in S^2 \times \matR$. One such plane is \emph{horizontal} if it is tangent to the $S^2$ factor and \emph{vertical} if it contains the line tangent to the $\matR$ factor.

\begin{prop}
The sectional curvatures of horizontal and vertical planes are $1$ and $0$, respectively.
\end{prop}
\begin{proof}
Let $\gamma \subset S^2$ be a closed geodesic. The surfaces $S^2 \times y$ and $\gamma\times \matR$ are totally geodesic, because they are fixed by some isometric reflections of $S^2\times \matR$. Therefore the sectional curvatures of horizontal and vertical planes equal the gaussian curvatures of these surfaces, which are $1$ and $0$.
\end{proof}

\begin{prop}
We have
$$\Iso(S^2\times \matR) = \Iso(S^2) \times \Iso(\matR).$$
\end{prop}
\begin{proof}
We certainly have the inclusion $\supset$, which gives to every point $x\in S^2\times \matR$ a stabiliser in $\Iso^+(S^2\times \matR)$ isomorphic to $\SO(2) \rtimes C_2$, a proper maximal subgroup of $\SO(3)$ by Proposition \ref{subgroups:SO3:prop}. 

If there were more isometries that that, there would be more fixing $p$ since they act transitively on $S^2\times \matR$ and the stabiliser would be the whole of $\SO(3)$, a contradiction because the sectional curvature of $S^2\times \matR$ is not constant.
\end{proof}
Since the isometry groups of $S^2$ and $\matR$ have two connected components each, the group $\Iso(S^2\times \matR)$ has four connected components, two of which are orientation-preserving.

\begin{prop}
An orientable manifold $M$ admits a finite-volume $S^2\times \matR$ geometry $\Longleftrightarrow$ $M$ is a closed Seifert manifold with $e= 0$ and $\chi>0$. 
\end{prop}
\begin{proof}
The closed Seifert manifolds with $e=0$ and $\chi>0$ are
just $S^2\times S^1$ and $\matRP^2 \timtil S^1$, and they are diffeomorphic to $(S^2\times \matR)/_{\Gamma}$ where $\Gamma$ is generated respectively by 
$$\big\{(\id, \tau)\big\}, \qquad \big\{(\iota, r), (\iota, r')\big\}$$
where $\tau$ is any translation, $\iota$ is the antipodal map and $r,r'$ are reflections with respect to distinct points in $\matR$.

Conversely, pick an orientable $M = (S^2\times S^1)/_\Gamma$. The discrete subgroup $\Gamma < \Iso(S^2) \times \Iso(\matR)$ preserves the foliation into spheres $S^2\times \{x\}$ which descends into a foliation into spheres and/or projective planes for $M$. Therefore $M$ decomposes into orientable interval bundles $S^2\times I$ and $\matRP^2\timtil I$, and is hence either $S^2\times S^1$ or $\matRP^2 \timtil S^1$.
\end{proof}

\subsection{$\matH^2 \times \matR$ geometry}
We give $\matH^2 \times \matR$ the product metric. The discussion of the previous section applies as is to this case, showing that horizontal and vertical planes in the tangent spaces have sectional curvature $-1$ and $0$. This in turn implies that
$$\Iso(\matH^2\times \matR) = \Iso(\matH^2) \times \Iso(\matR)$$
has four connected components, two being orientation-preserving. It is convenient to write the exact sequence
$$0 \longrightarrow \Iso(\matR) \longrightarrow \Iso(\matH^2\times \matR) \stackrel p \longrightarrow \Iso(\matH^2) \longrightarrow 0.$$
A discrete group $\Gamma<\Iso(X)$ is \emph{cofinite} if $X/_\Gamma$ has finite volume.

\begin{prop} \label{H2R:prop}
A discrete group $\Gamma < \Iso(\matH^2 \times \matR)$ is cofinite if and only if both $p(\Gamma)$ and $\Gamma \cap \ker p$ are discrete and cofinite.
\end{prop}
\begin{proof} 
If $p(\Gamma)$ is discrete we get
\begin{equation} \label{volume:product:eqn}
\Vol \big((\matH^2\times \matR)/_\Gamma \big) = \Area(\matH^2/_{p(\Gamma)}) \times {\rm Length}\big(\matR/ _{\Gamma \cap \ker p}\big). 
\end{equation}
This surprisingly simple formula is proved by picking a fundamental domain $D\subset \matH^2$ for $p(\Gamma)$ and noticing that 
$$\Vol \big((\matH^2\times \matR)/_\Gamma \big) = \Vol\big(p^{-1}(D)/_{\Gamma \cap \ker p}\big).$$ 
We deduce that $\Gamma$ is cofinite if and only if both $p(\Gamma)$ and $\Gamma \cap \ker p$ are. 

If $p(\Gamma)$ is not discrete, we prove that $\Gamma$ cannot be cofinite. Up to replacing $\Gamma$ with an index-four subgroup we may suppose 
$$\Gamma < \Iso^+(\matH^2) \times \Iso^+(\matR) = \Iso^+(\matH^2) \times \matR.$$
Pick a neighbourhood $U\subset \Iso^+(\matH^2) \times \matR$ of $e$ such that $[U,U] \subset U$ and $U\cap\Gamma = \{e\}$. Let $f,g \in \Gamma$ be two elements. 
We now prove that if $p(f), p(g) \in p(U)$ then $f$ and $g$ commute.
We note that the commutator $[f,g]$ depends only on the images $p(f)$ and $p(g)$, and since they lie in $p(U)$ we may suppose that $f,g\in U$ and get $[f,g]\in U\cap \Gamma$, which must be trivial.

The elements in $p(\Gamma)\cap p(U)$ commute. Two non-trivial elements in $\Iso^+(\matH^2)$ commute if and only if they are both hyperbolic, parabolic, or elliptic fixing the same line, point in $\partial \matH^2$, or point in $\matH^2$. Therefore all the isometries in $p(\Gamma)\cap p(U)$ are of the same type and fix the same line or point.

Analogously, for every $f \in \Gamma$ we pick a neighbourhood $U_f$ of $e$ such that $[f, U_f] \subset U$ and conclude that $f$ commutes with all the elements in $\Gamma $ projecting to $p(U_f)$. (There are non-trivial such elements since $p(\Gamma)$ is not discrete.) Therefore $p(f)$ also fixes the same line or point as above.

We have proved that $p(\Gamma)$ fixes a line, a horocycle, or a point in $\matH^2$ and hence fixes its inverse image in $\matH^2 \times \matR$ which is a line or a Euclidean plane. Moreover $\Gamma$ acts freely and proper discontinuously on it: hence $\Gamma =\matZ$ or $\matZ^2$ up to finite index and it is easy to deduce that $\Gamma$ is not cofinite.
\end{proof}

\begin{cor} \label{H2R:cor}
If the interior of a compact orientable manifold $M$ admits a finite-volume complete $\matH^2\times \matR$ geometry then $M$ is a Seifert manifold with $\chi<0$. If $M$ is closed then also $e = 0$.
\end{cor}
\begin{proof}
We have $\interior M = (\matH^2\times \matR)/_\Gamma$ with $\Gamma$ cofinite: by Proposition \ref{H2R:prop} the group $\Gamma \cap \ker p$ quotients every line $\{x\}\times \matR$ to a circle in $M$, giving a Seifert fibration $M\to S$ onto the finite-area orbifold $S=\matH^2/_{p(\Gamma)}$. We have $\chi(S)<0$, and either $e(M)=0$ or $\partial M \neq \emptyset$ because $\matH^2 \times y$ projects to a section for $M\to S$.
\end{proof}

We now prove the converse of Corollary \ref{H2R:cor}.
\begin{prop}
If $M$ is a Seifert manifold with  $\chi<0$ and either $\partial M \neq \emptyset$ or $e = 0$, the interior of $M$ admits a finite-volume complete $\matH^2\times S^1$ geometry.
\end{prop}
\begin{proof}
By hypothesis there is a section $\Sigma$ of $M\to S$, which is the fibre of a bundle $M\to O$ over a 1-orbifold $O$, see Section \ref{Seifert:fibers:subsection}. The two structures give two exact sequences
\begin{align*}
0 \longrightarrow K \longrightarrow \pi_1(M) & \stackrel f\longrightarrow \pi_1(S) \longrightarrow 0, \\
0 \longrightarrow \pi_1(\Sigma) \longrightarrow \pi_1(M) & \stackrel g\longrightarrow \pi_1(O) \longrightarrow 0.
\end{align*}
Since $\chi(S)<0$ we may write $S= \matH^2/_\Gamma$ and identify $\pi_1(S)$ with $\Gamma < \Iso(\matH^2)$. Analogously we consider $\pi_1(O)$ inside $\Iso(\matR)$.
The map 
$$(f,g) \colon \pi_1(M) \longrightarrow \Iso(\matH^2) \times \Iso(\matR)$$
is injective and its image is discrete and acts freely on $\matH^2\times \matR$, inducing a finite-volume $\matH^2 \times \matR$ structure on $M$.
\end{proof}

\section{Nil geometry}
There are still three non-isotropic geometries to analyse. These geometries are not products, but they have a reasonable bundle structure, so that many of the arguments of the previous section can be extended with not much effort.

We start with the geometry $\Nil$, which is a $\matR$-bundle over $\matR^2$. This geometry is fully carried by a Lie group called the \emph{Heisenberg group}.\index{group!Heisenberg group}\index{$\Nil$\ geometry}

\subsection{The Heisenberg group} 
The \emph{Heisenberg group} consists of all matrices
$$\begin{pmatrix} 1 & x & z \\ 0 & 1 & y \\ 0 & 0 & 1 \end{pmatrix}$$
with $x, y, z \in \matR$, with the multiplication operation. This is a nilpotent (but non-abelian) Lie group (see Exercise \ref{Heisenberg:ex}) and is hence also called Nil. It is clearly diffeomorphic to $\matR^3$.
There is a Lie groups exact sequence
$$ 0 \longrightarrow \matR \longrightarrow \Nil \longrightarrow \matR^2 \longrightarrow 0$$
where $\matR = [\Nil, \Nil] $ is the centre of $\Nil$ and consists of all matrices with $x=y=0$. Therefore Nil is naturally a line bundle over $\matR^2$.
We identify Nil with $\matR^3$ using the coordinates $(x,y,z)$. The product operation becomes
$$(x,y,z) \cdot (x',y',z') = (x+x', y+y', z+z'+xy').$$
\begin{ex} \label{Nil:commutator:ex}
The commutator of two elements is
$$\big[(x,y,z), (x',y',z')\big] = (0,0,xy'-x'y).$$
\end{ex}

\subsection{The geometry of Nil}
Every left-invariant Riemannian structure on a Lie group $G$ is constructed by fixing a scalar product at the tangent space $T_e G$ of $e\in G$ and then extending it by left-multiplication. Here $e=(0,0,0)$ and we give $T_e\Nil = \matR^3$ the Euclidean scalar product. This defines a left-invariant Riemannian structure on $\Nil$.

\begin{ex} \label{metric:tensor:Nil:ex}
The metric tensor at $(x,y,z)$ is
$$
\begin{pmatrix}
1 & 0 & 0 \\
0 & x^2 + 1 & -x \\
0 & -x & 1
\end{pmatrix}.
$$
\end{ex}
The following triple is an orthonormal basis at $(x,y,z)$:
\begin{equation} \label{orthonormal:eqn}
(1,0,0), \quad (0,1,x), \quad (0,0,1).
\end{equation}

\begin{figure}
\begin{center}
\includegraphics[width = 10 cm] {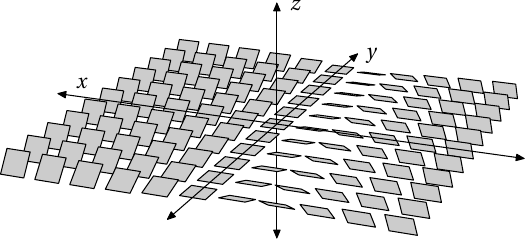}
\nota{The planes orthogonal to the $z$ axis in $\Nil$ form the standard contact structure on $\matR^3$.}
\label{contact:fig}
\end{center}
\end{figure}

The planes spanned by the first two vectors form the standard \emph{contact structure} in $\matR^3$ shown in Figure \ref{contact:fig}. It is a non-integrable distribution of planes: no surface is tangent to it at every point.\index{contact structure}

We remark that the metric tensor has unit determinant at every point: therefore the volume form on $\Nil$ is the standard one on $\matR^3$.

We can calculate the Christoffel symbols $\Gamma_{ij}^k$ by hand or using a computer code; these are all zero except the following:
$$\Gamma^1_{22} = -x, \quad \Gamma^1_{23} = \Gamma^1_{32} = \frac 12,$$  
$$\Gamma^2_{12} = \Gamma^2_{21} = \frac x 2, \quad \Gamma^2_{13} = \Gamma^2_{31} = -\frac 12, $$
$$\Gamma^3_{12} = \Gamma^3_{21} = \frac{x^2-1}2, \quad \Gamma^3_{13} = \Gamma^3_{31} = -\frac x2.$$
The Ricci tensor at $(x,y,z)$ is 
$$R_{ij} = \frac 12 \cdot \begin{pmatrix}
-1 & 0 & 0 \\
0 & {x^2-1} & -x \\
0 & -x & 1
\end{pmatrix}.
$$
When $x=0$ the Ricci tensor is a diagonal matrix with entries $-\frac 12, -\frac 12, \frac 12$. For a unit vector $v\in T_p\Nil$, recall that $R_{ij}v^iv^j$ equals twice the average value of the sectional curvatures of the planes containing $v$: this average value ranges here from $-\frac 14$ to $\frac 14$ and is maximal when $v=(0,0,\pm 1)$. This holds when $x=0$ and hence at all points $p\in \Nil$ by left-multiplication.

\subsection{The isometry group of Nil}
The group Nil acts on itself isometrically by left-multiplication. Left-multiplication by $(a,b,c)$ preserves the bundle $\Nil \to \matR^2$ and induces on $\matR^2$ a translation by the vector $(a,b)$. 

There are also more complicate isometries of $\Nil$ that preserve the bundle structure but induce rotations on $\matR^2$: one such isometry $\varphi$ sends $(x,y,z)$ to 
$$\Big(x\cos\theta -y\sin \theta, \  x\sin \theta + y \cos \theta,\ z+ \frac 12 (x^2-y^2)\sin \theta \cos \theta - xy\sin^2 \theta \Big).$$
This map preserves the bundle and induces a rotation on $\matR^2$.
\begin{prop} \label{varphi:isometry:prop}
The map $\varphi$ is an isometry of $\Nil$.
\end{prop}
\begin{proof}
Set $(x',y',z') = \varphi (x,y,z)$. The differential $d\varphi$ at $(x,y,z)$ is
$$
\begin{pmatrix} 
\cos \theta & -\sin \theta & 0 \\
\sin \theta & \cos \theta & 0 \\
x\sin\theta \cos\theta - y\sin^2 \theta & - y\sin\theta\cos \theta  - x\sin^2\theta & 1 
\end{pmatrix} 
$$
and may be rewritten as
$$
\begin{pmatrix} 
\cos \theta & -\sin \theta & 0 \\
\sin \theta & \cos \theta & 0 \\
x'\sin \theta & x'\cos \theta - x & 1.
\end{pmatrix} 
$$
The differential $d\varphi$ sends the orthonormal basis (\ref{orthonormal:eqn}) to
$$(\cos\theta, \sin \theta, x'\sin\theta), \quad (-\sin\theta, \cos\theta, x'\cos\theta), \quad (0,0,1).$$
These vectors at $(x',y',z')$ are also orthonormal.
\end{proof}

We deduce the following.

\begin{prop} \label{Nil:isometries:prop}
Every isometry of $\Nil$ preserves the line bundle and induces an isometry of $\matR^2$. We have
$$0 \longrightarrow \matR \longrightarrow \Iso^+(\Nil) \stackrel p \longrightarrow \Iso(\matR^2) \longrightarrow 0$$
where $\matR$ is the centre of $\Nil$.
The group $\Iso^+(\Nil)$ has two components.
\end{prop}
\begin{proof}
The rotational isometries $\varphi$ with angle $\theta$ introduced above form a $S^1$-subgroup of $\Iso^+(\Nil)$. The subgroups $S^1$ and $\Nil$ belong to the component $\Iso^+_\circ(\Nil) < \Iso^+(\Nil)$ containing the identity, which hence has dimension at least $1+3=4$ and acts transitively on $\Nil$. Stabilisers have dimension at least $4-3=1$ and cannot have bigger dimension by Proposition \ref{subgroups:SO3:prop}, otherwise $\Nil$ would have constant sectional curvature: therefore $\dim(\Iso^+_\circ(\Nil))=4$ and $\Iso^+_\circ(\Nil)$ is generated by $S^1$ and $\Nil$, and it preserves the line bundle since $S^1$ and $\Nil$ do.

The subgroup $\Iso^+_\circ(\Nil)$ preserves the orientation of the fibres and of $\matR^2$.
We leave as an exercise the existence of another component of $\Iso^+(\Nil)$ which inverts the orientation of the fibres and of $\matR^2$.
\end{proof}

\begin{figure}
\begin{center}
\includegraphics[width = 6 cm] {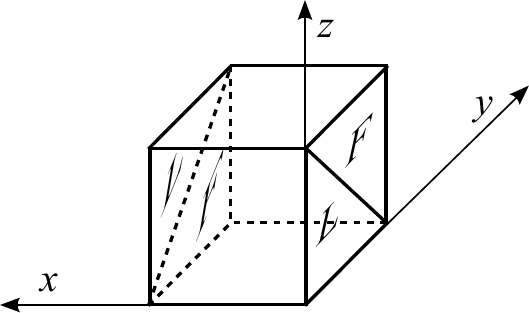}
\nota{The $\Nil$\ manifold $M=\Nil/_\Gamma$ is obtained by identifying the faces of this unit cube as follows: the triangular faces are glued via the affine maps as shown by the labels, and the pairs of unlabelled opposite square faces are identified by translations. The cube is a fundamental domain for $\Gamma$.}
\label{Nil_example:fig}
\end{center}
\end{figure}

\begin{example}  \label{Nil:example}
A manifold $M = \Nil/_\Gamma$ modelled on Nil is described in Figure \ref{Nil_example:fig}. The figure shows a fundamental domain for the group $\Gamma<\Iso^+(\Nil)$ generated by the isometries
$$(x,y,z) \longmapsto (x+1,y,z+y),$$
$$(x,y,z) \longmapsto (x, y+1, z), \quad (x,y,z) \longmapsto (x,y,z+1).$$
obtained by left-multiplication with the canonical basis of $\matR^3$.
Both the unit cube and $M$ have volume one. The manifold $M$ is clearly a torus bundle with monodromy $\matr 1101$ and is hence diffeomorphic to the Seifert manifold $\big(T,(1,1)\big)$ by Exercise \ref{Me:ex}.
\end{example}

A group $\Gamma< \Iso(X)$ is \emph{cocompact} if $X/_\Gamma$ is compact.
We prove that cocompact groups in $\Iso(\matR^2)$ do not lift.

\begin{prop} \label{Nil:no:lift:prop}
Let $\Gamma < \Iso(\matR^2)$ be discrete and cocompact. There is no homomorphism $f\colon \Gamma \to \Iso^+(\Nil)$ such that $p\circ f = \id$.
\end{prop}
\begin{proof}
Up to taking a finite-index subgroup we may suppose that $\Gamma$ is generated by two translations along independent vectors $(x,y)$ and $(x',y')$. We prove that two lifts $\varphi, \varphi' \in \Iso^+(\Nil)$ of these translations never commute, thus forbidding the existence of a homomorphism $f$.

If $\varphi$, $\varphi' \in \Nil$, Exercise \ref{Nil:commutator:ex} gives $[\varphi, \varphi'] = xy'-x'y \neq 0$. All the other lifts are of type $g\varphi, g'\varphi'$ for some $g,g'\in \matR$, hence we get the same commutator. 
\end{proof}

Let $\Iso_0^+(\Nil)$ be the connected component of $\Iso^+(\Nil)$ containing $e$.

\begin{ex} \label{isonil:nilpotente:ex}
We have 
$$\big[\Iso_0^+(\Nil), \Iso_0^+(\Nil)\big] = \matR.$$ 
Therefore $\Iso_0^+(\Nil)$ is nilpotent.
\end{ex}

\begin{cor} \label{Nil:is:Nil:cor}
If a three-manifold $M$ has a $\Nil$ geometry, then $\pi_1(M)$ is virtually nilpotent.
\end{cor}
\begin{proof}
It has a nilpotent subgroup of index at most two, because
the subgroup $\Iso_0^+(\Nil)< \Iso^+(\Nil)$ is nilpotent and has index two.
\end{proof}

\subsection{Nil geometry}
We now classify the manifolds modelled on $\Nil$. 

\begin{prop} \label{Nil:cofinite:prop}
A discrete group $\Gamma < \Iso^+ (\Nil)$ is cofinite if and only if both $p(\Gamma)$ and $\Gamma \cap \ker p$ are discrete and cofinite.
\end{prop}
\begin{proof}
We follow the proof of Proposition \ref{H2R:prop}. Up to taking finite-index subgroups we may suppose that $\Gamma<\Iso^+_0(\Nil)$.
If $p(\Gamma)$ is discrete the proof of Proposition \ref{H2R:prop} applies also here; note that the formula (\ref{volume:product:eqn}) holds because the volume form in $\Nil$ is the standard one on $\matR^3$.

If $p(\Gamma)$ is not discrete then that proof shows that there is a neighbourhood $U\subset \Iso^+(\Nil)$ of $e$ such that every two elements $f,g\in\Gamma$ projecting in $p(U)$ commute. (It is still true that $[f,g]$ depends only on $p(f)$ and $p(g)$ since $\ker p = \matR$ is central.) We deduce again that the isometries in $p(\Gamma)$ commute: hence they are either rotations fixing the same point or translations. 

In the former case $\Gamma$ is not cofinite. In the latter, we get $\Gamma < \Nil$ acting as left-multiplication. Every two $f,g \in \Gamma$ projecting to $p(U)$ commute and hence project to parallel translations by Exercise \ref{Nil:commutator:ex}. As in the proof of Proposition \ref{H2R:prop} we deduce that $p(\Gamma)$ preserves a line, and hence it is not cofinite.
\end{proof}

\begin{cor} If the interior of a compact orientable manifold $M$ admits a finite-volume complete $\Nil$ geometry then $M$ is a closed Seifert manifold with $\chi=0$ and $e \neq 0$.
\end{cor}
\begin{proof}
We have $\interior M = \Nil/_\Gamma$ with $\Gamma$ cofinite and Proposition \ref{Nil:cofinite:prop} provides a Seifert fibration $M\to S$ over a finite-area orbifold $S=\matR^2/_{p(\Gamma)}$. Finite-area flat orbifolds have $\chi(S)=0$ and are closed (there are no cusps in flat geometry), hence $M$ is closed. 

We have $e\neq 0$, otherwise up to finite-index we would get $M = T\times S^1$ contradicting Proposition \ref{Nil:no:lift:prop}.
\end{proof}

We prove the converse.

\begin{prop} \label{Nil:construct:prop}
If $M$ is a closed Seifert manifold with $\chi=0$ and $e\neq 0$, then $M$ admits a $\Nil$ geometry.
\end{prop}
\begin{proof}
We have
$$ M = \big(S, (p_1,q_1), \ldots, (p_k, q_k) \big).$$
Suppose $S$ is a closed orientable surface of genus $g\geqslant 0$. In the following $\pi_1(S)$ is the orbifold fundamental group. We have:
\begin{align*}
\pi_1(S) & = \langle a_1, b_1, \ldots, a_g, b_g, c_1, \ldots, c_k \ |\ [a_1,b_1]\cdots [a_g,b_g] c_1 \cdots c_k, c_i^{p_i} \rangle, \\
\pi_1(M) & = \langle a_1, b_1, \ldots, a_g, b_g, c_1, \ldots, c_k, l \ |\ [a_1,b_1]\cdots [a_g,b_g] c_1 \cdots c_k, c_i^{p_i}l^{q_i}, \\ 
& \qquad [a_i,l], [b_i,l] \rangle. 
\end{align*}
We fix any flat structure on the orbifold $S$ and get an injection $\pi_1(S) \to\Iso^+(\matR^2)$. We now want to lift this map to an injection $\pi_1(M) \to \Iso^+(\Nil)$, so that the resulting diagram commutes:
$$
\xymatrix{ 
\pi_1(M) \ar@{.>}[r] \ar[d] & \Iso^+(\Nil) \ar[d] \\
\pi_1(S)\ar[r] & \Iso^+(\matR^2)
}
$$
Recall that $\Iso_0^+(\Nil)$ contains the group $\matR$ of vertical translations: we use the multiplicative notation and indicate it as $\matR_{>0}$. We note that $\matR_{>0}$ is central in $\Iso_0^+(\Nil)$ because it commutes with $\Nil$ and the isometries $\varphi$ from Proposition \ref{varphi:isometry:prop}, which altogether generate $\Iso_0^+(\Nil)$.

We identify $\pi_1(S)$ with its image in $\Iso^+(\matR^2)$.
We lift arbitrarily $a_i,b_i,c_i$ inside $\Iso^+(\Nil)$ and pick $l\in \matR_{>0}$. We get
\begin{align*}
[a_1,b_1]\cdots [a_g,b_g] c_1 \cdots c_k & = e^\mu, \\
c_i^{p_i}l^{q_i} & = e^{\lambda_i}
\end{align*}
for some $\mu, \lambda_i \in \matR$. To get a homomorphism $\pi_1(M) \to \Iso^+(\Nil)$ we need $\mu = \lambda_i = 0$. To obtain that we change the lifts as $c_i' = e^{t_i} c_i$, $l' = e^u l$ to get
\begin{align*}
[a_1,b_1]\cdots [a_g,b_g] c'_1 \cdots c'_k & = e^{t_1 + \ldots + t_k + \mu}, \\
(c'_i)^{p_i}(l')^{q_i} & = e^{t_ip_i+ uq_i + \lambda_i}.
\end{align*}
We want
\begin{align*}
t_1+\ldots + t_k &= - \mu, \\
t_ip_i + uq_i &= - \lambda_i.
\end{align*}
The determinant of the $k\times k$ coefficient matrix was already calculated in the proof of Proposition \ref{Seifert:homology:prop}, and is
$$\pm\sum \frac{q_i}{p_i} (p_1\cdots p_k)= \pm e \cdot p_1\cdots p_k \neq 0.$$
Therefore the linear system has a unique solution $(t_1,\ldots, t_k, u)$. The lift $l'$ is necessarily non-trivial, otherwise we would get a lift $\pi_1(S) \to \Iso^+(\Nil)$ that is excluded by Proposition \ref{Nil:no:lift:prop}. Therefore the lift is injective.

If $S$ is non-orientable the proof is similar and left as an exercise.
\end{proof}

\section{$\widetilde{\SL_2}$ geometry}
The $\widetilde{\SL_2}$ geometry is similar to $\Nil$, but it is now a $\matR$-bundle over $\matH^2$. The geometry is again fully carried by the Lie group $\widetilde{\SL_2}$, which is the universal cover of both $\SL_2(\matR)$ and $\PSLR$. It is convenient to identify $\PSLR$ with the unit tangent bundle of $\matH^2$.\index{$\widetilde{\SL_2}$ geometry}

\subsection{The unit tangent bundle}
The tangent bundle $TM$ of a Riemannian $n$-manifold $M$ has a natural Riemannian structure, which we briefly introduce.

The tangent space at any point $x\in TM$ splits into a \emph{vertical} subspace $V_x  = \ker dp_x$ where $p\colon TM \to M$ is the projection, and a \emph{horizontal} subspace $H_x $ determined by the metric tensor on $M$ as follows: use the parallel transport to move the tangent vector $x$ along all geodesics exiting from $p(x)$; the result is a small $n$-surface in $TM$ containing $x$ and we set $H_x$ to be its tangent space at $x$. 

To define a metric on $TM$ we impose that $V_x$ and $H_x$ be orthogonal and we give to each space the metric of $T_xM$, via the natural identification $V_x = T_xM$ and via the isomorphism $dp_x\colon H_x \to T_xM$. 

The \emph{unit tangent bundle} $UM\subset TM$ consists of all unitary tangent vectors and inherits a Riemannian structure. Every isometry $f\colon M \to M$ induces an isometry $df\colon UM \to UM$.

\subsection{The space $U\matH^2$}
We now focus on the case $M=\matH^2$ we are interested in. Parallel transport was defined explicitly in Section \ref{parallel:subsection}. 
We represent $\matH^2$ using the upper half-plane model $H^2 = \{\Im z > 0\}$, so that $TH^2 = H^2 \times \matC$ and $UH^2 = H^2 \times S^1$ has coordinates $(z,\theta)$. The tangent space $T_{(z,\theta)}UH^2$ is naturally identified with $\matC\times \matR = \matR^3$.

\begin{lemma} \label{UH2:tensor:lemma}
The metric tensor of $UH^2$ at $(z,\theta)$ is
$$
\begin{pmatrix}
2y^{-2} & 0 & y^{-1} \\
0 & y^{-2} & 0 \\
y^{-1} & 0 & 1 
\end{pmatrix}
$$
where $z=x+iy$.
\end{lemma}
\begin{proof}
We first consider the disc model $D^2$ of $\matH^2$ with $TD^2 = D^2 \times \matC$ and $T_{(z,v)} TD^2 = \matC \times \matC$.  
We focus at a point $(0,v)$ and determine the decomposition 
$$T_{(0,v)} = H_{(0,v)} \oplus V_{(0,v)}.$$ 
We have $V_{(0,v)} = 0 \times \matC$ and we now determine $H_{(0,v)}$. Every geodesic through $0$ is a Euclidean line $l$ and the parallel transport of $v\in T_0D^2$ along $l$ forms a constant angle with $l$: the parallel transport of $v$ at $z\in l$ is just the rescaled vector $(1-|z|^2)v$. Parallel transports of $v$ along lines passing through $0$ form a surface, which is the graph of the function $f(z) = (1-|z|^2)v$, whose tangent plane at $0$ is $H_{(0,v)} = \matC \times 0$ since $\frac{\partial f}{\partial x} = \frac{\partial f}{\partial y} = 0$ at $0$.

We have discovered that the horizontal and vertical planes at $(0,v)$ are just the coordinate ones of $T_{(0,v)} = \matC \times \matC$. By definition these are orthogonal and inherit a metric tensor from that of $T_0D^2$, which is $4$ times the Euclidean one. Therefore the metric tensor at $T_{(0,v)} = \matC^2$ is $4$ times the Euclidean one. 

We now turn to the unitary sub-bundle $UD^2 = D^2 \times S^1$. We use the natural identifications $T_{(0,\theta)} UD^2 = \matC \times \matR = \matR^3$ and we deduce that the vectors $(1,0,0)$, $(0,1,0)$, $(0,0,1) \in T_{(0,\theta)}$ are orthogonal and have norm $2,2,1$ respectively. 

The M\"obius transformation
$$f(z) = \frac{z - i}{i z - 1}$$
is an isometry between the two models $D^2$ and $H^2$ that sends $0$ to $i$. 
It is holomorphic and its complex derivative is
$$f'(z) = \frac{-2}{(iz-1)^2}.$$
The isometry $f\colon D^2\to H^2$ induces an isometry $f_*\colon UD^2 \to UH^2$, that is $f_*\colon D^2\times S^1 \to H^2 \times S^1$, which is as follows:
\begin{align*}
f_*(z,\theta) & = \left(\frac{z-i}{iz-1}, \theta + {\rm arg} (f'(z))\right) \\
& = \left(\frac{z-i}{iz-1}, \theta +\pi - 2 {\rm arg} (iz-1)\right)  \\
& = \left(\frac{z-i}{iz-1}, \theta +\pi - 2 \Im \log (iz-1)\right).  
\end{align*}
In particular $f_*(0,\theta) = (i, \theta + \pi)$.
Recall that the Jacobians of a holomorphic $g$ and of $\log g$ are
$$
Jg = \begin{pmatrix} \Re g' & -\Im g' \\ \Im g' & \Re g' \end{pmatrix},
\quad
J(\log g) = \begin{pmatrix} \Re \frac{g'}g & -\Im \frac{g'}g \\ \Im \frac{g'}g & \Re \frac{g'}g \end{pmatrix}.
$$
The differential of $f_*$ at the point $(z,\theta)$ is hence
$$
(df_*)_{(z,\theta)} =
\begin{pmatrix}
\Re \frac{-2}{(iz-1)^2} & - \Im \frac{-2}{(iz-1)^2} & 0 \\
\Im \frac{-2}{(iz-1)^2} & \Re \frac{-2}{(iz-1)^2} & 0 \\
-2 \Im \frac{1}{z+i} & - 2 \Re \frac{1}{z+i} & 1 
\end{pmatrix}.
$$
In particular at $z=0$ we get
$$
(df_*)_{(0,\theta)} =
\begin{pmatrix}
-2 & 0 & 0 \\ 
0 & -2 & 0 \\
2 & 0 & 1 
\end{pmatrix}.
$$
Since this is an isometry, the three image vectors $(-2,0,2)$, $(0,-2,0)$, and $(0,0,1)$ are orthogonal with norm $2,2,1$, hence $(1,0,-1)$, $(0,1,0)$, and $(0,0,1)$ form an orthonormal basis at $T_{(i,\theta+\pi)}UH^2$ and the metric tensor there is
$$\begin{pmatrix}
2 & 0 & 1 \\
0 & 1 & 0 \\
1 & 0 & 1
\end{pmatrix}.
$$
The tensor is independent of $\theta$. To find its value at $(z,\theta)$ for a generic $z\in H^2$ it suffices to transport it via an isometry that sends $i$ to $z$. Every such isometry is a composition of a horizontal translation (which does not change the tensor) and a dilation $g(z) = \lambda z$ with $\lambda >0$. We have $g'(z) = \lambda$ and $g_*(z,\theta) = (\lambda z, \theta)$, therefore $dg_*$ sends the above orthonormal basis at $T_{(i,\theta)}$ to the basis $(\lambda, 0, -1)$, $(0,\lambda,0)$, $(0,0,1)$. The metric tensor at $(\lambda i, \theta)$ is hence
$$\begin{pmatrix}
2\lambda^{-2} & 0 & \lambda^{-1} \\
0 & \lambda^{-2} & 0 \\
\lambda^{-1} & 0 & 1
\end{pmatrix}.
$$
The proof is complete.
\end{proof}
At a point $(x,y,\theta) \in UH^2$, an orthonormal basis is
\begin{equation} \label{UH2:orthonormal:eqn}
(y,0,-1), \quad (0,y,0), \quad (0,0,1).
\end{equation}
The planes orthogonal to $(0,0,1)$ form a contact structure similar to the one drawn in Figure \ref{contact:fig}. 

We remark that the metric tensor at $(x,y,\theta)$ has determinant $y^{-4}$, like in the product metric $H^2 \times S^1$: therefore the volume form in $UH^2$ is the same as in the product metric (although the metric tensor is not).

The non-zero Christoffel symbols at $(x,y,\theta)$ are:
$$\Gamma^1_{12} = \Gamma^1_{21} = -\frac 3{2y}, \quad 
\Gamma^1_{23} = \Gamma^1_{32} = -\frac 12,$$
$$\Gamma^2_{11} = \frac 2y, \quad
\Gamma^2_{13} = \Gamma^2_{31} = \frac 12, \quad \Gamma^2_{22} = -\frac 1y,$$
$$\Gamma^3_{12} = \Gamma^3_{21} = \frac 1{y^2}, \quad
\Gamma^3_{23} = \Gamma^3_{32} = \frac 1{2y}.$$
The Ricci tensor is
$$
R_{ij} = \frac 12 \cdot
\begin{pmatrix}
-2 y^{-2} & 0 & y^{-1} \\
0 & - 3y^{-2} & 0 \\
y^{-1} & 0 & 1
\end{pmatrix}.
$$
When $y=1$ we can represent the Ricci tensor in the orthonormal basis (\ref{UH2:orthonormal:eqn}) and get a diagonal matrix with values $-\frac 32$, $-\frac 32$, $\frac 12$. For a vector $v\in T_{(z,\theta)} U\matH^2$, recall that $R_{ij}v_iv_j$ is twice the average value of the sectional curvatures of the planes containing $v$: this average value ranges from $-\frac 34$ to $\frac 14$ and is maximal when $v=(0,0,\pm 1)$. This holds when $y=1$ and hence for any $z\in U\matH^2$ by acting via isometries of $\matH^2$. 

\subsection{$\widetilde{\SL_2}$ geometry}
Let $\widetilde {\SL_2}$ be the universal cover of $\SL_2 = \SL_2(\matR)$. 
As a universal cover of a Lie group, it is also a Lie group.
We have coverings
$$\widetilde{\SL_2} \longrightarrow \SL_2 \longrightarrow \PSLR = \Iso^+(\matH^2).$$
The group $\PSLR$ acts freely and transitively on $U\matH^2$ and hence we can identify $\PSLR$ with $U\matH^2$. With this identification $\PSLR$ inherits a left-invariant Riemannian metric that lifts to a left-invariant Riemannian metric on the Lie group $\widetilde{\SL_2}$. 

In the previous section we have identified $U\matH^2$ with $H^2 \times S^1$ and we can likewise identify $\widetilde{\SL_2}$ with $H^2 \times \matR$, and get an explicit metric tensor from Lemma \ref{UH2:tensor:lemma}.

Since $U\matH^2$ is a circle bundle over $\matH^2$, likewise $\widetilde{\SL_2}$ is a line bundle over $\matH^2$. The group $\matR$ acts isometrically on $\widetilde \SL_2$ by translating every fibre.

\begin{prop}
Every isometry of $\widetilde{\SL_2}$ preserves the line bundle and induces an isometry of $\matH^2$. We have
$$0 \longrightarrow \matR \longrightarrow \Iso^+(\widetilde{\SL_2}) \stackrel p\longrightarrow \Iso(\matH^2) \longrightarrow 0.$$
The group $\Iso^+(\widetilde{\SL_2})$ has two components. 
\end{prop}
\begin{proof}
The groups $\matR$ and $\widetilde{\SL_2}$ belong to the component $\Iso^+_\circ (\widetilde{\SL_2})$ containing $e$ which has dimension at least $1+3=4$. We conclude as in the proof of Proposition \ref{Nil:isometries:prop}.
\end{proof}

\begin{ex}
The unit tangent bundle of a finite-area complete hyperbolic surface is naturally a manifold modelled on $\widetilde{\SL_2}$.
\end{ex}

As for $\Nil$, cocompact groups do no lift.

\begin{prop} \label{SL2:no:lift:prop}
Let $\Gamma < \Iso(\matH^2)$ be discrete and cocompact. There is no homomorphism $f\colon \Gamma \to \Iso^+(\widetilde{\SL_2})$ such that $p\circ f = \id$.
\end{prop}
\begin{proof}
Up to taking a finite index subgroup we may suppose that $\Gamma$ acts freely and hence $S=\matH^2/_\Gamma$ is a closed hyperbolic surface. If $\Gamma$ lifts, consider the group $G< \Iso^+(\widetilde{\SL_2})$ generated by $f(\Gamma)$ and $2 \pi\in \matR$, isomorphic to $\Gamma \times \matZ$. The quotient $\widetilde{\SL_2}/_G$ is naturally the unit tangent bundle of $S$ and is hence the Seifert manifold $\big(S, (1, \chi(S))\big)$ by Remark \ref{bundle:oss}. Its fundamental group is however not a product $\Gamma \times \matZ$ by Exercise \ref{bundle:ex}.
\end{proof}

\begin{prop} \label{SL2:cofinite:prop}
A discrete group $\Gamma < \Iso^+(\widetilde{\SL_2})$ is cofinite if and only if both $p(\Gamma)$ and $\Gamma \cap \ker p$ are discrete and cofinite.
\end{prop}
\begin{proof}
Same proof as Proposition \ref{H2R:prop}, with a minor variation: if $p(\Gamma)$ is not discrete, up to taking an index-two subgroup we suppose that $p(\Gamma) <  \Iso^+(\matH^2)$; the $\matR$-action commutes with $\Gamma$ and $[f,g]$ depends only on $p(f)$ and $p(g)$, so that proof applies.
\end{proof}

We now classify the manifolds having a $\widetilde{\SL_2}$ geometry.

\begin{cor}
If the interior of a compact orientable manifold $M$ admits a finite-volume complete $\widetilde{\SL_2}$ geometry then $M$ is a Seifert manifold with $\chi<0$. If $M$ is closed then $e\neq 0$.
\end{cor}
\begin{proof}
We have $\interior M= \widetilde{\SL_2}/_\Gamma$ with $\Gamma$ cofinite. Proposition \ref{SL2:cofinite:prop} furnishes a Seifert fibration $M\to S$ over the finite-area orbifold $S=\matH^2/_{p(\Gamma)}$. If $M$ is closed we get $e\neq 0$: if not, up to taking a finite index subgroup we would get $M=S\times S^1$ contradicting Proposition \ref{SL2:no:lift:prop}.
\end{proof}

We now prove the converse.

\begin{prop}
If $M$ is a Seifert manifold with $\chi<0$ and either $\partial M\neq \emptyset$ or $e\neq 0$, the interior of $M$ admits a finite-volume $\widetilde{\SL_2}$ geometry.
\end{prop}
\begin{proof}
If $M$ is closed we apply the proof of Proposition \ref{Nil:construct:prop}. If $\partial M \neq \emptyset$ the presentations of $\pi_1(S)$ and $\pi_1(M)$ are as described there, except that they do not contain the relator $c_i^{p_i}$ whenever $c_i$ represents a boundary component of $M$. So we have less constraints and we easily see that a solution to the final linear problem exists also in this case. 
\end{proof}

\section{Sol geometry}
The $\Sol$ geometry is the least symmetric one among the eight. It has a bundle structure, but with a one-dimensional basis: it is a $\matR^2$-bundle over $\matR$. Again, the geometry is fully governed by a Lie group $\Sol$.\index{$\Sol$ geometry}

\subsection{The Lie group Sol}
The Lie group $\Sol$ is the space $\matR^3$ equipped with the following operation
$$(x,y,z)\cdot(x',y',z') = (x+ e^{-z}x', y+e^zy', z+z').$$
\begin{ex} \label{Sol:conti:ex}
We have
\begin{align*}
(x,y,z)^{-1} & = (-xe^z, -ye^{-z}, -z), \\
[(x,y,z),(x',y',z')] & = \big(x(1\!-\!e^{-z'}) - x'(1\!-\!e^{-z}), y(1\!-\!e^{z'}) - y'(1\!-\!e^z), 0\big), \\
[(x,y,z), (x',y',0)] & = \big(- x'(1-e^{-z}), - y'(1 - e^z), 0\big). \
\end{align*}
\end{ex}
The subgroup $\matR^2$ consisting of all elements $(x,y,0)$ is the centre of $\Sol$ and by setting $p(x,y,z)=z$ we get an exact sequence
$$0 \longrightarrow \matR^2 \longrightarrow \Sol \stackrel p \longrightarrow \matR \longrightarrow 0.$$
Therefore $\Sol$ is a plane bundle over $\matR$. Exercise \ref{Sol:conti:ex} implies the following.
\begin{ex} We have $[\Sol, \Sol] = \matR^2$ and hence $\Sol$ is solvable. However $[\Sol, \matR^2] = \matR^2$ and hence $\Sol$ is not nilpotent.
\end{ex}

We define a Riemannian metric on $\Sol$ by assigning the scalar product
$$\begin{pmatrix}
e^{2z} & 0 & 0 \\
0 & e^{-2z} & 0 \\
0 & 0 & 1
\end{pmatrix}
$$
to the point $(x,y,z)$. The metric is left-invariant and every plane $z=k$ is isometric to the Euclidean $\matR^2$. This is the geometry with the smallest amount of symmetries.

We remark that the metric tensor has unit determinant at every point: therefore the volume form on $\Sol$ is the standard one on $\matR^3$.

The non-zero Christoffel symbols at $(x,y,z)$ are
$$\Gamma^1_{13} = \Gamma^1_{31} = 1,$$
$$\Gamma^2_{23} = \Gamma^2_{32} = -1,$$
$$\Gamma^3_{11} = -e^{2z}, \quad \Gamma^3_{22} = e^{-2z}.$$
The Ricci tensor is
$$\begin{pmatrix}
0 & 0 & 0 \\
0 & 0 & 0 \\
0 & 0 & -2
\end{pmatrix}.
$$
The average value of the sectional curvatures of the planes containing $v\in T_{(x,y,z)} \Sol$ ranges from $-1$ to $0$, and is minimal when $v=(0,0,\pm 1)$.

\subsection{Sol geometry}
We start with a simple exercise.
\begin{ex} \label{D8:ex}
The eight maps 
$$(x,y,z) \mapsto (\pm x, \pm y,z), \qquad (x,y,z) \mapsto (\pm y, \pm x, -z)$$
are isometries and form the dihedral group $D_8$. The orientation-preserving ones form the subgroup $D_4 = \matZ_2 \times \matZ_2$.
\end{ex}
Let $\Iso_*(\matR^2) < \Iso(\matR^2)$ be the subgroup consisting of all maps $v\mapsto \pm v +b$. It has two components, one being the translations $\matR^2$.
\begin{prop}
Every isometry of $\Sol$ preserves the plane bundle and induces an isometry of $\matR$. We have
$$0 \longrightarrow \Iso_*(\matR^2) \longrightarrow \Iso^+(\Sol) \stackrel p \longrightarrow \Iso(\matR) \longrightarrow 0.$$
The group $\Iso^+(\Sol)$ has four components, one of which is $\Sol$ acting by left-multiplication. 
\end{prop}
\begin{proof}
The group $\Sol$ acts transitively and freely on $\Sol$ itself, and to conclude it suffices to check that the stabiliser of the origin $0$ is the dihedral $D_4$ described in Exercise \ref{D8:ex}. 

The Ricci tensor tells us that an isometry fixing a point also fixes the vertical axis and the horizontal plane. Therefore the vertical unitary constant vector field $X=(0,0,1)$ is preserved up to sign by any isometry of $\Sol$. The covariant differentiation $v \mapsto \nabla_v X$ defines an endomorphism of $T_p \Sol$ for all $p \in \Sol$. We have $\nabla_{e_i} X = \Gamma_{i3}^ke_k$ and therefore the endomorphism is
$$\begin{pmatrix}
1 & 0 & 0 \\
0 & -1 & 0 \\
0 & 0 & 0 
\end{pmatrix}.
$$
The three coordinate axis are precisely the eigenvectors of the endomorphism, and being intrinsically defined they are preserved by every isometry. Therefore the orientation-preserving stabiliser of a point is $D_4$.
\end{proof}

\begin{prop}
A discrete group $\Gamma < \Iso^+(\Sol)$ is cofinite if and only if both $p(\Gamma)$ and $\Gamma \cap \ker p$ are discrete and cofinite.
\end{prop}
\begin{proof}
If $p(\Gamma)$ is discrete we get
$$\Vol \big(\Sol/_\Gamma \big) = {\rm Length}(\matR/_{p(\Gamma)}) \times \Area\big(\matR^2/ _{\Gamma \cap \ker p}\big). $$
This formula is proved as above by taking a fundamental domain for $p(\Gamma)$. We now prove that $p(\Gamma)$ is in fact discrete (since $\Gamma$ is).

Up to taking an index-four subgroup we suppose that $\Gamma < \Sol$. 
If $\gamma \in \Sol$ does not lie in $\matR^2$, then it fixes a vertical line. To prove that, note that $\gamma = (x,y,z)$ acts on $\matR^3$ as an affine transformation, which permutes the vertical lines and acts on $\matR^2$ as $(x',y') \mapsto (x+e^{-z}x', y+e^zy')$. If $z \neq 0$ this map has a fixed point (because $1$ is not an eigenvalue of its linearisation).

If $\Gamma$ is abelian, either it is contained in $\matR^2$ or it fixes a vertical line: in both cases we get a discrete $p(\Gamma)$. If $\Gamma$ is non-abelian, then $[\Gamma, \Gamma]<\matR^2$ is non-trivial so $\Gamma$ contains a non-trivial element $\gamma \in \matR^2$ and another $\eta \not \in \matR^2$. The elements $\gamma$ and $\eta\gamma\eta^{-1}$ are both in $\matR^2$ and independent, hence $\matR^2 /_{\Gamma \cap \ker p}$ is compact. This implies easily that $p(\Gamma)$ is discrete (since $\Gamma$ is).
\end{proof}

We now classify the manifolds modelled on $\Sol$. Recall that every semi-bundle is doubly covered by a canonical bundle. A torus bundle is \emph{of Anosov type} if its monodromy is Anosov. A torus semi-bundle is of Anosov type if its double-covering is.

\begin{prop} \label{Sol:bundle:prop}
The interior of a compact orientable manifold $M$ admits a finite-volume complete $\Sol$ geometry if and only if it is a torus (semi-)bundle of Anosov type.
\end{prop}
\begin{proof}
If $\interior M= \Sol/_\Gamma$ has finite volume, then both $\Gamma \cap \ker p$ and $p(\Gamma)$ are discrete and cofinite. The horizontal foliation of $\Sol$ into Euclidean planes $z=k$ is preserved by $\Gamma$ and projected to a surface bundle over the 1-orbifold $\matR /_{p(\Gamma)}$, which is either $S^1$ or an interval. The fibres are flat and have finite area, hence they are tori or Klein bottles. 

If $p(\Gamma)$ acts by translations, then $M$ is a torus bundle with an Anosov monodromy $\pm\matr {e^{-z}} 0 0 {e^z}$ for some $z>0$. If $M$ is a semi-bundle, this holds on its double-covering.

Conversely, let $M_A$ be a torus bundle with Anosov monodromy $A\in \SLZ$. We have $M_A = \matR^3/_\Gamma$ where $\Gamma$ is generated by the affine maps
$$(x,y,z) \mapsto (x+1,y,z), \quad (x,y,z) \mapsto (x,y+1,z), \quad (\bar x,z) \mapsto 
(A\bar x, z+h)$$
where $\bar x = (x,y)$ and $h>0$ is any positive number. A linear isomorphism $\psi$ of $\matR^2$ conjugates $A$ into a diagonal matrix $\pm\matr {e^{-a}} 0 0 {e^a}$ for some $a>0$ and transforms the generators into
$$(\bar x,z) \mapsto (\bar x +\bar v ,z), \ (\bar x,z) \mapsto (\bar x + \bar w, z), \ (x,y,z) \mapsto (\pm e^{-a}x, \pm e^{a}y, z+h)$$
where $\bar v = \psi(1,0)$ and $\bar w = \psi(0,1)$. If $\tr A>0$ and $h=a$ these maps are left-multiplications by the following elements of $\Sol$:
$$(\bar v,0), \quad (\bar w,0), \quad (\bar 0,a)$$
where $\bar 0 = (0,0)$, and we are done. If $\tr A <0$ we compose the third generator with the isometry $(x,y,z)\mapsto (-x,-y,z)$. 

If $M$ is a semi-bundle double-covered by $M_A$, we represent it as $\matR^3/_{\Gamma'}$ where $\Gamma'$ is generated by the group $\Gamma$ representing $M_A$ plus the element
$$(x,y,z) \mapsto \left(x+\frac 12, -y, -z\right)$$
which will be transformed via $\psi$ into an isometry of $\Sol$.
\end{proof}

\begin{cor} \label{Sol:not:Nil:cor}
If a closed three-manifold $M$ has a $\Sol$ geometry, then $\pi_1(M)$ is virtually solvable but not virtually nilpotent.
\end{cor}
\begin{proof}
It is virtually solvable because $\Sol$ has finite index in $\Iso^+(\Sol)$. The proof of Proposition \ref{Sol:bundle:prop} shows that
up to finite index $M=\Sol/_\Gamma$ and $\Gamma$ is generated by
$$(\bar v,0), \quad (\bar w,0), \quad (\bar 0,a).$$
Exercise \ref{Sol:conti:ex} implies that $[(\bar 0,a), (\bar v,0)] = (\bar v',0)$ for some $\bar v' \neq 0$ and hence $\Gamma$ is not nilpotent. This holds for any $\Gamma$ of this kind, and therefore $\Gamma$ is not virtually nilpotent, either.
\end{proof}

\section{Summary}
We have proudly completed the proof of Theorem \ref{Seifert:geometrisation:teo}. We state that theorem again here:

\begin{teo} 
A closed orentable 3-manifold has a geometric structure modelled on one of the following six geometries:
$$S^3,\ \matR^3,\ S^2\times \matR,\ \matH^2\times \matR,\ \Nil,\ \widetilde{\SL_2}$$
if and only if it is a Seifert manifold of the appropriate commensurability class, as prescribed by Table \ref{Seifert:geometry:table}. It has a $\Sol$ geometric structure if and only if it is a torus (semi-)bundle of Anosov type.
\end{teo}

We can deduce that the eight geometries are mutually exclusive:

\begin{prop}
Two closed 3-manifolds admitting different geometries are not diffeomorphic, and not even commensurable.
\end{prop}
\begin{proof}
We already know that the six Seifert geometries form six distinct commensurability classes. A manifold of type $\Sol$ is not Seifert by Proposition \ref{torus:bundle:prop}. A closed hyperbolic manifold is neither Seifert nor Sol because its fundamental group is infinite and does not contain an infinite normal abelian subgroup by Corollary \ref{normal:abelian:cor} (as opposite to the other geometries, see Proposition \ref{infinite:K:prop}).

Two manifolds having different geometries cannot be commensurable, since a common finite cover would inherit both geometries. 
\end{proof}

In fact one can tell the geometry of the manifold directly from its fundamental group. 

\begin{prop}
Let $M$ be a closed manifold modelled on one $\matX$ of the eight geometries; this flowchart shows how to determine $\matX$ from $\pi_1(M)$:
\begin{itemize}
\item if $\pi_1(M)$ is finite, then $\matX = S^3$; otherwise
\item if $\pi_1(M)$ is virtually cyclic, then $\matX = S^2\times \matR$; otherwise
\item if $\pi_1(M)$ is virtually abelian, then $\matX = \matR^3$; otherwise
\item if $\pi_1(M)$ is virtually nilpotent, then $\matX = \Nil$; otherwise
\item if $\pi_1(M)$ is virtually solvable, then $\matX = \Sol$; otherwise
\item if $\pi_1(M)$ contains a normal cyclic group $K$, then:
\begin{itemize}
\item if a finite-index subgroup of the quotient lifts, then $\matX = \matH^2 \times \matR$, 
\item otherwise $\matX = \widetilde{\SL_2}$; 
\end{itemize}
\item otherwise $\matX = \matH^3$.
\end{itemize}
\end{prop}
\begin{proof}
The virtually abelian cases were settled in Section \ref{virtually:abelian:subsection}. The $\Nil$ and $\Sol$ geometries were considered by Corollaries \ref{Nil:is:Nil:cor} and \ref{Sol:not:Nil:cor}. 

If $\matX = \matH^2\times \matR$ or $\widetilde{\SL_2}$ then 
$\pi_1(M)$ surjects onto $\pi_1(S)$ for some closed hyperbolic surface $S$. 
This implies that $\pi_1(M)$ is not solvable, for otherwise
(by Propositions \ref{subgroup:prop} and \ref{solvable:normal:prop})
the group $\pi_1(S)$ would also be, and hence it would contain a normal cyclic subgroup, contradicting Corollary \ref{normal:abelian:cor}.

If $\matX = \widetilde{\SL_2}$ then $\pi_1(S) = \pi_1(M)/_K$ and no finite-index subgroup of $\pi_1(S)$ can lift to $\pi_1(M)$ by Proposition \ref{SL2:no:lift:prop}.

Finally, if $\matX=\matH^3$ then $\pi_1(M)$ contains no normal cyclic group $K$ by Corollary \ref{normal:abelian:cor} again.
\end{proof}

The boundary case is slightly different. Of the eight geometries, only three produce non-compact finite-volume complete orientable manifolds:
$$\matH^3,\ \matH^2 \times \matR,\ \widetilde{\SL_2}.$$
These non-compact finite-volume manifolds are diffeomorphic to the interior of a compact manifold $M$ with boundary consisting of tori.

\begin{teo}
Let $M$ be a compact orientable 3-manifold with boundary consisting of tori. The following are equivalent:
\begin{itemize}
\item $M$ has a finite-volume complete $\matH^2 \times \matR$ structure,
\item $M$ has a finite-volume complete $\widetilde{\SL_2}$ structure,
\item $M$ is Seifert with $\chi<0$.
\end{itemize}
\end{teo}

In contrast with the closed case, a non-closed $M$ may admit two different geometric structures of Seifert type. On the other hand, a manifold $M$ cannot admit both a hyperbolic and a Seifert structure. 

\section{The geometrisation conjecture} \label{geometrisation:section}

We can finally state the famous Geometrisation Conjecture, proposed by Thurston in 1982 and proved twenty years later by Perelman in 2002.

\subsection{Statement and main consequences}
We say that a compact 3-manifold with (possibly empty) boundary consisting of tori is \emph{geometric} if its interior has a finite-volume complete geometric structure modelled on one of the eight geometries:\index{three-manifold!geometric three-manifold} 
$$S^3,\ \matR^3,\ \matH^3,\ S^2\times \matR,\ \matH^2\times \matR,\ \Nil,\ \Sol,\ \widetilde{\SL_2}.$$
The following conjecture was formulated by Thurston in 1982:\index{geometrisation conjecture}

\begin{conj}[Geometrisation Conjecture]
Let $M$ be an irreducible orientable compact 3-manifold with (possibly empty) boundary consisting of tori. Every block of the geometric decomposition of $M$ is geometric.
\end{conj}

The conjecture has been proved by Perelman in 2002 and its proof goes very very far from the scope of this book. It is however quite easy to deduce important consequences from it.

\begin{conj}[Poincar\'e conjecture] \label{Poincare:conj}
Every simply connected closed 3-manifold $M$ is diffeomorphic to $S^3$.
\end{conj}
\begin{proof}[Proof using geometrisation]
Via the prime decomposition we may restrict to the case $M$ is prime, hence irreducible. The group $\pi_1(M)$ is trivial and hence does not contain $\matZ \times \matZ$: every torus in $M$ is thus compressible and the geometric decomposition is trivial. By geometrisation $M$ is itself geometric. The only geometry with finite fundamental groups is $S^3$, and hence $M= S^3/_\Gamma$ is elliptic. Since $M$ is simply connected, the group $\Gamma = \pi_1(M)$ is trivial and hence $M=S^3$.
\end{proof}

\begin{conj}[Elliptisation]  \label{finite:conj}
Every closed 3-manifold $M$ with finite $\pi_1(M)$ is elliptic.
\end{conj}
\begin{proof}[Proof using geometrisation]
Same proof as above. Note that this is \emph{not} a consequence of Poincar\'e conjecture in general, for a manifold covered by $S^n$ needs not to be elliptic a priori.
\end{proof}

\begin{conj}[Hyperbolisation] \label{hyperbolisation:conj}
Every closed irreducible 3-manifold $M$ with infinite $\pi_1(M)$ not containing $\matZ \times \matZ$ is hyperbolic.
\end{conj}
\begin{proof}[Proof using geometrisation]
Since $\pi_1(M)$ does not contain $\matZ\times \matZ$ every torus is compressible and the geometric decomposition of $M$ is trivial. By geometrisation $M$ is geometric. Its geometry is not $S^3$ since $\pi_1(M)$ is infinite, and is not $S^2\times \matR$ since $M$ is irreducible. In the other Seifert geometries and in $\Sol$ the fundamental group $\pi_1(M)$ always contain a $\matZ\times \matZ$ (there is always a finite covering containing an incompressible torus).
\end{proof}

\begin{cor} 
Let $\tilde M \to M$ be a finite covering. If $\tilde M$ is geometric, then $M$ also is (with the same geometry).
\end{cor}
\begin{proof}[Proof using geometrisation]
Note that this is stronger than Conjecture \ref{finite:conj}. 
Since $\tilde M$ is geometric, it is irreducible and hence also $M$ is. Proposition \ref{JSJ:lifts:prop} implies that the geometric decomposition of $M$ is trivial, and by geometrisation $M$ is geometric.
\end{proof}

Concerning bounded manifolds, we get the following simple statement.

\begin{cor} \label{simple:bounded:hyperbolic}
Every simple compact manifold $M$ bounded by a non-empty collection of tori is hyperbolic.
\end{cor}
\begin{proof}[Proof using geometrisation]
Being simple, it is geometric. Seifert manifolds with boundary are never simple.
\end{proof}

That statement is not true for closed manifolds, because many Seifert manifolds fibering over $S^2$ with at most 3 singular fibres are simple. 

\subsection{Surface bundles}
There is a nice way to formulate geometrisation for surface bundles. 
We start with the much simpler torus case.

\begin{prop}
Let $M_A$ be a torus bundle with monodromy $A \neq \pm I$. The following holds:
\end{prop}
\begin{itemize}
\item if $|\tr A|<2$, \emph{i.e.}~$A$ has finite order, then $M_A$ is flat;
\item if $|\tr A|=2$, \emph{i.e.}~$A$ is reducible, then $M_A$ is $\Nil$;
\item if $|\tr A|>2$, \emph{i.e.}~$A$ is Anosov, then $M_A$ is $\Sol$.
\end{itemize}
\begin{proof}
Use Proposition \ref{torus:bundle:prop}.
\end{proof}

We now turn to the generic case. Let $\Sigma$ be a closed orientable surface with $\chi(\Sigma)<0$.
\begin{teo} \label{fibering:geometrisation:teo}
Let $M_\psi$ be a surface bundle with fibre $\Sigma$ and monodromy $\psi \in \MCG(\Sigma)$. The following holds:
\begin{itemize}
\item if $\psi$ has finite order, then $M_\psi$ has a $\matH^2 \times \matR$ geometry,
\item if $\psi$ is reducible, then $M_\psi$ contains an essential torus,
\item if $\psi$ is pseudo-Anosov, then $M_\psi$ is hyperbolic.
\end{itemize}
\end{teo}
\begin{proof}[Proof using geometrisation]
Use Proposition \ref{trichotomy:prop}.
\end{proof}

Theorem \ref{fibering:geometrisation:teo} was initially proved by Thurston in the 1980s. In the same years Thurston also proved the geometrisation conjecture for all Haken 3-manifolds. Before Perelman's proof the conjecture was open ``only'' in the non-Haken case, and it naturally split in three parts: the Conjectures \ref{Poincare:conj}, \ref{finite:conj}, and \ref{hyperbolisation:conj}. Perelman's proof certifies geometrisation in all cases with a unified technique.

\subsection{References}
This chapter contains many technical proofs. Most of them were taken from Scott \cite{S}, that is the standard reference on the geometrisation of Seifert manifolds, and Thurson's book \cite{Th_book}. We have also consulted a nice survey of Bonahon \cite{B3}, that contains in particular the suggestion of calculating the Ricci tensors in the non-product geometries. Thurson's geometrisation conjecture appears in \cite{Th_conj} in 1982, while Perelman's proof consists of three papers that he sent to the arXiv in 2003 \cite{Pe1, Pe2, Pe3}.

%% file: Mostow.tex
\chapter{Mostow rigidity theorem} \label{Mostow:chapter}
We have defined in Chapter \ref{Teichmuller:chapter} the
Teichm\"uller space $\Teich(S_g)$ of a genus-$g$ closed orientable surface $S_g$ as the space of all the hyperbolic metrics on $S_g$, considered up to isometries isotopic to the identity; we have then proved that $\Teich(S_g)$ is homeomorphic to $\matR^{6g-6}$ using the Fenchel--Nielsen coordinates.

This definition of $\Teich(M)$ actually applies to any closed hyperbolic manifold $M$, and we show here a striking difference between the dimensions two and three: if $\dim M = 3$ then $\Teich(M)$ is a single point. This strong result is known as the \emph{Mostow rigidity Theorem}.\index{Mostow--Prasad rigidity theorem} 

The impact of Mostow's rigidity on our understanding of three-dimensional topology cannot be overestimated. Thanks to this theorem every geometric information on a given closed hyperbolic three-manifold $M$ like its volume, geodesic spectrum, etc.~is promoted to a \emph{topological} invariant of $M$, that is it depends on the differentiable structure of $M$ only. In its strongest version, Mostow's rigidity says that the hyperbolic metric of $M$ is fully determined by the group $\pi_1(M)$ alone. 

We expose here Gromov's proof of Mostow's rigidity, which uses hyperbolic tetrahedra and introduces a nice invariant on closed manifolds of any dimension, called the \emph{simplicial volume}.\index{simplicial volume}

\section{Volume of tetrahedra} \label{volume:tetrahedra:section}
The volume of an ideal hyperbolic tetrahedron is a simple (but integral) formula that involves the \emph{Lobachevsky function}: we now prove this formula here. As a consequence we show that the regular ideal tetrahedron is the hyperbolic tetrahedron with maximal volume.

\subsection{The Lobachevsky function}
The \emph{Lobachevsky function} is\index{Lobachevsky function}
$$\Lambda(\theta) = -\int_0^\theta \log |2\sin t| dt.$$
The function $\log |2\sin t|$ is $-\infty$ on $\pi\matZ$ but is integrable, hence $\Lambda$ is well-defined and continuous on $\matR$. Its first derivatives are
$$\Lambda'(\theta) = -\log |2\sin \theta |, \qquad \Lambda''(\theta) = -\cot \theta.$$
The function $\Lambda$ has derivative $+\infty$ on $\pi\matZ$ and is an odd function, because its derivative is even. 

\begin{prop} \label{Loba:prop}
The function $\Lambda$ is $\pi$-periodic. We have $\Lambda (0) = \Lambda \left(\frac \pi 2 \right) = \Lambda (\pi) = 0$. The function $\Lambda$ is strictly positive on $\left(0,\frac \pi 2 \right)$, strictly negative on $\left(\frac \pi 2, \pi\right)$, and has absolute maximum and minimum at $\frac \pi 6 $ and $\frac 56 \pi$.
For all $m\in \matN$ the following holds:
$$\Lambda(m\theta) = m \sum_{k=0}^{m-1} \Lambda\left(\theta + \frac{k\pi}m\right).$$
\end{prop}
\begin{proof}
We prove the equality for $m=2$:
\begin{align*}
\frac{\Lambda(2\theta)}2 & = -\frac 12\int_0^{2\theta} \log |2\sin t|dt = - \int_0^\theta \log |2\sin 2t| dt \\
 & = - \int_0^\theta \log |2\sin t| dt - \int_0^\theta \log \left|2\sin \left(t+\frac \pi 2 \right) \right| dt \\
 & = \Lambda(\theta) - \int_{\frac \pi 2}^{\frac \pi 2 + \theta} \log |2 \sin t| dt \\
  & = \Lambda(\theta) + \Lambda\left(\theta + \frac \pi 2 \right) - \Lambda \left(\frac \pi 2 \right).
\end{align*}
By setting $\theta = \frac \pi 2$ we get $\Lambda(\pi)=0$. Since the derivative $\Lambda'$ is $\pi$-periodic and $\Lambda(\pi)=0$, also $\Lambda$ is $\pi$-periodic. Since $\Lambda$ is $\pi$-periodic and odd, we have $\Lambda \left(\frac \pi 2 \right) = 0$. We have also proved the formula for $m=2$.

To prove the formula for generic $m$ we use a generalisation of the duplication formula for the sinus. From the equality
$$z^m - 1 = \prod_{k=0}^{m-1} \left(z-e^{-\frac {2\pi ik}m}\right)$$
we deduce
$$2\sin(mt) = \prod_{k=0}^{m-1} 2\sin \left(t+\frac{k\pi}m\right)$$
and hence
\begin{align*}
\frac{\Lambda(m\theta)}m & = -\frac 1m \int_0^{m\theta} \log | 2\sin t| dt = - \int_0^\theta \log | 2 \sin (mt)| dt  \\
 & = - \sum_{k=0}^{m-1} \int_0^\theta \log \left|2 \sin \left( t + \frac {k \pi}m \right) \right| dt \\
 & = - \sum_{k=0}^{m-1} \left( \int_0^{\theta + \frac{k\pi}m} \log |2 \sin t| dt - \int_0^{\frac {k\pi}m} \log | 2\sin t | dt \right) \\
 & = - \sum_{k=0}^{m-1} \Lambda \left(\theta + \frac{k\pi}m \right) + C(m)
 \end{align*}
where $C(m)$ is a constant independent of $\theta$. By integrating both sides we get 
$$\frac 1m \int_0^\pi \Lambda(m\theta) = -\sum_{k=0}^{m-1}\int_0^\pi \Lambda\left(\theta + \frac{k\pi}m \right) + C(m)\pi.$$
Since $\Lambda$ is odd and $\pi$-periodic, we have
$$\int_0^\pi \Lambda (m\theta) = 0$$
for any integer $m$. Hence $C(m)=0$ and the formula is proved. Finally we note that $\Lambda''(\theta) = -\cot \theta $ is strictly negative in $(0,\frac \pi 2)$ and strictly positive in $(\frac \pi 2, \pi)$, hence $\Lambda$ is strictly positive in $(0, \frac \pi 2)$ and strictly negative in $(\frac \pi 2, \pi)$.
\end{proof}

\subsection{Volumes of ideal tetrahedra}
An ideal tetrahedron in $\matH^3$ is the convex hull of four non-planar ideal points. Quite surprisingly, every ideal tetrahedron has some non-trivial symmetries.\index{ideal tetrahedron} 

\begin{figure}
\begin{center}
\includegraphics[width = 4.5 cm] {\iftoggle{BW}{tetraedro4-BW}{tetraedro4}} 
\nota{Every pair of opposite edges in an ideal tetrahedron has an axis orthogonal to both which is a symmetry axis for the tetrahedron.}
\label{tetraedro4:fig}
\end{center}
\end{figure}

\begin{prop} \label{symmetries:ideal:tetrahedron:prop}
For any pair of opposite edges in an ideal tetrahedron $\Delta$ there is a unique line $r$ orthogonal to both as in Figure \ref{tetraedro4:fig} and $\Delta$ is symmetric with respect to a $\pi$-rotation around $r$.
\end{prop}
\begin{proof}
The opposite edges $e$ ed $e'$ are ultraparallel lines in $\matH^3$ and hence have a common perpendicular $r$. A $\pi$-rotation around $r$ inverts both $e$ and $e'$ but preserve the 4 ideal vertices of $\Delta$, hence $\Delta$ itself. 
\end{proof}

\begin{figure}
\begin{center}
\includegraphics[width = 10 cm] {\iftoggle{BW}{Loba-BW}{Loba}} 
\nota{The dihedral angles $\alpha, \beta, \gamma$ of an ideal tetrahedron. Opposite edges have the same angle and $\alpha+\beta+\gamma = \pi$ (left). To calculate the volume we use the half-space model, send a vertex to $\infty$, and divide the tetrahedron in six sub-tetrahedra (right).}
\label{Loba:fig}
\end{center}
\end{figure}

As a consequence, two opposite edges in $\Delta$ have coinciding dihedral angles as in Figure \ref{Loba:fig}-(left). Moreover, we have $\alpha + \beta + \gamma = \pi$ because a small horosphere based at a vertex intersects $\Delta$ into a Euclidean triangle with inner angles $\alpha, \beta$, and $\gamma$.
The regular ideal tetrahedron has of course equal angles $\alpha = \beta = \gamma = \frac \pi 3$.

\begin{teo} Let $\Delta$ be an ideal tetrahedron with dihedral angles $\alpha$, $\beta$ and $\gamma$. We have
$$\Vol(\Delta) = \Lambda(\alpha) + \Lambda(\beta) + \Lambda(\gamma).$$
\end{teo}
\begin{proof}
We represent $\Delta$ in the half-space model $H^3$ with one vertex $v_0$ at infinity and three vertices $v_1,v_2,v_3$ in $\matC$. Let $C$ be the circle containing $v_1$, $v_2$, and $v_3$: up to composing with elements in $\PSLC$ we can suppose that $C=S^1$. The Euclidean triangle $T\subset\matC$ with vertices $v_1$, $v_2$, and $v_3$ has interior angles $\alpha$, $\beta$, and $\gamma$. 

We first consider the case $0\in T$, that is $\alpha, \beta, \gamma \leqslant \frac \pi 2$. 
We decompose $T$ into six triangles as in Figure \ref{Loba:fig}: the tetrahedron $\Delta$ decomposes accordingly into six tetrahedra lying above them, and we prove that the one $\Delta_\alpha$ lying above the \iftoggle{BW}{light grey}{yellow} triangle has volume $\frac{\Lambda(\alpha)}2$. 
This proves the theorem.

The tetrahedron $\Delta_\alpha$ is the intersection of four half-spaces: three vertical ones bounded by the hyperplanes $y=0$, $x=\cos\alpha$, and $y=x\tan\alpha$, and one bounded by the half-sphere $z^2 = x^2+y^2$. Therefore
\begin{align*}
\Vol(\Delta_\alpha) & = \int_0^{\cos\alpha} dx \int_0^{x\tan \alpha} dy \int_{\sqrt{1-x^2-y^2}}^\infty \frac 1{z^3} dz \\
& = \int_0^{\cos\alpha} dx \int_0^{x\tan \alpha} dy \left[- \frac 1{2z^2} \right]_{\sqrt{1-x^2-y^2}}^\infty \\
& = \frac 12 \int_0^{\cos \alpha}dx \int_0^{x\tan\alpha}\frac 1{1-x^2-y^2}dy. 
\end{align*}
To solve this integral we use the relation 
$$\frac 1{1-x^2-y^2} = \frac 1{2\sqrt{1-x^2}} \left(\frac 1{\sqrt{1-x^2}-y} + \frac 1{\sqrt{1-x^2}+y}\right)$$
and hence $\Vol(\Delta_\alpha)$ equals
\begin{align*}
& \frac 14\int_0^{\cos\alpha} \!\!\!\!\!\! \frac {dx}{\sqrt{1-x^2}} 
\left(\left[-\log(\sqrt{1-x^2}-y)\right]_0^{x\tan \alpha} \!\!\!\!\! + \left[\log(\sqrt{1-x^2}+y)\right]_0^{x\tan\alpha}\right) \\
= & \frac 14 \int_0^{\cos\alpha} \!\!\!\!\! \frac{dx}{\sqrt{1-x^2}}
 \left(-\log(\sqrt{1-x^2}-x\tan\alpha)+\log(\sqrt{1-x^2}+x\tan\alpha)\right).
 \end{align*}
By writing $x=\cos t$ and hence $dx = -\sin t \, dt$ we obtain
\begin{align*} 
\Vol(\Delta_\alpha)  & = \frac 14 \int_{\frac\pi 2}^\alpha \frac{-\sin t}{\sin t} 
\left(-\log \frac{\sin t\cos\alpha - \cos t\sin \alpha}{\sin t \cos\alpha + \cos t \sin \alpha} \right) dt \\
 & = - \frac 14 \int_{\frac \pi 2}^\alpha \log \frac{\sin(t+\alpha)}{\sin(t-\alpha)}dt
 = - \frac 14 \int_{\frac \pi 2}^\alpha \log \frac{|2\sin(t+\alpha)|}{|2\sin(t-\alpha)|}dt \\
 & = \frac 14 \int_{2\alpha}^{\frac \pi 2+\alpha} \log | 2\sin t|dt - \frac 14 \int_0^{\frac \pi 2-\alpha} \log |2\sin t|dt \\
& = \frac 14\left(-\Lambda\left(\frac \pi 2 + \alpha \right) + \Lambda(2\alpha) + \Lambda\left(\frac \pi 2 - \alpha\right) \right) \\
& = \frac 14 \left(-\Lambda \left(\frac \pi 2 + \alpha \right) +2 \Lambda(\alpha) +2\Lambda\left(\frac \pi 2 + \alpha\right) - \Lambda\left(\frac \pi 2 + \alpha\right) \right) = \frac 12 \Lambda(\alpha)
\end{align*}
using Proposition \ref{Loba:prop}.

If $0\not \in T$ the triangle $T$ may be decomposed analogously into triangles, some of which contribute negatively to the volume, and we obtain the same formula.
\end{proof}

\begin{cor} \label{v3:cor}
The regular ideal tetrahedron is the hyperbolic tetrahedron of maximum volume.
\end{cor}
\begin{proof}
It is easy to prove that every hyperbolic tetrahedron is contained in an ideal tetrahedron: hence we may consider only ideal tetrahedra. Consider the triangle $T = \{0\leqslant \alpha, \beta, \alpha+\beta \leqslant \pi\}$ and 
\begin{align*}
f \colon \quad \ \ T & \longrightarrow  \matR \\ 
(\alpha, \beta) & \longmapsto \Lambda(\alpha) + \Lambda(\beta) + \Lambda(\pi-\alpha-\beta).
\end{align*}
The continuous function $f$ is null on $\partial T$ and strictly positive on the interior of $T$ because it measures the volume of the ideal tetrahedron of dihedral angles $\alpha, \beta, \gamma = \pi-\alpha-\beta$. Hence $f$ has at least a maximum on some interior point $(\alpha,\beta)$. The gradient 
$$
\nabla f = \begin{pmatrix}
\Lambda'(\alpha) - \Lambda'(\pi - \alpha-\beta) \\ \Lambda'(\beta) - \Lambda'(\pi-\alpha-\beta) 
\end{pmatrix}
 =
\begin{pmatrix}
-\log |2\sin \alpha| + \log |2\sin (\pi-\alpha-\beta)| \\ 
-\log |2\sin \beta| + \log |2\sin (\pi-\alpha-\beta) |
\end{pmatrix}
$$
must vanish there, and this holds if and only if $\sin\alpha = \sin (\pi-\alpha-\beta) = \sin \beta$, \emph{i.e.}~if and only if the tetrahedron has all dihedral angles equal to $\frac \pi 3$.
\end{proof}

\section{Simplicial volume}

Gromov has introduced a topological invariant on closed manifolds of any dimension called the \emph{simplicial volume}. This nice invariant can be used, among other things, to prove Mostow's rigidity theorem.\index{simplicial volume}

\subsection{Definition}
Gromov has introduced a measure of ``volume'' of a closed manifold $M$ which makes use only of the homology of $M$. Quite surprisingly, this notion of volume coincides (up to a factor) with the Riemannian one when $M$ is hyperbolic. 

Consider a topological space $X$ and its homology with ring $\matR$. We define the \emph{norm} of a cycle $\alpha = \lambda_1\alpha_1+\ldots +\lambda_h\alpha_h$ as follows:
$$|\alpha| = |\lambda_1|+\ldots +|\lambda_h|.$$

\begin{defn} The \emph{norm} of a class $a \in H_k(X,\matR)$ is the infimum of the norms of its elements: 
$$|a| = \inf \big\{|\alpha| \ \big|\ \alpha \in Z_k(X,\matR), [\alpha] = a \big\}.$$
\end{defn}

Recall that a \emph{seminorm} on a real vector space $V$ is a map $|\cdot|\colon V \to \matR_{\geqslant 0}$ such that
\begin{itemize}
\item $|\lambda v| = |\lambda||v|$ for any scalar $\lambda\in \matR$ and vector $v\in V$,
\item $|v+w| \leqslant |v|+|w|$ for any pair of vectors $v,w \in V$.
\end{itemize}
A \emph{norm} is a seminorm where $|v|=0$ implies $v=0$. The following is immediate.
\begin{prop} The norm $|\cdot |$ induces a seminorm on $H_k(X,\matR)$.
\end{prop}

Although it is only a seminorm, the function $|\cdot|$ is called a norm for simplicity. Let now $M$ be an oriented closed connected manifold: we know that $H_n(M,\matZ)\isom \matZ$ and the orientation of $M$ determines a \emph{fundamental class} $[M] \in H_n(M,\matZ)$ that generates the group. Moreover $H_n(M,\matR)\isom \matR$ and there is a natural inclusion\index{fundamental class}  
$$\matZ \isom H_n(M,\matZ) \hookrightarrow H_n(M,\matR) \isom \matR$$
hence the fundamental class $[M]$ is also an element of $H_n(M,\matR)$ and has a norm. 
\begin{defn} The \emph{simplicial volume} $\|M\|\in\matR_{\geqslant 0}$ of a closed oriented connected $M$ is the norm of its fundamental class:
$$\|M\| = |[M]|$$
\end{defn}
Since $|[M]| = |-[M]|$ the simplicial volume actually does not depend on the orientation. 
When $M$ is non-orientable we set $\|M\| = \big\|\tilde M\big\|/2$ where $\tilde M$ is the orientable double cover of $M$. The definition of $\|M\|$ is relatively simple but has various non-obvious consequences.

\subsection{Properties}
A continuous map $f\colon M \to N$ between closed oriented $n$-manifolds induces a homomorphism $f_*\colon H_n(M,\matZ) \to H_n(N,\matZ)$, and recall that the \emph{degree} of $f$ is the integer $\deg f$ such that
$$f_*([M]) = \deg f \cdot [N].$$

\begin{prop} \label{grado:prop}
Let $f\colon M\to N$ be a continuous map between closed oriented manifolds. The following inequality holds:
$$\|M\| \geqslant |\deg f| \cdot \|N\|.$$
\end{prop}
\begin{proof}
Every description of $[M]$ as a cycle $\lambda_1\alpha_1+\ldots +\lambda_h\alpha_h$ induces a description of $f_*([M]) = \deg f[N]$ as a cycle $\lambda_1f\circ\alpha_1+\ldots +\lambda_hf\circ\alpha_h$ with the same norm (or less, if there is some cancelation). 
\end{proof}

\begin{cor}
If $M$ and $N$ are closed orientable and homotopically equivalent $n$-manifolds then $\|M\|=\|N\|$.
\end{cor}
\begin{proof}
A homotopic equivalence consists of two maps $f\colon M\to N$ and $g\colon N \to M$ whose compositions are both homotopic to the identity. In particular both $f$ and $g$ have degree $\pm 1$.
\end{proof}

\begin{cor}
If $M$ admits a continuous self-map $f\colon M \to M$ of degree at least two then $\|M\| = 0$.
\end{cor}
\begin{cor}
Every sphere $S^n$ has norm zero. More generally we have $\|M\times S^n\|=0$ for every closed $M$ and any $n\geqslant 1$.
\end{cor}
\begin{proof}
A sphere $S^n$ admits self-maps of degree $\geqslant 2$, and hence also the product $M\times S^n$ does.
\end{proof}
Among the genus-$g$ surfaces $S_g$, we deduce that the sphere and the torus have simplicial volume zero. We will see soon that every surface of genus $g\geqslant 2$ has positive simplicial volume. When the continuous map is a covering the inequality from Proposition \ref{grado:prop} is promoted to an equality.

\begin{prop}
If $f\colon M \to N$ is a degree-$d$ covering we have 
$$\|M\| = d \cdot \|N \|.$$
\end{prop}
\begin{proof} 
The reason for this equality is that cycles can be lifted and projected through the covering. More precisely, we already know that $\|M\|\geqslant d\cdot \|N\|$. Conversely, let $\alpha = \lambda_1\alpha_1+\ldots +\lambda_h\alpha_h$ represent $[N]$; each $\alpha_i$ is a map $\Delta_n\to N$. Since $\Delta_n$ is simply connected, the map $\alpha_i$ lfts to $d$ distinct maps $\alpha_i^1,\ldots, \alpha_i^d\colon \Delta_n\to N$. The chain $\tilde\alpha = \sum_{ij} \lambda_i\alpha_i^j$ is a cycle in $M$ and $f_*(\tilde \alpha) = d\alpha$. Hence $\|M\|\leqslant d\cdot \|N\|$.
\end{proof}

We also note the following fact.

\begin{prop}
If $M$ is triangulated with $k$ simplices, then $\|M\| \leqslant k$.
\end{prop}
\begin{proof}
Up to taking a double cover we may suppose that $M$ is oriented. The closed $n$-manifold $M$ is triangulated into some simplices $\Delta_1, \ldots, \Delta_k$, and
we fix an orientation-preserving parametrisation $s_i\colon \Delta \to \Delta_i$ of each.
We would like to say that $s_1+\ldots + s_k$ is a fundamental cycle, however this singular chain is not necessarily a cycle because the restriction to adjacent faces coincide only up to the symmetries $S_{n+1}$ of $\Delta$. 

We can fix this problem easily by averaging each $s_i$ on all its permutations, that is we substitute each $s_i$ with $\frac{1}{(n+1)!}\sum_{\sigma \in S_{n+1}} (-1)^{{\rm sgn}(\sigma)} s_i \circ \sigma$. Now $s=s_1+\ldots + s_k$ is a fundamental cycle and $|s| = k$.
\end{proof}

\subsection{Seifert manifolds}
We can now calculate the simplicial volume of Seifert manifolds.

\begin{prop}
If $M$ is a closed Seifert 3-manifold then $\|M\|=0$.
\end{prop}
\begin{proof}
Every Seifert manifold is finitely covered by a product $S \times S^1$ (if $e=0$) or by a bundle $\big(S,(1,1)\big)$ with Euler number 1 (if $e\neq 0$). In the first case we are done since $\|S\times S^1\|=0$. 

In the second, if $S=S^2$ then $\big(S^2,(1,1)\big)=S^3$ and $\|S^3\|=0$, so we suppose $\chi(S)\leqslant 0$. There is a universal $K>0$ such that $\big(S,(1,1)\big)$ triangulates with at most $K|\chi (S)|+K$ tetrahedra (exercise) and therefore $\big\|\big(S,(1,1)\big)\big\| \leqslant K|\chi(S)|+K$. Exercise \ref{e:construct:ex} shows that for every $e>0$ there is a degree-$e$ covering $\tilde S \to S$ and two degree-$e$ coverings
$$\big(\tilde S, (1,1)\big) \longrightarrow \big(\tilde S, (1,e)\big) \longrightarrow \big(S,(1,1)\big)$$
which compose to a degree-$e^2$ covering 
$\big(\tilde S, (1,1)\big) \to \big(S,(1,1)\big)$. Therefore
$$\big\|\big(S, (1,1) \big) \big\| = \frac{\big\|\big(\tilde S, (1,1) \big) \big\|}{e^2} \leqslant \frac{K|\chi(\tilde S)| + K}{e^2} \leqslant \frac{Ke|\chi(S)| + K}{e^2} \to 0$$
as $e \to \infty$.
\end{proof}

\subsection{Simplicial and hyperbolic volume}
In the next pages we will prove the following theorem. Let $v_3$ be the volume of the regular ideal tetrahedron in $\matH^3$.

\begin{teo} \label{Gromov:teo}
Let $M$ be a closed hyperbolic $3$-manifold. We have 
$$\Vol(M) = v_3\|M\|.$$
\end{teo}
The theorem furnishes in particular some examples of manifolds with positive simplicial volume and shows that $\Vol(M)$ is a topological invariant of $M$, thus generalising the Gauss-Bonnet theorem to dimension $n=3$. Mostow rigidity will then strengthen this result in dimension $n = 3$, showing that the hyperbolic metric itself (not only its volume) is a topological invariant.

Both quantities $\Vol(M)$ and $\|M\|$ are multiplied by $d$ if we substitute $M$ with a degree-$d$ covering. In particular, up to substituting $M$ with its orientable double cover we can suppose that $M$ is orientable. 

\subsection{Cycle straightening}
The \emph{straight singular $k$-simplex} with vertices $v_0,\ldots, v_{k} \in \matH^n$ is the map
\begin{align*}
\alpha\colon \quad \ \Delta_k & \longrightarrow \matH^n \\
(t_0,\ldots, t_{k}) & \longmapsto t_0v_0 + \ldots + t_{k}v_{k}
\end{align*}
defined using convex combinations, see Section \ref{convex:subsection}. If the $k+1$ vertices are not contained in a $(k-1)$-plane the singular $k$-simplex is \emph{non-degenerate} and its image is a hyperbolic $k$-simplex.

The \emph{straightening} $\alpha^\st$ of a singular simplex $\alpha\colon \Delta_k \to \matH^n$ is the straight singular simplex with the same vertices of $\alpha$. The straightening $\alpha^\st$ of a singular simplex $\alpha\colon\Delta_k\to M$ in a hyperbolic manifold $M=\matH^n/_\Gamma$ is defined by lifting the singular simplex in $\matH^n$, straightening it, and projecting it back to $M$ by composition with the covering map.
Different lifts produce the same straightening in $M$ because they are related by isometries of $\matH^n$. 

The straightening extends by linearity to a homomorphism
$$\st\colon C_k(M,R) \to C_k(M,R)$$
which commutes with $\partial$ and hence induces a homomorphism in homology
$$\st_*\colon H_k(M,R) \to H_k(M,R).$$
\begin{prop} The map $\st_*$ is the identity.
\end{prop}
\begin{proof}
We may define a homotopy between a singular simplex $\sigma$ and its straightening $\sigma^{\rm st}$ using the convex combination
$$\sigma^t(x) =t\sigma(x) + (1-t)\sigma^{\rm st}(x).$$
This defines a chain homotopy between $\st_*$ and $\id$ via the same technique used to prove that homotopic maps induce the same maps in homology.
\end{proof}

The \emph{abstract volume} of a straightened singular $n$-simplex $\alpha\colon\Delta_n \to M$ is the volume of its lift in $\matH^n$ and may also be calculated as
$$\left|\int_\alpha \omega\right|$$
where $\omega$ is the volume form on $M$ pulled back along $\alpha$. 
If $\alpha$ is non-degenerate, we say that its \emph{sign} is \emph{positive} if $\alpha$ is orientation-preserving and \emph{negative} otherwise: equivalently, it is the sign of $\int_\alpha\omega$.

We now concentrate on the dimension $n=3$ where we know that a tetrahedron has maximum volume $v_3$ if and only if it is both regular and ideal: a straight simplex is compact and hence its abstract volume is strictly smaller than $v_3$. 
We can now easily prove one inequality. 

\begin{prop} \label{Gromov:1:prop}
Let $M$ be a closed hyperbolic $3$-manifold. We have
$$\Vol(M) \leqslant v_3 \|M\|.$$
\end{prop}
\begin{proof}
As we said above, we can suppose that $M$ is orientable. Take a cycle $\alpha = \lambda_1\alpha_1+\ldots+\lambda_k\alpha_k$ that represents $[M]$. We can suppose it is straightened, because the straightening preserves both the coefficients and the homology class. Let $\omega$ be the volume form on $M$. We get 
$$\Vol (M) = \int_M \omega = \int_\alpha \omega = \lambda_1 \int_{\alpha_1}\omega + \ldots + \lambda_k \int_{\alpha_k}\omega.$$
The quantity $\big|\int_{\alpha_i}\omega\big|$ is the abstract volume of $\alpha_i$. Hence $\big|\int_{\alpha_i}\omega\big|<v_3$ and
$$\Vol(M) < \big(|\lambda_1|+ \ldots +|\lambda_k|\big)v_3.$$
This holds for all $\alpha$, hence $\Vol(M) \leqslant v_3 \|M\|$.
\end{proof}

The proof of the converse inequality is less immediate.

\subsection{Efficient cycles} \label{efficienti:subsection}
Let $M=\matH^3/_\Gamma$ be a closed oriented hyperbolic 3-manifold. 
An \emph{$\varepsilon$-efficient cycle} for $M$ is a straightened cycle
$$\alpha = \lambda_1\alpha_1 + \ldots + \lambda_k\alpha_k$$
representing $[M]$ where the abstract volume of $\alpha_i$ if bigger than $v_3-\varepsilon$ and the sign of $\alpha_i$ is coherent with the sign of $\lambda_i$, for all $i$.

We will construct an $\varepsilon$-efficient cycle for every $\varepsilon >0$. This will conclude the proof of Theorem \ref{Gromov:teo} in virtue of the following:

\begin{lemma}
If for every $\varepsilon>0$ the manifold $M$ admits an $\varepsilon$-efficient cycle, then we have $\Vol(M) \geqslant v_3 \|M\|$.
\end{lemma}
\begin{proof}
Let $\alpha = \lambda_1\alpha_1 + \ldots \lambda_k\alpha_k$ be an $\varepsilon$-efficient cycle and $\omega$ be the volume form on $M$. Coherence of signs gives $\lambda_i\int_{\alpha_i}\omega>0$ for all $i$.
We get
\begin{align*}
\Vol(M) & = \int_M \omega = \int_\alpha \omega = \lambda_1 \int_{\alpha_1} \omega + \ldots + \lambda_k \int_{\alpha_k} \omega  \\
& \geqslant \big(|\lambda_1| + \ldots + |\lambda_k|\big) \cdot (v_3 - \varepsilon).
\end{align*}
Therefore $\Vol(M) \geqslant \|M\| \cdot (v_3 - \varepsilon)$ for all $\varepsilon > 0$.
\end{proof}

It remains to construct $\varepsilon$-efficient cycles.

\begin{ex} \label{volume:continuo:ex}
If $\Delta^i$ is a sequence of tetrahedra in $\matH^3$ whose vertices tend to the vertices of a regular ideal tetrahedron in $\partial\matH^3$, then
$$\Vol(\Delta^i)\to v_3.$$
\end{ex}

For any $t>0$, let $\Delta(t)$ be a regular tetrahedron obtained as in Section \ref{platonici:subsection} as follows.
Pick a point $x\in \matH^3$ and a regular tetrahedron in the Euclidean tangent space $T_x\matH^3$, centred at the origin with vertices at distance $t$ from it, project the vertices in $\matH^3$ via the exponential map, and pick their convex hull.

A \emph{$t$-simplex} is a tetrahedron isometric to $\Delta(t)$ equipped with an ordering of its vertices: the ordering allows us to consider it as a straightened singular simplex. Let $S(t)$ be the set of all $t$-simplices in $\matH^3$.

\begin{ex} The group $\Iso(\matH^3)$ acts on $S(t)$ freely and transitively.
\end{ex}

Recall from Corollary \ref{unimodular:cor} that $\Iso(\matH^3)$ is unimodular:
the Haar measure on $\Iso(\matH^3)$ induces an $\Iso(\matH^3)$-invariant measure on $S(t)$. 

Let $M=\matH^3/_\Gamma$ be a closed hyperbolic 3-manifold and $\pi\colon\matH^3 \to M$ the covering projection. Fix a base point $x_0 \in \matH^3$ and consider its orbit $O=\Gamma x_0$. Consider the set
$$\Sigma = \Gamma^{4}/_\Gamma$$
of the $4$-uples $(g_0,g_1,g_2, g_3)$ considered up to the diagonal action of $\Gamma$:
$$g \cdot (g_0,g_1,g_2, g_3) = (g g_0,gg_1, gg_2, g g_3).$$
An element $\sigma = (g_0,g_1,g_2, g_3) \in \Sigma$ determines a singular simplex $\tilde\sigma$ in $\matH^3$ with vertices $g_0(x_0), g_1(x_0), g_2(x_0), g_3(x_0)\in O$ only up to translations by $g\in\Gamma$, hence it gives a well-defined singular simplex in $M$, which we still denote by $\sigma$. We now introduce the chain
$$\alpha(t) = \sum_{\sigma \in \Sigma} \lambda_\sigma(t) \cdot \sigma$$
for some appropriate real coefficients $\lambda_\sigma(t)$ that we now define. The base point $x_0$ determines the Dirichlet tessellation of $\matH^3$ into domains $D(g(x_0))$, $g\in\Gamma$. For $\sigma = (g_0,g_1,g_2, g_3)$ we let $S^+_\sigma(t)\subset S(t)$ be the set of all positive $t$-simplices whose $i$-th vertex lies in $D(g_i(x_0))$ for all $i$. The number $\lambda^+_\sigma(t)$ is the measure of $S^+_\sigma(t)$. We define analogously $\lambda^-_\sigma(t)$ and set
$$\lambda_\sigma(t) = \lambda^+_\sigma(t) - \lambda^-_\sigma(t).$$

\begin{lemma} \label{ciclo:lemma}
The chain $\alpha(t)$ has finitely many addenda and is a cycle. If $t$ is sufficiently big the cycle $\alpha(t)$ represents a positive multiple of $[M]$ in the group $H_3(M,\matR)$.
\end{lemma}
\begin{proof}
We prove that the sum is finite. Let $d, T$ be the diameters of $D(x_0)$ and of a $t$-simplex. 
We write $\sigma = (\id,g_1,g_2,g_3)$ for all $\sigma \in \Sigma$: that is, all simplices have their first vertex at $x_0$. 
If $\lambda_\sigma(t)\neq 0$ then $d(g_ix_0, x_0) < 2d+T$ for all $i$: therefore $\alpha(t)$ has finitely many addenda (because $O$ is discrete).

We prove that $\alpha(t)$ is a cycle. The boundary $\partial \alpha(t)$ is a linear combination of straight $2$-simplices with vertices in $(g_0x_0, g_1x_0, g_2x_0)$ as $g_0, g_1$, and $g_2$ vary. The coefficient of one such $2$-simplex is
$$\sum_{g\in\Gamma} \big(- \lambda_{(g,g_0,g_1,g_2)}(t) + \lambda_{(g_0,g,g_1,g_2)}(t)
- \lambda_{(g_0,g_1,g,g_2)}(t) + \lambda_{(g_0,g_1,g_2,g)}(t)\big).$$
We prove that each addendum summed along $g\in \Gamma$ is zero; for simplicity we consider the last addendum and get
$$\sum_{g\in\Gamma} \lambda_{(g_0, g_1, g_{2}, g)}(t) = \sum_{g\in\Gamma} \lambda_{(g_0, g_1, g_{2}, g)}(t)^+ - \sum_{g\in\Gamma} \lambda_{(g_0, g_1, g_{2}, g)}(t)^-.$$
The first addendum measures the positive $t$-simplices whose first $3$ vertices lie in $D(g_0(x_0)), \ldots, D(g_2(x_0))$, the second measures the negative $t$-simplices with the same requirement. These two subsets have the same volume in $S(t)$ because they are related by the involution $r\colon S(t) \to S(t)$ that mirrors a simplex with respect to its first facet. 

We show that for sufficiently big $t$ the cycle is a positive multiple of $[M]$. Let $t$ be sufficiently big so that two vertices in a $t$-simplex have distance 
bigger than $2d$. This condition implies that if there is a positive $t$-simplex with vertices in $D(g_0(x_0)), \ldots, D(g_3(x_0))$, then any straight simplex with vertices in $D(g_0(x_0)), \ldots, D(g_3(x_0))$ is positive. 
Therefore in the expression
$$\alpha(t) = \sum_{\sigma \in \Sigma} \lambda_\sigma(t) \cdot \sigma$$
the signs of $\lambda_\sigma(t)$ and $\sigma$ are coherent and
$$\int_{\alpha(t)} \omega = \sum_{\sigma\in \Sigma} \lambda_\sigma(t) \cdot \int_{\sigma}\omega > 0.$$
Therefore $\alpha(t)$ is a positive multiple of $[M]$.
\end{proof}

For sufficiently big $t$ we have $\alpha(t) = k_t[M]$ in homology for some $k_t>0$. The rescaled $\bar\alpha(t) = \alpha(t)/k_t$ hence represents $[M]$.
We have found our $\varepsilon$-efficient cycles. 

\begin{lemma} \label{efficiente:lemma}
For any $\varepsilon>0$ there is a $t_0>0$ such that $\bar\alpha(t)$ is $\varepsilon$-efficient for all $t>t_0$.
\end{lemma}
\begin{proof}
Let $d$ be the diameter of the Dirichlet domain $D(x_0)$. Let a \emph{quasi $t$-simplex} be a simplex whose vertices are at distance $<d$ from those of a $t$-simplex. By construction $\bar\alpha(t)$ is a linear combination of quasi $t$-simplices.

We now show that for any $\varepsilon>0$ there is a $t_0>0$ such that for all $t>t_0$ every quasi $t$-simplex has volume bigger than $v_3-\varepsilon$. By contradiction, let $\Delta^t$ be a sequence of quasi $t$-simplices of volume smaller than $v_3-\varepsilon$ with $t\to\infty$. The vertices of $\Delta^t$ are $d$-closed to a $t$-simplex $\Delta_*^t$, and we move the pair $\Delta^t, \Delta^t_*$ isometrically so that the $t$-simplices $\Delta_*^t$ have the same barycenter. Now both the vertices of $\Delta^t$ and $\Delta^t_*$ tend to the vertices of an ideal regular tetrahedron and Exercise \ref{volume:continuo:ex} gives a contradiction.
\end{proof}

The previous lemmas together prove the second half of Theorem \ref{Gromov:teo}.
\begin{cor}
Let $M$ be a closed hyperbolic $3$-manifold. We have
$$\Vol(M) \geqslant v_3 \|M\|.$$
\end{cor}

Theorem \ref{Gromov:teo} has some non-trivial consequences.

\begin{cor} \label{degree:d:cor}
Let  $M, N$ be two closed orientable hyperbolic $3$-manifolds. If there is a map $f\colon M\to N$ of degree $d$ then $\Vol(M)\geqslant |d| \cdot \Vol(N)$.
\end{cor}

\begin{cor} \label{stesso:volume:cor}
Two homotopically equivalent and closed hyperbolic 3-manifolds have the same volume.
\end{cor}

\begin{rem}
If we are able to prove that the regular ideal simplex has the maximum volume $v_n$ among hyperbolic $n$-simplices in $\matH^n$, then the whole proof extends as is from $\matH^3$ to $\matH^n$ and shows that $\Vol(M) = v_n \cdot |M|$ for any closed hyperbolic $n$-manifold $M$. This is obviously true when $n=2$, and we get $v_2 = \pi$. 
\end{rem}

\section{Mostow rigidity}

\subsection{Introduction}
We want to prove the following.
\begin{teo}[Mostow rigidity] \label{Mostow:teo}
Let $M$ and $N$ be two closed connected orientable hyperbolic 3-manifolds. Every isomorphism $\pi_1(M) \stackrel \sim\to \pi_1(N)$ between fundamental groups is induced by a unique isometry $M\stackrel \sim\to N$.
\end{teo}
This very powerful theorem says that an algebraic isomorphism between fundamental groups alone is enough to produce and characterise an isometry. 

\begin{cor} Two closed orientable hyperbolic 3-manifolds with isomorphic fundamental groups are isometric.
\end{cor}

We prove Mostow's rigidity in this section. All the ingredients are already there, we only need to make a little last effort.

We note that closed hyperbolic manifolds are aspherical because their universal cover $\matH^n$ is contractible. For such manifolds every isomorphism $\pi_1(M)\to \pi_1(N)$ is induced by a homotopy equivalence $f\colon M \to N$, unique up to homotopy: see Corollary \ref{aspherical:cor}. To prove Mostow's theorem we need to promote this homotopy equivalence to an isometry in dimension $n=3$. We already know that $\Vol(M) = \Vol(N)$ by Corollary \ref{stesso:volume:cor}.

\subsection{Proof of Mostow's theorem}
Let $f\colon M \to N$ be a smooth homotopy equivalence. Recall from Theorem \ref{estensione:teo} that $f$lifts to a map $\tilde f\colon \matH^3 \to \matH^3$ which extends continuously to a homeomorphism $\tilde f\colon\partial \matH^3 \to \partial\matH^3$ of the boundary spheres.
We start with a lemma.

\begin{lemma} \label{regolare:lemma}
The extension $\tilde f\colon \partial\matH^3 \to \partial \matH^3$ sends the vertices of every regular ideal simplex to the vertices of some regular ideal simplex. 
\end{lemma}
\begin{proof}
Let $w_0,\ldots, w_3$ be vertices of a regular ideal simplex and suppose by contradiction that their images $\tilde f(w_0),\ldots ,\tilde f(w_3)$ span a non-regular ideal simplex, which has volume smaller than $v_3-2\delta$ for some $\delta >0$. 
By Exercise \ref{volume:continuo:ex} there are neighbourhoods $U_i$ of $w_i$ in $\overline{\matH^3}$ for $i=0,\ldots, 3$ such that the volume of the simplex with vertices $\tilde f(u_0), \ldots, \tilde f(u_3)$ is smaller than $v_3-\delta$ for any choice of $u_i\in U_i$. 

In Section \ref{efficienti:subsection} we have defined a cycle
$$\alpha(t) = \sum_{\sigma \in \Sigma} \lambda_\sigma(t) \cdot \sigma$$
where $t$ depends on $\varepsilon$. We say that a singular simplex $\sigma\in \Sigma$ is \emph{bad} if its $i$-th vertex is contained in $U_i$ for all $i$. 
Let $\Sigma^{\rm bad}\subset\Sigma$ be the subset of all bad singular simplices and define
$$\alpha(t)^{\rm bad} = \sum_{\sigma \in \Sigma^{\rm bad}} \lambda_\sigma(t) \cdot \sigma.$$
We want to estimate $|\alpha(t)|$ and $|\alpha(t)^{\rm bad}|$. We prove that  
$$|\alpha(t)| = \sum_{\sigma\in \Sigma} |\lambda_\sigma(t)|$$
is a real number independent of $t$: let $S_0\subset S(t)$ be the set of all $t$-simplices having the first vertex in the Dirichlet domain $D(x_0)$ of the fixed base point $x_0\in\matH^3$. It follows from the definitions that $|\alpha(t)|$ equals the measure of $S_0$ for sufficiently big $t$. Moreover the set $S_0$ is in natural correspondence with the set of all isometries that send $x_0$ to some point in $D(x_0)$: its volume does not depend on $t$. 

To estimate $|\alpha(t)^{\rm bad}|$ we fix $g_0\in \Gamma$ so that $D(g_0x_0)\subset U_0$. Let $S^{\rm bad}\subset S(t)$ be the set of all bad $t$-simplices with first vertex in $D(g_0x_0)$. If $t$ is sufficiently big, the volume of $S^{\rm bad}$ is bigger than a constant independent of $t$ (exercise). 

We have proved that $|\alpha(t)^{\rm bad}| / |\alpha(t)|>C>0$ independently of $t$. We may suppose that $\alpha(t)$ represents $[M]$ up to renormalising. The map $f\colon M \to N$ has degree $\pm 1$ and it sends $\alpha(t)$ to a class 
$$f_*(\alpha(t)) = \sum_{\sigma\in \Sigma} \lambda_\sigma(t) \cdot (f\circ\sigma)^{\rm st} $$
representing $\pm [N]$. Since a $C$-portion of $\alpha(t)$ is bad, a $C$-portion of simplices in $f_*(\alpha(t))$ has volume smaller than $v_3-\delta$ and hence
$$\Vol(N) = \left|\int_{f_*(\alpha(t))} \omega\right| < |\alpha(t)| ( (1-C)v_3 + C(v_3 - \delta)) = |\alpha(t)| (v_3 -\delta C).$$
Since this holds for all $t$ and $|\alpha(t)| \to \|M\|$ we get
$$\Vol(N) < \|M\| (v_3 - \delta C) = \Vol (M) - \delta C \cdot \|M\|.$$
Corollary \ref{stesso:volume:cor} gives $\Vol(M) = \Vol(N)$, a contradiction.
\end{proof}

\begin{prop} \label{riflessione:prop}
Every ideal triangle in $\matH^3$ is the face of precisely two regular ideal tetrahedra.
\end{prop}
\begin{proof}
Pick the line $l$ orthogonal to the barycenter of the triangle: the vertex of a regular ideal tetrahedron must be an endpoint of $l$. To prove that these vertices give regular ideal tetrahedra, note that all ideal triangles in $\matH^3$ are isometric, so one concrete example suffices. 
\end{proof}

We turn back to Mostow rigidity.

\begin{prop}
Let $f\colon M \to N$ be a smooth homotopic equivalence between closed hyperbolic orientable 3-manifolds.  
The restriction $\tilde f|_{\partial \matH^3}\colon \partial \matH^3 \to \partial \matH^3$ is the trace of an isometry $\psi\colon \matH^3 \to \matH^3$.
\end{prop}
\begin{proof}
Let $v_0,\ldots, v_3\in\partial \matH^n$ be vertices of a regular ideal tetrahedron $\Delta$. By Lemma \ref{regolare:lemma} the lift $\tilde f$ sends them to the vertices of some regular ideal tetrahedron, and let $\psi$ be the unique isometry of $\matH^3$ such that $\psi (v_i) = \tilde f(v_i)$ for all $i$. 

By Proposition \ref{riflessione:prop} there is a unique point $v_4\neq v_3$ such that $v_0, v_1, v_2, v_4$ are the vertices of an ideal regular tetrahedron, and $\psi (v_4)$ is the unique point other than $\psi(v_3)$ such that $v_0,v_1,v_2, \psi(v_4)$ are the vertices of an ideal regular tetrahedron. By Lemma \ref{regolare:lemma} we also have $\tilde f(v_4) = \psi(v_4)$.

If we mirror $\Delta$ along its faces iteratively we get a tessellation of $\matH^3$ via regular ideal tetrahedra, whose ideal vertices form a dense subset of $\partial \matH^3$. By iterating this argument in all directions the functions $\psi$ and $\tilde f$ coincide on this dense subset and hence on the whole of $\partial \matH^3$.
\end{proof}

We can finally prove Mostow's rigidity theorem.
\begin{teo}[Mostow Rigidity] Let $f\colon M \to N$ be a homotopic equivalence between closed orientable hyperbolic 3-manifolds. The map $f$ is homotopically equivalent to an isometry.
\end{teo}
\begin{proof}
Set $M = \matH^3/_\Gamma$ and $N = \matH^3/_{\Gamma'}$, and pick a lift $\tilde f$. We have
\begin{equation} \label{equivariance:eqn}
\tilde f\circ g = f_*(g)\circ \tilde f \quad \forall g \in \Gamma
\end{equation}
for an isomorphism $f_*\colon \Gamma \to \Gamma'$. We may suppose $f$ smooth. The boundary extension of $\tilde f$ is the trace of an isometry $\psi\colon \matH^3 \to \matH^3$ and hence
\begin{equation} \label{equivariance2:eqn}
\psi\circ g = f_*(g)\circ \psi \quad \forall g \in \Gamma
\end{equation}
holds at $\partial \matH^3$. Both terms in (\ref{equivariance2:eqn}) are isometries, and isometries are determined by their boundary traces: hence (\ref{equivariance2:eqn}) holds for all points in $\matH^3$. Therefore $\psi$ descends to an isometry 
$$\psi\colon M \to N.$$
A homotopy between $f$ and $\psi$ may be constructed from a convex combination of $\tilde f$ and $\psi$ in $\matH^n$, which is also $\Gamma$-equivariant and hence descends.
\end{proof}

\subsection{Consequences of Mostow rigidity}
The most important consequence is that the entire geometry of a closed hyperbolic $3$-manifold is a topological invariant:  numerical quantities like the volume of the manifold, its geodesic spectrum, etc.~depend only on the topology of the manifold.
We single out another application. 

\begin{cor} Let $M$ be a closed orientable hyperbolic 3-manifold. The natural map
$$\Iso (M) \to \Out(\pi_1(M))$$
is an isomorphism.
\end{cor}
\begin{proof}
We already know that it is injective by Proposition \ref{injective:prop}. We prove that it is surjective:
every automorphism of $\pi_1(M)$ is represented by a homotopy equivalence since $M$ is aspherical (see Corollary \ref{aspherical:cor}), which is in turn homotopic to an isometry by Mostow's rigidity.
\end{proof}
We note that this is false in dimension $n=2$, where $\Iso(S)$ is finite and $\Out(\pi_1(S))$ is infinite. 

\subsection{Orbifolds and cone manifolds}
While by Mostow's theorem an orientable closed 3-manifold $M$ can have at most one hyperbolic structure, it may have plenty of hyperbolic cone manifold or orbifold structures, typically distinguished by their singular sets and their cone angles. For instance $M=S^3$ has plenty of such structures, as we will see.

For the moment we note that the volumes of these hyperbolic orbifold structures may be arbitrarily big, but not arbitrarily small:

\begin{prop} \label{Vol:M:O:prop}
Let a closed orientable 3-manifold $M$ be the underlying space of a hyperbolic orbifold $O$. We have
$$\Vol(O) \geqslant v_3\|M\|.$$
\end{prop}
\begin{proof}
By Selberg's Lemma there is a degree-$d$ orbifold covering $N \to O$ which is a closed hyperbolic manifold. We have $\Vol(N) = d\Vol(O)$, and since the covering is a degree-$d$ map $N \to M$ we can apply Corollary \ref{degree:d:cor} and get $\|M\| \leqslant \frac{\|N\|}d = \frac{\Vol(N)}{dv_3} = \frac{\Vol(O)}{v_3}.$
\end{proof}

\subsection{References}
Most of the proofs presented in this chapter were taken from Benedetti -- Petronio \cite{BP}, and are also contained in Thurson's notes \cite{Th}. The proof of Mostow's rigidity presented here was proposed by Gromov and applies to every dimension $n\geqslant 3$, because ideal regular simplexes have indeed the maximum volume by Haagerup and Munkholm \cite{HM}. 

%% file: Three.tex
\chapter{Hyperbolic three-manifolds} \label{three:chapter}

We have studied and classified the three-manifolds having seven of the eight geometries, and we are now left with the most interesting ones: hyperbolic three-manifolds.

In dimension two every closed hyperbolic surface is constructed by gluing some geodesic pair-of-pants. In dimension three, although closed hyperbolic three-manifolds are everywhere, it is somehow harder to construct them explicitly: the most general procedure to determine a hyperbolic metric (if any) on a given closed 3-manifold consists of solving \emph{Thurston's equations}.

Thurston's equations arise naturally in the attempt of constructing a hyperbolic three-manifold by triangulating it into hyperbolic tetrahedra. 
The most relevant and unexpected aspect of the theory is that it is much easier to employ hyperbolic \emph{ideal} tetrahedra than compact ones. The combinatorial framework in the ideal case is so convenient, that we use it also for closed three-manifolds.

In this chapter we describe these equations and use them to determine various finite-volume hyperbolic 3-manifolds. We start with the cusped case and then turn to the slightly more complicate closed one.

\section{Cusped three-manifolds} \label{cusped:section}
We know that every cusped finite-volume complete hyperbolic surface is constructed by gluing isometrically finitely many ideal triangles along their edges (see Proposition \ref{straightens:prop}). Likewise, we now construct plenty of cusped hyperbolic three-manifolds by gluing isometrically finitely many ideal tetrahedra along their faces.

\subsection{Ideal tetrahedra}
Ideal triangles are all isometric, but ideal tetrahedra are not! We now show that they can be described up to isometry by a single complex parameter $z$ with $\Im z >0$.\index{ideal tetrahedron}

\begin{figure}
\begin{center}
\includegraphics[width = 7cm] {\iftoggle{BW}{tetrahedron-BW}{tetrahedron}} 
\nota{An ideal tetrahedron with three vertices in $0,1,\infty$ in the half-space mode is determined by the position $z\in \matC\cup\{\infty\}$ of the fourth vertex. A small horosphere centred at the ideal vertex intersects the tetrahedron in a Euclidean triangle uniquely determined up to similarities.}
\label{tetrahedron:fig}
\end{center}
\end{figure}

An ideal tetrahedron is determined by its four ideal vertices $v_1, v_2, v_3, v_4\in\partial \matH^3$. 
We use the half-space model $H^3$ and recall that $\partial H^3 = \matC \cup \{\infty\}$ and $\Iso^+(H^3) = \PSLC$, hence there is a unique orientation-preserving isometry of $H^3$ that sends the vertices $v_1, v_2, v_3, v_4$ respectively to $0, 1, \infty, z$ for some $z$. Up to mirroring with the orientation-reversing reflection $z\mapsto \bar z$ we can suppose that $\Im z > 0$. 

\begin{oss}
By definition the number $z$ is the \emph{cross-ratio} of the four complex numbers $v_1, v_2, v_3, v_4$.\index{cross-ratio}
\end{oss}

\begin{figure}
\begin{center}
\includegraphics[width = 9 cm] {\iftoggle{BW}{tetraedro2-BW}{tetraedro2}} 
\nota{At each ideal vertex we have a Euclidean triangle defined up to similarities: each vertex of the triangle has a well-defined complex angle (left). We can assign the complex angles directly to the edges of the tetrahedron. The argument is the dihedral angle of the edge (right).}
\label{tetraedro2:fig}
\end{center}
\end{figure}

A horosphere centred at the vertex $v_3=\infty$ is a horizontal Euclidean plane that intersects the ideal tetrahedron in a Euclidean triangle as in Figure \ref{tetrahedron:fig}. The oriented similarity class of the triangle depends only on the vertex $v_3$, because a horosphere change results in a dilation: we can represent it as a triangle in $\matC=\matR^2$ with vertices at $0$, $1$, and $z$ as in Figure \ref{tetraedro2:fig}-(left). The \emph{complex angle} of a vertex of the triangle is the ratio of the two adjacent sides, taken with clockwise order and seen as complex numbers. The three complex angles are shown in Figure \ref{tetraedro2:fig}-(left) and are:
$$z,\ \frac 1{1-z},\ \frac{z-1}z.$$
The argument is the usual angle, and the modulus is the ratio of the two lengths of the adjacent sides.

Proposition \ref{symmetries:ideal:tetrahedron:prop} shows that some symmetries of an ideal tetrahedron act on its vertices like the alternating group $A_4$ and hence transitively. Therefore every vertex has the same triangular section as in Figure \ref{tetraedro2:fig}-(left), and all the sections can be encoded by assigning the complex numbers directly to the edges of $\Delta$ as shown in Figure \ref{tetraedro2:fig}-(right). These labels on the edges are called \emph{moduli} and determine the ideal tetrahedron up to orientation-preserving isometries of $\matH^3$. The argument of a modulus is the dihedral angle of the edge.

On a regular ideal tetrahedron the triangular sections are equilateral and hence all edges have the same modulus $z = e^{\frac{\pi i}3}$.

\subsection{Ideal triangulations} \label{ideal:subsection}
Let $\Delta_1,\ldots, \Delta_n$ be identical copies of the standard oriented 3-simplex. As in the two-dimensional case, a \emph{triangulation} $\calT$ is a partition of the $4n$ faces of the tetrahedra into $2n$ pairs, and for each pair a simplicial isometry between the two faces. The triangulation is \emph{oriented} if the simplicial isometries are orientation-reversing. If we glue the tetrahedra along the simplicial isometries we get a topological space $X$, which is not necessarily a topological manifold. Let $M$ be $X$ minus the vertices of the triangulation: we say that $\calT$ is an \emph{ideal triangulation} for $M$.

\begin{figure}
\begin{center}
\includegraphics[width = 4 cm] {\iftoggle{BW}{Truncatedtetrahedron-BW}{269px-Truncatedtetrahedron}}
\nota{A truncated tetrahedron.}
\label{truncated:fig}
\end{center}
\end{figure}

\begin{prop} If $\calT$ is oriented then $M$ is a topological 3-manifold, homeomorphic to the interior of a compact oriented manifold with boundary.
\end{prop}
\begin{proof}
To prove that $M$ is a manifold we only need to check that a point $x\in e$ in an edge $e$ has a neighbourhood homeomorphic to an open ball. A cycle of tetrahedra is attached to $e$, and since $\calT$ is oriented we are certain that a neighbourhood of $x$ is a cone over a 2-sphere and not over a projective plane.

If we truncate the tetrahedra as in Figure \ref{truncated:fig} before gluing them, we get a compact manifold $N\subset M$ with boundary such that $M\setminus N \isom \partial N \times [0,1)$. Therefore $M$ is homeomorphic to $\interior N$.
\end{proof}

We will always suppose that $\calT$ is oriented and $M$ is connected, so $M=\interior N$ for some compact $N$ with boundary. Every ideal vertex $v$ in $\calT$ is locally a cone over a small triangulated closed surface $\Sigma\subset M$ obtained by truncating the tetrahedra incident to $v$. We call such a $\Sigma$ the \emph{link} of $v$. 

\subsection{Hyperbolic ideal triangulations}
Let $\calT$ be an oriented ideal triangulation with tetrahedra $\Delta_1,\ldots, \Delta_n$ of a 3-manifold $M$. We now substitute every $\Delta_i$ with an ideal hyperbolic tetrahedron and pair their faces with orientation-reversing isometries. As opposite to the two-dimensional case, the ideal hyperbolic tetrahedra are not unique (they depend on a complex modulus $z_i$) but the isometric pairing of their faces is uniquely determined in virtue of the following.

\begin{prop} Given two ideal triangles $\Delta$ and $\Delta'$, every bijection between the ideal vertices of $\Delta$ and of $\Delta'$ is realised by a unique isometry. 
\end{prop}
\begin{proof}
We see the ideal triangles in $\matH^2$ and recall that for any two triples of points in $\partial \matH^2$ there is a unique isometry sending pointwise the first triple to the second. Alternatively, we may use the barycentric decomposition shown in Figure \ref{triangoli_ideali:fig}-(left).
\end{proof}

\begin{figure}
\begin{center}
\includegraphics[width = 11 cm] {\iftoggle{BW}{tetrahedron3-BW}{tetrahedron3}} 
\nota{If we manage to glue all the tetrahedra incident to an edge $e$ inside $\matH^3$ as shown, the hyperbolic structure is defined also in $e$ (left). Let $z_1,\ldots, z_h$ be the complex moduli assigned to the sides of the $h$ incident tetrahedra (here $h=5$). We can glue everything in $\matH^3$ if and only if $z_1\cdots z_h=1$ and the arguments sum to $2\pi$ (right).}
\label{tetrahedron3:fig}
\end{center}
\end{figure}

If we substitute each $\Delta_i$ with an ideal hyperbolic tetrahedron, we immediately get a well-defined hyperbolic structure on $M$ minus the edges of $\calT$.
We now try to extend this hyperbolic structure to the edges: we can do this if we are able to glue all the $h$ tetrahedra around each edge $e$ inside $\matH^3$ as in Figure \ref{tetrahedron3:fig}. Let $z_1,\ldots, z_h$ be the complex moduli associated to the edges of the $h$ tetrahedra incident to $e$. As shown in the figure, if $z_1\cdots z_h = 1$ and the sum of their argument is $2\pi$ (and not some higher multiple of $2\pi$) then all tetrahedra can be glued simultaneously inside $\matH^3$ and the hyperbolic structure extends naturally to $e$. If this holds at every edge $e$ of $\calT$ then  $M$ inherits a hyperbolic structure and $\calT$ is said to be a \emph{hyperbolic} (or \emph{geometric}) \emph{ideal triangulation} for $M$.\index{hyperbolic ideal triangulation}

\subsection{Consistency equations}
We want to parametrize the hyperbolic structures on $M$ that may be constructed in this way from a fixed ideal triangulation $\calT$. We pick an arbitrary edge for every tetrahedron $\Delta_i$ and assign to it a complex variable $z_i$ with $\Im z_i>0$, and the other edges of $\Delta_i$ are automatically labeled by one of the variables $z_i$, $\frac{z_i-1}{z_i}$, or $\frac 1{1-z_i}$ as indicated in Figure \ref{tetraedro2:fig} (recall that the tetrahedra are oriented). As we have seen, for every edge $e$ in $\calT$ we obtain an equation of type 
$$w_1\cdots w_h = 1$$ 
(to which we must add the condition that the sum of the arguments is $2\pi$), where each $w_j$ equals $z_i$, $\frac{z_i-1}{z_i}$, or $\frac 1{1-z_i}$ for some $i$.

We have thus obtained a system of \emph{consistency equations}, with a variable $z_i$ for each tetrahedron and an equation $w_1\cdots w_h = 1$ for each edge. A solution $z = (z_1,\ldots, z_n)$ to these equations produces a hyperbolic ideal triangulation and hence a hyperbolic structure on $M$. Recall that we assume $\Im z_i>0$ $\ \forall i$.

As in the two-dimensional case, the resulting hyperbolic structure is not necessarily complete, and to get a complete hyperbolic manifold we must add more equations.

\subsection{Completeness equations} \label{complete:solutions:subsection}
Our aim is to construct a complete finite-volume hyperbolic metric on $M$. By Corollary \ref{tame:cor} if $M$ has such a metric the link of every ideal vertex of $\calT$ is a triangulated torus and identifies a cusp of $M$, so we will henceforth suppose that the links of all vertices are tori. In other words $M$ is the interior of a compact 3-manifold $N$ bounded by some $c>0$ tori.

\begin{figure}
\begin{center}
\includegraphics[width = 11 cm] {\iftoggle{BW}{complete_equations-BW}{complete_equations}} 
\nota{Every boundary torus $T\subset \partial N$ is triangulated, and each triangle has a Euclidean structure well-defined up to similarities and inherits three complex moduli $w_1, w_2, w_3$ at its vertices (left). At every vertex of the triangulation, the product of the adjacent moduli is $1$, for instance here $w_{12}w_{15}w_{16}w_{19}w_{23} = 1$. The \iftoggle{BW}{grey}{red} path contributes to $\mu(\gamma)$ with the factor $-w_{30}w_{29}w_{25}w_{24}w_{23}w_{19}w_{20}$ (right).}
\label{complete_equations:fig}
\end{center}
\end{figure}

Let $z=(z_1,\ldots, z_n)$ be a solution to the consistency equations, providing a hyperbolic structure to $M$. Every boundary torus $T\subset \partial N$ is triangulated by $\calT$: every triangle in $T$ is the truncation triangle of some $\Delta_i$ and hence inherits the complex moduli of the three adjacent edges of $\Delta_i$ as in Figure \ref{complete_equations:fig}-(left), thus it has a Euclidean structure well-defined up to similarities. We want to study the following problem: do these Euclidean structures on the triangles glue to form a Euclidean structure on $T$? We will see that $M$ is complete if and only if the answer is ``yes'' at each boundary torus $T$.

Pick $\gamma \in \pi_1(T)=H_1(T,\matZ) \isom \matZ^2$. We represent $\gamma$ as a simplicial path in the triangulation of $T$ and then define $\mu(\gamma)\in \matC^*$ to be $(-1)^{|\gamma|}$ times the product of all the complex moduli that $\gamma$ encounters at its right side, with $|\gamma|$ being the number of edges of $\gamma$, see Figure \ref{complete_equations:fig}-(right). 

\begin{figure}
\begin{center}
\includegraphics[width = 11 cm] {\iftoggle{BW}{complete_move-BW}{complete_move}} 
\nota{This move for $\gamma$ does not affect $\mu(\gamma)$.}
\label{complete_move:fig}
\end{center}
\end{figure}

\begin{prop} \label{mu:prop}
The element $\mu(\gamma)$ is well-defined and $\mu\colon\pi_1(T) \to \matC^*$ is a homomorphism.
\end{prop}
\begin{proof}
Two different paths for $\gamma$ are related by moves as in Figure \ref{complete_move:fig}. This move does not affect $\mu(\gamma)$ since $w_1w_2w_3=-1$ and the product of the moduli around a vertex is $+1$. The map $\mu$ is clearly a homomorphism.
\end{proof}

Let $C(T)\subset M$ be a closed collar of the torus $T$ in $N$, intersected with $M$. It is diffeomorphic to $T\times [0,+\infty)$.

\begin{prop} \label{mu:equivalent:prop}
The following facts are equivalent:
\begin{enumerate}
\item the homomorphism $\mu$ is trivial,
\item there is a Euclidean structure on $T$ that induces all the moduli,
\item the manifold $C(T)$ is complete and contains a truncated cusp.
\end{enumerate}
\end{prop}
\begin{proof}
The equivalence (1)$\Leftrightarrow$(2) is a simple exercise.
If (2) holds we may choose small horosections of all the hyperbolic tetrahedra incident to $v$ that match to give a Euclidean torus $T_*\subset C(T)$. The non-compact part of $C(T)$ bounded by $T_*$ is a truncated cusp with base $T_*$ and is thus complete.
On the other hand, if (1) does not hold, there is a path of horosections from a tetrahedron to itself which ends at a bigger height and we conclude that $C(T)$ is not complete as in the two-dimensional case, see the proof of Proposition \ref{completion:surface:prop}.
\end{proof}

\begin{cor} The hyperbolic manifold $M$ is complete $\Longleftrightarrow$ $\mu$ is trivial for every torus $T\subset \partial N$.
\end{cor}

Fix for every boundary torus $T$ two generators $m,l$ for $H_1(T, \matZ)$. The homomorphism $\mu$ is trivial $\Longleftrightarrow$ the following two equations are satisfied:
$$\mu(m) = 1, \quad \mu(l)=1.$$ 
Each equation is of some type $w_1\cdots w_k = 1$.
We get two equations for each of the $c$ boundary tori and hence $2c$ equations in total, called the \emph{completeness equations} for the triangulation $\calT$. We summarise our discussion:\index{Thurson's hyperbolicity equations}

\begin{prop}
Let $\calT$ be an ideal triangulation of $M=\interior N$ with $n$ tetrahedra and $\partial N$ consisting of $c$ tori. If a point $z = (z_1,\ldots, z_n)$ with $\Im z_i >0$ satisfies the $n$ consistency equations and the $2c$ completeness equations, then $M$ admits a finite-volume complete hyperbolic metric.
\end{prop}

\subsection{Examples}
On a triangulation $\calT$, the \emph{valence} of an edge $e$ is the number of tetrahedra incident to it, counted with multiplicity. The simplest kinds of solutions arise in the following construction. 

\begin{prop} \label{valence:6:prop}
Let $\calT$ be a triangulation where all edges have valence six. The point $(e^{\frac{\pi i}3}, \ldots, e^{\frac{\pi i}3})$ is a solution of both the consistency and completeness equations and defines a hyperbolic structure where all the ideal tetrahedra are regular.
\end{prop}
\begin{proof}
Note that when $z= e^{\frac{\pi i}3}$ the tetrahedron is regular and we get 
$$z = \frac{z-1}z = \frac 1{1-z} = e^{\frac{\pi i}3}.$$
Therefore the moduli are $e^{\frac{\pi i}3}$ everywhere. The consistency equations are satisfied since $(e^{\frac{\pi i}3})^6 = 1$ and also the completeness equations are, because by assigning length 1 to each edge of the triangulation of a torus $T$ we get a global Euclidean structure on $T$ tessellated into isometric equilateral triangles.
\end{proof}

\begin{figure}
\begin{center}
\includegraphics[width = 10 cm] {\iftoggle{BW}{triangulation_figure8-BW}{triangulation_figure8}} 
\nota{An oriented triangulation with $2$ tetrahedra: faces with the same letter $F, J, P, R$ are paired with a simplicial map that matches the letters. We get two edges with valence six (the \iftoggle{BW}{grey}{red} and white dots). }
\label{triangulation_figure8:fig}
\end{center}
\end{figure}

Figure \ref{triangulation_figure8:fig} shows an ideal triangulation with two tetrahedra and two edges, each with valence six: this defines a cusped finite-volume complete hyperbolic three-manifold $M$ that decomposes into two regular ideal tetrahedra. 

Every finite covering $\tilde M$ of $M$ is another example: the ideal triangulation $\calT$ lifts to an ideal triangulation $\tilde \calT$ where all the edges have valence six. Since the fundamental group $\pi_1(M)$ is residually finite (see Proposition \ref{RF:hyperbolic:prop}) there are plenty of such coverings.

\begin{figure}
\begin{center}
\includegraphics[width = 10 cm] {\iftoggle{BW}{octa-BW}{octa}} 
\nota{By pairing the faces of two ideal regular octahedra as shown we get a cusped complete hyperbolic three-manifold. Faces with the same letter are paired with the unique isometry that matches the letters. We get six edges with valence four (marked with coloured dots). }
\label{octa:fig}
\end{center}
\end{figure}

\begin{ex} Construct a triangulation $\calT$ with one tetrahedron and one edge. The edge has valence six and hence the construction of Proposition \ref{valence:6:prop} produces a cusped hyperbolic manifold. Note however that the triangulation is not orientable and a cusp section is a Klein bottle! This non-orientable hyperbolic manifold is called the \emph{Gieseking manifold}.\index{Gieseking manifold}
\end{ex}

Similar kinds of ``regular'' examples may be constructed using some other platonic solids.
The ideal regular octahedron is particularly interesting and useful because it is right-angled, see Section \ref{platonici:subsection}. If we pick finitely many oriented regular ideal octahedra $O_1, \ldots, O_h$ and pair orientation-reversingly their triangular faces in a way that every resulting edge has valence four, we construct a finite-volume cusped hyperbolic manifold $M$ because $4\times \frac \pi 2 = 2 \pi$. An example with two octahedra is shown in Figure \ref{octa:fig}.

We can further decompose each regular ideal octahedron into four non-regular ideal tetrahedra (there are three ways to do this: you must choose a diagonal connecting two opposite vertices) to get an ideal triangulation $\calT$ for $M$. The moduli of each such tetrahedron are 
$$i, \ \frac 1{1-i} = \frac{1+i}2, \ \frac{i-1}i = 1+i$$
and they satisfy the consistency and completeness equations in a less trivial way than before. The edges of the resulting triangulation have varying valences 4, 6, and 8.

\subsection{SnapPea}
It is of course hard to construct ideal triangulations by hand, except in some very symmetric and simple cases like the one just described. And it is even more difficult to solve the consistency and completeness equations.

There is a beautiful computer program that does all this for you! This is \emph{SnapPea}, written by Jeff Weeks in the 1980s. Using SnapPea you can draw any link diagram $L$ and the program immediately constructs an ideal triangulation for the complement of $L$ in $S^3$. Then it uses the Newton method to find a numerical solution to the consistency and completeness equations: the solution is only numerical, but then one can use some \emph{a posteriori} argument to confirm it rigorously.

For instance, if we draw the trefoil knot, SnapPea triangulates its complement and finds no solution to the equations: this is not surprising since its complement is Seifert and hence not hyperbolic, see Proposition \ref{torus:knot:complement:prop}. 

If we draw the figure-eight knot in Figure \ref{figure8:fig}-(left), SnapPea constructs precisely the ideal triangulation of Figure \ref{triangulation_figure8:fig}, and we have hence discovered that the complement of the figure-eight knot is hyperbolic and decomposes into two ideal regular tetrahedra. 

We may draw links with many components, and discover for instance that the complement of the Borromean rings from Figure \ref{figure8:fig}-(right) is hyperbolic and decomposes into two right-angled ideal regular octahedra. We can also Dehn-fill some boundary components to obtain more cusped manifolds that are not link complements in $S^3$.

\begin{figure}
\begin{center}
\includegraphics[width = 4.5 cm] {\iftoggle{BW}{Blue_Figure-Eight_Knot-BW}{422px-Blue_Figure-Eight_Knot}} 
\includegraphics[width = 5 cm] {\iftoggle{BW}{Borromean_Rings_Illusion-BW}{484px-Borromean_Rings_Illusion}}
\nota{The figure-eight knot and the Borromean rings are \emph{hyperbolic}, that is  their complements admit complete hyperbolic metrics.}
\label{figure8:fig}
\end{center}
\end{figure}

SnapPea can calculate numerically a wealth of geometrical invariants with some precision, including the volume and the first segment of the geodesic spectrum. It can  manipulate manifolds (Dehn-filling, drilling along short simple closed geodesics, finite covers) and change triangulations. It provides beautiful pictures of the Dirichlet domain and of the cusp shapes. 

Is SnapPea guaranteed to find a hyperbolic structure on $M$ if there is one? No, it is not, although in practice it succeeds most of the time. The fact that $M$ is hyperbolic does not guarantee that there is a solution to Thurston's equations on a given triangulation $\calT$, and the existence of a solution does not guarantee that SnapPea will be able to find it. 

We note that for most ideal triangulations $\calT$ of a hyperbolic $M$ there is no solution to Thurston's equations, and SnapPea is clever enough to modify $\calT$ to increase the probability to find one. In all known examples, a cusped hyperbolic $M$ admits at least \emph{one} geometric ideal triangulation, but whether this holds true for all cusped hyperbolic 3-manifolds $M$ is still an open question.

\subsection{Cusped census}
As with knots tabulations, topologists have used computers to list the cusped hyperbolic manifolds that can be triangulated with few ideal tetrahedra.

\begin{table}
\begin{center}
\begin{tabular}{c||ccccccc}
tetrahedra        & 1 & 2 & 3 & 4 & 5 & 6 & 7 \\
 \hline \hline
1 cusp &  & 2 & 9 & 52 & 223 & 913 & 3388 \\
2 cusps & & & & 4 & 11 & 48 & 162  \\
3 cusps & & & & & & 1 & 2
\end{tabular}
\vspace{.2 cm}
\nota{The number of cusped orientable hyperbolic three-manifolds that can be triangulated with at most 7 ideal tetrahedra.}
\label{cusped:census:table}
\end{center}
\end{table}

\begin{table}
\begin{center}
\begin{tabular}{c||ccccccc}
Name        & Volume & Homology & Symmetry & SG & C \\
 \hline \hline
$M2_1$ & $2.0298832128$ & $\matZ$ & $D_4$ & $1.09$ & a \\
$M2_2$ & $2.0298832128$ & $\matZ+\matZ_5$ & $\matZ_2+\matZ_4$ & $0.86$ & a \\
$M3_1$ & $2.5689706009$ & $\matZ+\matZ_5$ & $D_4$ & $0.65$ & c \\
$M3_2$ & $2.5689706009$ & $\matZ+\matZ_3$ & $D_4$ & $0.65$ & c \\
$M3_3$ & $2.6667447834$ & $\matZ+\matZ_2$ & $D_4$ & $0.63$ & c \\
$M3_4$ & $2.6667447834$ & $\matZ+\matZ_6$ & $D_4$ & $0.63$ & c \\
$M3_5$ & $2.7818339124$ & $\matZ$ & $\matZ_2$ & $0.51$ & c \\
$M3_6$ & $2.8281220883$ & $\matZ$ & $D_4$ & $0.56$ & c \\
$M3_7$ & $2.8281220883$ & $\matZ$ & $\matZ_2$ & $0.58$ & c \\
$M3_8$ & $2.8281220883$ & $\matZ+\matZ_7$ & $D_4$ & $0.56$ & c \\
$M3_9$ & $2.9441064867$ & $\matZ$ & $\matZ_2$ & $0.43$ & c 
\end{tabular}
\vspace{.2 cm}
\nota{The 1-cusped orientable hyperbolic manifolds that can be triangulated with at most 3 ideal tetrahedra. The name $Mi_j$ indicates that the manifold is the one with $j$-th smallest volume among those that can be triangulated with $i$ tetrahedra. The column SG shows the length of the shortest closed geodesic (volume and SG values are truncated after few digits). The column C indicates whether the manifold is achiral (a) or chiral (c), that is if it admits an orientation-reversing isometry or not.}
\label{cusped:census:2:table}
\end{center}
\end{table}

\begin{table}
\begin{center}
\begin{tabular}{c||ccccccc}
Name        & Volume & Homology & Symmetry & SG & C \\
 \hline \hline
$M4_1^2$ & $3.6638623767$ & $\matZ+\matZ$ & $D_8$ & $1.06$ & c \\
$M4_2^2$ & $3.6638623767$ & $\matZ+\matZ$ & $D_8$ & $0.96$ & c \\
$M4_3^2$ & $4.0597664256$ & $\matZ+\matZ$ & $D_{12}$ & $0.86$ & c \\
$M4_4^2$ & $4.0597664256$ & $\matZ+\matZ$ & $D_8$ & $0.86$ & a \\
\end{tabular}
\vspace{.2 cm}
\nota{The 2-cusped orientable hyperbolic manifolds that can be triangulated with 4 ideal tetrahedra. The name $Mi_j^k$ indicates that the manifold is the $k$-cusped one with $j$-th smallest volume among those that triangulate with $i$ tetrahedra.}
\label{cusped:census:3:table}
\end{center}
\end{table}

The number of cusped hyperbolic orientable manifolds that can be ideally triangulated with $n$ tetrahedra (and not less than $n$) is written in Table \ref{cusped:census:table} for all $n\leqslant 7$. The 1-cusped manifolds with $n=2,3$ are listed with more detail in Table \ref{cusped:census:2:table} and the 2-cusped ones with $n=4$ are in Table \ref{cusped:census:3:table}. These tables were produced by Callahan -- Hildebrand -- Weeks \cite{CaHiWe} in 1999 via a computer enumeration.

\begin{figure}
\begin{center}
\includegraphics[width = 12.5 cm] {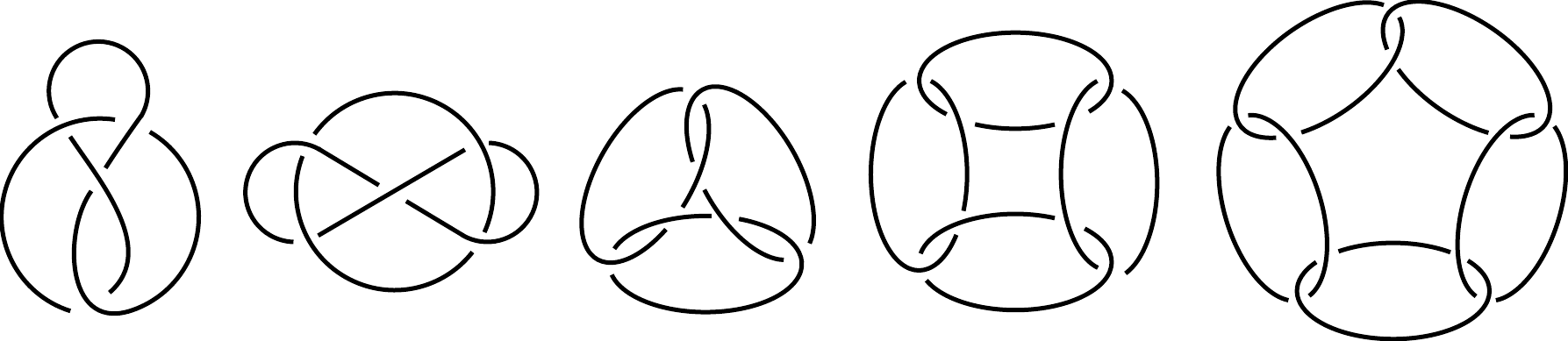}
\nota{A remarkable sequence of hyperbolic links. These are the \emph{figure eight knot}, the \emph{Whitehead link}, and some particular \emph{chain links} with 3,4,5 components. }
\label{sequence:fig}
\end{center}
\end{figure}

The figure-eight knot complement is $M2_1$, while $M2_2$ is another hyperbolic manifold that decomposes into two regular ideal tetrahedra: it is not a knot complement in $S^3$ since its homology is not $\matZ$, but it is  the complement of a knot $K$ in the lens space $L(5,1)$, called the \emph{figure-eight knot sibling}. This knot is obtained by performing a $(-5)$-surgery on one component of the \emph{Whitehead link} shown in Figure \ref{sequence:fig}: since this component is trivial, the surgered manifold is $L(5,1)$ and the other unsurgered component is our knot $K$. Both manifolds $M2_1$ and $M2_2$ have volume $2.0298832128\ldots$ since the volume of the ideal regular tetrahedron is $3\Lambda\big(\frac \pi 3\big) = 1.0149416064\ldots$

The manifold $M4_1^2$ is the Whitehead link complement, that can be obtained by pairing the faces of a single ideal regular octahedron $O$. Indeed the volume of $O$ is $8\Lambda\big(\frac {\pi} 4\big) = 3.6638623767\ldots $

After the figure eight knot and the Whitehead link, a notable sequence of hyperbolic links in $S^3$ with increasing number of components is shown in Figure \ref{sequence:fig}. Their complements are conjectured to be the smallest hyperbolic manifolds with $i=1,\ldots, 5$ cusps (see the references at the end of the chapter).

\subsection{Hyperbolic knots}
How many hyperbolic knots in $S^3$ are there? There are infinitely many hyperbolic knots, and infinitely many non-hyperbolic knots. Geometrisation translates hyperbolicity into an appealing topological condition:

\begin{teo}
A knot $K\subset S^3$ is either a torus knot, a satellite knot, of a hyperbolic knot.
\end{teo}
\begin{proof}[Proof using geometrisation]
Proposition \ref{trichotomy:knots:prop} says that $K$ is either a torus knot, a satellite knot, or has simple complement $M$. In the latter case $M$ is hyperbolic by Corollary \ref{simple:bounded:hyperbolic}.
\end{proof}

Note that the three cases are mutually exclusive. The number of prime knots with $c\leqslant 14$ crossings in each of the three classes is shown in Table \ref{torus:satellite:hyperbolic:knots:table}, taken from Hoste -- Thistlethwaite -- Weeks \cite{HTW}: the table shows a strong predominance of hyperbolic knots, but it is unknown whether it persists also when $c$ tends to infinity.

\begin{table}
\begin{center}
\begin{tabular}{c||cccccccccccc}
 $c$       &  
 $3$ & $4$ & $5$ & $6$ & $7$ & $8$
 & $9$ & $10$ & $11$ & $12$ & $13$ & $14$ \\
 \hline \hline
torus & 
1 & 0 & 1 & 0 & 1 & 1 & 1 & 1 & 1 & 0 & 1 & 1 \\
satellite & 0 & 0 & 0 & 0 & 0 & 0 & 0 & 0 & 0 & 0 & 2 & 2 \\
hyperbolic & 
0 & 1 & 1 & 3 & 6 & 20 & 48 & 164 & 551 & 2176 & 9985 & 46969 
 \end{tabular}
\vspace{.2 cm}
\nota{The number of torus, satellite, and hyperbolic prime knots with $c$ crossings, for all $c\leqslant 14$.}
\label{torus:satellite:hyperbolic:knots:table}
\end{center}
\end{table}

\section{Closed hyperbolic three-manifolds}
We have described in Section \ref{cusped:section} a method to construct cusped finite-volume hyperbolic three-manifolds. How can we now build some \emph{closed} hyperbolic three-manifolds? 

Following the same path, we could try to parametrize \emph{compact} hyperbolic tetrahedra via some variables, and then encode the isometric gluings of their faces via some equations. However, it is really hard to parametrize compact tetrahedra and their gluings, and nobody has ever constructed any closed hyperbolic three-manifold in this way! 

The usual procedure for building closed hyperbolic three-manifolds consists of using (again) ideal tetrahedra: we start with a cusped manifold, and then we slightly modify the completeness equations to allow some appropriate Dehn-filling of some (or all) cusps. 

\subsection{Dehn filling parameters} \label{DF:parameters:subsection}
We now consider an orientable three-manifold $M = \interior N$ where $N$ is compact and $\partial N = T_1 \sqcup \ldots \sqcup T_c$ consists of tori. For simplicity we will always assume that each $T_i$ is incompressible, so that $\pi_1(T_i)=\matZ\times \matZ$ injects in $\pi_1(N) = \pi_1(M)$ for all $i$. 

We fix once for all two generators $m_i, l_i$ of $\pi_1(T_i)$ for each $i=1,\ldots,c$. A \emph{Dehn filling parameter} $s = (s_1,\ldots, s_c)$ is a sequence where each $s_i$ is either a pair $(p,q)$ of coprime integers or the symbol $\infty$. It is useful to think of $s_i$ as a point in the two-sphere $S^2 = \matR^2 \cup \{\infty \}$.

A Dehn filling parameter $s$ determines a Dehn filling of $N$ as follows: for every $i$, if $s_i=(p,q)$ we fill $T_i$ by killing the slope $pm_i + ql_i$, while if $s_i=\infty$ we do nothing. The result is a new compact manifold, whose interior we denote by $M^{\rm fill}$, that may be closed (if there are no $\infty$ in $s$) or may be bounded by some tori.

We now fix a Dehn filling parameter $s$, producing a filled manifold $M^{\rm fill}$.

\subsection{The modified equations} \label{modified:equations:subsection}
Let $\calT$ be an oriented ideal triangulation for $M$. We now describe some equations similar to that of Section \ref{cusped:section}, whose solutions now identify hyperbolic structures on $M^{\rm fill}$ instead of $M$. The equations of course will depend on the Dehn filling parameter $s$, since the manifold $M^{\rm fill}$ does.

The triangulation $\calT$ consists of some $n$ tetrahedra $\Delta_1,\ldots, \Delta_n$; we choose an edge in each $\Delta_i$ and we assign it the variable $z_i\in \matC$ with $\Im z_i >0$ as we did in Section \ref{cusped:section}.

The edges of all tetrahedra are now coloured with moduli $z_i$, $\frac 1{1-z_i}$, $\frac{z_i-1}{z_i}$ that lie in the upper half-plane and as such may be written in polar coordinates as $\rho e^{i\theta}$ with $\theta \in (0,\pi)$. It is convenient to consider these variables as elements of the multiplicative group
$$\widetilde{\matC^*} = \big\{\rho e^{i\theta}\ \big| \ \rho \in \matR_{>0},\ \theta \in \matR\big\}$$
that covers $\matC^*$ and is isomorphic to $\matC$ via the map $\exp\colon \matC \to \widetilde{\matC^*}$. 

Every edge of $\calT$ furnishes a \emph{consistency equation} in $\widetilde {\matC^*}$ as described in Section \ref{cusped:section}, of the form
$$w_1\cdots w_h = e^{2 \pi i}.$$
Note that $e^{2 \pi i} \neq 1$ in $\widetilde{\matC^*}$, so by interpreting the consistency equations in $\widetilde{\matC^*}$ we have also incorporated the request that all angles sum to $2\pi$.

Let $z=(z_1,\ldots, z_n)$ be a solution to the consistency equations.
We have defined in Proposition \ref{mu:prop} a homomorphism $\mu\colon \pi_1(T_i) \to \matC^*$ for every boundary torus $T_i$, and we now lift it to $\widetilde{\matC^*}$. 

\begin{defn} \label{tilde:mu:defn}
Represent each non-trivial $\gamma \in \pi_1(T_i)$ as a simplicial path that lifts to an embedded path in the universal cover of $T_i$. Define $\tilde\mu(\gamma)$ to be $e^{-|\gamma| \pi i }$ times the product of all the moduli that $\gamma$ encounters at its right side, considered now as elements in the group $\widetilde{\matC^*}$. If $\gamma$ is trivial, set $\tilde\mu(\gamma) = 1$.
\end{defn}

\begin{ex} The element $\tilde\mu(\gamma)\in\widetilde{\matC^*}$ is well-defined and
$$\tilde\mu\colon \pi_1(T_i) \to \widetilde{\matC^*}$$
is a homomorphism.
\end{ex}

We now define some new \emph{completeness equations} at $T_i$ for each $i=1,\ldots, c$ that depend on the Dehn filling parameter $s_i$. There are two cases to consider for each $i=1,\ldots, c$. If $s_i = \infty$, the completeness equations relative to $T_i$ are the two equations 
$$\mu(m_i) = \mu(l_i) = 1$$ 
already considered in Section \ref{complete:solutions:subsection}. (If we wish, we can substitute them with the equations $\tilde\mu(m_i) = \tilde\mu(l_i) = 1$.
Both choices will work.) 

If $s_i = (p,q)$ we define one new completeness equation
$$\tilde\mu(m_i)^p \cdot \tilde\mu(l_i)^q = e^{2\pi i}.$$
The total number of completeness equations therefore varies from $c$ to $2c$, depending on the number of cusps that are left unfilled. The main goal of this section is to prove the following.

\begin{teo} \label{modified:equations:teo}
A solution $z = (z_1,\ldots, z_n)$ to the consistency and completeness equations determines a hyperbolic structure on $M$ whose completion $\overline M$ is a complete hyperbolic manifold diffeomorphic to $M^{\rm fill}$.
\end{teo}

If $s\neq (\infty, \ldots, \infty)$ the solution $z$ determines an incomplete hyperbolic metric on $M$, and the miracle here is that the completion $\overline M$ is another hyperbolic manifold! 

We note that both the consistency and completeness equations are of type $w_1\cdots w_h = 1$ or $e^{2\pi i}$.
The rest of this section is mainly devoted to the proof of Theorem \ref{modified:equations:teo}. We start with a short discussion that introduces (truncated) solid tori with cone angles: these objects will be crucial in the proof of the theorem.

\subsection{The infinite branched covering}
Let $l\subset \matH^3$ be any line and recall the universal cover
$$\pi \colon X \longrightarrow \matH^3 \setminus l$$ 
already considered in Section \ref{singularities:subsection} when we introduced manifolds with cone singularities. The manifold $X$ is incomplete and its completion $\overline X$ is obtained by adding a copy $\tilde l$ of $l$.

We use the half-space model and represent $l$ as the vertical axis, so that $\matH^3\setminus l = \matC^* \times \matR_{>0}$ and we can write
$$X = \widetilde{\matC^*} \times \matR_{>0}.$$ 

We note that $\widetilde{\matC^*}$ acts on $X$ via isometries: the element $w\in\widetilde{\matC^*}$ acts as
$$(z,t) \longmapsto (wz, |w| t).$$
When $w = e^{\alpha i}$ this is the rotation of angle $\alpha$ that was used to define hyperbolic cone manifolds. When $|w|\neq 1$ the map projects to a hyperbolic isometry of $\matH^3$ with axis $l$ and translation distance $\log |w|$.

\subsection{Tubes with cone angles}  \label{X:quotients:subsection}
A non-trivial discrete group $\Gamma <\widetilde{\matC^*}\isom \matC$ is either isomorphic to $\matZ$ or to $\matZ\times \matZ$; it acts freely on $X$ and (not necessarily freely!) on the line $\tilde l$. Concerning $\matZ\times \matZ$, there are two cases to consider.

If $\Gamma = \matZ\times \matZ$ contains a non-trivial rotation, it is generated by some maps
$$(z,t) \mapsto (e^{\alpha i}z, t), \quad (z,t) \mapsto e^\lambda(e^{\beta i}z, t)$$
with $\lambda\neq 0$, and it acts on $\tilde l$ as translations whose step is an integer multiple of $\lambda$. The quotient $\overline X/_{\Gamma}$ is naturally a complete hyperbolic manifold with cone angles: it is an open solid torus with singular locus a closed geodesic of length $|\lambda|$ and with cone angle $\alpha$. We call it a \emph{tube with cone angle $\alpha$}. When $\alpha = 2\pi$ we get an ordinary tube as in Section \ref{tubes:subsection}, with no singular points.

If $\Gamma = \matZ \times \matZ$ does not contain a rotation, it acts on $\tilde l$ as an indiscrete group of translations and hence $\overline X/_{\Gamma}$ is not a hyperbolic manifold with cone angles (it is not even Hausdorff). We also note that $\overline X/_{\Gamma}$ is not the completion of the hyperbolic manifold $X/_{\Gamma}$ in that case, but it only maps surjectively onto it: the completion adds a single point to $X/_{\Gamma}$.

Tubes have natural truncations. For every $R>0$ the $R$-neighbourhood $N_R(\tilde l)\subset \overline X$ of the singular line $\tilde l$ projects onto the $R$-neighbourhood $N_R(l)$ of $l$, which is a Euclidean cone with axis $l$. The group $\Gamma$ preserves $N_R(\tilde l)$ and in the first case its quotient is a truncated tube with cone angles.

\subsection{Incomplete solutions}
We now go back to our triangulation $\calT$, and we let $z = (z_1, \ldots, z_n)$ be a solution to the consistency equations for $\calT$. The solution $z$ furnishes a hyperbolic metric on $M$, which might not be complete: we want to understand the metric completion $\overline M$ of $M$. Recall that $M$ is diffeomorphic to the interior of a compact $N$ whose boundary consists of $c$ tori.

We pick a boundary torus $T\subset \partial N$ and define a \emph{collar} $C(T)\subset M$ of $T$ to be the intersection of a closed collar of $T$ in $N$ with $M$. We know from Proposition \ref{mu:equivalent:prop} that when $\mu\colon \pi_1(T) \to \matC^*$ is trivial there is a collar $C(T)$ isometric to a truncated cusp, that is the truncated quotient of $\matH^3$ by a discrete $\matZ \times \matZ$ of parabolic elements.

When $\mu$ is not trivial, we now show that a similar (but different) configuration arises: there is a collar $C(T)$ isometric to a truncated quotient of the hyperbolic manifold $X$ considered above by a discrete $\matZ \times \matZ$. The crucial difference is that $C(T)$ now is incomplete, because $X$ is. To prove this fact we use developing maps and holonomies.

The hyperbolic ideal triangulation $\calT$ of $M$ determined by $z$ lifts to a hyperbolic ideal triangulation $\tilde \calT$ for the universal cover $\tilde M \to M$ with infinitely many tetrahedra. The triangulated torus $T$ is the link of a vertex $v$ of $\calT$. We fix a lift $\tilde v$ of $v$ in $\tilde \calT$, and the link of $\tilde v$ is a triangulated surface $\tilde T$ that covers $T$. Since $T$ is incompressible, $\pi_1(T)$ injects in $\pi_1(M)$ and hence $\tilde T$ is a plane.

Recall from Section \ref{non-complete:subsection} that there is a developing map $D\colon \tilde {M} \to \matH^3$ and a holonomy $\rho\colon \pi_1(M) \to \Iso^+(\matH^3)$, and both are determined once we define $D$ on an ideal tetrahedron in $\tilde\calT$. We choose for our convenience one ideal tetrahedron incident to $\tilde v$ and we map it to the half-space model $H^3$ with an isometry $D$ that sends $\tilde v$ to $\infty$. 

The developing map $D$ sends all tetrahedra of $\tilde\calT$ incident to $\tilde v$ to vertical ideal tetrahedra in $H^3$ with one vertex at $\infty$. If $\gamma \in \pi_1(T)$, then $\rho(\gamma)$ permutes these tetrahedra and fixes $\infty$. Therefore the holonomy $\rho$ sends $\pi_1(T)=\matZ\times \matZ$ to a group of commuting elements in $\PSLC$ fixing $\infty$. Every such element may be written as $z\mapsto az+b$, and there are two possibilities:
\begin{itemize}
\item the group consists of translations $z \mapsto z+b$,
\item the group fixes a point $p\in\matC$.
\end{itemize}
In the second case we may suppose $p=0$ up to translating everything, so the maps are all of type $z\mapsto az$. Not surprisingly, the holonomy $\rho$ is tightly connected with the homomorphism $\mu\colon \pi_1(T) \to \matC^*$.
\begin{ex} \label{rho:derivative:ex}
The derivative $\rho(\gamma)'$, that is the coefficient $1$ or $a$ in the above examples, equals $\mu(\gamma)$. 
\end{ex}
We now suppose that $\mu$ is non-trivial and hence the first possibility is excluded: the group $\rho(\pi_1(T))$ consists of maps $z\mapsto az$ and we can identify $\rho$ with $\mu$. The map $D$ induces a developing map $D\colon \tilde T \to \matC$ with holonomy $\mu$.

\begin{prop}
The image $D(\tilde T)$ misses the origin.
\end{prop}
\begin{proof}
The vector field $v(z) = z$ on $\matC$ is $\mu(\pi_1(T))$-invariant and hence pulls-back via $D$ to a vector field on $T$. If $0\in D(\tilde T)$ then this vector field on $T$ has some zeroes, all with index $1$: a contradiction since $\chi(T)=0$.
\end{proof}

Let $\St_{\tilde v}\subset \tilde M$ be the open star of $\tilde v$, that is the union of all ideal tetrahedra incident to $\tilde v$, with the faces not incident to $\tilde v$ removed: it is homeomorphic to $\tilde T\times [0,+\infty)$ and hence simply connected. The proposition implies that the restriction $D|_{\St_{\tilde v}}$ of
the developing map misses the entire vertical line $l\subset \matH^3$ above the origin and hence lifts to a map 
$$\tilde {D} \colon \St_{\tilde v} \longrightarrow X.$$ 
Likewise the holonomy $\rho$ at $\pi_1(T)$ lifts to a holonomy 
$$\tilde \rho \colon \pi_1(T) \longrightarrow \widetilde{\matC^*} < \Iso^+(X)$$
which is of course related to the homomorphism $\tilde \mu\colon \pi_1(T) \to \widetilde{\matC^*}$.

\begin{ex} \label{tilde:rho:mu:ex}
We have $\tilde \rho = \tilde \mu$.
\end{ex}

After this long discussion, we can finally discover how a collar $C(T)$ of $T$ looks like.

\begin{figure}
\begin{center}
\includegraphics[width = 12.5 cm] {\iftoggle{BW}{Dehn_fill-BW}{Dehn_fill}} 
\nota{For sufficiently small $R>0$ the cone neighbourhood $N_R(\tilde l)$ of $\tilde l$ does not intersect the lower faces of the developed images of $\tilde\Delta_{i_1}, \ldots, \tilde\Delta_{i_h}$ (left). These tetrahedra intersect $\partial N_R(\tilde l)$ into triangles that glue up to determine a torus $T_R$ in $C(T)$ (right).}
\label{Dehn_fill:fig}
\end{center}
\end{figure}

\begin{prop} \label{DF:prop}
If $\mu$ is non-trivial then $\tilde\mu$ is injective and has discrete image. There is a collar $C(T)$ that is isometric to a truncation of $X/_{\img \tilde\mu}$.
\end{prop}
\begin{proof}
We have defined a lifted developing map $\tilde D\colon \St_{\tilde v} \to X$ with holonomy $\tilde\rho = \tilde \mu$.

Let $\Delta_{i_1}, \ldots, \Delta_{i_h}$ be the tetrahedra incident to $v$ (with multiplicities), and let $\tilde\Delta_{i_1}, \ldots, \tilde\Delta_{i_h}$ be any lifts in $\tilde\calT$ incident to $\tilde v$. We fix a sufficiently small $R_0>0$ such that the cone $N_{R_0}(\tilde l)\subset \overline X$ does not intersect the lower faces of the tetrahedra $\tilde D(\tilde\Delta_{i_1}), \ldots, \tilde D(\tilde\Delta_{i_h})$, see Figure \ref{Dehn_fill:fig}-(left). Since $N_{R_0}(\tilde l)$ is $\rho(\pi_1(T))$-invariant, it does not intersect the lower faces of any ideal tetrahedron in the image of $\tilde D$, and hence it intersects every ideal tetrahedron in a curved triangle as in Figure \ref{Dehn_fill:fig}-(right).

We define $\widetilde {C(T)}\subset \tilde M$ as the preimage of $N_{R_0}(\tilde l)$ along $\tilde D$. By construction it is a $\pi_1(T)$-invariant submanifold that projects to a submanifold $C(T) \subset M$. This submanifold intersects every tetrahedron incident to $v$ into a cone neighbourhood of $v$ and is hence a collar for $T$. We want to prove that the map
$$\tilde D\colon \widetilde {C(T)} \longrightarrow N_{R_0}(\tilde l) \setminus \tilde l$$
is an isometry. Being a local isometry, it suffices to prove that it is injective and surjective. We prove this using the natural foliation of $N_{R_0}(\tilde l)\setminus \tilde l$ into the sheets $\partial N_{R}(\tilde l)$ with $R\leqslant R_0$.

For every $R\leqslant R_0$ the map $\tilde D$ restricts to a local isometry of surfaces
$$\tilde D_R\colon \tilde {T}_R\to \partial N_R(\tilde l)$$
where $\tilde {T}_R = \tilde D^{-1}(\partial N_R(\tilde l))$ is $\pi_1(T)$-invariant and covers a torus $T_R\subset M$ parallel to $T$. Since $T_R$ is compact, the cover $\tilde T_R$ is complete. Therefore $\tilde D_R$ is a covering by Proposition \ref{local:isometry:prop} and is hence an isometry since $\partial N_R(\tilde l)$ is simply-connected. In particular it is a bijection for all $R<R_0$, and hence $\tilde D$ is a bijection.

Since $\tilde D$ is an isometry, its holonomy $\tilde\rho$ is discrete and injective and 
$$C(T) \isom \big(N_{R_0}(\tilde l) \setminus \tilde l \big)/_{\tilde\rho(\pi_1(T))}.$$
The proof is complete.
\end{proof}

\subsection{The completion}
We now have an isometric model for the collar $C(T)$ of $T$ and we study its completion. This analysis will lead to a proof of Theorem \ref{modified:equations:teo}.

We suppose that $\mu$ is non-trivial and hence $C(T)$ is not complete. By Proposition \ref{DF:prop}, the completion $\overline{C(T)}$ of $C(T)$ depends on the $\matZ\times \matZ$ group $\Img \tilde \mu< \widetilde{\matC^*}$ acting on $X$: see Section \ref{X:quotients:subsection} where we defined in particular the (truncated) tubes with cone angles. 

\begin{prop} \label{DF1:prop}
If $\Img \tilde \mu$ contains a non-trivial rotation then $\overline{C(T)}$ is a truncated tube with some cone angle; otherwise it is a one-point compactification of $C(T)$. 

In the first case, there are generators $\gamma, \eta$ for $\pi_1(T)$ such that 
$$\tilde \mu(\gamma) = e^{\alpha i}, \qquad \tilde \mu(\eta) = e^{\lambda + \beta i}$$
and the core geodesic of $\overline {C(T)}$ has cone angle $\alpha$ and length $|\lambda|$. The curve $\gamma$ is a meridian of the truncated tube.
\end{prop}
\begin{proof}
We know that $C(T)$ is isometric to a truncation of $X/_{\Img \tilde \mu}$ and we apply the discussion of Section \ref{X:quotients:subsection}.
\end{proof}

We now look at $M$ globally: every boundary torus $T_i\subset \partial N$ has its own homomorphisms $\mu$ and $\tilde \mu$.

\begin{cor} \label{DF2:cor}
Suppose that at every boundary torus $T_i\subset \partial N$ one of the following holds:
\begin{enumerate}
\item $\mu$ is trivial, or
\item there is a primitive $\gamma_i\in \pi_1(T_i)$ such that $\tilde\mu(\gamma_i) = e^{\alpha_i i}$.
\end{enumerate}
The completion $\overline M$ is a complete hyperbolic cone manifold obtained by Dehn filling the slopes $\gamma_i$. The core of the $i$-th Dehn filling is a closed geodesic with cone angle $\alpha_i$.
\end{cor}

Theorem \ref{modified:equations:teo} now follows immediately.

\dimo{modified:equations:teo}
The equation $\tilde\mu(m_i)^p\tilde\mu(l_i)^q = e^{2\pi i}$ implies that $\tilde\mu(pm_i + ql_i) = e^{2\pi i}$ and therefore the completion is a hyperbolic manifold (with no cone angles since $\alpha = 2\pi$) diffeomorphic to $M^{\rm fill}$.
\finedimo

\subsection{Generalised Dehn filling invariants} \label{generalized:DF:subsection}
We have proved Theorem \ref{modified:equations:teo} and we now make some comments that will be useful in the next chapter. Let $z$ be a solution of the consistency equations and consider a boundary torus $T_i\subset \partial N$ with fixed generators $m_i,l_i$ for $\pi_1(T)$.

\begin{prop} \label{pq:prop}
If $\mu$ is non-trivial, there is a unique $(p,q)\in\matR^2$ such that $\tilde\mu(m_i)^p \tilde\mu(l_i)^q = e^{2\pi i}$.
\end{prop}
\begin{proof}
We know from Proposition \ref{DF:prop} that $\tilde\mu$ is injective and has discrete image. Therefore $\tilde\mu(m_i)$ and $\tilde\mu(l_i)$ form a basis of $\widetilde{\matC^*}$ considered as a $\matR$-vector space.
\end{proof}

We now define the \emph{generalised Dehn filling invariant} $(p,q)\in S^2 = \matR^2 \cup \{\infty\}$ of the solution $z$ at the torus $T_i$ to be:
\begin{itemize}
\item $(p,q) = \infty$ if $\mu$ is trivial,
\item the $(p,q)\in\matR^2$ from Proposition \ref{pq:prop} if $\mu$ is non-trivial.
\end{itemize}

If $(p,q)=\infty$ there is a complete collar $C(T_i)$ of $T_i$ that is a truncated cusp; if $(p,q) \in \matR^2$ then Proposition \ref{DF:prop} provides a nice incomplete collar $C(T_i)$ and Proposition \ref{DF1:prop} tells us everything about its completion $\overline{C(T_i)}$. Namely, the following holds:

\begin{itemize}
\item if $\frac pq\in\matQ\cup \{\infty\}$ then $(p,q) = k(r,s)$ for a unique real number $k>0$ and coprime integers $(r,s)$, and $\overline{C(T_i)}$ is a solid torus with meridian $rm+sl$, isometric to a truncated tube with cone angle $\frac{2\pi}k$;
\item if $\frac pq \not\in \matQ \cup \{\infty\}$ then $\overline{C(T_i)}$ is the much less interesting one-point compactification.
\end{itemize}
We deduce the following:
\begin{itemize}
\item if $(p,q)$ are coprime integers then $\overline{C(T_i)}$ is a standard truncated tube; if this holds at all boundary tori $T_i$ then $\overline M$ is a hyperbolic manifold; 
\item if $(p,q)$ are integers then $k$ is a natural number and the cone angle divides $2\pi$; if this holds at all boundary tori $T_i$ then $\overline M$ may be interpreted as an orbifold thanks to Proposition \ref{from:cone:to:orbifold:prop};
\item if $\frac pq \in \matQ\cup \{\infty\}$ at all boundary tori $T_i$ then $\overline M$ is a hyperbolic cone manifold.
\end{itemize}

\begin{cor}
If the generalised Dehn filling invariants $(p,q)$ at each $T_i$ are either $\infty$ or coprime integers, the completion $\overline M$ is a complete hyperbolic manifold obtained by $(p,q)$-Dehn filling the cusps of the second type.
\end{cor}

We end this discussion by calculating $(p,q)$ explicitly. We define
$$u_i = \log \tilde \mu(m_i), \quad v_i = \log \tilde \mu(l_i)$$
and find the following.
\begin{prop} \label{pq:explicit:prop}
If $\mu$ is non trivial we have
$$p = - 2\pi\frac{\Re v_i}{\Im(\bar u_iv_i)}, \qquad 
q = 2\pi\frac{\Re u_i}{\Im(\bar u_iv_i)}.$$
\end{prop}
\begin{proof}
The pair $(p,q)$ is the unique solution of  
\begin{equation} \label{pquv:eqn}
pu_i + qv_i = 2\pi i
\end{equation}
when $u_i,v_i \neq 0$, or is $\infty$ otherwise. Note that $u_i=0 \Leftrightarrow v_i=0$.
The relation
$$i\Im (\bar u_iv) = -(\Re v_i)u_i + (\Re u_i)v_i$$
implies easily that the pair $(p,q)$ stated above is a solution to (\ref{pquv:eqn}).
\end{proof}

\subsection{Closed census}
The Lickorish-Wallace Theorem \ref{LW:teo} states that
every closed orientable three-manifold $M$ is the result of a Dehn surgery along some link $L\subset S^3$, and as such it can be easily presented to SnapPea. The program tries to solve numerically the consistency and completeness equations for $M$, based on some ideal triangulation for the complement of $L$.

If it succeeds to find a solution $z$, the closed manifold $M$ is hyperbolic and SnapPea calculates numerically various geometric invariants of $M$: the volume (which is just the sum of the volumes of the ideal tetrahedra), the length of the core geodesics of the filling solid tori (using Proposition \ref{DF1:prop}), a segment of the geodesic spectrum, a Dirichlet domain, etc.

Various closed manifolds $M$ have been tested and listed, and the ten closed hyperbolic three-manifolds of smallest volume known today are in Table \ref{closed:census:table}, taken from Hodgson -- Weeks \cite{Ho-We}.

\begin{table}
\begin{center}
\begin{tabular}{c||ccccccc}
Name & Volume & Homology & Symmetry & SG & C \\
 \hline \hline
Vol1 & $0.94270736$ & $\matZ_5+\matZ_5$ & $D_{12}$ & $0.5846$ & c \\
Vol2 & $0.98136883$ & $\matZ_5$ & $D_{4}$ & $0.5780$ & c \\
Vol3 & $1.01494161$ & $\matZ_3 + \matZ_6$ & $S_{16}$ & $0.8314$ & a \\
Vol4 & $1.26370924$ & $\matZ_5 + \matZ_5$ & $D_8$ & $0.5750$ & c \\
Vol5 & $1.28448530$ & $\matZ_6$ & $D_4$ & $0.4803$ & c \\
Vol6 & $1.39850888$ & $\{e\}$ & $D_4$ & $0.3661$ & c \\
Vol7 & $1.41406104$ & $\matZ_6$ & $D_4$ & $0.7941$ & c \\ 
Vol8 & $1.41406104$ & $\matZ_{10}$ & $D_4$ & $0.3648$ & c \\
Vol9 & $1.42361190$ & $\matZ_{35}$ & $D_4$ & $0.3523$ & c \\
Vol10 & $1.44069901$ & $\matZ_3$ & $D_4$ & $0.3615$ & c
\end{tabular}
\vspace{.2 cm}
\nota{The ten closed hyperbolic three-manifolds with smallest volume known. Here $S_{16}$ is the semidihedral group of order 16 with presentation $\langle x,y\ |\ x^8 = y^2 = 1, y^{-1}xy = x^3\rangle$, SG indicates the length of the shortest geodesics (volume and SG values are truncated after few digits), and C indicates whether the manifold is achiral (a) or chiral (c), that is whether it admits an orientation-reversing isometry or not.}
\label{closed:census:table}
\end{center}
\end{table}

\subsection{References}
The material contained in this chapter originated from Thurston's notes \cite{Th} and Neumann -- Zagier \cite{NZ}, see also Benedetti -- Petronio \cite{BP}. The program SnapPea is freely available and can now be used via a Python interface \cite{Sna}. We have used the computer censuses of Callahan -- Hildebrand -- Weeks \cite{CaHiWe}, Hoste -- Thistlethwaite -- Weeks \cite{HTW}, and Hodgson -- Weeks \cite{Ho-We}.

It is conjectured in \cite{A} that each link in Figure \ref{sequence:fig} with $c=1,\ldots, 5$ cusps has minimum volume among orientable hyperbolic manifolds with $c$ cusps. This has been proved for $c=1$ by Cao and Meyerhoff \cite{CM}, for $c=2$ by Agol \cite{A}, and for $c=4$ by Yoshida \cite{Yo}. The manifold Vol1 from Table \ref{closed:census:table} is called the \emph{Fomenko -- Matveev -- Weeks manifold} and has indeed smallest volume among all complete hyperbolic three-manifolds by Gabai -- R.~Meyerhoff -- P.~Milley \cite{GMM}.

%% file: Hyperbolic_Dehn.tex
\chapter{Hyperbolic Dehn filling} \label{Hyperbolic:Dehn:chapter}
We have completely classified in Chapter \ref{Seifert:chapter} the closed three-manifolds that belong to the six Seifert geometries, and now we long for a similar catalogue that displays all the closed hyperbolic three-manifolds that exist in nature. Is there something like a ``name'' to assign to each closed hyperbolic three-manifold, together with some reasonable tables that list all possible names?

There is not yet one such thing, and we do not know if there will ever be one: we still do not understand hyperbolic three-manifolds globally. The main difficulty is that there are really many hyperbolic three-manifolds around, so many that topologists often say informally that ``most three-manifolds are hyperbolic.''

This folk sentence is supported by a fundamental theorem that we prove in this chapter, the \emph{Hyperbolic Dehn filling Theorem}, which says roughly that by Dehn-filling generically a cusped hyperbolic three-manifold we still get a hyperbolic manifold.\index{hyperbolic Dehn filling theorem}

\section{Introduction}
It sometimes happens in three-dimensional topology, that a topological or geometric property of a manifold is preserved under Dehn fillings, with only few exceptions.
For instance, the Dehn filling of a Seifert manifold is again a Seifert manifold, with only one exception (a fibre-parallel Dehn filling gives a connected sum of Seifert manifolds, see Corollary~\ref{DF:Seifert:cor}).

The most striking appearance of this phenomenon is the Hyperbolic Dehn filling Theorem, proved by Thurston at the end of the 1970s. This theorem says roughly that ``most'' Dehn fillings of a cusped hyperbolic manifold are still hyperbolic.

\subsection{Generalised Dehn filling parameters} \label{generalised:subsection}
In all this chapter, we consider a compact oriented three-manifold $N$ with $\partial N = T_1 \sqcup \ldots \sqcup T_c$ consisting of tori, and we fix generators $m_i,l_i$ for every $\pi_1(T_i)$. 

Let a \emph{generalised Dehn filling parameter} $s=(s_1,\ldots, s_c)$ be a sequence where each $s_i$ is either the symbol $\infty$ or a rationally related pair of real numbers $(p,q) = (kp',kq') = k(p',q')$ where $k>0$ is real and $(p',q')$ are coprime integers. We think at $s_i$ as lying in the two-sphere $S^2 = \matR^2 \cup \{\infty\}$.

A generalised Dehn filling parameter $s$ determines a Dehn filling $N^{\rm fill}$ of $N$ as follows:
for every $i$, if $s_i=k(p',q')$ we fill $T_i$ by killing the slope $p'm_i + q'l_i$, while if $s_i=\infty$ we do nothing. The result is a new compact manifold, that may be closed (if there are no $\infty$ in $s$) or be bounded by some tori. We also mark the cores of the filled solid tori with the label $\alpha_i = \frac{2\pi} k>0$. 

\subsection{The Hyperbolic Dehn filling Theorem}
This chapter is mainly devoted to the proof of the following theorem. 

\begin{teo} \label{DF:teo} \label{hyperbolic:Dehn:filling:teo}
Let $M = \interior N$ be a complete orientable finite-volume cusped hyperbolic three-manifold. There is a neighbourhood $U$ of $(\infty,\ldots,\infty)$ in $S^2\times \ldots \times S^2$ such that for every generalised Dehn filling parameter $s \in U$ the interior $M^{\rm fill}$ of the manifold $N^{\rm fill}$ obtained by Dehn filling $N$ along $s$ admits a finite-volume complete hyperbolic structure with cone angles. 

The cores of the filling solid tori are closed geodesics with cone angles $\alpha_i$. The singular locus of $M^{\rm fill}$ consists of the core geodesics with $\alpha_i \neq 2\pi$.
\end{teo}

If $s\in U$ is an ordinary (not extended) Dehn filling parameter (see Section \ref{DF:parameters:subsection}), that is if $s_i$ is either $\infty$ or a coprime pair $(p,q)$ for all $i$, then $\alpha_i = 2\pi$ and $N^{\rm fill}$ is a hyperbolic manifold without cone angles: this is of course the case of most interest. 

\begin{cor}
Let $M = \interior N$ be a complete orientable finite-volume hyperbolic three-manifold. For every $i=1,\ldots, c$ there is a finite set $S_i$ of slopes in $T_i$ such that for every Dehn filling parameter $s$ with $s_i \not \in S_i$ for all $i$, the filled manifold $M^{\rm fill} = \interior{N^{\rm fill}}$ is hyperbolic.
\end{cor}

When $M$ has one cusp, the corollary may be stated simply as follows.

\begin{cor}
If $M$ is a complete orientable finite-volume hyperbolic three-manifold with one cusp, all but finitely many Dehn fillings $M^{\rm fill}$ are hyperbolic.
\end{cor}

To appreciate the power of these theorems, consider for instance a hyperbolic link $L\subset S^3$ (recall that $L\subset S^3$ is hyperbolic if $S^3\setminus L$ admits a finite-volume complete hyperbolic metric). By the hyperbolic Dehn filling theorem, there is a finite subset $S\subset \matQ$ such that every surgery on $L$ with coefficients in $\matQ\setminus S$ produces a closed hyperbolic manifold.

\subsection{Examples} We describe a couple of clarifying examples. We know that the figure eight knot $K\subset S^3$ is hyperbolic: by the Dehn filling Theorem, there is an open neighbourhood $U \subset S^2$ of $\infty$ such that every Dehn surgery on $K$ with parameter $s\in U$ produces a hyperbolic cone manifold. 

\begin{figure}
\begin{center}
\includegraphics[width = 9 cm] {\iftoggle{BW}{figure8_U-BW}{figure8_U}} 
\nota{Hyperbolic Dehn fillings on the figure-eight knot complement.}
\label{figure8_U:fig}
\end{center}
\end{figure}

We will prove in the next section that the subset $U$ shown in Figure \ref{figure8_U:fig} fulfils this requirement. The only coprime pairs $(p,q)$ that are not contained in $U$ are those with $\frac pq$ equal to one of the following
\begin{equation} \label{K:list:eqn}
-4,-3,-2,-1,0,1,2,3,4,\infty.
\end{equation}
Every other surgery on $K$ yields a closed hyperbolic three-manifold. For instance, the manifold Vol2 from Table \ref{closed:census:table} can be obtained with $\frac pq = \pm 5$ (there is a symmetry of $K$ sending $\frac pq$ to $-\frac pq$).

On the other hand, a $\frac pq$ surgery from the list (\ref{K:list:eqn}) does not produce a hyperbolic manifold: see Table \ref{non:hyp:table}, taken from Martelli -- Petronio \cite{MP}.

\begin{table}
\begin{center}
\begin{tabular}{c||c}
$\frac pq$ & manifold \phantom{\Big|}\!\! \\
\hline
\hline
$0$ & $T_{\matr 31{-1}0}$ \phantom{\Big|}\!\!  \\
$\pm 1$ & $\big(S^2, (2,1), (3,1), (7,-6)\big)$  \phantom{\Big|}\!\! \\
$\pm 2$ & $\big(S^2, (2,1), (4,1), (5,-4)\big)$  \phantom{\Big|}\!\! \\
$\pm 3$ & $\big(S^2, (3,1), (3,1), (4,-3)\big)$  \phantom{\Big|}\!\! \\
$\pm 4$ & $\big(D,(2,1), (2,1)\big) \bigcup_{\matr 0110} \big(D, (2,1), (3,1)\big) $\phantom{\Big|}\!\!
\end{tabular}
\vspace{.2 cm}
\nota{The non-hyperbolic Dehn surgeries of the figure-eight knot. We get a torus bundle with Anosov monodromy and hence of $\Sol$ geometry, three Seifert manifolds of $\widetilde {\SL_2}$ geometry, and a graph manifold that splits into two Seifert manifolds via a map that interchanges the fibres and the boundary sections of the two portions (the map is here denoted by a matrix in the section-fibre basis).}
\label{non:hyp:table}
\end{center}
\end{table}

More generally, a parameter $k(p,q) \in U$ produces a hyperbolic closed manifold with a core geodesic of cone angle $\frac{2\pi}k$, which may be interpreted as an orbifold when $k\in \matN$. In particular, if $(p,q) = (1,0)$ the closed manifold is just $S^3$, because an $\infty$-surgery on a knot in $S^3$ gives $S^3$ back. Therefore by looking at Figure \ref{figure8_U:fig} we discover that for every $k\geqslant 5$ there is a hyperbolic structure on $S^3$ with cone angle $\frac{2\pi}k$ on the figure-eight knot $K$. 

We now consider a link with two components.
The \emph{Whitehead link} $L$ shown in Figure \ref{sequence:fig} is hyperbolic: it is hard to determine an explicit open set $U\subset S^2\times S^2$ in this case, but a careful computer-aided analysis carried with SnapPy \cite{Sna} shows the following, see \cite{MP}.

\begin{teo}
A surgery on the Whitehead link with coefficients $(\frac pq, \frac rs)$ produces a closed hyperbolic manifold, unless one of the following holds:
\begin{itemize}
\item $\frac pq$ or $\frac rs$ belongs to the set $\{0, 1, 2, 3, 4, \infty\}$,
\item up to permutation the pair $\big( \frac pq, \frac rs \big)$ belongs to the set
$$\left\{ (-4,-1), (-3,-1), (-2,-2), (-2,-1), \left(\tfrac 32, 5\right), \left(\tfrac 43, 5\right) \left(\tfrac 52, \tfrac 72\right) \right\}.$$
\end{itemize}
\end{teo}
The proposition is symmetric in the two coefficients because the link itself has a symmetry that interchanges the two components. All the Dehn filling parameters listed in the theorem indeed produce non-hyperbolic manifolds. For instance, if $\frac pq\in \{1,2,3\}$ we get a Seifert manifold fibering over the disc with two singular fibres, see \cite{MP}.

We first prove Theorem \ref{hyperbolic:Dehn:filling:teo} for the figure-eight knot, where the combinatorics is so simple that everything can be verified by hand.
We will then prove the theorem in general in Section \ref{solution:space:section} using more sophisticated tools.

\subsection{The figure-eight knot example} \label{figure:8:example:subsection}

\begin{figure}
\begin{center}
\includegraphics[width = 11 cm] {\iftoggle{BW}{triangulation_figure8_2-BW}{triangulation_figure8_2}} 
\nota{We truncate the tetrahedra and flatten their boundary. We assign a modulus $z$ and $w$ to them. The edge numbered $i=1,2,3$ has modulus $z_i$ (on the left tetrahedron) or $w_i$ (on the right tetrahedron).}
\label{triangulation_figure8_2:fig}
\end{center}
\end{figure}

\begin{figure}
\begin{center}
\includegraphics[width = 9 cm] {\iftoggle{BW}{triangulation_figure8_3-BW}{triangulation_figure8_3}} 
\nota{The triangulated boundary is a torus obtained by identifying the opposite edges of this parallelogram. The triangulation $\calT$ has two edges, that we colour in \iftoggle{BW}{grey}{red} and white. Their endpoints are shown here.}
\label{triangulation_figure8_3:fig}
\end{center}
\end{figure}

We now discuss in detail the standard example: the figure-eight knot complement, one of the very few cases where the consistency and completeness equations can be solved by hand. This part will not be needed in the proof of Theorem \ref{hyperbolic:Dehn:filling:teo}, so the reader may wish to skip it and go directly to Section \ref{sketch:filling:subsection}.

Let $\calT$ be the triangulation with two tetrahedra shown in Figure \ref{triangulation_figure8:fig}. The truncated version is in Figure \ref{triangulation_figure8_2:fig}. We assign the complex variable $z$ and $w$ to the left and right tetrahedron, respectively. We set
$$z_1 = z, \quad z_2 = \frac{1}{1-z}, \quad z_3 = \frac {z-1}{z}$$
and define $w_1,w_2,w_3$ similarly. Recall that
\begin{equation} \label{m1:eq}
z_1z_2z_3 = w_1w_2w_3 = -1.
\end{equation}
The modulus of an edge in Figure \ref{triangulation_figure8_2:fig} labeled with $i\in \{1,2,3\}$ is $z_i$ or $w_i$ depending on the tetrahedron. With some patience one sees that the boundary triangulated surface is a torus $T$ as in Figure \ref{triangulation_figure8_3:fig}. The triangulation $\calT$ has two edges shown in Figure \ref{triangulation_figure8:fig}, each yielding a consistency equation. Figure \ref{triangulation_figure8_3:fig} shows that the consistency equations are
$$z_2^2z_3w_2^2w_3 = 1, \quad z_1^2z_3w_1^2w_3=1.$$
Using (\ref{m1:eq}) we see that the equations are both equivalent to 
$$z_2w_2 = z_1w_1,$$ 
that is
\begin{equation} \label{c:eqn}
z(1-z)w(1-w) = 1.
\end{equation}
We now look at the boundary torus $T$. Let $m$ and $l$ be the generators of $\pi_1(T)$ shown in Figure \ref{triangulation_figure8_3:fig}. We have
\begin{equation} \label{ml:eqn}
\mu(m) = - z_1z_3w_1 = w(1-z), \quad \mu(l) = z_2^2z_3^2w_1^4w_2^2w_3^2 = \frac{z^2}{w^2}.
\end{equation}
The completeness equations are
\begin{equation} \label{c2:eqn}
w(1-z) = 1,\quad z^2 = w^2.
\end{equation}

\begin{ex}
The only solution to (\ref{c:eqn}) and (\ref{c2:eqn}) is $z=w= e^{\frac{\pi i}3}$.
\end{ex}

We have confirmed that the figure-eight complement $M$ has a complete hyperbolic structure obtained by representing both tetrahedra with ideal regular hyperbolic tetrahedra. We now investigate the non-complete solutions of the consistency equations.

The two ideal tetrahedra have moduli $w$ and $z$ and the consistency equations reduce to one equation (\ref{c:eqn}) which we rewrite as
$$z^2 - z +\frac{1}{w(1-w)} = 0.$$
The solutions are
$$ z = \frac{1 \pm \sqrt{1-\frac 4{w(1-w)}}}{2}.$$
We are only interested in solutions with $\Im z, \Im w>0$. For every $w$ with $\Im w>0$ there is a unique solution $z$ with $\Im z >0$, except when $\Delta \in \matR_{\geqslant 0}$.

\begin{ex} 
We have $\Delta = 1-\frac 4{w(1-w)}\in\matR_{\geqslant 0}$ if and only if $w$ belongs to the half-line $s = \big\{\frac 12 + yi {\rm\ with\  } y\geqslant \frac{\sqrt{15}}2\big\}$, see Figure \ref{figure8_R:fig}-(left).
\end{ex}

\begin{figure}
\begin{center}
\includegraphics[width = 12.5 cm] {\iftoggle{BW}{figure8_R-BW}{figure8_R}} 
\nota{Every $w\in R$ determines a hyperbolic structure for $M$ (left). The Dehn filling generalised invariants $d$ map $R$ onto a neighbourhood of $\infty$ (right).}
\label{figure8_R:fig}
\end{center}
\end{figure}

We define the open region 
$$R =  \{w\in \matC\ |\ \Im w >0 \} \setminus s.$$
Every $w\in R$ determines a hyperbolic structure for $M$. The complete structure is obtained at $w_0=\frac 12 + \frac {\sqrt 3}2$, see Figure \ref{figure8_R:fig}-(left). 

The manifold $M$ is a knot complement and hence the boundary torus $T$ is equipped with its natural meridian/longitude basis $m',l'$, see Section \ref{canonical:longitudes:subsection}. Using SnapPea we find out that
$$m'=m, \qquad l'=l+2m.$$
From (\ref{ml:eqn}) we get
\begin{equation} \label{mlp:eqn}
 \tilde\mu(m') = w(1-z), \quad \tilde\mu(l') = \frac{z^2}{w^2} w^2(1-z)^2 = z^2(1-z)^2. 
\end{equation}
We fix $m',l'$ as generators for $\pi_1(T)$. We defined in Section \ref{generalized:DF:subsection} the generalised Dehn surgery invariant $(p,q) = d(w)$ for every non-complete solution $w\in R \setminus \{w_0\}$. This gives a continuous map
$$d \colon R \longrightarrow S^2$$
that sends $w_0$ to $\infty$.
\begin{prop}
The image $d(R)$ contains the coloured region shown in Figure \ref{figure8_R:fig}-(right).
\end{prop}
\begin{proof}
The domain $R$ is an open disc and its abstract closure $\bar R$ is homeomorphic to the closed disc. We show that $d$ extends to a continuous map $\bar R \to S^2$ which sends $\partial R$ to the almost-rectangle shown in Figure \ref{figure8_R:fig}. 

The invariants $d(w) = (p,q)$ are such that
\begin{equation} \label{pq:eqn}
w^p(1-z)^{p+2q}z^{2q} = e^{2\pi i}.
\end{equation}

The region $R$ has two involutions (determined in fact by isometries of $M$):
\begin{itemize}
\item the involution $\tau(w) = z = \frac{1 \pm \sqrt{1-\frac 4{w(1-w)}}}{2}$. This involution sends $\mu(m)$ to $\mu(m)^{-1} = \mu(m^{-1})$, hence $d(\tau(w)) = -d(w)$;
\item the involution $\sigma(w) = \overline{1-w}$: using (\ref{mlp:eqn}) we get $d(\sigma(w)) = \overline{d(w)}$.
\end{itemize}
Now we note that when $w$ tends to a point in the line $l_+ = [1, +\infty]$ then $z$ tends to a point in $l_- = [-\infty,0]$. Therefore the argument of $w,z,1-z$
tends respectively to $0, \pi, 0$. Equation (\ref{pq:eqn}) implies that $q \to 1$.

If $w\to 1$ then $z \to -\infty$ and (\ref{pq:eqn}) implies that $p+2q+2q=p+4q\to 0$, so $(p,q)\to (-4,1)$. If $w\to +\infty$ then $z \to 0$ and equation $z(1-z)w(1-w)=1$ implies that $|z||w|^2\to 1$; here (\ref{pq:eqn}) gives $|w|^{p-4q}\to 1$, hence $p-4q\to 0$, so $(p,q)\to (4,1)$. We have proved that $d$ maps $[1,+\infty]$ onto the segment 
$$\overline{(-4,+1), (+4,+1)}.$$
Using the involution $\tau$ we deduce that $d$ maps $[-\infty, 0]$ onto the segment
$$\overline{(+4,-1),(-4,-1)}.$$
When $w$ is near the right side of the half-line $s$, the number $z$ tends to the segment $(0,\frac 12]$. Therefore the arguments of $z$ and $1-z$ tend to $0$ and (\ref{pq:eqn}) implies that the argument of $w^p$ tends to $2\pi$. When $w\in s$ the argument of $w$ is at least $\arctan \sqrt{15} = 1.31811607\ldots $ and hence $p$ is at most $\frac{2\pi}{\arctan \sqrt{15}} = 4.374\ldots < 5$. This implies that we can connect $(4,1)$ and $(4,-1)$ while staying inside the image of $d$ with a curve as in Figure \ref{figure8_R:fig}.
\end{proof}

\subsection{A road map} \label{sketch:filling:subsection}
The rest of this chapter is mostly devoted to the proof of Theorem \ref{hyperbolic:Dehn:filling:teo}. Here is our plan: we first study the combinatorial properties of an arbitrary ideal triangulation $\calT$ for $M$, and we prove that the solutions to the consistency equations form a complex manifold $\Def(M,\calT)$ of (complex) dimension $c$ equal to the number of cusps of $M$. The manifold $\Def(M,\calT)$ is nicely parametrized by the holonomies on any $c$ fixed curves at different cusps. 

The manifold $\Def(M,\calT)$ is often empty, but since $M$ is hyperbolic we expect that there should be some ideal triangulation $\calT$ for which $\Def(M,\calT)$ is not empty and contains a solution $z$ that also satisfies the completeness equations (we are unfortunately not able to prove this in general, and this will be a technical issue). The generalised Dehn filling invariants defined in Section \ref{generalized:DF:subsection} furnish a local diffeomorphism $d\colon \Def(M,\calT) \to S^2 \times \ldots \times S^2$ that sends $z$ to $(\infty,\ldots,\infty)$. In particular this map is open and its image covers a neighbourhood $U$ of $(\infty,\ldots, \infty)$. We conclude thanks to the discussion of Section \ref{generalized:DF:subsection}. 

We start by defining and exploring the solution space $\Def(M,\calT)$.

\section{The solution space} \label{solution:space:section}
In this section we study the combinatorial properties of ideal triangulations of three-manifolds in general. We study the space of solutions to the consistency equations, we prove that it is a complex manifold and exhibit a concrete parametrisation via the holonomy of peripheral curves.

\subsection{Edges and tetrahedra}
Throughout this section, we let $N$ be a compact oriented three-manifold $N$ with $\partial N = T_1 \sqcup \ldots \sqcup T_c$ consisting of tori, and define $M=\interior N$. Let $\calT$ be any ideal triangulation for $M$. 

\begin{prop} \label{edges:prop}
The ideal triangulation $\calT$ has the same number $n$ of tetrahedra and edges.
\end{prop}
\begin{proof}
The total space $|\calT|$ of the triangulation $\calT$ has $v$ vertices, $e$ edges, $f$ faces, and $t$ tetrahedra. Since $\chi(\partial N) = 2\chi(N)$ for every compact 3-manifold $N$ and $\partial N$ consists of tori, we get $\chi(N)=0$. Therefore
$$v = v + \chi(N) = \chi(|\calT|) = v-e+f-t = v-e+t$$
since $f = 2t$. Then $e=t$.
\end{proof}

Recall that in the consistency equations we get one complex variable $z_i$ for each tetrahedron and one equation for each edge; since we have the same number of edges and tetrahedra, one would guess that the solutions form typically a discrete set of points, but this is surprisingly not the case: we now show that some equations are redundant, leaving some space for a higher-dimensional space of solutions. The origin of this redundancy lies in the combinatorial properties of three-dimensional triangulations.

\subsection{Incidence matrices} \label{incidence:matrices:subsection}
\begin{figure}
\begin{center}
\includegraphics[width = 5 cm] {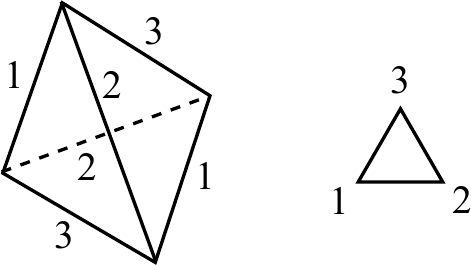} 
\nota{We label the edges of every \emph{oriented} tetrahedron as shown: the labelling is determined once we chose arbitrarily a pair of opposite edges and label them with 1 (left). The vertices of every \emph{oriented} triangle in the triangulated boundary inherit labels as in the figure (right).}
\label{tetrahedron_labels:fig}
\end{center}
\end{figure}

We denote respectively by $\Delta_1, \ldots, \Delta_n$ and $e_1,\ldots, e_n$ the tetrahedra and edges of $\calT$. We fix an orientation on $M$, which induces an orientation on each tetrahedron, and we assign numbers 1, 2, and 3 to pairs of opposite edges on each $\Delta_i$ in a way that is orientation-preservingly isomorphic to Figure \ref{tetrahedron_labels:fig}.

\begin{figure}
\begin{center}
\includegraphics[width = 11 cm] {\iftoggle{BW}{ABmatrices-BW}{ABmatrices}} 
\nota{The incidence matrices $A$ and $B$ count the incidences between edge pairs, tetrahedra, edges, and ideal vertices.}
\label{ABmatrices:fig}
\end{center}
\end{figure}

We now define a $2n \times 3n$ integral matrix $A$ that encodes some combinatorial adjacencies between tetrahedra. The matrix $A$ should be seen as a $2\times 3$ rectangle of $n\times n$ sub-matrices as in Figure \ref{ABmatrices:fig}-(left), where:
\begin{itemize}
\item each column of $A$ corresponds to a pair of opposite edges of one tetrahedron, more precisely the columns $i, i+n,$ and $i+2n$ correspond the opposite edges of type 1,2,3 of $\Delta_i$;
\item the first $n$ rows of $A$ correspond to the tetrahedra $\Delta_1,\ldots, \Delta_n$;
\item the last $n$ rows of $A$ correspond to the edges $e_1,\ldots, e_n$.
\end{itemize}
When $i\leqslant n$, the entry $A_{ij}$ is 1 or 0 depending on whether $\Delta_i$ contains the $j$-th edge pair. In other words, the top three $n\times n$ sub-matrices of $A$ are three identity matrices, see Figure \ref{ABmatrices:fig}-(left).

When $i=n+i'$ with $i'>0$, the entry $A_{ij}$ counts how many edges in the $j$-th edge pair are glued to the $i'$-th edge $e_{i'}$: this number is hence 0, 1, or 2.

\begin{ex}
The entries on each column of $A$ sum to $3$.
\end{ex}

Let now $v_1,\ldots, v_c$ be the ideal vertices of $\calT$. We define another $c \times 2n$ integral matrix $B$, which consists of two $c\times n$ sub-matrices as in Figure \ref{ABmatrices:fig}-(right), where:
\begin{itemize}
\item the rows of $B$ correspond to the ideal vertices $v_1,\ldots, v_c$;
\item the first $n$ columns of $B$ correspond to the $n$ tetrahedra $\Delta_1,\ldots, \Delta_n$;
\item the last $n$ columns correspond to the edges $e_1,\ldots, e_n$.
\end{itemize}

When $j \leqslant n$, the entry $B_{ij}$ is \emph{minus} the number of times $\Delta_j$ is incident to $v_i$, and when $j = n+j'$ with $j'>0$ the entry $B_{ij}$ is \emph{plus} the number of times $e_{j'}$ is incident to $v_i$. The possible entries are $0,-1,-2,-3,-4$ and $0,1,2$ respectively.

\begin{ex}
The entries on each column of $B$ sum to $-4$ or $2$.
\end{ex}

The matrices $A$ and $B$ are designed to be interpreted as linear maps, and to be composed.

\begin{prop} \label{AB:exact:prop}
The following sequence is exact:
$$ \matR^{3n} \stackrel{A}{\longrightarrow} \matR^{2n} \stackrel{B}{\longrightarrow} \matR^c \longrightarrow 0.$$
\end{prop}
The proof of this proposition splits into three lemmas. 

\begin{lemma} We have $BA = 0$.
\end{lemma}
\begin{proof}
The $i$-th row of $B$ corresponds to the ideal vertex $v_i$ and the $k$-th column of $A$ corresponds to the $k$-th pair of opposite edges. The sum $\sum_{j=1}^{2n} B_{ij} A_{jk}$ counts \emph{minus} the number of times $v_i$ is incident to the tetrahedron containing the edge pair, \emph{plus} the number of times it is incident to the two edges of the pair. These numbers are equal and hence we get 0.
\end{proof}

\begin{lemma} The $c$ rows in $B$ are independent vectors.
\end{lemma}
\begin{proof}
Suppose there is a vanishing linear combination $\sum_{i=1}^c \lambda_i B_i = 0$ of the rows $B_i$ of $B$. If the edge $e_j$ has endpoints at $v_a$ and $v_b$, the entries in the $(n+j)$-th column $B^{n+j}$ are 1 at the rows $a$ and $b$ and zero everywhere else: therefore we get $\lambda_a+\lambda_b=0$. If we apply this argument to the three edges of a triangle we get
$$\lambda_a+\lambda_b=0, \quad \lambda_b+\lambda_c=0, \quad \lambda_c+\lambda_a = 0$$
which implies $\lambda_a=\lambda_b=\lambda_c=0$. Then all coefficients $\lambda_a$ vanish.
\end{proof}

\begin{lemma} \label{B:surj:lemma}
The $c$ rows in $B$ generate $\ker(A^\intercal)$. 
\end{lemma}
\begin{proof}
Since $BA=0$ implies $A^\intercal B^\intercal=0$, we already know that the rows in $B$ are contained in $\ker(A^\intercal)$, and we must prove that they generate it. We pick a generic horizontal vector $q=(q_1, \ldots, q_{2n})\in \matR^{2n}$ with $qA=0$ and we need to prove that $q$ is generated by the rows $B_i$ of $B$. 

Recall that the columns $A^j, A^{j+n}, A^{j+2n}$ of $A$ correspond to the edge pairs of type 1, 2, and 3 of the tetrahedron $\Delta_j$. We label the vertices of $\Delta_j$ as $a,b,c,d$, the edges of $\Delta_j$ as pairs like $ab$, and the corresponding vertices and edges in $\calT$ as $j(a)$ and $j'(ab)$, interpreted as numbers in $\{1,\ldots,c\}$ and $\{1,\ldots, n\}$ respectively. Moreover we set $j(ab) = j'(ab)+n$. We suppose that $ab, ac, ad$ are of type 1, 2, 3.

The conditions $qA^j = qA^{j+n} = qA^{j+2n}=0$ say that
$$
\begin{aligned}
qA^j & = q_j + q_{j(ab)} + q_{j(cd)} = 0, \\
qA^{j+n} & = q_j + q_{j(ac)} + q_{j(bd)} = 0, \\
qA^{j+2n} & = q_j + q_{j(ad)} + q_{j(bc)} = 0. 
\end{aligned}
$$
This implies that
\begin{equation} \label{qq:eqn}
q_{j(ab)} + q_{j(cd)} = q_{j(ac)} + q_{j(bd)} = q_{j(ad)} + q_{j(bc)} = - q_j.
\end{equation}
We now prove that $q= \sum \lambda_i B_i$ where the coefficients $\lambda_i$ are defined as follows. Consider an ideal triangle in $\Delta_j$ with vertices $a,b,c$ and define
$$\lambda_{j(a)} = \frac{q_{j(ab)} + q_{j(ca)} - q_{j(bc)}}2.$$
We check that $\lambda_{j(a)}$ does not depend on the ideal triangle: the triangle in $\Delta_j$ with vertices $a,b,d$ gives 
$$\lambda_{j(a)} = \frac{q_{j(ab)} + q_{j(da)} - q_{j(bd)}}2 = \frac{q_{j(ab)} + q_{j(ca)} - q_{j(bc)}}2$$
using (\ref{qq:eqn}). Two triangles incident to the ideal vertex $v_{j(a)}$ are connected by a path of adjacent tetrahedra and hence $\lambda_{j(a)}$ is well-defined.

Finally, we need to check that indeed $q = \sum \lambda_iB_i$. We consider the edge $j(ab)$ and recall that the $(n+j(ab))$-th column of $B$ is everywhere zero except at rows $j(a)$ and $j(b)$; hence the $(n+j(ab))$-th entry of $\sum \lambda_iB_i$ is $\lambda_{j(a)} + \lambda_{j(b)}$, which by definition equals $q_{j(ab)}$, so we are done.

Analogously, the $j$-th column of $B$ is everywhere zero except at the rows $j(a),j(b),j(c),j(d)$ and hence the $j$-th entry of $\sum \lambda_iB_i$ is $-(\lambda_{j(a)} + \lambda_{j(b)} + \lambda_{j(c)} + \lambda_{j(d)})$ which equals $-q_{j(ab)} - q_{j(cd)} = q_j$. The proof is complete.
\end{proof}

These three lemmas imply Proposition \ref{AB:exact:prop}. In particular we get the following:

\begin{cor} \label{rk:A:cor}
We have $\rk A = 2n-c$.
\end{cor}

\subsection{A symplectic form} \label{symplectic:form:subsection}
We define an alternating bilinear form $\omega$ on $\matR^{3n}$ as follows: for every $i=1,\ldots, n$ consider the form
$$\omega_i = \begin{pmatrix}
0 & 1 & -1 \\
-1 & 0 & 1 \\
1 & -1 & 0
\end{pmatrix}
$$
on the subspace $V_i$ generated by $e_i, e_{i+n}$, and $e_{i+2n}$, and define $\omega = \oplus_i \omega_i$. Strictly speaking, the form $\omega$ is not symplectic because it is degenerate: its radical is generated by all vectors of type $e_i+e_{i+n} + e_{i+2n}$, that is by the first $n$ rows of the matrix $A$. The form $\omega$ becomes symplectic after quotienting $\matR^3$ by the radical.

\begin{prop} \label{rows:lagrangian:prop}
The rows of $A$ generate a lagrangian subspace.
\end{prop}
\begin{proof}
In other words, we need to prove that $\omega(A_j, A_k)=0$ for any rows $A_j, A_k$ of $A$. The first $n$ rows lie in the radical, hence we are left to check that $\omega(A_{n+j}, A_{n+k})=0$ for all $0 < j < k \leqslant n$ corresponding to the $j$-th and $k$-th edge of the triangulation. 

Each tetrahedron $\Delta_i$ contributes to the number $\omega(A_{n+j}, A_{n+k})$ as follows: if two edges contained in some face $f$ of $\Delta_i$ correspond to the edges $j$ and $k$ of the triangulation, they contribute with a $\pm 1$; the face $f$ is also contained in an adjacent tetrahedron where they contribute with opposite sign, and hence everything sums to zero.
\end{proof}

\subsection{The solution space} \label{solution:space:subsection}
We now employ the combinatorial results of the previous sections to study the consistency equations for $\calT$. 
We assign the variables $z_j$, $\frac 1{1-z_j}$, and $\frac{z_j-1}{z_j}$ to the edges of type 1, 2, and 3 in each tetrahedron $\Delta_j$. These three variables lie in the upper half-plane 
$$\matC_+ = \big\{\Im z > 0\big\}$$
which we see as the subset $\matC_+ \subset \widetilde{\matC^*}$ consisting of all $\rho e^{i\theta}$ with $0<\theta < \pi$.
The consistency equation at the edge $e_i$ is the following equation in $\widetilde{\matC^*}$:
\begin{equation} \label{gi:eqn}
g_i(z) = \prod_{j=1}^n z_j^{A_{n+i,j}} \left(\frac 1{1-z_j}\right)^{A_{n+i,n+j}} 
\left(\frac{z_j-1}{z_j}\right)^{A_{n+i,2n+j}} = e^{2\pi i}.
\end{equation}
The elements $\frac{1}{1-z_j}, \frac{z_j-1}{z_j}\in\matC_{+}\subset \widetilde \matC^*$ depend holomorphically on $z_j\in\matC_{+}$.
We have defined a complex function
$$g \colon \matC_{+}^n \longrightarrow \big(\widetilde{\matC^*}\big)^n$$
by setting $g=(g_1,\ldots, g_n)$. The \emph{solution space} is the set
$$\Def(M,\calT) = g^{-1}(e^{2\pi i},\ldots, e^{2\pi i}).$$
This space consists precisely of the solutions to the consistency equations of $\calT$. 
Each solution $z\in \Def(M,\calT)$ gives a hyperbolic metric on $M$, which ``deforms'' as $z$ varies: this justifies the use of the symbol ``$\Def$''. We want to prove that $\Def(M,\calT)$ is a complex manifold of dimension $c$. 

\subsection{Complex differentials} \label{complex:diff:subsection}
We use the biholomorphism $\log\colon \widetilde{\matC^*} \to \matC$ and define the function $G(z) = \log g(z)$ for all $z\in \matC_+^n$. Now 
$$\Def(M,\calT) = G^{-1}(2\pi i, \ldots, 2\pi i).$$ 
We may write
$$G_i(z) = \sum_{j=1}^n A_{i,j}' \log z_j + A'_{i,n+j}\log(1-z_j) + N_i\pi i$$
where $A'$ is the $n \times 2n$ integral matrix defined by
$$A'_{i,j} = A_{n+i,j} - A_{n+i,2n+j}, \quad A'_{i,n+j} = -A_{n+i,n+j} + A_{n+i,2n+j}$$ 
and $N_i$ is some fixed integer that depends on the way we interpret $\log(1-z_j)$. 
The complex differential of $G$ at $z$ splits as a composition of two linear maps
$$ dG_{ z}\colon \matC^{n} \stackrel{Z}{\longrightarrow} \matC^{2n} \stackrel{A'}{\longrightarrow} \matC^n$$
where
$$Z = \begin{pmatrix}
\frac 1{z_1} & 0 & \ldots & 0 \\
0 & \frac 1{z_2} & \ldots & 0 \\
\vdots & \vdots & \ddots & \vdots\\
0 & 0 & \ldots & \frac 1{z_n} \\
\frac {1}{z_1-1} & 0 & \ldots & 0 \\
0 & \frac {1}{z_2-1} & \ldots & 0 \\
\vdots & \vdots & \ddots & \vdots\\
0 & 0 & \ldots & \frac {1}{z_n-1} 
\end{pmatrix}
$$
is obtained by differentiating $\log z_j$ and $\log (1-z_j)$.
Let now
$$B'\colon \matC^n \longrightarrow \matC^c$$
be the linear map defined by the right $(c\times n)$-sub-matrix $B'$ of $B$, that is $$B'_{i,j} = B_{i.n+j}.$$

\begin{prop} \label{AB2:exact:prop}
The following sequence is exact:
$$ \matC^{2n} \stackrel{A'}{\longrightarrow} \matC^{n} \stackrel{B'}{\longrightarrow} \matC^c \longrightarrow 0.$$
\end{prop}
\begin{proof}
Lemma \ref{B:surj:lemma} shows that $B$ is surjective, and its proof actually shows that $B'$ is surjective. Recall that $A$ decomposes into six square matrices, and via Gauss moves on columns we get 
\begin{equation} \label{GaussA:eqn}
A=\begin{pmatrix}
I & I & I \\
A^1 & A^2 & A^3 
\end{pmatrix}
\longrightarrow
\begin{pmatrix}
0 & 0 & I \\
A^1-A^3 & -A^2 + A^3 & A^3
\end{pmatrix}
=
\begin{pmatrix}
0 & I \\
A' & A^3
\end{pmatrix}.
\end{equation}
Now $BA=0 \Rightarrow B'A'=0$, and $\rk A = 2n-c \Rightarrow \rk A' = n-c$.
\end{proof}

This shows in particular that the image of $G$ lies in a fixed affine subspace parallel to $\ker B'$ of complex codimension $c$. We now consider the standard symplectic form $\omega = \matr 0{I}{-I}0$ on $\matR^{2n}$. 
\begin{prop} \label{rows:lagrangian2:prop}
The rows of $A'$ generate a lagrangian subspace.
\end{prop}
\begin{proof}
We know that the rows of $A$ generate a lagrangian subspace with respect to the degenerate form on $\matR^{3n}$ by Proposition \ref{rows:lagrangian:prop}. The Gauss moves (\ref{GaussA:eqn}) imply the assertion.
\end{proof}

\begin{prop} \label{Choi:prop}
For every $z\in \matC_+^n$ the following sequence is exact:
$$ \matC^n \stackrel{dG_{z}} \longrightarrow \matC^n \stackrel {B'} \longrightarrow \matC^c \longrightarrow 0.$$
\end{prop}
\begin{proof}
We have $dG_{z} = A'Z$. Proposition \ref{AB2:exact:prop} gives an exact sequence
$$ \matC^{2n} \stackrel{A'}{\longrightarrow} \matC^{n} \stackrel{B'}{\longrightarrow} \matC^c \longrightarrow 0.$$
We only need to prove that $\Img (A'Z) = \Img A'$, that is $\rk (A'Z) = n-c$. We prefer to consider the transpose matrix $M= (A'Z)^\intercal =  Z^\intercal (A')^\intercal$. We have
$$ Z^\intercal = \begin{pmatrix}
\frac 1{z_1} & 0 & \ldots & 0 & \frac {1}{z_1-1} & 0 & \ldots & 0 \\
0 & \frac 1{z_2} & \ldots & 0 & 0 & \frac {1}{z_2-1} & \ldots & 0 \\
\vdots & \vdots & \ddots & \vdots & \vdots & \vdots & \ddots & \vdots \\
0 & 0 & \ldots & \frac 1{z_n} & 0 & 0 & \ldots & \frac {1}{z_n-1} 
\end{pmatrix}.
$$
The map $Z^\intercal\colon \matC^{2n} \to \matC^n$ is obviously not injective, but its restriction 
$$Z^\intercal |_{\matR^{2n}}\colon \matR^{2n} \longrightarrow \matC^n$$ 
is an isomorphism of real spaces because $\frac 1{z_i}$ and $\frac {1}{z_i-1}$ are independent over $\matR$ (we use here that $z_i\not\in \matR$).  

We define a symplectic form $\Omega$ on the real vector space $\matC^n$ by pushing-forward $\omega$ along $Z^\intercal |_{\matR^{2n}}$, that is we set $\Omega(Z^\intercal v,\, Z^\intercal w) = \omega(v,w)$. By construction $\Omega$ is a sum of symplectic forms on the coordinate complex lines of $\matC^n$. Note that a symplectic form on $\matC$ is unique up to rescaling.

Since $\Img z_i >0$, the complex numbers $\frac{1}{z_i}$ and $\frac{1}{z_i-1}$ form a negatively-oriented $\matR$-basis of $\matC$, and this implies that $\Omega (v,iv) < 0$ for every non-zero vector $v\in \matC^n$. In particular $\langle \cdot, \cdot \rangle = -\Omega(\cdot, i \cdot)$ is a positive-definite scalar product that induces a norm $\| \cdot \|$ on $\matC^n$. 

On each coordinate complex line of $\matC^n$, the form $\Omega$ is just a rescaling of the standard symplectic form and hence $\langle,\rangle$ is a rescaling of the standard scalar product. In particular we get $\Omega(iv,iw) = \Omega(v,w)$ and $\langle iv, iw \rangle  = \langle v,w \rangle$.

Finally, we suppose by contradiction that $\rk M < n-c$. Since $Z^\intercal |_{\matR^{2n}}$ is injective, the real subspace $M(\matR^{n})$ has real dimension $n-c$.
Since $\rk M < n-c$, there are vectors $v_1,v_2\in \matR^{n}$ such that $M(v_1+iv_2) = 0$ while $Mv_1 \neq 0$ and $Mv_2 \neq 0$. Proposition \ref{rows:lagrangian2:prop} gives  
$$\Omega(Mv_1, Mv_2)= \omega ((A')^\intercal v_1,(A')^\intercal v_2) = 0$$
and therefore
\begin{align*}
0 = \|M(v_1 + iv_2)\|^2 & = \|Mv_1\|^2 + \|iM v_2\|^2 + 2\langle M v_1, iM v_2 \rangle \\
& = \|Mv_1\|^2 + \|Mv_2\|^2 -2\Omega(M v_1, -Mv_2) \\
& = \|Mv_1\|^2 + \|Mv_2\|^2
\end{align*}
implies $Mv_1 = Mv_2=0$, a contradiction.
\end{proof}

\begin{cor} \label{smooth:cor}
The solution space $\Def(M,\calT) \subset \matC_+^n$ is either empty or a complex submanifold of dimension $c$.
\end{cor}
\begin{proof}
The image of $G\colon \matC_+^n \to \matC^n$ lies entirely in a $(n-c)$-dimensional complex affine space $S$ parallel to $\ker B' \subset \matC^n$ and by the previous proposition the holomorphic map $G\colon \matC_+^n \to S$ is a submersion onto its image. Therefore the counterimage of any point is a either empty or is a complex submanifold of dimension $n-(n-c)=c$.
\end{proof}

We now want to define a nice holomorphic parametrisation for the manifold $\Def(M,\calT)$, and to this purpose we need to investigate further the combinatorial properties of $\calT$ by looking more closely at its ideal vertices.

\subsection{Boundary curves} \label{boundary:curves:subsection}
Recall that $M=\interior N$ and the ideal vertices $v_1,\ldots, v_c$ of $\calT$ correspond to the boundary tori $T_1,\ldots, T_c$ of $N$, which are triangulated by $\calT$. Recall also that there are $n$ tetrahedra $\Delta_1,\ldots, \Delta_n$ in $\calT$ and each $\Delta_i$ has three edge pairs labeled with $i, n+i,$ and $2n+i$.

\begin{figure}
\begin{center}
\includegraphics[width = 12.5 cm] {\iftoggle{BW}{path_holonomy-BW}{path_holonomy}} 
\nota{We mark with a dot all the interior angles that $\gamma$ encounters on its right as in the picture: each dot determines an edge in some tetrahedron $\Delta_j$ and hence an edge pair (left). We also consider the ``degenerate path'' consisting of a single point ``run counterclockwise'', and recover the last $n$ rows of $A$ in this way (right).}
\label{path_holonomy:fig}
\end{center}
\end{figure}

Every simplicial closed oriented curve $\gamma$ in a triangulated torus $T_i$ determines an integer vector $v_\gamma \in \matR^{3n}$ that counts the (right-)incidences between $\gamma$ and the $3n$ edge pairs of $\calT$. More precisely, each interior angle of a triangle in $T_i$ determines an edge in some tetrahedron $\Delta_j$ and hence an edge pair: we define $(v_\gamma)_j$ to be the number of times $\gamma$ encounters the $j$-th edge pair on its right, as in Figure \ref{path_holonomy:fig}-(left).

We defined an alternating degenerate form $\omega$ on $\matR^{3n}$ and we now want to extend Proposition \ref{rows:lagrangian:prop}. We denote by $\gamma \cdot \gamma'$ the algebraic intersection of $\gamma$ and $\gamma'$ in $H_1(\partial M, \matZ)$.

\begin{prop} \label{A:gamma:gamma':prop}
Let $\gamma, \gamma'$ be simplicial closed curves in $\partial N$. We have:
\begin{itemize}
\item $\omega(v_\gamma, v_{\gamma'}) = 2 \gamma \cdot \gamma'$,
\item $\omega(v_\gamma, A_i) = 0 $ for every row $A_i$ of $A$.
\end{itemize}
\end{prop}
\begin{proof}
The row $A_i$ with $i\leqslant n$ lies in the radical of $\omega$, and $A_{n+i}$ records the incidences of the $i$-th edge $e_i$ of $\calT$ with the $3n$ edge pairs. We may write $A_{n+i} = v_{\gamma'}$ where $\gamma'$ is a ``degenerate'' closed path consisting of a single point (one of the endpoints of $e_i$), see Fig~\ref{path_holonomy:fig}-(right). Now the second assertion is just a special case of the first.

\begin{figure}
\begin{center}
\includegraphics[width = 7 cm] {\iftoggle{BW}{two_paths_holonomy-BW}{two_paths_holonomy}} 
\nota{We draw the paths $\gamma$, $\gamma'$ and their dots in \iftoggle{BW}{dark grey}{blue} and \iftoggle{BW}{grey}{red}.}
\label{two_paths_holonomy:fig}
\end{center}
\end{figure}

\begin{figure}
\begin{center}
\includegraphics[width = 10 cm] {\iftoggle{BW}{curves_dots_contribution-BW}{curves_dots_contribution}} 
\nota{Two dots in a triangle contribute to $\omega(\gamma, \gamma')$ with $+1$, $0$, or $-1$ according to their mutual position.}
\label{curves_dots_contribution:fig}
\end{center}
\end{figure}

\begin{figure}
\begin{center}
\includegraphics[width = 12.5 cm] {\iftoggle{BW}{curves_near_contribution-BW}{curves_near_contribution}} 
\nota{Two portions of \iftoggle{BW}{dark grey}{blue} and \iftoggle{BW}{grey}{red} paths that intersect the same triangle contribute to $\omega(\gamma, \gamma')$ as shown here. The contribution depends on their mutual configuration, and changes by a sign if we interchange the \iftoggle{BW}{dark grey}{blue}/\iftoggle{BW}{grey}{red} colours. If one of the two portions contains at least two edges (as in bottom-right) it determines a vector in the radical, and the contribution is always zero.}
\label{curves_near_contribution:fig}
\end{center}
\end{figure}

We draw the two paths $\gamma$ and $\gamma'$ in \iftoggle{BW}{dark grey}{blue} and \iftoggle{BW}{grey}{red} as in Figure \ref{two_paths_holonomy:fig}, marking with dots the angles they encounter on their right. A couple of \iftoggle{BW}{dark grey}{blue} and \iftoggle{BW}{grey}{red} dots in the same triangle contribute to $\omega(\gamma, \gamma')$ according to their mutual position as shown in Figure \ref{curves_dots_contribution:fig}. Therefore portions of \iftoggle{BW}{dark grey}{blue} and \iftoggle{BW}{grey}{red} paths intersecting the same triangle contribute to $\omega(\gamma, \gamma')$ as shown in Figure \ref{curves_near_contribution:fig}.

For every contribution $+1$ as in Figure \ref{curves_near_contribution:fig}-(top-left), there is an analogous but opposite contribution $-1$ provided by the other triangle adjacent to the edge joining the coloured endpoints: these all cancel. We are left with a $+1$ every time the \iftoggle{BW}{grey}{red} path exits leftward from an intersection with the \iftoggle{BW}{dark grey}{blue} path, a $-1$ every time it enters from the left, and the opposite contributions with the \iftoggle{BW}{dark grey}{blue}/\iftoggle{BW}{grey}{red} colours interchanged. We easily get a total of $2\gamma \cdot \gamma'$.

\begin{figure}
\begin{center}
\includegraphics[width = 12.5 cm] {\iftoggle{BW}{curves_far_contribution-BW}{curves_far_contribution}} 
\nota{When $\gamma$ and $\gamma'$ encounter in different truncated triangles of the same tetrahedron $\Delta_i$, they also contribute to $\omega(\gamma,\gamma')$.}
\label{curves_far_contribution:fig}
\end{center}
\end{figure}

The proof is not yet complete, because every tetrahedron $\Delta_i$ determines four triangles in $\partial N$ and some contributions to $\omega(\gamma, \gamma')$ arise also from paths intersecting distinct triangles of the same $\Delta_i$: we need to prove that these contributions all cancel.

Up to symmetries of the tetrahedron and up to interchanging the colours, the possible non-zero contributions are all shown in Figure \ref{curves_far_contribution:fig}. Each configuration as in the first row cancels with an opposite one produced by the tetrahedron adjacent to the grey face. The contributions in the second row arise when $\gamma$ and $\gamma'$ cross the endpoints of an edge $e$ on opposite sides, and one checks that the contributions of all the cycle of tetrahedra incident to $e$ sum to zero (exercise).
\end{proof}

\subsection{Holonomy parameters} \label{holonomy:subsection}
We now fix a non-trivial oriented simple closed curve $\gamma_i$ on each $T_i$ and a representation of $\gamma_i$ as a simplicial (possibly non-injective) path in the triangulated $T_i$.

\begin{figure}
\begin{center}
\includegraphics[width = 6 cm] {\iftoggle{BW}{ABmatrices_extended-BW}{ABmatrices_extended}} 
\nota{The enlarged matrix $\bar A$.}
\label{ABmatrices_extended:fig}
\end{center}
\end{figure}

We enlarge the incident $2n \times 3n$ matrix $A$ to a $(2n+c)\times 3n$ matrix $\bar A$ by adding the vectors $v_{\gamma_1}, \ldots, v_{\gamma_c}$ at its bottom, see Figure \ref{ABmatrices_extended:fig}.

\begin{prop} \label{A:bar:lagrangian:prop}
We have $\rk \bar A = 2n$. The rows of $\bar A$ generate a maximal lagrangian subspace in $\matR^{3n}$.
\end{prop}
\begin{proof}
Let $\mu_i \subset T_i$ be a simplicial closed curve with $\gamma_i\cdot \mu_i = 1$. We already know that the first $2n$ rows of $\bar A$ generate a lagrangian subspace of dimension $2n-c$, and to conclude it suffices to use Proposition \ref{A:gamma:gamma':prop} and deduce that $v_{\mu_i}$ is $\omega$-orthogonal to all the rows of $\bar A$ except $\bar A_{2n+i}=v_{\gamma_i}$. Therefore $\rk \bar A = \rk A +c = 2n$.
\end{proof}

Every point $z\in \Def(M,\calT)$ determines a hyperbolic structure on $M$ and a homomorphism
$$\tilde\mu_z \colon \pi_1(T_i) \to \widetilde{\matC^*}$$
for each $i=1,\ldots, c$, see Definition \ref{tilde:mu:defn}. We define
$$h_i(z) = \tilde\mu_z(\gamma_i) = e^{-|\gamma_i| \pi i } \prod_{j=1}^n z_j^{\bar A_{2n+i,j}} \left(\frac 1{1-z_j}\right)^{\bar A_{2n+i,n+j}} 
\left(\frac{z_j-1}{z_j}\right)^{\bar A_{2n+i,2n+j}} $$
and get a holomorphic function
$$h\colon \matC_+^n \longrightarrow \big(\widetilde{\matC^*}\big)^c.$$
We want to prove that $h$ furnishes a nice parametrisation for the $c$-dimensional space $\Def(M,\calT)$. The letter $h$ stands for ``holonomy'': recall that $h_i$ is the holonomy of $\gamma_i$, see Exercise \ref{tilde:rho:mu:ex}.

As in Section \ref{complex:diff:subsection}, we define $H(z) = \log h(z)$. 
We get a map
$$G \times H \colon \matC_+^n \longrightarrow \matC^n \times \matC^c.$$
As above, the differential splits as a composition of linear maps
$$ d(G\times H)_{z}\colon \matC^{n} \stackrel{Z}{\longrightarrow} \matC^{2n} \stackrel{\bar A'}{\longrightarrow} \matC^n \times \matC^c$$
where $\bar A'$ is constructed from $\bar A$ exactly as $A'$ from $A$.

\begin{prop} \label{Choi2:prop}
The differential $d(G\times H)_{z}$ is injective $\forall z \in \matC_+^n$.
\end{prop}
\begin{proof}
As above, Proposition \ref{A:bar:lagrangian:prop} implies that the rows of $\bar A'$ form a lagrangian subspace of dimension $n$ and we conclude with the same proof as Proposition \ref{Choi:prop}.
\end{proof}

\begin{cor}
The restriction
$$h|_{\Def(M,\calT)}\colon \Def(M,\calT) \longrightarrow \big(\widetilde{\matC^*}\big)^c$$ is a local biholomorphism.
\end{cor}

A \emph{complete solution} $z\in\Def(M,\calT)$ is one that determines a complete metric on $M$, that is such that $h(z)=(1,\ldots,1)$.

\begin{cor} \label{complete:isolated:cor}
Every complete solution in $\Def(M,\calT)$ is isolated.
\end{cor}

The local biholomorphism $h$ depends on the isotopy class of the curves $\gamma_i$ chosen: different choices give different parametrisations. Note however that the condition that $h_i(z)=1$ corresponds geometrically to the completeness of $T_i$ and hence does not depend on $\gamma_i$.

\section{Proof of the theorem}
We now conclude the proof of Theorem \ref{DF:teo}.
Through all this section, we let $M= \interior N$ be a complete orientable finite-volume cusped hyperbolic three-manifold. We have $\partial N = T_1 \cup \ldots \cup T_c$ for some $c\geqslant 1$. We also fix a positive pair $m_i, l_i$ of generators for $\pi_1(T_i)$ at each $T_i$. Recall that each $T_i$ has Euclidean structure, defined up to rescaling.

\subsection{Near the complete solution} \label{near:complete:subsection}
We now suppose that $M$ has a geometric decomposition $\calT$ into ideal hyperbolic tetrahedra, in other words the space $\Def(M,\calT)$ is non-empty and contains a complete solution $z^0$. We want to study $\Def(M,\calT)$ near $z^0$. 

We defined in Section \ref{generalized:DF:subsection} the generalised Dehn filling invariant of a solution $z\in\Def(M,\calT)$ at $T_i$ to be $\infty$ if $\mu(m_i) = \mu(l_i)=1$, or the unique point $(p,q)\in \matR^2$ such that $\tilde\mu(m_i)^p \tilde\mu(l_i)^q = e^{2\pi i}$ otherwise. We denote this invariant by $d_i(z)\in S^2 = \matR^2 \cup \{\infty\}$ and get a map 
$$d\colon \Def(M,\calT) \longrightarrow \underbrace{S^2 \times \ldots \times S^2}_c$$
such that $d(z^0) = (\infty,\ldots,\infty)$.

\begin{prop} \label{d:local:homeo:prop}
The map $d$ is a local homeomorphism at $z^0$.
\end{prop}
\begin{proof}
We set
$$u_i = \log \tilde\mu(m_i), \qquad v_i = \log \tilde \mu(l_i).$$
The sets $(u_1,\ldots,u_c)$ and $(v_1,\ldots, v_c)$ are both biholomorphic coordinates on $\Def(M,\calT)$ which vanish at $z^0$. Moreover $v_i$ vanishes precisely where $u_i$ does (on solutions that are complete at $T_i)$. 

Recall from Exercise \ref{rho:derivative:ex} that $\mu = \rho'$ is the derivative of the holonomy $\rho$ at $T_i$. When $u_i = v_i = 0$ the solution is complete at $T_i$ and up to conjugation we have
$$\rho(m_i)\colon z \mapsto z +1, \qquad \rho(l_i)\colon z \mapsto z + w_i$$
for some $w_i\in \matC$ with $\Im w_i>0$ which determines the cusp shape of $T_i$. When $u_i, v_i \neq 0$ the solution is incomplete and
$$\rho(m_i) \colon z \mapsto az, \qquad \rho(l_i) \colon z \mapsto bz$$
for some $a,b\neq 0$.
When $u_i, v_i \to 0$, a fundamental domain of $\Img\rho$ with vertices at $x,ax, bx, abx$ tends up to rescaling to one with vertices at $x, x+1, x+w_i, x+1+w_i$. Since $a,b \to 1$ we get
$$v_i = \log b \sim b-1 \sim w_i(a-1) \sim w_i \log a = w_iu_i$$
where $\sim$ means first-order equality. Therefore $\frac{v_i}{u_i} = w_i$ at the complete solution $v=u=0$. We deduce in particular that the function $\frac{v_i}{u_i}$ is analytic on $\Def(M,\calT)$.
In what follows we only use that $\Im w_i > 0$.

We consider the following map, defined whenever $(p_i,q_i) \neq (0,0)$ for all $i$:
$$\varphi\big((p_1,q_1), \ldots, (p_c, q_c)\big) = 2\pi i \left( \frac 1 {p_1+q_1w_1}, \ldots, \frac 1{p_c + q_cw_c}\right).$$
We extend it continuously to the case $(p_i, q_i) = \infty$ by setting the $i$-th component of the image to be zero.
By composing $\Phi = \varphi \circ d$ we get a map $\Phi\colon \Def(M,\calT) \to \matC^c$ that sends the complete solution $u=0$ to $(0,\ldots,0)$. We prove that $\Phi$ is differentiable and $d\Phi_0$ is invertible: this will conclude the proof. When $u_j, v_j \neq 0$ we have
$$p_ju_j + q_jv_j = 2\pi i$$
and Proposition \ref{pq:explicit:prop} gives
$$p_j = - 2 \pi \frac{\Re v_j}{\Im (\bar u_j v_j)}, \qquad q_j = 2 \pi \frac{\Re u_j}{\Im (\bar u_jv_j)}.$$
Near $u=v=0$ we have $v_j \sim w_ju_j$ and hence
$$p_j \sim - 2 \pi \frac{\Re w_j \Re u_j - \Im w_j \Im u_j}{|u_j|^2\Im w_j}, \qquad
q_j \sim 2 \pi \frac{\Re u_j}{|u_j|^2\Im w_j}
$$
and therefore
\begin{align*}
\varphi_j & = 2\pi i \frac 1{p_j + q_jw_j} \sim i\frac{|u_j|^2 \Im w_j}{-\Re w_j \Re u_j + \Im w_j \Im u_j + w_j \Re u_j} \\
& =
i\frac{|u_j|^2 \Im w_j}{\Im w_j \Im u_j + i\Im w_j \Re u_j} 
=
i\frac{|u_j|^2 \Im w_j}{i \bar u_j \Im w_j} = u_j
\end{align*}
giving $d\Phi_0 = \id$. 
\end{proof}

If $M$ has a geometric decomposition $\calT$ into ideal tetrahedra, we can easily conclude the proof of Theorem \ref{DF:teo} as follows: the geometric ideal triangulation $\calT$ determines a complete solution $z^0 \in \Def(M,\calT)$, and by Proposition \ref{d:local:homeo:prop} the image of the generalised Dehn filling map $d\colon \Def(M,\calT) \to S^2\times \ldots \times S^2$ is an open set $U$ containing $(\infty, \ldots, \infty)$. The discussion in Section \ref{generalized:DF:subsection} finishes the proof.

However, we still do not know if every complete cusped finite-volume $M$ decomposes into ideal tetrahedra! If that were the case, we would be done. Unfortunately we are not able to prove this, so more effort is needed to complete the proof.

\subsection{Triangulations with flat tetrahedra}
As we said, we are not able to prove that every cusped finite-volume $M$ decomposes into ideal tetrahedra. The best that we can do is to show that $M$ decomposes into ideal tetrahedra that are either ordinary or ``flat'', that is degenerated to an ideal quadrilateral. We show this using the Epstein--Penner canonical decomposition.

The Epstein-Penner canonical decomposition (see Theorem \ref{EP:teo}) partitions $M$ into some ideal polyhedra $P_1,\ldots, P_k$, glued isometrically along their faces. Can we further subdivide this partition to get a geometric triangulation of $M$? We start by triangulating each $P_i$ individually.

\begin{prop} \label{poly:tetra:prop}
Every ideal polyhedron $P\subset \matH^3$ subdivides into ideal tetrahedra.
\end{prop}
\begin{proof}
We can use the Klein model, so that ideal polyhedra in $\matH^3$ correspond to Euclidean polyhedra with vertices in $S^2$. We fix a vertex $v\in P$, and a vertex $v_i$ in every face $f_i$ of $P$ not incident to $P$. We triangulate each $f_i$ by coning along $v_i$, and then we triangulate $P$ by coning all the triangulated faces $f_i$ along $v$.
\end{proof}

\begin{figure}
\begin{center}
\includegraphics[width = 11 cm] {\iftoggle{BW}{layered-BW}{layered}} 
\nota{The triangulations of two isometric faces that do not match (left) because they are obtained by coning along distinct vertices $v\neq v'$ (centre). We can add some additional ``flat'' tetrahedra to correct this (right).}
\label{layered:fig}
\end{center}
\end{figure}

The problem is, that we are not able to guarantee that the triangulations of the polyhedra $P_1,\ldots, P_k$ match at each pair of glued faces. The typical situation is shown in Figure \ref{layered:fig}-(left): two isometric $n$-gons (with $n\geqslant 4$) are triangulated by coning on distinct vertices $v\neq v'$.

There is one thing that we can do, however: we can insert some \emph{flat} tetrahedra between the two $n$-gons that connect the two mismatching triangulations as in Figure \ref{layered:fig}-(right), and obtain a partially-flat ideal triangulation of $M$. The flat tetrahedra can be inserted as follows: if $v$ and $v'$ share an edge we insert $n-3$ flat tetrahedra as suggested by Figure \ref{layered:fig}-(right); otherwise, the diagonal joining $v$ and $v'$ is common to both triangulations and we proceed inductively on each of the two sub-polygons cut by this diagonal. 

The resulting ideal triangulation $\calT$ contains $n$ tetrahedra $\Delta_1,\ldots, \Delta_n$; each $\Delta_i$ is either \emph{fat} (it lies inside some polyhedron $P_j$) or \emph{flat} (it connects two faces of the $P_j$). We now fix arbitrarily a pair of opposite edges in each $\Delta_i$ and denote by $z_i\in \matC$ the modulus of $\Delta_i$ at these edges. If $\Delta_i$ is fat, then $\Im z_i>0$ as usual; if $\Delta_i$ is flat, the modulus $z_i$ still makes sense as a number in $\matR\setminus\{0,1\}$:
the vertex triangle of Figure \ref{tetraedro2:fig}-(left) is a degenerated flat triangle with three vertices $0,1,z_i$. The flat tetrahedron has three real moduli $z_i, \frac{z_i-1}{z_i}, \frac 1{1-z_i}$; two arguments are $0$ and one is $\pi$.

We set $z^0 = (z_1,\ldots, z_n)$. By construction, the point $z^0$ satisfies the consistency equations (\ref{gi:eqn}). 

\subsection{Perturbations}
We now let $U$ be a small contractible neighbourhood of $z^0$ in $(\matC \setminus \{0,1\})^n$. Note that $U$ necessarily contains some points $z$ with $\Im z_i <0$ for some $i$. The consistency equations (\ref{gi:eqn}) make sense in $U$, and we define $\Def_U \subset U$ to be the set of solutions. We still have $\Def_U = G^{-1}(2\pi i, \ldots, 2\pi i)$.

Does a solution $z\in \Def_U$ determine a hyperbolic structure on $M$? Intuitively, each modulus $z_j$ with $\Im z_j=0$ should represent a flat tetrahedron, each $z_j$ with $\Im z_j <0$ a negatively oriented fat tetrahedron, and positively and negatively oriented tetrahedra should overlap; although appealing, it is however hard (and sometimes impossible) to translate this idea into a concrete definition in general.

Our aim is to prove Theorem \ref{DF:teo} by perturbing $z^0$, so we only need here to furnish a geometric interpretation for the solutions $z$ that lie in a sufficiently small neighbourhood $U$. We do this as follows. 

If $z$ is sufficiently close to $z^0$, all the fat tetrahedra that triangulate a single polyhedron $P_i$ stay positively fat, and since they satisfy the consistency equations along the internal edges of $P_i$ they glue to form a possibly non-convex polyhedron $P_i(z)\subset \matH^3$. The original faces of $P_i$ are now bent in $P_i(z)$ along their diagonals.

We now consider the originally flat tetrahedra joining two faces of two polyhedra $P_i(z)$ and $P_j(z)$. We assign all of them to one of the two polyhedra, say $P_i(z)$, and forget about $P_j(z)$. These tetrahedra are attached to some $n$-gonal face $f$ of $P_i(z)$, bent along its diagonals. The effect of each attached tetrahedron $\Delta_j$ is to flip a diagonal, and we modify $P_i(z)$ accordingly by replacing the two ideal triangles in $f$ incident to the old diagonal with the two triangles sharing the new one. Geometrically, this corresponds to adding $\Delta_j$ if $\Im z_i >0$, or to cut $\Delta_j$ away from $P_j(z)$ is $\Im z_j <0$. At the end of the process, we get a new polyhedron $P_j'(z)$ with the same vertices as $P_j(z)$, whose faces are bent along different diagonals.

The new possibly non-convex polyhedra $P_j'(z)$ glue isometrically along their paired faces, that are now bent along the same diagonals. The consistency equations are satisfied along these diagonals, hence we get a hyperbolic structure for $M$. We summarise our discoveries:

\begin{prop}
If $U$ is sufficiently small, every solution $z \in \Def_U$ represents naturally a hyperbolic structure on $M$.
\end{prop}

\subsection{Local smoothness}
We have seen that every solution $z$ that is sufficiently close to $z^0$ determines a hyperbolic metric on $M$. We now need to extend Corollary \ref{smooth:cor}. 

\begin{prop} \label{smooth:flat:prop}
If $U$ is sufficiently small then $\Def_U$ is a complex manifold of dimension $c$.
\end{prop}
\begin{proof}
As in Corollary \ref{smooth:cor}, it suffices to prove that $dG_{z^0}\colon V^n \to \ker B'$ has maximum rank and is hence invertible.

The only place where we used $\Im z_j >0$ is the proof of Proposition \ref{Choi:prop}, and we now modify it so that it applies also to our partially flat solution $z^0$. We have $M = Z^\intercal (A')^\intercal$ and we need to prove that $\rk M = \rk A' = n-c$.

We needed $z_i \not \in \matR$ to prove that $Z^\intercal|_{\matR^{2n}}$ is injective: this is not the case here, since $\frac 1{z_i}$ and $\frac 1{z_i-1}$ are not $\matR$-independent when $\Delta_i$ is flat. Now we have
$$\ker Z^\intercal |_{\matR^{2n}} = \Span \{ z_ie_i +(1-z_i)e_{i+n}\ | \ \Delta_i\ {\rm is\ flat} \}.$$
For $\matK=\matR, \matC$ we write 
$$\matK^{2n} = \matK^{2n}_{\rm fat} \oplus \matK^{2n}_{\rm flat}, \qquad \matK^n = \matK^n_{\rm fat} \oplus \matK^n_{\rm flat}$$
where $\matK^{2n}_{\rm fat}$ ($\matK^{2n}_{\rm flat}$) is generated by the $e_i, e_{i+n}$ such that $\Delta_i$ is fat (flat), and similarly $\matK^{n}_{\rm fat}$ ($\matK^{n}_{\rm flat}$) is generated by the $e_i$ such that $\Delta_i$ is fat (flat).

The following restriction of $Z^\intercal$ is still an isomorphism:
$$Z^\intercal|_{\matR^{2n}_{\rm fat}}\colon \matR^{2n}_{\rm fat}
\longrightarrow \matC^n_{\rm fat}.$$
We use it to push-forward $\omega$ to a symplectic form $\Omega$ on $\matC^n_{\rm fat}$, and we extend it trivially to an alternating bilinear form $\Omega$ on $\matC^n$ with radical $\matC^n_{\rm flat}$.

It is still true that $M(\matR^{2n})$ has real dimension $n-c$, because $(A')^\intercal(\matR^{2n})$ intersects $\ker Z^\intercal$ trivially. To prove this, suppose that a vector $v'\in\matR^{2n}$ that is $\matR$-generated by the rows of $A'$ lies in the kernel of $Z^\intercal$, that is:
\begin{itemize}
\item $(v_i', v_{i+n}')=(0,0)$ for every fat $\Delta_i$,
\item $(v_i',v_{i+n}') = k_i(z_i, 1-z_i)$ for some $k_i \in \matR$, for every flat $\Delta_i$. 
\end{itemize}
We define $v \in \matR^{3n}$ by setting
$$v_i = v'_i, \quad v_{n+i} = - v'_{n+i}, \quad, v_{2n+i} = 0$$
for all $i=0,\ldots, n-1$. We note that each triple $(v_i, v_{n+i}, v_{2n+i})$ is either $(0,0,0)$, or it contains three pairwise distinct real numbers. The vector $v' \in\matR^{2n}$ is generated by the $n$ rows of $A'$ $\Longleftrightarrow$ the vector $v\in\matR^{3n}$ is generated by the $n$ rows of $A$ (exercise). 

We interpret $v$ as a colouring of the $3n$ edge pairs in $\calT$. The colouring $v$ is everywhere zero on all fat tetrahedra, and on each flat tetrahedron it is either $(0,0,0)$ or consists of three pairwise distinct numbers. We must prove that it is zero everywhere.

The vector $v$ is contained in the lagrangian subspace generated by the rows $A_1,\ldots, A_{2n}$ of $A$. In particular we have $\omega(v,A_{n+i})=0$ for all $i$. Every row $A_{n+i}$ corresponds to an edge $e$ of $\calT$, and this equality says that the colours on all tetrahedra incident to $e$ of the edge pairs not containing $e$ sum to zero, provided that they are counted with signs in each tetrahedron. This implies the following: if all such tetrahedra except one $\Delta_j$ have colour zero on all their edge pairs, then two edge pairs on $\Delta_j$ have the same colours, and by what just said on $v$ this implies that all colours are zero also on $\Delta_j$. 

Consider a sequence $\Delta_1, \ldots, \Delta_k$ of flat tetrahedra that connects two faces of two polyhedra $P_i$ and $P_j$ of the Epstein-Penner decomposition. There is one edge of $\Delta_1$ that is incident only to fat tetrahedra except $\Delta_1$, and the above discussion implies that all the colours on $\Delta_1$ are zero; by induction we deduce that all the flat $\Delta_i$ have zero colours. Therefore $v=0$, as required.

We have proved that $\dim M(\matR^{2n}) = n-c$, and we conclude as in the proof of Proposition \ref{Choi:prop}. The form $\Omega$ defines a positive semidefinite scalar product $\langle, \rangle$ on $\matC^n$ with radical $\matC^n_{\rm flat}$. If $\rk M < n-c$, there are two vectors $v_1,v_2\in \matR^n$ such that $M(v_1+iv_2) = 0$ while $Mv_1 \neq 0$ and $Mv_2 \neq 0$. We get $\|Mv_1\| = \|Mv_2\| = 0$ as in the proof of Proposition \ref{Choi:prop}, and this implies that $Mv_1, Mv_2$ both lie in the radical $\matC^n_{\rm flat}$.

If $Mv_1, Mv_2 \in \matC^n_{\rm flat}$ then both $Mv_1$ and $Mv_2$ are real since $\frac 1{z_i}, \frac{1}{z_i-1} \in \matR$ for every flat $\Delta_i$. Now $0 = M(v_1+iv_2) = Mv_1 + iMv_2$ implies that $Mv_1 = Mv_2 = 0$, a contradiction.
\end{proof}

\subsection{Conclusion of the proof}
The proof of Theorem \ref{DF:teo} finishes like in the geometric triangulation case.
All the discussion of Section \ref{holonomy:subsection} applies here: the holonomies are defined also in this context and may be used to parametrize the complex manifold $U(z^0)$. The strategy of Section \ref{near:complete:subsection} is still valid: the generalised Dehn filling invariant map $d$ sends $U(z^0)$ to an open neighbourhood of $(\infty,\ldots, \infty)$ in $S^2\times \ldots \times S^2$ and the discussion in Section \ref{generalized:DF:subsection} finishes the proof.
We have completed the proof of Theorem \ref{DF:teo}.

\subsection{Mostow--Prasad rigidity}
We can use the Dehn filling Theorem \ref{DF:teo} to extend Mostow Rigidity Theorem \ref{Mostow:teo} from the closed to the cusped case. The cusped case was proved by Prasad in 1973, hence the theorem is known as \emph{Mostow--Prasad Rigidity Theorem}.\index{Mostow--Prasad rigidity theorem}

Let $N$ a be compact orientable three-manifold bounded by some $c\geqslant 0$ tori, and let $M=\interior N$.

\begin{teo}[Mostow--Prasad rigidity]
The manifold $M$ admits at most one finite-volume complete hyperbolic metric up to isometries homotopic to the identity.
\end{teo}
\begin{proof}
The $c=0$ case is Theorem \ref{Mostow:teo}, so we suppose $c>1$. Roughly, the metric on $M$ is unique because it is the limit of hyperbolic metrics on its closed Dehn fillings, that are unique by Theorem \ref{Mostow:teo}. (We are using the hyperbolic Dehn filling theorem as a trick to extend Mostow's rigidity from the closed to the cusped case: this is certainly not Prasad's original proof, that was published before any hyperbolic Dehn filling was studied, and more importantly it applies to all dimensions $n\geqslant 3$.)

More precisely, let $M$ have two finite-volume hyperbolic metrics $m_1$ and $m_2$. The Dehn filling Theorem holds for both, and yields sequences $M_i^j$ of closed hyperbolic three-manifolds obtained from $M$ by Dehn filling with parameter $s^j$ such that  $s^j \to (\infty, \ldots, \infty)$, for $i=1,2$.

Mostow theorem furnishes an isometry $M_1^j \to M_2^j$ homotopic to the identity. The core closed geodesics of the filling solid tori are unique in their homotopy classes, hence by removing them we get isometric incomplete metrics. Since $m_1$ and $m_2$ are limits of these, they are also isometric.

This latter sentence can be made rigorous as follows: let $\calT$ be a partially flat ideal triangulation for $M$ in the metric $m_1$. Pick an ideal tetrahedron $\Delta$ in $\calT$ and note that it can be straightened uniquely also in the incomplete metrics for $M$ close to the complete ones, just by extending the construction in the proof of Proposition \ref{straighten:incomplete:prop} from dimension two to three (an orientation of the core geodesics is needed to define the straightening unambiguously). The modulus $z\in \matC_{>0} \cup (\matR \setminus \{0,1\})$ of the straightened $\Delta$ in $M$ is intrinsically determined as the limit of the moduli of the straightenings of $\Delta$ in the $M_i^j$, hence it is the same for both metrics $m_1$ and $m_2$. It follows that $m_2$ has the same geometric ideal decomposition of $m_1$, and hence $m_1=m_2$ after a homotopy. 
\end{proof}

\section{Volumes}
We have discovered that ``most'' Dehn fillings of a cusped finite-volume hyperbolic manifold are hyperbolic, and we now discuss their volumes.

Let $M$ be a complete orientable finite-volume cusped hyperbolic three-manifold. We say that  a Dehn filling of $M$ is \emph{non-trivial} if it fills at least one boundary torus; in other words, it is determined by some Dehn filling parameter $s$ different from $(\infty, \ldots, \infty)$, after fixing some arbitrary homology bases for the boundary tori. We prove in this section the following theorem. 

\begin{teo} \label{vol:vol:teo}
If a non-trivial Dehn filling $M^{\rm fill}$ of $M$ admits a complete finite-volume hyperbolic metric, then
$$\Vol(M^{\rm fill}) < \Vol(M).$$
If $s^i$ is any sequence of non-trivial Dehn filling parameters converging to $(\infty,\ldots, \infty)$, and $M^i$ is the manifold obtained by filling $M$ along $s^i$, then $M^i$ is eventually hyperbolic and
$$\Vol (M^i) \nearrow \Vol(M).$$
The cores of the filling solid tori are eventually simple closed geodesics and their lengths tend to zero.
\end{teo}

We already know that if $s^i \to (\infty, \ldots, \infty)$ then $\Vol(M^i) \to \Vol(M)$ and the cores of the solid tori are closed geodesics whose length tend to zero: this follows from the construction of the metric on $M^i$ via geometric ideal triangulations. In some sense (that we will not make precise), the metric on $M^i$ converges to that of $M$, and the core geodesics become so short that they disappear in the limit opening the filling solid torus to a cusp.

The only new fact to prove is that the volume strictly decreases under \emph{any} Dehn filling $M^{\rm fill}$, not only those that are close to $(\infty, \ldots, \infty)$: this is rather subtle, because the hyperbolic metric on $M^{\rm fill}$ might not be obtained by completing ideal triangulations of $M$. (For instance, the core tori in $M^{\rm fill}$ may not be isotopic to closed geodesics.)

We prove the theorem as follows. We first consider the case where $M$ has an ideal geometric triangulation $\calT$: in that case $\Def(M,\calT)$ is a nice smooth complex manifold and we investigate the volume function 
$$\Vol\colon \Def(M,\calT) \longrightarrow \matR_{>0}$$
that computes the volume of each hyperbolic structure on $M$. We show that the volume function makes sense in a bigger convex set containing $\Def(M,\calT)$ and is also easier to study there. This bigger convex set is the set of \emph{angle structures} on $\calT$. We prove that the volume function is strictly concave there, and the complete solution is a global maximum.\index{angle structure} 

We conclude with a couple of finite-covering tricks to pass from local (near $(\infty, \ldots, \infty)$) to global, and to deal with more general manifolds $M$ that do not decompose into hyperbolic ideal tetrahedra.

\subsection{Angle structures}
Let $M= \interior{N}$ where $N$ is an orientable compact three-manifold bounded by $c>0$ tori, equipped with an oriented ideal triangulation $\calT$ with $n$ tetrahedra $\Delta_1,\ldots, \Delta_n$ and $n$ edges $e_1,\ldots, e_n$. We identify orientation-preservingly each $\Delta_i$ with the one shown in Figure \ref{tetrahedron_labels:fig}, and we label the edge pairs of type 1, 2, 3, with the numbers $i$, $i+n$, and $i+2n$ respectively.

An \emph{angle structure} $\theta$ on $M$ is the assignment of an \emph{angle} $0< \theta_j < \pi$ to the $j$-th edge pair in $\calT$, with the following requirements:
\begin{enumerate}
\item at each tetrahedron $\Delta_i$ we have $\theta_i + \theta_{i+n} + \theta_{i+2n} = \pi$,
\item at each edge $e_i$, the sum of the incident angles is $2\pi$.
\end{enumerate}
The set of angle structures is thus a subset $\calA \subset \matR_{>0}^{3n}$. It is the intersection of the affine subspace in $\matR^{3n}$ determined by the conditions (1) and (2) with the cone $\matR_{>0}^{3n}$, so in particular it is convex. By definition this affine subspace is the solution space of the linear system
$$A\theta = \begin{pmatrix} \pi \\ 2\pi \end{pmatrix}$$
where $A$ is the incidence $2n \times 3n$ defined in Section \ref{incidence:matrices:subsection} and depicted in Figure \ref{ABmatrices:fig}-(left), and $\theta = (\theta_1, \ldots, \theta_{3n})\in \matR^{3n}$.

What is the geometric meaning of an angle structure $\theta$ on $\calT$? Thanks to condition (1), every $\Delta_i$ is realised by a unique ideal hyperbolic tetrahedron with dihedral angles $\theta_i, \theta_{i+n}$, and $\theta_{i+2n}$ whose complex modulus is
\begin{equation} \label{theta:z:eqn}
z_i = \frac{\sin \theta_{i+n}}{\sin \theta_{i+2n}}e^{i\theta_i}.
\end{equation}
Angle structures are then in natural 1-1 correspondence with the realisations of the tetrahedra $\Delta_1,\dots, \Delta_n$ as ideal hyperbolic tetrahedra, with the mild condition (2) that their dihedral angles around each edge $e_i$ must sum to $2\pi$. 

Condition (2) is only ``half'' of the consistency equation at each $e_i$, which requires the product of the complex moduli to be 1: condition (2) alone does \emph{not} guarantee that by gluing these tetrahedra we get a hyperbolic structure on $M$: some ``shearing'' may arise when you turn around an edge $e_i$.

With this geometric description, the deformation space $\Def(M,\calT)$ is naturally a smooth submanifold of $\calA$. 

\subsection{Dimensions}
We will henceforth suppose that $\Def(M,\calT)$ is non-empty and contains a complete solution $z^0$, that is $\calT$ can be realised as a geometric ideal triangulation. This implies in particular that $\calA$ is also non-empty and has dimension $\dim\ker A = 3n - (2n-c) =n+c$ by Corollary \ref{rk:A:cor}. The space $\Def(M,\calT)$ is a real smooth submanifold of dimension $2c$. 

The tangent space $T_\theta\calA$ at any point $\theta\in\calA$ is of course $\ker A$. We now construct an explicit set of generators for $\ker A$. Recall from Section \ref{boundary:curves:subsection} that every simplicial closed oriented curve $\gamma$ in a triangulated boundary torus $T_i$ defines an integer vector $v_\gamma \in \matR^{3n}$. 

Recall the alternating form $\omega$ on $\matR^{3n}$ from Section \ref{symplectic:form:subsection}: we consider it as a matrix and write $v^* = \omega v$ for all $v\in \matR^{3n}$, so that $\omega(v,w) = \langle v, w^* \rangle$ where $\langle, \rangle$ is the Euclidean scalar product. We denote by $A_i$ the $i$-th row of $A$, considered as a vector in $\matR^{3n}$.

We fix two arbitrary generators $m_i, l_i$ of $\pi_1(T_i)$ for all $i=1,\ldots, c$.

\begin{prop} \label{redundant:generators:prop}
The tangent space $T_\theta \calA = \ker A$ is generated by the (not independent) $n+2c$ vectors
$A_{n+1}^*,\ldots, A_{2n}^*, v_{m_1}^*, v_{l_1}^*, \ldots, v_{m_c}^*, v_{l_c}^*$.
\end{prop}
\begin{proof}
Propositions \ref{rows:lagrangian:prop} and \ref{A:gamma:gamma':prop} say that the rows of $A$ are $\omega$-orthogonal to themselves and to $v_\gamma$ for any $\gamma$, hence they are $\langle, \rangle$-orthogonal to all vectors $A_i^*$ and $v_\gamma^*$, in other words $A_i^*, v_\gamma^* \in \ker A$.

To prove that the $n+2c$ listed vectors generate $\ker A$ it suffices to show that they span a subspace of dimension $n+c$. Proposition \ref{A:bar:lagrangian:prop} says that if we add the rows $v_{m_1}, \ldots, v_{m_c}$ to $A$ we get a rank-$2n$ matrix. The same proof shows that by further adding the rows $v_{l_1},\ldots, v_{l_c}$ we get a rank-$(2n+c)$ matrix. 

The vectors $A_1,\ldots, A_{2n},v_{m_1}, v_{l_1}, \ldots, v_{m_c}, v_{l_c}$  span a dimension-$(2n+c)$ space. Since $\dim \ker \omega = n$, the vectors $A_1^*,\ldots, A_{2n}^*,v_{m_1}^*, v_{l_1}^*, \ldots, v_{m_c}^*, v_{l_c}^*$ span a space of dimension at least $2n+c-n = n+c$. Since $A_1^* = \ldots = A_n^* = 0$, we are done.
\end{proof}

\subsection{The volume function}
The volume function on $\Def(M,\calT)$ extends naturally to a function
$$\Vol\colon \calA \longrightarrow \matR_{>0}$$
which assigns to each angle structure $\theta$ the sum of the volumes of the ideal hyperbolic realisations of $\Delta_1,\ldots, \Delta_n$ determined by $\theta$. 
Recall from Section \ref{volume:tetrahedra:section} that the volume of $\Delta_i$ is
$$\Vol(\Delta_i) = \Lambda(\theta_i) + \Lambda(\theta_{i+n}) + \Lambda(\theta_{i+2n})$$
where $\Lambda$ is the Lobachevsky function
$$\Lambda(\alpha) = -\int_0^\alpha \log |2\sin t| dt.$$
We therefore get a simple-looking formula
$$\Vol(\theta) = \sum_{j=1}^{3n} \Lambda(\theta_j)$$
which implies immediately the following important fact.

\begin{prop} \label{derivative:Vol:prop}
For every $\theta \in \calA$ and $v\in T_\theta  \calA$ we have 
$$\frac{\partial \Vol}{\partial v} = \sum_{j=1}^{3n} -v_j \log \sin \theta_j, \qquad 
\frac{\partial^2 \Vol}{\partial v^2} < 0.$$
\end{prop}
\begin{proof}
We consider a single addendum $\Lambda(\theta_j)$ of $\Vol(\theta)$. We have
$$\frac{\partial \Lambda(\theta_j)}{\partial v} = 
\Lambda'(\theta_j) \frac{\partial \theta_j}{\partial v}
= - \log |2\sin \theta_j| v_j = -v_j\log 2 - v_j\log\sin \theta_j. $$
Since $v\in T_\theta \calA$ we have $v_i+v_{i+n}+v_{i+2n}=0$ for all $i$ and hence $\sum_j v_j\sin 2=0$. This proves the first equality. To estimate the second derivative, we may suppose up to symmetries that $\theta_i, \theta_{i+n} < \frac \pi 2$ and we get
\begin{align*}
-\frac {\partial^2 \Vol(\Delta_i)}{\partial v^2} 
& = v_i^2 \cot \theta_i + v_{i+n}^2 \cot \theta_{i+n} + v_{i+2n}^2 \cot \theta_{i+2n} \\
& = v_i^2 \cot \theta_i + v_{i+n}^2 \cot \theta_{i+n} + 
(v_i+v_{i+n})^2 \frac{1-\cot \theta_i\cot\theta_{i+n}}{\cot \theta_i + \cot \theta_{i+n}} \\
& = \frac{(v_i+v_{i+n})^2 + (v_i\cot\theta_i - v_{i+n}\cot\theta_{i+n})^2}
{\cot \theta_i + \cot \theta_{i+n}} > 0.
\end{align*}
Indeed the denominator is positive since $\theta_i, \theta_{i+n} < \frac \pi 2$, and if the numerator is zero we get $v_i=-v_{i+n}$ and hence $\cot \theta_i = -\cot \theta_{i+n}$, a contradiction.
\end{proof}

We have discovered that $\Vol$ is a smooth strictly concave function on $\calA$. 

\subsection{The complete solution}
We now remember that $\Def(M,\calT)$ contains a complete solution $z^0$ by assumption.

\begin{prop} \label{Vol:maximum:prop}
The function $\Vol\colon \calA \to \matR_{>0}$ has a unique global maximum at the complete solution $z^0$.
\end{prop}
\begin{proof}
The complete solution $z^0$ has an angle structure $\theta$. Proposition \ref{derivative:Vol:prop} says that
$$\frac{\partial \Vol}{\partial v} = \sum_{j=1}^{3n} -v_j \log \sin \theta_j$$
for all $v\in T_{\theta}\calA$. At every edge $e_i$, the consistency equation is of type $\sum_j \log w_j = 2\pi i$ where $w_1,\ldots, w_k$ are the moduli of the incident tetrahedra. The imaginary part of this equation is satisfied at every point of $\calA$, so we look at the real part and using (\ref{theta:z:eqn}) we discover that 
$$0 = \Re \sum_j \log w_j = \sum_j \log |w_j| = \sum_j \big(\log \sin \alpha_{j}^1 - \log \sin \alpha_j^2\big) $$
where $\alpha_j^1$ and $\alpha_j^2$ are the other two angles in the tetrahedron with modulus $w_j$, in counterclockwise order. The latter expression is in fact equivalent to 
$$-\sum_{j=1}^{3n} (A_{n+i}^*)_j\log\sin \theta_j = \frac{\partial \Vol}{\partial A_{n+i}^*}$$
which is hence zero for all $i=1,\ldots, n$.

We make some analogous considerations for the completeness equations. Since $z^0$ is complete, for every simplicial boundary curve $\gamma$ we have $\sum_j \log w_j = |\gamma|\pi i$ where $w_1,\ldots, w_k$ are the moduli encountered by $\gamma$ at its right, and the real part of this equation implies that $\frac{\partial \Vol}{\partial v_\gamma^*} = 0$.

Proposition \ref{redundant:generators:prop} shows that vectors of type $A_{n+i}^*$ and $v_\gamma^*$ generate $T_{\theta}\calA$, hence the gradient of $\Vol$ vanishes at $\theta$. Since $\Vol$ is strictly concave, this critical point is the unique global maximum.
\end{proof}

\begin{cor}
The space $\Def(M,\calT)$ contains at most one complete solution.
\end{cor}

\subsection{Volumes and Dehn filling}
We can now prove Theorem \ref{vol:vol:teo} on manifolds having a geometric ideal triangulation. We start with a ``local'' version which also considers generalised Dehn fillings, see Section \ref{generalised:subsection}.

\begin{lemma} \label{U:lemma}
If $M$ has a geometric ideal triangulation, there is a neighbourhood $U$ of $(\infty,\ldots, \infty)$ such that $\Vol(M^{\rm fill})< \Vol(M)$ for every non-trivial generalised Dehn filling parameter $s\in U$. Moreover $\Vol(M^{\rm fill}) \nearrow \Vol(M)$ as $s$ tends to $(\infty, \ldots, \infty)$.
\end{lemma}
\begin{proof}
The Dehn Filling Theorem \ref{DF:teo} says that all Dehn fillings in an open neighbourhood $U$ of $(\infty,\ldots,\infty)$ are hyperbolic and constructed by completing solutions in $\Def(M,\calT)$, and they have smaller volume than $M$ in virtue of Proposition \ref{Vol:maximum:prop}. (The completion adds a circle to $M$ which does not contribute to the volume.)

The volume of $M^{\rm fill}$ depends smoothly on $s\in U$ and has a global maximum at $(\infty,\ldots,\infty)$. In particular, when $s\to (\infty, \ldots, \infty)$ we get $\Vol(M^{\rm fill}) \to \Vol(M)$ from below.
\end{proof}

We now need to drop the open set $U$ and prove that volume decreases under \emph{any} Dehn filling of $M$. 

\begin{lemma} \label{strict:everywhere:lemma}
Let $M$ have a geometric ideal triangulation. If a non-trivial Dehn filling $M^{\rm fill}$ admits a complete finite-volume hyperbolic metric, then $\Vol(M^{\rm fill}) < \Vol(M)$.
\end{lemma}
\begin{proof}
Let $s=(s_1,\ldots, s_c)$ be the Dehn filling parameter giving $M^{\rm fill}$ and let
$U$ be the open neighbourhood of $(\infty,\ldots, \infty)$ furnished by Lemma \ref{U:lemma}. For sufficiently big integer $k>0$ we have $ks \in U$ and hence $M^{\rm fill}$ also admits a hyperbolic metric with cone angles $\frac{2\pi}k$ on the core geodesics of the filling solid tori, with volume strictly smaller than $\Vol(M)$. We may interpret this metric as a hyperbolic orbifold $O$, and we have $\Vol(O)<\Vol(M)$.

If $M^{\rm fill}$ is closed Proposition \ref{Vol:M:O:prop} gives $\Vol(M) > \Vol(O) \geqslant v_3\|M^{\rm fill}\| = \Vol(M^{\rm fill})$ and we are done.

If $M^{\rm fill}$ is not closed, we can find some closed fillings of $M^{\rm fill}$ with volume arbitrarily close to $\Vol(O)$ that can be interpreted as orbifolds, and get $\Vol(O)\geqslant \Vol(M^{\rm fill})$ as a limit.
\end{proof}

The proof of Lemma \ref{strict:everywhere:lemma} shows also the following.
\begin{lemma} \label{closed:decreases:lemma}
Let $M$ have a geometric ideal triangulation. We have $v_3\|M^{\rm fill}\| < \Vol(M)$ for every closed Dehn filling $M^{\rm fill}$ of $M$.
\end{lemma}

We prove a lemma that will be useful below.
\begin{lemma} 
Let $M$ have a geometric ideal triangulation. If $M_i$ is a hyperbolic Dehn filling of $M$ with parameter $s^i$ and $\Vol(M_i) \to \Vol(M)$, then $s^i \to (\infty,\ldots, \infty)$.
\end{lemma}
\begin{proof}
We may suppose that the $M_i$ are closed, since every sequence of Dehn fillings can be approximated by a sequence of closed ones. If $s^i$ lies definitely in the open set $U$ furnished by Theorem \ref{DF:teo}, the assertion follows because the volume function there has $\Vol(M)$ as a unique maximum (thanks to Proposition \ref{Vol:maximum:prop}). 

Otherwise, there is a fixed big $k\in\matN$ such that $ks^i\in U$ for all $i$, but $ks^i$ does not converge to $(\infty,\ldots, \infty)$. Therefore $M_i$ admits a hyperbolic orbifold structure $O_i$ with $\Vol(M) - \varepsilon > \Vol(O_i) $ and Proposition \ref{Vol:M:O:prop} gives $\Vol(M) -\varepsilon > \Vol(O) \geqslant v_3\|M_i\| = \Vol(M_i)$ for all $i$. 
\end{proof}

Again, the proof shows also the following.
\begin{lemma} \label{sse:converge:lemma}
Let $M$ have a geometric ideal triangulation. If $M_i$ is a closed Dehn filling of $M$ with parameter $s^i$ and $v_3\|M_i\| \to \Vol(M)$, then $s^i \to (\infty, \ldots, \infty)$.
\end{lemma}

To conclude the proof of Theorem \ref{vol:vol:teo} we only need to consider the unlucky case where $M$ has no geodesic ideal triangulations. We prove that geometric triangulations exist on some finite cover of $M$, and this will suffice.

\subsection{Geometric triangulations exist virtually}
We prove here the following.

\begin{teo}
Every cusped complete finite-volume three-manifold $M$ has a finite-sheeted cover $\tilde M$ that has a geometric ideal triangulation $\calT$.
\end{teo}
\begin{proof}
The Epstein-Penner decomposition subdivides $M$ into ideal polyhedra $P_1,\ldots, P_k$. If distinct vertices of $P_i$ lie in distinct cusps of $M$ for all $i$, then the decomposition can be easily subdivided into an ideal triangulation: we order the cusps, so the vertices of each $P_i$ inherit an ordering which is preserved along the matching faces; we subdivide each $P_i$ as in the proof of Proposition \ref{poly:tetra:prop} by coning on the smallest vertex at each step: the triangulations on the paired faces match and we get a geometric ideal triangulation for $M$ itself.

We now show that some cover $\tilde M$ has a decomposition of this type. The decomposition of $M$ into $P_1,\ldots, P_k$ lifts in any degree-$d$ cover $\tilde M$ to a decomposition of $\tilde M$ into $P_1^i, \ldots, P_k^i$ with $i=1,\ldots, d$. We now prove that there is a $\tilde M$ such that distinct vertices of $P_j^i$ lie in distinct cusps, for all $i,j$.

To do so, we prove that for every pair of vertices $v \neq v'$ of some $P_i$ there is a covering where the lifts $\tilde v, \tilde v'$ in $P_i^j$ for some $j$ lie in distinct cusps. This is enough to conclude: pick a \emph{regular} covering $\tilde M$ that covers all these finitely many coverings; since its desk transformation group acts transitively on lifts, all distinct vertices on every $P_i^j$ lie in distinct cusps.

We have $M=\matH^3/_\Gamma$. We fix an arbitrary lift $\tilde P_i \subset \matH^3$ of $P_i$. Two vertices $v \neq v'$ of $P_i$ lie in the same cusp $T$ $\Leftrightarrow$ there is a $\gamma \in \Gamma$ such that $\gamma(\tilde v) = \tilde v'$, where $\tilde v, \tilde v' \in \tilde P_i$ are the lifts of $v,v'$.
The subgroup $\Stab_\Gamma(\tilde v) = \pi_1(T)$ is separable by Corollary \ref{separable:cor} and hence there is a finite-index $H< \Gamma$ that contains $\Stab_\Gamma(\tilde v)$ but avoids $\gamma$. 

No element $\varphi\in H$ is such that $\varphi(\tilde v) = \tilde v'$, otherwise we would have $\gamma^{-1}\varphi \in \Stab_\Gamma(v) < H$ and hence $\gamma \in H$, a contradiction. We have constructed a finite-cover $\matH^3/_H$ where in some lift $P_i^j$ of $P_i$ the two lifted vertices $\tilde v,\tilde v'$ lie in distinct cusps.
\end{proof}

\subsection{Conclusion of the proof}
We can finally prove Theorem \ref{vol:vol:teo}.

\dimo{vol:vol:teo}
Some degree-$d$ cover $\tilde M$ of $M$ has an ideal geometric triangulation, hence the theorem holds for $\tilde M$. 

For every Dehn filling $M^{\rm fill}$ of $M$, there is a Dehn filling $\tilde M^{\rm fill}$ of $\tilde M$ such that the cover $\tilde M \to M$ extends to a degree-$d$ map $\tilde M^{\rm fill} \to M^{\rm fill}$. 

If $M^{\rm fill}$ is closed hyperbolic, we have 
$$d\Vol(M) =\Vol(\tilde M) > v_3\| (\tilde M^{\rm fill})\|  \geqslant dv_3\|M^{\rm fill}\| =d\Vol(M^{\rm fill}) $$
using Lemma \ref{closed:decreases:lemma}, hence $\Vol(M) > \Vol(M^{\rm fill})$. 

If $M^{\rm fill}$ is cusped hyperbolic, it can be approximated by closed hyperbolic fillings $M^{\rm fill}_i$, and there are closed fillings $\tilde M^{\rm fill}_i$ of $\tilde M^{\rm fill}$ with degree-$d$ maps $\tilde M^{\rm fill}_i \to M^{\rm fill}_i$.
By Lemma \ref{sse:converge:lemma} there is an $\varepsilon > 0$ such that 
$$\Vol(\tilde M) - \varepsilon > v_3\|\tilde M^{\rm fill}_i\| \geqslant d v_3\|M^{\rm fill}_i\| = d\Vol(M^{\rm fill}_i)$$
for all $i$, hence $d\Vol(M) - \varepsilon = \Vol(\tilde M) - \varepsilon \geqslant d\Vol(M^{\rm fill})$. 
\finedimo

\subsection{Bounded volume}
How many hyperbolic manifolds are there with bounded volume? In general, infinitely many. 

We have just discovered that there are infinitely many hyperbolic manifolds with volume smaller than $2v_3 = 2.0298832128 \ldots$ because the figure-eight knot complement $M$ has $\Vol(M) = 2v_3$ and every hyperbolic Dehn filling of $M$ has volume smaller than $\Vol(M)$. Moreover, the volumes of the Dehn fillings of $M$ form an indiscrete set that tends to $2v_3$ from below.

Theorem \ref{vol:vol:teo} shows that this is the typical situation: the hyperbolic Dehn filling produces infinite sequences of manifolds with bounded volume. Despite this variety of manifolds, we can still control topologically the hyperbolic manifolds having bounded volume, at least in principle.

\begin{teo} \label{V:teo}
For every $V>0$ there is a compact 3-manifold $N$ bounded by tori, such that every complete finite-volume orientable hyperbolic 3-manifold $M$ with $\Vol(M)<V$ is diffeomorphic to the interior of some Dehn filling of $N$.
\end{teo}
\begin{proof}
This is a consequence of the thick-thin decomposition. Let $\varepsilon>0$ be a fixed Margulis constant: every complete finite-volume $M$ decomposes along disjoint embedded tori as $M=M^{\rm thick} \cup M^{\rm thin}$ where $M^{\rm thick}$ is compact and has injectivity radius $>\varepsilon$, and $M^{\rm thin}$ consists of truncated cusps and tube neighbourhoods of short geodesics. In particular $M$ is a Dehn filling of $M^{\rm thick}$.

We now show that there are only finitely many possible diffeomorphism types for $M^{\rm thick}$ with volume $<V$. This concludes the proof: by Proposition \ref{LW:bordo:prop}
every such diffeomorphism type is realised from a link $L\cup C\subset S^3$ by surgerying along $L$ and drilling along $C$; we take a disjoint union of all these links in $S^3$, and the complement of the resulting link in $S^3$ is our $N$.

Let $X\subset M^{\rm thick}$ be a \emph{maximal $\frac{\varepsilon}2$-net}, that is a finite set of points that stay at pairwise distance $\geqslant \frac{\varepsilon}2$, such that every other point in $M^{\rm thick}$ stays at distance $< \frac {\varepsilon}2$ from $X$ ($X$ exists because $M^{\rm thick}$ is compact). The $\frac{\varepsilon}4$-balls centred at the points in $X$ are embedded and disjoint: hence $X$ has cardinality at most $C(V) = V/\Vol\big(B(x, \frac{\varepsilon}4)\big)$.

The set $X$ determines a cellularisation of $M^{\rm thick}$ as follows. Lift $X$ to $\tilde X\subset \matH^3$, take the Voronoi tessellation (see Section \ref{Voronoi:subsection}) determined by $\tilde X$ and project it back to $M$: this furnishes a subdivision of $M^{\rm thick}$ into polyhedra. Since $X$ is maximal, every polyhedron has at most $C'(\varepsilon)$ faces with $C'$ depending only on $\varepsilon$. Hence there are finitely many combinatorial types of polyhedra, depending only on $\varepsilon$. With at most $C(V)$ of them we obtain only finitely many manifolds. (The proof of Proposition \ref{X:triangulation:prop} applies to polyhedra and shows that the combinatorial subdivision determines the smooth manifold.)
\end{proof}

\subsection{Volumes of hyperbolic three-manifolds}
We finish this chapter by furnishing some qualitative information on the set of volumes of hyperbolic three-manifolds. We start with the following.

\begin{prop}
Let $M_i$ be a sequence of pairwise non-diffeomorphic complete orientable hyperbolic three-manifolds of uniformly bounded volume. After passing to a subsequence, we may suppose that there is a hyperbolic manifold $M$ such that each $M_i$ is a Dehn filling of $M$ with parameter $s^i\to (\infty, \ldots, \infty)$ . 
\end{prop}
\begin{proof}[Proof using geometrisation]
We know from Theorem \ref{V:teo} that there is a manifold $N$ such that each $M_i$ is a Dehn filling of $N$ with parameter $s^i$. After passing to a subsequence, we may suppose that at each boundary torus $T_j$ of $N$ the slope $s^i_j$ is either constant or goes to infinity. Let $M$ be obtained from $N$ by Dehn filling the slopes that stay constant.

If $M$ is hyperbolic, we are done. Unfortunately, this may not be the case, and we use geometrisation to solve this annoying possibility (the original proof of J\o rgensen and Thurston \cite{Th} does not need geometrisation and follows from a discussion on the geometric convergence of hyperbolic manifolds: we employ geometrisation only for simplicity). If $M$ is not prime, it contains some essential spheres that always bound balls in the filled $M_i$ and hence we may cut $M$ along them and cup off with balls (reducing the components of $\partial M$). Now $M$ is prime and its geometric decomposition must contain some hyperbolic block, otherwise $M$ would be a graph manifold and no Dehn filling of $M$ would be hyperbolic (it would be a graph manifold again). The tori of the decomposition of $M$ become compressible in $M_i$, so up to passing to a subsequence we may suppose that each $M_i$ is a Dehn filling of one fixed hyperbolic block of $M$. If $s_i$ does not tend to $(\infty, \ldots, \infty)$ in this new setting, we restart from the beginning and proceed by induction on the number of boundary components of $M$.
\end{proof}

\begin{cor}
The set of volumes of all complete orientable hyperbolic three-manifolds is well-ordered. For every value $V$ there are only finitely many such manifolds with volume $V$.
\end{cor}
\begin{proof}
Suppose that there is a sequence of manifolds $M_i$ with $\Vol(M_1) > \Vol (M_2) > \ldots $ By the previous proposition these are all obtained by Dehn filling the same hyperbolic $M$ with $s^i \to (\infty,\ldots,\infty)$, and Theorem \ref{vol:vol:teo} says they should converge to $\Vol(M)$ from below, a contradiction. The same argument show that there are finitely many manifolds of any fixed volume.
\end{proof}

\begin{figure}
\begin{center}
\includegraphics[width = 6 cm] {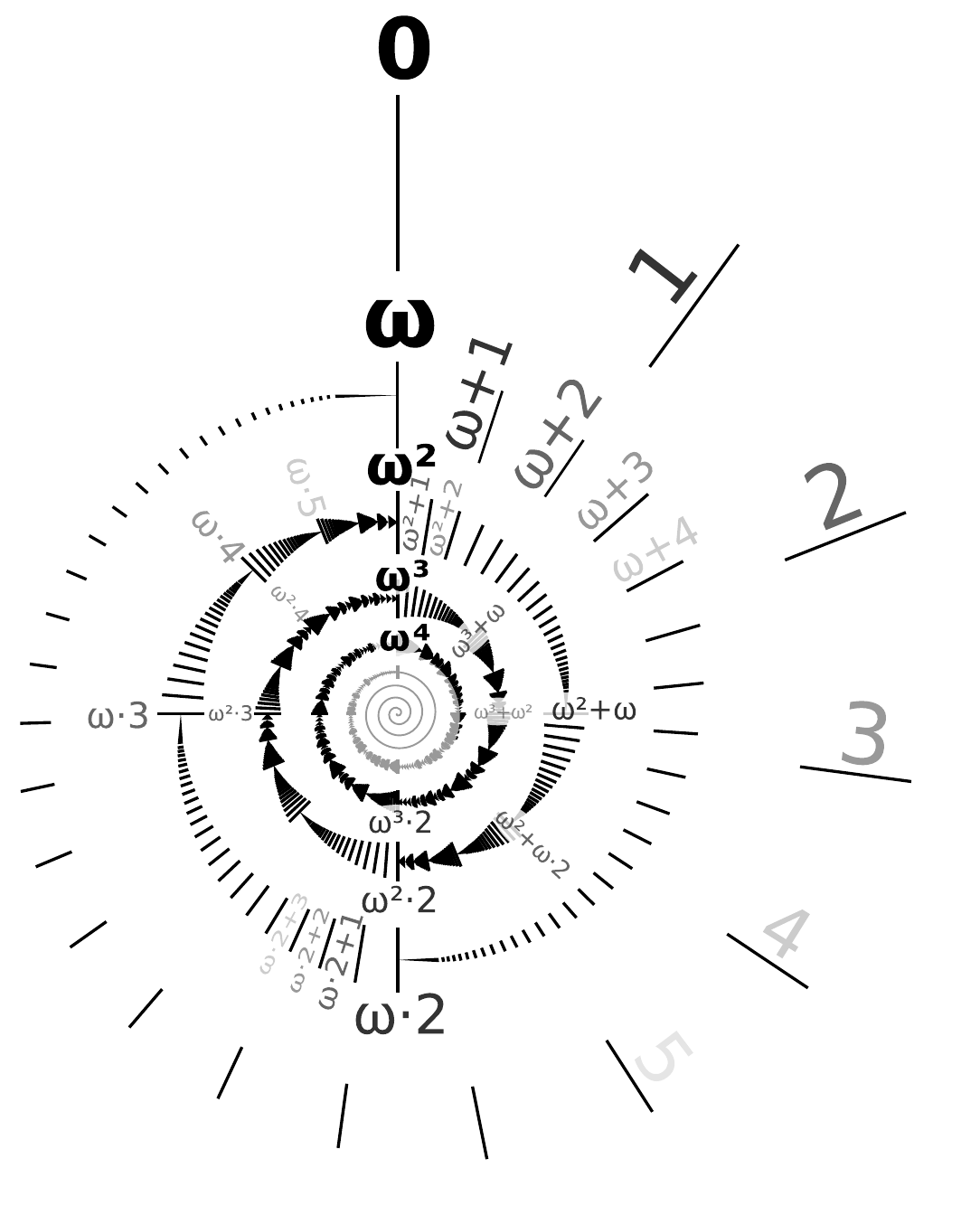} 
\nota{The volumes of all the complete orientable hyperbolic three-manifolds form a well-ordered subset of $\matR$ of type $\omega^\omega$.}
\label{ordinal:fig}
\end{center}
\end{figure}

We deduce that the volumes of all the complete orientable hyperbolic three-manifolds are indexed by countable ordinals. The volumes $v_1< v_2< \ldots, $ of the smallest closed manifolds converge to the volume $v_\omega$ of the smallest 1-cusped manifold; the volumes $v_\omega < v_{2\omega} < \ldots$ of the smallest 1-cusped manifolds converge to the volume $v_{\omega^2}$ of the smallest 2-cusped manifold, and so on. The set of all volumes has order type $\omega^\omega$, see Figure \ref{ordinal:fig}.

We can also deduce the following.

\begin{cor}
For every $V>0$ and $R>0$ there are only finite many closed hyperbolic three-manifolds $M$ with $\Vol(M)<V$ and $\inj(M)>R$.
\end{cor}
\begin{proof}
On a sequence of closed hyperbolic three-manifolds with bounded volume, the length of some core geodesic tends to zero.
\end{proof}

\subsection{References}
To furnish a complete proof of Thurston's Dehn filling Theorem we have patched together a number of different sources. The most relevant ones are as usual Thurston's notes \cite{Th}, Benedetti -- Petronio \cite{BP}, and Ratcliffe \cite{R}, together with Neumann -- Zagier \cite{NZ}. 

Some arguments have been updated here in light of some new results that were discovered more recently by some authors. The list of non-exceptional slopes for the figure-eight knot and the Whitehead link can be found in Martelli -- Petronio \cite{MP}. The smoothness of the deformation space stated in Corollary \ref{smooth:cor} has been proved by Choi \cite{Ch} in 2004. Many of the arguments of this chapter, including all the discussion on angle structures, have been taken from a nice paper of Futer -- Gu\'eritaud \cite{FG}. The first complete proof of Thurson's Dehn filling Theorem, that takes care of the annoying case of triangulations with flat tetrahedra, is due to Petronio -- Porti \cite{PP}. We diverge from that proof in the last arguments, where we extend Choi's smoothness theorem to this partially flat context in Proposition \ref{smooth:flat:prop} (this proposition is probably the only original result of this book). The virtual existence of geometric triangulations was proved by Luo -- Schleimer -- Tillmann \cite{LST} in 2008.